%% file: Manuskript.tex
\documentclass[a4paper,11pt,twoside]{scrartcl}
\usepackage{german,amssymb,amsmath,epsfig,graphics,graphicx,mathrsfs,xcolor,amsthm,cite} 
\usepackage[latin1]{inputenc}
\usepackage{fancyhdr}
\usepackage{makeidx}  

\newtheorem{theorem}{Theorem}[section]

\newtheorem{definition}[theorem]{Definition}
\newtheorem{lemma}[theorem]{Lemma}

\newtheorem{folgerung}[theorem]{Folgerung}
\newtheorem{beispiel}[theorem]{Beispiel}
\newtheorem{satz}[theorem]{Satz}
\newtheorem{bemerkung}[theorem]{Bemerkung}
\newtheorem*{definit}{Definition}

\DeclareMathOperator*{\esssup}{ess\,sup}
\DeclareMathOperator*{\Var}{{\bf V}}

\DeclareOldFontCommand{\it}{\normalfont\itshape}{\mathit}
\DeclareOldFontCommand{\rm}{\normalfont\rmfamily}{\mathrm}
\DeclareOldFontCommand{\bf}{\normalfont\bfseries}{\mathbf}

\numberwithin{equation}{section}
\newcommand{\R}{{\mathbb R}}
\newcommand{\N}{{\mathbb N}}

\newcommand{\Zt}{{\mathscr{Z}^k([t_0,t_1])}}

\makeindex

\begin{document}
\thispagestyle{empty}
\vspace*{10mm}
{\bf{\LARGE Nico Tauchnitz}} \\[2cm]
{\bf{\Huge Optimalitätsprinzipien \\[5mm] in der Steuerungstheorie}} \\[2cm]
{\bf{\Large Zu schwachen und starken lokalen Minima in Steuerungs- \\[5mm]
            problemen mit endlichem und unendlichem Zeithorizont}} \\[10mm]            
\begin{center}\fbox{\includegraphics[width=0.95\textwidth]{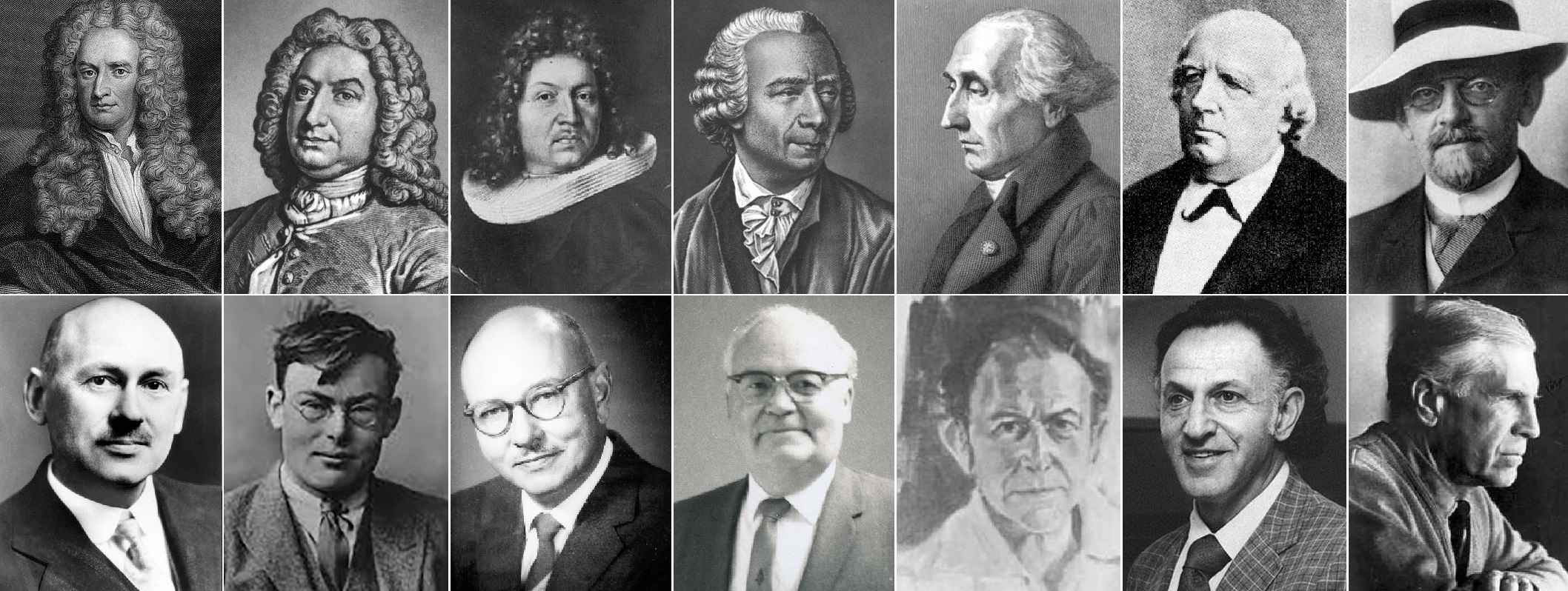}}\end{center}
\vspace*{15mm}

\cleardoublepage
\thispagestyle{empty}

\begin{minipage}{0.2\textwidth}
\includegraphics[width=2cm]{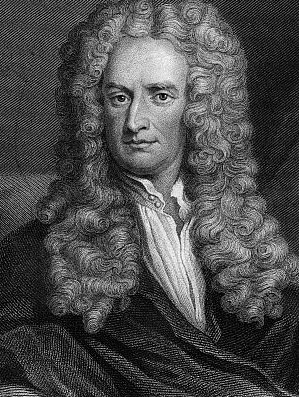}
\end{minipage}
\begin{minipage}{0.79\textwidth}
{\bfseries Isaac Newton} (1643--1727) präsentierte in seinem Jahrtausendwerk {\em Principia Mathematica}
(1687 veröffentlicht) mit der Berechnung eines Rotationskörpers mit geringstem Luftwiderstand
das erste vollständig gelöste Variationsproblem.
\vspace*{8.5mm}
\end{minipage} \\[2mm]
\begin{minipage}{0.2\textwidth}
\includegraphics[width=2cm]{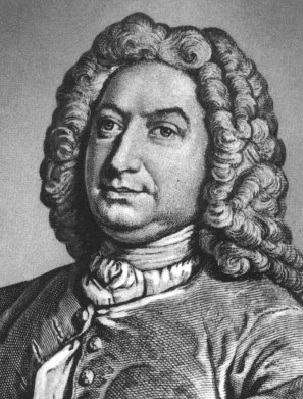}
\end{minipage}
\begin{minipage}{0.79\textwidth}
{\bfseries Johann Bernoulli} (1667--1748) veröffentlichte im Jahr 1696 eine ``Einladung zur Lösung eines neuen Problems''
an die ``scharfsinnigsten Mathematiker der Welt''.
Das gestellte Problem der {\em Brachistochrone} gilt als Geburtsstunde der Variationsrechnung.
\vspace*{8.5mm}
\end{minipage} \\[2mm]
\begin{minipage}{0.2\textwidth}
\includegraphics[width=2cm]{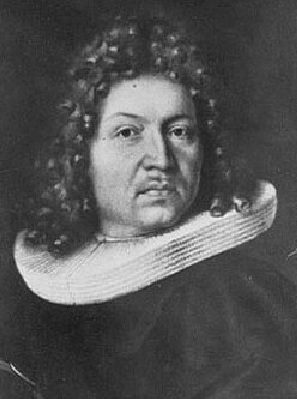}
\end{minipage}
\begin{minipage}{0.79\textwidth}
{\bfseries Jakob Bernoulli} (1655--1705), älterer Bruder von Johann,
schlug einen neuartigen Weg zur Lösung des Problems der Brachistochrone vor.
Dieser bemerkenswerte Ansatz bildete die Grundlage einer Methode,
die Euler später als Variationsrechnung bezeichnete.
\vspace*{8.5mm}
\end{minipage} \\[2mm]
\begin{minipage}{0.2\textwidth}
\includegraphics[width=2cm]{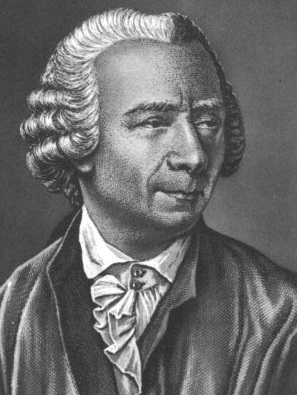}
\end{minipage}
\begin{minipage}{0.79\textwidth}
Die Variationsrechnung verdankt {\bfseries Leonhard Euler} (1707--1783),
dem bedeutendsten Mathematiker des 18.\,Jahrhunderts, nicht nur ihren Namen
({\em Calculus variationum}, 1756), sondern auch ihr erstes Lehrbuch im Jahr 1744.
\vspace*{9mm}
\end{minipage} \\[2mm]
\begin{minipage}{0.2\textwidth}
\includegraphics[width=2cm]{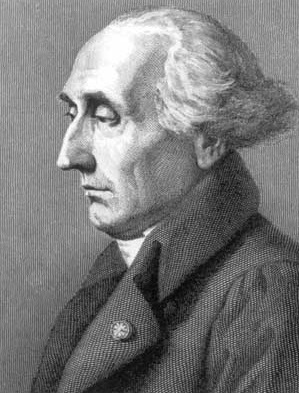}
\end{minipage}
\begin{minipage}{0.79\textwidth}
{\bfseries Joseph-Louis Lagrange} (1736--1813) erweiterte mit Hilfe von ersten Variationen von Funktionalen
({\em Richtungsvariationen}) die analytischen Methoden der
Klassischen Variationsrechnung zur Bestimmung schwacher lokaler Optimalstellen grundlegend.
\vspace*{8mm}
\end{minipage} \\[2mm]
\begin{minipage}{0.2\textwidth}
\includegraphics[width=2cm]{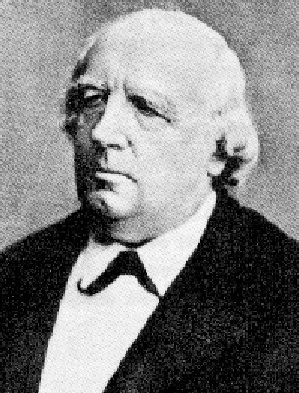}
\end{minipage}
\begin{minipage}{0.79\textwidth}
{\bfseries Karl Weierstraß} (1815--1897) führte die {\em Nadelvariationsmethode} ein,
aus denen neue Optimalitätsbedingungen und der allgemeinere Optimalitätsbgeriff einer
starken lokalen Optimalstelle resultieren.
\vspace*{12mm}
\end{minipage} \\[2mm]
\begin{minipage}{0.2\textwidth}
\includegraphics[width=2cm]{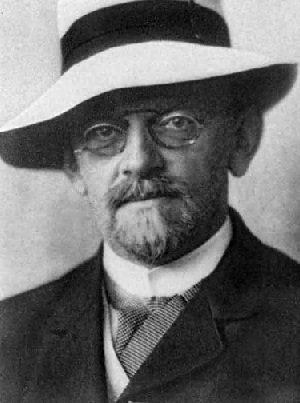}
\end{minipage}
\begin{minipage}{0.79\textwidth}
Im Jahr 1900 stellte {\bfseries David Hilbert} (1862--1943) auf dem Internationalen Mathematiker-Kongress 
eine Liste von 23 Problemen vor, die zu diesem Zeitpunkt ungelöst waren.
Als eine Anwort auf das 23.\,Hilbertsche Problem, der Weiterentwicklung der Variationsrechnung,
gilt die Optimale Steuerung.
\vspace*{4mm}
\end{minipage}

\newpage
\thispagestyle{empty}

\begin{minipage}{0.2\textwidth}
\includegraphics[width=2cm]{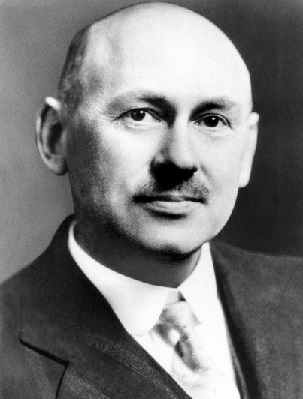}
\end{minipage}
\begin{minipage}{0.79\textwidth}
Die ersten Anwendungen von Steuerungsproblemen entwickelten sich durch die Luft- und Raumfahrt.
Im Bereich der Raumfahrt publizierte der Raketenpionier {\bfseries Robert Goddard} (1882--1945)
zu ``Methoden zum Erreichen extremer Höhen'' bereits 1920 eine Abhandlung.
\vspace*{7.5mm}
\end{minipage} \\[2mm]
\begin{minipage}{0.2\textwidth}
\includegraphics[width=2cm]{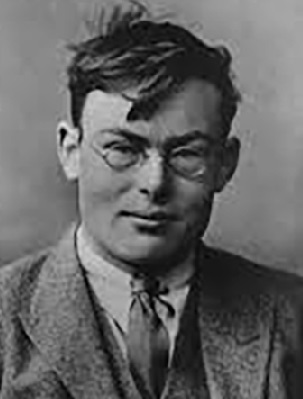}
\end{minipage}
\begin{minipage}{0.79\textwidth}
Mit der Fragestellung nach der ``optimalen Sparquote einer Ökonomie'' schuf
{\bfseries Frank Ramsey} (1903--1930) im Jahr 1928 einen neuen Zweig der neoklassischen Wachstumstheorie.
Diese Betrachtung eines Variationsproblemes mit einem unbeschränkten Planungszeitraum führte zu den
Steuerungsproblemen mit unendlichem Zeithorizont.
\vspace*{3mm}
\end{minipage} \\[2mm]
\begin{minipage}{0.2\textwidth}
\includegraphics[width=2cm]{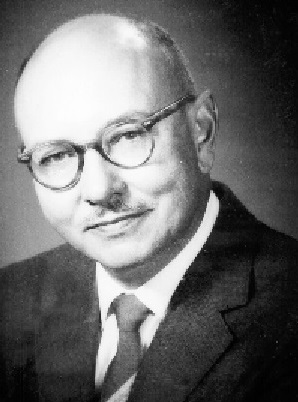}
\end{minipage}
\begin{minipage}{0.79\textwidth}
Im Zuge der Bestrebungen der Chicagoer Schule um Oskar Bolza (1857--1942) und Gilbert Bliss (1876--1951) um einer
Verallgemeinerung der Variationsrechnung
führte {\bfseries Edward McShane} (1904--1989) bereits 1939 diejenige Form der Nadelvariationsmethode ein,
auf der der spätere Beweis des Maximumprinzips um Pontrjagin beruhte.
\vspace*{3mm}
\end{minipage} \\[2mm]
\begin{minipage}{0.2\textwidth}
\includegraphics[width=2cm]{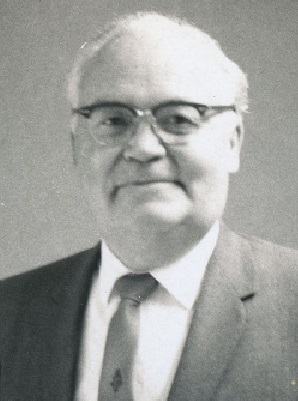}
\end{minipage}
\begin{minipage}{0.79\textwidth}
Im RAND Research Memorandum No. 100 formulierte zu Beginn der 1950er Jahre
{\bfseries Magnus Hestenes} (1906--1991) eine frühere Version des Maximumprinzips. 
Insbesondere gab er dabei den Steuerungsproblemen die heutige Form mit Zustands- und Steuerungsvariable.
\vspace*{8mm}
\end{minipage} \\[2mm]
\begin{minipage}{0.2\textwidth}
\includegraphics[width=2cm]{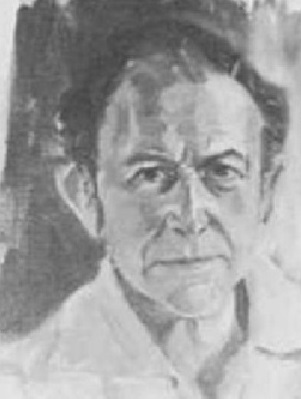}
\end{minipage}
\begin{minipage}{0.79\textwidth}
{\bfseries Rufus Isaacs} (1914--1977) ist in den 1950er Jahren Vorreiter in der Untersuchung von {\em Differentialspielen},
die wegen dem gleichzeitigen Auftreten mehrerer Zielkriterien als Verallgemeinerung der Standardaufgabe der Steuerungstheorie
gelten.
\vspace*{7.5mm}
\end{minipage} \\[2mm]
\begin{minipage}{0.2\textwidth}
\includegraphics[width=2cm]{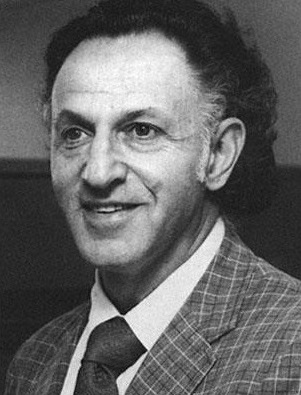}
\end{minipage}
\begin{minipage}{0.79\textwidth}
{\bfseries Richard Bellman} (1920--1984) beschäftigte sich in 1950er Jahren mit Entscheidungsproblemen 
und entwickelte das {\em Prinzip der Dynamischen Programmierung},
das unter zusätzlichen Annahmen in direktem Bezug zum Pontrjaginschen Maximumprinzip steht.
\vspace*{6.5mm}
\end{minipage} \\[2mm]
\begin{minipage}{0.2\textwidth}
\includegraphics[width=2cm]{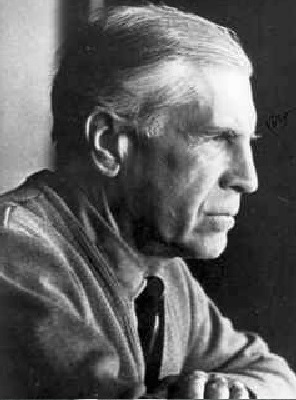}
\end{minipage}
\begin{minipage}{0.79\textwidth}
Das {\em Pontrjaginsche Maximumprinzip},
ein kompletter Satz notwendiger Optimalitätsbedingungen für ein starkes lokales Minimum in Steuerungsproblemen,
formulierten {\bfseries Lew Pontrjagin} (1908--1988), Rewas Gamkrelidze (1927--) und 
Wladimir Boltjanski (1925--2019) zunächst in den 1950er Jahren als Hypothese und veröffentlichten den Nachweis 1961.
\end{minipage}

\cleardoublepage


\newpage
\setcounter{page}{1}
\pagenumbering{roman}

\lhead[\thepage \hspace*{1mm} Inhaltsverzeichnis]{}
\rhead[]{Inhaltsverzeichnis \hspace*{1mm} \thepage}
\tableofcontents
\cleardoublepage




\newpage
\setcounter{page}{1}
\pagenumbering{arabic}
\addcontentsline{toc}{section}{Einleitung}
\lhead[\thepage \hspace*{1mm} Einleitung]{}
\rhead[]{Einleitung \hspace*{1mm} \thepage}
\input{0-Einleitung}
\cleardoublepage

\newpage
\addcontentsline{toc}{section}{Wichtige Bezeichnungen}
\lhead[\thepage \hspace*{1mm} Wichtige Bezeichnungen]{}
\rhead[]{Wichtige Bezeichnungen \hspace*{1mm} \thepage}
\input{0-Bezeichnungen}
\cleardoublepage

\newpage
\addcontentsline{toc}{section}{Stückweise Stetigkeit und Differenzierbarkeit}
\lhead[\thepage \hspace*{1mm} Stückweise Stetigkeit]{}
\rhead[]{Stückweise Stetigkeit \hspace*{1mm} \thepage}
\input{0-Stueckweise}
\cleardoublepage

\newpage
\lhead[\thepage \hspace*{1mm} Richtungsvariationen]{}
\rhead[]{ }
\input{1-0-KlassischeVariationsrechnung}

       \rhead[]{Euler-Lagrangesche Gleichung \hspace*{1mm} \thepage}
       \input{1-1-EulerLagrange}
       \input{1-11-Beispiele}
       \input{1-12-BolzaAufgabe}

       \rhead[]{Lagrange-Aufgabe \hspace*{1mm} \thepage}
       \input{1-2-LagrangeAufgabe}
       \input{1-21-Beweis}
       \input{1-22-Isoperimetrisch}

       \rhead[]{Schwaches Optimalit\"atsprinzip \hspace*{1mm} \thepage}
       \input{1-3-Richtungsvariationen}
       \input{1-31-Aufgabenstellung}
       \input{1-32-SchwachesPrinzip}
       \input{1-33-Beweis}
       \input{1-34-Hinreichend}
       \input{1-35-FreieZeit}
       \input{1-36-Zustandsaufgabe}
       \newpage
       \input{1-37-Beweis}

\cleardoublepage


\newpage
\lhead[\thepage \hspace*{1mm} Standardaufgabe der Optimalen Steuerung]{}
\rhead[]{ }
\input{2-0-Nadelvariationen}
      
       \newpage
       \rhead[]{Klassische Nadelvariation \hspace*{1mm} \thepage}
       \input{2-1-Nadelvariationen}
       
       \newpage
       \rhead[]{Elementare Aufgabe \hspace*{1mm} \thepage}
       \input{2-2-PMPeinfach}
       \input{2-21-Beweis}
       \newpage
       \input{2-22-Interpretation}
       
       \newpage
       \rhead[]{Aufgabenstellung \hspace*{1mm} \thepage}
       \input{2-3-Aufgabenstellung}
      
       \newpage
       \rhead[]{Pontrjaginsches Maximumprinzip \hspace*{1mm} \thepage}
       \input{2-4-PontrjaginAufgabe}
       \input{2-41-Beweis}
       \input{2-42-Hinreichend}   
      
       \rhead[]{Zustandsbeschr\"ankungen \hspace*{1mm} \thepage}
       \input{2-5-Zustandsaufgabe}
       \input{2-51-Beweis}
       \input{2-52-HinreichendZB}
       \newpage
       \input{2-53-FreieZeit}    
\cleardoublepage
\addtocontents{toc}{\protect\newpage} 


\newpage
\lhead[\thepage \hspace*{1mm} Erweiterungen der Standardaufgabe]{}
\rhead[]{ }
\input{3-0-Erweiterungen}
       \newpage
       \rhead[]{Zielfunktionale in gemischter Form \hspace*{1mm} \thepage}
       \input{3-10-GemischtesZF}
       \input{3-11-PMP}
       \input{3-12-Produktion}
       \input{3-13-Instandhaltung}
       
       \newpage
       \rhead[]{Zeitschranken \hspace*{1mm} \thepage}
       \input{3-20-Zeitfenster}

       \newpage
       \rhead[]{Multiprozesse \hspace*{1mm} \thepage}
       \input{3-30-Multiprozesse}
       \input{3-31-Zerlegungen}
       \input{3-32-PMPeinfach}
       \input{3-33-Maximumprinzip}
       \input{3-34-Beweis}
       \input{3-35-FreieZeit}
       \input{3-36-Investitionsmodell}
       \input{3-37-Wassercontainer}
       
       \newpage
       \rhead[]{Zeitverz\"ogerte Systeme \hspace*{1mm} \thepage} 
       \input{3-40-RetardierteSysteme}
       \input{3-41-PMPeinfach}
       \input{3-42-Aufgabenstellung}
       \input{3-43-Maximumprinzip}
       \input{3-44-Beweis}
       \input{3-45-Arrow}
       \newpage       
       \input{3-46-Chemoimmuntherapie}         
       
       \newpage
       \rhead[]{Differentialspiele \hspace*{1mm} \thepage} 
       \input{3-50-Differentialspiele}
       \input{3-51-Kapitalismusspiel}
       \input{3-52-Fischerei}
\cleardoublepage


\newpage
\lhead[\thepage \hspace*{1mm} Unendlicher Zeithorizont]{}
\rhead[]{ }
\input{4-0-StarkeVariationen}
       \newpage
       \rhead[]{Elementare Aufgabe \hspace*{1mm} \thepage}
       \input{4-1-EinfacheNV}
       
       \newpage
       \rhead[]{Aufgabenstellung \hspace*{1mm} \thepage}
       \input{4-2-Aufgabenstellung}

       \rhead[]{Pontrjaginsches Maximumprinzip \hspace*{1mm} \thepage}
       \input{4-3-PontrjaginAufgabe}
       \input{4-31-Beweis}
       \input{4-32-Normalform}
       \input{4-33-Hinreichend}  

       \rhead[]{Zustandsbeschr\"ankungen \hspace*{1mm} \thepage}
       \input{4-4-Zustandsaufgabe}
       \input{4-41-Beweis}
       \input{4-42-Hinreichend} 
       \newpage
       \rhead[]{Einordnung der Aufgabe \hspace*{1mm} \thepage}
       \input{4-5-Einordnung}
       \input{4-51-Approximation}
       \input{4-52-Approximationssatz}
       \input{4-53-Zeittransformation}
       \input{4-54-Einordnung}

\cleardoublepage
\addtocontents{toc}{\protect\newpage} 


\newpage
\lhead[\thepage \hspace*{1mm} Volterrasche Integralgleichungen]{}
\rhead[]{ }
\input{5-0-Integralgleichungen}\newpage
       \rhead[]{Elementare Aufgabe \hspace*{1mm} \thepage}
       \input{5-1-PMPeinfach}       
       
       \newpage
       \rhead[]{Aufgabenstellung \hspace*{1mm} \thepage}
       \input{5-2-Aufgabenstellung} 
       
       \newpage
       \rhead[]{Pontrjaginsches Maximumprinzip \hspace*{1mm} \thepage}
       \input{5-3-PontrjaginAufgabe}
       \input{5-31-Beweis}
       \input{5-32-Hinreichend}
       \input{5-33-Werbestrategien}  
      
       \newpage
       \rhead[]{Zustandsbeschr\"ankungen \hspace*{1mm} \thepage}
       \input{5-4-Zustandsaufgabe}
       \input{5-41-Beweis}
       \input{5-42-HinreichendZB}
       \input{5-43-Werbestrategien}
       
       \newpage
       \rhead[]{Freier Anfangs- und Endzeitpunkt \hspace*{1mm} \thepage}
       \input{5-5-FreieZeit}
       \input{5-51-Substitution}
       \input{5-52-Beweis}
       \input{5-53-Maximumprinzip}
       \input{5-54-Werbestrategien}
       \input{5-55-Standardaufgabe}

\cleardoublepage

\newpage
\addcontentsline{toc}{part}{Anhang}
\begin{appendix}
\lhead[\thepage \hspace*{1mm} Maß- und Integrationstheorie]{}
\rhead[]{Maß- und Integrationstheorie \hspace*{1mm} \thepage}
\input{I-Masstheorie}
\newpage
\lhead[\thepage \hspace*{1mm} Funktionalanalytische Hilfsmittel]{}
\rhead[]{Funktionalanalytische Hilfsmittel \hspace*{1mm} \thepage}
\input{II-Hilfsmittel}
\newpage
\lhead[\thepage \hspace*{1mm} Differential- und Integralgleichungen]{}
\rhead[]{Differential- und Integralgleichungen \hspace*{1mm} \thepage}
\input{III-Differentialgleichung}
\newpage
\addtocontents{toc}{\protect\newpage} 
\lhead[\thepage \hspace*{1mm} Konvexe Analysis]{}
\rhead[]{Konvexe Analysis \hspace*{1mm} \thepage}
\input{IV-KonvexeAnalysis}
\newpage
\lhead[\thepage \hspace*{1mm} Mehrfache Nadelvariationen]{}
\rhead[]{Mehrfache Nadelvariationen \hspace*{1mm} \thepage}
\input{V-Nadelvariation}
      \input{V-1-Konstruktion}
      \input{V-2-Definition}
      \input{V-3-Eigenschaften}
      \input{V-4-Multiprozesse}
      \input{V-5-UnendlicherHorizont}
      \input{V-6-Integralgleichung}
\newpage
\lhead[\thepage \hspace*{1mm} Theorie der Extremalaufgaben]{}
\rhead[]{Theorie der Extremalaufgaben \hspace*{1mm} \thepage}
\input{VI-Extremalprinzip}
      \input{VI-1-GlattesExtremalprinzip}
      \input{VI-2-SchwachesExtremalprinzip}
      \input{VI-3-StarkesExtremalprinzip}
       
\cleardoublepage
\end{appendix}
\cleardoublepage


\newpage
\addcontentsline{toc}{part}{Abbildungs-, Literatur- und Sachverzeichnis}
\addcontentsline{toc}{section}{Abbildungsverzeichnis}
\lhead[\thepage \hspace*{1mm} Abbildungsverzeichnis]{}
\rhead[]{Abbildungsverzeichnis \hspace*{1mm} \thepage}
\listoffigures
\cleardoublepage

\newpage
\addcontentsline{toc}{section}{Literatur}
\lhead[\thepage \hspace*{1mm}  Literatur]{}
\rhead[]{Literatur \hspace*{1mm} \thepage}
\input{Literatur}
\cleardoublepage




\newpage
\addcontentsline{toc}{section}{Sachverzeichnis}
\printindex
\end{document}

%% file: 0-Einleitung.tex
\section*{Einleitung}
Die Theorie der Optimalen Steuerung entwickelte sich aus der Klassischen Variationsrechnung.
Notwendig wurde die Weiterentwicklung der Methoden der Klassischen Variationsrechnung durch Problemstellungen aus verschiedenen 
Fachbereichen, z.\,B. der Physik, der Biologie, der Chemie, den Ingenieurwissenschaften und vorallem
der Luft- und Raumfahrt.
Die Aufgabenstellungen haben gemein,
dass sie steuerbar, d.\,h. durch den Menschen beeinflussbar, sind.
Es ergibt sich damit die Frage nach derjenigen Steuerung,
die im gew"unschten Sinn die beste, wir sagen die optimale, ist.
Diese optimale Steuerung zeigt dem Anwender Handlungsweisen
auf, wie auf den Zustand des steuerbaren Objektes einzuwirken ist, damit sich
dieser entsprechend der Zielvorstellungen auf die beste Art und Weise verh"alt. \\[2mm]
Historisch gilt als Geburtsstunde der Klassischen Variationsrechnung das Jahr 1696,
als Johann Bernoulli (1667--1748) seine Zeitgenossen mit der Aufgabe der Brachistochrone\index{Brachistochrone} herausforderte.
Aber es gibt weitaus "altere Variationsprobleme. \\
Vergil (70--19 v.\,Chr.) berichtet in seinem antiken Epos der Aeneis die Sage von der Gr"undung Karthagos etwa
900 v.\,Chr. durch die ph"onizische Prinzessin Dido.
Aus Furcht vor ihrem Bruder Pygmalion, der ihren Gemahl erschlagen hatte,
floh Dido und gelangte an die K"uste Nordafrikas. 
Man wollte ihr dort jedoch nur soviel Land "uberlassen,
wie sie mit einer Ochsenhaut begrenzen k"onne.
Dido schnitt die Haut in feine Streifen, n"ahte diese zusammen und konnte damit ein
nicht unerhebliches Landst"uck an der K"uste eingrenzen,
auf dem sie die Burg Byrsa errichten lie"s. \\
Das Problem der Dido\index{Problem der Dido} wirft die Frage nach der Gestalt derjenigen geschlossenen,
ebenen Kurve gegebener L"ange auf,
die den Bereich mit gr"o"stem Fl"acheninhalt berandet. \\
Der (orientierte) Inhalt $I$ des Bereiches $B$, der von einer geschlossenen, ebenen und stetig differenzierbaren Kurve
mit der Parameterdarstellung
$$t \to x(t)=\big(x_1(t),x_2(t)\big), \qquad t \in [a,b],$$
umrandet wird, kann durch
$$I=\frac{1}{2} \int_a^b [ x_1(t) \dot{x}_2(t)-x_2(t) \dot{x}_1(t) ] \, dt$$
bestimmt werden (das Vorzeichen der Zahl $I$ h"angt vom Durchlaufsinn der Kurve ab).
Dabei muss die Kurve die isoperimetrische Beschr"ankung der vorgegebenen L"ange $l$,
$$\int_a^b \sqrt{\dot{x}_1^2(t)+\dot{x}_2^2(t)} \, dt =l,$$
erf"ullen.
Zudem betrachten wir ausschlie"slich geschlossene Kurven,
was wir durch die Angaben $x_1(a)=x_1(b)$ und $x_2(a)=x_2(b)$ erreichen.
Wir beschreiben der Form halber die isoperimetrische Beschr"ankung durch die Funktion $x_3(\cdot)$ mit den Eigenschaften
$$\dot{x}_3(t)=\sqrt{\dot{x}_1^2(t)+\dot{x}_2^2(t)}, \qquad x_3(a)=0, \quad x_3(b)=l.$$
Dann f"uhrt das Problem der Dido auf der Menge der geschlossenen und glatten Kurven gegebener L"ange zu der Maximierungsaufgabe
$$J\big(x_1(\cdot),x_2(\cdot),x_3(\cdot)\big) =  \frac{1}{2}\int_a^b [ x_1(t) \dot{x}_2(t)-x_2(t) \dot{x}_1(t) ] \, dt \to \sup$$
unter den Nebenbedingungen
$$\dot{x}_3(t)=\sqrt{\dot{x}_1^2(t)+\dot{x}_2^2(t)}, \qquad x_1(a)=x_1(b), \quad x_2(a)=x_2(b), \quad x_3(a)=0, \; x_3(b)=l.$$
Die L"osung des Problems der Dido wird durch die isoperimetrische Eigenschaft des Kreises geliefert:
Unter allen m"oglichen geschlossenen Kurven gegebener L"ange ist diejenige,
deren Innengebiet den gr"o"stm"oglichen Fl"acheninhalt besitzt, die Kreislinie. \\
Das Problem der Dido zeigt, dass sich Variationsaufgaben grundlegend zu der Extremwertsuche einer Funktion unterscheiden.
Während man bei Extremwertaufgaben einer Funktion diejenige Stelle sucht,
an der die Funktion einen kleinsten oder gr"o"sten Wert annimmt,
ist man im Gegensatz dazu im Problem der Dido bestrebt eine Funktion zu finden,
die dem Zielfunktional $J$ den optimalen Wert zuordnet.
Die Maximierung des Funktionals $J$ dr"ucken wir durch die Suche nach einem Supremum aus,
da wir "uber die Existenz eines Maximum a priori keine Kenntnis besitzen. \\[2mm]
Der\index{Goddard-Problem} Raketenpionier Robert H. Goddard (1882--1945) formulierte in der ersten
H"alfte des 20. Jahrhunderts die Aufgabe,
dass eine H"ohenrakete mit einer vorgegebenen Treibstoffmenge eine m"oglichst gro"se Flugh"ohe erreichen soll.
Das besondere Merkmal daran ist,
dass die Schubleistung gewissen Schranken unterliegt und Grenzlagen annehmen darf.
Diese besondere Charakteristik ist in der Klassischen Variationsrechnung nicht vorgesehen.
Damit wird die Frage nach einer passenden Verallgemeinerung der Klassischen Variationsrechnung aufgeworfen,
die zur Antwort die Entwicklung der Steuerungstheorie erhielt. \\
In einer ``einfachen'' Formulierung einer Rakete maximaler Flugh"ohe (Maurer \cite{Maurer}) seien 
\begin{center}\begin{tabular}{lcl}
$h(t)$ &--& die H"ohe der Rakete zur Zeit $t$, \\
$v(t)$ &--& die Geschwindigkeit der Rakete zur Zeit $t$, \\
$m(t)$ &--& die Masse der Rakete zur Zeit $t$, \\
$u(t)$ &--& Kraftstoffverbrauch zur Zeit $t$, \\
$D(h,v)$ &--& der Luftwiderstand, \\
$g(h)$ &--& die Gravitation, \\
$c$ &--& der spezifische Impuls pro Einheit Treibstoff. \hspace*{3cm}
\end{tabular} \end{center}
Nach dem Luftwiderstands- und dem Gravitationsgesetz gelten
$$D(h,v) = \alpha v^2 e^{-\beta h}, \qquad g(h)=g_0 \frac{r_0^2}{r_0^2+h^2},$$
wobei $\alpha,\beta$ spezifische Konstanten des Modells,
$r_0$ den Erdradius und $g_0$ die Gravitationskonstante auf der Erdoberfl"ache bezeichnen. \\
Die Dynamik des Systems wird durch die Differentialgleichungen
$$\dot{h}(t) = v(t), \qquad \dot{v}(t) = \frac{c u(t)-D\big(h(t),v(t)\big)}{m(t)} - g\big(h(t)\big), \qquad \dot{m}(t) = -u(t)$$
gegeben. Die Randbedingungen lauten
$$h(0)=0, \qquad v(0)=0,\qquad m(0)=m_0\mbox{ (Startmasse)}, \qquad m(T)=m_T \mbox{ (Leermasse)}$$
und es treten Steuerungsbeschr"ankungen der Form
$$u(t) \in [0,u_{\max}]$$
auf. Schlie"slich ist das Optimierungskriterium die Maximierung der Steigh"ohe:
$$h(T)= \int_0^T v(t)\, dt \to \sup.$$
Als naheliegende L"osung k"onnte man die Steuerung
$$u(t)=\left\{ \begin{array}{ll} u_{\max}, & t \in [0,\tau), \\[1mm] 0, & t \in [\tau,T], \end{array} \right.$$
erwarten, die im Goddard-Problem ohne Luftwiderstand optimal w"are.
Jedoch wird die dadurch erreichte Flugh"ohe um fast 40\% "ubertroffen,
\begin{figure}[h]
	\centering
	\fbox{\includegraphics[width=6cm]{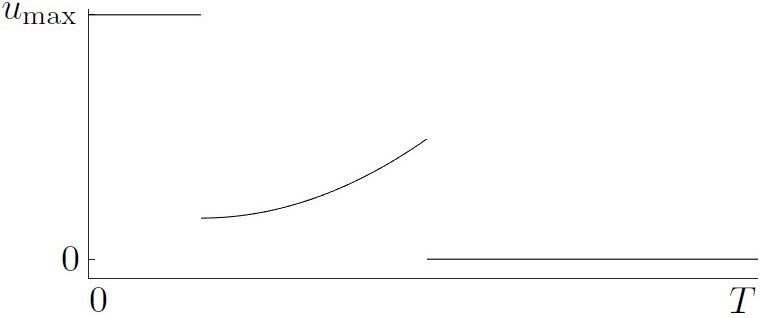}}
	\caption[Optimale Steuerung im Goddard Problem mit Luftwiderstand]{Optimale Steuerung im Goddard Problem mit Luftwiderstand.}
	\label{AbbGoddard}
\end{figure}
wenn man eine Steuerung anwendet, deren charakteristischer Verlauf in Abbildung \ref{AbbGoddard} dargestellt ist. \\
Weitere Einsatzm"oglichkeiten der mathematischen Methoden f"ur die Optimierung einer H"ohenrakete ergeben sich,
wenn gegens"atzliche Ziele gleichzeitig beachtet werden m"ussen.
Gewicht und Sicherheit sind Beispiele daf"ur.
Einerseits sollte das Gesamtgewicht so gering wie m"oglich sein,
um weniger Treibstoff zu verbrauchen.
Andererseits ist die Explosionsgefahr bei dickeren W"anden der Brennstoffkammern geringer,
dabei aber das Gewicht der Rakete h"oher.
Die optimale Abstimmung dieser gegenl"aufigen Bestrebungen ist gesucht.

\newpage
Wir geben nun eine allgemeine Gestalt f"ur die Variations- und Steuerungsprobleme mit endlichem Zeithorizont an,
die wir in dieser Arbeit unter verschiedenen Voraussetzungen und Rahmenbedingungen betrachten werden: \\
Ohne Einschr"ankung sei das Problem stets als Minimierungsaufgabe gestellt.
Es bezeichne ferner $x(\cdot):\R \to \R^n$ den Zustand des steuerbaren Objektes und
$u(\cdot):\R \to \R^m$ die Steuerung,
mit der auf den Zustand Einfluss genommen werden kann.
Die Einschr"ankung, dass die Steuerungsparameter $u(t)$ nur Werte aus einer vorgegebenen Menge $U$ (dem Steuerungsbereich)
annehmen d"urfen,
nennen wir Steuerungsbeschr"ankungen. \\
Das Optimierungskriterium wird durch das Zielfunktional
$$J\big(x(\cdot),u(\cdot)\big) = \displaystyle \int_{t_0}^{t_1} f\big(t,x(t),u(t)\big) \, dt$$
mit dem Integranden $f:\R \times \R^n \times \R^m \to \R$ beschrieben.
Dabei bezeichnet die Variable $t$ oft die Zeit.
Die Dynamik, d.\,h. das dynamische Verhalten des Zustandes $x(\cdot)$ "uber dem Zeitraum $[t_0,t_1]$,
ist durch das Differentialgleichungssystem
$$\dot{x}(t)= \varphi\big(t,x(t),u(t)\big)$$
mit der rechten Seite $\varphi:\R \times \R^n \times \R^m \to \R^n$ beschrieben.
Schlie"slich liegen durch
$$h_0\big(x(t_0)\big)=0, \qquad h_1\big(x(t_1)\big)=0$$
Randbedingungen an die Start- und Zielwerte f"ur den Zustand $x(\cdot)$ im Anfangs- und Endpunkt
des Betrachtungszeitraumes $[t_0,t_1]$ vor.
Zus"atzliche Beschr"ankungen des Zustandes $x(\cdot)$ der Form
$$g_j\big(t,x(t)\big) \leq 0 \quad \mbox{ f"ur alle } t \in [t_0,t_1], \quad j=1,...,l,$$
nennen wir Zustandsbeschr"ankungen. \\
Zusammenfassend ergibt sich folgende Aufgabenklasse,
in die sich das Problem der Dido und das Goddard-Problem einordnen lassen:
\begin{equation}\tag{1} \label{Einleitung2}
\left. \begin{array}{l}
J\big(x(\cdot),u(\cdot)\big) = \displaystyle \int_{t_0}^{t_1} f\big(t,x(t),u(t)\big) \, dt \to \inf, \\[2mm]
\dot{x}(t)= \varphi\big(t,x(t),u(t)\big), \quad t \in [t_0,t_1], \\[2mm]
h_0\big(x(t_0)\big)=0, \qquad h_1\big(x(t_1)\big)=0, \\[2mm]
u(t) \in U, \qquad U \not=\emptyset, \\[2mm]
g_j\big(t,x(t)\big) \leq 0 \quad \mbox{ f"ur alle } t \in [t_0,t_1], \quad j=1,...,l.
\end{array} \right\}
\end{equation}
Dabei sind die Zeitpunkte $t_0<t_1$ fest vorgegeben.
Wird das Problem (\ref{Einleitung2}) mit freiem Anfangs- und Endzeitpunkt $t_0,t_1$ betrachtet, so lauten die Randbedingungen
$$h_0\big(t_0,x(t_0)\big)=0, \qquad h_1\big(t_1,x(t_1)\big)=0.$$
Im Jahr 1928 stellte Frank P. Ramsey (1903--1930) die Frage nach der optimalen Sparquote \cite{Ramsey},
die einer "Okonomie langfristiges und wohlfahrtsoptimiertes Wachstum garantiert.
Das besondere an der Modellierung des Variationsproblems war die Einf"uhrung des unendlichen Zeithorizontes.
Die Philosophie dahinter ist die Vorstellung,
dass kein nat"urliches Ende f"ur den Betrachtungszeitraum existiert.
M"ochte man s"amtlichen nachfolgenden Generationen Beachtung schenken, dann ist die Idealisierung in Form des
zeitlich unbeschr"ankten Rahmens die einzig m"ogliche Konsequenz. \\
In dieser Arbeit enthält die Aufgabe mit unendlichem Zeithorizont
im Zielfunktional eine nichtnegative und integrable Dichtefunktion $\omega(\cdot)$,
welche häufig die Diskontierung $\omega(t)=e^{-\varrho t}$ beschreibt.
Damit besitzt die Aufgabe die Gestalt:
\begin{equation}\tag{2} \label{Einleitung3}
\left. \begin{array}{l}
J\big(x(\cdot),u(\cdot)\big) = \displaystyle \int_0^\infty \omega(t) f\big(t,x(t),u(t)\big) \, dt \to \inf, \\[2mm]
\dot{x}(t)= \varphi\big(t,x(t),u(t)\big), \quad t \in [0,\infty), \\[2mm]
h_0\big(x(0)\big)=0, \qquad \lim\limits_{t \to \infty} h_1\big(t,x(t)\big)=0, \\[2mm]
u(t) \in U, \qquad U \not=\emptyset, \\[2mm]
g_j\big(t,x(t)\big) \leq 0 \quad \mbox{ f"ur alle } t \in \R_+, \quad j=1,...,l.
\end{array} \right\}
\end{equation}
Viele Phänomene in Natur, Technik oder Ökonomie lassen sich durch eine Differentialgleichung nicht adäquat beschreiben.
Stattdessen treten Integralgleichungen in Erscheinung.
Neben der Standardaufgabe (\ref{Einleitung2}) und der Aufgabe mit unendlichem Zeithorizont (\ref{Einleitung3})
gehen wir auf die Steuerung Volterrascher Integralgleichungen ausführlich ein.
Diese Aufgabenklasse besitzt im vorliegenden Manuskript die folgende Gestalt: 
\begin{equation}\tag{3} \label{Einleitung4}
\left. \begin{array}{l}
J\big(x(\cdot),u(\cdot)\big) = \displaystyle \int_{t_0}^{t_1} f\big(t,x(t),u(t)\big) \, dt \to \inf, \\[3mm]
\displaystyle x(t)= x(t_0)+ \int_{t_0}^t \varphi\big(t,s,x(s),u(s)\big), \quad t \in [t_0,t_1], \\[3mm]
h_0\big(x(t_0)\big)=0, \qquad h_1\big(x(t_1)\big)=0, \\[2mm]
u(t) \in U, \qquad U \not=\emptyset, \\[2mm]
g_j\big(t,x(t)\big) \leq 0 \quad \mbox{ f"ur alle } t \in [t_0,t_1], \quad j=1,...,l.
\end{array} \right\}
\end{equation}
Im Rahmen dieser Arbeit befassen wir uns vorrangig mit notwendigen Bedingungen f"ur die Standardaufgabe (\ref{Einleitung2}).
Ausgehend von der Klassischen Variationsrechnung werden wir die Standardaufgabe sukzessive zu einem Steuerungsproblem erweitern und
Verallgemeinerungen wie Multiprozesse oder zeitverzögerte Systeme behandeln. \\
Zudem bilden im Hinblick auf Methoden und Resultate die Aufgabe mit unendlichem Zeithorizont (\ref{Einleitung3}) und
die Steuerung Volterrascher Integralgleichungen (\ref{Einleitung4}) zwei grundlegende Erweiterungen der Standardaufgabe (\ref{Einleitung2}).
Deshalb sind sie in dieser Einleitung besonders hervorgehoben und werden außerdem in eigenständigen Kapiteln behandelt.

\newpage
Die Literatur zur Theorie der Optimalen Steuerung und zur Herleitung notwendiger Optimalit"atsbedingungen ist diesbez"uglich 
sehr umfangreich.
Unser Verzeichnis spiegelt nur diejenigen wider,
die die vorliegende Monographie am st"arksten beeinflusst haben.
Eine vollst"andige Angabe der Quellen und die W"urdigung der einzelnen Beitr"age,
die ma"sgeblich zur Entwicklung der Optimalen Steuerung beigetragen haben,
ist im Rahmen der vorliegenden Arbeit nicht m"oglich und ist auch nicht angestrebt. \\
Hervorheben m"ochten wir zun"achst den fundamentalen Beitrag von Pontrjagin et al. \cite{Pontrjagin}.
Wesentlichen Einfluss auf die Entwicklung der Optimalen Steuerung haben au"serdem u.\,a.
Bellman \cite{Bellman}, Isaacs \cite{Isaacs}, Gamkrelidze \cite{Gamkrelidze}, McShane \cite{McShane} und Young \cite{Young} 
genommen.
Weiterhin beziehen sich unsere Darstellungen zur Theorie der Extremalaufgaben auf die Arbeiten von
Dubovickii \& Milyutin \cite{DuboMil}, Girsanov \cite{Girsanov}, Halkin \cite{Halkin2},
Hestenes \cite{Hestenes}, Ioffe \& Tichomirov \cite{Ioffe}, Kurcyusz \cite{Kurcyusz} und Neustadt \cite{Neustadt}. \\
Aus dem umfassenden Vorrat an Lehrb"uchern und "Ubersichtsartikeln z"ahlen wir 
Boltjanski \cite{Boltjanski}, Carlson et al. \cite{CarHauLei}, Hartl et al. \cite{Hartl},
Sethi \& Thompson \cite{Sethi} und Vinter \cite{Vinter} auf.
Eine Vielzahl an Beispielen sind in den Lehrb"uchern Feichtinger \& Hartl \cite{Feichtinger},
Grass et al. \cite{Grass}, Kamien \& Schwartz \cite{Kamien} und Seierstad \& Syds\ae ter \cite{Seierstad} angegeben.
Abschlie"send möchten wir auf die Arbeiten Pesch \& Plail \cite{PeschPlail} und Plail \cite{Plail} zur historischen Entwicklung
der Optimalen Steuerung zu einem eigenst"andigen mathematischen Fachgebiet hinweisen.
Der Einfluss des Ramsey-Modells auf die Ökonomische Wachstumstheorie ist in den Lehrbüchern von Arnold \cite{Arnold} und
Barro \& Sala-i-Martin \cite{Barro} ausführlich dargestellt. \\[2mm]
Im Aufbau der Arbeit mussten zwei wesentliche Kriterien ber"ucksichtigt werden:
die "ubersichtliche Darstellung und die mathematische Beweisf"uhrung.
Zur besseren Lesbarkeit sind in den Kapiteln die einzelnen Aufgabenklassen und Resultate systematisch aufgef"uhrt,
es werden einfachere Beweise vorgestellt und es werden die Resultate an Beispielen demonstriert.
Demgegenüber sind die Beweise der notwendigen Optimalit"atsprinzipien sehr anspruchsvoll und "au"serst aufwendig.
Aus diesem Grund sind diese Beweise in den Kapiteln lediglich mit Verweis auf die entsprechenden Stellen im Anhang
skizziert.
Der Anhang f"ur sich stellt in seinem Umfang einen eigenen Teil dieses Manuskriptes dar. \\[2mm]
Im ersten Kapitel betrachten wir die klassischen Richtungsvariationen zur Herleitung der Euler-Lagrangeschen Gleichung.
F"ur die Standardaufgabe (\ref{Einleitung2}) als Lagrange-Aufgabe der Klassischen Variationsrechnung 
bilden diese Variationen die Grundlage zur Auswertung notwendiger Optimalit"atsbedingungen.
Der wesentliche Punkt in der Auswertung der Lagrange-Aufgabe ist dabei das Verst"andnis des Problems als eine Extremalaufgabe 
in Funktionenr"aumen und der Nachweis der G"ultigkeit des Lagrangeschen Prinzips in diesem abstrakten Rahmen.
Mit der Erweiterung der Standardaufgabe um Steuerungsbeschr"ankungen verlassen wir die Klassische Variationsrechnung.
Wiederum auf der Basis von Richtungsvariationen erhalten wir in der Standardaufgabe notwendige Bedingungen in Form eines
Schwachen Optimalit"atsprinzips. \\[2mm]
Im zweiten Kapitel steht die Nadelvariationsmethode im Fokus.
Der Unterschied zu den Richtungsvariationen besteht darin,
dass sich ein Kandidat innerhalb einer gr"o"seren Konkurrentenmenge behaupten muss.
Folglich beziehen sich notwendige Bedingungen in Form des Pontrjaginschen Maximumprinzips,
hergeleitet auf der Basis von Nadelvariationen,
auf optimale Kandidaten ``h"oherer Qualit"at''.
Dementsprechend spricht man im Zusammenhang mit Richtungsvariationen von schwachen lokalen Optimalstellen,
während die Nadelvariationsmethode zum Begriff der starken lokalen Optimalität führt. \\[2mm]
Das dritte Kapitel behandelt verschiedene Erweiterungen der Standardaufgabe.
Zunächst die Zielfunktionale in gemischter Form und die Aufgabe mit freiem Anfangs- und Endzeitpunkt unter Zeitschranken.
Die weiteren Aufgabenklassen, d.\,h. die optimalen Multiprozesse, die Steuerung zeitverzögerter Systeme und die Differentialspiele,
besitzen ihre eigenen Charakteristiken und sind methodisch von eigenständigem Interesse. \\[2mm]
Im vierten Kapitel behandeln wir die Standardaufgabe (\ref{Einleitung3}) mit unendlichem Zeithorizont.
Dabei m"ussen wir im Gegensatz zur Standardaufgabe (\ref{Einleitung2}) den unbeschr"ankten Zeithorizont beachten,
der die Problemstellung grundlegend ändert.
Im Zuge der Untersuchungen ergibt sich schließlich die Frage,
wie sich die Steuerungsprobleme mit endlichem und unendlichem Zeithorizont einander zuordnen.
Wir werden zeigen,
dass sich die Ergebnisse für die Standardaufgabe (\ref{Einleitung2}) vollständig aus den Resultaten für die
Aufgabe (\ref{Einleitung3}) ableiten lassen.
Somit stellen die Betrachtungen im vierten Kapitel direkte Verallgemeinerungen der angewandten Methoden im zweiten Kapitel dar. \\[2mm]
Im fünften Kapitel diskutieren die Steuerung von Volterraschen Integralgleichungen.
Diese Aufgabenklasse zeichnet sich durch das Vorhandensein von zwei Zeitvariablen aus.
Deswegen muss die Nadelvariationsmethode auf einen zweidimensionalen Zeitbereich erweitert werden.
Weiterhin treffen wir bei der Behandlung der Aufgabe mit freiem Anfangs- und Endzeitpunkt auf Zustände,
die von jeweils einer der beiden Zeitvariablen beeinflusst werden.
Dieser Umstand ist noch kein Bestandteil in der Aufgabenstellung (\ref{Einleitung4}).
Im Vergleich zur Standardaufgabe (\ref{Einleitung2}) bestehen die Herausforderungen dieses Kapitels darin,
den Einfluss der beiden Zeitvariablen aufzuzeigen und die Nadelvariationsmethode auf den zweidimensionalen Zeitbereich zu erweitern. \\[2mm]
In dem Streben nach einer höheren Übersichtlichkeit und besseren Lesbarkeit der einzelnen Kapitel
sind im Rahmen der Untersuchungen der verschiedenen Aufgabenklassen die komplexeren und aufwendigeren Beweise lediglich skizziert.
Aber ohne Beweise und dem Prüfen der Richtigkeit der Aussagen nimmt man der Mathemtik ihre Natur.
Um diesen Anspruch doch gerecht zu werden,
enthält der sehr umfangreiche Anhang wichtige Grundlagen und Resultate,
die in die komplexen Beweise der Optimalitätsprinzipien ma"sgeblich einflie"sen und die Beweisskizzen vervollständigen.
Weite Teile des Anhangs sind deswegen nur für diejenigen Leser von Interesse,
die die Beweise genauer nachvollziehen möchten.

%% file: 0-Bezeichnungen.tex
\section*{Wichtige Bezeichnungen}
Die Arbeit ist untergliedert in die Untersuchungen von Aufgabenstellungen der 
Klassischen Variationsrechnung und der Optimalen Steuerung in Bezug auf
\begin{enumerate}
\item[(1)] schwache lokale Minimalstellen,
\item[(2)] starke lokale Minimalstellen,
\item[(3)] verschiedene Erweiterungen des Standardproblems,
\item[(4)] die elementare Aufgabe mit freiem rechten Endpunkt.
\end{enumerate}
Insbesondere der Aufbau der Steuerungsprobleme, die wir behandeln werden,
enthält die wiederkehrenden Elemente von Zielfunktional, Dynamik, Randbedingungen,
Zustands- und Steuerungsbeschränkungen.
Dafür verwenden wir meistens die Bezeichnungen \\[2mm]
\begin{tabular}{lll}
$J$ &--& für das Zielfunktional, \\[1mm]
$f$ &--& für den Integranden im Zielfunktional, \\[1mm]
$S$ &--& für den Wiedergewinnungswert, \\[1mm]
$\varphi$ &--& für die rechte Seite des dynamischen Systems, \\[1mm]
$h$ &--& für Randbedingungen an die Zustandstrajektorie, \\[1mm]
$g$ &--& für Zustandsbeschränkungen, \\[1mm]
$U$ &--& für den Steuerbereich, \\[1mm]
$x(\cdot)$ &--& für den Zustand, \\[1mm]
$u(\cdot)$ &--& für die Steuerung. \\[2mm]
\end{tabular}

Die verschiedenen Problemklassen besitzen ihre eigenen Charakteristiken und erfordern teils individuelle Voraussetzungen.
Das betrifft einerseits die Art der lokalen Optimalität und andererseits die Eigenheiten der spezifischen Problemklassen.
Deswegen verwenden wir unterschiedliche Bezeichnungen für die Abbildungen und die Annahmen in den Untersuchungen von
schwachen und starken Optimalstellen,
und versehen sie zudem mit einem zusätzlichen, zur spezifischen Aufgabenstellung gehörenden Index. \\[2mm]
Im Rahmen der Klassischen Variationsrechnung und der Anwendung der Richtungsvariationen verwenden wir für die Annahmen an die Aufgabenklasse die Schreibweise $\mathscr{A}$,
z.\,B. $\mathscr{A}_{\rm adm}$ für die Menge der zul"assigen Steuerungsprozesse. \\
Um im Gegensatz dazu den Rahmen der Optimalen Steuerung und insbesondere der Nadelvariationsmethode hervorzuheben,
verwenden wir z.\,B. für die Menge der zul"assigen Steuerungsprozesse die Abkürzung $\mathscr{B}_{\rm adm}$. \\
Ferner werden die Pontrjagin-Funktion $H$ und die Hamilton-Funktion $\mathscr{H}$ eingeführt.

\newpage
Zussammenfassend treten bei der Methode der Richtungsvariation die Bezeichnungen \\[2mm]
\begin{tabular}{lll}
$\mathcal{L}\,$agrange-Aufgabe &--& $\mathscr{A}^{\mathcal{L}}_{\rm adm}$, $\mathscr{A}^{\mathcal{L}}_{\rm Lip}$; \\[2mm]
$\mathcal{I}\,$soperimetrische Aufgabe &--& $\mathscr{A}^{\mathcal{I}}_{\rm adm}$, $\mathscr{A}^{\mathcal{I}}_{\rm Lip}$; \\[2mm]
$\mathcal{S}\,$tandardaufgabe &--& $\mathscr{A}^{\mathcal{S}}_{\rm adm}$, $\mathscr{A}^{\mathcal{S}}_{\rm Lip}$, $H^{\mathcal{S}}$; \\[2mm]
$\mathcal{F}\,$reie Zeit &--& $\mathscr{A}^{\mathcal{F}}_{\rm adm}$, $\mathscr{A}^{\mathcal{F}}_{\rm Lip}$; \\[2mm]
\end{tabular}

auf.
Für die Nadelvariationsmethode verwenden wir die Symbole \\[2mm]
\begin{tabular}{lll}
$\mathcal{S}\,$tandardaufgabe &--&$\mathscr{B}^{\,\mathcal{S}}_{\rm adm}$, $\mathscr{B}^{\,\mathcal{S}}_{\rm Lip}$,
                                  $H^{\mathcal{S}}$, $\mathscr{H}^{\mathcal{S}}$; \\[2mm]
$\mathcal{M}\,$ultiprozesse &--& $\mathscr{B}^{\,\mathcal{M}}_{\rm adm}$, $\mathscr{B}^{\,\mathcal{M}}_{\rm Lip}$,
                                 $H^{\mathcal{M}}$, $\mathscr{H}^{\mathcal{M}}$; \\[2mm]
$\mathcal{U}\,$nendlicher Zeithorizont\hspace*{1.9mm} &--& $\mathscr{B}^{\,\mathcal{U}}_{\rm adm}$, $\mathscr{B}^{\,\mathcal{U}}_{\rm Lip}$,
                                          $\mathscr{B}^{\,\mathcal{U}}_{\lim}$, $H^{\mathcal{U}}$, $\mathscr{H}^{\,\mathcal{U}}$; \\[2mm]
$\mathcal{Z}\,$eitverzögerte Systeme &--& $\mathscr{B}^{\,\mathcal{Z}}_{\rm adm}$, $\mathscr{B}^{\,\mathcal{Z}}_{\rm Lip}$,
                                          $H^{\mathcal{Z}}$, $\mathscr{H}^{\mathcal{Z}}$; \\[2mm]
$\mathcal{I}\,$ntegralgleichungen &--& $\mathscr{B}^{\,\mathcal{I}}_{\rm adm}$, $\mathscr{B}^{\,\mathcal{I}}_{\rm Lip}$,
                                       $H^{\mathcal{I}}$, $\mathscr{H}^{\mathcal{I}}$; \\[2mm]
$\mathcal{F}\,$reie Zeit &--& $\mathscr{B}^{\mathcal{F}}_{\rm adm}$, $\mathscr{B}^{\mathcal{F}}_{\rm Lip}$. \\[2mm]
\end{tabular}

In der Behandlung von Problemen der Klassischen Varationsrechnung und der Optimalen Steuerung wandeln sich die Rahmenbedingungen bezüglich
der Art der Funktionen und wir treffen die folgenden Funktionsklassen an: \\[2mm]
\begin{tabular}{lll}
$C([a,b],\R)$ &--& Raum der stetigen Funktionen; \\[1mm]
$C_0([a,b],\R)$ &--& Raum der stetigen Funktionen $x(\cdot)$ mit $x(a)=0$; \\[1mm]
$C_1([a,b],\R)$ &--& Raum der stetig differenzierbaren Funktionen; \\[1mm]
$C_0(\R_+,\R)$ &--& Raum der stetigen Funktionen, die im Unendlichen verschwinden; \\[1mm]
$C_{\lim}(\R_+,\R)$ &--& Raum der stetigen Funktionen, die im Unendlichen konvergieren; \\[1mm]
$L_\infty([a,b],\R)$ &--& Raum der messbar und beschränkten Funktionen; \\[1mm]
$W^1_\infty([a,b],\R)$ &--& Raum der absolutstetigen Funktionen. \\[2mm]
\end{tabular}

Hierbei bezeichnet $\R_+$ das halboffene Intervall $[0,\infty)$.
Außerdem können wir den Raum $C_{\lim}(\R_+,\R)$ mit dem Raum $C(\overline{\R}_+,\R)=C([0,\infty],\R)$ identifizieren. \\[2mm]
In unseren Untersuchungen treten einseitige Grenzwerte wiederholt auf.
Für den links- bzw. rechtseitigen Grenzwert einer Abbildung $\psi$ an der Stelle $x_0$ schreiben wir
$$\lim_{x \to x_0^-} \psi(x) \quad\mbox{bzw.}\quad \lim_{x \to x_0^+} \psi(x).$$ 

%% file: 0-Stueckweise.tex
\section*{Stückweise Stetigkeit und Differenzierbarkeit}
Eine gewisse Sonderstellung haben die elementaren Aufgaben mit freiem rechten Endpunkt inne.
Diese betrachten wir im Zuge
der Standardaufgabe, der Multipozesse, der zeitverzögerten Systeme, des unendlichen Zeithorizontes
und der Volterraschen Integralgleichung
in den Abschnitten \ref{AbschnittPMPeinfach}, \ref{AbschnittPMPeinfachMP}, \ref{AbschnittPMPeinfachRet}, \ref{AbschnittPMPeinfachUH} und
\ref{AbschnittPMPeinfachIGL}. \\
Für die Aufgaben mit freiem rechten Endpunkt eignet sich die Methode der einfachen Nadelvariation.
Sie erlaubt eine vergleichsweise einfache Beweisführung und die Wahl von stückweise stetigen und
stückweise stetig differenzierbaren Funktionen.
Dieser  Rahmen gestattet zudem die auftretenden Differentationen und Integrationen im klassischen Sinn aufzufassen,
d.\,h. es bezeichnen $\dot{x}(t)$ die klassische Zeitableitung und $\displaystyle \int_a^b \,dt$ das (für $b=\infty$ uneigentliche) Riemann-Integral.
Da wir die elementaren Aufgaben mit freiem rechten Endpunkt über den stückweise stetigen Funktionen betrachten,
verwenden für die Menge der zulässigen Steuerungsprozesse die Bezeichnungen
$\mathscr{D}^{\,\mathcal{S}}_{\rm adm}$,..., $\mathscr{D}^{\,\mathcal{I}}_{\rm adm}$.

\begin{definit}[Stückweise Stetigkeit]
Die Funktion $f(\cdot): [a,b] \to \R$ heißt stückweise stetig, \index{Funktion, stückweise stetige}
wenn sie in endlich vielen Stellen $a<s_1<s_2<...<s_N<b$ Sprünge besitzt,
d.\,h. in diesen Stellen existieren beide einseitigen Grenzwerte von $f$ im eigentlichen Sinn.
In den Stellen $s_1,...,s_N$ wählen die Werte der Funktion $f$ so, dass $f$ rechtsseitig stetig ist.
\end{definit}

\begin{definit}[Stückweise stetige Differenzierbarkeit]
Die Funktion $f(\cdot): [a,b] \to \R$ heißt stückweise stetig differenzierbar,
\index{Funktion, stückweise stetige!ststetigdiff@--, stückweise stetig differenzierbare}
wenn sie auf $[a,b]$ stetig und in den endlich vielen Teilintervallen $(a,s_1),(s_1,s_2),...,(s_N,b)$ stetig differenzierbar ist,
sowie ihre Ableitung eine stückweise stetige und in den Stellen $s_1,...,s_N$ rechtsseitig stetige Funktion über $[a,b]$ ist.
\end{definit}

Diese Klassen stückweise stetiger Funktionen kennzeichnen wir mit \\[2mm]
\begin{tabular}{lll}
$PC([a,b],\R)$ &--& Raum der stückweise stetigen Funktionen, \\[1mm]
$PC_1([a,b],\R)$ &--& Raum der stückweise stetig differenzierbaren Funktionen. \\[2mm]
\end{tabular}

In der Aufgabe mit unendlichem Zeithorizont wird das Intervall $[a,b]$ durch $[0,\infty)$ ersetzt.
Es ergeben sich die Funktionenräume  \\[2mm]
\begin{tabular}{lll}
$PC([0,\infty),\R^n)$ &--& Raum der stückweise stetigen Funktionen, \\[1mm]
$PC_1([0,\infty),\R^n)$ &--& Raum der stückweise stetig differenzierbaren Funktionen \\[2mm]
\end{tabular}

über $[0,\infty)$.
In dem Fall des unbeschränkten Intervalls $[0,\infty)$ definieren wir:
\begin{definit}[Stückweise Stetigkeit]
Die Funktion $f(\cdot): [0,\infty) \to \R$ heißt stückweise stetig, 
wenn sie über $[0,\infty)$ beschränkt und über jedem endlichen Intervall $[0,T]$ stückweise stetig ist.
\end{definit}

\begin{definit}[Stückweise stetige Differenzierbarkeit]
Die Funktion $f(\cdot): [0,\infty) \to \R$ heißt stückweise stetig differenzierbar,
wenn sie über $[0,\infty)$ beschränkt und über jedem endlichen Intervall $[0,T]$ stückweise stetig differenzierbar ist.
\end{definit}

%% file: 1-0-KlassischeVariationsrechnung.tex
\section{Richtungsvariationen in der Variationsrechnung und Optimalen Steuerung} \label{KapitelKlassischeVariationsrechnung}
Die Entdeckungen der Variationsprinzipien der Optik und der Mechanik im 17. und 18. Jahrhundert stehen f"ur den Beginn eines neuen
Kapitels im Verst"andnis grundlegender Vorg"ange in Natur und Technik.
Angeregt durch das Problem der Brachistochrone\index{Brachistochrone},
n"amlich der Frage nach der Kurve k"urzester Fallzeit, entstand die Variationsrechnung. \\[2mm]
Im Problem der Brachistochrone ist diejenige Kurve gesucht,
auf der ein Massenpunkt allein unter der Wirkung der Schwerkraft unter Vernachl"assigung der Reibung
am schnellsten von einem Punkt
\begin{figure}[h]
	\centering
	\fbox{\includegraphics[width=4cm]{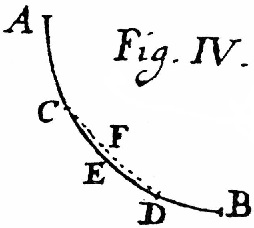}}
	\caption[Brachistochrone]{Kurve k"urzester Fallzeit von $A$ nach $B$.}
\end{figure}
$A$ zu einem tieferliegenden Punkt $B$ gelangt. \\
Die L"osung ist eine nach oben ge"offnete Zykloide mit der Parameterdarstellung
$$x(t)=r(t-\sin t), \qquad y(t)=r(1-\cos t), \qquad 0 \leq t \leq \pi.$$
Den bedeutendsten Einfluss auf die Entwicklung der Klassischen Variationsrechnung nahmen im 18. Jahrhundert
Leonhard Euler (1707--1783) und Joseph-Louis Lagrange (1736--1813),
die die grundlegende notwendige Bedingung in Form einer Differentialgleichung,
der Euler-Lagrange-Gleichung, entdeckten. \\[2mm]
Wir werden zu Beginn dieses Abschnitts die Euler-Lagrangesche Gleichung ausf"uhrlich diskutieren.
Als n"achstes gehen wir zur Untersuchung der Lagrange-Aufgabe\index{Lagrange-Aufgabe} der Klassischen Variationsrechnung "uber.
Dabei folgen wir Ioffe \& Tichomirov \cite{Ioffe}. \\
Die Lagrange-Aufgabe ist die Standardaufgabe
\begin{eqnarray*}
&& J\big(x(\cdot),u(\cdot)\big) = \int_{t_0}^{t_1} f\big(t,x(t),u(t)\big) \, dt \to \inf, \\
&& \dot{x}(t)= \varphi\big(t,x(t),u(t)\big), \\
&&h_0\big(x(t_0)\big)=0, \qquad h_1\big(x(t_1)\big)=0, \\
&& u(t) \in U, \qquad U=\R^m.
\end{eqnarray*}
f"ur Paare $\big(x(\cdot),u(\cdot)\big) \in C_1([t_0,t_1],\R^n) \times C([t_0,t_1],\R^m)$,
die die Nebenbedingungen erf"ullen.
Als Optimalit"atsbegriff verwenden wir das schwache lokale Minimum.
Darunter verstehen wir ein zul"assiges Paar $\big(x_*(\cdot),u_*(\cdot)\big)$ der Lagrange-Aufgabe,
zu dem eine Zahl $\varepsilon > 0$ derart existiert, dass die Ungleichung 
$$J\big(x(\cdot),u(\cdot)\big) \geq J\big(x_*(\cdot),u_*(\cdot)\big)$$
f"ur alle zul"assigen Paare $\big(x(\cdot),u(\cdot)\big)$ mit 
$\|x(\cdot)-x_*(\cdot) \|_{C_1} < \varepsilon$, $\|u(\cdot)-u_*(\cdot)\|_\infty < \varepsilon$ gilt. \\[2mm]
Der Umgebungsbegriff in der Definition eines schwachen lokalen Minimum erfordert geeignete Annahmen an die Standardaufgabe.
Wir fordern in unseren Betrachtungen "uber schwache lokale Minimalstellen,
dass die Abbildungen $f(t,x,u)$, $\varphi(t,x,u)$ und die Randbedingungen $h_0(x_0)$, $h_1(x_1)$
stetig und stetig differenzierbar auf einer gewissen, gleichmäßigen Umgebungsmenge des Referenzpaares $\big(x_*(\cdot),u_*(\cdot)\big)$ sind. \\[2mm]
Die Zykloide als L"osung der Brachistochrone mit der oben angegebenen Parameterdarstellung besitzt im Anfangspunkt eine unbeschr"ankte
Ableitung.
Wir stehen damit vor der Frage, ob die angegebene Zykloide im Problem der Brachistochrone\index{Brachistochrone} variierbar ist.
Im Rahmen der vorliegenden Arbeit k"onnen wir diese Frage nicht kl"aren. \\
Obwohl dieser Umstand stets in der Literatur auftaucht und eingehend bekannt ist,
liegt unserer Kenntnis nach für das Problem der Brachistochrone in Gestalt der Standardaufgabe noch keine einwandfreie Analyse vor.

%% file: 1-1-EulerLagrange.tex
\subsection{Die Euler-Lagrangesche Gleichung der Klassischen Variationsrechnung} \label{AbschnittELG}
\subsubsection{Die Herleitung der Euler-Lagrange-Gleichung}
Die Grundaufgabe der Klassischen Variationsrechnung\index{Grundaufgabe der Variationsrechnung}
besteht in der Minimierung des Zielfunktionals
\begin{equation}\label{ELG1} J\big(x(\cdot)\big) = \int_{t_0}^{t_1} f\big(t,x(t),\dot{x}(t)\big) \, dt \to \inf \end{equation}
unter den Nebenbedingungen
\begin{equation}\label{ELG2} x(t_0)=x_0, \qquad x(t_1)=x_1. \end{equation}
Dabei sind in (\ref{ELG1}) und (\ref{ELG2}) die Punkte $t_0<t_1$ und $x_0,x_1$ fest vorgegeben.
Ferner setzen wir voraus,
dass in (\ref{ELG1}) die Funktion $f:\R \times \R \times \R \to \R$ 
einmal stetig differenzierbar ist. \\[2mm]
Wir nennen die Funktion $x(\cdot):[t_0,t_1]\to \R$ zul"assig in der Aufgabe (\ref{ELG1})--(\ref{ELG2}),
wenn $x(\cdot)$ stetig differenzierbar ist und die Randbedingungen (\ref{ELG2}) erf"ullt.
Die zul"assige Funktion $x_*(\cdot) \in C_1([t_0,t_1],\R)$ ist eine schwache lokale Minimalstelle
\index{Minimum, schwaches lokales!Variation@-- Variationsrechnung}
in der Aufgabe
(\ref{ELG1})--(\ref{ELG2}),
falls ein $\varepsilon > 0$ derart existiert, dass die Ungleichung 
$$J\big(x(\cdot)\big) \geq J\big(x_*(\cdot)\big)$$
f"ur alle zul"assigen Funktionen $x(\cdot)$ mit $\| x(\cdot)-x_*(\cdot) \|_{C_1} < \varepsilon$ gilt.

\newpage
Zur Bestimmung einer notwendigen Optimalit"atsbedingung f"ur eine schwache lokale Minimalstelle $x_*(\cdot)$
in der Aufgabe (\ref{ELG1})--(\ref{ELG2}) betrachten wir die Funktion einer Ver"anderlichen
\begin{equation}\label{ELG3}
\varphi(\lambda)= J\big(x_*(\cdot) + \lambda x(\cdot)\big)
                = \int_{t_0}^{t_1} f\big(t,x_*(t)+\lambda x(t),\dot{x}_*(t)+\lambda \dot{x}(t)\big) \, dt,
\end{equation}
die durch die Variation $x_\lambda(\cdot) =x_*(\cdot) + \lambda x(\cdot)$ von $x_*(\cdot)$
bez"uglich der Richtung $x(\cdot)$ erzeugt wird.
Dies wirft die naheliegende Frage auf, ob $x_\lambda(\cdot)$ die Randbedingungen (\ref{ELG2}) erf"ullt.
Fordern wir f"ur die Richtung $x(\cdot)$ die Bedingung
$$x(\cdot) \in M_0 = \{ y(\cdot) \in C_1([t_0,t_1],\R) \,|\, y(t_0)=y(t_1)=0 \},$$
so ist die Variation $x_\lambda(\cdot)$ stets ein zul"assiges Element in der Aufgabe (\ref{ELG1})--(\ref{ELG2}). \\[2mm]
Unter unseren Voraussetzungen "uber die Funktionen $f$, $x_*(\cdot)$, $x(\cdot)$ ist in (\ref{ELG3}) die Differentiation unter dem
Integralzeichen zul"assig, und es gilt mit den Funktionen
$$q(t)=f_x\big(t,x_*(t),\dot{x}_*(t)\big), \qquad p(t)=f_{\dot{x}}\big(t,x_*(t),\dot{x}_*(t)\big)$$
f"ur die erste Variation
$$\varphi'(0)= \delta J\big(x_*(\cdot)\big)x(\cdot)
             = \int_{t_0}^{t_1} \big( q(t) x(t) + p(t) \dot{x}(t) \big) \, dt.$$
Da zu $x(\cdot) \in M_0$ f"ur hinreichend kleine $\lambda>0$ die Variation $x_\lambda(\cdot)$ die Bedingung
$$\| x_\lambda(\cdot)-x_*(\cdot) \|_{C_1} = \| \lambda x(\cdot) \|_{C_1} < \varepsilon$$
erf"ullt,
erhalten wir folgende notwendige Bedingung f"ur eine schwache lokale Minimalstelle:
\begin{equation}\label{ELG4}
\varphi'(0)= \delta J\big(x_*(\cdot)\big)x(\cdot) = 0 \qquad \mbox{f"ur alle } x(\cdot) \in M_0.
\end{equation}

\begin{satz} \label{SatzELG}
Es sei die Funktion $f$ stetig differenzierbar.
Dann gen"ugt ein schwaches lokales Minimum $x_*(\cdot)$ der Aufgabe (\ref{ELG1})--(\ref{ELG2}) auf $[t_0,t_1]$
der Euler-Lagrange-Gleichung\index{Euler-Lagrangesche Gleichung}
\begin{equation}\label{ELG5}
\bigg(-\frac{d}{dt} f_{\dot{x}} + f_x \bigg)\bigg|_{x_*(t)}
=-\frac{d}{dt} f_{\dot{x}}\big(t,x_*(t),\dot{x}_*(t)\big) + f_x\big(t,x_*(t),\dot{x}_*(t)\big) =0.
\end{equation}
\end{satz}

\begin{definition}
Jede L"osung $x_*(\cdot) \in C_1([t_0,t_1],\R)$ von (\ref{ELG5}) nennen wir eine Extremale\index{Extremale}.
\end{definition}

{\bf Beweis} Wir zeigen Satz \ref{SatzELG} auf zweierlei Weisen:
Im ersten Vorgehen f"uhren wir in (\ref{ELG4}) bzgl. dem Term $p(t) \dot{x}(t)$ eine partielle Integration durch
und wenden das Fundamentallemma von Lagrange an.
Dies erfordert die zus"atzliche Annahme, dass die Funktion $f$ zweimal stetig differenzierbar ist.
Im Gegensatz dazu erfordert die partielle Integration bez"uglich dem Term $q(t) x(t)$ in (\ref{ELG4}) keine zus"atzlichen Annahmen.
Das Fundamentallemma von du Bois-Reymond liefert dann die Euler-Lagrangesche Gleichung in Satz \ref{SatzELG}. \\[2mm] 
Wegen der zweimaligen stetigen Differenzierbarkeit von $f$ und wegen $x(\cdot) \in M_0$
erhalten wir nach partieller Integration des Terms $p(t) \dot{x}(t)$ f"ur (\ref{ELG4}) die Darstellung
$$\int_{t_0}^{t_1} \big( q(t) - \dot{p}(t) \big) x(t) \, dt  = 0 \qquad \mbox{f"ur alle } x(\cdot) \in M_0.$$
Daraus schlie"sen wir mit dem Fundamentallemma von Lagrange,
dass die stetige Funktion $\varphi(t)= q(t) - \dot{p}(t)$ auf $[t_0,t_1]$ identisch verschwinden muss.
Das bedeutet ausf"uhrlich geschrieben, dass die Gleichung (\ref{ELG5}) gilt. \\[2mm]
F"uhren wir demgegen"uber in (\ref{ELG4}) die partielle Integration bez"uglich dem Term $q(t) x(t)$ durch,
erhalten wir f"ur alle $x(\cdot) \in M_0$ die Darstellung
$$\int_{t_0}^{t_1} \big( -Q(t)+ p(t) \big) \dot{x}(t) \, dt  = 0, \qquad Q(t) = \int_{t_0}^t q(s) \, ds.$$
Es folgt dann aus dem Fundamentallemma von du Bois-Reymond,
dass die stetige Funktion $\psi(t)= -Q(t) +p(t)$ auf $[t_0,t_1]$ konstant sein muss.
Dies bedeutet ferner, dass $\psi(\cdot)$ stetig differenzierbar ist und die Ableitung auf $[t_0,t_1]$ verschwindet.
Da nun $Q(\cdot)$ eine stetige Ableitung besitzt, muss das auch f"ur $p(\cdot)$ gelten.
Damit ist (\ref{ELG5}) gezeigt. \hfill $\blacksquare$

\begin{lemma}[Fundamentallemma von Lagrange]\index{Fundamentallemma!Lagrange@-- von Lagrange}
Es sei $\varphi(\cdot)$ auf $[t_0,t_1]$ stetig mit
$$\int_{t_0}^{t_1} \varphi(t) x(t) \, dt  = 0 \qquad \mbox{f"ur alle } x(\cdot) \in M_0.$$
Dann gilt $\varphi(t)\equiv 0$ auf $[t_0,t_1]$.
\end{lemma}

{\bf Beweis} Angenommen, es ist $\varphi(\tau)\not=0$ f"ur ein $\tau \in (t_0,t_1)$.
Es sei $\varphi(\tau)>0$.
Dann gibt es ein $\varepsilon >0$ mit $\varphi(t)>0$ auf $\Delta=[\tau-\varepsilon,\tau+\varepsilon] \subseteq [t_0,t_1]$.
Wir w"ahlen
$$z(t) = \left\{ \begin{array}{ll}
        (t-\tau+\varepsilon)^2(t-\tau-\varepsilon)^2, & t \in \Delta, \\[1mm] 0,& t \not\in \Delta. \end{array}\right.$$ 
Dann geh"ort $z(\cdot)$ offenbar zur Menge $M_0$ und es gilt
$$\int_{t_0}^{t_1} \varphi(t) z(t) \, dt = \int_{\Delta} \varphi(t) z(t) \, dt > 0,$$
im Widerspruch zur Annahme des Fundamentallemma von Lagrange. \hfill $\blacksquare$

\begin{lemma}[Fundamentallemma von du Bois-Reymond]\index{Fundamentallemma!du Bois@-- von du Bois-Reymond}
Es sei $\psi(\cdot)$ eine stetige Funktion mit
$$\int_{t_0}^{t_1} \psi(t) \dot{x}(t) \, dt  = 0 \qquad \mbox{f"ur alle } x(\cdot) \in M_0.$$
Dann ist $\psi(\cdot)$ auf $[t_0,t_1]$ konstant.
\end{lemma}

{\bf Beweis} Wir betrachten die Funktion
$$z(t)=\int_{t_0}^t \big(\psi(s)-c_0\big) \, ds, \qquad c_0:= \frac{1}{t_1-t_0} \int_{t_0}^{t_1} \psi(t) \, dt.$$
Dann gelten $z(t_0)=z(t_1)=0$, also $z(\cdot) \in M_0$, $\dot{z}(t)=\psi(t)-c_0$ und wir erhalten
$$0=\int_{t_0}^{t_1} \psi(t) \dot{z}(t) \, dt
   =\int_{t_0}^{t_1} \psi(t) \dot{z}(t) \, dt - c_0 \int_{t_0}^{t_1} \dot{z}(t) \, dt
   = \int_{t_0}^{t_1} \big(\psi(t) - c_0 \big)^2 \, dt.$$
Daher muss $\psi(\cdot)$ auf $[t_0,t_1]$ konstant sein. \hfill $\blacksquare$

\begin{bemerkung}{\rm
Statt $x(\cdot) \in M_0$ wird das Lemma von du Bois-Reymond oft f"ur stetige Funktionen $y(\cdot)$ mit Integralmittelwert Null,
d.\,h.
$$\int_{t_0}^{t_1} y(t) \, dt =0,$$
formuliert.
Offensichtlich besitzt $\dot{x}(\cdot)$ f"ur jedes $x(\cdot) \in M_0$ den Integralmittelwert Null.
\hfill $\square$}
\end{bemerkung}

\begin{beispiel} \label{BeispielELGRegler}
{\rm Wir betrachten eine Aufgabe mit quadratischem Zielfunktional:
$$J\big(x(\cdot)\big) = \int_0^T \big[\big(x(t)-1\big)^2+\dot{x}^2(t)\big] \, dt \to \inf, \qquad x(0)=0,\quad x(T)=2,
  \qquad T>0.$$
\begin{minipage}{0.57\textwidth}
Mit $f(t,x,\dot{x})=(x-1)^2+\dot{x}^2$ erhalten wir aus der Euler-Lagrange-Gleichung (\ref{ELG5}),
dass
eine schwache lokale Minimalstelle $x_*(\cdot)$ der Differentialgleichung
$$\ddot{x}(t)-x(t)=-1$$
gen"ugen muss.
Die Anpassung an die Randbedingungen $x(0)=0$, $x(T)=2$ liefert den Kandidaten
$$x_*(t)=1+\frac{1+e^{-T}}{e^T-e^{-T}}e^t- \frac{1+e^T}{e^T-e^{-T}}e^{-t}$$
\end{minipage}
\begin{minipage}{0.42\textwidth}
	\centering
	\includegraphics[width=4.5cm]{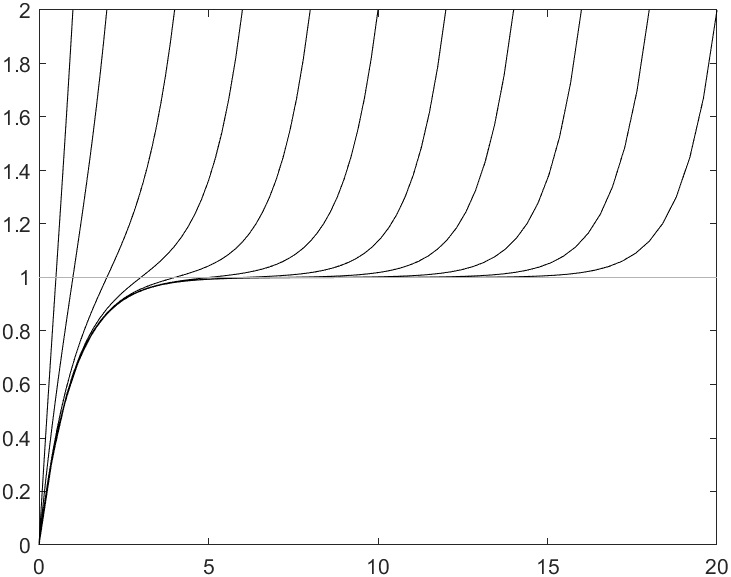}
	\captionof{figure}[Optimale Trajektorien in einem quadratischen Regler]{Optimale Trajektorien für $T=2,4,...,20$.}
	\label{AbbRegler}
\end{minipage} \\[1mm]
f"ur ein schwaches lokales Minimum. \hfill $\square$}
\end{beispiel}

\begin{folgerung} \label{FolgerungELG}
Wir geben Spezialf"alle der Euler-Lagrange-Gleichung\index{Euler-Lagrangesche Gleichung} (\ref{ELG5}) an,
die sich f"ur eine L"osung $x_*(\cdot)$ des Variationsproblems ergeben.
\begin{enumerate}
\item[(a)] H"angt die Funktion $f$ nicht von $\dot{x}$ ab, dann muss auf $[t_0,t_1]$ gelten:
           $$f_x\big(t,x_*(t)\big) =0.$$
\item[(b)] H"angt die Funktion $f$ nicht von $x$ ab, dann folgt auf $[t_0,t_1]$:
           $$p(t)=f_{\dot{x}}\big(t,\dot{x}_*(t)\big) = \mbox{konstant}.$$           
\item[(c)] H"angt die Funktion $f$ nicht von $t$ ab,
           so erh"alt man mit (\ref{ELG5}),
           dass die Ableitung nachstehender Funktion $H(\cdot)$ auf $[t_0,t_1]$ identisch verschwindet; oder gleichbedeutend
           $$H(t)=f_{\dot{x}}\big(x_*(t),\dot{x}_*(t)\big) \cdot \dot{x}_*(t)-f\big(x_*(t),\dot{x}_*(t)\big)
                  = \mbox{konstant} \quad \mbox{auf } [t_0,t_1].$$
\end{enumerate}
\end{folgerung}

An die Betrachtungen zur Grundaufgabe (\ref{ELG1})--(\ref{ELG2}) schlie"st sich die Frage nach notwendigen 
Optimalit"atsbedingungen f"ur eine schwache lokale Minimalstelle in der Aufgabe mit freiem Endpunkt an:
\begin{equation} 
\label{ELGFE1} J\big(x(\cdot)\big) = \int_{t_0}^{t_1} f\big(t,x(t),\dot{x}(t)\big) \, dt \to \inf, \qquad
               x(t_0)=x_0, \quad x(t_1) \mbox{ frei}.
\end{equation}
Unter den gleichen Voraussetzungen  und mit dem gleichen Vorgehen wie in der Grundaufgabe (\ref{ELG1})--(\ref{ELG2}) erhalten wir 
f"ur eine schwache lokale Minimalstelle
\begin{equation}\label{ELGFE2}
\delta J\big(x_*(\cdot)\big)x(\cdot) = \int_{t_0}^{t_1} \big( q(t) x(t) + p(t) \dot{x}(t) \big) \, dt
   = 0 \qquad \mbox{f"ur alle } x(\cdot) \in M_1.
\end{equation}
Dabei bezeichnen $M_1=\{ y(\cdot) \in C_1([t_0,t_1],\R) \,|\, y(t_0)=0 \}$ und $q(\cdot),p(\cdot)$ die Funktionen
$$q(t)=f_x\big(t,x_*(t),\dot{x}_*(t)\big), \qquad p(t)=f_{\dot{x}}\big(t,x_*(t),\dot{x}_*(t)\big).$$
Integrieren wir in (\ref{ELGFE2}) zun"achst partiell den Term $q(t) x(t)$, so folgt
$$Q(t_1)x(t_1)+\int_{t_0}^{t_1} \big( -Q(t)+ p(t) \big) \dot{x}(t) \, dt  = 0 \quad \mbox{f"ur alle } x(\cdot) \in M_1,
  \quad Q(t) = \int_{t_0}^t q(s) \, ds.$$
F"ur $x(\cdot) \in M_0$ erhalten wir nach dem Lemma von du Bois-Reymond die Euler-Lagrange-Gleichung (\ref{ELG5}).
Au"serdem liefert die Argumentation nach du Bois-Reymond,
dass die Funktion $p(\cdot)$ stetig differenzierbar ist.
Daher d"urfen wir ohne weitere Annahmen in (\ref{ELGFE2}) den Term $p(t) \dot{x}(t)$ partiell integrieren.
Damit folgt
$$p(t_1)x(t_1)+\int_{t_0}^{t_1} \big( q(t) - \dot{p}(t) \big) x(t) \, dt  = 0 \qquad \mbox{f"ur alle } x(\cdot) \in M_1.$$
Darin verschwindet $\varphi(t)=q(t) - \dot{p}(t)$ identisch auf $[t_0,t_1]$.
Damit diese Gleichung auf $M_1$ erf"ullt ist, muss also $p(t_1)=0$ gelten.

\begin{satz} \label{SatzELGFE}
Es sei $f$ stetig differenzierbar.
Ist $x_*(\cdot)$ ein schwaches lokales Minimum in der Aufgabe (\ref{ELGFE1}),
dann ist notwendig, dass auf $[t_0,t_1]$ die Euler-Lagrange-Gleichung
\begin{equation}\label{SatzELGFE1}
\bigg(-\frac{d}{dt} f_{\dot{x}} + f_x \bigg)\bigg|_{x_*(t)}
=-\frac{d}{dt} f_{\dot{x}}\big(t,x_*(t),\dot{x}_*(t)\big) + f_x\big(t,x_*(t),\dot{x}_*(t)\big) =0
\end{equation}
zur Randbedingung
\begin{equation}\label{SatzELGFE2} p(t_1)= f_{\dot{x}}\big(t_1,x_*(t_1),\dot{x}_*(t_1)\big)=0 \end{equation}
erf"ullt ist.
\end{satz}

\begin{folgerung} \label{FolgerungELGFE}
Es folgt auf die gleiche Weise,
dass in der Aufgabe freier Randpunkte, d.\,h. $x(t_0),x(t_1)$ sind frei, folgende Randbedingungen gelten:
$$p(t_0)= f_{\dot{x}}\big(t_0,x_*(t_0),\dot{x}_*(t_0)\big)=0, \qquad p(t_1)= f_{\dot{x}}\big(t_1,x_*(t_1),\dot{x}_*(t_1)\big)=0.$$
\end{folgerung}

\begin{bemerkung}{\rm
S"amtliche Ausf"uhrungen im Abschnitt \ref{AbschnittELG} lassen sich unmittelbar auf den Vektorfall $x(\cdot) \in C_1([t_0,t_1],\R^n)$
erweitern. \hfill $\square$}
\end{bemerkung}

%% file: 1-11-Beispiele.tex
\subsubsection{Die Euler-Lagrangesche Gleichung in der Wirtschaftstheorie}\index{Kapitalakkumulation} \label{AbschnittKapitalakkumulation}
\begin{beispiel}[Isoelastische Nutzenfunktion]\index{Nutzenfunktion, isoelastische}
\index{Funktion, absolutstetige!isoelastisch@--, isoelastische Nutzenfunktion} \label{BeispielIsoNutzen}
{\rm Es bezeichnen $K(\cdot)$ das Verm"ogen einer "Okonomie,
$Y(\cdot)$ das Nationaleinkommen und $C(\cdot)$ die Konsumption.
Mit einer streng wachsenden und konkaven Produktionsfunktion $f$ sei $Y(t)=f\big(K(t)\big)$.
Zu jedem Zeitpunkt wird die Produktion $Y$ in Konsumption und Investition aufgeteilt, d.\,h.
$$f\big(K(t)\big)=C(t)+\dot{K}(t).$$
Das Wohlbefinden der Gesellschaft wird durch die zweimal stetig differenzierbare,
streng wachsende und streng konkave Nutzenfunktion $U$ dargestellt. \\[1mm]
Die Aufgabe der "Okonomischen Wachstumstheorie lautet dann:
Man bestimme "uber dem Zeitrahmen $[0,T]$ die optimale Konsumption $C(\cdot)$,
die den Gesamtnutzen der Gesellschaft,
$$J\big(K(\cdot)\big) = \int_0^T e^{-\varrho t} U\big(C(t)\big) \, dt
                      = \int_0^T e^{-\varrho t} U\big(f\big(K(t)\big)-\dot{K}(t)\big) \, dt,$$
unter den gegebenen Randbedingungen $K(0)=K_0$, $K(T)=K_T$ maximiert. \\[2mm]
In dieser Form l"asst sich keine qualitative Aussage "uber eine optimale Konsumptionspolitik finden.
Stattdessen leiten wir die Form der sogenannten isoelastischen Nutzenfunktionen mit Hilfe der Euler-Lagrangeschen Gleichung ab
und gewinnen daraus die Fundamentalgleichung f"ur die Wachstumsraten einer optimalen Konsumption:
Es gelten
\begin{eqnarray*}
     \frac{\partial}{\partial K} \big[e^{-\varrho t} U\big(f(K)-\dot{K}\big)\big]
&=&  e^{-\varrho t} U'\big(f(K)-\dot{K}\big) f'(K) =e^{-\varrho t} U'(C) f'(K), \\
     \frac{\partial}{\partial \dot{K}} \big[e^{-\varrho t} U\big(f(K)-\dot{K}\big)\big]
&=&  -e^{-\varrho t} U'\big(f(K)-\dot{K}\big) =-e^{-\varrho t} U'(C).
\end{eqnarray*}
Die Anwendung der Euler-Lagrange-Gleichung (\ref{ELG5}) liefert
\begin{eqnarray*}
0 &=& -\frac{d}{dt} \big[-e^{-\varrho t} U'\big(C_*(t)\big)\big]+ e^{-\varrho t} U'\big(C_*(t)\big) f'\big(K_*(t)\big) \\
  &=& e^{-\varrho t} \Big[ U''\big(C_*(t)\big) \dot{C}_*(t)+ \big[f'\big(K_*(t)\big)-\varrho\big]U'\big(C_*(t)\big)\Big].
\end{eqnarray*}
Nach einfachen Umformungen erhalten wir daraus die optimale Wachstumsrate
$$\frac{\dot{C}_*(t)}{C_*(t)}= \frac{U'\big(C_*(t)\big)}{C_*(t) U''\big(C_*(t)\big)} \big[\varrho - f'\big(K_*(t)\big)\big].$$
Eine isoelastische Nutzenfunktion gen"ugt der Differentialgleichung der konstanten relativen Risikoaversion
$$-\sigma= \frac{C U''(C)}{U'(C)}, \qquad C>0.$$
Dabei gilt $\sigma>0$, da die Funktion $U$ streng monoton wachsend und streng konkav ist. \\
Mit den Integrationskonstanten $k_0,k_1$ lautet die L"osung der Differentialgleichung
$$U(C)= k_0\frac{C^{1-\sigma}}{1-\sigma}+k_1 \quad\mbox{f"ur } \sigma\not=1, \qquad
  U(C)=k_0\ln(C)+k_1 \quad\mbox{f"ur } \sigma=1.$$
Die Festlegung $k_0=1,k_1=-1/(1-\sigma)$ liefert f"ur alle $\sigma >0$
$$U(C)= \frac{C^{1-\sigma}-1}{1-\sigma},$$
denn der Spezialfall $U(C)=\ln(C)$ f"ur $\sigma=1$ ergibt sich darin im Grenz"ubergang $\sigma \to 1$
mit Hilfe der Regel von l'Hospital. \\[2mm]
Damit erhalten wir bei einer isoelastischen Nutzenfunktion die Fundamentalgleichung
\begin{equation} \label{Fundamentalgleichung}
\frac{\dot{C}_*(t)}{C_*(t)}=\frac{f'\big(K_*(t)\big) - \varrho}{\sigma}
\end{equation}
f"ur die Wachstumsrate einer optimalen Konsumptionspolitik. \hfill $\square$}
\end{beispiel}

\begin{beispiel} \label{BeispielWachstumstheorie1}
{\rm Wir verwenden die Bezeichnungen aus Beispiel \ref{BeispielIsoNutzen} und es seien
$$f(K)=rK+W, \qquad r>0, \qquad K(T) \mbox{ frei}.$$
Darin bezeichnet $W$ den Lohn f"ur geleistete Arbeit.
Dies f"uhrt zu der Aufgabe
$$J\big(K(\cdot)\big) = \int_0^T e^{-\varrho t} U\big(rK(t)+W(t)-\dot{K}(t)\big) \, dt \to \sup, \qquad
  K(0)=K_0, \quad K(T) \mbox{ frei}.$$
Die Nutzenfunktion sei isoelastisch.
Daher gelten $U \in C_2$, $U'>0$, $U''<0$ f"ur $C>0$. \\[1mm]
Die Fundamentalgleichung (\ref{Fundamentalgleichung}) und die Randbedingung (\ref{SatzELGFE2}) liefern
$$\frac{\dot{C}_*(t)}{C_*(t)}=\frac{r - \varrho}{\sigma}, \qquad 0=-e^{-\varrho T} U'\big(C_*(T)\big).$$
Da $U'>0$ gilt, gibt es keinen Parameter $C>0$, der die Gleichung $U'(C)=0$ l"ost.
Deswegen erhalten wir keinen Kandidaten f"ur eine schwache lokale Minimalstelle. \\[2mm]
Bei n"aherer Betrachtung ist das wenig "uberraschend:
Da keine Beschr"ankungen an die Entwicklung des Kapitalstocks vorliegen, liefert die Maximierung des Zielfunktionals
$$J\big(K(\cdot)\big) = \int_0^T e^{-\varrho t} U\big(C(t)\big) \, dt$$
die wenig brauchbare L"osung $C_*(t)\equiv \infty$. \hfill $\square$}
\end{beispiel}

\begin{beispiel}{\rm
Wir fahren mit dem vorherigen Beispiel fort,
indem wir dem Gesamtnutzen der Gesellschaft den R"uckkaufswert zum Abschluss der Planungsperiode gegen"uberstellen:
$$\int_0^T e^{-\varrho t} U\big(C(t)\big) \, dt+ e^{-\varrho T}K(T) \to \sup.$$
Da die zul"assigen Trajektorien $K(\cdot)$ stetig differenzierbar sind,
k"onnen wir dies in "aquivalenter Form so ausdr"ucken:
$$J\big(K(\cdot)\big) = \int_0^T e^{-\varrho t} \big[U\big(rK(t)+W(t)-\dot{K}(t)\big)
                                 +\dot{K}(t)-\varrho K(t)  \big] \, dt \to \sup.$$
Die Euler-Lagrange-Gleichung f"uhrt wieder zur Konsumption mit Wachstumsrate
$$\frac{\dot{C}_*(t)}{C_*(t)}=\frac{r - \varrho}{\sigma}.$$
Ferner wird durch die Randbedingung (\ref{SatzELGFE2}) bei freiem Endpunkt,
$$0= e^{-\varrho T}\big[1- U'\big(C_*(T)\big)\big] \qquad\Leftrightarrow\qquad U'\big(C_*(T)\big)=1,$$
die optimale Konsumption eindeutig festgelegt. \hfill $\square$}
\end{beispiel}

%% file: 1-12-BolzaAufgabe.tex
\subsubsection{Notwendige Bedingungen in der einfachsten Bolza-Aufgabe} \label{BolzaAufgabe}
Wir geben nun den Weg zur Herleitung von notwendigen Bedingungen f"ur eine schwache lokale Minimalstelle in der 
einfachsten Bolza-Aufgabe\index{Bolza-Aufgabe},
\begin{equation} \label{ELGBA1}
B\big(x(\cdot)\big) = \psi_0\big(x(t_0)\big) + \psi_1\big(x(t_1)\big) + \int_{t_0}^{t_1} f\big(t,x(t),\dot{x}(t)\big) \, dt \to \inf
\end{equation}
an,
wo im Unterschied zur Aufgabe (\ref{ELG1})--(\ref{ELG2}) keine Randbedingungen auftreten,
aber daf"ur die Terminalfunktionale $\psi_0,\psi_1$ vorkommen. \index{Zielfunktional!Terminal@-- Terminalfunktional}
Wir setzen die stetige Differenzierbarkeit der Abbildungen $\psi_0,\psi_1$ und $f$ voraus.
Dann erhalten wir f"ur die erste Variation des Funktionals $B$ die Darstellung
\begin{eqnarray*}
    \delta B\big(x_*(\cdot)\big) x(\cdot)
&=& \psi_0'\big(x_*(t_0)\big) x(t_0) + \psi_1'\big(x_*(t_1)\big)x(t_1) + \delta J\big(x_*(\cdot)\big)x(\cdot) \\
&=& \psi_0'\big(x_*(t_0)\big) x(t_0) + \psi_1'\big(x_*(t_1)\big)x(t_1)
    +  \int_{t_0}^{t_1} \big( q(t) x(t) + p(t) \dot{x}(t) \big) \, dt
\end{eqnarray*}
mit den vertrauten Bezeichnungen
$$q(t)=f_x\big(t,x_*(t),\dot{x}_*(t)\big), \qquad p(t)=f_{\dot{x}}\big(t,x_*(t),\dot{x}_*(t)\big).$$
Ist $x_*(\cdot)$ ein schwaches lokales Minimum der Bolza-Aufgabe (\ref{ELGBA1}), dann muss gelten:
\begin{equation}\label{ELGBA3}
\delta B\big(x_*(\cdot)\big)x(\cdot) = 0 \qquad \mbox{f"ur alle } x(\cdot) \in C_1([t_0,t_1],\R).
\end{equation}
Auf dem Teilraum $M_0 \subset C_1([t_0,t_1],\R)$ stimmen $\delta B\big(x_*(\cdot)\big)$ und $\delta J\big(x_*(\cdot)\big)$ "uberein.
Damit ist nach Lemma \ref{SatzELG} die Euler-Lagrangesche Gleichung (\ref{ELG5}) erf"ullt.
Das Vorgehen nach du Bois-Reymond
liefert zudem die stetige Differenzierbarkeit der Funktion $p(\cdot)$.
Integrieren wir im Ausdruck f"ur die erste Variation des Funktionals $B$ den Term $p(t) \dot{x}(t)$ partiell und 
ber"ucksichtigen die Euler-Lagrange-Gleichung, so erhalten wir
$$\delta B\big(x_*(\cdot)\big)x(\cdot)
 = \big[\psi_0'\big(x_*(t_0)\big)-p(t_0)\big] x(t_0) + \big[\psi_1'\big(x_*(t_1)\big)+p(t_1)\big]x(t_1).$$
Wegen (\ref{ELGBA3}) gelten somit die zus"atzlichen Bedingungen
$$p(t_0)= \psi_0'\big(x_*(t_0)\big), \qquad p(t_1)=-\psi_1'\big(x_*(t_1)\big).$$

\begin{satz} \label{SatzELGBA}
Es seien $\psi_0,\psi_1$ und $f$ stetig differenzierbar.
Ist $x_*(\cdot)$ ein schwaches lokales Minimum in der Aufgabe (\ref{ELGBA1}),
dann ist notwendig, dass auf $[t_0,t_1]$ die Euler-Lagrange-Gleichung
\begin{equation}\label{SatzELGBA1}
\bigg(-\frac{d}{dt} f_{\dot{x}} + f_x \bigg)\bigg|_{x_*(t)}
=-\frac{d}{dt} f_{\dot{x}}\big(t,x_*(t),\dot{x}_*(t)\big) + f_x\big(t,x_*(t),\dot{x}_*(t)\big) =0
\end{equation}
mit den Randbedingungen
\begin{equation}\label{SatzELGBA2}
\left. \begin{array}{l}
p(t_0)= f_{\dot{x}}\big(t_0,x_*(t_0),\dot{x}_*(t_0)\big)=\psi_0'\big(x_*(t_0)\big), \\[1mm]
p(t_1)= f_{\dot{x}}\big(t_1,x_*(t_1),\dot{x}_*(t_1)\big)=-\psi_1'\big(x_*(t_1)\big)
\end{array} \right\}
\end{equation}
erf"ullt ist.
\end{satz}


\begin{bemerkung}{\rm 
Mit Hilfe der Funktion $y(\cdot)$, die der Differentiagleichung
$$\dot{y}(t)=f\big(t,x(t),\dot{x}(t)\big), \qquad y(t_0)=0$$
genügt,
entsteht aus einem Integralfunktional ein "aquivalentes Terminalfunktional: 
$$\int_{t_0}^{t_1} f\big(t,x(t),\dot{x}(t)\big) \, dt \to \inf \qquad\Leftrightarrow\qquad y(t_1) \to \inf.$$
Auf diese Weise lässt sich ein gemischtes Zielfunktional in der Form nach Oskar Bolza (1857--1942) in die Gestalt
$$\psi\big(x(t_0),x(t_1)\big) + \int_{t_0}^{t_1} f\big(t,x(t),\dot{x}(t)\big) \, dt= \psi\big(x(t_0),x(t_1)\big) + y(t_1)
  =\tilde{\psi}\big(x(t_0),x(t_1)\big)$$
eines Terminalfunktionals mit der zusätzlichen Nebenbedingung an die Funktion $y(\cdot)$ überführen.
Aufgaben der Variationsrechnung, die ausschließlich Terminalfunktionale enthalten,
werden zu Ehren von Adolph Mayer (1839--1907) als Mayer-Aufgabe bezeichnet. \\
In Abschnitt \ref{AbschnittGemischtesZF} gehen wir im Rahmen der Optimalen Steuerung auf Zielfunkionale in gemischter Form
bestehend aus Integral- und Terminalfunktionalen nochmals ein.
Vorerst beschränken wir uns jedoch auf Integralfunktionale und untersuchen im nächsten Abschnitt die Lagrange-Aufgabe. \hfill $\square$}
\end{bemerkung}

%% file: 1-2-LagrangeAufgabe.tex
\subsection{Die Lagrange-Aufgabe und das Isoperimetrische Problem} \label{AbschnittLagrangeAufgabe}
\subsubsection{Notwendige Bedingungen in der Lagrange-Aufgabe}
Als Lagrange-Aufgabe\index{Lagrange-Aufgabe} der Klassischen Variationsrechnung bezeichnen wir das Problem
\begin{eqnarray}
&& \label{LA1} J\big(x(\cdot),u(\cdot)\big) = \int_{t_0}^{t_1} f\big(t,x(t),u(t)\big) \, dt \to \inf, \\
&& \label{LA2} \dot{x}(t) = \varphi\big(t,x(t),u(t)\big), \\
&& \label{LA3} h_0\big(x(t_0)\big)=0, \qquad h_1\big(x(t_1)\big)=0, \\
&& \label{LA4} u(t) \in U=\R^m.
\end{eqnarray}
In der Lagrange-Aufgabe seien 
$\big(x(\cdot),u(\cdot)\big) \in C_1([t_0,t_1],\R^n) \times C([t_0,t_1],\R^m)$ und
$$f:\R \times \R^n \times \R^m \to \R, \quad \varphi:\R \times \R^n \times \R^m \to \R^n, \quad
  h_i:\R^n \to \R^{s_i}, \; i=0,1.$$
Zur Menge $\mathscr{A}_{\rm Lip}^{\mathcal{L}}$ gehören diejenigen
$\big(x(\cdot),u(\cdot)\big) \in  C_1([t_0,t_1],\R^n) \times C([t_0,t_1],\R^m)$,
f"ur die es eine Zahl $\gamma>0$ derart gibt,
dass die Abbildungen $f(t,x,u)$, $\varphi(t,x,u)$, $h_0(x_0),h_1(x_1)$ auf der Menge aller
$(t,x,x_0,x_1,u) \in \R \times \R^n \times \R^n \times \R^n \times \R^m$ mit
$$t_0 \leq t\leq t_1, \quad \|x-x(t)\| < \gamma, \quad \|x_0-x(t_0)\| < \gamma, \quad \|x_1-x(t_1)\| < \gamma, \quad \|u-u(t)\| < \gamma$$
stetig in allen Variablen und stetig differenzierbar bezüglich $x,x_0,x_1,u$ sind. \\[2mm]
Das Paar $\big(x(\cdot),u(\cdot)\big) \in C_1([t_0,t_1],\R^n) \times C([t_0,t_1],\R^m)$
nennen wir zul"assig in der Aufgabe (\ref{LA1})--(\ref{LA4}),
falls $\big(x(\cdot),u(\cdot)\big)$ dem System (\ref{LA2}) gen"ugt und die Randbedingungen (\ref{LA3}) erf"ullt.
Die Menge $\mathscr{A}_{\rm adm}^{\mathcal{L}}$ bezeichnet die Menge der zul"assigen Paare $\big(x(\cdot),u(\cdot)\big)$. \\[2mm]
Ein zul"assiges Paar $\big(x_*(\cdot),u_*(\cdot)\big)$ ist in der Aufgabe (\ref{LA1})--(\ref{LA4})
ein schwaches lokales Minimum\index{Minimum, schwaches lokales!Variation@-- Variationsrechnung},
falls eine Zahl $\varepsilon > 0$ derart existiert, dass die Ungleichung 
$$J\big(x(\cdot),u(\cdot)\big) \geq J\big(x_*(\cdot),u_*(\cdot)\big)$$
f"ur alle $\big(x(\cdot),u(\cdot)\big) \in \mathscr{A}_{\rm adm}^{\mathcal{L}}$ mit 
$\|x(\cdot)-x_*(\cdot)\|_{C_1} < \varepsilon$ und $\|u(\cdot)-u_*(\cdot)\|_\infty < \varepsilon$ gilt. \\[2mm]
F"ur die Aufgabe (\ref{LA1})--(\ref{LA4}) bezeichnet $L: \R \times \R^n \times \R^n \times \R^m \times \R^n \times \R \to \R$ die Funktion
$$L(t,x,\dot{x},u,p,\lambda_0) = \lambda_0 f(t,x,u) + \langle p, \dot{x}-\varphi(t,x,u) \rangle.$$

\begin{theorem} \label{SatzLA}
Es sei $\big(x_*(\cdot),u_*(\cdot)\big) \in \mathscr{A}_{\rm adm}^{\mathcal{L}} \cap \mathscr{A}_{\rm Lip}^{\mathcal{L}}$. 
Ist $\big(x_*(\cdot),u_*(\cdot)\big)$ ein schwaches lokales Minimum der Aufgabe (\ref{LA1})--(\ref{LA4}),
dann existieren nicht gleichzeitig verschwindende Multiplikatoren $\lambda_0 \geq 0$, $l_0 \in \R^{s_0},l_1 \in \R^{s_1}$
und $p(\cdot) \in C_1([t_0,t_1],\R^n)$ derart, dass auf $[t_0,t_1]$
\begin{enumerate}
\item[(a)] die Euler-Lagrangesche Gleichung f"ur die Funktion $L$ bez"uglich $x$,
          \begin{equation}\label{SatzLA1}
          \bigg(-\frac{d}{dt} L_{\dot{x}} + L_x \bigg)\bigg|_{(x_*(t),u_*(t))}=0,
          \end{equation}
\item[(b)] die Randbedingungen
          \begin{equation}\label{SatzLA2} 
          L_{\dot{x}}\big|_{(x_*(t_0),u_*(t_0))}={h_0'}^T\big(x_*(t_0)\big)l_0, \quad
          L_{\dot{x}}\big|_{(x_*(t_1),u_*(t_1))}=-{h_1'}^T\big(x_*(t_1)\big)l_1
          \end{equation}
\item[(c)] und die Euler-Lagrangesche Gleichung f"ur die Funktion $L$ bez"uglich $u$,
          \begin{equation}\label{SatzLA3}
          L_u\big|_{(x_*(t),u_*(t))}=0,
          \end{equation}
\end{enumerate}
erfüllt sind.
\end{theorem}

\begin{bemerkung} \label{BemerkungSatzLA} {\rm
Ausf"uhrlich geschrieben nehmen die Bedingungen (\ref{SatzLA1})--(\ref{SatzLA3}) in Hamiltonscher Form die folgende Gestalt an:
Die Gleichung (\ref{SatzLA1}) stellt eine Differentialgleichung, die sogenannte adjungierte Gleichung,
\index{adjungierte Gleichung!Lagrange@-- Lagrange-Aufgabe}
$$\dot{p}(t)=-\varphi_x^T\big(t,x_*(t),u_*(t)\big) p(t) + \lambda_0 f_x\big(t,x_*(t),u_*(t)\big) \quad
  \mbox{f"ur alle } t \in [t_0,t_1],$$
f"ur den Multiplikator $p(\cdot) \in C_1([t_0,t_1],\R^n)$ dar.
Die Randbedingungen (\ref{SatzLA2}) treten als die Transversalit"atsbedingungen
\index{Transversalitätsbedingungen!Lagrange@-- Lagrange-Aufgabe}
$$p(t_0)= {h_0'}^T\big(x_*(t_0)\big)l_0, \qquad p(t_1)=-{h_1'}^T\big(x_*(t_1)\big)l_1$$
f"ur die Adjungierte $p(\cdot)$ auf.
Schlie"slich ist die Gleichung
$$\lambda_0 f_u\big(t,x_*(t),u_*(t)\big)-\varphi_u^T\big(t,x_*(t),u_*(t)\big) p(t) = 0 \quad \mbox{f"ur alle } t \in [t_0,t_1]$$
die ausf"uhrlich geschriebene Darstellung der Bedingung (\ref{SatzLA3}). \hfill $\square$}
\end{bemerkung}

\begin{bemerkung} \label{BemerkungFreierEndpunkt}
{\rm In der Lagrange-Aufgabe geben die Randbedingungen (\ref{LA3}) Einschr"ankungen der zul"assigen Trajektorie $x(\cdot)$ in den Punkten
$t=t_0$ und $t=t_1$ vor.
Wir schlie"sen dabei den Fall nicht aus, dass diese Einschr"ankungen nicht f"ur alle Koordinaten der Punkte
$x(t_0)$ und $x(t_1)$ gelten.
Sind die Komponenten $x_i(t_0)$, $x_j(t_1)$ f"ur gewisse $i,j \in \{1,...,n\}$ frei,
so gelten die Transversalit"atsbedingungen $p_i(t_0)=0$ und $p_j(t_1)=0$. \\
Liegt die Lagrange-Aufgabe mit freiem Anfangs- bzw. Endpunkt vor,
d.\,h. durch wenigstens eine der Abbildung $h_0(x_0),h_1(x_1)$ sind keine Einschr"ankungen an die Punkte $x(t_0)$ bzw. $x(t_1)$ definiert,
dann ergeben sich die Transversalit"atsbedingungen $p(t_0)=0$ bzw. $p(t_1)=0$.
Wegen der Nichttrivialit"at der Multiplikatoren folgt in diesem Fall aus der adjungierten Gleichung unmittelbar der normale Fall,
d.\,h. $\lambda_0=1$. \hfill $\square$}
\end{bemerkung}

W"ahrend die Euler-Lagrangesche Gleichung (\ref{ELG5}) f"ur die Aufgabe (\ref{ELG1})--(\ref{ELG2}) der Variationsrechnung
das Gegenstück zur Fermatschen Gleichung $J'\big(x_*(\cdot)\big)=0$ darstellt,
findet im Satz \ref{SatzLA} das Lagrangesche Prinzip seine Best"atigung. \\
Dazu versehen wir die Nebenbedingungen mit geeigneten Lagrangeschen Multiplikatoren und bilden die Lagrange-Funktion $\mathscr{L}$,
$$\mathscr{L}:C_1([t_0,t_1],\R^n) \times C([t_0,t_1],\R^m) \times C_1([t_0,t_1],\R^n) \times \R \times \R^{s_0} \times \R^{s_1},$$
mit Hilfe der zur Aufgabe (\ref{LA1})--(\ref{LA4}) definierten Funktion $L$ auf folgende Weise:
\begin{eqnarray*}
\mathscr{L} &=& \mathscr{L}\big(x(\cdot),u(\cdot),p(\cdot),\lambda_0,l_0,l_1\big) \\
&=& \int_{t_0}^{t_1} L\big(t,x(t),\dot{x}(t),u(t),p(t),\lambda_0\big) \, dt +
    \big\langle l_0 , h_0\big(x(t_0)\big) \big\rangle + \big\langle l_1 , h_1\big(x(t_1)\big) \big\rangle.
\end{eqnarray*}
Gem"a"s dem Lagrangeschen Prinzip\index{Lagrangesches Prinzip} betrachten wir die Aufgabe ohne Nebenbedingungen:
$$\inf_{(x(\cdot),u(\cdot))} \mathscr{L}.$$
Halten wir $u_*(\cdot)$ fest, so f"uhrt dies zu einer Bolza-Aufgabe (\ref{ELGBA1}) in Vektorform,
$$\int_{t_0}^{t_1} L\big(t,x(t),\dot{x}(t),u_*(t),p(t),\lambda_0\big) \, dt +
    \big\langle l_0 , h_0\big(x(t_0)\big) \big\rangle + \big\langle l_1 , h_1\big(x(t_1)\big) \big\rangle \to \inf,$$
und (\ref{SatzLA1}), (\ref{SatzLA2}) ergeben sich in voller "Ubereinstimmung mit Satz \ref{SatzELGBA}. \\
Wird in der Lagrange-Funktion $\mathscr{L}$ die Trajektorie $x_*(\cdot)$ festgehalten,
dann ergibt sich bez"uglich $u(\cdot) \in C([t_0,t_1],\R^m)$ die Aufgabe
$$\int_{t_0}^{t_1} g\big(t,u(t)\big) \, dt := \int_{t_0}^{t_1} L\big(t,x_*(t),\dot{x}_*(t),u(t),p(t),\lambda_0\big) \, dt \to \inf.$$
Setzen wir darin $u(t)=\dot{w}(t)$ mit der Funktion $w(\cdot) \in C_1([t_0,t_1],\R^m)$,
dann entsteht die Grundaufgabe mit freien Endpunkten
$$\int_{t_0}^{t_1} g\big(t,\dot{w}(t)\big) \, dt \to \inf.$$
Aufgrund Teil (b) in Folgerung \ref{FolgerungELG} ist
$$p(t)=g_{\dot{w}}\big(t,\dot{w}_*(t)\big)= \mbox{konstant}.$$
Mit Folgerung \ref{FolgerungELGFE} gelten weiterhin
$$p(t_0)= g_{\dot{w}}\big(t_0,\dot{x}_*(t_0)\big)=0, \qquad p(t_1)= g_{\dot{w}}\big(t_1,\dot{x}_*(t_1)\big)=0,$$
und wir erhalten zusammen (\ref{SatzLA3}).

\begin{beispiel} \label{BeispielVorsorgemodell}
{\rm Wir betrachten ein einfaches Vorsorgemodell "uber dem gegebenen Zeitrahmen $[0,T]$ und
wollen das optimale Konsumptionsverhalten charakterisieren.
Dabei bezeichne $K(\cdot)$ den Kapitalstock, $W(\cdot)$ den Lohn f"ur die erbrachte Arbeitsleistung und $C(\cdot)$ die Konsumption.
Weiterhin seien $U$ eine zweimal stetig differenzierbar, monoton wachsende und streng konkave Nutzenfunktion,
d.\,h. $U \in C_2$, $U'>0$, $U''<0$, und es seien $r,\varrho>0$ die Zins- bzw. Diskontrate.
Die Lohnfunktion $W(\cdot)$,
sowie der Kapitalbestand $K_0$ zu Beginn und das angestrebte Verm"ogen $K_T$ zum Zeitpunkt $T$ sind vorgegeben.\\[2mm]
Das Modell lautet mit diesen Gr"o"sen
\begin{eqnarray*}
&& J\big(K(\cdot),C(\cdot)\big) = \int_0^T e^{-\varrho t} U\big(C(t)\big) \, dt \to \sup, \\
&& \dot{K}(t) = rK(t)+W(t)-C(t), \qquad K(0)=K_0, \quad K(T)=K_T, \qquad C(t) \in \R.
\end{eqnarray*}
Als zul"assige Paare fassen wir $\big(K(\cdot),C(\cdot)\big) \in C_1([0,T],\R) \times C([0,T],\R)$ auf. \\[2mm]
Wir wenden die notwendigen Bedingungen (\ref{SatzLA1})--(\ref{SatzLA3}) in der Form nach Bemerkung \ref{BemerkungSatzLA} an:
Die adjungierte Gleichung ergibt 
$$\dot{p}(t)=-rp(t)$$
und die Transversalit"atsbedingungen liefern
$$p(t)=l_0e^{-rt}.$$
Weiterhin erhalten wir aus der Euler-Lagrange-Gleichung bez"uglich $C$ die Bedingung
$$p(t)=\lambda_0 e^{-\varrho t} U'\big(C_*(t)\big).$$
Der Fall $\lambda_0=0$ kann ausgeschlossen werden, da ansonsten s"amtliche Multiplikatoren verschwinden w"urden.
Ferner muss $l_0>0$ gelten.
Das Konsumptionsverhalten reduziert sich damit auf die Auswertung der Bedingung 
$$U'\big(C_*(t)\big)= l_0e^{(\varrho-r)t}.$$
\begin{enumerate}
\item[(a)] $r=\varrho$: In diesem Fall ist $C_*(t)$ "uber $[0,T]$ konstant.
\item[(b)] $r>\varrho$: Da $U$ konkav ist, ist $C_*(t)$ "uber $[0,T]$ monoton wachsend.
\item[(c)] $r<\varrho$: Es ist $C_*(t)$ "uber $[0,T]$ monoton fallend.
\end{enumerate}
Ist daher die Teuerungsrate $\varrho$ kleiner bzw. gr"oßer als die Zinsrate $r$,
so wird in stetig h"oher bzw. geringer werdendem Ma"se konsumiert. 
Das eindeutige Konsumptionsverhalten $C_*(\cdot)$ ergibt sich abschlie"send aus den Randbedingungen
$K_*(0)=K_0$ und $K_*(T)=K_T$. \hfill $\square$}
\end{beispiel}

%% file: 1-21-Beweis.tex
\subsubsection{Der Nachweis der notwendigen Optimalit\"atsbedingungen}
Wir definieren f"ur $\big(x(\cdot),u(\cdot)\big) \in C_1([t_0,t_1],\R^n) \times C([t_0,t_1],\R^m)$ die Abbildungen
\begin{eqnarray*}
J\big(x(\cdot),u(\cdot)\big) &=& \int_{t_0}^{t_1} f\big(t,x(t),u(t)\big) \, dt, \\
F\big(x(\cdot),u(\cdot)\big)(t) &=& \dot{x}(t) - \varphi\big(t,x(t),u(t)\big), \quad t \in [t_0,t_1], \\
H_i\big(x(\cdot)\big) &=& h_i\big(x(t_i)\big), \quad i=0,1.
\end{eqnarray*}
Dabei fassen wir sie als Abbildungen zwischen folgenden Funktionenr"aumen auf:
\begin{eqnarray*}
J &:& C_1([t_0,t_1],\R^n) \times C([t_0,t_1],\R^m) \to \R, \\
F &:& C_1([t_0,t_1],\R^n) \times C([t_0,t_1],\R^m) \to C([t_0,t_1],\R^n), \\
H_i &:& C_1([t_0,t_1],\R^n) \to \R^{s_i}, \quad i=0,1.
\end{eqnarray*}
Wir setzen $\mathscr{F}=(F,H_0,H_1)$ und pr"ufen f"ur die Extremalaufgabe
\begin{equation} \label{LAExtremalaufgabe}
J\big(x(\cdot),u(\cdot)\big) \to \inf, \qquad \mathscr{F}\big(x(\cdot),u(\cdot)\big)=0
\end{equation}
im Punkt $\big(x_*(\cdot),u_*(\cdot)\big) \in \mathscr{A}_{\rm Lip}^{\mathcal{L}}$
die Voraussetzungen des Extremalprinzips \ref{SatzExtremalprinzipGlatt}:
\begin{enumerate}
\item[(A$_1$)] Beachten wir f"ur $x(\cdot) \in C_1([t_0,t_1],\R^n)$ und $u(\cdot) \in C([t_0,t_1],\R^m)$ die Relationen
               $$\|x(\cdot) \|_{L_\infty} \leq \|x(\cdot) \|_{C_1}, \qquad \|u(\cdot) \|_{L_\infty} = \|u(\cdot) \|_\infty,$$
               so ist $J$ nach Beispiel \ref{DiffZielfunktional1}
               im Punkt $\big(x_*(\cdot),u_*(\cdot)\big)$ Fr\'echet-differenzierbar.
\item[(A$_2$)] Die Abbildung $F$ ist die Summe der stetigen linearen Abbildung $x(\cdot) \to \dot{x}(t)$ und
               der Abbildung $\big(x(\cdot),u(\cdot)\big) \to -\varphi\big(t,x(t),u(t)\big)$,
               deren stetige Fr\'echet-Differenzierbarkeit im Beispiel \ref{DiffAbbildung} nachgewiesen ist.
               F"ur die Abbildungen $H_i$ ist die stetige Differenzierbarkeit offensichtlich.
\item[(B)] Wir setzen $A(t)=\varphi_x\big(t,x_*(t),u_*(t)\big)$.
           Dann ist die Surjektivit"at des Operators $F_x\big(x_*(\cdot),u_*(\cdot)\big)$ "aquivalent zur Aussage,
           dass f"ur jedes $y(\cdot) \in C([t_0,t_1],\R^n)$ die Gleichung
           $$\dot{x}(t) - A(t)x(t)=y(t), \qquad t \in [t_0,t_1],$$
           l"osbar ist.            
           Dies besagt gerade Lemma \ref{LemmaDGL5}.
           Somit besitzt der Operator $\mathscr{F}_x\big(x_*(\cdot),u_*(\cdot)\big)$ eine endliche Kodimension
           und ${\rm Im\,}\mathscr{F}'\big(x_*(\cdot),u_*(\cdot)\big)$ ist abgeschlossen.
\end{enumerate}

Zur Extremalaufgabe (\ref{LAExtremalaufgabe}) definieren wir auf
$$C_1([t_0,t_1],\R^n) \times C([t_0,t_1],\R^m) \times \R \times C^*([t_0,t_1],\R^n) \times \R^{s_0} \times \R^{s_1}$$
die Lagrange-Funktion $\mathscr{L}=\mathscr{L}\big(x(\cdot),u(\cdot),\lambda_0,y^*,l_0,l_1\big)$,
$$\mathscr{L}= \lambda_0 J\big(x(\cdot),u(\cdot)\big)+ \big\langle y^*, F\big(x(\cdot),u(\cdot)\big) \big\rangle
                         +l_0^T H_0\big(x(\cdot)\big)+l_1^T H_1\big(x(\cdot)\big).$$
Ist $\big(x_*(\cdot),u_*(\cdot)\big)$ ein lokales Minimum der Aufgabe (\ref{LAExtremalaufgabe}),
dann existieren nach Theorem \ref{SatzExtremalprinzipGlatt}
nicht gleichzeitig verschwindende Lagrangesche Multiplikatoren $\lambda_0 \geq 0$, $l_i \in \R^{s_i}$ und $y^* \in C^*([t_0,t_1],\R^n)$ derart,
dass
\begin{equation}\label{LAExtremalaufgabe1}
0 = \mathscr{L}_x\big(x_*(\cdot),u_*(\cdot),\lambda_0,y^*,l_0,l_1\big), \quad 0 = \mathscr{L}_u\big(x_*(\cdot),u_*(\cdot),\lambda_0,y^*,l_0,l_1\big)
\end{equation}
gelten.
Wir folgern daraus nun die Beziehungen (\ref{SatzLA1})--(\ref{SatzLA3}): \\
Nach dem Rieszschen Darstellungssatz in der Version der Folgerung \ref{FolgerungRiesz}
existiert eine Vektorfunktion $\mu(\cdot)= \big(\mu_1(\cdot),...,\mu_n(\cdot)\big)$,
wobei die $\mu_i(\cdot)$ von beschr"ankter Variation und rechtsseitig stetig sind,
mit
$$\big\langle y^*, F\big(x(\cdot),u(\cdot)\big) \big\rangle
  = \int_{t_0}^{t_1} \Big[\dot{x}(t)- \varphi\big(t,x(t),u(t)\big)\Big]^T \, d\mu(t).$$
Dann folgt aus der ersten Gleichung in (\ref{LAExtremalaufgabe1}):
\begin{eqnarray}
0 &=& \lambda_0 \int_{t_0}^{t_1} \big\langle f_x\big(t,x_*(t),u_*(t)\big),x(t) \big\rangle \, dt
      + \int_{t_0}^{t_1} \Big[\dot{x}(t)- \varphi_x\big(t,x_*(t),u_*(t)\big)x(t)\Big]^T \, d\mu(t) \nonumber \\
  & & \label{BeweisschlussLA1}
     + \big\langle l_0, h_0'\big(x_*(t_0)\big) x(t_0) \big\rangle + \big\langle l_1, h_1'\big(x_*(t_1)\big) x(t_1) \big\rangle
\end{eqnarray}
f"ur alle $x(\cdot) \in C_1([t_0,t_1],\R^n)$.
Wir setzen der "Ubersichtlichkeit halber
$$a(t)= f_x\big(t,x_*(t),u_*(t)\big), \qquad A(t)= \varphi_x\big(t,x_*(t),u_*(t)\big), \qquad \Gamma_i=h_i'\big(x_*(t_i)\big).$$
Unter Verwendung von $\displaystyle x(t)=x(t_0) + \int_{t_0}^t \dot{x}(s) \, ds$ erhalten wir
$$\int_{t_0}^{t_1} \langle \lambda_0 a(t),x(t) \rangle \, dt
  = \int_{t_0}^{t_1} \lambda_0 a^T(t) \cdot \bigg[x(t_0)+\int_{t_0}^t \dot{x}(s) \, ds \bigg] \, dt$$
und es ergibt sich weiter durch Vertauschen der Integrationsreihenfolge  
$$\int_{t_0}^{t_1} \langle \lambda_0 a(t),x(t) \rangle \, dt
  = \int_{t_0}^{t_1} \bigg[\int_t^{t_1} \lambda_0 a(s) \, ds \bigg]^T \dot{x}(t) \, dt +
    \bigg[\int_{t_0}^{t_1} \lambda_0 a(t) \, dt \bigg]^T x(t_0).$$
Bezüglich der Abbildung $A(\cdot)$ verfahren wir auf die gleiche Weise und erhalten
$$\int_{t_0}^{t_1} [A(t)x(t)]^T \, d\mu(t)
  = \int_{t_0}^{t_1} \bigg[\int_t^{t_1} A^T(s) \, d\mu(s) \bigg]^T \dot{x}(t) \, dt
    + \bigg[\int_{t_0}^{t_1} A^T(t) \, d\mu(t) \bigg]^T x(t_0).$$
Mit diesen Beziehung k"onnen wir (\ref{BeweisschlussLA1}) in folgende Form bringen:
\begin{eqnarray}
0 &=& \int_{t_0}^{t_1} \bigg[\int_t^{t_1} \lambda_0 a(s) \, ds
                               - \int_t^{t_1} A^T(s) \, d\mu(s) + \Gamma_1^Tl_1 \bigg]^T \dot{x}(t) \, dt
      + \int_{t_0}^{t_1} [\dot{x}(t)]^T \, d\mu(t) \nonumber \\
  & & \label{BeweisschlussLA2} + \bigg[\Gamma_0^T l_0 + \Gamma_1^T l_1 +\int_{t_0}^{t_1} \lambda_0 a(t) \, dt
              - \int_{t_0}^{t_1} A^T(t) \, d\mu(t) \bigg]^T x(t_0).
\end{eqnarray}
Die rechte Seite in (\ref{BeweisschlussLA2}) ist ein stetiges lineares Funktional im Raum $C_1([t_0,t_1],\R^n)$:
$$\langle x^*, x(\cdot) \rangle= \int_{t_0}^{t_1} [\dot{x}(t)]^T \, d\nu(t) + \langle a,x(t_0) \rangle$$
mit
\begin{eqnarray}
\nu(t) &=& \label{BeweisschlussLA3} \mu(t)+ \int_{t_0}^t \bigg[ \int_s^{t_1} \lambda_0 a(\tau) \, d\tau \bigg] \, ds
                - \int_{t_0}^t \bigg[ \int_s^{t_1} A^T(\tau) \, d\mu(\tau) \bigg] \, ds + \Gamma_1^Tl_1 t, \\
a &=& \label{BeweisschlussLA4}
      \Gamma_0^T l_0 + \Gamma_1^T l_1 +\int_{t_0}^{t_1} \lambda_0 a(t) \, dt - \int_{t_0}^{t_1} A^T(t) \, d\mu(t).
\end{eqnarray}
Wegen der eindeutigen Darstellung (Folgerung \ref{FolgerungRiesz2}) eines stetigen linearen Funktionals im Raum $C_1([t_0,t_1],\R^n)$
und aus der Gleichung $\mathscr{L}_x=0$ ergeben sich
\begin{equation}\label{BeweisschlussLA5} \nu(t) \equiv 0, \qquad a=0.\end{equation}
Aus (\ref{BeweisschlussLA3}) und (\ref{BeweisschlussLA5}) ist ersichtlich,
dass die Vektorfunktion $\mu(\cdot)$ absolutstetig ist.
Setzen wir nun $p(t)=\dot{\mu}(t)$, so folgt mit (\ref{BeweisschlussLA3}) und (\ref{BeweisschlussLA5}):
\begin{equation}\label{BeweisschlussLA6}
p(t)+ \int_t^{t_1} \lambda_0 a(s) \, ds -\int_t^{t_1} A^T(s) p(s) \, d(s) + \Gamma_1^Tl_1=0.
\end{equation}
F"ur $t=t_0$ bzw. f"ur $t=t_1$ liefert diese Gleichung zusammen mit (\ref{BeweisschlussLA4}) die Beziehungen
$$p(t_0)= \Gamma_0^T l_0= {h_0'}^T\big(x_*(t_0)\big)l_0, \qquad p(t_1)= -\Gamma_1^T l_1= -{h_1'}^T\big(x_*(t_1)\big)l_1.$$
Differenzieren wir schlie"slich (\ref{BeweisschlussLA6}),
so erhalten wir die adjungierte Gleichung
$$\dot{p}(t)=-A^T(t) p(t)+\lambda_0 a(t) = -\varphi_x^T\big(t,x_*(t),u_*(t)\big) p(t) + \lambda_0 f_x\big(t,x_*(t),u_*(t)\big).$$
Wegen der Stetigkeit der Daten ergibt sich $p(\cdot) \in C_1([t_0,t_1],\R^n)$.
Damit sind (\ref{SatzLA1}) und (\ref{SatzLA2}) in Theorem \ref{SatzLA} gezeigt. \\[2mm]
Setzen wir in der Lagrange-Funktion $\mathscr{L}$ anstelle von $d\mu(t)$ jetzt $p(t)\,dt$ ein,
so liefert die zweite Gleichung in (\ref{LAExtremalaufgabe1}) f"ur alle $u(\cdot) \in C([t_0,t_1],\R^m)$:
\begin{equation} \label{BeweisschlussLA7}
0=\int_{t_0}^{t_1} \Big[ \lambda_0 f_u\big(t,x_*(t),u_*(t)\big)- \varphi^T_u\big(t,x_*(t),u_*(t)\big)p(t)\Big]^T u(t) \, dt. 
\end{equation}
Die rechte Seite definiert im Raum $C([t_0,t_1],\R^m)$ ein stetiges lineares Funktional.
Nach dem Rieszschen Darstellungssatzes (Satz \ref{SatzRiesz}) folgt hieraus
$$\varphi^T_u\big(t,x_*(t),u_*(t)\big)p(t) = \lambda_0 f_u\big(t,x_*(t),u_*(t)\big).$$
Damit ist (\ref{SatzLA3}) gezeigt und Theorem \ref{SatzLA} vollst"andig bewiesen. \hfill $\blacksquare$

%% file: 1-22-Isoperimetrisch.tex
\subsubsection{Das Isoperimetrische Problem}
Das Isoperimetrische Problem\index{Isoperimetrisches Problem}
der Klassischen Variationsrechnung ist folgende Aufgabe:
\begin{eqnarray}
&& \label{IP1} J\big(x(\cdot)\big) = \int_{t_0}^{t_1} f_0\big(t,x(t),\dot{x}(t)\big) \, dt \to \inf, \\
&& \label{IP2} \int_{t_0}^{t_1} f_i\big(t,x(t),\dot{x}(t)\big) \, dt = \alpha_j, \quad j=1,...,m, \\
&& \label{IP3} h_0\big(x(t_0)\big) = 0, \qquad h_1\big(x(t_1)\big)=0.
\end{eqnarray}
Dabei seien $f_j:\R \times \R^n \times \R^n \to \R$, $j=0,...,m$, und $h_i:\R^n \to \R^{s_i}$, $i=0,1$. \\[2mm]
In der Aufgabe (\ref{IP1})--(\ref{IP3}) sind diejenigen Funktionen $x(\cdot) \in C_1([t_0,t_1],\R^n)$ zul"assig,
die die isoperimetrischen Beschr"ankungen (\ref{IP2}) und die Randbedingungen (\ref{IP3}) erf"ullen.
Die Menge $\mathscr{A}_{\rm adm}^{\mathcal{I}}$ bezeichnet wieder die Menge der zul"assigen Trajektorien $x(\cdot)$.
Ferner ist $\mathscr{A}_{\rm Lip}^{\mathcal{I}}$ die Menge aller $x(\cdot) \in C_1([t_0,t_1],\R^n)$,
f"ur die die Abbildungen $f_j(t,x,\dot{x})$, $j=0,...,m$, und $h_0(x_0),h_1(x_1)$ auf der Menge
der $(t,x,x_0,x_1,\dot{x}) \in \R \times \R^n \times \R^n \times \R^n \times \R^n$ mit
$$t_0 \leq t\leq t_1, \quad \|x-x(t)\| < \gamma, \quad \|x_0-x(t_0)\| < \gamma, \quad \|x_1-x(t_1)\| < \gamma, \quad \|\dot{x}-\dot{x}(t)\| < \gamma$$ 
stetig in allen Variablen und stetig differenzierbar bezüglich $x,x_0,x_1,\dot{x}$ sind. \\[2mm]
In der Aufgabe (\ref{IP1})--(\ref{IP3}) bezeichnet
$L: \R \times \R^n \times \R^n \times \R \times \R^m \to \R$ die Funktion
$$L(t,x,\dot{x},\lambda_0,\lambda) = \lambda_0 f_0(t,x,\dot{x}) + \sum_{j=1}^m \lambda_j f_j(t,x,\dot{x})
  =\sum_{j=0}^m \lambda_j f_j(t,x,\dot{x}).$$

\begin{satz} \label{SatzIP}
Es sei $x_*(\cdot) \in \mathscr{A}_{\rm adm}^{\mathcal{I}} \cap \mathscr{A}_{\rm Lip}^{\mathcal{I}}$. 
Ist $x_*(\cdot)$ ein schwaches lokales Minimum der Aufgabe (\ref{IP1})--(\ref{IP3}),
dann existieren nicht gleichzeitig verschwindende Multiplikatoren $\lambda_j \in \R$, $j=0,...,m$, und $l_i \in \R^{s_i}$, $i=0,1$,
derart, dass f"ur die Funktion $L$ die Euler-Lagrange-Gleichung
\begin{equation}\label{SatzIP1}
\bigg(-\frac{d}{dt} L_{\dot{x}} + L_x \bigg)\bigg|_{x_*(t)} =0
\end{equation}
auf $[t_0,t_1]$ erf"ullt ist und folgende Randbedingungen gelten:
\begin{equation}\label{SatzIP2} 
 L_{\dot{x}}\big|_{x_*(t_0)}={h_0'}^T\big(x_*(t_0)\big)l_0, \qquad L_{\dot{x}}\big|_{x_*(t_1)}=-{h_1'}^T\big(x_*(t_1)\big)l_1.
\end{equation}
\end{satz}

{\bf Beweis} Im Isoperimetrischen Problem (\ref{IP1})--(\ref{IP3}) setzen wir
$$f=(f_1,...,f_m)^T, \qquad \alpha=(\alpha_1,...,\alpha_m)^T, \qquad u=\dot{x}, \qquad \dot{y}=f(t,x,u).$$
Damit ensteht die Lagrange-Aufgabe
\begin{eqnarray*}
&& \tilde{J}\big(x(\cdot),y(\cdot),u(\cdot)\big) = \int_{t_0}^{t_1} f_0\big(t,x(t),u(t)\big) \, dt \to \inf, \\
&& \dot{x}(t) = u(t), \qquad \dot{y}(t)=f\big(t,x(t),u(t)\big), \\
&& h_0\big(x(t_0)\big)=0, \quad h_1\big(x(t_1)\big)=0, \quad y(t_0)=0, \quad y(t_1)-\alpha=0, \quad u(t) \in \R^n.
\end{eqnarray*}
Die dazugeh"orende Funktion $\tilde{L}$ lautet
$$\tilde{L}(t,x,y,\dot{x},\dot{y},u,p,q,\lambda_0) = \lambda_0 f_0(t,x,u) + \langle p, \dot{x}-u \rangle
  + \langle q, \dot{y}-f(t,x,u) \rangle.$$
Betrachten wir zun"achst (\ref{SatzLA1}) bez"uglich $y$,
so erhalten wir f"ur $q(\cdot)\in C_1([t_0,t_1],\R^n)$:
$$-\frac{d}{dt} q(t)=0 \qquad\Rightarrow\qquad
  q(t)\equiv -\lambda=-(\lambda_1,...,\lambda_m)^T \in \R^m \quad\mbox{auf } [t_0,t_1].$$
Damit folgt unmittelbar f"ur  (\ref{SatzLA1}) bez"uglich $x$ die Bedingung
$$\frac{d}{dt} \tilde{L}_{\dot{x}}\big|_{(x_*(t),y_*(t),u_*(t))}
  = \frac{d}{dt} p(t)= \sum_{j=0}^m \lambda_j f_{jx}\big(t,x_*(t),u_*(t))=\tilde{L}_{x}\big|_{(x_*(t),y_*(t),u_*(t))}.$$
Aus den Randbedingungen (\ref{SatzLA2}) ergibt sich 
$$p(t_0)={h_0'}^T\big(x_*(t_0)\big)l_0, \quad p(t_1)=-{h_1'}^T\big(x_*(t_1)\big)l_1, \quad
  q(t_0)=q(t_1)=-\lambda=\tilde{l}_0=-\tilde{l}_1.$$
Abschlie"send folgt mit $\dot{x}_*(t)=u_*(t)$ aus (\ref{SatzLA3}) der Zusammenhang
$$\tilde{L}_u\big|_{(x_*(t),y_*(t),u_*(t))}=  - p(t)+L_{\dot{x}}\big|_{x_*(t)} = 0.$$
Beachten wir in s"amtlichen Bedingungen die Beziehungen
$$\dot{x}_*(t)=u_*(t), \qquad p(t)=L_{\dot{x}}\big|_{x_*(t)},$$
dann ist Satz \ref{SatzIP} aus Satz \ref{SatzLA} vollst"andig abgeleitet. \hfill $\blacksquare$

\begin{beispiel}[Kettenlinie\index{Kettenlinie}] {\rm
Eine homogene schwere Kette der L"ange $l$ ist an ihren Enden in den Punkten $(t_0,x_0)$ und $(t_1,x_1)$ aufgeh"angt.
Welche Form besitzt sie?

\begin{figure}[h]
	\centering
	\fbox{\includegraphics[height=3.5cm]{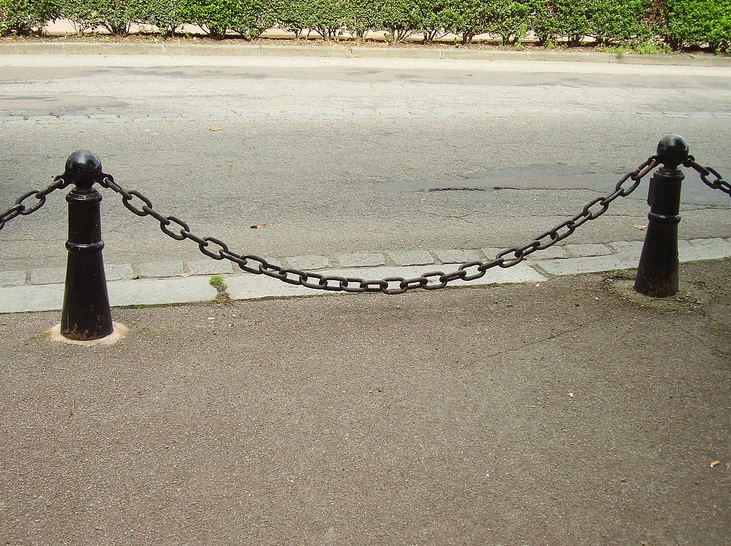}} \hspace*{5mm}
	\fbox{\includegraphics[height=3.5cm]{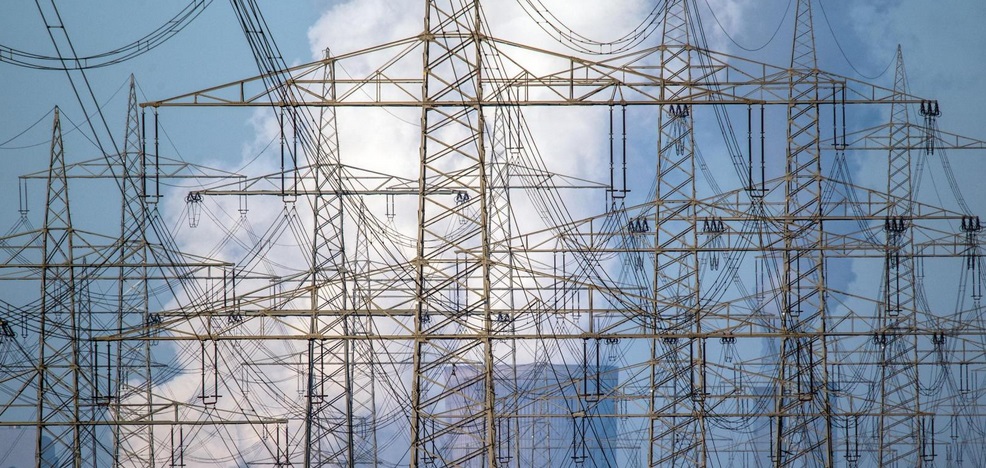}}
	\caption[Kettenlinien im Alltag]{Seilkurven, die durch ihr Eigengewicht belastet sind.}
\end{figure}

\begin{figure}[h]
	\centering
	\fbox{\includegraphics[width=12cm]{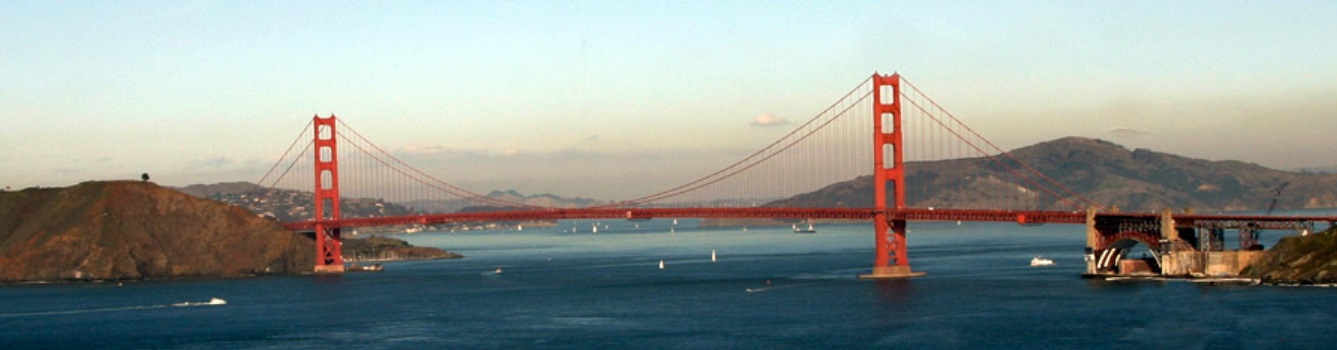}}
	\caption[Kettenlinie an der Golden Gate Bridge]{Golden Gate Bridge.}
\end{figure}

Der Einfachheit halber besitzen die Aufh"angepunkte die gleiche H"ohe $x_0=x_1=h$ und es seien $t_0=-a,t_1=a$ mit einem $a >0$.
Damit die Aufgabe sinnvoll gestellt ist,
d"urfen die Punkte nicht weiter voneinander entfernt sein als die L"ange $l$ der Kette, d.\,h. $l \geq 2a$. \\[2mm]
Die Aufgabe der Kettenlinie f"uhrt zur Minimierung der gesamten potentiellen Energie,
$$J\big(x(\cdot)\big) = \int_{-a}^a x(t) \sqrt{1+\dot{x}^2(t)} \, dt \to \inf,$$
unter der isoperimetrischen Beschr"ankung
$$\int_{-a}^a \sqrt{1+\dot{x}^2(t)} \, dt =l, \qquad l \geq 2a,$$
und den Randbedingungen
$$x(-a)= x(a)=h.$$
Die Euler-Lagrange-Gleichung (\ref{SatzIP1}) liefert f"ur die Funktion
$$L(x,\dot{x},\lambda_0,\lambda) = \lambda_0 x \sqrt{1+\dot{x}^2} + \lambda \sqrt{1+\dot{x}^2},$$
bei Beachtung des Spezialfalles (c) in Folgerung \ref{FolgerungELG}, die Differentialgleichung
$$\lambda_0 x_*(t)+\lambda = \frac{1}{k}\sqrt{1+\dot{x}_*^2(t)}.$$
Ist dabei $\lambda_0=0$, so muss $\dot{x}_*(t)$ konstant sein.
Dies kann nur f"ur die direkte Verbindungsstrecke gelten, also im Fall $l=2a$.
Es sei nun $l>2a$ und $\lambda_0=1$.
Dann erhalten wir aus der Differentialgleichung durch Differentation nach $t$:
$$\ddot{x}_*(t)= k \sqrt{1+\dot{x}_*^2(t)} \qquad\Rightarrow\qquad x_*(t)+\lambda=\frac{1}{k} \cosh(kt+D).$$
Aus physikalischen Gr"unden ist dabei $k>0$.
Wegen der Symmetrie der Randbedingungen gilt ferner $D=0$.
Die "L"ange der Kurve $x_*(t)$ "uber $[-a,a]$ ist gleich $2 \sinh(ka) /k$
und es gibt genau eine positive Konstante $k_0$ mit $2 \sinh(k_0a) /k_0=l$.
Dies ergibt abschlie"send
$$x_*(t)=\frac{1}{k_0} \cosh(k_0t) - \lambda,$$
wobei sich $\lambda$ aus der Randbedingung $x_*(a)=h$ ermitteln l"asst. \hfill $\square$}
\end{beispiel}

%% file: 1-3-Richtungsvariationen.tex
\subsection{Richtungsvariationen in der Optimalen Steuerung} \label{KapitelSchwach}
In diesem Abschnitt widmen wir uns einer direkten Verallgemeinerung der Lagrange-Aufgabe,
indem wir den zul"assigen Paaren die Anforderung
$$\big(x(\cdot),u(\cdot)\big) \in W^1_\infty([t_0,t_1],\R^n) \times L_\infty([t_0,t_1],U)$$
auferlegen.
Dabei erscheint der "Ubergang von der Lagrange-Aufgabe der Klassischen Variationsrechnung zu einem Steuerungsproblem einen rein formalen Charakter zu besitzen,
denn wir behalten den grundlegenden Aufbau der Standardaufgabe:
\begin{eqnarray*}
&& J\big(x(\cdot),u(\cdot)\big) = \displaystyle \int_{t_0}^{t_1} f\big(t,x(t),u(t)\big) \, dt \to \inf, \\
&& \dot{x}(t)= \varphi\big(t,x(t),u(t)\big), \\
&& h_0\big(x(t_0)\big)=0, \qquad h_1\big(x(t_1)\big)=0, \\
&& u(t) \in U, \qquad U=\R^m.
\end{eqnarray*}
Jedoch enthält die Untersuchung der Standardaufgabe auf ein schwaches lokales Minimum im Rahmen der Optimalen Steuerung
signifikante Unterschiede,
die durch die Steuerungsbeschr"ankungen $u(\cdot) \in L_\infty([t_0,t_1],U)$ bzw. $u(t) \in U$ hervorgerufen werden. \\[2mm]
Im Vergleich zur Lagrange-Aufgabe bestehen die wesentlichen Verallgemeinerungen darin,
dass erstens die Steuerung unstetig sein kann und damit die Zustandstrajektorie nicht stetig differenzierbar sein muss.
Sie braucht also nur st"uckweise der Euler-Lagrangeschen Gleichung gen"ugen.
Eine derartige Trajektorie hei"st geknickte Extremale\index{Extremale!--, geknickte}. \\
Aber insbesondere kann zweitens der Fall eines kompakten Steuerungsbereiches vorliegen,
d.\,h. die Steuerungsparameter $u \in U$ d"urfen gewisse Grenzlagen annehmen.
Der Umstand,
dass der Steuerungsbereich $U$ kompakt sein darf,
hebt die nachfolgende Aufgabenklasse aus dem Rahmen der Klassischen Variationsrechnung heraus.
In den Anwendungen sind jedoch die Aufgaben mit kompaktem Steuerungsbereich und mit Steuerungen,
die Werte auf dem Rand des Steuerungsbereiches annehmen, von zentraler Bedeutung. \\[2mm]
Neben der grundlegenden Verallgemeinerung der Aufgabenklasse erweitern wir au"serdem das Spektrum an Methoden und Resultaten.
Einerseits gehen wir n"amlich auf hinreichende Bedingungen nach Mangasarian ein,
die sich in der Standardaufgabe unter zus"atzlichen Konkavit"atsannahmen direkt einbinden lassen.
Ferner zeigen wir,
wie man mit Hilfe der Substitution der Zeit die Aufgaben mit freiem Anfangs- und Endzeitpunkt behandelt.
Zum Abschluss behandeln wir die Aufgaben mit Zustandsbeschr"ankungen,
d.\,h. Aufgaben mit punktweisen Ungleichungsrestriktionen für die Zustandstrajektorie. \\[2mm]
Dieser Abschnitt behandelt die Standardaufgabe der Steuerungstheorie mit Hilfe der Richtungsvariationen.
Allerdings greifen die Beweise der Optimalitätsprinzipien auf grundlegende funktionalanalytische und ma"stheoretische Methoden zurück,
die im Anhang aufgeführt sind.

%% file: 1-31-Aufgabenstellung.tex
\subsubsection{Die Aufgabenstellung f\"ur ein schwaches lokales Minimum}
Wir betrachten als Steuerungsproblem die Aufgabe
\begin{eqnarray}
&& \label{SOP1} J\big(x(\cdot),u(\cdot)\big) = \int_{t_0}^{t_1} f\big(t,x(t),u(t)\big) \, dt \to \inf, \\
&& \label{SOP2} \dot{x}(t) = \varphi\big(t,x(t),u(t)\big), \\
&& \label{SOP3} h_0\big(x(t_0)\big)=0, \qquad h_1\big(x(t_1)\big)=0, \\
&& \label{SOP4} u(t) \in U \subseteq \R^m, \quad U\not= \emptyset, \quad U \mbox{ konvex}, \\
&& \label{SOP5} g_j\big(t,x(t)\big) \leq 0 \quad \mbox{f"ur alle } t \in [t_0,t_1], \quad j=1,...,l.
\end{eqnarray}
Dabei gelten $f:\R \times \R^n \times \R^m \to \R$, $\varphi:\R \times \R^n \times \R^m \to \R^n$, $h_i:\R^n \to \R^{s_i}$
f"ur $i=0,1$, sowie $g_j:\R \times \R^n \to \R$ f"ur $j=1,...,l$ in den Zustandsbeschränkungen (\ref{SOP5}). \\
Die Aufgabe (\ref{SOP1})--(\ref{SOP5}) betrachten wir bez"uglich der Paare
$$\big(x(\cdot),u(\cdot)\big) \in W^1_\infty([t_0,t_1],\R^n) \times L_\infty([t_0,t_1],U).$$
Wir nennen die Trajektorie $x(\cdot)$ eine L"osung der Gleichung (\ref{SOP2}) zur Steuerung $u(\cdot)$,
falls $x(\cdot)$ auf $[t_0,t_1]$ definiert ist und die Dynamik im Sinn von Carath\'eodory l"ost. \\[2mm]
Da die Steuerung lediglich messbar und beschränkt ist,
gestaltet sich im Vergleich zur Lagrange-Aufgabe die Festlegung eines Umgebungsstreifens für einen Steuerungsprozess $\big(x(\cdot),u(\cdot)\big)$ umständlicher.
Deswegen schränken wir die Differenzierbarkeit bezüglich der Steuervariable $u$ nicht auf eine gewisse Umgebung ein. \\[2mm]
Zur Menge $\mathscr{A}^{\mathcal{S}}_{\rm Lip}$ gehören diejenigen
$\big(x(\cdot),u(\cdot)\big) \in W^1_\infty([t_0,t_1],\R^n) \times L_\infty([t_0,t_1],U)$,
f"ur die es eine Zahl $\gamma>0$ derart gibt,
dass die Abbildungen $f(t,x,u)$, $\varphi(t,x,u)$, $h_i(x_i)$ und $g_j(t,x)$ auf der Menge aller
$(t,x,x_0,x_1,u) \in \R \times \R^n \times \R^n \times \R^n \times \R^m$ mit
$$t_0 \leq t\leq t_1, \quad \|x-x(t)\| < \gamma, \quad \|x_0-x(t_0)\| < \gamma, \quad \|x_1-x(t_1)\| < \gamma, \quad u \in \R^m$$
stetig in allen Variablen und stetig differenzierbar bezüglich $x,x_0,x_1,u$ sind. \\[2mm]
Das Paar $\big(x(\cdot),u(\cdot)\big) \in W^1_\infty([t_0,t_1],\R^n) \times L_\infty([t_0,t_1],U)$
hei"st ein zul"assiger Steuerungsprozess in der Aufgabe (\ref{SOP1})--(\ref{SOP5}),
falls $\big(x(\cdot),u(\cdot)\big)$ dem System (\ref{SOP2}) gen"ugt, sowie die Randbedingungen (\ref{SOP3}) und
Zustandsbeschr"ankungen (\ref{SOP5}) erf"ullt.
Die Menge $\mathscr{A}^{\mathcal{S}}_{\rm adm}$ bezeichnet die Menge der zul"assigen Steuerungsprozesse. \\[2mm]
Ein zul"assiger Steuerungsprozess $\big(x_*(\cdot),u_*(\cdot)\big)$ ist eine schwache lokale
Minimalstelle\index{Minimum, schwaches lokales!Optimale@-- Optimale Steuerung}
der Aufgabe (\ref{SOP1})--(\ref{SOP5}),
falls eine Zahl $\varepsilon > 0$ derart existiert, dass die Ungleichung 
$$J\big(x(\cdot),u(\cdot)\big) \geq J\big(x_*(\cdot),u_*(\cdot)\big)$$
f"ur alle $\big(x(\cdot),u(\cdot)\big) \in \mathscr{A}^{\mathcal{S}}_{\rm adm}$ mit 
$\|x(\cdot)-x_*(\cdot)\|_\infty < \varepsilon$ und $\|u(\cdot)-u_*(\cdot)\|_{L_\infty} < \varepsilon$ gilt. \\[2mm]
Ferner bezeichnet $H^{\mathcal{S}}: \R \times \R^n \times \R^m \times \R^n \times \R \to \R$ die
Pontrjagin-Funktion \index{Pontrjagin-Funktion!Standard@-- $H^{\mathcal{S}}$ der Standardaufgabe}
$$H^{\mathcal{S}}(t,x,u,p,\lambda_0) = \langle p, \varphi(t,x,u) \rangle - \lambda_0 f(t,x,u)$$
der Standardaufgabe (\ref{SOP1})--(\ref{SOP5}).

%% file: 1-32-SchwachesPrinzip.tex
\subsubsection{Ein Schwaches Optimalit\"atsprinzip} \label{AbschnittSchwachesPrinzip}
\begin{theorem} \label{SatzSOP} \index{Schwaches Optimalitätsprinzip}
Es sei $\big(x_*(\cdot),u_*(\cdot)\big) \in \mathscr{A}^{\mathcal{S}}_{\rm adm} \cap \mathscr{A}^{\mathcal{S}}_{\rm Lip}$. 
Ist $\big(x_*(\cdot),u_*(\cdot)\big)$ ein schwaches lokales Minimum der Aufgabe (\ref{SOP1})--(\ref{SOP4}),
dann existieren nicht gleichzeitig verschwindende Multiplikatoren $\lambda_0 \geq 0$, $l_0 \in \R^{s_0},l_1 \in \R^{s_1}$
und $p(\cdot) \in W^1_\infty([t_0,t_1],\R^n)$ derart, dass
\begin{enumerate}
\item[(a)] für fast alle $t \in [t_0,t_1]$ die adjungierte Gleichung
           \index{adjungierte Gleichung!Schwach@-- Schwaches Optimalitätsprinzip}
           \begin{equation}\label{SatzSOP1}
           \dot{p}(t)=-H^{\mathcal{S}}_x\big(t,x_*(t),u_*(t),p(t),\lambda_0\big),
           \end{equation}
\item[(b)] die Transversalit"atsbedingungen
           \index{Transversalitätsbedingungen!Schwach@-- Schwaches Optimalitätsprinzip}
           \begin{equation}\label{SatzSOP2} 
           p(t_0)= {h_0'}^T\big(x_*(t_0)\big)l_0, \qquad p(t_1)=-{h_1'}^T\big(x_*(t_1)\big)l_1
           \end{equation}
\item[(c)] und in fast allen Punkten $t \in [t_0,t_1]$ f"ur alle $u \in U$ die Variationsungleichung
           \begin{equation}\label{SatzSOP3}
           \big\langle H^{\mathcal{S}}_u\big(t,x_*(t),u_*(t),p(t),\lambda_0\big),\big(u-u_*(t)\big) \big\rangle\leq 0
           \end{equation}
\end{enumerate}
erfüllt sind.
\end{theorem}

\begin{bemerkung} \label{BemerkungSOP} {\rm
Ohne Einschr"ankung l"asst sich an dieser Stelle die Bemerkung \ref{BemerkungFreierEndpunkt} "uber die Transversalit"atsbedingungen und
"uber das Eintreten des normalen Falles bei fehlenden Randbedingungen "ubernehmen. \hfill $\square$}
\end{bemerkung}

\begin{bemerkung} \label{BemerkungSOP2} {\rm
Die adjungierte Gleichung (\ref{SatzSOP1}) und die Variationsungleichung (\ref{SatzSOP3}) besitzen in Hamiltonscher Form das Aussehen
\begin{eqnarray*}
\dot{p}(t) &=& -\varphi_x^T\big(t,x_*(t),u_*(t)\big) p(t) + \lambda_0 f_x\big(t,x_*(t),u_*(t)\big), \\
0 &\geq& \big\langle \varphi_u^T\big(t,x_*(t),u_*(t)\big) p(t) + \lambda_0 f_u\big(t,x_*(t),u_*(t)\big), u-u_*(t) \big\rangle.
\end{eqnarray*}
Darin erkennt man die Optimalitätsbedingungen der Lagrange-Aufgabe in Theorem \ref{SatzLA},
deren Hamiltonsche Form in Bemerkung \ref{BemerkungSatzLA} angegeben sind.
Wegen dem gewählten Rahmen $x_*(\cdot) \in W^1_\infty([t_0,t_1],\R^n)$ und $u_*(\cdot) \in L_\infty([t_0,t_1],U)$
gelten die Bedingungen (\ref{SatzSOP1}) und (\ref{SatzSOP3}) jedoch nur fast überall.
Ferner führt die Steuerungsbeschränkung $u \in U$ mit der konvexen Menge $U \not= \emptyset$ zur Ungleichung (\ref{SatzSOP3})
anstelle der Gleichung (\ref{SatzLA3}).  \hfill $\square$}
\end{bemerkung}

\begin{beispiel}{\rm Wir suchen eine L"osung des Variationsproblems
$$\int_0^1 \big(\dot{x}^2(t)-1\big)^2 \, dt \to \inf, \qquad x(0)=x(1)=0.$$
Die L"osungen dieser Aufgabe bestehen aus allen ``Zick-Zack-Kurven'',
die nur die Anstiege $\pm 1$ besitzen und die Punkte $(0,0)$ und $(1,0)$ miteinander verbinden.
Zwischen den Stellen, in denen der Anstieg das Vorzeichen wechselt,
gen"ugen diese geknickten Extremalen\index{Extremale!--, geknickte} der Euler-Lagrangeschen Gleichung. \hfill $\square$}
\end{beispiel}

Die nachstehenden Investitionsmodelle sind Seierstad \& Syds\ae ter \cite{Seierstad} entnommen:
\begin{beispiel} \label{BeispielLinInv} \index{Kapitalakkumulation}
{\rm Es bezeichne $K(t)$ das Kapital und
$u(t)$ die Investitions- bzw. $\big(1-u(t)\big)$ die Konsumptionsrate zum Zeitpunkt $t$.
Wir betrachten das lineare Investitionsmodell
\begin{eqnarray*}
&& J\big(K(\cdot),u(\cdot)\big) = \int_0^T \big( 1-u(t)\big) \cdot K(t) \, dt \to \sup, \\
&& \dot{K}(t) = u(t) \cdot K(t), \quad K(0)=K_0 >0, \\
&& u(t) \in [0,1], \quad T > 1 \mbox{ fest}.
\end{eqnarray*}
Den Integranden mulitiplizieren wir mit $-1$ und gehen zu einem Minimierungsproblem "uber.
Wir erhalten damit in der Aufgabe die Abbildungen
$$f(t,K,u) = -( 1-u) \cdot K, \quad \varphi(t,K,u) = u \cdot K.$$
Nach Bemerkung \ref{BemerkungSOP} d"urfen wir Theorem \ref{SatzSOP} mit $\lambda_0=1$ anwenden.
Damit gelten für alle $u \in [0,1]$ die Variationsungleichung
$$\big\langle H^{\mathcal{S}}_u\big(t,x_*(t),u_*(t),p(t),1\big),\big(u-u_*(t)\big) \big\rangle
  = \langle (p(t) -1) \cdot K_*(t), u-u_*(t) \rangle \leq 0$$
und $p(\cdot)$ ist die L"osung der adjungierten Gleichung
$$\dot{p}(t) = - u_*(t) \cdot p(t) - \big( 1 - u_*(t)\big) = -1 - u_*(t) \big( p(t) - 1 \big), \quad p(T)=0.$$
Aus der Variationsungleichung k"onnen wir unmittelbar
$$u_*(t) = 0 \mbox{ f"ur } p(t) < 1 \quad\mbox{und}\quad u_*(t) = 1 \mbox{ f"ur } p(t) > 1$$
entnehmen. 
Für die adjungierte Gleichung folgt damit sofort
$$\dot{p}(t) = \left\{ \begin{array}{ll}
               -1 ,   & p(t) < 1 \\
               -p(t), & p(t) > 1 \end{array} \right\} = - \max \{1, p(t) \} \leq -1 <0 \mbox{ für alle } t \in (0,T).$$
Die Stelle $\tau_1$ mit $p(\tau_1)=1$ ist wegen $\dot{p}(t)<-1$ und $p(T)=0$ eindeutig bestimmt. \\[2mm]
Aus formalen Gründen betrachten wir zur Bestimmung von $\tau_1$ die Funktion $g$ mit
$$\dot{g}(t) = -1 \mbox{ und } g(T) = 0 \quad\Rightarrow\quad g(t)=T-t \mbox{ für } t \leq T.$$
Also gibt es genau eine L"osung $\tau_1 \in (0,T)$ der Gleichung $g(t)=1$ und es gilt $\tau_1 = T-1$. \\[2mm]
Da die adjungierte Funktion $p(\cdot)$ streng monoton fallend ist, folgt f"ur diese
$$\dot{p}(t) = \left\{ \begin{array}{ll}
               -p(t), & t \in (0,T-1), \\
               -1,    & t \in (T-1,T), \end{array} \right.
  \qquad
  p(t) = \left\{ \begin{array}{ll}
               e^{T-(1+t)}, & t \in [0,T-1), \\
               T-t,         & t \in [T-1,T]. \end{array} \right.$$
Mit der Steuerung
$$u_*(t) = \left\{ \begin{array}{ll}
               1 & \mbox{ f"ur } t \in [0,T-1), \\
               0 & \mbox{ f"ur } t \in [T-1,T], \end{array} \right.$$
erhalten wir f"ur das Kapital
$$\dot{K}_*(t) = \left\{ \begin{array}{ll}
                K_*(t), & t \in (0,T-1), \\
                0,     & t \in (T-1,T), \end{array} \right.
  \qquad
  K_*(t) = \left\{ \begin{array}{ll}
               K_0 \cdot e^t     & t \in [0,T-1), \\
               K_0 \cdot e^{T-1} & t \in [T-1,T], \end{array} \right.$$
und weiterhin f"ur den Wert des Zielfunktionals
$$J\big(K_*(\cdot),u_*(\cdot)\big) = \int_0^T \big( 1-u_*(t)\big) \cdot K_*(t) \, dt = \int_{T-1}^{T} K_0 \cdot e^{T-1} \, dt = K_0 \cdot e^{T-1}.$$
In diesem Beispiel ist die optimale Steuerung $u_*(\cdot)$ unstetig und nimmt ausschlie"slich Werte auf dem Rand des 
Steuerungsbereiches $U=[0,1]$ an. \hfill $\square$}
\end{beispiel}

\begin{beispiel} \label{BeispielKonInv}
{\rm Im Gegensatz zum vorherigen linearen Modell wird die Produktion durch die Cobb-Douglas-Funktion $f(K)=K^\alpha$ mit
$\alpha \in (0,1)$ angegeben.
Damit erhalten wir folgende Aufgabe 
\begin{eqnarray*}
&& J\big(K(\cdot),u(\cdot)\big) = \int_0^T \big( 1-u(t)\big) \cdot K^\alpha(t) \, dt \to \sup, \\
&& \dot{K}(t) = u(t) \cdot K^\alpha(t), \quad K(0)=K_0 >0, \\
&& u(t) \in [0,1], \quad \alpha \in (0,1) \mbox{ konstant}, \quad T \mbox{ fest mit } \alpha T -K_0^{1-\alpha} > 0.
\end{eqnarray*}
Der Term $K^\alpha$ f"uhrt zur Einschr"ankung $\gamma<K_0$ bei der Festlegung des Umgebungsstreifens $|x-x_*(t)| < \gamma$
einer zulässigen Zustandstrajektorie.
Wir multiplizieren den Integranden mit $-1$ und wenden Theorem \ref{SatzSOP} mit $\lambda_0=1$ an:
Es gelten für alle $u \in [0,1]$ die Variationsungleichung
$$\big\langle H^{\mathcal{S}}_u\big(t,x_*(t),u_*(t),p(t),1\big),\big(u-u_*(t)\big) \big\rangle
  = \langle (p(t) -1) \cdot K^\alpha_*(t), u-u_*(t) \rangle \leq 0$$
und $p(\cdot)$ ist die L"osung der adjungierten Gleichung
$$\dot{p}(t) = - \alpha K_*^{\alpha-1}(t) \big[ 1 + u_*(t) \big( p(t) - 1 \big) \big], \quad p(T)=0.$$
Aus der Variationsungleichung ergibt sich unmittelbar
$$u_*(t) = 0 \mbox{ f"ur } p(t) < 1 \quad\mbox{und}\quad u_*(t) = 1 \mbox{ f"ur } p(t) > 1.$$
Für die adjungierte Gleichung folgt weiterhin
$$\dot{p}(t) = \left\{ \begin{array}{ll}
               -\alpha \, K_*^{\alpha-1}(t),            & p(t) < 1 \\
               -\alpha \, K_*^{\alpha-1}(t) \cdot p(t), & p(t) > 1 \end{array} \right\}
             = -\alpha \, K_*^{\alpha-1}(t) \cdot \max \{ 1,p(t) \}<0 \mbox{ über } (0,T).$$
Die Adjungierte $p(\cdot)$ ist stetig und streng monoton fallend über $[0,T]$ mit $p(T)=0$. \\[2mm]
Wir zeigen nun die Existenz einer Zahl $\tau_2 \in (0,T)$ mit $p(\tau_2)=1$.
Es sei $\tau_2 \in [0,T]$ mit
$$\tau_2 = \inf \{ t \in [0,T] \,|\, p(t)<1 \}.$$
F"ur $t > \tau_2$ gilt $u_*(t)=0$ und es ist $K_*(t)$ auf $[\tau_2,T]$ konstant.
Es bezeichne $K=K_*(\tau_2)$ diesen konstanten Wert.
Anstelle der Adjungieren $p(\cdot)$ betrachten wir mit Hilfe der Konstanten $K$ zun"achst die Funktion $g$ mit
$$\dot{g}(t) = -\alpha \, K^{\alpha-1}, \quad g(T) = 0, \quad \alpha \in (0,1), \quad K>0 \mbox{ konstant}.$$
Es ergibt sich $g(t) = \alpha K^{\alpha-1} (T-t)$.
Nun bestimmen wir $\tau_2$.
F"ur $t < \tau_2$ ist $p(t) >1$ und damit $u_*(t) = 1$.
In der Dynamik der Aufgabe sei $u(t) \equiv 1$.
Es ergibt sich die Funktion $h(\cdot)$ mit $\dot{h}(t) = h^\alpha(t)$, $\alpha \in (0,1)$ und $h(0) = K_0$.
Wir erhalten nach Trennung der Ver"anderlichen und anschlie"sender Integration
$$\int_{K_0}^h x^{-\alpha}\, dx = \int_0^t 1 \, ds \quad\Rightarrow\quad 
  \frac{1}{1-\alpha}\big[h^{1-\alpha}- K_0^{1-\alpha}\big] = t\quad\Rightarrow\quad 
  h(t) = \big[ (1-\alpha) t + K_0^{1-\alpha} \big]^{\frac{1}{1-\alpha}}.$$
Aus $K=h(\tau_2)$ und $g(\tau_2)=1$ ergibt sich f"ur $\tau_2$ die Gleichung
$$1 = \underbrace{\alpha K^{\alpha-1} (T-\tau_2)}_{=g(\tau_2)}
    = \alpha \cdot \underbrace{\big[ (1-\alpha) \tau_2 + K_0^{1-\alpha} \big]^{\frac{\alpha-1}{1-\alpha}}}_{=h^{\alpha-1}(\tau_2)} \cdot (T-\tau_2)
    = \frac{\alpha (T-\tau_2)}{(1-\alpha) \tau_2 + K_0^{1-\alpha}}.$$
Aus dieser Gleichung folgt zusammen mit der Annahme $\alpha T -x_0^{1-\alpha} > 0$ der Zeitpunkt $\tau_2 = \alpha T - K_0^{1-\alpha}>0$. \\[2mm]
Wir  erhalten als optimalen Steuerungsprozess
$$u_*(t) = \left\{ \begin{array}{ll}
               1 & \mbox{ f"ur } t \in [0,\tau_2), \\
               0 & \mbox{ f"ur } t \in [\tau_2,T], \end{array} \right. \quad
K_*(t) = \left\{ \begin{array}{ll}
               \big[ (1-\alpha)t +  K_0^{1-\alpha} \big]^\frac{1}{1-\alpha}            & t \in [0,\tau_2), \\
               \big[ \alpha(1-\alpha)T + \alpha K_0^{1-\alpha} \big]^\frac{1}{1-\alpha} & t \in [\tau_2,T], \end{array} \right.$$
mit dem Umschaltpunkt
$$\tau_2 = \alpha T - K_0^{1-\alpha}.$$
Die adjungierte Funktion lautet dabei
$$p(t)=\left\{ \begin{array}{ll}  
       \displaystyle \left[ \frac{\alpha(1-\alpha)T + \alpha K_0^{1-\alpha}}{(1-\alpha)t + K_0^{1-\alpha}} \right]^\frac{\alpha}{1-\alpha}
                                                        & \mbox{ f"ur } t \in [0,\tau_2), \\[3mm]
       \displaystyle \frac{T-t}{(1-\alpha)T + K_0^{1-\alpha}} & \mbox{ f"ur } t \in [\tau_2,T]. \end{array} \right.$$
Im Beispiel mit Cobb-Douglas-Produktionsfunktion ergibt sich
$$J\big(x_*(\cdot),u_*(\cdot)\big) = \alpha^{\frac{\alpha}{1-\alpha}} \cdot \big[ (1-\alpha)T + K_0^{1-\alpha} \big]^\frac{1}{1-\alpha}$$
f"ur den Wert des Zielfunktionals. \hfill $\square$}
\end{beispiel}

%% file: 1-33-Beweis.tex
\subsubsection{Der Nachweis der notwendigen Optimalit\"atsbedingungen} \label{AbschnittBeweisSOP}
Wir betrachten f"ur $\big(x(\cdot),u(\cdot)\big) \in C([t_0,t_1],\R^n) \times L_\infty([t_0,t_1],\R^m)$ die Abbildungen
\begin{eqnarray*}
J\big(x(\cdot),u(\cdot)\big) &=& \int_{t_0}^{t_1} f\big(t,x(t),u(t)\big) \, dt, \\
F\big(x(\cdot),u(\cdot)\big)(t) &=& x(t) -x(t_0) -\int_{t_0}^t \varphi\big(s,x(s),u(s)\big) \, ds, \quad t \in [t_0,t_1],\\
H_i\big(x(\cdot)\big) &=& h_i\big(x(t_i)\big), \quad i=0,1.
\end{eqnarray*}
Da $x(\cdot)$ zu $C([t_0,t_1],\R^n)$ geh"ort, gilt f"ur diese Abbildungen
\begin{eqnarray*}
J &:& C([t_0,t_1],\R^n) \times L_\infty([t_0,t_1],\R^m) \to \R, \\
F &:& C([t_0,t_1],\R^n) \times L_\infty([t_0,t_1],\R^m) \to C_0([t_0,t_1],\R^n), \\
H_i &:& C([t_0,t_1],\R^n) \to \R^{s_i}, \quad i=0,1.
\end{eqnarray*}

Wir setzen $\mathscr{F}=(F,H_0,H_1)$ und pr"ufen f"ur die Extremalaufgabe
\begin{equation} \label{ExtremalaufgabeSOP}
J\big(x(\cdot),u(\cdot)\big) \to \inf, \qquad \mathscr{F}\big(x(\cdot),u(\cdot)\big)=0, \qquad
u(\cdot) \in L_\infty([t_0,t_1],U)
\end{equation}
im $\big(x_*(\cdot),u_*(\cdot)\big) \in \mathscr{A}^{\mathcal{S}}_{\rm Lip}$,
wobei wir $x_*(\cdot)$ als Element des Raumes $C([t_0,t_1],\R^n)$ auffassen,
die Voraussetzungen von Theorem \ref{SatzExtremalprinzipSchwach}:
\begin{enumerate}
\item[(A$_1$)] Beachten wir f"ur $x(\cdot) \in C([t_0,t_1],\R^n)$ die Gleichung
               $\|x(\cdot) \|_{L_\infty} = \|x(\cdot) \|_\infty$,
               so ist $J$ nach Beispiel \ref{DiffZielfunktional1}
               im Punkt $\big(x_*(\cdot),u_*(\cdot)\big)$ Fr\'echet-differenzierbar.
\item[(A$_2$)] Die Abbildung $F$ ist die Summe der Abbildung $x(\cdot) \to x(t)$ und der Abbildung
               $$\big(x(\cdot),u(\cdot)\big) \to -\int_{t_0}^t \varphi\big(s,x(s),u(s)\big) \, ds,$$
               deren stetige Fr\'echet-Differenzierbarkeit im Beispiel \ref{DiffDynamik1} nachgewiesen ist.
               F"ur die Abbildungen $H_i$ ist die stetige Differenzierbarkeit offensichtlich.
\item[(B)] Wir setzen $A(t)=\varphi_x\big(t,x_*(t),u_*(t)\big)$.
           Dann ergibt sich die Surjektivit"at des Operators $F_x\big(x_*(\cdot),u_*(\cdot)\big)$ aus der Aussage,
           dass f"ur jedes $y(\cdot) \in C([t_0,t_1],\R^n)$ die Gleichung
           $$x(t) - \int_{t_0}^t A(s)x(s) \, ds =y(t), \qquad t \in [t_0,t_1],$$
           l"osbar ist.            
           Dies besagen Lemma \ref{LemmaDGL3} und Bemerkung \ref{BemDGL}.
           Somit besitzt der Operator $\mathscr{F}_x\big(x_*(\cdot),u_*(\cdot)\big)$ eine endliche Kodimension.
\end{enumerate}

Zur Extremalaufgabe (\ref{ExtremalaufgabeSOP}) definieren wir auf
$$C([t_0,t_1],\R^n) \times L_\infty([t_0,t_1],\R^m) \times \R \times C_0^*([t_0,t_1],\R^n) \times \R^{s_0} \times \R^{s_1}$$
die Lagrange-Funktion $\mathscr{L}=\mathscr{L}\big(x(\cdot),u(\cdot),\lambda_0,y^*,l_0,l_1\big)$,
$$\mathscr{L}= \lambda_0 J\big(x(\cdot),u(\cdot)\big)+ \big\langle y^*, F\big(x(\cdot),u(\cdot)\big) \big\rangle
                         +l_0^T H_0\big(x(\cdot)\big)+l_1^T H_1\big(x(\cdot)\big).$$
Ist $\big(x_*(\cdot),u_*(\cdot)\big)$ eine schwache lokale Minimalstelle der Aufgabe (\ref{ExtremalaufgabeSOP}),
dann existieren nach Theorem \ref{SatzExtremalprinzipSchwach}
nicht gleichzeitig verschwindende Lagrangesche Multiplikatoren $\lambda_0 \geq 0$, $y^* \in C_0^*([t_0,t_1],\R^n)$ und $l_i \in \R^{s_i}$ derart,
dass gelten:
\begin{enumerate}
\item[(a)] Die Lagrange-Funktion besitzt bez"uglich $x(\cdot)$ in $x_*(\cdot)$ einen station"aren Punkt, d.\,h.
          \begin{equation}\label{SatzSOPLMR1}
          \mathscr{L}_x\big(x_*(\cdot),u_*(\cdot),\lambda_0,y^*,l_0,l_1\big)=0;
          \end{equation}         
\item[(b)] Die Lagrange-Funktion erf"ullt bez"uglich $u(\cdot)$ in $u_*(\cdot)$ die Variationsungleichung
          \begin{equation}\label{SatzSOPLMR2}
           \big\langle \mathscr{L}_u\big(x_*(\cdot),u_*(\cdot),\lambda_0,y^*,l_0,l_1\big), u(\cdot)-u_*(\cdot) \big\rangle \geq 0
          \end{equation}
          f"ur alle $u(\cdot) \in L_\infty([t_0,t_1],U)$.
\end{enumerate}
Zu $y^* \in C_0^*([t_0,t_1],\R^n)$ gibt es ein eindeutig bestimmtes regul"ares Borelsches Vektorma"s $\mu$ derart,
dass nach (\ref{SatzSOPLMR1}) f"ur alle $x(\cdot) \in C([t_0,t_1],\R^n)$ die Variationsgleichung
\begin{eqnarray}
0 &=& \lambda_0 \int_{t_0}^{t_1} \big\langle f_x\big(t,x_*(t),u_*(t)\big),x(t) \big\rangle \, dt
      + \big\langle l_0, h_0'\big(x_*(t_0)\big) x(t_0) \big\rangle + \big\langle l_1, h_1'\big(x_*(t_1)\big) x(t_1) \big\rangle
      \nonumber \\
  & & \label{BeweisschlussSOP1}
      + \int_{t_0}^{t_1} \bigg[ x(t)-x(t_0) - \int_{t_0}^{t} \varphi_x\big(s,x_*(s),u_*(s)\big) x(s)\, ds \bigg]^T d\mu(t)
\end{eqnarray}
erf"ullt ist.
Der K"urze halber seien $a(t)=f_x\big(t,x_*(t),u_*(t)\big)$ und $A(t)= \varphi_x\big(t,x_*(t),u_*(t)\big)$.
Durch vertauschen der Integrationsreihenfolge,
$$\int_{t_0}^{t_1} \bigg[\int_{t_0}^t A(s)x(s) \, ds \bigg]^T \, d\mu(t)
  = \int_{t_0}^{t_1} \bigg[\int_t^{t_1} d\mu(s) \bigg]^T [A(t)x(t)] \, dt,$$
bringen wir (\ref{BeweisschlussSOP1}) in die Gestalt
\begin{eqnarray}
0 &=& \int_{t_0}^{t_1} \big\langle \lambda_0 a(t) - A^T(t) \int_{t}^{t_1} d\mu(s) , x(t) \big\rangle \, dt + \int_{t_0}^{t_1} [x(t)]^T \, d\mu(t) \nonumber \\
  & & \label{BeweisschlussSOP2} 
      + \Big\langle {h_0'}^T\big(x_*(t_0)\big)l_0 - \int_{t_0}^{t_1} d\mu(t) , x(t_0) \Big\rangle
      + \langle {h_1'}^T\big(x_*(t_1)\big)l_1 , x(t_1) \rangle.
\end{eqnarray}
Setzen wir $p(t)=\displaystyle \int_t^{t_1} d\mu(s)$ in (\ref{BeweisschlussSOP2}),
so folgen wegen der eindeutigen Darstellung eines stetigen linearen Funktionals im Raum $C([t_0,t_1],\R^n)$
\begin{eqnarray*}
p(t) &=& -{h_1'}^T\big(x_*(t_1)\big)l_1
         + \int_{t}^{t_1} \Big( \varphi_x^T\big(s,x_*(s),u_*(s)\big) p(s) - \lambda_0 f_x\big(s,x_*(s),u_*(s)\big) \Big) \, ds, \\
p(t_0) &=& {h_0'}^T\big(x_*(t_0)\big)l_0, \qquad p(t_1)= -{h_1'}^T\big(x_*(t_1)\big)l_1.
\end{eqnarray*}
Die Vektorfunktion $p(\cdot)$ ist der Setzung nach "uber $[t_0,t_1]$ messbar, beschr"ankt und besitzt die verallgemeinerte Ableitung
$$\dot{p}(t)= -\varphi_x^T\big(t,x_*(t),u_*(t)\big) p(t) + \lambda_0 f_x\big(t,x_*(t),u_*(t)\big)
            = - H^{\mathcal{S}}_x\big(t,x_*(t),u_*(t),p(t),\lambda_0\big).$$
Wegen der wesentlichen Beschr"anktheit der Abbildungen
$$t \to f_x\big(t,x_*(t),u_*(t)\big), \qquad t \to \varphi_x\big(t,x_*(t),u_*(t)\big)$$
geh"ort $\dot{p}(\cdot)$ dem Raum $L_\infty([t_0,t_1],\R^n)$ an.
Damit sind (\ref{SatzSOP1}) und (\ref{SatzSOP2}) gezeigt. \\[2mm]
Gem"a"s (\ref{SatzSOPLMR2}) gilt f"ur alle $u(\cdot) \in L_\infty([t_0,t_1],U)$ die Ungleichung
$$\int_{t_0}^{t_1} \big\langle H^{\mathcal{S}}_u\big(t,x_*(t),u_*(t),p(t),\lambda_0\big),u(t)-u_*(t) \big\rangle \, dt \leq 0.$$
Es seien $u \in U$, $\lambda>0$ und $\tau \in (t_0,t_1)$ ein Lebesguescher Punkt der Abbildungen $t \to u_*(t)$
und $t \to H^{\mathcal{S}}_u\big(t,x_*(t),u_*(t),p(t),\lambda_0\big)$.
Wir setzen f"ur $0<\lambda<t_1-\tau$
$$u_\lambda(t)=\left\{\begin{array}{ll} u,& t \in [\tau,\tau+\lambda], \\[1mm]
                                u_*(t), & t \not\in [\tau,\tau+\lambda]. \end{array}\right.$$
Dann gilt $u_\lambda(\cdot) \in L_\infty(\R_+,U)$.
Wir erhalten damit
\begin{eqnarray*}
\lefteqn{\frac{1}{\lambda}
               \int_{t_0}^{t_1} \big\langle H^{\mathcal{S}}_u\big(t,x_*(t),u_*(t),p(t),\lambda_0\big),u_\lambda(t)-u_*(t) \big\rangle\, dt}\\
&=& \frac{1}{\lambda} \int_\tau^{\tau+\lambda} \big\langle H^{\mathcal{S}}_u\big(t,x_*(t),u_*(t),p(t),\lambda_0\big),   
                  \big(u-u_*(t)\big)\big\rangle \, dt \leq 0.
\end{eqnarray*}
Da $\tau \in [t_0,t_1)$ ein beliebiger Lebesguescher Punkt ist, folgt daraus f"ur fast alle $t \in [t_0,t_1]$:
$$\big\langle H^{\mathcal{S}}_u\big(t,x_*(t),u_*(t),p(t),\lambda_0\big),\big(u-u_*(t)\big)\big\rangle \leq 0.$$
Ferner, weil $u \in U$ willk"urlich gew"ahlt war, gilt diese Ungleichung f"ur alle $u \in U$.
Damit ist (\ref{SatzSOP3}) gezeigt und der Beweis von Theorem \ref{SatzSOP} ist abgeschlossen. \hfill $\blacksquare$

%% file: 1-34-Hinreichend.tex
\subsubsection{Hinreichende Bedingungen nach Mangasarian} \label{AbschnittMangasarianSOP}
Die\index{hinreichende Bedingungen, Mangasarian} Darstellung der hinreichenden Bedingungen nach Mangasarian \cite{Mangasarian}
f"ur das folgende Steuerungsproblem ist Seierstad \& Syds\ae ter \cite{Seierstad} entnommen:
\begin{eqnarray}
&& \label{HBSOP1} J\big(x(\cdot),u(\cdot)\big) = \int_{t_0}^{t_1} f\big(t,x(t),u(t)\big) \, dt \to \inf, \\
&& \label{HBSOP2} \dot{x}(t) = \varphi\big(t,x(t),u(t)\big), \\
&& \label{HBSOP3} x(t_0)=x_0, \qquad x(t_1)=x_1, \\
&& \label{HBSOP4} u(t) \in U \subseteq \R^m, \quad U\not= \emptyset, \quad U \mbox{ konvex}.
\end{eqnarray}

Wir schlie"sen in der Aufgabenstellung nicht aus,
dass durch die Randbedingungen (\ref{HBSOP3}) gewisse Komponenten der Punkte $x_0$ und $x_1$ nicht eingeschr"ankt werden.
Genauer formuliert sind folgende Varianten von Randbedingungen enthalten:
\begin{enumerate}
\item[(1)] freie Randwerte, falls $x(t_0)$ und $x(t_1)$ ohne Einschr"ankung sind;
\item[(2)] feste Randwerte, falls $x(t_0)=x_0$ und $x(t_1)=x_1$ mit gegebenen $x_0,x_1 \in \R^n$ gelten;
\item[(3)] gemischt freie und feste Randwerte, falls die einzelnen Komponenten von $x(t_0)$ und $x(t_1)$ mit festen bzw. freien Randwerten vorliegen.
\end{enumerate}

Ferner setzen wir $U_\gamma(t)=\{ (x,u) \in \R^n \times \R^m \,|\, \|x-x_*(t)\| < \gamma,\; \|u-u_*(t)\| < \gamma\}$.

\begin{theorem} \label{SatzHBSOP}
Sei $\big(x_*(\cdot),u_*(\cdot)\big) \in \mathscr{A}^{\mathcal{S}}_{\rm Lip} \cap \mathscr{A}^{\mathcal{S}}_{\rm adm}$. Ferner gelte:
\begin{enumerate}
\item[(a)] Das Tripel $\big(x_*(\cdot),u_*(\cdot),p(\cdot)\big)$
           erf"ullt (\ref{SatzSOP1})--(\ref{SatzSOP3}) in Theorem \ref{SatzSOP} mit $\lambda_0=1$.        
\item[(b)] F"ur jedes $t \in [t_0,t_1]$ ist $H^{\mathcal{S}}\big(t,x,u,p(t),1\big)$ konkav bez"uglich $(x,u)$ auf $U_\gamma(t)$.
\end{enumerate}
Dann ist $\big(x_*(\cdot),u_*(\cdot)\big)$ ein schwaches lokales Minimum der Aufgabe (\ref{HBSOP1})--(\ref{HBSOP4}).
\end{theorem}

{\bf Beweis} Es sei $\big(x(\cdot),u(\cdot)\big) \in \mathscr{A}^{\mathcal{S}}_{\rm adm}$ mit
$\|x(\cdot)-x_*(\cdot)\|_\infty < \gamma$, $\|u(\cdot)-u_*(\cdot)\|_{L_\infty} < \gamma$.
Wir betrachten die Differenz
\begin{eqnarray*}
    \Delta
&=& \int_{t_0}^{t_1} \big[f\big(t,x(t),u(t)\big)-f\big(t,x_*(t),u_*(t)\big)\big] dt \\
&=& \int_{t_0}^{t_1} \big[H^{\mathcal{S}}\big(t,x_*(t),u_*(t),p(t),1\big)-H^{\mathcal{S}}\big(t,x(t),u(t),p(t),1\big)\big] \, dt \\
& & \hspace*{10mm}    + \int_{t_0}^{t_1} \langle p(t), \dot{x}(t)-\dot{x}_*(t) \rangle \, dt.
\end{eqnarray*}
Aus den elementaren Eigenschaften konkaver Funktionen folgt
\begin{eqnarray*}
\lefteqn{H^{\mathcal{S}}\big(t,x_*(t),u_*(t),p(t),1\big)-H^{\mathcal{S}}\big(t,x(t),u(t),p(t),1\big)} \\
&\geq& -\big\langle H^{\mathcal{S}}_x\big(t,x_*(t),u_*(t),p(t),1\big) , x(t)-x_*(t) \big\rangle \\
& &  -\big\langle H^{\mathcal{S}}_u\big(t,x_*(t),u_*(t),p(t),1\big) , u(t)-u_*(t) \big\rangle.
\end{eqnarray*}
Mit der adjungierten Gleichung (\ref{SatzSOP1}) und der Ungleichung (\ref{SatzSOP3}) ergibt sich daraus
$$H^{\mathcal{S}}\big(t,x_*(t),u_*(t),p(t),1\big)-H^{\mathcal{S}}\big(t,x(t),u(t),p(t),1\big) \geq \langle \dot{p}(t),x(t)-x_*(t)\rangle.$$
Somit gilt die Beziehung
\begin{eqnarray*}
\Delta &\geq&  \int_{t_0}^{t_1} \big[\langle \dot{p}(t),x(t)-x_*(t)\rangle + \langle p(t), \dot{x}(t)-\dot{x}_*(t) \rangle \big] \, dt \\
       &=& \langle p(t_1),x(t_1)-x_*(t_1)\rangle-\langle p(t_0),x(t_0)-x_*(t_0)\rangle.
\end{eqnarray*}
Im Fall fester Endbedingungen verschwinden die Differenzen $x(t_i)-x_*(t_i)$.
Sind jedoch gewisse Komponenten im Anfangs- oder Endpunkt $x_0$ bzw. $x_1$ frei,
dann liefern die Transversalit"atsbedingungen,
dass die entsprechenden Komponenten der Adjungierten $p(\cdot)$ zum Zeitpunkt $t_0$ bzw. $t_1$ verschwinden.
Daher folgt die Beziehung $\Delta \geq 0$ f"ur alle zul"assigen $\big(x(\cdot),u(\cdot)\big)$ mit 
$\|x(\cdot)-x_*(\cdot)\|_\infty < \gamma$ und $\|u(\cdot)-u_*(\cdot)\|_{L_\infty} < \gamma$. \hfill $\blacksquare$

\begin{beispiel} {\rm Das Beispiel \ref{BeispielELGRegler} "uberf"uhren wir in die Form
\begin{eqnarray*}
&& J\big(x(\cdot),u(\cdot)\big) = \int_0^T \big[\big(x(t)-1\big)^2+u^2(t)\big] \, dt \to \inf, \\
&& \dot{x}(t)=u(t), \qquad x(0)=0,\quad x(T)=2, \qquad u(t)\in \R, \qquad T>0.
\end{eqnarray*}
Wir bilden die zugeh"orige Pontrjagin-Funktion: $\quad H^{\mathcal{S}}(t,x,u,p,1)=pu-(x-1)^2-u^2$. \\
Sie ist offenbar streng konkav in $(x,u)$.
Daher liefert die Extremale
$$x_*(t)=1+\frac{1+e^{-T}}{e^T-e^{-T}}e^t- \frac{1+e^T}{e^T-e^{-T}}e^{-t}$$
ein schwaches lokales Minimum. \hfill $\square$}
\end{beispiel}

\begin{beispiel}
{\rm Im einfachen Vorsorgemodell \ref{BeispielVorsorgemodell} mit konkaver Nutzenfunktion $U$,
\begin{eqnarray*}
&& J\big(K(\cdot),C(\cdot)\big) = \int_0^T e^{-\varrho t} U\big(C(t)\big) \, dt \to \sup, \\
&& \dot{K}(t) = rK(t)+W(t)-C(t), \qquad K(0)=K_0, \quad K(T)=K_T, \qquad C(t) \in \R,
\end{eqnarray*}
wird die optimale Konsumption $C_*(\cdot)$ durch die Bedingung $U'\big(C_*(t)\big)= l_0e^{(\varrho-r)t}$
vorgeschlagen.
Die Pontrjagin-Funktion
$H^{\mathcal{S}}(t,K,C,p,1)=p(rK+W(t)-C)+ e^{-\varrho t} U(C)$
konkav in $(K,C)$.
Daher ist der eindeutig bestimmte Steuerungsprozess $\big(K_*(\cdot),C_*(\cdot)\big)$,
der die Bedingungen
$$U'\big(C_*(t)\big)= l_0e^{(\varrho-r)t}$$
und
$$\dot{K}_*(t) = rK_*(t)+W(t)-C_*(t), \qquad K_*(0)=K_0,\; K_*(T)=K_T,$$
erf"ullt, ein schwaches lokales Maximum. \hfill $\square$}
\end{beispiel}

%% file: 1-35-FreieZeit.tex
\subsubsection{Freier Anfangs- und Endzeitpunkt} \label{AbschnittFreieZeitSOP}
Wir betrachten in diesem Abschnitt die Aufgabe mit freiem Anfangs- und Endzeitpunkt.
Im Gegensatz zur Aufgabe (\ref{SOP1})--(\ref{SOP4}) bedeutet dies,
dass die Abbildungen $h_0$ bzw. $h_1$, die die Start- und Zielmannigfaltigkeiten definieren, zeitabh"angig sind:
$$h_i (t_i,x_i) : \R \times \R^n \to \R^{s_i }, \quad i = 0,1.$$
Zur Behandlung dieser Aufgaben gehen wir analog zu Ioffe \& Tichomirov \cite{Ioffe}
vor und "uberf"uhren das Problem mit freiem Anfangs- und Endzeitpunkt mittels einer Substitution der Zeit
in eine Aufgabe auf festem Zeitintervall. \\[2mm]
Wir untersuchen nun die Aufgabe mit freiem Anfangs- und Endzeitpunkt:
\begin{eqnarray}
&& \label{FZSOP1} J\big(x(\cdot),u(\cdot)\big) = \int_{t_0}^{t_1} f\big(t,x(t),u(t)\big) \, dt \to \inf, \\
&& \label{FZSOP2} \dot{x}(t) = \varphi\big(t,x(t),u(t)\big), \\
&& \label{FZSOP3} h_0\big(t_0,x(t_0)\big)=0, \qquad h_1\big(t_1,x(t_1)\big)=0, \\
&& \label{FZSOP4} u(t) \in U \subseteq \R^m, \quad U\not= \emptyset, \quad U \mbox{ konvex}.
\end{eqnarray}
Die Aufgabe (\ref{FZSOP1})--(\ref{FZSOP4}) betrachten wir f"ur Tripel $\big([t_0,t_1],x(\cdot),u(\cdot)\big)$ mit
$$[t_0,t_1] \subset \R, \quad x(\cdot) \in W^1_\infty([t_0,t_1],\R^n), \quad u(\cdot) \in L_\infty([t_0,t_1],U).$$
Zur Menge $\mathscr{A}^{\,\mathcal{F}}_{\rm Lip}$ geh"oren diejenigen Tripel
$\big([t_0,t_1],x(\cdot),u(\cdot)\big)$,
f"ur die es eine Zahl $\gamma>0$ derart gibt,
dass die Abbildungen $f(t,x,u)$, $\varphi(t,x,u)$ und $h_i(\tau_i,x_i)$
auf der Menge aller Punkte $(t,\tau_0,\tau_1,x,x_0,x_1,u) \in \R \times \R \times \R \times \R^n \times \R^n \times \R^n \times \R^m$ mit
\begin{eqnarray*}
&& t_0-\gamma < t < t_1+\gamma,\quad t_0-\gamma < \tau_0< t_0+\gamma,\quad t_1-\gamma < \tau_1 < t_1+\gamma, \\
&& \|x-x(t)\| < \gamma, \quad \|x_0-x(t_0)\| < \gamma, \quad \|x_1-x(t_1)\| < \gamma, \quad u \in \R^m
\end{eqnarray*}
stetig und stetig differenzierbar in allen Variablen sind.
(Zur unmissverständlichen Angabe der Punktemenge treten $\tau_0,\tau_1$ in $h_i$ anstelle der Zeitvariablen $t_0,t_1$ auf.) \\[2mm]
In der Aufgabe (\ref{FZSOP1})--(\ref{FZSOP4}) nennen wir ein Tripel $\big([t_0,t_1],x(\cdot),u(\cdot)\big)$ mit
$[t_0,t_1] \subset \R$, $x(\cdot) \in W^1_\infty([t_0,t_1],\R^n)$ und $u(\cdot) \in L_\infty([t_0,t_1],U)$ einen
Steuerungsprozess.
Ein Steuerungsprozess hei"st in der Aufgabe (\ref{FZSOP1})--(\ref{FZSOP4}) zul"assig,
wenn auf dem Intervall $[t_0,t_1]$ die Funktion $x(\cdot)$ fast "uberall der Gleichung (\ref{FZSOP2}) gen"ugt und
die Randbedingungen (\ref{FZSOP3}) erf"ullt.
Die Menge $\mathscr{A}^{\mathcal{F}}_{\rm adm}$ bezeichnet die Menge der zul"assigen Steuerungsprozesse. \\[2mm]
Einen zul"assigen Steuerungsprozess $\big([t_{0*},t_{1*}],x_*(\cdot),u_*(\cdot)\big)$ nennen wir ein
schwaches lokales Minimum\index{Minimum, schwaches lokales!Optimale@-- Optimale Steuerung},
wenn eine Zahl $\varepsilon > 0$ derart existiert,
dass f"ur jeden anderen zul"assigen Steuerungsprozess $\big([t_0,t_1],x(\cdot),u(\cdot)\big)$ mit
$|t_0 - t_{0*}| < \varepsilon$, $|t_1 - t_{1*}| < \varepsilon$ und 
$$\max_{t \in I} \| x(t)-x_*(t) \| < \varepsilon, \qquad \esssup_{t \in I} \| u(t)-u_*(t) \| < \varepsilon, \qquad
  I=[t_{0*},t_{1*}] \cap [t_0,t_1]$$
die Ungleichung $J\big(x(\cdot),u(\cdot)\big) \geq J\big(x_*(\cdot),u_*(\cdot)\big)$ erfüllt ist.

\newpage
Wir werden nun zeigen, wie man mit Hilfe des Theorems \ref{SatzSOP} f"ur Aufgaben mit fester Zeit
unter Verwendung der Substitution der Zeit die notwendigen Optimalit"atsbedingungen f"ur die Aufgabe
mit freiem Anfangs- und Endzeitpunkt erh"alt: \\[2mm]
Es sei $\big([t_{0*},t_{1*}],x_*(\cdot),u_*(\cdot)\big)$ ein schwaches lokales Minimum der Aufgabe (\ref{FZSOP1})--(\ref{FZSOP4}). 
Wir f"uhren eine neue unabh"angige Ver"anderliche $s$ ein, die das Intervall $[0,1]$ durchl"auft,
und betrachten folgendes Differentialgleichungssystem:
\begin{equation} \label{FZSOP6}
t'(s) = v(s), \qquad y'(s) = v(s) \cdot \varphi\big(t(s),y(s),w(s)\big).
\end{equation}
Es sei $v(s) > 0$ f"ur alle $s \in [0,1]$.
Ist $\big( t(\cdot),y(\cdot) \big)$ eine L"osung dieses Systems,
die der Steuerung $\big( v(\cdot),w(\cdot) \big)$ entspricht,
so ist $s \to t(s)$ eine streng wachsende stetige Funktion.
Die zu ihr inverse Funktion $t \to s(t)$ ist ebenfalls stetig und streng wachsend.
In diesem Fall ist $x(\cdot)$ mit $x(t) = y \big( s(t) \big)$ f"ur $t \in [t(0),t(1)]$ L"osung der Gleichung (\ref{FZSOP2}) zur Steuerung
$u(\cdot)$ mit $u(t) = w\big(s(t)\big)$ f"ur $t \in [t(0),t(1)]$, wobei
\begin{equation} \label{FZSOP7}
\int_{t(0)}^{t(1)} f\big(t,x(t),u(t)\big) \, dt = \int_0^1 v(s) \cdot f\big(t(s),y(s),w(s)\big) \, ds
\end{equation}
gilt.
Weil $v(\cdot)$ auf $[0,1]$ stets positiv ist, sind diese Behauptungen klar. \\[2mm]
Ist umgekehrt $x(\cdot)$ eine auf dem Intervall $[t_0,t_1]$ definierte L"osung der Gleichung (\ref{FZSOP2}),
die der Steuerung $u(\cdot)$ entspricht,
so ist $\big(t(\cdot),y(\cdot)\big)$ mit $t(s) = t_0 + (t_1 - t_0) s$, $y(s) = x\big(t(s)\big)$ eine L"osung des Systems (\ref{FZSOP6}),
die den Steuerungen $\big(v(\cdot),w(\cdot)\big)$ mit $v(s) \equiv t_1-t_0$, $w(s) = u\big(t(s)\big)$ zugeordnet ist.
Dabei gilt die Beziehung (\ref{FZSOP7}). \\[2mm]
Somit ist $\big( t_*(\cdot), y_*(\cdot), w_*(\cdot), v_*(\cdot)\big)$ mit
$$\begin{array}{ll}
t_*(s) = t_{0*} + (t_{1*} - t_{0*}) s, &\quad y_*(s) = x_*\big(t_*(s)\big), \\
v_*(s) \equiv v_* = t_{1*} - t_{0*},   &\quad w_*(s) = u_*\big(t_*(s)\big),
\end{array}$$
für $s \in [0,1]$, ein schwaches lokales Minimum der Aufgabe:
\begin{eqnarray}
&& \label{FZSOP8} \int_0^1 v \cdot f\big(t(s),y(s),w(s)\big) \, ds \to \inf, \\
&& \label{FZSOP9} t'(s) = v, \quad y'(s) = v \cdot \varphi\big(t(s),y(s),w(s)\big), \\
&& \label{FZSOP10} h_0\big( t(0),y(0) \big) = 0, \quad h_1\big( t(1),y(1) \big) = 0, \\
&& \label{FZSOP11} v > 0, \quad w(s) \in U, \quad U\not= \emptyset, \quad U \mbox{ konvex}.
\end{eqnarray}
Dies ist bereits eine Aufgabe mit fester Zeit, auf die wir Theorem \ref{SatzSOP} anwenden k"onnen.
Als Zustände agieren $t(\cdot) \in  W^1_\infty([0,1],\R)$ und $y(\cdot) \in  W^1_\infty([0,1],\R^n)$.
Die Steuerungen werden durch $w(\cdot) \in L_\infty([0,1],U)$ und $v >0$ widergegeben.

\newpage
Die Pontrjaginsche Funktion der Aufgabe (\ref{FZSOP8})--(\ref{FZSOP11}) lautet
\begin{eqnarray}
\tilde{H}(t,x,u,v,p,q,\lambda_0) &=& \langle p ,v \cdot \varphi(t,x,u) \rangle + qv - \lambda_0 v \cdot f(t,x,u) \nonumber \\
 \label{FZSOP13} &=& v \cdot [H^{\mathcal{S}}(t,x,u,p,\lambda_0) + q],
\end{eqnarray}
wobei $H^{\mathcal{S}}(t,x,u,p,\lambda_0) = \langle p , \varphi(t,x,u) \rangle - \lambda_0 f(t,x,u)$
die Pontrjagin-Funktion zu der Aufgabe (\ref{FZSOP1})--(\ref{FZSOP4}) ist.
Dann existieren nach Theorem \ref{SatzSOP} nicht gleichzeitig verschwindende Multiplikatoren
$\lambda_0 \geq 0$, $l_0 \in \R^{s_0}$, $l_1 \in \R^{s_1}$,
$\tilde{p}(\cdot) \in W^1_\infty([0,1],\R^n)$ und $\tilde{q}(\cdot) \in W^1_\infty([0,1],\R)$ derart,
dass die Beziehungen
\begin{eqnarray*} 
\tilde{p}'(s) &=& - \tilde{H}_x\big(t_*(s),y_*(s),w_*(s),v_*,\tilde{p}(s),\tilde{q}(s),\lambda_0\big), \\
\tilde{p}(0) &=& h^T_{0,x_0} \big(t_*(0),y_*(0)\big) l_0, \qquad \tilde{p}(1)=- h^T_{1,x_1}\big(t_*(1),y_*(1)\big) l_1
\end{eqnarray*}
und
\begin{eqnarray*} 
\tilde{q}'(s) &=& - \tilde{H}_t\big(t_*(s),y_*(s),w_*(s),v_*,\tilde{p}(s),\tilde{q}(s),\lambda_0\big), \\             
\tilde{q}(0) &=& \big\langle h_{0,t_0} \big(t_*(0),y_*(0)\big) , l_0 \big\rangle,
                 \qquad \tilde{q}(1) = - \big\langle h_{1,t_1}\big(t_*(1),y_*(1)\big),l_1 \big\rangle
\end{eqnarray*}
gelten.
Weiterhin sind die Variationsungleichungen
\begin{eqnarray*}
&& \big\langle \tilde{H}_u\big(t_*(s),y_*(s),w_*(s),v_*,\tilde{p}(s),\tilde{q}(s),\lambda_0\big),\big(w-w_*(s)\big) \big\rangle\leq 0, \\
&& \big\langle \tilde{H}_v\big(t_*(s),y_*(s),w_*(s),v_*,\tilde{p}(s),\tilde{q}(s),\lambda_0\big),(v-v_*) \big\rangle\leq 0
\end{eqnarray*}
f"ur alle $v>0$, $w \in U$ und f"ur fast alle $s \in [0,1]$ erf"ullt. \\[2mm]
Es sei $s_*(\cdot)$ die zu $t_*(\cdot)$ inverse Funktion, d.\,h.
$$s_*(t) = \frac{t-t_{0*}}{t_{1*}-t_{0*}}, \qquad t \in [t_{0*},t_{1*}].$$
Wir f"uhren die Bezeichnungen
$$p(t) = \tilde{p}\big(s_*(t)\big), \quad q(t) = \tilde{q}\big(s_*(t)\big)$$
ein.
So lassen sich die obigen Beziehungen unter Ber"ucksichtigung von (\ref{FZSOP13}) auf folgende Form bringen:
\begin{eqnarray} 
\dot{p}(t) &=& -H^{\mathcal{S}}_x\big(t,x_*(t),u_*(t),p(t),\lambda_0\big), \nonumber \\
p(t_{0*}) &=& h^T_{0,x_0} \big(t_{0*},x_*(t_{0*})\big) l_0, \qquad p(t_{1*})=-h^T_{1,x_1}\big(t_{1*},x_*(t_{1*})\big) l_1, \nonumber\\
\dot{q}(t) &=& \label{FZSOP16} - H^{\mathcal{S}}_t\big(t,x_*(t),u_*(t),p(t),\lambda_0\big), \\
q(t_{0*}) &=& \label{FZSOP17} \big\langle h_{0,t_0}\big(t_{0*},x_*(t_{0*})\big) , l_0 \big\rangle,
              \qquad q(t_{1*})=- \big\langle h_{1,t_1}\big(t_{1*},x_*(t_{1*})\big) , l_1 \big\rangle
\end{eqnarray}
und wir erhalten die Variationsungleichungen
\begin{eqnarray*}
&&  v_* \cdot \big\langle H^{\mathcal{S}}_u\big(t,x_*(t),u_*(t),p(t),\lambda_0\big),\big(u-u_*(t)\big) \big\rangle\leq 0, \\
&& \big\langle H^{\mathcal{S}} \big(t,x_*(t),u_*(t),p(t),\lambda_0\big)+q(t),(v-v_*) \big\rangle\leq 0.
\end{eqnarray*}
Wegen $0 < v_* < \infty$ folgen daraus f"ur alle $u \in U$ die Variationsungleichung
$$\big\langle H^{\mathcal{S}}_u\big(t,x_*(t),u_*(t),p(t),\lambda_0\big),\big(u-u_*(t)\big) \big\rangle\leq 0$$
und die Beziehung
$$H^{\mathcal{S}}\big(t,x_*(t),u_*(t),p(t),\lambda_0\big) = -q(t).$$
Durch Vergleich mit (\ref{FZSOP16}) und (\ref{FZSOP17}) erhalten wir weiter
\begin{eqnarray*}
&& \frac{d}{dt} H^{\mathcal{S}}\big(t,x_*(t),u_*(t),p(t),\lambda_0\big) = H^{\mathcal{S}}_t\big(t,x_*(t),u_*(t),p(t),\lambda_0\big), \\
&& H^{\mathcal{S}}\big(t_{0*},x_*(t_{0*}),u_*(t_{0*}),p(t_{0*}),\lambda_0\big) = - \big\langle h_{0,t_0} \big(t_{0*},x_*(t_{0*})\big) , l_0 \big\rangle, \\
&& H^{\mathcal{S}}\big(t_{1*},x_*(t_{1*}),u_*(t_{1*}),p(t_{1*}),\lambda_0\big) = \big\langle h_{1,t_1} \big(t_{1*},x_*(t_{1*})\big) , l_1 \big\rangle.
\end{eqnarray*}
Wir fassen alle Bedingungen zusammen und formulieren das Schwache Optimalit"atsprinzip f"ur die
Aufgabe mit freiem Anfangs- und Endzeitpunkt: 

\begin{theorem} \label{SatzFZSOP} \index{Schwaches Optimalitätsprinzip}
In der Aufgabe (\ref{FZSOP1})--(\ref{FZSOP4}) sei 
$\big([t_{0*},t_{1*}],x_*(\cdot),u_*(\cdot)\big) \in \mathscr{A}^{\mathcal{F}}_{\rm adm} \cap \mathscr{A}^{\mathcal{F}}_{\rm Lip}$.
Ist $\big([t_{0*},t_{1*}],x_*(\cdot),u_*(\cdot)\big)$ ein schwaches lokales Minimum der Aufgabe (\ref{FZSOP1})--(\ref{FZSOP4}),
dann existieren nicht gleichzeitig verschwindende Multiplikatoren
$\lambda_0 \geq 0$, $l_0 \in \R^{s_0}$, $l_1 \in \R^{s_1}$ und $p(\cdot) \in W^1_\infty([t_{0*},t_{1*}],\R^n)$
derart, dass
\begin{enumerate}
\item[(a)] für fast alle $t \in [t_{0*},t_{1*}]$ die adjungierte Gleichung
           \index{adjungierte Gleichung!Schwach@-- Schwaches Optimalitätsprinzip}
           \begin{equation} \label{SatzFZSOP1}
           \dot{p}(t) = -H^{\mathcal{S}}_x\big(t,x_*(t),u_*(t),p(t),\lambda_0\big),
           \end{equation}
\item[(b)] die Transversalit"atsbedingungen
           \index{Transversalitätsbedingungen!Schwach@-- Schwaches Optimalitätsprinzip}
           \begin{equation}\label{SatzFZSOP2}
           p(t_{0*})= h^T_{0,x_0} \big(t_{0*},x_*(t_{0*})\big) l_0, \qquad p(t_{1*})=-h^T_{1,x_1}\big(t_{1*},x_*(t_{1*})\big) l_1,
           \end{equation}
\item[(c)] in fast allen Punkten $t\in [t_{0*},t_{1*}]$ für alle $u \in U$ die Variationsungleichung
          \begin{equation}\label{SatzFZSOP3}
          \big\langle H^{\mathcal{S}}_u\big(t,x_*(t),u_*(t),p(t),\lambda_0\big),\big(u-u_*(t)\big) \big\rangle\leq 0
          \end{equation}
\item[(d)] und f"ur die Funktion $t \to H^{\mathcal{S}}\big(t,x_*(t),u_*(t),p(t),\lambda_0\big)$ die Bedingungen
          \begin{eqnarray}
          \frac{d}{dt} H^{\mathcal{S}}\big(t,x_*(t),u_*(t),p(t),\lambda_0\big) &=& H^{\mathcal{S}}_t\big(t,x_*(t),u_*(t),p(t),\lambda_0\big), \\
          \sup_{u \in U} H^{\mathcal{S}}\big(t_{0*},x_*(t_{0*}),u,p(t_{0*}),\lambda_0\big)
             &=& - \big\langle h_{0,t_0} \big(t_{0*},x_*(t_{0*})\big) , l_0 \big\rangle, \\
          \sup_{u \in U} H^{\mathcal{S}}\big(t_{1*},x_*(t_{1*}),u,p(t_{1*}),\lambda_0\big)
             &=& \big\langle h_{1,t_1} \big(t_{1*},x_*(t_{1*})\big) , l_1 \big\rangle
          \end{eqnarray}
\end{enumerate}
erfüllt sind.
\end{theorem}

\begin{beispiel} \label{BeispielMassenpunkt}
{\rm Wir betrachten die Aufgabe einer zeitminimalen Bewegung eines Massenpunktes entlang einer horizontalen Geraden
bei vernachl"assigter Reibung (\!\!\cite{Feichtinger}, S.\,75--77):
\begin{eqnarray*}
&& J\big(s(\cdot),v(\cdot),a(\cdot)\big) = \int_0^T 1 \, dt \to \inf, \\
&& \dot{s}(t) = v(t), \qquad \dot{v}(t) = a(t), \\
&& s(0)=0, \quad s(T)=s_T>0, \quad v(0)=v(T)=0, \quad a(t) \in [-1,1].
\end{eqnarray*}
Darin bezeichnen $s=s(t)$ den Ort, $v=v(t)$ die Geschwindigkeit und $a=a(t)$ die Beschleunigung zum Zeitpunkt $t$.
Der Steuerungsparameter wird durch ausge"ubte Beschleunigung gegeben.
Die Zustandsgr"oßen sind der Ort und die Geschwindigkeit. \\[2mm]
F"ur diese Aufgabe lautet die Pontrjagin-Funktion
$$H^{\mathcal{S}}(t,s,v,a,p_1,p_2,\lambda_0) = p_1 v + p_2 a - \lambda_0.$$
Die notwendigen Bedingungen in Theorem \ref{SatzFZSOP} liefern f"ur die Adjungierten
\begin{eqnarray*}
\dot{p}_1(t) \equiv 0 \quad&\Rightarrow&\quad p_1(t)\equiv \pi_1 \in \R, \\
\dot{p}_2(t) = -p_1(t)=-\pi_1 \quad&\Rightarrow&\quad p_2(t) = \pi_2 -\pi_1 t,  \; \pi_2 \in \R.
\end{eqnarray*}
Au"serdem erhalten wir aus der Maximumbedingung
$$a_*(t) = -1 \quad\mbox{f"ur } p_2(t)<0, \qquad a_*(t) = 1 \quad\mbox{f"ur } p_2(t)>0.$$
Da der Massenpunkt aus einer Ruhelage zeitoptimal in eine andere Ruhelage "uberf"uhrt wird,
folgt aus der Darstellung von $p_2(\cdot)$ ein Vorzeichenwechsel zum Zeitpunkt $\tau=T_*/2$. \\
Weiterhin liegt eine autonome Aufgabe vor und es ergibt sich
$$H^{\mathcal{S}}\big(t,s_*(t),v_*(t),a_*(t),p_1(t),p_2(t),\lambda_0\big) \equiv 0.$$
Daraus ergibt sich, dass der anormale Fall $\lambda_0=0$ ausgeschlossen werden kann. \\
Zusammenfassend liefern die notwendigen Bedingungen
$$\lambda_0=1, \qquad \pi_1=\frac{2}{T_*}, \qquad \pi_2=1.$$
\begin{minipage}{0.54\textwidth}
Die Bewegung des Massenpunktes besitzt damit die Darstellung
\begin{eqnarray*} 
a_*(t) &=& \left\{ \begin{array}{ll} \;\;\;1, & t \in [0,\tau), \\ -1, \quad & t \in [\tau,T_*], \end{array} \right. \\[2mm]
v_*(t) &=& \left\{ \begin{array}{ll} t, & t \in [0,\tau), \\ T_*-t, \quad & t \in [\tau,T_*], \end{array} \right. \\[2mm]
s_*(t) &=& \left\{ \begin{array}{ll}
                   \frac{1}{2} t^2, & t \in [0,\tau), \\[1mm]
                   \tau^2-\frac{1}{2}(T_*-t)^2, \quad & t \in [\tau,T_*]. \end{array} \right.
\end{eqnarray*}
\end{minipage}
\begin{minipage}{0.45\textwidth}
\centering
\includegraphics[width=5cm]{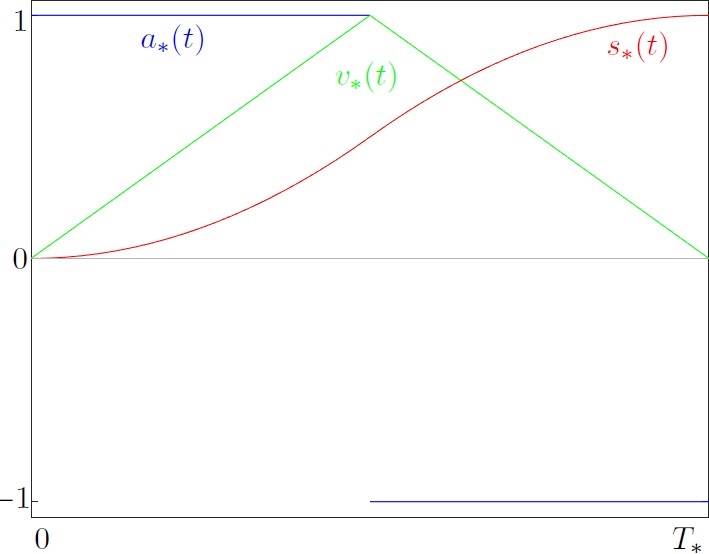}
\captionof{figure}[Bewegung eines Massenpunktes]{Bewegung für $s(T)=1$.}
\label{AbbMassenpunkt}
\end{minipage} \\[2mm]
Abschlie"send erhalten wir aus $s_*(T_*)=s_T$ unmittelbar $T_*=2 \sqrt{s_T}$. \hfill $\square$}
\end{beispiel}

%% file: 1-36-Zustandsaufgabe.tex
\subsubsection{Ein Schwaches Optimalit\"atsprinzip unter Zustandsbeschr\"ankungen} \label{AbschnittZustandSOP}
\begin{theorem} \label{SatzSOPZA} \index{Schwaches Optimalitätsprinzip}
Es sei $\big(x_*(\cdot),u_*(\cdot)\big) \in \mathscr{A}^{\mathcal{S}}_{\rm adm} \cap \mathscr{A}^{\mathcal{S}}_{\rm Lip}$. 
Ist $\big(x_*(\cdot),u_*(\cdot)\big)$ ein schwaches lokales Minimum der Aufgabe (\ref{SOP1})--(\ref{SOP5}),
dann existieren eine Zahl $\lambda_0 \geq 0$, eine Vektorfunktion $p(\cdot):[t_0,t_1] \to \R^n$, Vektoren $l_0 \in \R^{s_0}$, $l_1 \in \R^{s_1}$
und auf den Mengen
$$T_j=\big\{t \in [t_0,t_1] \,\big|\, g_j\big(t,x_*(t)\big)=0\big\}, \quad j=1,...,l,$$
konzentrierte nichtnegative regul"are Borelsche Ma"se $\mu_j$ endlicher Totalvariation
(wobei s"amtliche Gr"o"sen nicht gleichzeitig verschwinden) derart,
dass die Vektorfunktion $p(\cdot)$ von beschr"ankter Variation und rechtsseitig stetig ist, und
\begin{enumerate}
\item[(a)] die adjungierte Gleichung
           \index{adjungierte Gleichung!Schwach@-- Schwaches Optimalitätsprinzip}
           \begin{eqnarray}
           p(t)&=&-{h_1'}^T\big(x_*(t_1)\big)l_1 + \int_t^{t_1} H^{\mathcal{S}}_x\big(s,x_*(s),u_*(s),p(s),\lambda_0\big) \, ds \nonumber \\
           \label{SatzSOPZA1} & & -\sum_{j=1}^l \int_t^{t_1} g_{j,x}\big(s,x_*(s)\big)\, d\mu_j(s),
           \end{eqnarray}
\item[(b)] die Transversalit"atsbedingung
           \index{Transversalitätsbedingungen!Schwach@-- Schwaches Optimalitätsprinzip}
           \begin{equation}\label{SatzSOPZA2}
           p(t_0)= {h_0'}^T\big(x_*(t_0)\big)l_0, \qquad p(t_1)=-{h_1'}^T\big(x_*(t_1)\big)l_1,
           \end{equation}
           in $t=t_1$ die Sprung-Transversalit"atsbedingung
           \begin{equation} \label{ZBEndbedingung}
           p(t_1^-) - p(t_1) = - \sum_{j=1}^l \mu_j(\{t_1\}) \, g_{j,x}\big(t_1,x_*(t_1)\big)
           \end{equation}
\item[(c)] und in fast allen Punkten $t \in [t_0,t_1]$ f"ur alle $u \in U$ die Variationsungleichung
           \begin{equation}\label{SatzSOPZA3}
           \big\langle H^{\mathcal{S}}_u\big(t,x_*(t),u_*(t),p(t),\lambda_0\big),\big(u-u_*(t)\big) \big\rangle\leq 0
           \end{equation}
\end{enumerate}
erfüllt sind.
\end{theorem}

Als Funktion beschr"ankter Variation existieren an jeder Stelle $t \in [t_0,t_1]$ die einseitigen Grenzwerte der Adjungierten $p(\cdot)$. 
Wegen der rechtsseitigen Stetigkeit kann die Adjungierte an der Stelle $t=t_1$ unstetig sein.
In diesem Fall gilt im Endpunkt $t_1$ die Sprung-Transversalit"atsbedingung (\ref{ZBEndbedingung}). \\[2mm]
Als ein anschauliches Beispiel,
in dem die Zustandsbeschränkung eine Sprungbedingung im Endpunkt hervorruft,
betrachten wir nach Seierstad \& Syds\ae ter \cite{Seierstad} das folgende Lagerhaltungsmodell:

\begin{beispiel} \label{BspLagerhaltungZB} \index{Lagerhaltung}
{\rm Ein Gesch"aft m"ochte die Nachfrage $a$ des Produktes $x$ gegen"uber seinen Kunden gew"ahrleisten.
Es entstehen dadurch zum Zeitpunkt $t$ die Einkaufs- und Lagerhaltungskosten $\alpha u(t)$ bzw. $\beta x(t)$.
Daraus ergibt sich das folgende Modell:
\begin{eqnarray*}
&& J\big(x(\cdot),u(\cdot)\big) =\int_0^T \big(\alpha u(t)+\beta x(t)\big) \, dt \to \inf, \\
&& \dot{x}(t)=u(t)-a, \quad x(0)=1, \quad x(T) \mbox{ frei}, \quad x(t) \geq 0, \quad u(t) \geq 0, \quad T>1/a.
\end{eqnarray*}
Der Fall $\lambda_0=0$ kann ausgeschlossen werden und die Pontrjagin-Funktion lautet
$$H^{\mathcal{S}}(t,x,u,p,1)=p(u-a)-(\alpha u + \beta x)=(p-\alpha)u-(pa+\beta x).$$
In diesem Modell k"onnen wir als optimalen Kandidaten den Steuerungsprozess
$$u_*(t)= \left\{\begin{array}{ll} 0, & t \in [0,1/a), \\[1mm] a, & t \in [1/a,T], \end{array}\right. \quad
  x_*(t)= \left\{\begin{array}{ll} 1-at, & t \in [0,1/a), \\[1mm] 0, & t \in [1/a,T], \end{array}\right.$$
ableiten.
F"ur $t \in [0,1/a)$ sind die Zustandsbeschr"ankungen nicht aktiv und es gilt
$$\dot{p}(t)=\beta , \qquad p(t)=p(0)+\beta t.$$
F"ur $t \in [1/a,T]$ ist $u_*(t)=a$ ein innerer Punkt des Steuerungsbereiches $U = \R_+$.
Also folgt aus der Variationsungleichung (\ref{SatzSOPZA3}), dass $p(t)=\alpha$ gelten muss.
Damit erhalten wir
$$p(t)= \left\{\begin{array}{ll} \alpha+ \beta (t-1/a), & t \in [0,1/a), \\[1mm] \alpha, & t \in [1/a,T), \\[1mm] 0, & t=T. \end{array}\right.$$
Mit der Zustandsbeschr"ankung $g\big(x(t)\big)=-x(t) \leq 0$ lautet die adjungierte Gleichung
$$\dot{p}(t)=0= \beta + \lambda(t)g'\big(x_*(t)\big) = \beta - \lambda(t), \qquad t \in (1/a,T),$$
und wir erhalten den absolutstetigen Anteil $\lambda(t)=\beta$ f"ur das positive Ma"s $\mu$.
Ferner ergibt sich in der Sprungbedingung (\ref{ZBEndbedingung}) im Endpunkt
$$p(T^-)-p(T)= p(T^-)= -g_x\big(x_*(T)\big)\, \mu(\{T\}) = \mu(\{T\}) = \alpha$$
f"ur den atomaren Anteil des Ma"ses $\mu$ im Punkt $t=T$. \\[2mm]
Wir zeigen nun, dass der Kandidat $\big(x_*(\cdot),u_*(\cdot)\big)$ in der Tat global optimal ist: \\
Für den Zustand $x_*(\cdot)$ gilt,
dass keine zul"assige Trajektorie $x(\cdot)$ mit $x(t)<x_*(t)$ f"ur ein $t \in [0,T]$ existiert.
Damit ergibt sich direkt:
\begin{eqnarray*}
    J\big(x(\cdot),u(\cdot)\big)
&=& \int_0^T \big(\alpha u(t)+\beta x(t)\big) \, dt = \int_0^T \big(\alpha [a+\dot{x}(t)] +\beta x(t)\big) \, dt \\
&=& \alpha a T + \alpha [x(T)-x(0)] + \int_0^T \beta x(t) \, dt \\
&\geq& \alpha (aT-1)  + \int_0^T \beta x_*(t) \, dt = J\big(x_*(\cdot),u_*(\cdot)\big).
\end{eqnarray*}
Diese Ungleichung zeigt die globale Optimalit"at des Steuerungsprozesses $\big(x_*(\cdot),u_*(\cdot)\big)$. \hfill $\square$}
\end{beispiel}

%% file: 1-37-Beweis.tex
\subsubsection{Der Nachweis der notwendigen Optimalit\"atsbedingungen}
Wir betrachten f"ur $\big(x(\cdot),u(\cdot)\big) \in C([t_0,t_1],\R^n) \times L_\infty([t_0,t_1],\R^m)$ die Abbildungen
\begin{eqnarray*}
J\big(x(\cdot),u(\cdot)\big) &=& \int_{t_0}^{t_1} f\big(t,x(t),u(t)\big) \, dt, \\
F\big(x(\cdot),u(\cdot)\big)(t) &=& x(t) -x(t_0) -\int_{t_0}^t \varphi\big(s,x(s),u(s)\big) \, ds, \quad t \in [t_0,t_1],\\
H_i\big(x(\cdot)\big) &=& h_i\big(x(t_i)\big), \quad i=0,1, \\
G_j\big(x(\cdot)\big)(t) &=& g_j\big(t,x(t)\big), \quad t \in [t_0,t_1], \quad j=1,...,l.
\end{eqnarray*}
Da $x(\cdot)$ zu $C([t_0,t_1],\R^n)$ geh"ort, gilt f"ur diese Abbildungen
\begin{eqnarray*}
J &:& C([t_0,t_1],\R^n) \times L_\infty([t_0,t_1],\R^m) \to \R, \\
F &:& C([t_0,t_1],\R^n) \times L_\infty([t_0,t_1],\R^m) \to C_0([t_0,t_1],\R^n), \\
H_i &:& C([t_0,t_1],\R^n) \to \R^{s_i}, \quad i=0,1, \\
G_j &:& C([t_0,t_1],\R^n) \to C([t_0,t_1],\R), \quad j=1,...,l.
\end{eqnarray*}
Wir setzen $\mathscr{F}=(F,H_0,H_1)$ und $G=(G_1,...,G_l)$.
Au"serdem f"uhren wir f"ur die Elemente $x(\cdot)= \big(x_1(\cdot),...,x_l(\cdot)\big)$ des Raumes $C([t_0,t_1],\R^l)$
folgende Halbordnung ``$\preceq$'' ein:
$$x(\cdot) \preceq y(\cdot) \qquad\Leftrightarrow\qquad x_j(t) \leq y_j(t) \mbox{ f"ur alle } t \in [t_0,t_1], \; j=1,...,l.$$
Im Raum $C([t_0,t_1],\R^l)$ bezeichnen wir mit $\mathscr{K}$ den Kegel
$$\mathscr{K} =\{ x(\cdot) \in C([t_0,t_1],\R^l) \,|\, x(\cdot) \preceq 0\}.$$  
$\mathscr{K}$ ist ein konvexer, abgeschlossener Kegel mit Spitze in Null.
Zudem besitzt der Kegel $\mathscr{K}$ ein nichtleeres Inneres, d.\,h. ${\rm int\,}\mathscr{K} \not= \emptyset$. \\[1mm]
Mit diesen Setzungen pr"ufen wir f"ur die Extremalaufgabe
\begin{equation} \label{ExtremalaufgabeSOPZA}
J\big(x(\cdot),u(\cdot)\big) \to \inf, \quad \mathscr{F}\big(x(\cdot),u(\cdot)\big)=0, \quad
G\big(x(\cdot)\big) \in \mathscr{K}, \quad  u(\cdot) \in L_\infty([t_0,t_1],U)
\end{equation}
in $\big(x_*(\cdot),u_*(\cdot)\big) \in \mathscr{A}^{\mathcal{S}}_{\rm Lip}$,
wobei wir $x_*(\cdot)$ als Element des Raumes $C([t_0,t_1],\R^n)$ auffassen,
die Voraussetzungen von Theorem \ref{SatzExtremalprinzipSchwach}:

\begin{enumerate}
\item[(A$_2$)] Mit Verweis auf Abschnitt \ref{AbschnittBeweisSOP} ist nur noch die Abbildung $G$ zu diskutieren,
               deren stetige Fr\'echet-Differenzierbarkeit im Beispiel \ref{DiffAbbildung} nachgewiesen ist.              
\end{enumerate}

Zur Extremalaufgabe (\ref{ExtremalaufgabeSOPZA}) definieren wir auf
$$C([t_0,t_1],\R^n) \times L_\infty([t_0,t_1],\R^m) \times \R \times C_0^*([t_0,t_1],\R^n) \times \R^{s_0} \times \R^{s_1} \times \big(C^*([t_0,t_1],\R)\big)^l$$
die Lagrange-Funktion $\mathscr{L}=\mathscr{L}\big(x(\cdot),u(\cdot),\lambda_0,y^*,l_0,l_1,z_1^*,...,z_l^*\big)$,
$$\mathscr{L}= \lambda_0 J\big(x(\cdot),u(\cdot)\big)+ \big\langle y^*, F\big(x(\cdot),u(\cdot)\big) \big\rangle
                         +l_0^T H_0\big(x(\cdot)\big)+l_1^T H_1\big(x(\cdot)\big) + \sum_{j=1}^l \big\langle z_j^*, G_j\big(x(\cdot)\big) \big\rangle.$$
Ist $\big(x_*(\cdot),u_*(\cdot)\big)$ eine schwache lokale Minimalstelle der Aufgabe (\ref{ExtremalaufgabeSOPZA}),
dann existieren nach Theorem \ref{SatzExtremalprinzipSchwach}
nicht gleichzeitig verschwindende Lagrangesche Multiplikatoren $\lambda_0 \geq 0$, $y^* \in C_0^*([t_0,t_1],\R^n)$, $l_i \in \R^{s_i}$
und $z_j^* \in C^*([t_0,t_1],\R)$ derart,
dass gelten:
\begin{enumerate}
\item[(a)] Die Lagrange-Funktion besitzt bez"uglich $x(\cdot)$ in $x_*(\cdot)$ einen station"aren Punkt, d.\,h.
           \begin{equation}\label{SatzSOPZALMR1}
           \mathscr{L}_x\big(x_*(\cdot),u_*(\cdot),\lambda_0,y^*,l_0,l_1,z_1^*,...,z_l^*\big)=0;
           \end{equation}         
\item[(b)] Die Lagrange-Funktion erf"ullt bez"uglich $u(\cdot)$ in $u_*(\cdot)$ die Variationsungleichung
           \begin{equation}\label{SatzSOPZALMR2}
           \big\langle \mathscr{L}_u\big(x_*(\cdot),u_*(\cdot),\lambda_0,y^*,l_0,l_1,z_1^*,...,z_l^*\big), u(\cdot)-u_*(\cdot) \big\rangle \geq 0
           \end{equation}
           f"ur alle $u(\cdot) \in L_\infty([t_0,t_1],U)$;
\item[(c)] Die komplement"aren Schlupfbedingungen gelten, d.\,h.
           \begin{equation}\label{SatzSOPZALMR3}
           0 = \big\langle z_i^*, G_i\big(x_*(\cdot)\big) \big\rangle, \quad
           \langle z_i^*,z(\cdot) \rangle \leq 0 \quad\mbox{f"ur alle } z(\cdot) \preceq 0, \quad i=1,...,l.
           \end{equation}
\end{enumerate}

Die zweite Bedingung in (\ref{SatzSOPZALMR3}) liefert,
dass die Funktionale $z_j^* \in C^*([t_0,t_1],\R)$ positiv sind.
Dann ergibt sich aus der ersten Bedingung in (\ref{SatzSOPZALMR3}),
dass nur diejenigen stetigen linearen Funktionale $z_j^*$ von Null verschieden sein k"onnen,
f"ur welche die Zustandsbeschr"ankungen aktiv sind.
Au"serdem folgt in diesem Fall aus der ersten Bedingung,
dass die zugeh"origen endlichen nichtnegativen regul"aren Borelschen Ma"se $\mu_j$
auf den Mengen
$$T_j=\big\{t \in [t_0,t_1] \,\big|\, g_j\big(t,x_*(t)\big)=0\big\}, \quad j=1,...,l,$$
konzentriert sind.
Daher kann man ohne Einschr"ankung annehmen,
dass alle Ma"se $\mu_j$ auf den Mengen $T_j$ konzentriert sind. \\[2mm]
Aufgrund (\ref{SatzSOPZALMR1}) ist folgende Variationsgleichung f"ur alle $x(\cdot) \in C([t_0,t_1],\R^n)$ erf"ullt: 
\begin{eqnarray}
0 &=& \lambda_0 \int_{t_0}^{t_1} \big\langle f_x\big(t,x_*(t),u_*(t)\big),x(t) \big\rangle\, dt
      + \big\langle l_0, h_0'\big(x_*(t_0)\big) x(t_0) \big\rangle + \big\langle l_1, h_1'\big(x_*(t_1)\big) x(t_1) \big\rangle
      \nonumber \\
  & & + \int_{t_0}^{t_1} \bigg[ x(t)-x(0) - \int_0^t \varphi_x\big(s,x_*(s),u_*(s)\big) x(s) \,ds \bigg]^T d\mu(t) \nonumber \\
  & & \label{BeweisschlussSOPZA1}
      + \sum_{j=1}^l \int_{t_0}^{t_1}\big\langle g_{j,x}\big(t,x_*(t)\big),x(t) \big\rangle \,d\mu_j(t).
\end{eqnarray}
In der Gleichung (\ref{BeweisschlussSOPZA1}) "andern wir die Integrationsreihenfolge im zweiten Summanden und bringen sie in die Form
\begin{eqnarray}
0 &=& \int_{t_0}^{t_1} \big\langle \lambda_0 f_x\big(t,x_*(t),u_*(t)\big) - \varphi_x^T\big(t,x_*(t),u_*(t)\big) \int_{t}^{t_1} d\mu(s) ,
                       x(t) \big\rangle \, dt \nonumber \\
  & & + \int_{t_0}^{t_1} [x(t)]^T \, d\mu(t) + \Big\langle {h_0'}^T\big(x_*(t_0)\big)l_0 - \int_{t_0}^{t_1} d\mu(t) , x(t_0) \Big\rangle
      + \langle {h_1'}^T\big(x_*(t_1)\big)l_1 , x(t_1) \rangle \nonumber \\
  & & \label{BeweisschlussSOPZA2}
      + \sum_{j=1}^l \int_{t_0}^{t_1} \big\langle g_{j,x}\big(t,x_*(t)\big),x(t) \big\rangle \,d\mu_j(t).
\end{eqnarray}
Wir setzen $p(t)=\displaystyle \int_t^{t_1} \, d\mu(s)$,
so ist $p(\cdot)$ eine Funktion von beschränkter Variation und gemäß den Eigenschaften einer Verteilungsfunktion rechtsseitig stetig. \\[1mm]
Die rechte Seite in (\ref{BeweisschlussSOPZA2}) definiert ein stetiges lineares Funktional im Raum $C([t_0,t_1],\R^n)$.
Wenden wir den Darstellungssatz von Riesz an,
so folgen aus der eindeutigen Darstellung eines stetigen linearen Funktionals im Raum $C([t_0,t_1],\R^n)$ die Beziehungen
\begin{eqnarray*}
p(t) &=& -{h_1'}^T\big(x_*(t_1)\big)l_1
         + \int_t^{t_1} \big[ \varphi^T_x\big(s,x_*(s),u_*(s)\big)p(s)-\lambda_0 f_x\big(s,x_*(s),u_*(s)\big)\big]\, ds \\
     & & - \sum_{j=1}^l \int_t^{t_1} g_{j,x}\big(s,x_*(s)\big) \,d\mu_j(s) \\
     &=&-{h_1'}^T\big(x_*(t_1)\big)l_1 + \int_t^{t_1} H^{\mathcal{S}}_x\big(s,x_*(s),u_*(s),p(s),\lambda_0\big) \, ds \\
     & & -\sum_{j=1}^l \int_t^{t_1} g_{j,x}\big(s,x_*(s)\big)\, d\mu_j(s), \\
p(t_0) &=& {h_0'}^T\big(x_*(t_0)\big)l_0.
\end{eqnarray*}
Damit sind (\ref{SatzSOPZA1}) und (\ref{SatzSOPZA2}) gezeigt.
Die Beziehung (\ref{SatzSOPZALMR2}) ist "aquivalent zu
$$\int_{t_0}^{t_1} \big\langle H^{\mathcal{S}}_u\big(t,x_*(t),u_*(t),p(t),\lambda_0\big),u(t)-u_*(t) \big\rangle \, dt \leq 0.$$
Wir bemerken, dass die Funktion $p(\cdot)$ von beschr"ankter Variation ist.
Daher kann man $p(\cdot)$ als Differenz zweier monotoner Funktionen schreiben und
$p(\cdot)$ besitzt h"ochstens abz"ahlbar viele Unstetigkeiten.
Mit Hilfe dieser Anmerkung folgt abschlie"send die Variationsungleichung (\ref{SatzSOPZA3}) via Standardtechniken
f"ur Lebesguesche Punkte. \hfill $\blacksquare$

%% file: 2-0-Nadelvariationen.tex
\section{Nadelvariationen in der Optimalen Steuerung} \label{KapitelStark}
Die Methoden der Optimalen Steuerung liefern geeignete Werkzeuge f"ur viele Problemklassen.
In den meisten praktischen Anwendungen sind allerdings Steuerungselemente enthalten,
die nur gewisse Positionen annehmen d"urfen.
Eine Richtungsvariation des letzten Kapitels ist bez"uglich dieser Elemente nicht zul"assig.
Deswegen wollen wir in diesem Kapitel eine neue Klasse von Variationen betrachten.
Die Grundlage daf"ur wurde im Rahmen der Klassischen Variationsrechnung von Karl Weistra"s (1815--1897)
durch die Einf"uhrung der Nadelvariationsmethode gelegt.
Im Vergleich zur schwachen lokalen Optimalstelle entsteht dadurch ein allgemeinerer Optimalit"atsbegriff. \\[2mm]
Der "Ubergang von der Klassischen Variationsrechung zur Optimalen Steuerung 
wurde durch die Identifizierung der Aufgabenklassen der Steuerungstheorie und durch die Entwicklung von geeigneten Methoden realisiert.
Wesentliche Beitr"age zu diesen Entdeckungen haben
Constantin Carath\'eodory (1873--1950), Edward James McShane (1904--1989), Magnus Hestenes (1906--1991),
Richard Bellman (1920--1984),
Rufus Isaacs (1914--1977),
sowie Lew Semjonowitsch Pontrjagin (1908--1988) in Zusammenarbeit mit Wladimir Boltjanski (1925--2019) und
Rewas Gamkrelidze (1927--) geleistet.
Die bedeutendsten Ergebnisse in der Entwicklung von Optimalit"atsbedingungen sind das Pontrjaginsche Maximumprinzip und das
Bellmansche Prinzip der Dynamischen Programmierung. \\[2mm]
Im Fokus dieses Kapitels steht das Pontrjaginsche Maximumprinzip.
Der Errungenschaft des Maximumprinzips gehen wichtige Entwicklungen der Mathematik in der ersten H"alfte des 20. Jahrhunderts voraus.
Ohne die Methoden der Funktionalanalysis oder der Theorie von Differentialgleichungen mit unstetigen rechten Seiten
ist dessen vollst"andiger Nachweis kaum denkbar.
Dieses vielschichtige und umfassende theoretische Fundament,
der komplexe Beweis des Maximumprinzips und die eigenst"andigen Methoden, die zum Nachweis ausgearbeitet wurden,
bilden einen ma"sgeblichen Beitrag zur Etablierung der Optimalen Steuerung als eigenst"andige mathematische Disziplin. \\[2mm]
Das Kapitel zu Nadelvariationen und starken lokalen Minimalstellen in Steuerungsproblemen beginnen wir mit der Untersuchung
der Aufgabe mit freiem Endpunkt.
Da neben der Dynamik zu einem Anfangswert keine weiteren Beschr"ankungen an die Zustandstrajektorie vorliegen,
ist die Methode der einfachen Nadelvariation ein geeignetes Mittel zur Herleitung von notwendigen Optimalt"atsbedingungen. \\[2mm]
Der "Ubergang zu allgemeineren Aufgabenklassen ist mit der einfachen Nadelvariationsmethode nicht m"oglich.
Das Pontrjaginsche Maximumprinzip f"ur die Standardaufgabe und die Aufgabe mit Zustandsbeschr"ankungen beweisen wir deswegen 
auf der Grundlage der mehrfachen Nadelvariationen nach Ioffe \& Tichomirov \cite{Ioffe},
die im Anhang ausf"uhrlich behandelt werden.
Au"serdem behandeln wir hinreichende Bedingungen nach Arrow und die Aufgabe mit freiem Anfangs- und Endzeitpunkt.

%% file: 2-1-Nadelvariationen.tex
\subsection{Die Nadelvariation der Klassischen Variationsrechnung}
Zur Einstimmung auf die Untersuchung von starken lokalen Minimalstellen betrachten wir zun"achst wieder die einfachste Aufgabe der
Klassischen Variationsrechnung:
\begin{eqnarray}
&& \label{WNV1} J\big(x(\cdot)\big) = \int_{t_0}^{t_1} f\big(t,x(t),\dot{x}(t)\big) \, dt \to \inf, \\
&& \label{WNV2} x(t_0)=x_0, \qquad x(t_1)=x_1.
\end{eqnarray}
Ebenso wie in Abschnitt \ref{AbschnittELG} seien die Punkte $t_0<t_1$, $x_0,x_1$ fest vorgegeben und es sei
die Funktion $f:\R \times \R \times \R \to \R$ in (\ref{WNV1}) 
stetig differenzierbar. \\[2mm]
Im Gegensatz zur Analyse der Aufgabe (\ref{ELG1})--(\ref{ELG2}) betrachten wir an dieser Stelle nicht die
Richtungsvariation
$x_\lambda(\cdot) = x_*(\cdot) + \lambda x(\cdot)$
mit $x_*(\cdot),x(\cdot) \in C_1([t_0,t_1],\R)$ und entsprechender Normkonvergenz
$\| x_\lambda(\cdot) - x_*(\cdot) \|_{C_1} = |\lambda| \| x(\cdot) \|_{C_1} \to 0$ f"ur $\lambda \to 0$.
Sondern wir definieren zu $v \in \R$, $\tau \in (t_0,t_1)$ und $0 < \lambda<\varepsilon$ mit $\tau+\varepsilon < t_1$
folgende st"uckweise lineare Funktion:\index{Nadelvariation, einfache}
\begin{equation}\label{WNV3}
\xi_\lambda(t)= \left\{\begin{array}{ll} 0, & t \not\in [\tau,\tau+\varepsilon], \\[1mm]
                                 v(t-\tau), & t \in [\tau,\tau+\lambda), \\[1mm]
                                 \mbox{linear}, & t \in [\tau+\lambda,\tau+\varepsilon], \end{array}\right. \quad
  \dot{\xi}_\lambda(t)= \left\{\begin{array}{rl}
  0, & t \not\in [\tau,\tau+\varepsilon], \\[1mm]
  v, & t \in [\tau,\tau+\lambda), \\[1mm]
  -\frac{\lambda v}{\varepsilon-\lambda}, & t \in [\tau+\lambda,\tau+\varepsilon]. \end{array} \right.
\end{equation}
Die Funktion $\xi_\lambda(\cdot)$ besitzt die Form einer Nadelspitze,
die bei Verkleinerung des Parameters $\lambda$ zwar schmaler wird,
aber deren Anstieg sich auf dem Intervall $[\tau,\tau+\lambda)$ nicht verringert.

\begin{figure}[h]
	\centering
	\includegraphics[width=6.5cm]{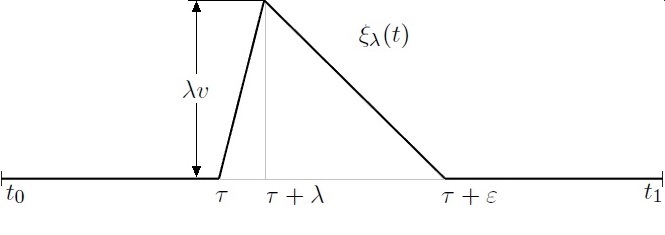} \hspace*{5mm}
          \includegraphics[width=6.5cm]{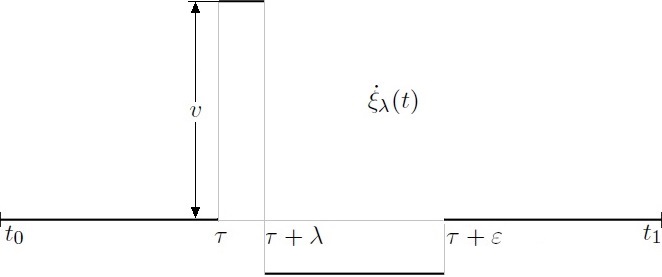}
	\caption[Nadelvariation der Klassischen Variationsrechnung]{Nadelvariation der Klassischen Variationsrechnung.}
\end{figure}

Die daraus resultierende Variation $x_\lambda(\cdot) = x_*(\cdot) + \xi_\lambda(\cdot)$ besitzt dann die Eigenschaften
\begin{eqnarray*}
\|x_\lambda(\cdot) - x_*(\cdot) \|_\infty &=& \| \xi_\lambda(\cdot) \|_\infty = \lambda |v| \to 0 \qquad \mbox{f"ur } \lambda \to 0, \\
\|x_\lambda(\cdot) - x_*(\cdot)\|_{C_1} &\geq& \|\dot{\xi}_\lambda(\cdot)\|_\infty \geq \;\;|v| \not\to 0 \qquad \mbox{f"ur } \lambda \to 0.
\end{eqnarray*}
Der "Ubergang von den Richtungs- zu den Nadelvariationen vergr"o"sert die Menge der konkurrierenden Trajektorien.
Entsprechend muss der Optimalitätsbegriff angepasst werden.
Wir nennen $x_*(\cdot)$ ein starkes lokales Minimum\index{Minimum, starkes lokales!Variation@-- Variationsrechnung}
in der Aufgabe (\ref{WNV1})--(\ref{WNV2}), wenn sie die Randbedingungen (\ref{WNV2}) erf"ullt und
falls ein $\varepsilon > 0$ derart existiert, dass die Ungleichung
$J\big(x(\cdot)\big) \geq J\big(x_*(\cdot)\big)$
f"ur alle zul"assigen Funktionen $x(\cdot)$ mit $\| x(\cdot)-x_*(\cdot) \|_\infty < \varepsilon$ gilt. \\[2mm]
Einen anschaulichen Vergleich von starken und schwachen lokalen Minimalstellen liefert die folgende Aufgabe
der Klassischen Variationsrechnung (\!\!\cite{Ioffe}, S.\,107/108):
$$J\big(x(\cdot)\big) = \int_0^1 \dot{x}^3(t) \, dt \to \inf, \qquad x(0)=0, \quad x(1)=1.$$
Wir werden nun zeigen,
dass es in dieser Aufgabe eine eindeutig bestimmte L"osung der Euler-Lagrange-Gleichung gibt,
diese ein schwaches lokales Minimum darstellt, aber kein starkes.
Demgegen"uber existiert kein starkes lokales Minimum, denn $\inf J\big(x(\cdot)\big) = -\infty$
in jeder gleichm"a"sigen Umgebung der schwachen lokalen Minimalstelle. \\[2mm]
Zum schwachen lokalen Minimum:
Die Euler-Lagrangesche Gleichung als notwendige Bedingung f"ur ein schwaches lokales Minimum besitzt die Gestalt
$$-\frac{d}{dt} f_{\dot{x}}\big(t,x_*(t),\dot{x}_*(t)\big) + f_x\big(t,x_*(t),\dot{x}_*(t)\big) = -\frac{d}{dt}3\dot{x}_*^2(t)=0.$$
Eindeutige L"osung dieser Gleichung, die den Randbedingungen gen"ugt, ist $x_*(t)=t$. \\
Wir weisen jetzt nach, dass $x_*(\cdot)$ tats"achlich ein schwaches lokales Minimum liefert.
Es sei $x(\cdot) \in C_1([0,1],\R)$ mit $x(0)=x(1)=0$.
Dann ist $x_*(\cdot)+x(\cdot)$ zul"assig.
Wir erhalten
$$J\big(x_*(\cdot)+x(\cdot)\big) = \int_0^1 \bigg[\frac{d}{dt} \big(t+x(t)\big)\bigg]^3 \, dt
  = J\big(x_*(\cdot)\big) + \int_0^1 [3 \dot{x}^2(t)+\dot{x}^3(t)] \, dt.$$
Wenn daher $3 \dot{x}^2(t)+\dot{x}^3(t) \geq 0$ auf $[0,1]$ gilt, was insbesondere bei $\|x(\cdot)\|_{C_1}\leq 3$
wegen $|\dot{x}^3(t)|=|\dot{x}(t)| \cdot \dot{x}^2(t) \leq 3 \dot{x}^2(t)$ erf"ullt ist, so ist
$$J\big(x_*(\cdot)+x(\cdot)\big) \geq J\big(x_*(\cdot)\big).$$
D.\,h., dass $x_*(t)=t$ ein schwaches lokales Minimum in der Aufgabe liefert. \\[2mm]
Zum starken lokalen Minimum:
Es ergeben sich f"ur die Funktionen $x_n(\cdot)=x_*(\cdot)+\xi_n(\cdot)$,
$$\xi_n(0)=\xi_n(1)=0, \qquad
  \dot{\xi}_n(t)= \left\{\begin{array}{ll} -\sqrt{n}, & t \not\in [0,1/n], \\[1mm]
                                       \frac{\sqrt{n}}{n-1}, & t \in (1/n,1], \end{array}\right. \qquad
  n=2,3,...,$$
dass das Funktional $J$ f"ur gro"se $n$ asymptotisch gleich $-\sqrt{n}$ ausf"allt und deswegen
$$J\big(x_n(\cdot)\big) \to -\infty \qquad\mbox{f"ur } n \to \infty$$
gilt.
Da nun bez"uglich der Supremumsnorm aber
$$\|x_n(\cdot)-x_*(\cdot)\|_\infty = \|\xi_n(\cdot)\|_\infty = 1/\sqrt{n}$$
ist, liegt mit $x_*(\cdot)$ kein starkes lokales Minimum vor.
\newpage

Wir leiten nun die notwendige Weierstra"ssche Bedingung \index{Weierstra"ssche Bedingung} für ein starkes lokales Minimum in der
Grundaufgabe (\ref{WNV1})--(\ref{WNV2}) der Klassischen Variationsrechnung, \index{Grundaufgabe der Variationsrechnung}
$$J\big(x(\cdot)\big) = \int_{t_0}^{t_1} f\big(t,x(t),\dot{x}(t)\big) \, dt \to \inf; \qquad x(t_0)=x_0, \quad x(t_1)=x_1,$$
her. Wir folgen dabei \cite{Ioffe}.
Es seien die Voraussetzungen von Satz \ref{SatzELG} erfüllt.
Ferner bezeichne die stetig differenzierbare Funktion $x_*(\cdot)$ ein starkes lokales Minimum der Grundaufgabe.
Zu $\tau \in (t_0,t_1)$ und $0 < \lambda < \varepsilon$ mit $\tau+\varepsilon<t_1$
betrachten wir die Schar $x_\lambda(\cdot)=x_*(\cdot)+\xi_\lambda(\cdot)$,
wobei $\xi_\lambda(\cdot)$ nach (\ref{WNV3}) definiert ist.
Wir bilden die Funktion
\begin{eqnarray*}
\varphi(\lambda) &=& J\big(x_\lambda(\cdot)\big) = \int_{t_0}^{t_1} f\big(t,x_\lambda(t),\dot{x}_\lambda(t)\big) \, dt \\
                 &=& \int_{[t_0,t_1]\setminus [\tau,\tau+\varepsilon]} f\big(t,x_*(t),\dot{x}_*(t)\big) \, dt
                     + \int_{\tau}^{\tau+\varepsilon} f\big(t,x_\lambda(t),\dot{x}_\lambda(t)\big) \, dt \\
                 &=& J\big(x_*(\cdot)\big) - \int_{\tau}^{\tau+\varepsilon} f\big(t,x_*(t),\dot{x}_*(t)\big) \, dt
                     + \int_{\tau}^{\tau+\lambda} f\big(t,x_*(t)+(t-\tau)v,\dot{x}_*(t)+v\big) \, dt \\
                 & & + \int_{\tau+\lambda}^{\tau+\varepsilon} f\Big(t,x_*(t)+\lambda v-(t-\tau-\lambda)\frac{\lambda v}{\varepsilon-\lambda},
                        \dot{x}_*(t)-\frac{\lambda v}{\varepsilon-\lambda}\Big) \, dt.
\end{eqnarray*}
Deren Differenzenquotient besitzt in $\lambda_0=0$ für hinreichend kleine $\lambda >0$ die Gestalt
\begin{eqnarray*}
\frac{\varphi(\lambda)-\varphi(0)}{\lambda}
&=& \frac{1}{\lambda} \bigg[ \int_{\tau}^{\tau+\lambda} \Big[f\big(t,x_*(t)+(t-\tau)v,\dot{x}_*(t)+v\big)- f\big(t,x_*(t),\dot{x}_*(t)\big)\Big] \, dt \\
& &\hspace*{-25mm} + \int_{\tau+\lambda}^{\tau+\varepsilon} \Big[f\Big(t,x_*(t)+\lambda v-(t-\tau-\lambda)\frac{\lambda v}{\varepsilon-\lambda},
                        \dot{x}_*(t)-\frac{\lambda v}{\varepsilon-\lambda}\Big)-f\big(t,x_*(t),\dot{x}_*(t)\big)\Big] \, dt. \bigg]
\end{eqnarray*}
Im Grenzübergang $\lambda \to 0^+$ ergeben sich für die einzelnen Terme auf der rechten Seite
\begin{eqnarray*}
&& \lim_{\lambda \to 0^+} \frac{1}{\lambda} \int_{\tau}^{\tau+\lambda} \Big[f\big(t,x_*(t)+(t-\tau)v,\dot{x}_*(t)+v\big)- f\big(t,x_*(t),\dot{x}_*(t)\big)\Big] \, dt \\
&& \hspace*{15mm} = f\big(\tau,x_*(\tau),\dot{x}_*(\tau)+v\big)- f\big(\tau,x_*(\tau),\dot{x}_*(\tau)\big)
\end{eqnarray*}
und
\begin{eqnarray*}
& & \lim_{\lambda \to 0^+} \frac{1}{\lambda} \int_{\tau+\lambda}^{\tau+\varepsilon} \Big[f\Big(t,x_*(t)+\lambda v-(t-\tau-\lambda)\frac{\lambda v}{\varepsilon-\lambda},
                        \dot{x}_*(t)-\frac{\lambda v}{\varepsilon-\lambda}\Big)-f\big(t,x_*(t),\dot{x}_*(t)\big)\Big] \, dt \\
&& \hspace*{10mm} = \int_{\tau}^{\tau+\varepsilon} f_x\big(t,x_*(t),\dot{x}_*(t)\big) \Big(v -(t-\tau)\frac{v}{\varepsilon}\Big)  \, dt
   - \int_{\tau}^{\tau+\varepsilon} f_{\dot{x}}\big(t,x_*(t),\dot{x}_*\big) \frac{v}{\varepsilon} \, dt.
\end{eqnarray*}
In der letzten Gleichung gilt
$$\int_{\tau}^{\tau+\varepsilon} f_x\big(t,x_*(t),\dot{x}_*(t)\big) \cdot (t-\tau)\frac{v}{\varepsilon}  \, dt = O(\varepsilon).$$
Damit ergibt sich zusammenfassend:
\begin{eqnarray*}
\varphi'(0^+) &=& f\big(\tau,x_*(\tau),\dot{x}_*(\tau)+v\big)- f\big(\tau,x_*(\tau),\dot{x}_*(\tau)\big) \\
              & & + v \int_{\tau}^{\tau+\varepsilon} f_x\big(t,x_*(t),\dot{x}_*(t)\big) \, dt 
                  - \frac{v}{\varepsilon} \int_{\tau}^{\tau+\varepsilon} f_{\dot{x}}\big(t,x_*(t),\dot{x}_*(t)\big)  \, dt + O(\varepsilon).
\end{eqnarray*}
Der Beweis von Satz \ref{SatzELG} mit Hilfe des Fundamentallemma von du Bois-Reymond führte auf die Euler-Lagrangesche Gleichung in Integralform:
$$Q(t) - p(t) = \int_{t_0}^t f_x\big(s,x_*(s),\dot{x}_*(s)\big) \, ds - f_{\dot{x}}\big(t,x_*(t),\dot{x}_*(t)\big) = \mbox{konstant}.$$
Es sei $x_*(\cdot)$ eine Extremale der Grundaufgabe,
so ergibt sich unter Benutzung der Euler-Lagrangeschen Gleichung
$$\int_{\tau}^{\tau+\varepsilon} f_x\big(t,x_*(t),\dot{x}_*(t)\big) \, dt = 
  f_{\dot{x}}\big(\tau+\varepsilon,x_*(\tau+\varepsilon),\dot{x}_*(\tau+\varepsilon)\big) - f_{\dot{x}}\big(\tau,x_*(\tau),\dot{x}_*(\tau)\big)$$
und demzufolge
\begin{eqnarray*}
\varphi'(0^+) &=& f\big(\tau,x_*(\tau),\dot{x}_*(\tau)+v\big)- f\big(\tau,x_*(\tau),\dot{x}_*(\tau)\big) \\
              & & + v \big[f_{\dot{x}}\big(\tau+\varepsilon,x_*(\tau+\varepsilon),\dot{x}_*(\tau+\varepsilon)\big)
                         - f_{\dot{x}}\big(\tau,x_*(\tau),\dot{x}_*(\tau)\big) \big] \\
              & & - \frac{v}{\varepsilon} \int_{\tau}^{\tau+\varepsilon} f_{\dot{x}}\big(t,x_*(t),\dot{x}_*(t)\big)  \, dt + O(\varepsilon).
\end{eqnarray*}
Weiterhin muss $\varphi'(0^+) \geq 0$ gelten, da $x_*(\cdot)$ ein starkes lokales Minimum der Grundaufgabe darstellt.
Führen wir nun abschließend den Grenzübergang $\varepsilon \to 0^+$ durch,
so erhalten wir die Ungleichung
$$f\big(\tau,x_*(\tau),\dot{x}_*(\tau)+v\big)- f\big(\tau,x_*(\tau),\dot{x}_*(\tau)\big)- v  f_{\dot{x}}\big(\tau,x_*(\tau),\dot{x}_*(\tau)\big) \geq 0.$$
Als Weierstra"ssche $\mathscr{E}-$Funktion \index{Funktion, absolutstetige!Weier@--, Weierstra"ssche $\mathscr{E}-$Funktion} bezeichnet man die Funktion
$$\mathscr{E}(t,x,\dot{x},v)= f(t,x,v)-f(t,x,\dot{x})- (v-\dot{x})f_{\dot{x}}(t,x,\dot{x}).$$
Dann lautet die notwendige Weierstra"ssche Bedingung für ein starkes lokales Minimum der Grundaufgabe:

\begin{satz} \label{SatzWNV}
Unter den genannten Voraussetzungen an die Grundaufgabe und an $x_*(\cdot)$ sei $x_*(\cdot)$ eine Extremale der Grundaufgabe.
Liefert $x_*(\cdot)$ ein starkes lokales Minimum in der Grundaufgabe der Klassischen Variationsrechnung,
so ist für jeden Punkt $t \in (t_0,t_1)$ und für jedes $v \in \R$ die Ungleichung
$$\mathscr{E}\big(t,x_*(t),\dot{x}_*(t),v\big)= f\big(t,x_*(t),v\big)-f\big(t,x_*(t),\dot{x}_*(t)\big)
  - \big(v-\dot{x}_*(t)\big)f_{\dot{x}}\big(t,x_*(t),\dot{x}_*(t)\big) \geq 0$$
erfüllt.
\end{satz}

%% file: 2-2-PMPeinfach.tex
\subsection{Die elementare Aufgabe mit freiem Endpunkt} \label{AbschnittPMPeinfach}
\subsubsection{Formulierung der Aufgabe und das Maximumprinzip}
Über dem gegebenen Intervall $[t_0,t_1]$ betrachten wir zu $x_0 \in \R^n$ die Aufgabe 
\begin{eqnarray}
&&\label{PMPeinfach1} J\big(x(\cdot),u(\cdot)\big) = \int_{t_0}^{t_1} f\big(t,x(t),u(t)\big) \, dt \to \inf, \\
&&\label{PMPeinfach2} \dot{x}(t) = \varphi\big(t,x(t),u(t)\big), \quad x(t_0)=x_0,\\
&&\label{PMPeinfach3} u(t) \in U \subseteq \R^m, \quad U \not= \emptyset.
\end{eqnarray}
Die Aufgabe untersuchen wir bez"uglich $\big(x(\cdot),u(\cdot)\big) \in PC_1([t_0,t_1],\R^n) \times PC([t_0,t_1],U)$. \\[2mm]
Mit $\mathscr{D}^{\,\mathcal{S}}_{\rm Lip}$ bezeichnen wir die Menge aller Paare $\big(x(\cdot),u(\cdot)\big)$,
für die es ein $\gamma>0$ derart gibt,
dass die Abbildungen $f(t,x,u)$, $\varphi(t,x,u)$ auf der Menge aller $(t,x,u) \in \R \times \R^n \times \R^m$ mit
$$t \in [t_0,t_1], \qquad \|x-x(t)\| < \gamma, \qquad u \in \R^m$$
stetig in der Gesamtheit aller Variablen und stetig differenzierbar bezüglich $x$ sind. \\[2mm]
Das Paar $\big(x(\cdot),u(\cdot)\big) \in PC_1([t_0,t_1],\R^n) \times PC([t_0,t_1],U)$
hei"st ein zul"assiger Steuerungsprozess in der Aufgabe (\ref{PMPeinfach1})--(\ref{PMPeinfach3}),
falls $\big(x(\cdot),u(\cdot)\big)$ dem System (\ref{PMPeinfach2}) zu $x(t_0)=x_0$ gen"ugt.
Mit $\mathscr{D}^{\,\mathcal{S}}_{\rm adm}$ bezeichnen wir die Menge der zul"assigen Steuerungsprozesse. \\[2mm]
Ein zul"assiger Steuerungsprozess $\big(x_*(\cdot),u_*(\cdot)\big)$ ist eine
starke lokale Minimalstelle\index{Minimum, starkes lokales!elementar@-- elementare Aufgabe}
der Aufgabe (\ref{PMPeinfach1})--(\ref{PMPeinfach3}),
falls eine Zahl $\varepsilon > 0$ derart existiert, dass die Ungleichung 
$$J\big(x(\cdot),u(\cdot)\big) \geq J\big(x_*(\cdot),u_*(\cdot)\big)$$
f"ur alle $\big(x(\cdot),u(\cdot)\big) \in \mathscr{D}^{\,\mathcal{S}}_{\rm adm}$ mit $\|x(\cdot)-x_*(\cdot)\|_\infty < \varepsilon$ gilt. \\[2mm]
Weiterhin bezeichnet $H^{\mathcal{S}}: \R \times \R^n \times \R^m \times \R^n \times \R \to \R$ die Pontrjagin-Funktion
$$H^{\mathcal{S}}(t,x,u,p,\lambda_0) = \langle p, \varphi(t,x,u) \rangle - \lambda_0 f(t,x,u).$$

\begin{theorem}[Pontrjaginsches Maximumprinzip] \label{SatzPMPeinfach}
\index{Pontrjaginsches Maximumprinzip!Standard@-- Standardaufgabe} 
\index{Pontrjaginsches Maximumprinzip!elementar@-- elementare Aufgabe} 
Es sei $\big(x_*(\cdot),u_*(\cdot)\big) \in \mathscr{D}^{\,\mathcal{S}}_{\rm adm} \cap \mathscr{D}^{\,\mathcal{S}}_{\rm Lip}$. 
Ist $\big(x_*(\cdot),u_*(\cdot)\big)$ ein starkes lokales Minimum der Aufgabe (\ref{PMPeinfach1})--(\ref{PMPeinfach3}),
dann existiert eine Vektorfunktion $p(\cdot) \in PC_1([t_0,t_1],\R^n)$ derart, dass
\begin{enumerate}
\item[(a)] die adjungierten Gleichung
           \index{adjungierte Gleichung!Standard@-- Standardaufgabe}
           \index{adjungierte Gleichung!elementar@-- elementare Aufgabe}
           \begin{equation}\label{PMPeinfach4} 
           \dot{p}(t) =  -H_x^{\mathcal{S}}\big(t,x_*(t),u_*(t),p(t),1\big),
           \end{equation} 
\item[(b)] in $t=t_1$ die Transversalitätsbedingung
           \index{Transversalitätsbedingungen!Standard@-- Standardaufgabe}
           \index{Transversalitätsbedingungen!elementar@-- elementare Aufgabe}
           \begin{equation}\label{PMPeinfach5} p(t_1)=0 \end{equation}
\item[(c)] und in fast allen Punkten $t \in [t_0,t_1]$ die Maximumbedingung
           \index{Maximumbedingung!Standard@-- Standardaufgabe}
           \index{Maximumbedingung!elementar@-- elementare Aufgabe}
           \begin{equation}\label{PMPeinfach6} 
           H^{\mathcal{S}}\big(t,x_*(t),u_*(t),p(t),1\big) = \max_{u \in U}H^{\mathcal{S}}\big(t,x_*(t),u,p(t),1\big)
           \end{equation}
\end{enumerate}
erfüllt sind.
\end{theorem}


%% file: 2-21-Beweis.tex
\subsubsection{Der Beweis des Maximumprinzips} \label{AbschnittPMPBeweiseinfach}
Der nachstehende Beweis ist \cite{Ioffe} entnommen.
Da s"amtliche Abbildungen stetig und stetig differenzierbar bezüglich $x$ sind, und die Funktion $u_*(\cdot)$ dem Raum $PC([t_0,t_1],U)$ angeh"ort,
erf"ullt die adjungierte Gleichung (\ref{PMPeinfach4}) die Voraussetzungen von Lemma \ref{LemmaDGL3} und Lemma \ref{LemmaDGL5}.
Daher gibt es eine eindeutige L"osung $p(\cdot)$ der Gleichung (\ref{PMPeinfach4}) zur Randbedingung $p(t_1)=0$,
die dem Raum $PC_1([t_0,t_1],\R^n)$ angehören muss. \\[1mm]
Das Pontrjaginsche Maximumprinzip \ref{SatzPMPeinfach} ist mit $\lambda_0=1$ in der normalen Form angegeben,
der für die Aufgabe mit freiem rechten Endpunkt stets vorliegt.
Angenommen, es ist $\lambda_0=0$.
Dann besitzt die adjungierte Gleichung die Form
$$\dot{p}(t) =  -H_x^{\mathcal{S}}\big(t,x_*(t),u_*(t),p(t),0\big)= -\varphi_x^T\big(t,x_*(t),u_*(t)\big) p(t)$$
und es würde wegen der Transversalitätsbedingung $p(t_1)=0$ au"serdem $p(t) \equiv 0$ im Widerspruch zur
Nichttrivialität der Multiplikatoren $\big(\lambda_0,p(\cdot)\big)$ gelten. \\[1mm]
Es sei $\tau \in (t_0,t_1)$ ein Stetigkeitspunkt der Steuerung $u_*(\cdot)$.
Dann ist $u_*(\cdot)$ auch in einer gewissen hinreichend kleinen Umgebung von $\tau$ stetig und wir w"ahlen ein festes
$\lambda$ positiv und hinreichend klein, so dass sich $\tau-\lambda$ in dieser Umgebung befindet.
Weiter sei nun $v$ ein beliebiger Punkt aus $U$.
Wir setzen\index{Nadelvariation, einfache}\\
\begin{minipage}{0.59\textwidth}
$$u(t;v,\tau,\lambda) = u_{\lambda}(t) = 
  \left\{ \begin{array}{ll}
          u_*(t) & \mbox{ f"ur } t \not\in [\tau-\lambda,\tau), \\
          v      & \mbox{ f"ur } t     \in [\tau-\lambda,\tau), 
          \end{array} \right.$$
und es bezeichne $x_\lambda(\cdot)$, $x_\lambda(t)=x(t;v,\tau,\lambda)$, die eindeutige L"osung der Gleichung
$$\dot{x}(t) = \varphi\big(t,x(t),u_\lambda(t)\big), \qquad x(t_0)=x_0.$$
\end{minipage}
\begin{minipage}{0.4\textwidth}
\centering
\includegraphics[width=4cm]{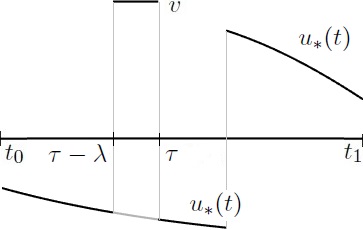}
\captionof{figure}[Einfache Nadelvariation der Steuerungstheorie]{Variation $u(t;v,\tau,\lambda)$.}
\end{minipage} \\[2mm]
Dann ist $x_{\lambda}(t) = x_*(t)$ f"ur $t_0 \leq t \leq  \tau - \lambda$.
F"ur $t \geq \tau$ betrachten wir den Grenzwert
$$y(t)=\lim_{\lambda \to 0^+}\frac{x_{\lambda}(t) - x_*(t)}{\lambda},$$
und werden zeigen, dass dieser existiert, die Funktion $y(\cdot)$ der Integralgleichung
\begin{equation} \label{BeweisPMPeinfach1}
y(t)=y(\tau) + \int_{\tau}^{t} \varphi_x\big(s,x_*(s),u_*(s)\big)\,y(s) \, ds
\end{equation}
zur Anfangsbedingung
\begin{equation} \label{BeweisPMPeinfach2}
y(\tau)=\varphi\big(\tau,x_*(\tau),v\big) - \varphi\big(\tau,x_*(\tau),u_*(\tau)\big)
\end{equation} 
gen"ugt und die Beziehung
\begin{equation} \label{BeweisPMPeinfach3}
\langle p(\tau) , y(\tau) \rangle = - \int_{\tau}^{t_1} \big\langle f_x\big(t,x_*(t),u_*(t)\big) , y(t) \big\rangle \, dt
\end{equation}
erf"ullt ist. 

\newpage
Herleitung von (\ref{BeweisPMPeinfach2}):
Nach Wahl des Punktes $\tau$ ist $\dot{x}_*(\cdot)$ in einer Umgebung dieses Punktes
stetig und f"ur hinreichend kleine positive $\lambda$ gelten
\begin{eqnarray*}
    x_*(\tau)
&=& x_*(\tau - \lambda) + \lambda 
    \varphi\big(\tau - \lambda,x_*(\tau - \lambda),u_*(\tau - \lambda)\big) + o(\lambda), \\
    x_{\lambda}(\tau)
&=& x_*(\tau - \lambda) + \lambda
    \varphi\big(\tau - \lambda,x_*(\tau - \lambda),v\big) + o(\lambda). \hspace{13,8mm}
\end{eqnarray*}
Hieraus folgt
$$\frac{x_{\lambda}(\tau) - x_*(\tau)}{\lambda}
= \varphi\big(\tau - \lambda,x_*(\tau - \lambda),v\big) 
- \varphi\big(\tau - \lambda,x_*(\tau - \lambda),u_*(\tau - \lambda)\big)
+ \frac{o(\lambda)}{\lambda}.$$
D.\,h., dass der Grenzwert
$$y(\tau)=\lim_{\lambda \to 0^+}\frac{x_{\lambda}(\tau) - x_*(\tau)}{\lambda}$$
existiert und gleich (\ref{BeweisPMPeinfach2}) ist. \\[2mm]
Herleitung von (\ref{BeweisPMPeinfach1}):
Auf dem Intervall $[\tau,t_1]$ gen"ugen sowohl $x_*(\cdot)$ als auch $x_{\lambda}(\cdot)$ der Gleichung
$$\dot{x}(t) = \varphi\big(t,x(t),u_*(t)\big).$$
Aus den Sätzen \ref{SatzEEglobal} und \ref{SatzDGLDifferenzierbarkeit} "uber die Stetigkeit und Differenzierbarkeit
der L"osung eines Differentialgleichungssystems in Abh"angigkeit
von den Anfangsdaten folgt, dass f"ur hinreichend kleine positive $\lambda$ die Vektorfunktionen
$x_{\lambda}(\cdot)$ auf $[\tau,t_1]$ definiert sind,
dass sie f"ur $\lambda \to 0^+$ gleichm"a"sig gegen $x_*(\cdot)$ konvergieren
und dass der Grenzwert
$$y(t) = \lim_{\lambda \to 0^+}\frac{x_{\lambda}(t) - x_*(t)}{\lambda}$$
f"ur jedes $t \in [\tau,t_1]$ existiert.
Weiterhin liefert Satz \ref{SatzDGLDifferenzierbarkeit}, dass
$y(\cdot)$ für alle $t \in [\tau,t_1]$ der Gleichung
$$y(t)=y(\tau) + \int_{\tau}^{t} \varphi_x\big(s,x_*(s),u_*(s)\big)\,y(s) \, ds$$
genügt, d.\,h. $y(\cdot)$ die Gleichung (\ref{BeweisPMPeinfach1}) zum Anfangswert (\ref{BeweisPMPeinfach2}) erfüllt. \\[2mm]
Herleitung von (\ref{BeweisPMPeinfach3}):
Mit (\ref{PMPeinfach4}) und (\ref{BeweisPMPeinfach1}) erhalten wir f"ur $t \geq \tau$
$$\frac{d}{dt} \langle p(t),y(t) \rangle = \big\langle f_x\big(t,x_*(t),u_*(t)\big) , y(t) \big\rangle.$$
Daher gilt für $t \in [\tau,t_1]$ die Gleichung
$$\langle p(t_1) ,y(t_1) \rangle - \langle p(t) ,y(t) \rangle = \int_{t}^{t_1} \big\langle f_x\big(s,x_*(s),u_*(s)\big) , y(s) \big\rangle \, ds.$$
F"ur $t = \tau$ folgt daraus mit $p(t_1) = 0$ insbesondere (\ref{BeweisPMPeinfach3}):
$$\langle p(\tau) ,y(\tau) \rangle = - \int_{\tau}^{t_1} \big\langle f_x\big(s,x_*(s),u_*(s)\big) , y(s) \big\rangle \, ds.$$
Beweisschluss:
Da $\big(x_*(\cdot),u_*(\cdot)\big)$ ein starkes lokales Minimum ist, ist f"ur alle hinreichend kleine und positive $\lambda$
$$\frac{J\big(x_\lambda(\cdot),u_\lambda(\cdot)\big) - J\big(x_*(\cdot),u_*(\cdot)\big)}{\lambda} \geq 0.$$
In diesem Ausdruck ergibt sich im Grenzwert $\lambda \to 0^+$:
\begin{eqnarray*}
&& \lim_{\lambda \to 0^+}\frac{1}{\lambda} \int_{\tau - \lambda}^{\tau} 
    \Big[ f\big(t,x_\lambda(t),u_\lambda(t)\big) - f\big(t,x_*(t),u_*(t)\big) \Big] \, dt \\
&& + \lim_{\lambda \to 0^+} \frac{1}{\lambda} \int_{\tau}^{t_1} 
    \Big[ f\big(t,x_\lambda(t),u_\lambda(t)\big) - f\big(t,x_*(t),u_*(t)\big) \Big] \, dt \\
&=& f\big(\tau,x_*(\tau),v\big) - f\big(\tau,x_*(\tau),u_*(\tau)\big)
    + \int_{\tau}^{t_1} \big\langle f_x\big(t,x_*(t),u_*(t)\big) , y(t) \big\rangle \, dt.
\end{eqnarray*}
Zusammenfassend haben wir
\begin{eqnarray*}
0 &\leq& \lim_{\lambda \to 0^+} \frac{J\big(x_\lambda(\cdot),u_\lambda(\cdot)\big)- J\big(x_*(\cdot),u_*(\cdot)\big)}{\lambda} \\
  &=   & f\big(\tau,x_*(\tau),v\big) - f\big(\tau,x_*(\tau),u_*(\tau)\big)
         + \int_{\tau}^{t_1} \big\langle f_x\big(t,x_*(t),u_*(t)\big),y(t) \big\rangle \, dt
\end{eqnarray*}
gezeigt.
Nach (\ref{BeweisPMPeinfach3}) besteht in dieser Ungleichung der Zusammenhang
$$-\langle p(\tau) , y(\tau) \rangle = \int_{\tau}^{t_1} \big\langle f_x\big(t,x_*(t),u_*(t)\big),y(t) \big\rangle \, dt$$
und die Ungleichung erhält die Gestalt
$$0 \leq f\big(\tau,x_*(\tau),v\big) - f\big(\tau,x_*(\tau),u_*(\tau)\big) -\langle p(\tau) , y(\tau) \rangle.$$
Bei Anwendung der Gleichungen (\ref{BeweisPMPeinfach2}) und (\ref{BeweisPMPeinfach3}) ergibt sich weiter
\begin{eqnarray*}
\lefteqn{\big\langle p(\tau) , \varphi\big(\tau,x_*(\tau),u_*(\tau)\big) \big\rangle - f\big(\tau,x_*(\tau),u_*(\tau)\big)} \\
&\geq& \big\langle p(\tau) , \varphi\big(\tau,x_*(\tau),v\big) \big\rangle - f\big(\tau,x_*(\tau),v\big).
\end{eqnarray*}
Damit folgt aus der Definition der Pontrjagin-Funktion mit $\lambda_0=1$:
\begin{equation*}
     H^{\mathcal{S}}\big(\tau,x_*(\tau),u_*(\tau),p(\tau),1\big) \geq H^{\mathcal{S}}\big(\tau,x_*(\tau),v,p(\tau),1\big).
\end{equation*}
Nun ist $\tau$ ein beliebiger Stetigkeitspunkt von $u_*(\cdot)$ und $v$ ein beliebiger Punkt der Menge $U$. 
Demzufolge ist die Beziehung (\ref{PMPeinfach6}) in allen Stetigkeitspunkten von $u_*(\cdot)$ wahr und damit ist das
Maximumprinzip bewiesen. \hfill $\blacksquare$

%% file: 2-22-Interpretation.tex
\subsubsection{\"Okonomische Deutung des Maximumprinzips} \label{AbschnittDeutung}
Die "okonomische Interpretation des Pontrjaginschen Maximumprinzips
\index{Pontrjaginsches Maximumprinzip!oekonomische@--, ökonomische Interpretation}
 wurde von Dorfman \cite{Dorfman} angeregt.
Wir geben eine Variante an,
die sich auf den Beweis von Theorem \ref{SatzPMPeinfach} im letzten Abschnitt bezieht. \\[2mm]
Der Kern des Beweises im letzten Abschnitt war,
die Auswirkungen der Nadelvariation
$$u_{\lambda}(t) = 
  \left\{ \begin{array}{ll}
          u_*(t) & \mbox{ f"ur } t \not\in [\tau-\lambda,\tau), \\
          v      & \mbox{ f"ur } t     \in [\tau-\lambda,\tau), 
          \end{array} \right.$$
auf das Steuerungsproblem (\ref{PMPeinfach1})--(\ref{PMPeinfach3}) zu untersuchen.
Die erzeugte marginale "Anderung
$$y(t)=\lim_{\lambda \to 0^+}\frac{x_{\lambda}(t) - x_*(t)}{\lambda}$$
des optimalen Zustandes $x_*(\cdot)$ bewirkt im Zielfunktional
\begin{equation} \label{Interpretation1}
\langle p(\tau) , y(\tau) \rangle = - \int_{\tau}^{t_1} \big\langle f_x\big(t,x_*(t),u_*(t)\big) , y(t) \big\rangle \, dt,
\end{equation}
wobei $p(\cdot)$ die L"osung der adjungierten Gleichung ist.
Au"serdem gilt f"ur den marginalen Profit die Beziehung
\begin{eqnarray}
\lefteqn{\lim_{\lambda \to 0^+} \frac{J\big(x_\lambda(\cdot),u_\lambda(\cdot)\big)- J\big(x_*(\cdot),u_*(\cdot)\big)}{\lambda}}
         \nonumber \\
  &=   & f\big(\tau,x_*(\tau),v\big) - f\big(\tau,x_*(\tau),u_*(\tau)\big)
         + \int_{\tau}^{t_1} \big\langle f_x\big(t,x_*(t),u_*(t)\big),y(t) \big\rangle \, dt \nonumber \\
  &=   & \label{Interpretation2} H^{\mathcal{S}}\big(\tau,x_*(\tau),u_*(\tau),p(\tau),1\big) - H^{\mathcal{S}}\big(\tau,x_*(\tau),v,p(\tau),1\big).
\end{eqnarray}
Die Entscheidung zum Zeitpunkt $\tau$ mit dem Parameter $v$ statt $u_*(\tau)$ zu steuern hat einen direkten und indirekten Effekt
(Feichtinger \& Hartl \cite{Feichtinger}, S.\,29):
\begin{enumerate}
\item[--] Der unmittelbare Effekt besteht darin, dass die Profitrate $f\big(\tau,x_*(\tau),v\big)$ erzielt wird.
\item[--] Die indirekte Wirkung manifestiert sich in der "Anderung der Kapitalbest"ande um $y(\tau)$.
          Der Kapitalstock wird durch die in $\tau$ getroffene Entscheidung transformiert,
          was bemessen mit der Adjungierten $p(\cdot)$ in (\ref{Interpretation1}) ausgedr"uckt wird.
\end{enumerate}
Die Gleichung (\ref{Interpretation1}) beschreibt damit die Opportunit"atskosten,
die sich durch die Entscheidung $v$ zum Zeitpunkt $\tau$ ergeben.
Wegen der Darstellung der Opportunit"atskosten in der Form $\langle p(\tau) , y(\tau) \rangle$ in (\ref{Interpretation1}),
wird die Adjungierte $p(\cdot)$ in der "Okonomie h"aufig als Schattenpreis\index{Schattenpreis} bezeichnet. \\
Ferner bemisst die Funktion $v \to H^{\mathcal{S}}\big(\tau,x_*(\tau),v,p(\tau),1\big)$ nach Gleichung (\ref{Interpretation2}) die Profitrate aus
direktem und indirektem Gewinn,
die sich aus der "Anderung der Kapitalbest"ande ergibt.
Die Maximumbedingung sagt daher aus,
dass die Instrumente zu jedem Zeitpunkt so eingesetzt werden sollen,
dass die totale Profitrate maximal wird.

%% file: 2-3-Aufgabenstellung.tex
\subsection{Die Aufgabenstellung f\"ur ein starkes lokales Minimum}
Wir betrachten als Steuerungsproblem die Aufgabe
\begin{eqnarray}
&& \label{PMP1} J\big(x(\cdot),u(\cdot)\big) = \int_{t_0}^{t_1} f\big(t,x(t),u(t)\big) \, dt \to \inf, \\
&& \label{PMP2} \dot{x}(t) = \varphi\big(t,x(t),u(t)\big), \\
&& \label{PMP3} h_0\big(x(t_0)\big)=0, \qquad h_1\big(x(t_1)\big)=0, \\
&& \label{PMP4} u(t) \in U \subseteq \R^m, \quad U\not= \emptyset, \\
&& \label{PMP5} g_j\big(t,x(t)\big) \leq 0 \quad \mbox{f"ur alle } t \in [t_0,t_1], \quad j=1,...,l.
\end{eqnarray}
Darin gelten f"ur die Abbildungen
$$f:\R \times \R^n \times \R^m \to \R, \quad \varphi:\R \times \R^n \times \R^m \to \R^n, \quad
  h_i:\R^n \to \R^{s_i}, \; i=0,1,$$
sowie
$$g_j:\R \times \R^n \to \R, \quad j=1,...,l.$$
Die Aufgabe (\ref{PMP1})--(\ref{PMP5}) betrachten wir bez"uglich der Paare
$$\big(x(\cdot),u(\cdot)\big) \in W^1_\infty([t_0,t_1],\R^n) \times L_\infty([t_0,t_1],U).$$
Wir nennen die Trajektorie $x(\cdot)$ eine L"osung der Gleichung (\ref{PMP2}) zur Steuerung $u(\cdot)$,
falls $x(\cdot)$ auf $[t_0,t_1]$ definiert ist und die Dynamik im Sinn von Carath\'eodory l"ost. \\[2mm]
Mit $\mathscr{B}^{\,\mathcal{S}}_{\rm Lip}$ bezeichnen wir die Menge aller Paare $\big(x(\cdot),u(\cdot)\big)$,
für die es ein $\gamma>0$ derart gibt, dass die Abbildungen
$f(t,x,u)$, $\varphi(t,x,u)$, $h_i(x_i)$ und $g_j(t,x)$ auf der Menge aller Punkte
$(t,x,x_0,x_1,u) \in \R \times \R^n \times \R^n \times \R^n \times \R^m$ mit
$$t_0 \leq t\leq t_1, \quad \|x-x(t)\| < \gamma, \quad \|x_0-x(t_0)\| < \gamma, \quad \|x_1-x(t_1)\| < \gamma, \quad u \in \R^m$$
stetig in der Gesamtheit der Variablen und stetig differenzierbar bez"uglich $x,x_0,x_1$ sind. \\[2mm]
Das Paar $\big(x(\cdot),u(\cdot)\big) \in W^1_\infty([t_0,t_1],\R^n) \times L_\infty([t_0,t_1],U)$
hei"st ein zul"assiger Steuerungsprozess in der Aufgabe (\ref{PMP1})--(\ref{PMP5}),
falls $\big(x(\cdot),u(\cdot)\big)$ dem System (\ref{PMP2}) gen"ugt, sowie die Randbedingungen (\ref{PMP3}) und
Zustandsbeschr"ankungen (\ref{PMP5}) erf"ullt.
Die Menge $\mathscr{B}^{\,\mathcal{S}}_{\rm adm}$ bezeichnet die Menge der zul"assigen Steuerungsprozesse. \\[2mm]
Ein zul"assiger Steuerungsprozess $\big(x_*(\cdot),u_*(\cdot)\big)$ ist eine starke lokale
Minimalstelle\index{Minimum, starkes lokales!Standardl@-- Standardaufgabe}
der Aufgabe (\ref{PMP1})--(\ref{PMP5}),
falls eine Zahl $\varepsilon > 0$ derart existiert, dass die Ungleichung 
$$J\big(x(\cdot),u(\cdot)\big) \geq J\big(x_*(\cdot),u_*(\cdot)\big)$$
f"ur alle $\big(x(\cdot),u(\cdot)\big) \in \mathscr{B}^{\,\mathcal{S}}_{\rm adm}$ mit 
$\|x(\cdot)-x_*(\cdot) \|_\infty < \varepsilon$ gilt. \\[2mm]
F"ur die Aufgabe (\ref{PMP1})--(\ref{PMP5}) bezeichnet
$H^{\mathcal{S}}: \R \times \R^n \times \R^m \times \R^n \times \R \to \R$ die Pontrjagin-Funktion,
die wir wie folgt definieren:
$$H^{\mathcal{S}}(t,x,u,p,\lambda_0) = \langle p, \varphi(t,x,u) \rangle - \lambda_0 f(t,x,u).$$

%% file: 2-4-PontrjaginAufgabe.tex
\subsection{Das Pontrjaginsche Maximumprinzip} \label{AbschnittPMPStandard}
\subsubsection{Notwendige Optimalit\"atsbedingungen}
\begin{theorem}[Pontrjaginsches Maximumprinzip] \label{SatzPMP}
\index{Pontrjaginsches Maximumprinzip!Standard@-- Standardaufgabe} 
Es sei $\big(x_*(\cdot),u_*(\cdot)\big) \in \mathscr{B}^{\,\mathcal{S}}_{\rm adm} \cap \mathscr{B}^{\,\mathcal{S}}_{\rm Lip}$. 
Ist $\big(x_*(\cdot),u_*(\cdot)\big)$ ein starkes lokales Minimum der Aufgabe (\ref{PMP1})--(\ref{PMP4}),
dann existieren nicht gleichzeitig verschwindende Multiplikatoren $\lambda_0 \geq 0$,
$p(\cdot) \in W^1_\infty([t_0,t_1],\R^n)$ und $l_i \in \R^{s_i}$, $i=0,1$, derart, dass
\begin{enumerate}
\item[(a)] die adjungierte Gleichung
           \index{adjungierte Gleichung!Standard@-- Standardaufgabe}
           \begin{equation}\label{SatzPMP1}
           \dot{p}(t)=-H_x^{\mathcal{S}}\big(t,x_*(t),u_*(t),p(t),\lambda_0\big),
           \end{equation}
\item[(b)] die Transversalitätsbedingungen
           \index{Transversalitätsbedingungen!Standard@-- Standardaufgabe}
           \begin{equation}\label{SatzPMP2} 
           p(t_0)= {h_0'}^T\big(x_*(t_0)\big)l_0, \qquad p(t_1)=-{h_1'}^T\big(x_*(t_1)\big)l_1
           \end{equation}
\item[(c)] und in fast allen Punkten $t \in [t_0,t_1]$ die Maximumbedingung
           \index{Maximumbedingung!Standard@-- Standardaufgabe}
           \begin{equation}\label{SatzPMP3}
           H^{\mathcal{S}}\big(t,x_*(t),u_*(t),p(t),\lambda_0\big) = \max_{u \in U} H^{\mathcal{S}}\big(t,x_*(t),u,p(t),\lambda_0\big)
           \end{equation}
\end{enumerate}
erfüllt sind.
\end{theorem}

\begin{bemerkung} \label{BemerkungPMP} {\rm
F"ur das Pontrjaginsche Maximumprinzip l"asst sich an dieser Stelle die Bemerkung \ref{BemerkungFreierEndpunkt} "uber die
Transversalit"atsbedingungen und "uber das Eintreten des normalen Falles bei fehlenden Randbedingungen "ubernehmen. \hfill $\square$}
\end{bemerkung}

Das Pontrjaginsche Maximumprinzip enthält die vorgestellten notwendigen Optimalitätsbedingungen der Klassischen Variationsrechnung
und des Schwachen Optimalitätsprinzips. \\[2mm]
In der Grundaufgabe der Klassischen Variationsrechnung über dem Raum $C_1([t_0,t_1],\R^n)$,
$$J\big(x(\cdot)\big) = \int_{t_0}^{t_1} f\big(t,x(t),\dot{x}(t)\big) \, dt \to \inf; \qquad x(t_0)=x_0, \quad x(t_1)=x_1,$$
setzen wir $\dot{x}(t)=u(t)$ und schreiben sie in der Form
$$J\big(x(\cdot),u(\cdot)\big) = \int_{t_0}^{t_1} f\big(t,x(t),u(t)\big) \, dt \to \inf; \quad \dot{x}(t)=u(t), \quad x(t_0)=x_0, \; x(t_1)=x_1.$$
Die Pontrjagin-Funktion lautet dann $H^{\mathcal{S}}(t,x,u,p,\lambda_0)=pu-\lambda_0 f(t,x,u)$,
die adjungierte Gleichung besitzt die Gestalt
$$\dot{p}(t)=\lambda_0 f_x\big(t,x_*(t),u_*(t)\big)=-H^{\mathcal{S}}_x\big(t,x_*(t),u_*(t),p(t),\lambda_0\big)$$
und die Maximumbedingung erhält für $u \in \R$ die Form
$$H^{\mathcal{S}}_u\big(t,x_*(t),u_*(t),p(t),\lambda_0\big)=p(t)-\lambda_0 f_u\big(t,x_*(t),u_*(t)\big)=0.$$
Es sei nun die stetig differenzierbare Funktion $x_*(\cdot)$ ein lokales Minimum in der Grundaufgabe.
Da in diesem Fall $u_*(\cdot)$ stetig ist, so würde für $\lambda_0=0$ nach der eben umgeformten Maximumbedingung die Adjungierte $p(\cdot)$ 
im Widerspruch zur Nichttrivialität der Multiplikatoren ebenfalls identisch verschwinden.
Daher kann man $\lambda_0=1$ annehmen und gelangt zur Euler-Lagrangeschen Gleichung
$$\dot{p}(t)-\frac{d}{dt} f_u\big(t,x_*(t),u_*(t)\big)
  =f_x\big(t,x_*(t),u_*(t)\big)=\frac{d}{dt} f_u\big(t,x_*(t),u_*(t)\big)=0.$$
Aus der Maximumbedingung ergibt sich weiterhin 
$$\max_{u \in U} \big[ p(t)u- f\big(t,x_*(t),u\big)\big] = p(t)u_*(t)- f\big(t,x_*(t),u_*(t)\big)$$
für alle $t$, da $u_*(\cdot)=\dot{x}_*(\cdot)$ stetig ist.
Beachten wir nun $p(t) = f_u\big(t,x_*(t),u_*(t)\big)$ wegen $H^{\mathcal{S}}_u=0$,
so gewinnen wir die Gültigkeit der Weierstra"sschen Bedingung
$$\mathscr{E}\big(t,x_*(t),u_*(t),u\big)= f\big(t,x_*(t),u\big)-f\big(t,x_*(t),u_*(t)\big)
  - \big(u-u_*(t)\big)f_u\big(t,x_*(t),u_*(t)\big) \geq 0$$
für alle $t \in (t_0,t_1)$ und alle $u \in \R^n$. \\[2mm]
Der Bezug des Pontrjaginschen Maximumprinzips zu den Optimalitätsbedingungen in Satz \ref{SatzLA} für die Lagrange-Aufgabe ist der
Bemerkung \ref{BemerkungSatzLA} unmittelbar zu entnehmen.
Entsprechend ist der Zusammenhang zum Schwachen Optimalitätsprinzip in Theorem \ref{SatzSOP} direkt einsichtig.

%% file: 2-41-Beweis.tex
\subsubsection{Der Nachweis der notwendigen Optimalit\"atsbedingungen} \label{AbschnittBeweisPMP}
Wir betrachten f"ur $\big(x(\cdot),u(\cdot)\big) \in C([t_0,t_1],\R^n) \times L_\infty([t_0,t_1],\R^m)$ die Abbildungen
\begin{eqnarray*}
J\big(x(\cdot),u(\cdot)\big) &=& \int_{t_0}^{t_1} f\big(t,x(t),u(t)\big) \, dt, \\
F\big(x(\cdot),u(\cdot)\big)(t) &=& x(t) -x(t_0) -\int_{t_0}^t \varphi\big(s,x(s),u(s)\big) \, ds, \quad t \in [t_0,t_1],\\
H_i\big(x(\cdot)\big) &=& h_i\big(x(t_i)\big), \quad i=0,1.
\end{eqnarray*}
Da $x(\cdot)$ zu $C([t_0,t_1],\R^n)$ geh"ort, gilt f"ur diese Abbildungen
\begin{eqnarray*}
J &:& C([t_0,t_1],\R^n) \times L_\infty([t_0,t_1],\R^m) \to \R, \\
F &:& C([t_0,t_1],\R^n) \times L_\infty([t_0,t_1],\R^m) \to C_0([t_0,t_1],\R^n), \\
H_i &:& C([t_0,t_1],\R^n) \to \R^{s_i}, \quad i=0,1.
\end{eqnarray*}

Wir setzen $\mathscr{F}=(F,H_0,H_1)$ und pr"ufen f"ur die Extremalaufgabe
\begin{equation} \label{ExtremalaufgabePMP}
J\big(x(\cdot),u(\cdot)\big) \to \inf, \qquad \mathscr{F}\big(x(\cdot),u(\cdot)\big)=0, \qquad
u(\cdot) \in L_\infty([t_0,t_1],U)
\end{equation}
im Punkt $\big(x_*(\cdot),u_*(\cdot)\big) \in \mathscr{B}^{\,\mathcal{S}}_{\rm Lip}$,
wobei wir $x_*(\cdot)$ als Element des Raumes $C([t_0,t_1],\R^n)$ auffassen,
die Voraussetzungen von Theorem \ref{SatzExtremalprinzipStark}:

\begin{enumerate}
\item[(A$_1$)] F"ur jedes $u(\cdot) \in L_\infty([t_0,t_1],U)$ ist die Abbildung $x(\cdot) \to J\big(x(\cdot),u(\cdot)\big)$
               nach Beispiel \ref{DiffZielfunktional2} im Punkt $x_*(\cdot)$ Fr\'echet-differenzierbar.
\item[(A$_2$)] Die Abbildung $F$ ist die Summe der Abbildung $x(\cdot) \to x(t)$ und der Abbildung
               $$\big(x(\cdot),u(\cdot)\big) \to -\int_{t_0}^t \varphi\big(s,x(s),u(s)\big) \, ds.$$
               Im Beispiel \ref{DiffDynamik2} ist die Fr\'echet-Differenzierbarkeit der zweiten Abbildung im Punkt $x_*(\cdot)$
               f"ur jedes $u(\cdot) \in L_\infty([t_0,t_1],U)$ nachgewiesen.
               F"ur die Abbildungen $H_i$ ist die stetige Differenzierbarkeit offensichtlich.
\item[(B)] Die endliche Kodimension des Operators $\mathscr{F}_x\big(x_*(\cdot),u_*(\cdot)\big)$ folgt wie in
           Abschnitt \ref{AbschnittBeweisSOP}.
\item[(C)] Der Nachweis dieser Voraussetzungen ist im Anhang \ref{AnhangNV},
           speziell in Lemma \ref{LemmaEigenschaftNadelvariation1} und Lemma \ref{LemmaEigenschaftNadelvariation2}
           im Abschnitt \ref{AbschnittEigenschaftenNadelvariation}, dargestellt.
\end{enumerate}

Zur Extremalaufgabe (\ref{ExtremalaufgabePMP}) definieren wir auf
$$C([t_0,t_1],\R^n) \times L_\infty([t_0,t_1],\R^m) \times \R \times C_0^*([t_0,t_1],\R^n) \times \R^{s_0} \times \R^{s_1}$$
die Lagrange-Funktion $\mathscr{L}=\mathscr{L}\big(x(\cdot),u(\cdot),\lambda_0,y^*,l_0,l_1\big)$,
$$\mathscr{L}= \lambda_0 J\big(x(\cdot),u(\cdot)\big)+ \big\langle y^*, F\big(x(\cdot),u(\cdot)\big) \big\rangle
                         +l_0^T H_0\big(x(\cdot)\big)+l_1^T H_1\big(x(\cdot)\big).$$
Ist $\big(x_*(\cdot),u_*(\cdot)\big)$ eine starke lokale Minimalstelle der Aufgabe (\ref{ExtremalaufgabePMP}),
dann existieren nach Theorem \ref{SatzExtremalprinzipStark}
nicht gleichzeitig verschwindende Lagrangesche Multiplikatoren 
$$\lambda_0 \geq 0, \qquad y^* \in C_0^*([t_0,t_1],\R^n), \qquad l_i \in \R^{s_i}$$
derart,
dass folgende Bedingungen gelten:
\begin{enumerate}
\item[(a)] Die Lagrange-Funktion besitzt bez"uglich $x(\cdot)$ in $x_*(\cdot)$ einen station"aren Punkt, d.\,h.
          \begin{equation}\label{SatzPMPLMR1}
          \mathscr{L}_x\big(x_*(\cdot),u_*(\cdot),\lambda_0,y^*,l_0,l_1\big)=0;
          \end{equation}         
\item[(b)] Die Lagrange-Funktion erf"ullt bez"uglich $u(\cdot)$ in $u_*(\cdot)$ die Minimumbedingung
           \begin{equation}\label{SatzPMPLMR2}
           \mathscr{L}\big(x_*(\cdot),u_*(\cdot),\lambda_0,y^*,l_0,l_1\big)
           = \min_{u(\cdot) \in L_\infty([t_0,t_1],U)} \mathscr{L}\big(x_*(\cdot),u(\cdot),\lambda_0,y^*,l_0,l_1\big).
           \end{equation}
\end{enumerate}
Zu $y^* \in C_0^*([t_0,t_1],\R^n)$ gibt es ein eindeutig bestimmtes regul"ares Borelsches Vektorma"s $\mu$ derart,
dass nach (\ref{SatzPMPLMR1}) f"ur alle $x(\cdot) \in C([t_0,t_1],\R^n)$ die Variationsgleichung
\begin{eqnarray}
0 &=& \lambda_0 \int_{t_0}^{t_1} \big\langle f_x\big(t,x_*(t),u_*(t)\big),x(t) \big\rangle \, dt
      + \big\langle l_0, h_0'\big(x_*(t_0)\big) x(t_0) \big\rangle + \big\langle l_1, h_1'\big(x_*(t_1)\big) x(t_1) \big\rangle
      \nonumber \\
  & & \label{BeweisschlussPMP1}
      + \int_{t_0}^{t_1} \bigg[ x(t)-x(t_0) - \int_{t_0}^{t} \varphi_x\big(s,x_*(s),u_*(s)\big) x(s)\, ds \bigg]^T d\mu(t)
\end{eqnarray}
erf"ullt ist.
In der Gleichung (\ref{BeweisschlussPMP1}) "andern wir die Integrationsreihenfolge im zweiten Summanden und bringen sie in 
die Form
\begin{eqnarray}
0 &\!=&\! \int_{t_0}^{t_1} \big\langle \lambda_0 f_x\big(t,x_*(t),u_*(t)\big)
                          - \varphi_x^T\big(t,x_*(t),u_*(t)\big) \int_{t}^{t_1} d\mu(s) , x(t) \big\rangle \, dt
      + \int_{t_0}^{t_1} [x(t)]^T \, d\mu(t) \nonumber \\
  & & \label{BeweisschlussPMP2} 
      + \Big\langle {h_0'}^T\big(x_*(t_0)\big)l_0 - \int_{t_0}^{t_1} d\mu(t) , x(t_0) \Big\rangle
      + \langle {h_1'}^T\big(x_*(t_1)\big)l_1 , x(t_1) \rangle.
\end{eqnarray}
Setzen wir $p(t)=\displaystyle \int_t^{t_1} d\mu(s)$ in (\ref{BeweisschlussPMP2}),
so folgen wegen der eindeutigen Darstellung eines stetigen linearen Funktionals im Raum $C([t_0,t_1],\R^n)$
\begin{eqnarray*}
p(t) &=& -{h_1'}^T\big(x_*(t_1)\big)l_1
         + \int_{t}^{t_1} \Big( \varphi_x^T\big(s,x_*(s),u_*(s)\big) p(s) - \lambda_0 f_x\big(s,x_*(s),u_*(s)\big) \Big) \, ds, \\
p(t_0) &=& {h_0'}^T\big(x_*(t_0)\big)l_0, \qquad p(t_1)= -{h_1'}^T\big(x_*(t_1)\big)l_1.
\end{eqnarray*}
Damit sind (\ref{SatzPMP1}) und (\ref{SatzPMP2}) gezeigt.
Die Beziehung (\ref{SatzPMPLMR2}) ist "aquivalent zu
\begin{eqnarray*} 
\lefteqn{\int_{t_0}^{t_1} \lambda_0 f\big(t,x_*(t),u_*(t)\big) \, dt
    - \int_{t_0}^{t_1}  \left[ \int_{t_0}^{t} \varphi\big(s,x_*(s),u_*(s)\big) \, ds \right]^T d\mu(t)} \\
&=& \min_{u(\cdot) \in L_\infty ([t_0,t_1],U)}
    \left\{\int_{t_0}^{t_1} \lambda_0 f\big(t,x_*(t),u(t)\big) \, dt
            - \int_{t_0}^{t_1} \bigg[ \int_{t_0}^{t} \varphi\big(s,x_*(s),u(s)\big)\, ds \bigg]^T d\mu(t) \right\}.
\end{eqnarray*}
Es gilt aber
\begin{eqnarray*}
    \int_{t_0}^{t_1} \bigg[ \int_{t_0}^t \varphi\big(s,x_*(s),u(s)\big) \, ds \bigg]^T d\mu(t)
&=& \int_{t_0}^{t_1} \bigg( \big[\varphi\big(t,x_*(t),u(t)\big)\big]^T \int_t^{t_1} d\mu(s) \bigg) \, dt \\
&=& \int_{t_0}^{t_1} \big\langle p(t) , \varphi\big(t,x_*(t),u(t)\big) \big\rangle \, dt.
\end{eqnarray*}
Aus diesen beiden Beziehungen folgt
$$\int_{t_0}^{t_1} H^{\mathcal{S}}\big(t,x_*(t),u_*(t),p(t),\lambda_0\big) \, dt \geq \int_{t_0}^{t_1} H^{\mathcal{S}}\big(t,x_*(t),u(t),p(t),\lambda_0\big) \, dt$$
f"ur alle $u(\cdot) \in L_\infty([t_0,t_1],U)$.
Daraus folgt abschlie"send die Maximumbedingung (\ref{SatzPMP3}) via Standardtechniken f"ur Lebesguesche Punkte. \hfill $\blacksquare$

%% file: 2-42-Hinreichend.tex
\subsubsection{Hinreichende Bedingungen nach Arrow} \label{AbschnittHBPMP}
Unser\index{hinreichende Bedingungen, Arrow!Standard@-- Standardaufgabe}
Vorgehen zur Herleitung der hinreichenden Bedingungen basiert auf der Darstellung in 
Seierstad \& Syds\ae ter \cite{Seierstad},
die wir um Ausf"uhrungen in Aseev \& Kryazhimskii \cite{AseKry} erg"anzt haben.
Wir betrachten das Steuerungsproblem
\begin{eqnarray}
&& \label{HBPMP1} J\big(x(\cdot),u(\cdot)\big) = \int_{t_0}^{t_1} f\big(t,x(t),u(t)\big) \, dt \to \inf, \\
&& \label{HBPMP2} \dot{x}(t) = \varphi\big(t,x(t),u(t)\big), \\
&& \label{HBPMP3} x(t_0)=x_0, \qquad x(t_1)=x_1, \\
&& \label{HBPMP4} u(t) \in U \subseteq \R^m, \quad U\not= \emptyset.
\end{eqnarray}

In der Aufgabenstellung schlie"sen wir wieder den Fall nicht aus,
dass durch die Randbedingungen (\ref{HBPMP3}) gewisse Komponenten der Punkte $x_0$ und $x_1$ nicht fest vorgegeben,
sondern ohne Einschr"ankung sind. \\[2mm]
Wir definieren die Menge $V^{\mathcal{S}}_\gamma(t)=\{ x \in \R^n \,|\, \|x-x_*(t)\| < \gamma\}$.
Au"serdem bezeichnet 
\begin{equation} \label{PMPHamilton}
\mathscr{H}^{\mathcal{S}}(t,x,p) = \sup_{u \in U} H^{\mathcal{S}}(t,x,u,p,1)
\end{equation}
die Hamilton-Funktion $\mathscr{H}^{\mathcal{S}}$ im normalen Fall.

\begin{theorem} \label{SatzHBPMP}
In der Aufgabe (\ref{HBPMP1})--(\ref{HBPMP4}) sei
$\big(x_*(\cdot),u_*(\cdot)\big) \in \mathscr{B}^{\mathcal{S}}_{\rm Lip} \cap \mathscr{B}^{\mathcal{S}}_{\rm adm}$ 
und es sei $p(\cdot) \in W^1_\infty([t_0,t_1],\R^n)$. Ferner gelte:
\begin{enumerate}
\item[(a)] Das Tripel $\big(x_*(\cdot),u_*(\cdot),p(\cdot)\big)$
           erf"ullt (\ref{SatzPMP1})--(\ref{SatzPMP3}) in Theorem \ref{SatzPMP} mit $\lambda_0=1$.        
\item[(b)] F"ur jedes $t \in [t_0,t_1]$ ist die Funktion $\mathscr{H}^{\mathcal{S}}\big(t,x,p(t)\big)$ konkav in $x$ auf $V^{\mathcal{S}}_\gamma(t)$.
\end{enumerate}
Dann ist $\big(x_*(\cdot),u_*(\cdot)\big)$ ein starkes lokales Minimum der Aufgabe (\ref{HBPMP1})--(\ref{HBPMP4}).
\end{theorem}

{\bf Beweis} Es sei $t \in [t_0,t_1]$ gegeben.
Da die Abbildung $x \to -\mathscr{H}^{\mathcal{S}}\big(t,x,p(t)\big)$ auf $V^{\mathcal{S}}_\gamma(t)$ konvex ist, ist die Menge
$Z=\big\{ (\alpha,x) \in \R \times \R^n \,\big|\, x \in V^{\mathcal{S}}_\gamma(t), \alpha \geq -\mathscr{H}^{\mathcal{S}}\big(t,x,p(t)\big) \big\}$
konvex und besitzt ein nichtleeres Inneres.
Wir setzen $\alpha_*= -\mathscr{H}^{\mathcal{S}}\big(t,x_*(t),p(t)\big)$.
Dann ist $\big(\alpha_*,x_*(t)\big)$ mit ein Randpunkt der Menge $Z$.
Daher existiert nach dem Trennungssatz für konvexe Mengen ein nichttrivialer Vektor $\big(a_0(t),a(t)\big) \in \R \times \R^n$ mit
\begin{equation} \label{BeweisHBPMP1}
a_0(t) \alpha + \langle a(t),x\rangle \geq a_0(t) \alpha_* + \langle a(t),x_*(t)\rangle \quad \mbox{ f"ur alle } (\alpha,x) \in Z.
\end{equation}
Es ist $x_*(t)$ ein innerer Punkt der Menge $V^{\mathcal{S}}_\gamma(t)$.
Weiterhin folgt aus den elementaren Eigenschaften konvexer Funktionen,
dass $x \to -\mathscr{H}^{\mathcal{S}}\big(t,x,p(t)\big)$ in $V^{\mathcal{S}}_\gamma(t)$ stetig ist,
da sie auf $V^{\mathcal{S}}_\gamma(t)$ konvex und nach oben durch $x \to -H^{\mathcal{S}}\big(t,x,u_*(t),p(t)\big)$ beschr"ankt ist. \\
Deswegen existiert ein $\delta>0$ mit $x_*(t)+\xi \in V^{\mathcal{S}}_\gamma(t)$ und $\big(\alpha_*+1,x_*(t)+\xi\big) \in Z$
f"ur alle $\|\xi\| < \delta$.
Aus (\ref{BeweisHBPMP1}) folgt daher $a_0(t) +\langle a(t),\xi\rangle \geq 0$ f"ur alle $\|\xi\| < \delta$.
Dies zeigt $a_0(t) >0$ und wir k"onnen ohne Einschr"ankung $a_0(t)=1$ annehmen.
Wiederum (\ref{BeweisHBPMP1}) liefert damit
\begin{equation} \label{BeweisHBPMP2}
\langle a(t),x-x_*(t)\rangle \geq  \mathscr{H}^{\mathcal{S}}\big(t,x,p(t)\big) - \mathscr{H}^{\mathcal{S}}\big(t,x_*(t),p(t)\big)
  \quad \mbox{ f"ur alle } x \in V^{\mathcal{S}}_\gamma(t).
\end{equation}
Es sei nun $t \in [t_0,t_1]$ so gew"ahlt,
dass die Maximumbedingung (\ref{PMPeinfach6}) zu diesem Zeitpunkt erf"ullt ist.
Dann folgt aus (\ref{BeweisHBPMP2}), dass
\begin{eqnarray*}
\langle a(t),x-x_*(t)\rangle
&\geq& \sup_{u \in U} H^{\mathcal{S}}\big(t,x,u,p(t)\big) - H^{\mathcal{S}}\big(t,x_*(t),u_*(t),p(t),1\big) \\
&\geq& H^{\mathcal{S}}\big(t,x,u_*(t),p(t),1\big) - H^{\mathcal{S}}\big(t,x_*(t),u_*(t),p(t),1\big)
\end{eqnarray*}
f"ur alle $x \in V^{\mathcal{S}}_\gamma(t)$ gilt.
Wir setzen
\begin{eqnarray*}
\Phi(x) &=& H^{\mathcal{S}}\big(t,x,u_*(t),p(t),1\big) - H^{\mathcal{S}}\big(t,x_*(t),u_*(t),p(t),1\big) - \langle a(t),x-x_*(t)\rangle.
\end{eqnarray*}
Die Funktion $\Phi(x)$ ist stetig differenzierbar auf $V^{\mathcal{S}}_\gamma(t)$.
Ferner gelten $\Phi(x) \leq 0$ f"ur alle $x \in V^{\mathcal{S}}_\gamma(t)$ und $\Phi(x_*(t))=0$.
Damit nimmt die Funktion $\Phi$ in dem inneren Punkt $x_*(t)$ der Menge $V^{\mathcal{S}}_\gamma(t)$ ihr globales Maximum an.
Also gilt $0=\Phi'(x_*(t))$, d.\,h.
\begin{equation} \label{BeweisHBPMP3}
-a(t)= -\varphi_x^T\big(t,x_*(t),u_*(t)\big) p(t) + f_x\big(t,x_*(t),u_*(t)\big).
\end{equation}
Die Gleichung (\ref{BeweisHBPMP3}) wurde unter der Annahme erzielt,
dass die Maximumbedingung (\ref{PMPeinfach6}) in dem Zeitpunkt $t \in [t_0,t_1]$ erf"ullt ist.
Da (\ref{PMPeinfach6}) f"ur fast alle $t \in [t_0,t_1]$ gilt,
stimmt $-a(t)$ mit der Ableitung $\dot{p}(t)$ "uberein.
Also gilt auf $V^{\mathcal{S}}_\gamma$ die Ungleichung
\begin{equation} \label{BeweisHBPMP4}
\langle \dot{p}(t),x-x_*(t)\rangle \leq \mathscr{H}^{\mathcal{S}}\big(t,x_*(t),p(t)\big)- \mathscr{H}^{\mathcal{S}}\big(t,x,p(t)\big)
\end{equation}
f"ur fast alle $t \in [t_0,t_1]$.
Für $\big(x(\cdot),u(\cdot)\big) \in \mathscr{B}^{\,\mathcal{S}}_{\rm adm}$ mit $\|x(\cdot)-x_*(\cdot)\|_\infty < \gamma$ erhalten wir:
\begin{eqnarray*}
\lefteqn{J\big(x(\cdot),u(\cdot)\big)-J\big(x_*(\cdot),u_*(\cdot)\big)
          = \int_{t_0}^{t_1} \big[f\big(t,x(t),u(t)\big)-f\big(t,x_*(t),u_*(t)\big)\big] \, dt} \\
&\geq& \int_{t_0}^{t_1} \big[\mathscr{H}^{\mathcal{S}}\big(t,x_*(t),p(t)\big)-\mathscr{H}^{\mathcal{S}}\big(t,x(t),p(t)\big)\big] \, dt 
    + \int_{t_0}^{t_1} \langle p(t), \dot{x}(t)-\dot{x}_*(t) \rangle dt \\
&\geq& \int_{t_0}^{t_1} \big[\langle \dot{p}(t),x(t)-x_*(t)\rangle + \langle p(t), \dot{x}(t)-\dot{x}_*(t) \rangle \big] \, dt \\
&=& \langle p(t_1),x(t_1)-x_*(t_1)\rangle-\langle p(t_0),x(t_0)-x_*(t_0)\rangle.
\end{eqnarray*}
Im Fall fester Anfangs- und Endbedingungen verschwinden die Differenzen $x(t_i)-x_*(t_i)$.
Sind jedoch gewissen Komponenten im Anfangs- oder Endpunkt $x_0$ bzw. $x_1$ frei,
dann liefern die Transversalit"atsbedingungen,
dass die entsprechenden Komponenten der Adjungierten $p(\cdot)$ zum Zeitpunkt $t_0$ bzw. $t_1$ verschwinden.
Daher folgt die Beziehung
$$J\big(x(\cdot),u(\cdot)\big)-J\big(x_*(\cdot),u_*(\cdot)\big) \geq 0$$
f"ur alle zul"assigen $\big(x(\cdot),u(\cdot)\big)$ mit $\|x(\cdot)-x_*(\cdot)\|_\infty < \gamma$. \hfill $\blacksquare$

\begin{beispiel} {\rm \index{Kapitalakkumulation}
Wir betrachten nach Seierstad \& Syds\ae ter \cite{Seierstad} ein Zwei-Sektoren-Modell,
das aus der Produktion von Investitons- und Konsumg"uter besteht.
Wir bezeichnen mit $x(t)$ bzw. $y(t)$ die Rate der G"uterproduktion im Investitions- bzw. Konsumsektor
und es beschreibe $u(t) \in [0,1]$ die Aufteilung der Investitionsg"uter auf beide Sektoren zur Produktion.
Es entsteht damit die Aufgabe
\begin{equation}
\left.\begin{array}{l}
\hspace*{-5mm}
\displaystyle J\big(x(\cdot),y(\cdot),u(\cdot)\big) =\int_0^T y(t) \, dt \to \sup,  \\[3mm]
\hspace*{-5mm}   
\displaystyle \dot{x}(t)=u(t)x(t), \quad x(0)=x_0>0,  \\[2mm]
\hspace*{-5mm}   
\displaystyle \dot{y}(t)=\big(1-u(t)\big)x(t), \quad y(0)=y_0>0,  \\[2mm]
\hspace*{-5mm}   
u \in [0,1], \qquad T>2.
\end{array} \right\}
\end{equation}
Wir gehen zu einem Minimierungsproblem über und wenden Theorem \ref{SatzPMP} an. \\
Die Pontrjagin-Funktion lautet $H^{\mathcal{S}}(t,x,y,u,p,q) = uxp+(1-u)xq+y$.
F"ur die Adjungierten $p(\cdot)$ und $q(\cdot)$ ergeben sich die Gleichungen
$$\dot{p}(t)=-u_*(t)p(t)-\big(1-u_*(t)\big)q(t),\quad p(T)=0,\qquad \dot{q}(t)=-1, \quad q(T)=0.$$
Es folgt unmittelbar $q(t)=T-t$.
Die Maximumbedingung ist äquivalent zu der Bedingung $\max\limits_{u \in [0,1]} u \cdot \big[p(t)-q(t)\big] \cdot x(t)$
und wir erhalten $u_*(t) = 1$ f"ur $p(t) > q(t)$ und $u_*(t) = 0$ f"ur $p(t) < q(t)$ für die optimale Steuerung. \\
Wegen $p(T)=0$ und $q(t)=T-t$ ist $\dot{p}(t)>-1$ über einem gewissen Intervall $(\tau,T)$.
Demzufolge gilt $p(t)<q(t)$ und $\dot{p}(t)=-(T-t)$ über $(\tau,T)$.
Daraus ergibt sich $\tau=T-2$. 
Andererseits ist $\dot{p}(t)<-2$ für $t \in (0,T-2)$ und $p(t)> q(t)$ über $[0,T-2)$. \\
Wir erhalten zusammenfassend für die adjungierten Funktionen 
$$p(t)= \left\{\begin{array}{ll} 2e^{T-(t+2)}, & t \in [0,T-2), \\[1mm] \frac{1}{2}(T-t)^2, & t \in [T-2,T], \end{array}\right.
  \qquad q(t)= T-t,$$
für die optimale Steuerung 
$$u_*(t) = 1 \mbox{ f"ur } t \in [0,T-2) \quad\mbox{und}\quad u_*(t) = 0 \mbox{ f"ur } t \in [T-2,T],$$
für die zugeh"origen Zustandstrajektorien
$$\big(x_*(t),y_*(t)\big) = \left\{\begin{array}{ll}
  \big(x_0 e^t,y_0\big), & t \in [0,T-2), \\[1mm]
  \big(x_0e^{T-2},y_0+(t-T+2)x_0e^{T-2}\big), & t \in [T-2,T], \end{array}\right. $$
und für den Wert des Zielfunktionals
$$J\big(x_*(\cdot),y_*(\cdot),u_*(\cdot)\big) = \int_0^T y_*(t) \, dt= y_0 T+ \frac{1}{2}(T-t-2)^2x_0e^{T-2}.$$
Die Pontrjagin-Funktion $H^{\mathcal{S}}(t,x,y,u,p,q)$ ist nicht konkav bez"uglich $(x,y,u)$ und
die hinreichenden Bedingungen nach Mangasarian im Abschnitt \ref{AbschnittMangasarianSOP} sind nicht anwendbar. \\
Beachten wir, dass f"ur alle zul"assigen Trajektorien die Bedingung $x(t) \geq x_0>0$ gilt,
dann erhalten wir f"ur die Hamilton-Funktion
$$\mathscr{H}^{\mathcal{S}}(t,x,y,p,q) = \sup_{u \in U} H^{\mathcal{S}}(t,x,y,u,p,q,1)
  =\left\{\begin{array}{ll} x p+y, & p>q, \\ x q+y, & p \leq q. \end{array}\right.$$
F"ur jedes $t \in [0,T]$ ist $\mathscr{H}^{\mathcal{S}}$ linear in $(x,y)$, also eine konkave Funktion in $(x,y)$.
Da der Kandidat $\big(x_*(\cdot),y_*(\cdot),u_*(\cdot)\big)$ alle Bedingungen des Maximumprinzips erf"ullt,
ist er eine starke lokale Minimalstelle. \hfill $\square$}
\end{beispiel}

%% file: 2-5-Zustandsaufgabe.tex
\subsection{Aufgaben mit Zustandsbeschr\"ankungen} \label{AbschnittZustandsaufgabe}
\subsubsection{Notwendige Optimalit\"atsbedingungen} \label{AbschnittZustandPMP}
\begin{theorem}[Pontrjaginsches Maximumprinzip] \label{SatzPMPZA}
\index{Pontrjaginsches Maximumprinzip!Standard@-- Standardaufgabe} 
Sei $\big(x_*(\cdot),u_*(\cdot)\big) \in \mathscr{B}^{\,\mathcal{S}}_{\rm adm} \cap \mathscr{B}^{\,\mathcal{S}}_{\rm Lip}$.
Ist $\big(x_*(\cdot),u_*(\cdot)\big)$ ein starkes lokales Minimum der Aufgabe (\ref{PMP1})--(\ref{PMP5}),
dann existieren eine Zahl $\lambda_0 \geq 0$, Vektoren $l_0 \in \R^{s_0}$ und $l_1 \in \R^{s_1}$, eine Vektorfunktion $p(\cdot):[t_0,t_1] \to \R^n$
und auf den Mengen
$$T_j=\big\{t \in [t_0,t_1] \,\big|\, g_j\big(t,x_*(t)\big)=0\big\}, \quad j=1,...,l,$$
konzentrierte nichtnegative regul"are Borelsche Ma"se $\mu_j$ endlicher Totalvariation
(wobei s"amtliche Gr"o"sen nicht gleichzeitig verschwinden) derart,
dass die Vektorfunktion $p(\cdot)$ von beschr"ankter Variation und rechtsseitig stetig ist, und
\begin{enumerate}
\item[(a)] die adjungierte Gleichung
           \index{adjungierte Gleichung!Standard@-- Standardaufgabe}
           \begin{eqnarray}
           p(t)&=&-{h_1'}^T\big(x_*(t_1)\big)l_1 + \int_t^{t_1} H^{\mathcal{S}}_x\big(s,x_*(s),u_*(s),p(s),\lambda_0\big) \, ds \nonumber \\
           \label{SatzPMPZA1} & & -\sum_{j=1}^l \int_t^{t_1} g_{j,x}\big(s,x_*(s)\big)\, d\mu_j(s),
           \end{eqnarray}
\item[(b)] die Transversalit"atsbedingungen
           \index{Transversalitätsbedingungen!Standard@-- Standardaufgabe}
           \begin{equation}\label{SatzPMPZA2}
           p(t_0)= {h_0'}^T\big(x_*(t_0)\big)l_0, \qquad p(t_1)=-{h_1'}^T\big(x_*(t_1)\big)l_1,
           \end{equation}
           in $t=t_1$ die Sprung-Transversalit"atsbedingung
           \index{Transversalitätsbedingungen!Standard@-- Standardaufgabe}
           \begin{equation} \label{SatzPMPZA4}
           p(t_1^-) - p(t_1) = - \sum_{j=1}^l \mu_j(\{t_1\}) \, g_{j,x}\big(t_1,x_*(t_1)\big)
           \end{equation}
\item[(c)] und in fast allen Punkten $t \in [t_0,t_1]$ die Maximumbedingung
           \index{Maximumbedingung!Standard@-- Standardaufgabe}
           \begin{equation}\label{SatzPMPZA3}
           H^{\mathcal{S}}\big(t,x_*(t),u_*(t),p(t),\lambda_0\big) = \max_{u \in U} H^{\mathcal{S}}\big(t,x_*(t),u,p(t),\lambda_0\big)
           \end{equation}
\end{enumerate}
erfüllt sind.
\end{theorem}

\newpage
Wie im Abschnitt \ref{AbschnittZustandSOP} existieren an jeder Stelle $t \in [t_0,t_1]$ die einseitigen Grenzwerte der Adjungierten $p(\cdot)$.
Jedoch kann die Adjungierte an der Stelle $t=t_1$ unstetig sein und es gilt (\ref{SatzPMPZA4}).

\begin{beispiel}
{\rm Wir\index{Wirtschaftswachstum und Umweltschutz} betrachten nach Seierstad \& Syds\ae ter \cite{Seierstad} eine Aufgabe,
in welcher der Ertragsmaximierung eine Umweltschutzrichtlinie gegen"uber steht.
Es sei zu jedem Zeitpunkt $t$ die Produktion proportional zum Kapital $k(t)$ mit einer positiven Proportionalit"atskonstanten $a>0$.
Ein Teil $s(t)$ vom Kapital wird wieder in die weitere Produktion investiert und der verbleibende Teil $(1-s(t))$ dient zur
Konsumption.
Ferner besteht die M"oglichkeit die Produktionskapazit"aten "uber die Auslastungsrate $u(t)$ zu steuern. \\
Bei der Produktion f"allt Abfall an.
Es beschreibt $z(t)$ das gesamte M"ullvolumen zur Zeit $t$ und $\dot{z}(t)$ den bei der Produktion zum Zeitpunkt $t$ entstehenden Abfall.
Wir nehmen an, dass $\dot{z}(t)$ proportional zur Produktion und somit proportional zum Kapital ist,
d.\,h. $\dot{z}(t)=ck(t)$ mit einer Konstanten $c>0$. \\
In diesem Modell hat die Verschmutzung $z(t)$ keinen direkten Einfluss auf die Ertragsmaximierung.
Die Planungsperiode $[0,T]$ und die Startwerte $k(0), z(0)$ f"ur das Kapital und das M"ullvolumen sind bekannt.
An den Endwert $k(T)$ wird keine Bedingung gestellt,
aber das Abfallvolumen $z(t)$ soll "uber $[0,T]$ die vorgegebene Richtlinie $Z$ nicht "ubersteigen. \\[2mm]
Wir betrachten demnach die folgende Aufgabe:
\begin{eqnarray}
&& \label{BspZB1} J\big(k(\cdot),z(\cdot),u(\cdot),s(\cdot)\big) = \int_0^T \big(1-s(t)\big)u(t)k(t) \, dt \to \sup, \\
&& \label{BspZB2} \dot{k}(t) = s(t)u(t)k(t), \quad k(0)=k_0, \quad k(T) \mbox{ frei}, \\
&& \label{BspZB3} \dot{z}(t) = u(t)k(t)-az(t), \quad z(0)=z_0, \quad z(T) \mbox{ frei}, \\
&& \label{BspZB4} (u,s) \in U = [0,1] \times [0,1], \\
&& \label{BspZB5} g(t,k,z)=z(t)-Z \leq 0.
\end{eqnarray}
Weiterhin sind $k_0, z_0, Z, T$ und $a$ positive Konstanten, f"ur die die Relationen
\begin{equation} \label{BspZB6}
z_0 < Z, \qquad T > 1, \qquad \left(z_0-\frac{k_0}{a}\right)e^{-aT} + \frac{k_0}{a} < Z
\end{equation}
gelten.
Die letzte Annahme wird dadurch erkl"art, dass zul"assige Handlungsm"oglichkeiten existieren.
Werden n"amlich keine Investitionen get"atigt, d.\,h. $s(t) \equiv 0$ auf $[0,T]$,
so folgt $k(t) \equiv k_0$ und wir erhalten mit der Auslastungsrate $u(t) \equiv 1$:
$$z(T)=\bigg( z_0 - \frac{k_0}{a} \bigg) e^{-aT} + \frac{k_0}{a}.$$
Das Modell legt es nahe folgende drei F"alle zu betrachten:
\begin{enumerate}
\item[(A)] Bei der Produktion wird in der Planungsperiode niemals die Grenze $Z$
           f"ur die Abfallmenge erreicht und es kann "uber dem gesamten Planungshorizont
           das maximal m"ogliche Produktionsvolumen ausgenutzt werden, d.\,h.
           $$z(t) < Z \quad\mbox{ f"ur alle } t \in [0,T].$$
\item[(B)] Die Richtlinie $Z$ wird genau im Endzeitpunkt der Planungsperiode erreicht, d.\,h.
           $$z(t) < Z \quad \mbox{ f"ur alle } t \in [0,T), \qquad z(T)=Z.$$
\item[(C)] Innerhalb der Planungsperiode wird bereits die Grenze $Z$ erreicht.
           Anschlie"send l"auft die Produktion auf der h"ochstm"oglichen Stufe, so dass
           die Abfallmenge $z(t)$ die Richtlinie $Z$ nicht "uberschreitet, d.\,h.
           $$z(t) < Z \quad\mbox{ f"ur alle } t \in [0,t'), \; 0<t'<T, \qquad z(t)=Z \quad\mbox{ f"ur alle } t \in [t',T].$$
\end{enumerate}
Im Weiteren geben wir f"ur jeden Fall nach \cite{Seierstad} optimale Steuerungsprozesse und die wesentlichen Merkmale der Adjungierten an.
Auf die Details der Analyse, die wir nicht ausf"uhrlich belegen, verweisen wir auf \cite{Seierstad}.
\begin{enumerate}
\item[(A)] Ist die Zustandsbeschr"ankung nicht aktiv, so erhalten wir die Steuerungen
           $$s_*(t)=\left\{ \begin{array}{ll} 1, & t \in [0,T-1), \\ 0, & t \in [T-1,T], \end{array} \right. \qquad
             u_*(t)\equiv 1 \mbox{ auf } [0,T].$$
           F"ur die Dynamiken der Zustandstrajektorien $k_*(\cdot)$, $z_*(\cdot)$ ergeben sich daraus
           $$\dot{k}_*(t)=\left\{ \begin{array}{ll} k_*(t), & t \in [0,T-1), \\ 0, & t \in [T-1,T], \end{array} \right. \qquad
             \dot{z}_*(t)=\left\{ \begin{array}{ll} k_0e^t-az_*(t), & t \in [0,T-1), \\
                                                    k_0e^{T-1}-az_*(t), & t \in [T-1,T]. \end{array} \right.$$
           Die Adjungierten $p_1(\cdot)$, $p_2(\cdot)$, die den Zust"anden $k_*(\cdot)$ bzw. $z_*(\cdot)$ zugeordnet sind,
           sind absolutstetig und erf"ullen die Transversalit"atsbedingungen $p_1(T)=0$ und $p_2(T)=0$.
\item[(B)] Ist die Zustandsbeschr"ankung in $t=T$ aktiv, dann gilt f"ur die Steuerungen
           $$s_*(t)=\left\{ \begin{array}{ll} 1, & t \in [0,t_*), \\ 0, & t \in [t_*,T], \end{array} \right. \qquad
             u_*(t)\equiv 1 \mbox{ auf } [0,T],$$
           wobei $t_*$ der Bedingung
           $$e^{t_*}-\frac{e^{-aT}}{1+a}e^{(1+a)t_*} = \frac{a}{k_0} \bigg[Z-\bigg(z_0-\frac{k_0}{1+a}\bigg)e^{-aT}\bigg]$$
           gen"ugt. Die Zustandstrajektorien verhalten sich entsprechend dem Fall (A),
           wobei der Umschaltpunkt $t=T-1$ durch $t=t_*$ zu ersetzen ist.
           F"ur die Adjungierte $p_2(\cdot)$ ergibt sich in $t=T$ die Transversalit"atsbedingung
           $$p_2(T)=-\mu(\{T\})=a\frac{t_*-(T-1)}{1-e^{-a(T-t_*)}} \in [-1,0].$$
           Daher muss $t_* \leq T-1$ gelten. Im Vergleich zum Fall (A) werden die Investitionen gegebenenfalls fr"uher beendet.
           Ist dabei $t_*=T-1$, dann stimmt die L"osung mit derjenigen L"osung im Fall (A) "uberein.
\item[(C)] Es sei $t' = \sup\{ t \in [0,T] \,|\, z_*(t)<Z\}$.
           Ist die Zustandsbeschr"ankung auf dem Intervall $[t',T]$ aktiv,
           dann gibt es einen Zeitpunkt $t''$ mit $0<t''<t'<T$ und
           $$s_*(t)=\left\{ \begin{array}{ll} 1, & t \in [0,t''), \\ 0, & t \in [t'',T], \end{array} \right. \qquad
            u_*(t)=\left\{ \begin{array}{ll} 1, & t \in [0,t'), \\ \displaystyle\frac{aZ}{k_*(t'')}, & t \in [t',T]. \end{array} \right.$$
           Dabei erf"ullen die Zeitpunkte  $t',t''$ die gekoppelten Bedingungen
           \begin{eqnarray*}
           1+\frac{1}{a} &=& t'-t''+\frac{1}{a} e^{-a(t'-t'')}, \\
           Z &=& \bigg[z_0-\frac{k_0}{1+a}-\frac{k_0}{a(1+a)}e^{(1+a)t''}\bigg]e^{-at'}+\frac{k_0}{a}e^{t''}.
           \end{eqnarray*}
           F"ur die Dynamiken der Zustandstrajektorien $k_*(\cdot)$, $z_*(\cdot)$ ergeben sich daraus
           $$\dot{k}_*(t)=\left\{ \begin{array}{ll} k_*(t), & t \in [0,t''), \\ 0, & t \in [t'',T], \end{array} \right. \qquad
             \dot{z}_*(t)=\left\{ \begin{array}{ll} k_*(t)-az_*(t), & t \in [0,t'), \\ 0, & t \in [t',T]. \end{array} \right.$$
           Die Adjungierte $p_2(t)$ gen"ugt auf $(t',T)$ der Differentialgleichung
           $$\dot{p}_2(t)=ap_2(t)+\lambda(t), \qquad \lambda(t)=a$$
           und es ergibt sich in $t=T$ die Transversalit"atsbedingung $p_2(T)=-1$.  \hfill $\square$
\end{enumerate}
}
\end{beispiel}

%% file: 2-51-Beweis.tex
\subsubsection{Der Nachweis der notwendigen Optimalit\"atsbedingungen} \label{AbschnittBeweisPMPZA}
Wir betrachten f"ur $\big(x(\cdot),u(\cdot)\big) \in C([t_0,t_1],\R^n) \times L_\infty([t_0,t_1],\R^m)$ die Abbildungen
\begin{eqnarray*}
J\big(x(\cdot),u(\cdot)\big) &=& \int_{t_0}^{t_1} f\big(t,x(t),u(t)\big) \, dt, \\
F\big(x(\cdot),u(\cdot)\big)(t) &=& x(t) -x(t_0) -\int_{t_0}^t \varphi\big(s,x(s),u(s)\big) \, ds, \quad t \in [t_0,t_1],\\
H_i\big(x(\cdot)\big) &=& h_i\big(x(t_i)\big), \quad i=0,1, \\
G_j\big(x(\cdot)\big) &=& \max_{t \in [t_0,t_1]} g_j\big(t,x(t)\big), \quad j=1,...,l.
\end{eqnarray*}
Da $x(\cdot)$ zu $C([t_0,t_1],\R^n)$ geh"ort, gilt f"ur diese Abbildungen
\begin{eqnarray*}
J &:& C([t_0,t_1],\R^n) \times L_\infty([t_0,t_1],\R^m) \to \R, \\
F &:& C([t_0,t_1],\R^n) \times L_\infty([t_0,t_1],\R^m) \to C_0([t_0,t_1],\R^n), \\
H_i &:& C([t_0,t_1],\R^n) \to \R^{s_i}, \quad i=0,1, \\
G_j &:& C([t_0,t_1],\R^n) \to \R, \quad j=1,...,l.
\end{eqnarray*}
Wir setzen $\mathscr{F}=(F,H_0,H_1)$ und pr"ufen f"ur die Extremalaufgabe
\begin{equation} \label{ExtremalaufgabePMPZA}
J\big(x(\cdot),u(\cdot)\big) \to \inf, \quad \mathscr{F}\big(x(\cdot),u(\cdot)\big)=0, \quad G_j\big(x(\cdot)\big) \leq 0,
\quad u(\cdot) \in L_\infty([t_0,t_1],U)
\end{equation}
im Punkt $\big(x_*(\cdot),u_*(\cdot)\big) \in \mathscr{B}^{\,\mathcal{S}}_{\rm Lip}$,
wobei wir $x_*(\cdot)$ als Element des Raumes $C([t_0,t_1],\R^n)$ auffassen,
die Voraussetzungen von Theorem \ref{SatzExtremalprinzipStark}:

\begin{enumerate}
\item[(A$_2$)] Mit Verweis auf Abschnitt \ref{AbschnittBeweisPMP} sind nur noch die Abbildungen $G_j$ zu diskutieren.
               In den Beispielen \ref{SubdifferentialMaximum1} und \ref{SubdifferentialMaximum2} wird gezeigt,
               dass die Funktionen $G_j$ Hintereinanderausf"uhrungen einer stetigen, konvexen, eigentlichen Funktion und
               einer Fr\'echet-differenzierbaren Abbildung sind.
               Daher sind nach Lemma \ref{LemmaRichtungsableitung} die Funktionen $G_j$ in $x_*(\cdot)$ lokalkonvex und bez"uglich jeder Richtung
               gleichm"a"sig differenzierbar.              
\end{enumerate}

Zur Extremalaufgabe (\ref{ExtremalaufgabePMPZA}) definieren wir auf
$$C([t_0,t_1],\R^n) \times L_\infty([t_0,t_1],\R^m) \times \R \times C_0^*([t_0,t_1],\R^n) \times \R^{s_0} \times \R^{s_1} \times \R^l$$
die Lagrange-Funktion $\mathscr{L}=\mathscr{L}\big(x(\cdot),u(\cdot),\lambda_0,y^*,l_0,l_1,\lambda\big)$,
$$\mathscr{L}= \lambda_0 J\big(x(\cdot),u(\cdot)\big)+ \big\langle y^*, F\big(x(\cdot),u(\cdot)\big) \big\rangle
                         +l_0^T H_0\big(x(\cdot)\big)+l_1^T H_1\big(x(\cdot)\big) + \sum_{j=1}^l \lambda_j G_j\big(x(\cdot)\big).$$
Ist $\big(x_*(\cdot),u_*(\cdot)\big)$ eine starke lokale Minimalstelle der Extremalaufgabe (\ref{ExtremalaufgabePMPZA}),
dann existieren nach Theorem \ref{SatzExtremalprinzipStark}
nicht gleichzeitig verschwindende Lagrangesche Multiplikatoren $\lambda_0 \geq 0$, $y^* \in C_0^*([t_0,t_1],\R^n)$, $l_i \in \R^{s_i}$
und $\lambda_1 \geq 0,...,\lambda_l \geq 0$ derart,
dass gelten:
\begin{enumerate}
\item[(a)] Die Lagrange-Funktion besitzt bez"uglich $x(\cdot)$ in $x_*(\cdot)$ einen station"aren Punkt, d.\,h.
           \begin{equation}\label{SatzPMPZALMR1}
           0 \in \partial_x \mathscr{L}\big(x_*(\cdot),u_*(\cdot),\lambda_0,y^*,l_0,l_1,\lambda\big);
           \end{equation}         
\item[(b)] Die Lagrange-Funktion erf"ullt bez"uglich $u(\cdot)$ in $u_*(\cdot)$ die Minimumbedingung
           \begin{equation}\label{SatzPMPZALMR2}
           \hspace*{-3mm} \mathscr{L}\big(x_*(\cdot),u_*(\cdot),\lambda_0,y^*,l_0,l_1,\lambda\big)
           = \min_{u(\cdot) \in L_\infty([t_0,t_1],U)} \mathscr{L}\big(x_*(\cdot),u(\cdot),\lambda_0,y^*,l_0,l_1,\lambda\big);
           \end{equation}
\item[(c)] Die komplement"aren Schlupfbedingungen gelten, d.\,h.
           \begin{equation}\label{SatzPMPZALMR3}
           0 = \lambda_j G_j\big(x(\cdot)\big), \qquad i=1,...,l.
           \end{equation}
\end{enumerate}
Aufgrund (\ref{SatzPMPZALMR1}) ist folgende Variationsgleichung f"ur alle $x(\cdot) \in C([t_0,t_1],\R^n)$ erf"ullt: 
\begin{eqnarray}
0 &=& \lambda_0 \int_{t_0}^{t_1} \big\langle f_x\big(t,x_*(t),u_*(t)\big),x(t) \big\rangle\, dt
      + \big\langle l_0, h_0'\big(x_*(t_0)\big) x(t_0) \big\rangle + \big\langle l_1, h_1'\big(x_*(t_1)\big) x(t_1) \big\rangle
      \nonumber \\
  & & + \int_{t_0}^{t_1} \bigg[ x(t)-x(0) - \int_0^t \varphi_x\big(s,x_*(s),u_*(s)\big) x(s) \,ds \bigg]^T d\mu(t) \nonumber \\
  & & \label{BeweisschlussPMPZA1}
      + \sum_{j=1}^l \lambda_j \int_{t_0}^{t_1}\big\langle g_{j,x}\big(t,x_*(t)\big),x(t) \big\rangle \,d\tilde{\mu}_j(t).
\end{eqnarray}
In (\ref{BeweisschlussPMPZA1}) verschwinden die Lagrangeschen Multiplikatoren $\lambda_j \geq 0$,
falls die entsprechende Zustandsbeschr"ankung nichtaktiv ist.
Andernfalls besitzen die nichtnegativen regul"aren Borelschen Ma"se $\tilde{\mu}_j$ die Totalvariation $\|\tilde{\mu}_j\|=1$ und sind 
jeweils auf der Menge
$$T_j = \big\{ t \in [t_0,t_1] \big| g_j\big(t,x_*(t)\big) = 0 \big\}$$
konzentriert.
Wir schreiben $\mu_j=\lambda_j \tilde{\mu}_j$.
Dann k"onnen unter den Ma"sen $\mu_j$ nur diejenigen von Null verschieden sein,
f"ur die die Menge $T_j \not= \emptyset$ ist.
Daher kann man o.\,B.\,d.\,A. annehmen, alle Ma"se $\mu_j$ seien auf den Mengen $T_j$ konzentriert. \\[2mm]
In der Gleichung (\ref{BeweisschlussPMPZA1}) "andern wir die Integrationsreihenfolge im zweiten Summanden und bringen sie in die Form
\begin{eqnarray}
0 &=& \int_{t_0}^{t_1} \big\langle \lambda_0 f_x\big(t,x_*(t),u_*(t)\big) - \varphi_x^T\big(t,x_*(t),u_*(t)\big) \int_{t}^{t_1} d\mu(s) ,
                       x(t) \big\rangle \, dt \nonumber \\
  & & + \int_{t_0}^{t_1} [x(t)]^T \, d\mu(t) + \Big\langle {h_0'}^T\big(x_*(t_0)\big)l_0 - \int_{t_0}^{t_1} d\mu(t) , x(t_0) \Big\rangle
      + \langle {h_1'}^T\big(x_*(t_1)\big)l_1 , x(t_1) \rangle \nonumber \\
  & & \label{BeweisschlussPMPZA2}
      + \sum_{j=1}^l \int_{t_0}^{t_1} \big\langle g_{j,x}\big(t,x_*(t)\big),x(t) \big\rangle \,d\mu_j(t).
\end{eqnarray}
Wir setzen $p(t)=\displaystyle \int_t^{t_1} \, d\mu(s)$,
so ist $p(\cdot)$ eine Funktion von beschränkter Variation und gemäß den Eigenschaften einer Verteilungsfunktion rechtsseitig stetig. \\[1mm]
Die rechte Seite in (\ref{BeweisschlussPMPZA2}) definiert ein stetiges lineares Funktional im Raum $C([t_0,t_1],\R^n)$.
Wenden wir den Darstellungssatz von Riesz an,
so folgen aus der eindeutigen Darstellung eines stetigen linearen Funktionals im Raum $C([t_0,t_1],\R^n)$ die Beziehungen
\begin{eqnarray*}
p(t) &=& -{h_1'}^T\big(x_*(t_1)\big)l_1
         + \int_t^{t_1} H^{\mathcal{S}}_x\big(s,x_*(s),u_*(s),p(s),\lambda_0\big) \, ds \\
     & & - \sum_{j=1}^l \int_t^{t_1} g_{j,x}\big(s,x_*(s)\big) \,d\mu_j(s), \\
p(t_0) &=& {h_0'}^T\big(x_*(t_0)\big)l_0.
\end{eqnarray*}
Damit sind (\ref{SatzPMPZA1}) und (\ref{SatzPMPZA2}) gezeigt. \\
Genauso wie im Abschnitt \ref{AbschnittBeweisPMP} ist die Beziehung (\ref{SatzPMPZALMR2}) "aquivalent zu
\begin{eqnarray*} 
\lefteqn{\int_{t_0}^{t_1} \lambda_0 f\big(t,x_*(t),u_*(t)\big) \, dt
    - \int_{t_0}^{t_1}  \left[ \int_{t_0}^{t} \varphi\big(s,x_*(s),u_*(s)\big) \, ds \right]^T d\mu(t)} \\
&&\hspace*{-5mm} = \min_{u(\cdot) \in L_\infty ([t_0,t_1],U)}
    \left\{\int_{t_0}^{t_1} \lambda_0 f\big(t,x_*(t),u(t)\big) \, dt
            - \int_{t_0}^{t_1} \bigg[ \int_{t_0}^{t} \varphi\big(s,x_*(s),u(s)\big)\, ds \bigg]^T d\mu(t) \right\}
\end{eqnarray*}
und es ergibt sich durch Vertauschen der Integrationsreihenfolge die Ungleichung
\begin{eqnarray} 
\lefteqn{\int_{t_0}^{t_1} \big[ \lambda_0 f\big(t,x_*(t),u_*(t)\big) - \big\langle p(t) , \varphi\big(t,x_*(t),u_*(t)\big) \big\rangle \big] \, dt} \nonumber \\
\label{BeweisschlussPMPZA3}
&\leq & \int_{t_0}^{t_1} \big[ \lambda_0 f\big(t,x_*(t),u(t)\big) - \big\langle p(t) , \varphi\big(t,x_*(t),u(t)\big) \big\rangle \big] \, dt
\end{eqnarray}
f"ur alle $u(\cdot) \in L_\infty([t_0,t_1],U)$.
Die Herleitung der Maximumbedingung (\ref{SatzPMPZA3}) aus (\ref{BeweisschlussPMPZA3})
basiert auf den Lebesgueschen Punkten\index{Lebesguescher Punkt}:
Ist $f(t)$ über $[a,b]$ integrierbar, so ist fast jeder Punkt $t \in (a,b)$, insbesondere jeder Stetigkeitspunkt,
ein Lebesguescher Punkt (vgl. \cite{Natanson}) und es gilt in diesen Stellen
$$\lim_{h \to 0^+} \frac{1}{h} \int_0^h |f(t_0 \pm t)-f(t_0)| \, dt =0.$$
In (\ref{BeweisschlussPMPZA3}) sind $t \to f\big(t,x_*(t),u_*(t)\big)$ und $t \to \varphi\big(t,x_*(t),u_*(t)\big)$ über
$[t_0,t_1]$ integrierbar.
Demnach besitzen sie fast überall in $(t_0,t_1)$ Lebesguesche Punkte.
Mit $T_f$ bzw. $T_\varphi$ bezeichnen wir die Stellen, in denen für diese Abbildungen keine Lebesgueschen Punkte vorliegen.\\
Ferner ist in (\ref{BeweisschlussPMPZA3}) die Adjungierte $p(\cdot)$ eine Funktion beschränkter Variation.
Damit besitzt $p(\cdot)$ als Differenz zweier monotoner Funktionen (\!\!\cite{Natanson}, S.\,229) nur abzählbar viele
Unstetigkeitsstellen.
Mit $T_p$ bezeichnen wir die Menge der Unstetigkeitsstellen.
Ferner sind für jedes $u \in U$ die Abbildungen $t \to f\big(t,x_*(t),u\big)$ und $t \to \varphi\big(t,x_*(t),u\big)$ über
$[t_0,t_1]$ integrierbar und stetig.
Also sind für jedes $u \in U$ alle $t \in (t_0,t_1)$ Lebesguesche Punkte. \\[1mm]
Es seien $u \in U$, $\lambda>0$ und $\tau \in (t_0,t_1)\setminus (T_f \cup T_\varphi \cup T_p)$.
Wir setzen f"ur $0<\lambda<t_1-\tau$
$$u_\lambda(t)=\left\{\begin{array}{ll} u,& t \in [\tau,\tau+\lambda], \\[1mm]
                                u_*(t), & t \not\in [\tau,\tau+\lambda]. \end{array}\right.$$
Dann gilt $u_\lambda(\cdot) \in L_\infty(\R_+,U)$ und ist zulässig.
In (\ref{BeweisschlussPMPZA3}) ergibt sich nun
\begin{eqnarray*}
\lefteqn{\frac{1}{\lambda}
               \int_{t_0}^{t_1} \big[ H^{\mathcal{S}}\big(t,x_*(t),u_\lambda(t),p(t),\lambda_0\big)
                                -H^{\mathcal{S}}\big(t,x_*(t),u_*(t),p(t),\lambda_0\big)\big]\, dt}\\
&=& \frac{1}{\lambda} \int_\tau^{\tau+\lambda} \big[ H^{\mathcal{S}}\big(t,x_*(t),u,p(t),\lambda_0\big)
                                -H^{\mathcal{S}}\big(t,x_*(t),u_*(t),p(t),\lambda_0\big)\big]\, dt \leq 0.
\end{eqnarray*}
Da $\tau \in (t_0,t_1)\setminus (T_f \cup T_\varphi \cup T_p)$ ein Lebesguescher Punkt aller eingehender Abbildungen ist,
folgt $H^{\mathcal{S}}\big(\tau,x_*(\tau),u,p(\tau),\lambda_0\big)-H^{\mathcal{S}}\big(\tau,x_*(\tau),u_*(\tau),p(\tau),\lambda_0\big) \leq 0$.
Zusammenfassend gilt diese Ungleichung für alle $\tau \in (t_0,t_1)\setminus (T_f \cup T_\varphi \cup T_p)$ und alle $u \in U$.
Somit ist für fast alle $t \in [t_0,t_1]$ die Maximumbedingung (\ref{SatzPMPZA3}) erfüllt. \hfill $\blacksquare$

%% file: 2-52-HinreichendZB.tex
\subsubsection{Hinreichende Bedingungen nach Arrow} \label{AbschnittArrowZB}
Die\index{hinreichende Bedingungen, Arrow!Standard@-- Standardaufgabe}
Herleitung der hinreichenden Bedingungen basiert wieder auf Seierstad \& Syds\ae ter \cite{Seierstad}
mit den erg"anzenden Ausf"uhrungen in Aseev \& Kryazhimskii \cite{AseKry}. \\[2mm]
Wir betrachten das Steuerungsproblem
\begin{eqnarray}
&& \label{HBZA1} J\big(x(\cdot),u(\cdot)\big) = \int_{t_0}^{t_1} f\big(t,x(t),u(t)\big) \, dt \to \inf, \\
&& \label{HBZA2} \dot{x}(t) = \varphi\big(t,x(t),u(t)\big), \\
&& \label{HBZA3} x(t_0)=x_0, \qquad x(t_1)=x_1, \\
&& \label{HBZA4} u(t) \in U \subseteq \R^m, \quad U\not= \emptyset, \\
&& \label{HBZA5} g_j\big(t,x(t)\big) \leq 0 \quad \mbox{f"ur alle } t \in [t_0,t_1], \quad j=1,...,l.
\end{eqnarray}

In der Aufgabenstellung schlie"sen wir wieder den Fall nicht aus,
dass durch die Randbedingungen (\ref{HBZA3}) gewisse Komponenten der Punkte $x_0$ und $x_1$ nicht fest vorgegeben,
sondern ohne Einschr"ankung sind. \\[2mm]
Wir definieren die Menge $V^{\mathcal{S}}_\gamma(t)=\{ x \in \R^n \,|\, \|x-x_*(t)\| < \gamma\}$.
Weiterhin bezeichnet
\begin{equation} \label{ZAHamilton}
\mathscr{H}^{\mathcal{S}}(t,x,p) = \sup_{u \in U} H^{\mathcal{S}}(t,x,u,p,1)
\end{equation}
die Hamilton-Funktion $\mathscr{H}^{\mathcal{S}}$ im normalen Fall.

\begin{theorem} \label{SatzHBZA}
In der Aufgabe (\ref{HBZA1})--(\ref{HBZA5}) sei $\big(x_*(\cdot),u_*(\cdot)\big) \in \mathscr{B}^{\,\mathcal{S}}_{\rm Lip} \cap \mathscr{B}^{\,\mathcal{S}}_{\rm adm}$.
Au"serdem sei die Vektorfunktion $p(\cdot):[t_0,t_1] \to \R^n$ st"uckweise stetig,
besitze h"ochstens endlich viele Sprungstellen (d.\,h. die einseitigen Grenzwerte existieren) $s_k \in (t_0,t_1)$
und sei zwischen diesen Spr"ungen stetig differenzierbar.
Ferner gelte:
\begin{enumerate}
\item[(a)] Das Tripel $\big(x_*(\cdot),u_*(\cdot),p(\cdot)\big)$
           erf"ullt (\ref{SatzPMPZA1})--(\ref{SatzPMPZA3}) in Theorem \ref{SatzPMPZA} mit $\lambda_0=1$.        
\item[(b)] F"ur jedes $t \in [t_0,t_1]$ ist die Funktion $\mathscr{H}^{\mathcal{S}}\big(t,x,p(t)\big)$ konkav 
           und es sind die Funktionen $g_j(t,x)$, $j=1,...,l$, konvex bez"uglich $x$ auf $V^{\mathcal{S}}_\gamma(t)$.
\end{enumerate}
Dann ist $\big(x_*(\cdot),u_*(\cdot)\big)$ ein starkes lokales Minimum der Aufgabe (\ref{HBZA1})--(\ref{HBZA5}).
\end{theorem}

\begin{bemerkung}\label{BemHBSOPZB} {\rm
Der Teil (a) in Theorem \ref{SatzHBZA} bedarf einer detaillierteren Diskussion: \\[1mm]
Da wir von einer st"uckweise, nach Theorem \ref{SatzPMPZA} rechtsseitig stetigen und zwischen den Sprungstellen stetig differenzierbaren Adjungierten $p(\cdot)$ ausgehen,
k"onnen wir die adjungierte Gleichung (\ref{SatzPMPZA1}) in Integraldarstellung in die Form einer st"uckweise definierten Differentialgleichung mit
Sprungbedingungen "uberf"uhren. \\[1mm]
Es bezeichnen $t_0<s_1<...<s_d < t_1$ die Unstetigkeitsstellen der Adjungierten $p(\cdot)$ im Intervall $(t_0,t_1)$.
Dann gelten die Sprungbedingungen
$$p(s_k)-p(s_k^-)= \sum_{j=1}^l \beta_j^k g_{j,x}\big(s_k,x_*(s_k)\big),\qquad \beta_j^k = \mu_j(\{s_k\}) \geq 0, \quad k=1,...,d.$$
Ferner gibt es eine st"uckweise stetige Vektorfunktion $\lambda(\cdot):[t_0,t_1] \to \R^l$ derart,
dass
$$\dot{p}(t)=- H^{\mathcal{S}}_x\big(t,x_*(t),u_*(t),p(t),1\big) + \sum_{j=1}^l \lambda_j(t) g_{j,x}\big(t,x_*(t)\big)$$
st"uckweise auf $(s_{k-1},s_k)$, $k=1,...,d+1$, gilt. Dabei haben wir $s_0=t_0$, $s_{d+1}=t_1$ gesetzt. \\[1mm]
Abschließend halten wir fest,
dass wegen der Positivit"at der Ma"se $\mu_j$ und der Konzentration dieser Ma"se auf den Mengen
$T_j=\big\{t \in [t_0,t_1] \,\big|\, g_j\big(t,x_*(t)\big)=0\big\}$ die Komplementaritätsbedingungen
$\lambda_j(t) \geq 0$ und $\lambda_j(t)g_j\big(t,x_*(t)\big)=0$ auf $[t_0,t_1]$,
in den Sprungstellen $\beta_j^k \geq 0$ und $\beta_j^k g_j\big(s_k,x_*(s_k)\big)=0$ f"ur $j=1,...,l$ und $k=1,...,d$,
sowie $\beta_j =\mu_j(\{t_1\}) \geq 0$ und $\beta_j g_j\big(t_1,x_*(t_1)\big)=0$ im Endpunkt $t_1$ gelten. \hfill $\square$}
\end{bemerkung}

{\bf Beweis} Wie im Beweis von Theorem \ref{SatzHBPMP} in Abschnitt \ref{AbschnittHBPMP} ergibt sich 
$$\langle a(t),x-x_*(t)\rangle \geq  \mathscr{H}^{\mathcal{S}}\big(t,x,p(t)\big) - \mathscr{H}^{\mathcal{S}}\big(t,x_*(t),p(t)\big)
  \quad \mbox{ f"ur alle } x \in V_\gamma(t),$$
wobei $a(t)=H^{\mathcal{S}}_x\big(t,x_*(t),u_*(t),p(t),1\big)$ f"ur fast alle $t \in [t_0,t_1]$ erfüllt ist,
vgl. (\ref{BeweisHBPMP2}) und (\ref{BeweisHBPMP3}).
Damit und mit Bemerkung \ref{BemHBSOPZB} gilt f"ur fast alle $t \in [t_0,t_1]$ die Ungleichung
$$\mathscr{H}^{\mathcal{S}}\big(t,x_*(t),p(t)\big) - \mathscr{H}^{\mathcal{S}}\big(t,x,p(t)\big) 
  \geq \Big\langle \dot{p}(t) - \sum_{j=1}^l \lambda_j(t) g_{j,x}\big(t,x_*(t)\big) , x-x_*(t) \Big\rangle.$$
Ferner gelten die Komplementarit"atsbedingungen $\lambda_j(t) \geq 0$ und $\lambda_j(t)g_j\big(t,x_*(t)\big)=0$. \\[2mm]
Wir betrachten zun"achst den Fall, dass $p(\cdot)$ stetig ist.
Die Konvexit"at der Abbildungen $g_j(t,x)$ bez"uglich $x$ auf $V_\gamma(t)$ liefert die Ungleichungen
$$g_j(t,x) \geq g_j\big(t,x_*(t)\big) + \big\langle g_{j,x}\big(t,x_*(t)\big) , x-x_*(t) \big\rangle.$$
Gen"ugt $x(\cdot)$ den Zustandsbeschr"ankungen $g_j(t,x)\leq 0$,
dann ergeben sich mit den Komplementarit"atsbedingungen $\lambda_j(t)g_j\big(t,x_*(t)\big)=0$ und $\lambda_j(t) \geq 0$ die Ungleichungen
$$-\lambda_j(t) \big\langle g_{j,x}\big(t,x_*(t)\big) , x(t)-x_*(t) \big\rangle \geq - \lambda_j(t) g_j\big(t,x(t)\big)\geq 0, \quad j=1,...,l.$$
Die gleiche Argumentation führt an der Stelle $t=t_1$ zu der Ungleichung
\begin{equation} \label{BeweisHBZA}
-\sum_{j=1}^l \big\langle \beta_j  g_{j,x}\big(t_1,x_*(t_1)\big) ,x(t_1)-x_*(t_1) \big\rangle \geq 0.
\end{equation}

Damit erhalten wir f"ur $\big(x(\cdot),u(\cdot)\big) \in \mathscr{B}^{\,\mathcal{S}}_{\rm adm}$ mit $\|x(\cdot)-x_*(\cdot)\|_\infty < \gamma$:
\begin{eqnarray*}
\lefteqn{J\big(x(\cdot),u(\cdot)\big)-J\big(x_*(\cdot),u_*(\cdot)\big)
          = \int_{t_0}^{t_1} \big[f\big(t,x(t),u(t)\big)-f\big(t,x_*(t),u_*(t)\big)\big] \, dt} \\
&\geq& \int_{t_0}^{t_1} \big[\mathscr{H}^{\mathcal{S}}\big(t,x_*(t),p(t)\big)-\mathscr{H}^{\mathcal{S}}\big(t,x(t),p(t)\big)\big] \, dt 
    + \int_{t_0}^{t_1} \langle p(t), \dot{x}(t)-\dot{x}_*(t) \rangle dt \\
&\geq& \int_{t_0}^{t_1} \big[ \langle \dot{p}(t),x(t)-x_*(t)\rangle + \langle p(t), \dot{x}(t)-\dot{x}_*(t) \rangle \big] \, dt \\
&&    \hspace*{10mm} - \int_{t_0}^{t_1} \sum_{j=1}^l \lambda_j(t) \big\langle g_{j,x}\big(t,x_*(t)\big) , x(t)-x_*(t) \big\rangle \, dt \\
&\geq& \int_{t_0}^{t_1} \big[\langle \dot{p}(t),x(t)-x_*(t)\rangle + \langle p(t),\dot{x}(t)-\dot{x}_*(t) \rangle\big] \, dt \\
&=& \langle p(t_1),x(t_1)-x_*(t_1)\rangle -\langle p(t_0),x(t_0)-x_*(t_0)\rangle \\
&\geq& -\langle l_1,x(t_1)-x_*(t_1)\rangle -\langle l_0,x(t_0)-x_*(t_0)\rangle.
\end{eqnarray*}
Die letzte Zeile folgt mit Hilfe (\ref{BeweisHBZA}) und der Transversalitätsbedingungen (\ref{SatzPMPZA2}), (\ref{SatzPMPZA4}). \\[2mm]
Es besitze nun $p(\cdot)$ im Intervall $(t_0,t_1)$ die Unstetigkeitsstellen $t_0<s_1<...<s_d \leq t_1$.
Dann gelten nach Bemerkung \ref{BemHBSOPZB} die Sprungbedingungen
$$p(s_k)-p(s_k^-)= \sum_{j=1}^l \beta_j^k g_{j,x}\big(s_k,x_*(s_k)\big),\qquad \beta_j^k \geq 0, \qquad k=1,...,d.$$
Erneut mit Hilfe der Konvexit"at der Abbildungen $g_j(t,x)$ und der Komplementaritätsbedingungen erhalten wir die Ungleichungen
$$-\langle p(s_k)-p(s_k^-),x-x_*(s_k)\rangle \geq 0, \qquad x\in V_\gamma(s_k), \qquad k=1,...,d.$$
F"ur $\big(x(\cdot),u(\cdot)\big) \in \mathscr{B}^{\,\mathcal{S}}_{\rm adm}$ mit $\|x(\cdot)-x_*(\cdot)\|_\infty < \gamma$ ergibt sich im Fall $s_d < t_1$
\begin{eqnarray*}
\lefteqn{\int_{t_0}^{t_1} \langle \dot{p}(t),x(t)-x_*(t)\rangle + \langle p(t), \dot{x}(t)-\dot{x}_*(t) \rangle \, dt} \\
\!\!&=&\!\! \langle p(t_1),x(t_1)-x_*(t_1)\rangle - \sum_{k=1}^d \langle p(s_k)-p(s_k^-),x(s_k)-x_*(s_k)\rangle
  -\langle p(t_0),x(t_0)-x_*(t_0)\rangle \\
\!\!&\geq&\!\! -\langle l_1,x(t_1)-x_*(t_1)\rangle -\langle l_0,x(t_0)-x_*(t_0)\rangle.
\end{eqnarray*}
Mit der gleichen Argumentationn wie im Abschnitt \ref{AbschnittHBPMP} folgt damit 
$$J\big(x(\cdot),u(\cdot)\big)-J\big(x_*(\cdot),u_*(\cdot)\big) \geq 0$$
f"ur alle zul"assigen $\big(x(\cdot),u(\cdot)\big)$ mit $\|x(\cdot)-x_*(\cdot)\|_\infty < \gamma$. \hfill $\blacksquare$

\begin{beispiel}
{\rm Im\index{Speichermanagement eines Wasserkraftwerkes}
Folgenden untersuchen wir ein optimales Speichermanagement eines Wasserkraftwerkes.
\begin{figure}[h]
	\centering
	\fbox{\includegraphics[width=12cm]{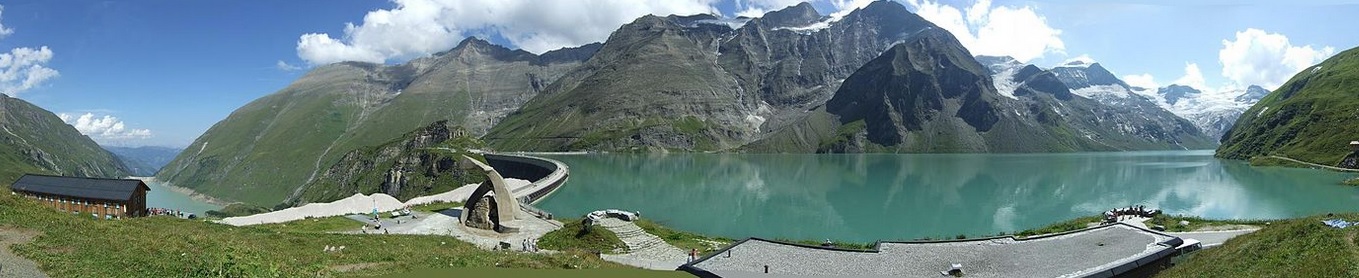}}
	\caption[Panorama-Ansicht des Stausees Mooserboden]{Panorama-Ansicht des Stausees Mooserboden.}
\end{figure}
Es bezeichne $x(t)$ die H"ohe des Stausees innerhalb der Stauanlage.
Die nat"urliche Zulaufmenge, die in den Stausee einflie"st, wird durch die Funktion $w(t)>0$ beschrieben.
Gesteuert wird die Stauh"ohe $x(t)$ durch die Ablaufmenge $v(t)$.
Weiterhin bezeichnen $x_{\min},x_{\max}$ die minimale bzw. maximale Stauh"ohe und
$v_{\min}, v_{\max}$ die minimale bzw. maximale Ablaufmenge.
Die erzeugte elektrische Energie bei der Energieumwandlung bestimmt sich aus physikalischen Gr"unden mit Hilfe der streng
konkaven Produktionsfunktion $U$.
Au"serdem sei $\pi$ der Preis, der pro Einheit elektrischer Energie am Markt erzielt wird.
Damit stellt sich die Aufgabe wie folgt dar:
\begin{eqnarray*}
&& \hspace*{-5mm} J\big(x(\cdot),v(\cdot)\big) = \int_0^T \pi U\big(x(t)\big)v(t) \, dt \to \sup,  \\
&& \hspace*{-5mm} \dot{x}(t) = w(t)-v(t), \qquad x(0)=x_0, \quad x(T) \mbox{ frei},  \\
&& \hspace*{-5mm} \displaystyle v \in [v_{\min},v_{\max}], \qquad 0\leq v_{\min} < v_{\max},  \\
&& \hspace*{-5mm} x_{\min} \leq x(t) \leq x_{\max}, \quad v_{\min} < w(t) < v_{\max} \quad\mbox{f"ur alle } t \in [0,T].
\end{eqnarray*}
Im normalen Fall lautet die Pontrjagin-Funktion $H^{\mathcal{S}}(t,x,v,p,1)= p[w(t)-v] + \pi U(x)v$.
Da die Produktionsfunktion $U$ streng konkav und die Zustandsbeschr"ankungen linear in $x$ sind,
sind die notwendigen Optimalit"atsbedingungen (\ref{SatzPMPZA1})--(\ref{SatzPMPZA3}) nach Theorem \ref{SatzHBZA} auch
hinreichend.
Die Maximumbedingung (\ref{SatzPMPZA3}) liefert
$$H^{\mathcal{S}}\big(t,x_*(t),v_*(t),p(t),1\big) = \max_{v \in [v_{\min},v_{\max}]} v \cdot \big[\pi U\big(x_*(t)\big)-p(t)\big] +p(t)w(t).$$
Daraus ergibt sich:
$\qquad \displaystyle v_*(t) \left\{\begin{array}{ll} = v_{\max}, & \pi U\big(x_*(t)\big)>p(t), \\
                                  \in [v_{\min},v_{\max}], & \pi U\big(x_*(t)\big)=p(t), \\
                                  = v_{\min}, & \pi U\big(x_*(t)\big)<p(t). \end{array}\right.$ \\[1mm]
Ist keine Zustandsbeschr"ankung aktiv, dann lautet die adjungierte Gleichung (\ref{SatzPMPZA1})
$$\dot{p}(t)=-\pi U'\big(x_*(t)\big)v_*(t).$$
Die Adjungierte ist in diesem Fall monoton fallend und wir erhalten unmittelbar
\begin{equation} \label{BeispielStaudamm}
\pi \frac{d}{dt} U\big(x_*(t)\big) = -\pi U'\big(x_*(t)\big) \cdot \big(v_*(t)-w(t)\big) > -\pi U'\big(x_*(t)\big)v_*(t) =\dot{p}(t).
\end{equation}
D.\,h., dass f"ur $x_*(t) \in (x_{\min},x_{\max})$ keine singul"are Steuerung $v_*(t) \in (v_{\min},v_{\max})$ auftritt.
Ist auf einem Teilst"uck die Zustandsbeschr"ankung
$g_1\big(x(t)\big)=x(t)-x_{\max} \leq 0$
aktiv, dann erhalten wir aus $\pi U\big(x_*(t)\big)=p(t)$ f"ur die st"uckweise stetige Funktion $\lambda_1(\cdot)$ mit
$$\dot{p}(t)=- H^{\mathcal{S}}_x\big(t,x_*(t),v_*(t),p(t),1\big) + \lambda_1(t) g_1'\big(x_*(t)\big)$$
die Darstellung $\lambda_1(t)=\pi U'(x_{\max}) w(t) \geq 0$.
Im Gegensatz dazu darf die Zustandsbeschr"ankung
$g_2\big(x(t)\big)=x_{\min}-x(t) \leq 0$
auf keinem Teilst"uck aktiv sein, denn dies h"atte $\lambda_2(t)=-\pi U'(x_{\min}) w(t) < 0$ zur Folge.
Andererseits kann $p(\cdot)$ in $t=T$ den Wert
$$p(T)=-\beta_2 g_2'\big(x_*(T)\big) = \beta_2 \geq 0$$
aufweisen, wenn die Zustandsbeschr"ankung $g_2\big(x(t)\big)$ in $t=T$ aktiv ist.
Damit ergeben sich auf gewissen Teilintervallen von $[0,T]$ folgende Verhaltensweisen, die auftreten k"onnen:
\begin{enumerate}
\item[(A)] F"ur $\pi U\big(x_*(t)\big)<p(t)$ wird das Becken mit der Rate $v_*(t)=w(t)-v_{\min}>0$ gef"ullt.
\item[(B)] F"ur $\pi U\big(x_*(t)\big)=\pi U(x_{\max})=p(t)$ wird das maximale Stauvolumen in dem Wasserkraftwerk gehalten.
           Die zugeh"orige Steuerung ist $v_*(t)=w(t)$.
\item[(C)] F"ur $\pi U\big(x_*(t)\big)>p(t)$ wird das Becken mit der Rate $v_*(t)=v_{\max}-w(t)>0$ geleert.
\end{enumerate}
In Abh"angigkeit des Planungszeitraumes ergeben sich folgende Szenarien:
\begin{enumerate}
\item[(i)] Ist der Zeitraum $[0,T]$ kurz, dann tritt nur (C) ein.
\item[(ii)] Bei einer mittelfristigen Planung ergibt sich der "Ubergang (A)$\to$(C).
\item[(iii)] Im Fall eines langfristigen Zeithorizontes erhalten wir (A)$\to$(B)$\to$(C).
\end{enumerate}
Wir demonstrieren die Szenarien anhand folgender gew"ahlter Daten:
\begin{eqnarray*}
&& \hspace*{-5mm} J\big(x(\cdot),v(\cdot)\big) = \int_0^T \ln\big(x(t)\big)v(t) \, dt \to \sup,  \\
&& \hspace*{-5mm} \dot{x}(t) = 1-v(t), \quad x(0)=5, \quad x(T) \mbox{ frei},  \\
&& \hspace*{-5mm} v \in [0,2], \quad 1 \leq x(t) \leq 20, \quad \pi = 1.
\end{eqnarray*}

In Abh"angigkeit vom Planungszeitraum $[0,T]$ ergeben sich:
\begin{enumerate}
\item[(1)] F"ur $0<T\leq 5-\sqrt{5}$ tritt nur (C) ein. Durch
           $$v_*(t)= 2, \qquad x_*(t)=5-t, \qquad x_*(T) \in [\sqrt{5},5)$$
           wird der optimale Steuerungsprozess gegeben.
           Die Adjungierte lautet
           $$p(t)=2 [\ln(5-t)-\ln(5-T)], \qquad p(T)=0.$$
           Unter Beachtung von (\ref{BeispielStaudamm}) und $U\big(x_*(0)\big) \geq p(0)$ erhalten wir $T\leq 5-\sqrt{5}$.
\item[(2)] F"ur $5-\sqrt{5}<T\leq 15-\sqrt{10}$ liegt der Fall (A)$\to$(C) vor.
           Es gelten
           $$v_*(t)= \left\{\begin{array}{ll} 0, & t \in [0,\tau), \\ 2,& t \in [\tau,T], \end{array}\right. \qquad
             x_*(t)= \left\{\begin{array}{ll} 5+t, & t \in [0,\tau), \\ 5+2\tau-t,& t \in [\tau,T], \end{array}\right.$$
           $x_*(T) \in (\sqrt{5},\sqrt{10}]$ und die zugeh"orige Adjungierte ist
           $$p(t)= \left\{\begin{array}{ll} 2 [\ln(5+t)-\ln(5+2\tau-T)], & t \in [0,\tau), \\ 
                                      2 [\ln(5+2\tau-t)-\ln(5+2\tau-T)], & t \in [\tau,T], \end{array}\right. \qquad p(T)=0.$$
           Die Beziehungen $\tau \leq 5$ und $T-\tau= 5+\tau-\sqrt{5+\tau}$ ergeben sich aus den Bedingungen
           $$x_*(t)<10 \quad\mbox{f"ur alle } t \not=\tau, \qquad U\big(x_*(\tau)\big) = p(\tau).$$
\item[(3)] Im Fall eines langfristigen Horizontes $T>15-\sqrt{10}$ erhalten wir (A)$\to$(B)$\to$(C).
           Der optimale Steuerungsprozess ist
           $$v_*(t)= \left\{\begin{array}{ll} 0, & t \in [0,\tau_1), \\ 1,& t \in [\tau_1,\tau_2),\\
                                              2,& t \in [\tau_2,T], \end{array}\right. \qquad
             x_*(t)= \left\{\begin{array}{ll} 5+t, & t \in [0,\tau_1), \\ 10,& t \in [\tau_1,\tau_2),\\
                                              10-(t-\tau_2),& t \in [\tau_2,T], \end{array}\right.$$
           es gilt $x_*(T)=\sqrt{10}$ und die zugeh"orige Adjungierte lautet
           $$p(t)= \left\{\begin{array}{ll} \ln(10), & t \in [0,\tau_2), \\ 
                                            2 [\ln(10+\tau_2-t)-\ln(10+\tau_2-T)], & t \in [\tau_2,T], \end{array}\right. \qquad p(T)=0.$$
           Es ergeben sich $\tau_1=5$ und $T-\tau_2=10-\sqrt{10}$ aus den Beziehungen
           $$x_*(t)<10,\; t \in [0,\tau_1),\quad x_*(\tau_1)=10, \quad \ln(10) = p(\tau_2), \quad p(T)=0.$$
\end{enumerate}
Eine interessante Erweiterung ergibt sich im letzten Fall,
wenn der Preis $\pi$ nicht konstant ist,
sondern im Zeitpunkt $\sigma \in (\tau_1,\tau_2)$ eine Preiserh"ohung $\pi^+$ bzw. -verringerung $\pi^-$ erwartet wird.
Im Fall der Preiserh"ohung wird die interne Bewertung (Schattenpreis\index{Schattenpreis}) der Ressource im Zeitpunkt $t =\sigma$ durch
eine Sprungstelle der Adjungierten angepasst:
$$p(\sigma^-)-p(\sigma)=(\pi-\pi^{+}) U(x_{\max}) = -\beta_1 g_1'\big(x_*(\sigma)\big) = -\beta_1 < 0.$$
Wird eine Preissenkung erwartet,
dann wird die Ressource ab dem Zeitpunkt $\sigma_1$ vor der "Anderung zum Preis $\pi$ mit $v_*(t)=v_{\max}$ ver"au"sert und danach
mit der Steuerung $v_*(t)=v_{\min}$ eingelagert.
Und zwar in der Form,
dass die Bewertung der Ressource durch die Bedingung $p(\sigma_2)=\pi^{-} U(x_{\max})$ dem aktuellen Marktpreis $\pi^{-}$
angepasst wird.
Dieser Argumentation ist die Annahme $[\sigma_1,\sigma_2] \subseteq [\tau_1,\tau_2]$ auferlegt.
In dem spezifizierten Modell erhalten wir im Fall der Preissenkung mit $\pi=2$ und $\pi^{-}=1$ f"ur die Zeitpunkte den Zusammenhang
$\sigma-\sigma_1=\sigma_2-\sigma=10-\sqrt{10}$. \hfill $\square$}
\end{beispiel}

%% file: 2-53-FreieZeit.tex
\subsubsection{Freier Anfangs- und Endzeitpunkt} \label{AbschnittFreieZeit}
Wir betrachten in diesem Abschnitt die zustandsbeschr"ankte Aufgabe mit freiem Anfangs- und Endzeitpunkt.
Ebenso wie im Abschnitt \ref{AbschnittFreieZeitSOP} seien daher die Abbildungen $h_0$ bzw. $h_1$,
die die Start- und Zielmannigfaltigkeit definieren, zeitabh"angig:
$$h_i (t_i,x_i) : \R \times \R^n \to \R^{s_i }, \quad i = 0,1.$$
Mit den getroffenen Bezeichnungen und Festlegungen untersuchen wir die Aufgabe:
\begin{eqnarray}
&& \label{FZ1} J\big(x(\cdot),u(\cdot)\big) = \int_{t_0}^{t_1} f\big(t,x(t),u(t)\big) \, dt \to \inf, \\
&& \label{FZ2} \dot{x}(t) = \varphi\big(t,x(t),u(t)\big), \\
&& \label{FZ3} h_0\big(t_0,x(t_0)\big)=0, \qquad h_1\big(t_1,x(t_1)\big)=0, \\
&& \label{FZ4} u(t) \in U \subseteq \R^m, \quad U \not= \emptyset, \\
&& \label{FZ5} g_j\big(t,x(t)\big) \leq 0 \quad \mbox{f"ur alle } t \in [t_0,t_1], \quad j=1,...,l.
\end{eqnarray}
Die Aufgabe (\ref{FZ1})--(\ref{FZ5}) betrachten wir f"ur Tripel $\big([t_0,t_1],x(\cdot),u(\cdot)\big)$ mit
$$[t_0,t_1] \subset \R, \qquad x(\cdot) \in W^1_\infty([t_0,t_1],\R^n), \qquad u(\cdot) \in L_\infty([t_0,t_1],U).$$
Zur Menge $\mathscr{B}^{\,\mathcal{F}}_{\rm Lip}$ geh"oren diejenigen Tripel
$\big([t_0,t_1],x(\cdot),u(\cdot)\big)$,
f"ur die es eine Zahl $\gamma>0$ derart gibt,
dass die Abbildungen $f(t,x,u)$, $\varphi(t,x,u)$, $g_j(t,x)$ und $h_i(\tau_i,x_i)$
auf der Menge aller Punkte $(t,\tau_0,\tau_1,x,x_0,x_1,u) \in \R \times \R \times \R \times \R^n \times \R^n \times \R^n \times \R^m$ mit
\begin{eqnarray*}
&& t_0-\gamma < t < t_1+\gamma,\quad t_0-\gamma < \tau_0< t_0+\gamma,\quad t_1-\gamma < \tau_1 < t_1+\gamma, \\
&& \|x-x(t)\| < \gamma, \quad \|x_0-x(t_0)\| < \gamma, \quad \|x_1-x(t_1)\| < \gamma, \quad u \in \R^m
\end{eqnarray*}
stetig in der Gesamtheit aller Variablen und
stetig differenzierbar bez"uglich der Variablen $t, t_0, t_1, x, x_0, x_1$ sind.
(Zur unmissverständlichen Angabe der Punktemenge treten $\tau_0,\tau_1$ in $h_i$ anstelle der Zeitvariablen $t_0,t_1$ auf.) \\[2mm]
In der Aufgabe (\ref{FZ1})--(\ref{FZ5}) mit freiem Anfangs- und Endzeitpunkt nennen wir ein Tripel
$\big([t_0,t_1],x(\cdot),u(\cdot)\big)$ mit
$[t_0,t_1] \subset \R$, $x(\cdot) \in W^1_\infty([t_0,t_1],\R^n)$ und $u(\cdot) \in L_\infty([t_0,t_1],U)$ einen Steuerungsprozess.
Ein Steuerungsprozess hei"st zul"assig in dem Steuerungsproblem (\ref{FZ1})--(\ref{FZ5}),
wenn auf dem Intervall $[t_0,t_1]$ die Funktion $x(\cdot)$ fast "uberall der Gleichung (\ref{FZ2}) gen"ugt,
die Randbedingungen (\ref{FZ3}) und die Zustandsbeschr"ankungen (\ref{FZ5}) erf"ullt.
Die Menge $\mathscr{B}^{\,\mathcal{F}}_{\rm adm}$ bezeichnet die Menge der zul"assigen Steuerungsprozesse. \\[2mm]
Einen zul"assigen Steuerungsprozess $\big([t_{0*},t_{1*}],x_*(\cdot),u_*(\cdot)\big)$ nennen wir ein starkes lokales
Minimum\index{Minimum, starkes lokales!Standardl@-- Standardaufgabe},
wenn eine Zahl $\varepsilon > 0$ derart existiert,
dass f"ur jeden zul"assigen Steuerungsprozess $\big([t_0,t_1],x(\cdot),u(\cdot)\big)$ mit den Eigenschaften
$$|t_0 - t_{0*}| < \varepsilon, \qquad |t_1 - t_{1*}| < \varepsilon, \qquad
  \| x(t)-x_*(t) \| < \varepsilon \quad \mbox{ f"ur alle } t \in [t_{0*},t_{1*}] \cap [t_0,t_1]$$
die Ungleichung $J\big(x(\cdot),u(\cdot)\big) \geq J\big(x_*(\cdot),u_*(\cdot)\big)$ gilt.
\newpage
Wir wenden die Methode der Substitution der Zeit aus Abschnitt \ref{AbschnittFreieZeitSOP} an und erhalten folgende Aufgabe
mit fester Zeit, auf die wir Theorem \ref{SatzPMPZA} anwenden k"onnen:
\begin{eqnarray}
&& \label{FZ8} \int_0^1 v \cdot f\big(t(s),y(s),w(s)\big) \, ds \to \inf, \\
&& \label{FZ9} t'(s) = v, \quad y'(s) = v \cdot \varphi\big(t(s),y(s),w(s)\big), \\
&& \label{FZ10} h_0\big( t(0),y(0) \big) = 0, \quad h_1\big( t(1),y(1) \big) = 0, \\
&& \label{FZ11} v > 0, \quad w(s) \in U, \\
&& \label{FZ12} g_j\big(t(s),y(s)\big) \leq 0, \quad s \in [0,1], \quad j=1,...,l.
\end{eqnarray}
Die Pontrjagin-Funktion der Aufgabe (\ref{FZ8})--(\ref{FZ12}) lautet
\begin{eqnarray}
\tilde{H}(t,x,u,v,p,q,\lambda_0) &=& \langle p ,v \cdot \varphi(t,x,u) \rangle + qv - \lambda_0 v \cdot f(t,x,u) \nonumber \\
 \label{FZ13} &=& v \cdot [H^{\mathcal{S}}(t,x,u,p,\lambda_0) + q],
\end{eqnarray}
wobei $H^{\mathcal{S}}(t,x,u,p,\lambda_0) = \langle p , \varphi(t,x,u) \rangle - \lambda_0 f(t,x,u)$
die Pontrjagin-Funktion zu der Aufgabe (\ref{FZ1})--(\ref{FZ5}) ist.
Weiterhin seien $\Delta_j = \big\{ s \in [0,1] \, \big| \, g_j\big( t_*(s),y_*(s) \big) = 0 \big\}$ für $j=1,...,l$.
Dann existieren nach Theorem \ref{SatzPMPZA} eine Zahl $\lambda_0 \geq 0$, Vektoren $l_0 \in \R^{s_0},l_1 \in \R^{s_1},$
eine Vektorfunktion $\tilde{p}(\cdot)$,
eine Funktion $\tilde{q}(\cdot)$ und nichtnegative, regul"are,
auf den Mengen $\Delta_j$ konzentrierte Borelsche Ma"se
$\tilde{\mu}_j, j=1,...,l,$ (diese Gr"o"sen verschwinden nicht gleichzeitig) derart, dass die Beziehungen
\begin{eqnarray*} 
\tilde{p}(s) &=& - h^T_{1,x_1}\big(t_*(1),y_*(1)\big) l_1
                 +\int_s^1 \tilde{H}_x\big(t_*(\sigma),y_*(\sigma),w_*(\sigma),v_*,\tilde{p}(\sigma),\tilde{q}(\sigma),\lambda_0\big) \, d\sigma \\
             & & - \sum_{j=1}^l \int_s^1 g_{j,x}\big(t_*(\sigma),y_*(\sigma)\big)\, d\tilde{\mu}_j(\sigma), \\
\tilde{p}(0) &=& h^T_{0,x_0} \big(t_*(0),y_*(0)\big) l_0, \\
\tilde{q}(s) &=& - \big\langle h_{1,t_1}\big(t_*(1),y_*(1)\big),l_1 \big\rangle
                 +\int_s^1 \tilde{H}_t\big(t_*(\sigma),y_*(\sigma),w_*(\sigma),v_*,\tilde{p}(\sigma),\tilde{q}(\sigma),\lambda_0\big) \, d\sigma \\
             & & - \sum_{j=1}^l \int_s^1 g_{j,t}\big( t_*(\sigma),y_*(\sigma)\big) \, d\tilde{\mu}_j(\sigma), \\
\tilde{q}(0) &=& \big\langle h_{0,t_0} \big(t_*(0),y_*(0)\big) , l_0 \big\rangle 
\end{eqnarray*}
gelten und die Maximumbedingung
$$\tilde{H}\big(t_*(s),y_*(s),w_*(s),v_*,\tilde{p}(s),\tilde{q}(s),\lambda_0\big)
  = \max_{u \in U,\, v> 0} v \cdot \big[ H^{\mathcal{S}}\big(t_*(s),y_*(s),u,\tilde{p}(s),\lambda_0\big) + \tilde{q}(s) \big]$$
f"ur fast alle $s \in [0,1]$ erf"ullt ist.
Es sei $s_*(\cdot)$ die zu $t_*(\cdot)$ inverse Funktion, d.\,h.
$$s_*(t) = \frac{t-t_{0*}}{t_{1*}-t_{0*}}, \qquad t \in [t_{0*},t_{1*}].$$
Wir f"uhren die Bezeichnungen
\begin{eqnarray*}
&& p(t) = \tilde{p}\big(s_*(t)\big), \quad q(t) = \tilde{q}\big(s_*(t)\big), \\
&& T_j = \big\{ t \in [t_{0*},t_{1*}] \, \big| \, g_j\big( t,x_*(t) \big) = 0 \big\} = t_*(\Delta_j), \quad j=1,...,l,
\end{eqnarray*}
ein und verwenden die Bildma"se $\mu_j$ der Ma"se $\tilde{\mu}_j$ bei der Abbildung $t_*(\cdot)$.
Dann gilt
$$\int_{t_{0*}}^{t_{1*}} \psi(t) \, d\mu_j(t) = \int_0^1 \psi\big(t_*(s)\big) \, d\tilde{\mu}_j(s) \qquad \mbox{ f"ur alle }
  \psi(\cdot) \in C([t_{0*},t_{1*}],\R)$$
und die Ma"se $\mu_j$ sind offensichtlich auf $T_j$ konzentriert.
So lassen sich die obigen Beziehungen unter Ber"ucksichtigung von (\ref{FZ13}) auf folgende Form bringen:
\begin{eqnarray} 
p(t) &=& -h^T_{1,x_1}\big(t_{1*},x_*(t_{1*})\big) l_1
         + \int_t^{t_{1*}} H^{\mathcal{S}}_x\big(s,x_*(s),u_*(s),p(s),\lambda_0\big) \, ds \nonumber \\
     & & - \sum_{j=1}^l \int_t^{t_{1*}} g_{j,x}\big( s,x_*(s)\big) \, d\mu_j(s), \nonumber \\
p(t_{0*}) &=& h^T_{0,x_0} \big(t_{0*},x_*(t_{0*})\big) l_0, \nonumber \\
q(t) &=& - \big\langle h_{1,t_1}\big(t_{1*},x_*(t_{1*})\big) , l_1 \big\rangle
         + \int_t^{t_{1*}} H^{\mathcal{S}}_t\big(s,x_*(s),u_*(s),p(s),\lambda_0\big) \, ds \nonumber \\
     & & \label{FZ16} - \sum_{j=1}^l \int_t^{t_{1*}} g_{j,t}\big( s,x_*(s)\big) \, d\mu_j(s), \\
q(t_{0*}) &=& \label{FZ17} \big\langle h_{0,t_0}\big(t_{0*},x_*(t_{0*})\big) , l_0 \big\rangle
\end{eqnarray}
und wir erhalten die Maximumbedingung
\begin{equation} \label{FZ18}
v_* \cdot \Big( H^{\mathcal{S}}\big(t,x_*(t),u_*(t),p(t),\lambda_0\big) + q(t) \Big)
=  \max_{u \in U,\, v> 0} v \cdot \big[ H^{\mathcal{S}}\big(t,x_*(t),u,p(t),\lambda_0\big)  + q(t) \big].
\end{equation}
Aus (\ref{FZ18}) folgt
$$H^{\mathcal{S}}\big(t,x_*(t),u_*(t),p(t),\lambda_0\big) = \max_{u \in U} H^{\mathcal{S}}\big(t,x_*(t),u,p(t),\lambda_0\big).$$
Ferner erhalten wir wegen $0 < v_* < \infty$ aus (\ref{FZ18}) die Beziehung
$$\max_{u \in U} H^{\mathcal{S}}\big(t,x_*(t),u,p(t),\lambda_0\big) = -q(t).$$
Nachstehend verwenden wir die Hamilton-Funktion $\mathscr{H}^{\mathcal{S}}(t,x,p,\lambda_0) = \max\limits_{u \in U} H^{\mathcal{S}}(t,x,u,p,\lambda_0)$:
Durch Vergleich mit (\ref{FZ16}) und (\ref{FZ17}) erhalten wir weiter
\begin{eqnarray*}
&& \mathscr{H}^{\mathcal{S}}\big(t,x_*(t),p(t),\lambda_0\big)
   = \big\langle h_{1,t_1}\big(t_{1*},x_*(t_{1*})\big) , l_1 \big\rangle
   - \int_t^{t_{1*}} H^{\mathcal{S}}_t\big(s,x_*(s),u_*(s),p(s),\lambda_0\big) \, ds \\
&& \hspace*{45mm} + \sum_{j=1}^l \int_t^{t_{1*}} g_{j,t}\big( s,x_*(s)\big) \, d\mu_j(s), \\
&& \mathscr{H}^{\mathcal{S}}\big(t_{0*},x_*(t_{0*}),p(t_{0*}),\lambda_0\big)
   = - \big\langle h_{0,t_0} \big(t_{0*},x_*(t_{0*})\big) , l_0 \big\rangle.
\end{eqnarray*}
Zusammenfassend lautet damit das Maximumprinzip f"ur Aufgaben mit freier Zeit: 

\begin{theorem}[Pontrjaginsches Maximumprinzip] \label{SatzFZ}
\index{Pontrjaginsches Maximumprinzip!Standard@-- Standardaufgabe} 
In der Aufgabe (\ref{FZ1})--(\ref{FZ5}) sei der Steuerungsprozess
$\big([t_{0*},t_{1*}],x_*(\cdot),u_*(\cdot)\big) \in \mathscr{B}^{\,\mathcal{F}}_{\rm adm} \cap \mathscr{B}^{\,\mathcal{F}}_{\rm Lip}$.
Ist $\big([t_{0*},t_{1*}],x_*(\cdot),u_*(\cdot)\big)$ ein starkes lokales Minimum der Aufgabe (\ref{FZ1})--(\ref{FZ5}),
dann existieren eine Zahl $\lambda_0 \geq 0$, Vektoren $l_0 \in \R^{s_0}$ und $l_1 \in \R^{s_1}$, eine Vektorfunktion $p(\cdot):[t_{0*},t_{1*}] \to \R^n$
und auf den Mengen
$$T_j=\big\{t \in [t_{0*},t_{1*}] \,\big|\, g_j\big(t,x_*(t)\big)=0\big\}, \quad j=1,...,l,$$
konzentrierte nichtnegative regul"are Borelsche Ma"se $\mu_j$ endlicher Totalvariation
(wobei s"amtliche Gr"o"sen nicht gleichzeitig verschwinden) derart,
dass die Vektorfunktion $p(\cdot)$ von beschr"ankter Variation und rechtsseitig stetig ist, und
\begin{enumerate}
\item[(a)] die adjungierte Gleichung
           \index{adjungierte Gleichung!Standard@-- Standardaufgabe}
           \begin{eqnarray}
           p(t)&=&-h^T_{1,x_1}\big(t_{1*},x_*(t_{1*})\big) l_1 
           + \int_t^{t_{1*}} H^{\mathcal{S}}_x\big(s,x_*(s),u_*(s),p(s),\lambda_0\big) \, ds \nonumber \\
           \label{SatzFZ1} & & -\sum_{j=1}^l \int_t^{t_{1*}} g_{j,x}\big(s,x_*(s)\big)\, d\mu_j(s),
           \end{eqnarray}
\item[(b)] die Transversalit"atsbedingungen
           \index{Transversalitätsbedingungen!Standard@-- Standardaufgabe}
           \begin{equation}\label{SatzFZ2}
            \left. \begin{array}{lcl} p(t_{0*}) &=& h^T_{0,x_0} \big(t_{0*},x_*(t_{0*})\big) l_0, \\[2mm]
                  p(t_{1*}^-)-p(t_{1*}) &=& \displaystyle - \sum_{j=1}^l \mu_j(\{t_{1*}\}) \, g_{j,x}\big(t_{1*},x_*(t_{1*})\big)
                             \end{array} \right\}
           \end{equation} 
\item[(c)] und in fast allen Punkten $t\in [t_{0*},t_{1*}]$ die Maximumbedingung 
           \index{Maximumbedingung!Standard@-- Standardaufgabe}
           \begin{equation}\label{SatzFZ3}
           H^{\mathcal{S}}\big(t,x_*(t),u_*(t),p(t),\lambda_0\big) = \max_{u \in U} H^{\mathcal{S}}\big(t,x_*(t),u,p(t),\lambda_0\big)
           \end{equation}
\end{enumerate}
erfüllt sind und weiterhin
\begin{enumerate}
\item[(d)] die Beziehungen 
           \begin{eqnarray}
           && \hspace*{-10mm} \mathscr{H}^{\mathcal{S}}\big(t,x_*(t),p(t),\lambda_0\big) 
              = \big\langle h_{1,t_1}\big(t_{1*},x_*(t_{1*})\big) , l_1 \big\rangle
               - \int_t^{t_{1*}} H^{\mathcal{S}}_t\big(s,x_*(s),u_*(s),p(s),\lambda_0\big) \, ds \nonumber \\
           && \hspace*{30mm} + \sum_{j=1}^l \int_t^{t_{1*}} g_{j,t}\big(s,x_*(s)\big) \, d\mu_j(s), \\
           && \hspace*{-10mm}  \mathscr{H}^{\mathcal{S}}\big(t_{0*},x_*(t_{0*}),p(t_{0*}),\lambda_0\big)
              = - \big\langle h_{0,t_0}\big(t_{0*},x_*(t_{0*})\big) , l_0 \big\rangle
           \end{eqnarray}
\end{enumerate}
gelten.
\end{theorem}

%% file: 3-0-Erweiterungen.tex
\section{Erweiterungen der Standardaufgabe}
Im ersten Kapitel wurde ausgehend von den Richtungsvariationen die Euler-Lagrangesche Gleichung und
die Schwachen Optimalitätsprinzipien hergeleitet.
Nachfolgend wurden im zweiten Kapitel auf der Basis von Nadelvariationen das Pontrjaginsche Maximumprinzip für die
Standardaufgabe über endlichem Zeithorizont entwickelt.
Dieses Kapitel befasst sich nun mit verschiedenen Erweiterungen der Standardaufgabe.
Die Charakteristiken der nachstehenden Themen sind grundverschiedenen;
folglich von eigenständigem Interesse. \\[2mm]
In der Literatur gibt es häufig abweichende Formulierungen der Standardaufgabe,
nämlich in der Gestalt des Zielfunktionals oder in der Beschreibung der Randbedingungen.
Aus diesem Grund gehen wir auf gemischte Zielfunktionale bestehend aus Integral- und Terminalfunktionalen ein.
Weiterhin ergeben sich in den nachfolgenden Abschnitten Randbedingungen,
auf welche die bisherige Form nicht passt.
Deswegen verallgemeinern wir außerdem die Formulierungen der Randbedingungen in Gleichungs- und Ungleichungsform. \\[2mm]
Im Rahmen der Multiprozesse besteht der Ausbau der bisherigen Betrachtungen darin,
dass eine spezifische Steuervariable ``Wechselstrategie'' eingeführt wird.
Diese Variable darf ausschließlich Werte annehmen,
die Eckpunkte eines $k$-dimensionalen Simplexes darstellen.
Die mathematische Modellierung eines ``Springen'' zwischen verschiedenen Eckpunkten ist mittels Richtungsvariationen
nicht realisierbar.
An dieser Stelle zeigt sich der Vorteil der Nadelvariationsmethode,
mit der das Springen zwischen den Eckpunkten beschrieben und in die Modellierung aufgenommen werden kann. \\
Dies offenbart neben der Verallgemeinerung des Optimalitätsbegriffes vom schwachen zum starken lokalen Minimum
einen neuen Aspekt der Nadelvariationsmethode:
Sie ermöglicht eine flexiblere und realistischere Modellierung der Problemstellung. \\[2mm]
Differentialgleichungen mit Zeitverzögerungen von fester Zeitspanne unterscheiden sich vom Standardtyp dynamischer Systeme.
Jedoch tritt mittels eines Rekursionsansatzes die Eigenschaft hervor,
dass sich eine zeitverzögerte Differentialgleichung auf die sukzessive Betrachtung der Differentialgleichung über
vorangegangene Teilintervalle gegebener Länge zurückführen lässt.
Auf diese Weise kann man eine verzögerte Differentialgleichung in ein System üblicher Differentialgleichungen überführen.
Das Zurückführen des rekursiven Ersatzsystems induziert letztlich zeitverzögerte Terme im Pontrjaginschen Maximumprinzip
für die Aufgabe mit verzögerten Steuerungsprozessen. \\[2mm]
In der Spieltheorie untersucht man Entscheidungssitutationen,
in denen Interessenkonflikte durch das Auftreten verschiedener Zielkriterien zustande kommen.
Charakteristisch für Differentialspiele sind dabei dynamische Nebenbedingungen,
die durch ein Differentialgleichungssystem gegeben sein können.
Die Spielsituation entsteht durch das Zuordnen der verschiedenen Zielkriterien zu einzelnen ``Spielern'',
die in der ``Spielsituation'' Konkurrenten um das bestmögliche Ergebnis darstellen.
Lösungsstrategien für Differentialspiele sind Gleichgewichtskonzepte wie zum Beispiel das Nash-Gleichgewicht,
welches wir vorstellen und in Beispielen anwenden werden.

%% file: 3-10-GemischtesZF.tex
\subsection{Zielfunktionale in gemischter Form und allgemeine Randbedingungen} \label{AbschnittGemischtesZF}
\index{Zielfunktional!gemischt@-- in gemischter Form}
Bisher beschränkten wir unsere Betrachtungen auf Zielfunktionale in Integralform.
Eine unmittelbare Erweiterung ist der Einbezug von Terminalfunktionalen. \index{Zielfunktional!Terminal@-- Terminalfunktional}
Terminalfunktionale waren bereits Bestandteil der einfachsten Bolza-Aufgabe in Abschnitt \ref{BolzaAufgabe}.
In diesem Abschnitt binden wir Terminalfunktionale in die Standardaufgabe (\ref{PMP1})--(\ref{PMP5}) ein. \\
Gekoppelte Randbedingungen der Form $x(t_0)=x(t_1)$ traten in der Standardaufgabe bisher nicht auf,
da diese Bedingungen für Anfangs- und Endzeitpunkt getrennt wurden.
Wir werden nun allgemeinere Randbedingungen in Gleichungs- und Ungleichungsform in die Aufgabe aufnehmen. \\[1mm]
Die Aufgabe (\ref{PMP1})--(\ref{PMP5}) mit Terminalfunktionalen und Randbedingungen in allgemeiner
Gleichungs- und Ungleichungsform besitzt die Gestalt
\begin{eqnarray}
&& \label{GZF1} \tilde{J}\big(x(\cdot),u(\cdot)\big)
                = \int_{t_0}^{t_1} f\big(t,x(t),u(t)\big) \, dt + S\big(x(t_0),x(t_1)\big) \to \inf, \\
&& \label{GZF2} \dot{x}(t) = \varphi\big(t,x(t),u(t)\big), \\
&& \label{GZF3} h\big(x(t_0),x(t_1)\big)=0, \qquad  r\big(x(t_0),x(t_1)\big) \leq 0,\\
&& \label{GZF4} u(t) \in U \subseteq \R^m, \quad U\not= \emptyset, \\
&& \label{GZF5} g_j\big(t,x(t)\big) \leq 0 \quad \mbox{f"ur alle } t \in [t_0,t_1], \quad j=1,...,l,
\end{eqnarray}
wobei wir sämtliche Annahmen und Bezeichnungen aus Kapitel \ref{KapitelStark} übernehmen.
Das Terminalfunktional $S$ und die Randbegingungen fassen wir als Funktionen der Variablen $x_0,\, x_1$ auf,
d.\,h. $S=(x_0,x_1)$, $h=h(x_0,x_1)$ und $r=r(x_0,x_1)$.
Ferner setzen wir die Abbildungen
$S:\R^n \times \R^n \to \R$, $h:\R^n \times \R^n \to \R^s$ und $r: \R^n \times \R^n \times \R$ als stetig differenzierbar voraus.
Zwischen (\ref{PMP1}) und (\ref{GZF1}) gilt außerdem der Zusammenhang
\begin{equation} \label{GFZZusammenhang}
\tilde{J}\big(x(\cdot),u(\cdot)\big)= J\big(x(\cdot),u(\cdot)\big) + S\big(x(t_0),x(t_1)\big).
\end{equation}

Die Bezeichnung der Terminalfunktionale mit $S$ ist durch ökonomische Anwendungen motiviert:
Bei der Betrachtung von Strategien zur optimalen Instandhaltung einer Maschine bezeichnet $S\big(x(t_1)\big)$ den
Schrotterlös (``scrap value'') oder den Wert von wiederverwendbaren Materialien oder Teilen (``salvage value'') einer nicht
weiter nutzbaren Maschine.

%% file: 3-11-PMP.tex
\subsubsection{Das Pontrjaginsche Maximumprinzip}
Wir formulieren in diesem Abschnitt das Pontrjaginsche Maximumprinzip für verschiedene Varianten
der Aufgabe (\ref{GZF1})--(\ref{GZF5}).
Da sich die Aufgabe (\ref{GZF1})--(\ref{GZF5}) von der Standardaufgabe (\ref{PMP1})--(\ref{PMP5}) lediglich im Zielfunktional
durch die Beziehung (\ref{GFZZusammenhang}) und in den Randbedingungen (\ref{GZF3}) unterscheidet,
stimmen die Kernelemente des Pontrjaginschen Maximumprinzips und dessen Herleitung mit den Ausführungen in
Kapitel \ref{KapitelStark} überein.
Die Unterschiede manifestieren sich hier ausschließlich in den Transversalitätsbedingungen und der komplementären Schlupfbedingung
bezüglich der Randbedingung $r \leq 0$. \\[1mm]
F"ur die Aufgabe (\ref{GZF1})--(\ref{GZF5}) bezeichnet $H^{\mathcal{S}}$ die Pontrjagin-Funktion
$$H^{\mathcal{S}}(t,x,u,p,\lambda_0) = \langle p, \varphi(t,x,u) \rangle - \lambda_0 f(t,x,u).$$

\begin{theorem}[Pontrjaginsches Maximumprinzip] \label{SatzPMPGZFZA}
\index{Pontrjaginsches Maximumprinzip!gemischt@-- gemischtes Zielfunktional} 
Sei $\big(x_*(\cdot),u_*(\cdot)\big) \in \mathscr{B}^{\,\mathcal{S}}_{\rm adm} \cap \mathscr{B}^{\,\mathcal{S}}_{\rm Lip}$.
Ist $\big(x_*(\cdot),u_*(\cdot)\big)$ ein starkes lokales Minimum der Aufgabe (\ref{GZF1})--(\ref{GZF5}),
dann existieren Zahlen $\lambda_0 \geq 0$ und $\alpha \geq 0$, ein Vektor $\nu  \in \R^s$, eine Vektorfunktion $p(\cdot):[t_0,t_1] \to \R^n$
und auf den Mengen
$$T_j=\big\{t \in [t_0,t_1] \,\big|\, g_j\big(t,x_*(t)\big)=0\big\}, \quad j=1,...,l,$$
konzentrierte nichtnegative regul"are Borelsche Ma"se $\mu_j$ endlicher Totalvariation
(wobei s"amtliche Gr"o"sen nicht gleichzeitig verschwinden) derart,
dass die Vektorfunktion $p(\cdot)$ von beschr"ankter Variation und rechtsseitig stetig ist, und
\begin{enumerate}
\item[(a)] die adjungierte Gleichung
           \index{adjungierte Gleichung!gemischt@-- gemischtes Zielfunktional}
           \begin{eqnarray}
           p(t)&=& -h_{x_1}^T\big(x_*(t_0),x_*(t_1)\big) \nu - \alpha r_{x_1}\big(x_*(t_0),x_*(t_1)\big) \nonumber \\
               & &   + \int_t^{t_1} H^{\mathcal{S}}_x\big(s,x_*(s),u_*(s),p(s),\lambda_0\big) \, ds - \lambda_0 S_{x_1}\big(x_*(t_0),x_*(t_1)\big) \nonumber \\
           \label{SatzPMPGZFZA1} & & -\sum_{j=1}^l \int_t^{t_1} g_{j,x}\big(s,x_*(s)\big)\, d\mu_j(s),
           \end{eqnarray}
\item[(b)] die Transversalit"atsbedingungen
           \index{Transversalitätsbedingungen!gemischt@-- gemischtes Zielfunktional}
           \begin{equation} \label{SatzPMPGZFZA2}
            \left. \begin{array}{lcl}
            p(t_0) &=& h_{x_0}^T\big(x_*(t_0),x_*(t_1)\big) \nu  + \lambda_0 S_{x_0}\big(x_*(t_0),x_*(t_1)\big) \\[2mm]
                   & & + \alpha r_{x_0}\big(x_*(t_0),x_*(t_1)\big), \\[2mm]
            p(t_1^-)-p(t_1) &=& \displaystyle - \sum_{j=1}^l \mu_j(\{t_1\}) \, g_{j,x}\big(t_1,x_*(t_1)\big),
           \end{array} \right\}
           \end{equation}
\item[(c)] in fast allen Punkten $t \in [t_0,t_1]$ die Maximumbedingung
           \index{Maximumbedingung!gemischt@-- gemischtes Zielfunktional}
           \begin{equation}\label{SatzPMPGZFZA3}
           H^{\mathcal{S}}\big(t,x_*(t),u_*(t),p(t),\lambda_0\big) = \max_{u \in U} H^{\mathcal{S}}\big(t,x_*(t),u,p(t),\lambda_0\big)
           \end{equation}
\item[(d)] und die komplementäre Schlupfbedingung
           \begin{equation}\label{SatzPMPGZFZA4}
           \alpha r\big(x_*(t_0),x_*(t_1)\big)=0
           \end{equation}
\end{enumerate}
erfüllt sind.
\end{theorem}

{\bf Beweis} Die Randbedingungen (\ref{GZF3}) verbinden wir mit den Abbildungen
\begin{eqnarray*}
&& H\big(x(\cdot)\big)= h\big(x(t_0),x(t_1)\big), \qquad H:C([t_0,t_1],\R^n) \to \R^s, \\
&& R\big(x(\cdot)\big)= r\big(x(t_0),x(t_1)\big), \qquad R:C([t_0,t_1],\R^n) \to \R.
\end{eqnarray*}
Die Abbildungen $H$ und $R$ sind stetig Fr\'echet-differenzierbar und es gelten
\begin{eqnarray*}
H'\big(x(\cdot)\big)\xi(\cdot) &=& h_{x_0}\big(x(t_0),x(t_1)\big) \xi(t_0)+ h_{x_1}\big(x(t_0),x(t_1)\big) \xi(t_1), \\
R'\big(x(\cdot)\big)\xi(\cdot) &=& \big\langle r_{x_0}\big(x(t_0),x(t_1)\big), \xi(t_0) \big \rangle
                                   + \big\langle r_{x_1}\big(x(t_0),x(t_1)\big), \xi(t_1) \big\rangle.
\end{eqnarray*}
Zur Aufgabe (\ref{GZF1})--(\ref{GZF5}) besitzt die Lagrange-Funktion
$\tilde{\mathscr{L}}= \tilde{\mathscr{L}}\big(x_*(\cdot),u_*(\cdot),\lambda_0,y^*,\nu,\alpha,\lambda\big)$
die Gestalt
$$\tilde{\mathscr{L}}=
  \lambda_0 \tilde{J}\big(x(\cdot),u(\cdot)\big)+ \big\langle y^*, F\big(x(\cdot),u(\cdot)\big) \big\rangle
  +\nu^T H\big(x(\cdot)\big) + \alpha R\big(x(\cdot)\big) + \sum_{j=1}^l \lambda_j G_j\big(x(\cdot)\big).$$

Die Bedingungen des Extremalprinzips \ref{SatzExtremalprinzipStark} in Form von (\ref{SatzPMPZALMR1})--(\ref{SatzPMPZALMR3}) lauten:
\begin{enumerate}
\item[(a)] Die Lagrange-Funktion besitzt bez"uglich $x(\cdot)$ in $x_*(\cdot)$ einen station"aren Punkt:
           \begin{equation}\label{SatzGZFLMR1}
           0 \in \partial_x \tilde{\mathscr{L}}\big(x_*(\cdot),u_*(\cdot),\lambda_0,y^*,\nu,\alpha,\lambda\big).
           \end{equation}         
\item[(b)] Die Lagrange-Funktion erfüllt bez"uglich $u(\cdot)$ in $u_*(\cdot)$ die Minimumbedingung:
           \begin{equation}\label{SatzGZFLMR2}
           \hspace*{-3mm} \tilde{\mathscr{L}}\big(x_*(\cdot),u_*(\cdot),\lambda_0,y^*,\nu,\alpha,\lambda\big)
           = \min_{u(\cdot) \in L_\infty([t_0,t_1],U)}
           \tilde{\mathscr{L}}\big(x_*(\cdot),u(\cdot),\lambda_0,y^*,\nu,\alpha,\lambda\big).
           \end{equation}
\item[(c)] Die komplementären Schlupfbedingungen sind gültig:
           \begin{equation}\label{SatzGZFLMR3}
           \alpha R\big(x_*(\cdot)\big) = 0, \quad \lambda_1 G_1\big(x_*(\cdot)\big) = 0 \,,...,\, \lambda_l G_l\big(x_*(\cdot)\big) = 0.
           \end{equation}         
\end{enumerate}
Diese Bedingungen führen auf die gleiche Weise wie im Abschnitt \ref{AbschnittBeweisPMPZA} auf folgende Beziehungen:
Die Gleichung (\ref{SatzGZFLMR1}) kann in die Form
\begin{eqnarray*}
0 \!\!&=& \!\! \int_{t_0}^{t_1} \Big\langle \lambda_0 f_x\big(t,x_*(t),u_*(t)\big) 
                            - \varphi_x^T\big(t,x_*(t),u_*(t)\big) \int_{t}^{t_1} d\mu(s) ,  x(t) \Big\rangle \, dt + \int_{t_0}^{t_1} [x(t)]^T \, d\mu(t)\\
  & & \!\! + \big\langle h_{x_1}^T\big(x_*(t_0),x_*(t_1)\big) \nu + \alpha r_{x_1}\big(x_*(t_0),x_*(t_1)\big)
           + \lambda_0 S_{x_1}\big(x_*(t_0),x_*(t_1)\big) , x(t_1) \big\rangle \\
  & & \!\! + \big\langle h_{x_0}^T\big(x_*(t_0),x_*(t_1)\big) \nu + \alpha  r_{x_0}\big(x_*(t_0),x_*(t_1)\big)
           + \lambda_0 S_{x_0}\big(x_*(t_0),x_*(t_1)\big), x(t_0) \big\rangle \\
  & & \!\! - \Big\langle \int_{t_0}^{t_1} d\mu(t) , x(t_0) \Big\rangle + \sum_{j=1}^l \int_{t_0}^{t_1} \big\langle g_{j,x}\big(t,x_*(t)\big),x(t) \big\rangle \,d\mu_j(t)
\end{eqnarray*}
überführt werden.
Die eindeutige Darstellung eines stetigen linearen Funktionals im Raum $C_0([t_0,t_1],\R^n)$ liefert für
$p(t)=\displaystyle \int_t^{t_1} \, d\mu(s)$ die Gleichungen
\begin{eqnarray*}
p(t) &=& -h_{x_1}^T\big(x_*(t_0),x_*(t_1)\big) \nu - \alpha r_{x_1}\big(x_*(t_0),x_*(t_1)\big) - \lambda_0 S_{x_1}\big(x_*(t_0),x_*(t_1)\big) \\
     & & + \int_t^{t_1} \big[ \varphi^T_x\big(s,x_*(s),u_*(s)\big)p(s)-\lambda_0 f_x\big(s,x_*(s),u_*(s)\big)\big]\, ds \\
     & & - \sum_{j=1}^l \int_t^{t_1} g_{j,x}\big(s,x_*(s)\big) \,d\mu_j(s), \\
p(t_0) &=& h_{x_0}^T\big(x_*(t_0),x_*(t_1)\big) \nu + \alpha  r_{x_0}\big(x_*(t_0),x_*(t_1)\big) + \lambda_0 S_{x_0}\big(x_*(t_0),x_*(t_1)\big).
\end{eqnarray*}
Damit sind die Optimalitätsbedingungen (\ref{SatzPMPGZFZA1}) und (\ref{SatzPMPGZFZA2}) gezeigt. \\[2mm]
Aufgrund von (\ref{SatzGZFLMR2}) gilt f"ur alle $u(\cdot) \in L_\infty([t_0,t_1],U)$ die Ungleichung
\begin{eqnarray*} 
\lefteqn{\int_{t_0}^{t_1} \big[
    \lambda_0 f\big(t,x_*(t),u_*(t)\big) - \big\langle p(t) , \varphi\big(t,x_*(t),u_*(t)\big) \big\rangle \big] \, dt} \\
&\leq & \int_{t_0}^{t_1} \big[
    \lambda_0 f\big(t,x_*(t),u(t)\big) - \big\langle p(t) , \varphi\big(t,x_*(t),u(t)\big) \big\rangle \big] \, dt,
\end{eqnarray*}
woraus sich mit Hilfe der Eigenschaft Lebesguescher Punkte die Maximumbedingung (\ref{SatzPMPGZFZA3}) ableitet.
Abschließend ergibt sich (\ref{SatzPMPGZFZA4}) unmittelbar aus (\ref{SatzGZFLMR3}). \hfill $\blacksquare$


\begin{theorem}[Pontrjaginsches Maximumprinzip] \label{SatzPMPGZF}
\index{Pontrjaginsches Maximumprinzip!gemischt@-- gemischtes Zielfunktional} 
Es sei $\big(x_*(\cdot),u_*(\cdot)\big) \in \mathscr{B}^{\,\mathcal{S}}_{\rm adm} \cap \mathscr{B}^{\,\mathcal{S}}_{\rm Lip}$. 
Ist $\big(x_*(\cdot),u_*(\cdot)\big)$ ein starkes lokales Minimum der Aufgabe (\ref{GZF1})--(\ref{GZF4}),
dann existieren nicht gleichzeitig verschwindende Multiplikatoren $\lambda_0 \geq 0$, $\alpha \geq0$,
$p(\cdot) \in W^1_\infty([t_0,t_1],\R^n)$ und $\nu \in \R^s$ derart, dass
\begin{enumerate}
\item[(a)] die adjungierte Gleichung
           \index{adjungierte Gleichung!gemischt@-- gemischtes Zielfunktional}
           \begin{equation}\label{SatzPMPGZF1}
           \dot{p}(t)=-H^{\mathcal{S}}_x\big(t,x_*(t),u_*(t),p(t),\lambda_0\big),
           \end{equation}
\item[(b)] die Transversalit"atsbedingungen
           \index{Transversalitätsbedingungen!gemischt@-- gemischtes Zielfunktional}
           \begin{equation}\label{SatzPMPGZF2}
           \left. \begin{array}{l}
\hspace*{-9mm} p(t_0) \,=\, h_{x_0}^T\big(x_*(t_0),x_*(t_1)\big) \nu + \alpha r_{x_0}\big(x_*(t_0),x_*(t_1)\big) + \lambda_0 S_{x_0}\big(x_*(t_0),x_*(t_1)\big), \\[2mm]
\hspace*{-9mm} p(t_1) \,=\, - h_{x_1}^T\big(x_*(t_0),x_*(t_1)\big) \nu - \alpha r_{x_1}\big(x_*(t_0),x_*(t_1)\big) - \lambda_0 S_{x_1}\big(x_*(t_0),x_*(t_1)\big),
           \end{array}\right\}
           \end{equation}
\item[(c)] in fast allen Punkten $t \in [t_0,t_1]$ die Maximumbedingung
           \index{Maximumbedingung!gemischt@-- gemischtes Zielfunktional}
           \begin{equation}\label{SatzPMPGZF3}
           H^{\mathcal{S}}\big(t,x_*(t),u_*(t),p(t),\lambda_0\big) = \max_{u \in U} H^{\mathcal{S}}\big(t,x_*(t),u,p(t),\lambda_0\big)
           \end{equation}
\item[(d)] und die komplementäre Schlupfbedingung
           \begin{equation}\label{SatzPMPGZF4}
           \alpha r\big(x_*(t_0),x_*(t_1)\big)=0
           \end{equation}           
\end{enumerate}
erfüllt sind.
\end{theorem}

\newpage
Weiterhin lassen sich hinreichende Bedingungen nach Arrow
\index{hinreichende Bedingungen, Arrow!gemischt@-- gemischtes Zielfunktional}
direkt anfügen.
Wir übernehmen die im Abschnitt \ref{AbschnittArrowZB} getroffenen Bezeichnungen und formulieren die hinreichenden Bedingungen
bei gemischtem Zielfunktional:

\begin{theorem} \label{SatzHBGZFZA}
In der Aufgabe (\ref{GZF1})--(\ref{GZF5}) mit den festen Randbedingungen $x(t_0)=x_0$ und $x(t_1)=x_1$
sei $\big(x_*(\cdot),u_*(\cdot)\big) \in \mathscr{B}^{\,\mathcal{S}}_{\rm Lip} \cap \mathscr{B}^{\,\mathcal{S}}_{\rm adm}$.
Au"serdem sei $p(\cdot):[t_0,t_1] \to \R^n$ st"uckweise stetig,
besitze h"ochstens endlich viele Sprungstellen (d.\,h. die einseitigen Grenzwerte existieren) $s_k \in (t_0,t_1)$
und sei zwischen diesen Spr"ungen stetig differenzierbar.
Ferner:
\begin{enumerate}
\item[(a)] Das Tripel $\big(x_*(\cdot),u_*(\cdot),p(\cdot)\big)$
           erf"ullt (\ref{SatzPMPGZFZA1})--(\ref{SatzPMPGZFZA3}) in Theorem \ref{SatzPMPGZFZA} mit $\lambda_0=1$.        
\item[(b)] F"ur jedes $t \in [t_0,t_1]$ ist die Funktion $\mathscr{H}^{\mathcal{S}}\big(t,x,p(t)\big)$ konkav 
           und es sind die Funktionen $g_j(t,x)$, $j=1,...,l$, konvex bez"uglich $x$ auf $V^{\mathcal{S}}_\gamma(t)$.
\item[(c)] Die Funktion $S(x_0,x_1)$ ist auf $V^{\mathcal{S}}_\gamma(t_0) \times V^{\mathcal{S}}_\gamma(t_1)$ konvex.
\end{enumerate}
Dann ist $\big(x_*(\cdot),u_*(\cdot)\big)$ ein starkes lokales Minimum der Aufgabe (\ref{GZF1})--(\ref{GZF5}).
\end{theorem}

{\bf Beweis} Die Beweise der Arrow-Bedingungen in den Abschnitten \ref{AbschnittHBPMP} und \ref{AbschnittArrowZB} liefern
zusammen mit der Konvexität der Funktion $S$:
\begin{eqnarray*}
    \Delta
&=& \tilde{J}\big(x(\cdot),u(\cdot)\big) - \tilde{J}\big(x_*(\cdot),u_*(\cdot)\big) \\
&\geq& \langle p(t_1),x(t_1)-x_*(t_1)\rangle-\langle p(t_0),x(t_0)-x_*(t_0)\rangle \\
& & + \langle S_{x_0}\big(x_*(t_0),x_*(t_1),x(t_0)-x_*(t_0)\rangle + \langle S_{x_1}\big(x_*(t_0),x_*(t_1),x(t_1)-x_*(t_1)\rangle \\
&=& \langle p(t_1)+S_{x_1}\big(x_*(t_0),x_*(t_1),x(t_1)-x_*(t_1)\rangle \\
& & -\langle p(t_0)-S_{x_0}\big(x_*(t_0),x_*(t_1),x(t_0)-x_*(t_0)\rangle.
\end{eqnarray*}
Zusammen mit den Transversalitätsbedingungen,
die auf die gleiche Weise wie im Abschnitt \ref{AbschnittArrowZB} um die Anteile der Zustandsbeschränkungen bereinigt sind,
ergibt sich $\Delta \geq 0$ in den Fällen fester und freier Randwerte. \hfill $\blacksquare$ \\[2mm]
In der Aufgabe (\ref{GZF1})--(\ref{GZF5}) mit freier Zeit seien neben der Randbedingungen au"serdem
im Zielfunktional das Terminalfunktional eine Funktion der Anfangs- und Endzeit:
$$S=S(t_0,t_1,x_0,x_1), \quad h= h(t_0,t_1,x_0,x_1), \quad r=r(t_0,t_1,x_0,x_1).$$
Die Methode der Substitution der Zeit in Abschnitt \ref{AbschnittFreieZeit} führt dabei im Zielfunktional zu
$S\big(t(0),t(1),y(0),y(1)\big)$, sowie in den Randbedingungen zu $h\big(t(0),t(1),y(0),y(1)\big)$ und $r\big(t(0),t(1),y(0),y(1)\big)$. 
Dem weiteren Ablauf in Abschnitt \ref{AbschnittFreieZeit} folgend ergibt sich:

\begin{theorem}[Pontrjaginsches Maximumprinzip] \label{SatzGZFFZ}
\index{Pontrjaginsches Maximumprinzip!gemischt@-- gemischtes Zielfunktional}
In der Aufgabe (\ref{GZF1})--(\ref{GZF5}) mit freier Zeit sei
$\big([t_{0*},t_{1*}],x_*(\cdot),u_*(\cdot)\big) \in \mathscr{B}^{\,\mathcal{F}}_{\rm adm} \cap \mathscr{B}^{\,\mathcal{F}}_{\rm Lip}$.
Ist $\big([t_{0*},t_{1*}],x_*(\cdot),u_*(\cdot)\big)$ ein starkes lokales Minimum in der Aufgabe mit freier Zeit,
dann existieren Zahlen $\lambda_0 \geq 0$ und $\alpha \geq 0$, ein $\nu \in \R^s$, eine Vektorfunktion $p(\cdot):[t_{0*},t_{1*}] \to \R^n$
und auf den Mengen
$$T_j=\big\{t \in [t_{0*},t_{1*}] \,\big|\, g_j\big(t,x_*(t)\big)=0\big\}, \quad j=1,...,l,$$
konzentrierte nichtnegative regul"are Borelsche Ma"se $\mu_j$ endlicher Totalvariation
(wobei s"amtliche Gr"o"sen nicht gleichzeitig verschwinden) derart,
dass die Vektorfunktion $p(\cdot)$ von beschr"ankter Variation und rechtsseitig stetig, und
\begin{enumerate}
\item[(a)] die adjungierte Gleichung
           \index{adjungierte Gleichung!gemischt@-- gemischtes Zielfunktional}
           \begin{eqnarray}
           p(t)&=& -h^T_{x_1}\big(t_{0*},t_{1*},x_*(t_{0*}),x_*(t_{1*})\big) \nu - \alpha r_{x_1}\big(t_{0*},t_{1*},x_*(t_{0*}),x_*(t_{1*})\big) \nonumber \\
               & & - \lambda_0 S_{x_1}\big(t_{0*},t_{1*},x_*(t_{0*}),x_*(t_{1*})\big)
                   + \int_t^{t_{1*}} H^{\mathcal{S}}_x\big(s,x_*(s),u_*(s),p(s),\lambda_0\big) \, ds \nonumber \\
           \label{SatzGZFFZ1} & & -\sum_{j=1}^l \int_t^{t_{1*}} g_{j,x}\big(s,x_*(s)\big)\, d\mu_j(s),
           \end{eqnarray}
\item[(b)] die Transversalit"atsbedingungen
           \index{Transversalitätsbedingungen!gemischt@-- gemischtes Zielfunktional}
           \begin{equation} \label{SatzGZFFZ2}
           \left. \begin{array}{l}
           p(t_{0*}) \,=\, h^T_{x_0}\big(t_{0*},t_{1*},x_*(t_{0*}),x_*(t_{1*})\big) \nu + \lambda_0 S_{x_0}\big(t_{0*},t_{1*},x_*(t_{0*}),x_*(t_{1*})\big) \\[2mm]
                     \hspace*{10mm} + \alpha r_{x_0}\big(t_{0*},t_{1*},x_*(t_{0*}),x_*(t_{1*})\big), \\[2mm]
           p(t_{1*}^-)-p(t_{1*}) \,=\, \displaystyle - \sum_{j=1}^l \mu_j(\{t_{1*}\}) \, g_{j,x}\big(t_{1*},x_*(t_{1*})\big)
           \end{array} \right\}
           \end{equation}
\item[(c)] in fast allen Punkten $t \in [t_0,t_1]$ die Maximumbedingung
           \index{Maximumbedingung!gemischt@-- gemischtes Zielfunktional}
           \begin{equation}\label{SatzGZFFZ3}
           H^{\mathcal{S}}\big(t,x_*(t),u_*(t),p(t),\lambda_0\big) = \max_{u \in U} H^{\mathcal{S}}\big(t,x_*(t),u,p(t),\lambda_0\big),
           \end{equation}
\item[(d)] die komplementäre Schlupfbedingung
           \begin{equation}
           \alpha r\big(t_{0*},t_{1*},x_*(t_{0*}),x_*(t_{1*})\big) = 0,
           \end{equation}
\item[(e)] die Hamilton-Funktion eine Funktion beschr"ankter Variation ist und
           \begin{eqnarray}
           && \hspace*{-20mm} \mathscr{H}^{\mathcal{S}}\big(t,x_*(t),p(t),\lambda_0\big) 
              = \big\langle h_{t_1}\big(t_{0*},t_{1*},x_*(t_{0*}),x_*(t_{1*})\big) , \nu \big\rangle \nonumber \\
           && \hspace*{-10mm} + \lambda_0 S_{t_1}\big(t_{0*},t_{1*},x_*(t_{0*}),x_*(t_{1*})\big)+ \alpha r_{t_1}\big(t_{0*},t_{1*},x_*(t_{0*}),x_*(t_{1*})\big) \nonumber \\
           && \hspace*{-10mm} - \int_t^{t_{1*}} H^{\mathcal{S}}_t\big(s,x_*(s),u_*(s),p(s),\lambda_0\big) \, ds
                              + \sum_{j=1}^l \int_t^{t_{1*}} g_{j,t}\big(s,x_*(s)\big) \, d\mu_j(s), \\
           && \hspace*{-20mm}  \mathscr{H}^{\mathcal{S}}\big(t_{0*},x_*(t_{0*}),p(t_{0*}),\lambda_0\big)
              = - \big\langle h_{t_0}\big(t_{0*},t_{1*},x_*(t_{0*}),x_*(t_{1*})\big) , \nu \big\rangle  \nonumber \\[1mm]
           && \hspace*{-10mm} -\alpha r_{t_0}\big(t_{0*},t_{1*},x_*(t_{0*}),x_*(t_{1*})\big)- \lambda_0 S_{t_0}\big(t_{0*},t_{1*},x_*(t_{0*}),x_*(t_{1*})\big)
           \end{eqnarray}
\end{enumerate}
erfüllt sind.
\end{theorem}

%% file: 3-12-Produktion.tex
\subsubsection{Investition in Produktionskapazitäten}\index{Investition in Produktionskapazitäten}
Ein Investor hat die Möglichkeit eine Unternehmung aufzubauen.
Die Produktion schafft Werte und Güter in Höhe von $f\big(x(t)\big)$ Einheiten wenn die Produktionsrate $x(t)$ beträgt.
Beim Aufbau der Unternehmung und der Investition in die Produktionsstruktur entstehen Kosten $c$ pro Einheit für die
Anschaffungen,
die am Ende der Planungsperiode zum Preis $d$ wieder veräu"sert werden können.
Über den gesamten Planungszeitraum kann ein Teil $\big(1-u(t)\big)f\big(x(t)\big)$ der Einnahmen konsumiert und der verbliebene
Anteil $u(t)f\big(x(t)\big)$ in die weitere Produktion mit Effizienz $\alpha$ eingesetzt werden. \\[2mm]
Es entsteht ein interessantes Beispiel (Seierstad \cite{Seierstad}, S.\,187),
in dem der Anfangszustand $x(0)$ frei angepasst werden kann,
um eine adäquate Startsituation zu schaffen:
\begin{eqnarray}
&&\label{Produktion1} J\big(x(\cdot),u(\cdot)\big) = \int_0^T \big(1-u(t)\big)f\big(x(t)\big)e^{-rt} \, dt
                      -c x(0) + dx(T)e^{-rT} \to \sup, \hspace*{5mm} \\
&&\label{Produktion2} \dot{x}(t) = \alpha u(t)f\big(x(t)\big), \qquad x(0) \mbox{ frei }, \quad x(T) \mbox{ frei},
                    \qquad u(t) \in [0,1].
\end{eqnarray}
Dabei sei $f$ eine zweimal stetig differenzierbare und streng konkave Funktion mit $f'(x) \to 0$ für $x \to \infty$.
Ferner gelte $d<1/\alpha<c$.
In der Aufgabe liegen keine Randbedingungen vor,
so dass die Transversalitätsbedingungen in den notwendigen Optimalitätsbedingungen einzig durch die
Terminalfunktionale festgelegt werden. \\[2mm]
Wir wenden auf die Aufgabe (\ref{Produktion1})--(\ref{Produktion2}) mit freien Randbedingungen 
das Pontrjaginsche Maximumprinzip in der Form von Theorem \ref{SatzPMPGZF} an: 
Die Pontrjagin-Funktion $H^{\mathcal{S}}$ der Aufgabe lautet im normalen Fall
$$H^{\mathcal{S}}(t,x,u,p,1)=p \alpha u f(x) +(1-u)f(x)e^{-rt}=[p\alpha u+(1-u)e^{-rt}]f(x).$$
Die Optimalitätsbedingungen in Theorem \ref{SatzPMPGZF} liefern für die Adjungierte
$$\dot{p}(t)=-\big[p(t)\alpha u_*(t)+\big(1-u_*(t)\big)e^{-rt}\big]f'\big(x_*(t)\big),
  \qquad p(0)=c, \quad p(T)=de^{-rT},$$
und mittels der Maximumbedingung an eine optimale Steuerung die Kriterien
$$u_*(t)= 1 \quad\mbox{für } \alpha p(t) > e^{-rt}, \qquad u_*(t)= 0 \quad\mbox{für }\alpha p(t) < e^{-rt}.$$
Weiter ergibt sich $\dot{p}(t)=-\max\{e^{-rt},\alpha p(t)\} \cdot f'\big(x_*(t)\big) < 0$
und $p(t)$ fällt streng monoton vom Wert $p(0)=c$ zum Wert $p(T)=de^{-rT}$.
Zu dem gelten $\alpha p(0)=\alpha c >1$ und $\alpha p(T)=\alpha de^{-rT} < e^{-rT}$.
Daher existiert ein eindeutiger Schaltzeitpunkt $\tau \in (0,T)$,
in welchem von Investition in Konsumption gewechselt wird. \\[2mm]
Zusammenfassend ergibt sich $u_*(t)=1$ für $t \in [0,\tau)$ und $u_*(t)=0$ für $t \in [\tau,T]$ als die optimale
Investitionsstrategie.
Die zugehörige Trajektorie $x_*(t)$ erhält man als die Lösung von $\dot{x}_*(t) = \alpha f\big(x_*(t)\big)$ zum 
Anfangswert $x_*(0)$ für $t \in [0,\tau)$ und ist konstant $x_*(t)=x_*(\tau)$ für $t \in [\tau,T]$.
Die Parameter $x_*(0)$ und $\tau$ sind dabei so zu bestimmen, dass die Transversalitätsbedingungen $p(0)=c$, $p(T)=de^{-rT}$
und ferner $\alpha p(\tau)=  e^{-r\tau}$ erfüllt sind. \\
Die Aufgabe der optimalen Investition in Produktionskapazitäten genügt
den Anforderungen der Arrow-Bedingungen in Theorem \ref{SatzHBGZFZA}.
Damit ist der Steuerungsprozess $\big(x_*(\cdot),u_*(\cdot)\big)$ ein starkes lokales Minimum der Aufgabe. \hfill $\square$

%% file: 3-13-Instandhaltung.tex
\subsubsection{Instandhaltungsmanagement} \index{Instandhaltung}
Vorbeugende Instandhaltung wirkt sich positiv auf Lebensdauer und Funktionalität von Produktionsanlagen aus,
woraus sich au"serdem Einsparungen beim Ersatz von Anlagen ergeben.
Die Aufgabe besteht hier in der Ermittlung eines optimierten Zusammenspiels von Instandhaltungsintensität im Laufe der
Produktionsdauer und dem Verkaufswert bzw. dem Schrotterlös der Anlage am Ende der Produktionsdauer. \\
Das nachstehende Modell wurde durch die Betrachtungen in Feichtinger \cite{Feichtinger} angeregt. \\[2mm]
Es sei $x(t)$ der Zustand einer Maschine,
die im Produktionseinsatz den Erlös $E \cdot x(t)$ pro Zeiteinheit einbringt.
Der Einsatz der Maschine führt zum Verschleiß mit der Rate $\delta \cdot x(t)$,
welchem durch instandhaltende Maßnahmen entgegengewirkt werden kann.
Die Kosten für die Instandhaltung wird durch $u(t)$ angegeben und die Effizienz der Maßnahmen durch $C\big(u(t)\big)$.
Zusammenfassend bezeichnen:
\begin{center}\begin{tabular}{rcl}
$x(t)$ &--& den Zustand der Maschine zur Zeit $t$, \\
$E \cdot x(t)$ &--& den Erlös der Maschine pro Zeiteinheit aus dem Betrieb der Maschine, \\
$\delta \cdot x(t)$ &--& den Verschleiß der Maschine bei Produktion, \\
$W \cdot x(T)$ &--& den Wiederverkaufswert am Ende der Planungsperiode, \\
$u(t)$ &--& die Instandshaltungskosten, \\
$C(u)$ &--& die Effizienz der Instandhaltung zur Rate $u$.
\end{tabular}\end{center}
Der zu erwartende Nettoerlös aus dem Betrieb der Maschine, den anfallenden Instandhaltungskosten und dem Barwert des Verkaufes
berechnet sich gemäß
$$\int_0^T e^{-\varrho t}\big[E x(t)-u(t)\big] \, dt + e^{-\varrho T} \cdot W x(T).$$
Die instandhaltenden Maßnahmen $C(u)$ seien mit wachsendem Kostenaufwand $u$ weniger effizient.
Daher besitze die Funktion $C$ die Eigenschaften: 
$$C(0)=0, \quad C(\infty)=\infty, \quad C'(u) >0, \quad C'(0)=\infty, \quad C'(\infty)=0, \quad C''(u)<0, \quad u \geq 0.$$
Diese Eigenschaften spiegelt zum Beispiel die Funktion $C(u)=u^\sigma$ mit $\sigma \in (0,1)$ wider.
Über dem Planungszeitraum werden außerdem die instandhaltenden Maßnahmen durch den betrieblichen Verschleiß und durch wiederkehrende Reparaturen weniger effektiv.
Wir machen für diese Beobachtung den Ansatz $e^{-\alpha t} \cdot C(u)$. \\[2mm]  
Zusammenfassend ergibt sich für unser Instandhaltungsmodell die Aufgabe
\begin{equation} \label{Instand} 
\left.\begin{array}{l}
\hspace*{-5mm}
\displaystyle J\big(x(\cdot),u(\cdot)\big) = \int_0^T e^{-\varrho t}\big[E x(t)-u(t)\big] \, dt + e^{-\varrho T} \cdot W x(T) \to \sup,  \\[4mm]
\hspace*{-5mm}   
\displaystyle \dot{x}(t)=e^{-\alpha t} \cdot C\big(u(t)\big)- \delta x(t), \qquad x(0)=x_0>0,  \\[4mm]
\hspace*{-5mm}   
u \geq 0, \qquad \varrho,\,\delta,\, \alpha,\, E,\, W >0, \qquad\displaystyle W < \frac{E}{\varrho + \delta}.
\end{array} \right\}
\end{equation}
Wir wenden Theorem \ref{SatzPMPGZF} an:
Die Pontrjagin-Funktion $H^{\mathcal{S}}$ der Aufgabe (\ref{Instand}) hat die Form
$$H^{\mathcal{S}}(t,x,u,p)=p (e^{-\alpha t}C(u)-\delta x) + e^{-\varrho t}(E x -u).$$
Sie ist linear in $x$, womit die Hamilton-Funktion $\mathscr{H}^{\mathcal{S}}(t,x,p)$ konkav ist.
Daher sind die Bedingungen des Maximumprinzips hinreichend für ein starkes lokales Maximum. 
Die Maximumbedingung (\ref{PMPeinfach6}) lautet
$$H^{\mathcal{S}}\big(t,x_*(t),u_*(t),p(t)\big)=\max_{u \geq 0} \Big[ p(t) \big(e^{-\alpha t}C(u)-\delta x_*(t)\big) + e^{-\varrho t}\big(E x_*(t) -u\big)\Big].$$
Aufgrund der Eigenschaften der Funktion $C(u)$ können nur positive Instandhaltungskosten $u_*(t)$ optimal sein.
Deswegen führt die Maximumbedingung nach Ausschluss der Randlösung $u=0$ auf die Gleichung
$$H^{\mathcal{S}}_u\big(t,x_*(t),u_*(t),p(t)\big)=p(t)e^{-\alpha t}C'\big(u_*(t)\big) - e^{-\varrho t} =0
  \Leftrightarrow p(t)= e^{-(\varrho - \alpha) t} \frac{1}{C'\big(u_*(t)\big)}.$$
Wir beachten den Vorzeichenwechsel beim Übergang zu einem Minimierungsproblem in der Transversalitätsbedingung.
Damit erhalten (\ref{SatzPMPGZF1}) und (\ref{SatzPMPGZF2}) die Form
$$\dot{p}(t)=\delta p(t)-E \cdot e^{-\varrho t}, \quad p(T)=W \cdot  e^{-\varrho T}.$$
Für die adjungierte Funktion $p(\cdot)$ erhalten wir die Abbildung
$$p(t)=W e^{-(\varrho + \delta)T} e^{\delta t} + \frac{E}{\varrho + \delta}\big[e^{-\varrho t} - e^{-(\varrho + \delta)T} e^{\delta t}\big].$$
Den Ausdruck in der letzten Klammer formen wir um und erhalten
$$e^{-\varrho t} - e^{-(\varrho + \delta)T} e^{\delta t} = e^{-\varrho t} - e^{-\varrho T} \cdot e^{-\delta (T-t)} > e^{-\varrho t} - e^{-\varrho T} >0$$
für $t \in [0,T)$.
Daher nimmt die Adjungierte $p(t)$ nur positive Werte über $[0,T]$ an und
die optimale Strategie $u_*(\cdot)$ ist durch $p(t)= e^{-(\varrho - \alpha) t} / C'\big(u_*(t)\big)$ sinnvoll festgelegt. \\[2mm]
In $t=0$ ergibt sich für die Adjungierte $p(\cdot)$ der Wert
$$p(0)= W e^{-(\varrho + \delta)T} + \frac{E}{\varrho + \delta}\big[1 - e^{-(\varrho + \delta)T}\big] \approx \frac{E}{\varrho + \delta}.$$
Im Sinn des Schattenpreises bemisst $p(\cdot)$ zu Beginn den gesamten Erlös unter Beachtung von Diskontierung und Verschleiß,
$$p(0) \approx \int_0^T E e^{-(\varrho + \delta)t} \, dt \approx \int_0^\infty E e^{-(\varrho + \delta)t} \, dt = \frac{E}{\varrho + \delta},$$
und nimmt am Ende des Planungszeitraumes den Barwert $p(T)=W e^{-\varrho T}$ des Verkaufspreises an. \hfill $\square$

%% file: 3-20-Zeitfenster.tex
\subsection{Freier Anfangs- und Endzeitpunkt unter Zeitschranken}
In dem Lehrbuch von Feichtinger \& Hartl \cite{Feichtinger} werden unter der Rubrik
``Optimale Wahl des Endzeitpunktes'' für die Aufgabe mit freiem Anfangs- und Endzeitpunkt
notwendige Optimalitätsbedingungen in dem Fall formuliert,
wenn sich der Zeitraum $[t_0,t_1]$ innerhalb eines gewissen Zeitfensters $[T_0,T_1]$ befinden soll.
Die Argumentation in \cite{Feichtinger} basiert auf einer getrennten Herleitung der bekannten Optimalitätsbedingungen und der anschließenden Analyse
der Zeitpunkte $t_{0*}$ und $t_{1*}$.
Da ein optimaler Steuerungsprozess und damit die Gestalt von Optimalitätsbedingungen wesentlich durch das optimale Zeitintervall $[t_{0*},t_{1*}]$
beeinflusst werden,
binden wir in den nachstehenden Untersuchungen das vorgegebene Zeitfenster unmittelbar ein.
Die Aufgabe mit einem gegebenen Zeitfenster gehört nicht zum Standard in den Darstellungen von Optimalitätsbedingungen.
Deswegen widmen wir diesem Detail diesen eigenständigen Abschnitt. \\[2mm]
Wir kehren zurück zur Standardaufgabe (\ref{FZ1})--(\ref{FZ5})
mit freiem Anfangs- und Endzeitpunkt im Abschnitt \ref{AbschnittFreieZeit}
und fügen dieser Aufgabe gewisse Zeitschranken hinzu:
 
\begin{equation} \label{FZZF1}
\left. \begin{array}{l}
\displaystyle  J\big(x(\cdot),u(\cdot)\big) = \int_{t_0}^{t_1} f\big(t,x(t),u(t)\big) \, dt \to \inf, \\[3mm]
\dot{x}(t) = \varphi\big(t,x(t),u(t)\big), \\[1mm]
h_0\big(t_0,x(t_0)\big)=0, \quad h_1\big(t_1,x(t_1)\big)=0, \\[1mm]
[\Theta_0,\Theta_1] \subseteq [t_0,t_1] \subseteq [T_0,T_1],\quad T_0 \leq \Theta_0 < \Theta_1  \leq T_1, \\[1mm]
u(t) \in U \subseteq \R^m, \quad U \not= \emptyset, \\[1mm]
g_j\big(t,x(t)\big) \leq 0 \quad \mbox{f"ur alle } t \in [t_0,t_1], \quad j=1,...,l.
\end{array}\right\}
\end{equation}
Die Methode der Substitution der Zeit führt auf das Steuerungsproblem
\begin{equation} \label{FZZF2}
\left. \begin{array}{l}
\displaystyle \int_0^1 v \cdot f\big(t(s),y(s),w(s)\big) \, ds \to \inf, \\[3mm]
t'(s) = v, \quad y'(s) = v \cdot \varphi\big(t(s),y(s),w(s)\big), \\[1mm]
h_0\big( t(0),y(0) \big) = 0, \quad h_1\big( t(1),y(1) \big) = 0, \\[1mm]
T_0-t(0) \leq 0, \quad t(1)-T_1 \leq 0, \quad t(0) - \Theta_0 \leq 0, \quad \Theta_1-t(1) \leq 0, \\[1mm]
v > 0, \quad w(s) \in U, \\[1mm]
g_j\big(t(s),y(s)\big) \leq 0, \quad s \in [0,1], \quad j=1,...,l.
\end{array} \right\}
\end{equation}
Da die Zustandsvariable $t(\cdot)$ wegen $t'(s)=v >0$ streng monoton wachsend ist,
bewirken die Randbedingungen $r_0\big(t(0)\big)=T_0-t(0) \leq 0$ und $r_1\big(t(1)) = t(1)-T_1 \leq 0$
die Eingrenzung der Zeit $t(\cdot)$ in das Zeitfenster $[T_0,T_1]$ bzw.
die Randbedingungen $\varrho_0\big(t(0)\big)=t(0) - \Theta_0 \leq 0$ und $\varrho_1\big(t(1)) = \Theta_1 -t(1) \leq 0$
die Berücksichtigung der Zeitdauer $[\Theta_0,\Theta_1]$. \\
Die Randbedingungen $r_0(t_0)=T_0-t_0,r_1(t_1)=t_1-T_1, \varrho_0(t_0)=t_0-\Theta_0,\varrho_1(t_1)=\Theta_1-t_1$ haben eine denkbar einfache Gestalt
und besitzen die einfachen Ableitungen $r_0'=\varrho_1'=-1$ und $r_1'=\varrho_0'=1$.
Deswegen fließen lediglich die Lagrangeschen Multiplikatoren $\alpha_0,\alpha_1$ bezüglich $r_0,r_1$ bzw.
$\beta_0,\beta_1$ bezüglich $\varrho_0,\varrho_1$ in die Optimalitätsbedingungen ein.

\newpage
Die Aufgabe (\ref{FZZF2}) gliedert sich in den Rahmen der Aufgabe (\ref{GZF1})--(\ref{GZF5}) mit allgemeinen Randbedingungen ein.
Daher existieren nach Theorem \ref{SatzPMPGZFZA} Zahlen $\lambda_0 \geq 0$, $\alpha_0 \geq 0$, $\alpha_1 \geq 0$, $\beta_0 \geq 0$ und $\beta_1 \geq 0$,
Vektoren $l_0 \in \R^{s_0},l_1 \in \R^{s_1},$ eine Vektorfunktion $\tilde{p}(\cdot)$,
eine Funktion $\tilde{q}(\cdot)$, sowie nichtnegative, regul"are und
auf den Mengen $\Delta_j$ konzentrierte Borelsche Ma"se
$\tilde{\mu}_j, j=1,...,l,$ (diese Gr"o"sen verschwinden nicht gleichzeitig) derart, dass die Beziehungen
\begin{eqnarray*} 
\tilde{p}(s) &=& - h^T_{1,x_1}\big(t_*(1),y_*(1)\big) l_1
                 +\int_s^1 \tilde{H}^{\mathcal{S}}_x\big(t_*(\sigma),y_*(\sigma),w_*(\sigma),v_*,\tilde{p}(\sigma),\tilde{q}(\sigma),\lambda_0\big) \, d\sigma \\
             & & - \sum_{j=1}^l \int_s^1 g_{j,x}\big(t_*(\sigma),y_*(\sigma)\big)\, d\tilde{\mu}_j(\sigma), \\
\tilde{p}(0) &=& h^T_{0,x_0} \big(t_*(0),y_*(0)\big) l_0, \\
\tilde{q}(s) &=& - \big\langle h_{1,t_1}\big(t_*(1),y_*(1)\big),l_1 \big\rangle
                 +\int_s^1 \tilde{H}^{\mathcal{S}}_t\big(t_*(\sigma),y_*(\sigma),w_*(\sigma),v_*,\tilde{p}(\sigma),\tilde{q}(\sigma),\lambda_0\big) \, d\sigma \\
             & & - \sum_{j=1}^l \int_s^1 g_{j,t}\big( t_*(\sigma),y_*(\sigma)\big) \, d\tilde{\mu}_j(\sigma) - \alpha_1 + \beta_1, \\
\tilde{q}(0) &=& \big\langle h_{0,t_0} \big(t_*(0),y_*(0)\big) , l_0 \big\rangle - \alpha_0 + \beta_0,
\end{eqnarray*}
für fast alle $s \in [0,1]$ die Maximumbedingung
$$\tilde{H}^{\mathcal{S}}\big(t_*(s),y_*(s),w_*(s),v_*,\tilde{p}(s),\tilde{q}(s),\lambda_0\big)
  = \max_{u \in U,\, v> 0} v \cdot \big[ H^{\mathcal{S}}\big(t_*(s),y_*(s),u,\tilde{p}(s),\lambda_0\big) + \tilde{q}(s) \big]$$
und die komplementären Schlupfbedingungen
$$\alpha_0 \big(T_0-t(0)\big)=0, \quad \alpha_1 \big(t(1)-T_1\big)=0, \quad \beta_0 \big(t(0)-\Theta_0\big)=0, \quad \beta_1 \big(\Theta_1 - t(1)\big)=0$$
erfüllt sind. \\[2mm]
Die Randbedingungen $r_0,...,\varrho_1$ wirken sich lediglich bezüglich der Zeitvariablen und damit nur auf die Adjungierte $\tilde{q}(\cdot)$
in Form der Multiplikatoren $\alpha_0,...,\beta_1$ aus.
Dieser Einfluss wird durch die Beziehung $H^{\mathcal{S}}\big(t_*(s),y_*(s),u,\tilde{p}(s),\lambda_0\big) = -\tilde{q}(s)$
an die Hamilton-Funktion
$\mathscr{H}^{\mathcal{S}}(t,x,p,\lambda_0) = \max\limits_{u \in U} H^{\mathcal{S}}(t,x,u,p,\lambda_0)$
übergeben. \\
Für die Auswirkung der Zeitschranken auf die ursprüngliche Aufgabe (\ref{FZZF1}) bedeutet dies,
dass in den notwendigen Optimalitätsbedingungen
nur die zeitabhängigen Elemente durch das Auftreten der Multiplikatoren $\alpha_0,..., \beta_1$ beeinflusst werden. \\[2mm]
Die Umkehrung der Substitution der Zeit erfolgt wieder auf die gleiche Weise wie im Abschnitt \ref{AbschnittFreieZeit}.
Zusammenfassend lautet damit das Pontrjaginsche Maximumprinzip wir für die Aufgabe mit freiem Anfangs- und Endzeitpunkt unter Zeitschranken wie folgt:

\newpage
\begin{theorem}[Pontrjaginsches Maximumprinzip]
\index{Pontrjaginsches Maximumprinzip!Standard@-- Standardaufgabe} 
In der Aufgabe (\ref{FZZF1}) sei der Steuerungsprozess
$\big([t_{0*},t_{1*}],x_*(\cdot),u_*(\cdot)\big) \in \mathscr{B}^{\,\mathcal{F}}_{\rm adm} \cap \mathscr{B}^{\,\mathcal{F}}_{\rm Lip}$.
Ist $\big([t_{0*},t_{1*}],x_*(\cdot),u_*(\cdot)\big)$ ein starkes lokales Minimum der Aufgabe (\ref{FZZF1}),
dann existieren Zahlen $\lambda_0 \geq 0$, $\alpha_0 \geq 0$, $\alpha_1 \geq 0$, $\beta_0 \geq 0$ und $\beta_1 \geq 0$,
Vektoren $l_0 \in \R^{s_0}$ und $l_1 \in \R^{s_1}$, eine Vektorfunktion $p(\cdot):[t_{0*},t_{1*}] \to \R^n$
und auf den Mengen
$$T_j=\big\{t \in [t_{0*},t_{1*}] \,\big|\, g_j\big(t,x_*(t)\big)=0\big\}, \quad j=1,...,l,$$
konzentrierte nichtnegative regul"are Borelsche Ma"se $\mu_j$ endlicher Totalvariation
(wobei s"amtliche Gr"o"sen nicht gleichzeitig verschwinden) derart,
dass die Vektorfunktion $p(\cdot)$ von beschr"ankter Variation und rechtsseitig stetig ist, und
\begin{enumerate}
\item[(a)] die adjungierte Gleichung
           \index{adjungierte Gleichung!Standard@-- Standardaufgabe}
           \begin{eqnarray}
           p(t)&=&-h^T_{1,x_1}\big(t_{1*},x_*(t_{1*})\big) l_1 
           + \int_t^{t_{1*}} H^{\mathcal{S}}_x\big(s,x_*(s),u_*(s),p(s),\lambda_0\big) \, ds \nonumber \\
           & & -\sum_{j=1}^l \int_t^{t_{1*}} g_{j,x}\big(s,x_*(s)\big)\, d\mu_j(s),
           \end{eqnarray}
\item[(b)] die Transversalit"atsbedingungen
           \index{Transversalitätsbedingungen!Standard@-- Standardaufgabe}
           \begin{equation}
            \left. \begin{array}{lcl} p(t_{0*}) &=& h^T_{0,x_0} \big(t_{0*},x_*(t_{0*})\big) l_0, \\[2mm]
                  p(t_{1*}^-)-p(t_{1*}) &=& \displaystyle - \sum_{j=1}^l \mu_j(\{t_{1*}\}) \, g_{j,x}\big(t_{1*},x_*(t_{1*})\big),
                             \end{array} \right\}
           \end{equation} 
\item[(c)] in fast allen Punkten $t\in [t_{0*},t_{1*}]$ die Maximumbedingung 
           \index{Maximumbedingung!Standard@-- Standardaufgabe}
           \begin{equation}
           H^{\mathcal{S}}\big(t,x_*(t),u_*(t),p(t),\lambda_0\big) = \max_{u \in U} H^{\mathcal{S}}\big(t,x_*(t),u,p(t),\lambda_0\big)
           \end{equation}
\item[(d)] die komplementären Schlupfbedingungen
           \begin{equation}
           \alpha_0 \big(t_{0*}-T_0\big)=0, \quad \alpha_1 \big(t_{1*}-T_1\big)=0, \quad
           \beta_0 \big(t_{0*}-\Theta_0\big)=0, \quad \beta_1 \big(t_{1*}-\Theta_1\big)=0
           \end{equation}
\end{enumerate}
erfüllt sind und weiterhin
\begin{enumerate}
\item[(e)] die Beziehungen 
           \begin{eqnarray}
           && \hspace*{-10mm} \mathscr{H}^{\mathcal{S}} \big(t,x_*(t),p(t),\lambda_0\big) 
              = \big\langle h_{1,t_1}\big(t_{1*},x_*(t_{1*})\big) , l_1 \big\rangle
               - \int_t^{t_{1*}} H^{\mathcal{S}}_t\big(s,x_*(s),u_*(s),p(s),\lambda_0\big) \, ds \nonumber \\
           && \hspace*{30mm} + \sum_{j=1}^l \int_t^{t_{1*}} g_{j,t}\big(s,x_*(s)\big) \, d\mu_j(s) + \alpha_1 - \beta_1, \\
           && \hspace*{-10mm}  \mathscr{H}^{\mathcal{S}}\big(t_{0*},x_*(t_{0*}),p(t_{0*}),\lambda_0\big)
              = - \big\langle h_{0,t_0}\big(t_{0*},x_*(t_{0*})\big) , l_0 \big\rangle + \alpha_0 - \beta_0
           \end{eqnarray}
\end{enumerate}
gelten.
\end{theorem}

\newpage
\begin{theorem} \label{SatzFZZF}
\index{Pontrjaginsches Maximumprinzip!Standard@-- Standardaufgabe} 
In der Aufgabe (\ref{FZZF1}) sei 
$\big([t_{0*},t_{1*}],x_*(\cdot),u_*(\cdot)\big) \in \mathscr{B}^{\,\mathcal{F}}_{\rm adm} \cap \mathscr{B}^{\,\mathcal{F}}_{\rm Lip}$.
Ist der Steuerungsprozess $\big([t_{0*},t_{1*}],x_*(\cdot),u_*(\cdot)\big)$ ein starkes lokales Minimum der Aufgabe (\ref{FZZF1}),
dann existieren nicht gleichzeitig verschwindende Multiplikatoren
$\lambda_0 \geq 0$, $\alpha_0 \geq 0$, $\alpha_1 \geq 0$, $\beta_0 \geq 0$, $\beta_1 \geq 0$,
$l_0 \in \R^{s_0}$, $l_1 \in \R^{s_1}$ und $p(\cdot) \in W^1_\infty([t_{0*},t_{1*}],\R^n)$
derart, dass
\begin{enumerate}
\item[(a)] für fast alle $t \in [t_{0*},t_{1*}]$ die adjungierte Gleichung
           \index{adjungierte Gleichung!Standard@-- Standardaufgabe}
           \begin{equation} \label{SatzFZZF1}
           \dot{p}(t) = -H^{\mathcal{S}}_x\big(t,x_*(t),u_*(t),p(t),\lambda_0\big),
           \end{equation}
\item[(b)] die Transversalit"atsbedingungen
           \index{Transversalitätsbedingungen!Standard@-- Standardaufgabe}
           \begin{equation}\label{SatzFZZF2}
           p(t_{0*})= h^T_{0,x_0} \big(t_{0*},x_*(t_{0*})\big) l_0, \qquad p(t_{1*})=-h^T_{1,x_1}\big(t_{1*},x_*(t_{1*})\big) l_1,
           \end{equation}
\item[(c)] in fast allen Punkten $t\in [t_{0*},t_{1*}]$ die Maximumbedingung 
           \index{Maximumbedingung!Standard@-- Standardaufgabe}
           \begin{equation}
           H^{\mathcal{S}}\big(t,x_*(t),u_*(t),p(t),\lambda_0\big) = \max_{u \in U} H^{\mathcal{S}}\big(t,x_*(t),u,p(t),\lambda_0\big),
           \end{equation}
\item[(d)] die komplementären Schlupfbedingungen
           \begin{equation}
           \alpha_0 \big(t_{0*}-T_0\big)=0, \qquad \alpha_1 \big(t_{1*}-T_1\big)=0, \quad
           \beta_0 \big(t_{0*}-\Theta_0\big)=0, \quad \beta_1 \big(t_{1*}-\Theta_1\big)=0
           \end{equation}
\item[(e)] und f"ur die Funktion $t \to \mathscr{H}^{\mathcal{S}}\big(t,x_*(t),p(t),\lambda_0\big)$ die Bedingungen
          \begin{eqnarray}
          \frac{d}{dt} \mathscr{H}^{\mathcal{S}}\big(t,x_*(t),p(t),\lambda_0\big) &=& H^{\mathcal{S}}_t\big(t,x_*(t),u_*(t),p(t),\lambda_0\big), \\
          \mathscr{H}^{\mathcal{S}}\big(t_{0*},x_*(t_{0*}),p(t_{0*}),\lambda_0\big)
             &=& - \big\langle h_{0,t_0} \big(t_{0*},x_*(t_{0*})\big) , l_0 \big\rangle + \alpha_0 - \beta_0, \\
          \mathscr{H}^{\mathcal{S}}\big(t_{1*},x_*(t_{1*}),p(t_{1*}),\lambda_0\big)
             &=& \big\langle h_{1,t_1} \big(t_{1*},x_*(t_{1*})\big) , l_1 \big\rangle + \alpha_1-\beta_1
          \end{eqnarray}
\end{enumerate}
erfüllt sind.
\end{theorem}

\begin{beispiel}
{\rm Wir modifizieren das Beispiel \ref{BeispielMassenpunkt} zur Bewegung eines Massenpunktes indem wir
der zeitminimalen Bewegung den Beschleunigungsaufwand gegenüber stellen:
\begin{equation} \label{MassenpunktZF1}
\left. \begin{array}{l}
\displaystyle J\big(s(\cdot),v(\cdot),a(\cdot)\big) = \int_0^T \bigg(1 + \frac{1}{2}a^2(t) \bigg) \, dt \to \inf, \\[3mm]
\dot{s}(t) = v(t), \qquad \dot{v}(t) = a(t), \\[2mm]
s(0)=0, \quad s(T)=s_T>0, \quad v(0)=v(T)=0, \quad a(t) \in \R.
\end{array} \right\}
\end{equation}
F"ur diese Aufgabe lautet die Pontrjagin-Funktion
$$H^{\mathcal{S}}(t,s,v,a,p_1,p_2,\lambda_0) = p_1 v + p_2 a - \lambda_0 \bigg(1+\frac{1}{2}a^2\bigg).$$
Die notwendigen Bedingungen in Theorem \ref{SatzFZZF} liefern f"ur die Adjungierten wieder
\begin{eqnarray*}
\dot{p}_1(t) \equiv 0 \quad&\Rightarrow&\quad p_1(t)\equiv \pi_1 \in \R, \\
\dot{p}_2(t) = -p_1(t)=-\pi_1 \quad&\Rightarrow&\quad p_2(t) = \pi_2 -\pi_1 t,  \; \pi_2 \in \R.
\end{eqnarray*}
Mit Hilfe der Maximumbedingung ergibt sich $\lambda_0=1$ und wir erhalten für $a_*(\cdot)$:
$$\max_{a \in \R} \bigg( p_2(t) a - \frac{1}{2} a^2 \bigg) \quad\Rightarrow\quad  a_*(t)=p_2(t).$$
Weiterhin liegt eine autonome Aufgabe vor und es ergibt sich
\begin{eqnarray*}
0 &\equiv& \mathscr{H}^{\mathcal{S}}\big(t,s_*(t),v_*(t),p_1(t),p_2(t),\lambda_0\big) = p_1(t) v_*(t) + p_2(t) a_*(t)-1 - \frac{1}{2}a_*^2(t) \\
  &=& \pi_1 \bigg(\pi_2 t - \frac{1}{2} \pi_1 t^2\bigg) +(\pi_2 - \pi_1 t)^2 - 1 - \frac{1}{2}(\pi_2- \pi_1 t)^2
      =\frac{1}{2} \pi_2^2 -1.
\end{eqnarray*}
Damit gilt $\pi_2^2=2$.
Der Parameter $\pi_1$ und die optimale Stoppzeit $T_*$ ergeben sich aus
$$v_*(T_*) = \pi_2 T_*- \frac{1}{2} \pi_1 T_*^2=0, \qquad s_*(T_*) = \frac{1}{2}\pi_2 T_*^2 - \frac{1}{6} \pi_1 T_*^3 = S_T.$$
Es sei nun $0 <T_1 < T_*$ und wir fügen der Aufgabe \ref{MassenpunktZF1} das Zeitfenster $[0,T_1]$ hinzu.
Dadurch entsteht die Aufgabe
\begin{equation} \label{MassenpunktZF2}
\left. \begin{array}{l}
\displaystyle J\big(s(\cdot),v(\cdot),a(\cdot)\big) = \int_{t_0}^{t_1} \bigg(1 + \frac{1}{2}a^2(t) \bigg) \, dt \to \inf, \\[3mm]
\dot{s}(t) = v(t), \qquad \dot{v}(t) = a(t), \quad [t_0,t_1] \subseteq [0,T_1], \\[2mm]
s(t_0)=0, \quad s(t_1)=s_T>0, \quad v(t_0)=v(t_1)=0, \quad a(t) \in \R.
\end{array} \right\}
\end{equation}
Während die adjungierten Gleichungen und die Maximumbedingung mit denen der Aufgabe \ref{MassenpunktZF1} übereinstimmen,
liegen nun die Randbedingungen 
$$\mathscr{H}^{\mathcal{S}}\big(t_{0*},x_*(t_{0*}),p(t_{0*}),\lambda_0\big)= \alpha_0, \quad
  \mathscr{H}^{\mathcal{S}}\big(t_{1*},x_*(t_{1*}),p(t_{1*}),\lambda_0\big) = \alpha_1$$
an die Hamilton-Funktion vor.
Wegen $H^{\mathcal{S}}_t\big(t,x_*(t),u_*(t),p(t),\lambda_0\big) \equiv 0$ muss darin $\alpha_0 = \alpha_1$ gelten
und weiterhin sind beide Zahlen positiv, da $T_1 < T_*$ vorausgesetzt wurde.
Da $\alpha_0 >0$ und $\alpha_1 >0$ ausfallen, gelten $t_{0*}=0$ und $t_{1*}=T_1$. 
Ferner ergibt sich
$$\alpha_0 \equiv \mathscr{H}^{\mathcal{S}}\big(t,s_*(t),v_*(t),p_1(t),p_2(t),\lambda_0\big) =\frac{1}{2} \pi_2^2 -1,$$
also $\pi_2^2=2(\alpha_0+1)$.
Anhand der Beziehungen
$$v_*(T_1) = \pi_2 T_1- \frac{1}{2} \pi_1 T_1^2=0, \qquad s_*(T_1) = \frac{1}{2}\pi_2 T_1^2 - \frac{1}{6} \pi_1 T_1^3 = S_T$$
können die Größen $\pi_1=12 S_T /T_1^3$ und $\pi_2 = 6 S_T / T_1^2$,
z.\,B. mit Hilfe der Cramerschen Regel, ermittelt werden. \hfill $\square$}
\end{beispiel}

%% file: 3-30-Multiprozesse.tex
\subsection{Multiprozesse} \label{KapitelMultiprozess}
Der Begriff des Multiprozesses\index{Multiprozess} l"asst sich am Beispiel des Autofahrens illustrieren:
Neben den kontinuierlichen Steuerungen ``beschleunigen'' und ``bremsen''
bildet die Wahl des konkreten Ganges einen wesentlichen Beitrag zur Minimierung des Treibstoffverbrauchs oder
zum Erreichen des Zielortes in k"urzester Zeit.
Denn die Auswahl des jeweiligen Ganges beeinflusst ma"sgeblich das dynamische Beschleunigungsverhalten, den momentanen Benzinverbrauch
und die Geschwindigkeit.
Das Fahrverhalten und der Treibstoffverbrauch wird in jedem einzelnen Gang durch eine eigene Dynamik und Verbrauchsfunktion beschrieben.
Dementsprechend setzt sich dieses Optimierungsproblem aus verschiedenen,
dem jeweilig ausgew"ahlten Gang zugeordneten, Steuerungsproblemen zusammen.
Die Schaltfolge zwischen den einzelnen G"angen wird durch eine Wechselstrategie\index{Wechselstrategie}
\index{Multiprozess!Wechselstrategie@--, Wechselstrategie} beschrieben.
Das Beispiel des Autofahrens verdeutlicht dabei den speziellen Charakter der Wechselstrategie:
Im Vergleich zu den Pedalen, die stufenlos gesteuert werden k"onnen,
ist die Auswahl des Ganges eine rein diskrete Gr"o"se. \\[2mm]
Ein Multiprozess besteht aus einer gewissen endlichen Anzahl von einzelnen Steuerungssystemen,
welche sich jeweils aus individuellen Dynamiken, Zielfunktionalen und Steuerungsbereichen zusammensetzen.
Neben der Suche nach der optimalen Steuerung f"ur das jeweils gew"ahlte Steuerungssystem
liegt das Hauptaugenmerk bei der Untersuchung von Multiprozessen auf der Bestimmung der optimalen Wechselstrategie.
Dabei wird eine Wechselstrategie durch die Anzahl der Wechsel zwischen den einzelnen Steuerungssystemen,
durch die konkrete Auswahl des jeweiligen Steuerungssystems und durch diejenigen Zeitpunkte,
zu denen diese Wechsel stattfinden, beschrieben.
Wir betrachten ausschlie"slich Multiprozesse mit stetigen Zust"anden.
Au"serdem d"urfen keine Wechselkosten anfallen. \\[2mm]
Eine anschauliche Beschreibung eines Multiprozesses wird durch die Auflistung der jeweilig ausgew"ahlten Steuerungssysteme gegeben.
Dazu zerlegt man das Zeitintervall $[t_0,t_1]$ zu den Wechselzeitpunkten in Teilintervalle,
$$t_0=s_0<s_1<...<s_N=t_1,$$
und beschreibt f"ur $j=1,...,N$ auf den Teilabschnitten das Steuerungsproblem:
$$\int_{s_{j-1}}^{s_j} f_{i_j}\big(t,x(t),u(t)\big) \, dt, \quad
  \dot{x}(t) = \varphi_{i_j}\big(t,x(t),u(t)\big) \mbox{ f"ur } t \in (s_{j-1},s_j),\quad u(t) \in U_{i_j} \subseteq \R^{m_{i_j}}.$$
Dabei geh"ort die Zahl $i_j$ zu der gegebenen Menge $\{1,...,k\}$ und gibt das ausgew"ahlte Steuerungssystem an.
Der Multiprozess erh"alt dadurch die Gestalt
\begin{eqnarray*}
&& J\big(x(\cdot),u(\cdot)\big) = \sum_{j=1}^N \int_{s_{j-1}}^{s_j} f_{i_j}\big(t,x(t),u(t)\big) \, dt \to \inf, \\
&& \dot{x}(t) = \varphi_{i_j}\big(t,x(t),u(t)\big) \quad \mbox{ f"ur } t \in (s_{j-1},t_j), \quad j=1,...,N, \\
&& u(t) \in U_{i_j} \subseteq \R^{m_{i_j}},\quad {i_j} \in \{1,...,k\}, \quad j=1,...,N.
\end{eqnarray*}
In dieser Darstellung eines Multiprozesses ist anzumerken,
dass durch die Folge der Indizes $\{i_j\}$ die Wechselstrategie vordefiniert ist.
D.\,h., dass die Anzahl der Wechsel und die jeweilige Auswahl des Steuerungssystems implizit bereits determiniert werden.
Die Beschreibung einer Wechselstrategie als ``echte'' Steuerungsvariable wird in dieser Form nicht erreicht. \\
Allerdings basiert auf dieser Formulierung weitestgehend die Strategie zur Untersuchung von allgemeinen Multiprozessen.
N"amlich auf dem Vergleich von Multiprozessen, die auf der gleichen Wechselstrategie beruhen.
Die fundamentale Arbeit von Sussmann \cite{Sussmann} baut auf diesem Ansatz auf.
Jedoch bezieht sie sich auf allgemeinere Aufgaben als wir sie betrachten werden,
da sie ebenso unstetige Zust"ande und Wechselkosten beinhaltet.
Weitere grundlegenden Beitr"age zu notwendigen Optimalit"atsbedingungen f"ur Multiprozesse, z.\,B. bei
Clarke \& Vinter \cite{ClarkeVinter,ClarkeVinterII}, Dmitruk \& Kaganovich \cite{Dmitruk},
Galbraith \& Vinter \cite{Galbraith}, Garavello \& Piccoli \cite{Garavello}, Shaikh \& Caines \cite{Shaikh04,Shaikh07},
basieren ebenso auf dieser Philosophie. \\[2mm]
Wir werden im Folgenden f"ur Multiprozesse einen Weg aufzeigen,
die Wechselstrategien als ``echte'' Steuerungsvariable einzuf"uhren.
Damit sind wir in der Lage,
Multiprozesse auf starke lokale Optimalstellen zu untersuchen.
Das beinhaltet,
dass wir zur Bestimmung von notwendigen Optimalit"atsbedingungen Multiprozesse mit beliebigen
Wechseltstrategien vergleichen (vgl. Tauchnitz \cite{TauchnitzDiss}). \\
Unser grundlegendes Element zur Beschreibung von Wechselstrategien sind gewisse Zerlegungen des Zeitintervalls und
die Verkn"upfung dieser Zerlegungen mit vektorwertigen Abbildungen.
Die Wechselstrategien in Form von Zerlegungen dienen als Steuerungsvariable in der Beschreibung eines Multiprozesses
und flie"sen so in die konkrete Aufgabestellung und in das Pontrjaginsche Maximumprinzip ein.

%% file: 3-31-Zerlegungen.tex
\subsubsection{$k$-fache Zerlegungen und die Aufgabenstellung} \label{AbschnittZerlegungen}
In diesem Abschnitt f"uhren wir auf der Grundlage einfacher "Uberlegungen die $k$-fachen Zerlegungen eines Intervalls
ein.\index{kfacheZerlegung@$k$-fache Zerlegung}\index{Multiprozess!kfacheZerlegung@--, $k$-fache Zerlegung}
Anschlie"send legen wir dar,
wie diese Gr"o"sen eine Wechselstrategie zwischen einer gegebenen Anzahl von Abbildungen realisiert. \\[2mm]
Es seien das Intervall $I \subseteq \R$ und $k \in \N$ gegeben.
\begin{definition} \label{DefinitionZerlegung} 
Unter einer $k$-fachen Zerlegung des Intervalls $I \subseteq \R$ verstehen wir ein endliches System
$\mathscr{A} = \{\mathscr{A}_1,...,\mathscr{A}_k\}$ von Lebesgue-messbaren Teilmengen von $I$ mit
$$\bigcup_{1\leq s\leq k} \mathscr{A}_s = I, \qquad
  \mathscr{A}_{s} \cap \mathscr{A}_{s'} = \emptyset \quad\mbox{f"ur } s \not=s'.$$
Es bezeichnet $\mathscr{Z}^k(I) = \{ \mathscr{A} \} = \big\{\{ \mathscr{A}_s\}_{1\leq s\leq k}\big\}$ die Menge der
$k$-fachen Zerlegungen von $I$.
\end{definition}

Wir identifizieren ein Element $\mathscr{A} \in \mathscr{Z}^k(I)$ durch die Vektorfunktion
$$\chi_{\mathscr{A}}(t) = \big(\chi_{\mathscr{A}_1}(t),...,\chi_{\mathscr{A}_k}(t) \big), \quad t \in I,$$
der charakteristischen Funktionen der Mengen ${\mathscr A}_s$.
Jede Funktion $\chi_{\mathscr{A}}(\cdot)$ geh"ort der Menge
$$\mathscr{Y}^k(I) = \bigg\{ y(\cdot) \in L_\infty(I,\R^k)\,\bigg|\, y(t)=\big(y_1(t),...,y_k(t)\big), y_s(t) \in \{0,1\},
                             \sum_{s=1}^k y_s(t)=1\bigg\}$$
an.
Die Menge $\mathscr{Z}^k(I)$ der $k$-fachen Zerlegungen bzw. $\mathscr{Y}^k(I)$ der charakteristischen Vektorfunktionen
repr"asentieren bei geeigneter Verkn"upfung s"amtliche Wechselstrategien in Multiprozessen mit $k$ Steuerungssystemen. \\[2mm]
Die Grundlage bilden daf"ur die folgenden Bezeichnungen: \\
Wir betrachten $k$ Funktionen $h_s(t,x,u_s): I \times \R^n \times \R^{m_s} \to \R^m$ und fassen diese zur
Vektorfunktion $h=(h_1,...,h_k)$ zusammen.
Dabei setzen wir $u=(u_1,...,u_k) \in \R^{m_1+...+m_k}$.
Die Verkn"upfung der Funktion $h$ mit einer $k$-fachen Zerlegung legen wir wie folgt fest:
$$\mathscr{A} \circ h(t,x,u) = \chi_{\mathscr{A}}(t) \circ h(t,x,u)
                             = \sum_{s=1}^k \chi_{\mathscr{A}_s}(t) \cdot h_s(t,x,u_s).$$
Sind die Funktionen $h^{i}$ differenzierbar, dann setzen wir f"ur die partiellen Ableitungen:
\begin{eqnarray*}
\mathscr{A} \circ h_t(t,x,u) &=& \chi_{\mathscr{A}}(t) \circ h_t(t,x,u) = \sum_{s=1}^k \chi_{\mathscr{A}_s}(t) \cdot h_{s,t}(t,x,u_s), \\
\mathscr{A} \circ h_x(t,x,u) &=& \chi_{\mathscr{A}}(t) \circ h_x(t,x,u) = \sum_{s=1}^k \chi_{\mathscr{A}_s}(t) \cdot h_{s,x}(t,x,u_s).
\end{eqnarray*}

In diesen Bezeichnungen w"ahlt eine $k$-fache Zerlegung zum Zeitpunkt $t \in I$ einerseits in eindeutiger Weise eine der
Funktionen $h_s(t,x,u_s)$ und
au"serdem die entsprechende Komponente $u_s$ des Vektors $u=(u_1,...,u_k)$ aus. \\[2mm]
Es sei $k \in \N$ die Anzahl der verschiedenen gegebenen Steuerungssysteme,
die individuelle Integranden $f_s(t,x,u_s): \R \times \R^n \times \R^{m_s} \to \R$,
Dynamiken $\varphi_s(t,x,u_s): \R \times \R^n \times \R^{m_s} \to \R^n$
und Steuerbereiche $U_s \subseteq \R^{m_s}$ besitzen k"onnen.
Ferner seien
$$h_0(x_0):\R^n \to \R^{s_0}, \quad h_1(x_1):\R^n \to \R^{s_1}, \qquad g_j(t,x):\R \times \R^n \to \R, \quad j=1,...,l.$$
Wir fassen die Steuerungssysteme durch die Setzungen
$$f=(f_1,...,f_k),\qquad \varphi=(\varphi_1,...,\varphi_k), \qquad U=U_1 \times ... \times U_k$$
zusammen.
Mit den Vorbetrachtungen des letzten Abschnitts formulieren wir "uber $[t_0,t_1]$
die Aufgabe des optimalen Multiprozesses mit Zustandsbeschr"ankungen wie folgt:
\begin{eqnarray}
&& \label{HA1} J\big(x(\cdot),u(\cdot),\mathscr{A}\big)
               = \int_{t_0}^{t_1} \chi_{\mathscr{A}}(t) \circ f\big(t,x(t),u(t)\big) dt \to \inf, \\
&& \label{HA2} \dot{x}(t) = \chi_{\mathscr{A}}(t) \circ \varphi\big(t,x(t),u(t) \big), \\
&& \label{HA3} h_0\big( x(t_0) \big) = 0, \quad h_1\big( x(t_1) \big) = 0, \\
&& \label{HA4} u(t) \in U= U_1 \times ... \times U_k, \quad U_s \not= \emptyset, \quad \mathscr{A} \in \Zt, \\
&& \label{HA5} g_j\big(t,x(t)\big) \leq 0, \quad t \in [t_0,t_1], \quad j=1,...,l.
\end{eqnarray}
Die Aufgabe (\ref{HA1})--(\ref{HA5}) mit Zustandsbeschr"ankungen betrachten wir bez"uglich
$$\big(x(\cdot),u(\cdot),\mathscr{A}\big) \in W^1_\infty([t_0,t_1],\R^n) \times L_\infty([t_0,t_1],U) \times \Zt.$$
Mit $\mathscr{B}^{\,\mathcal{M}}_{\rm Lip}$ bezeichnen wir die Menge aller Tripel $\big(x(\cdot),u(\cdot),\mathscr{A}\big)$,
für die es ein $\gamma>0$ derart gibt, dass die Abbildungen
$f_s(t,x,u_s)$, $\varphi_s(t,x,u_s)$, $h_i(x_i)$ und $g_j(t,x)$ auf der Menge aller
$(t,x,x_0,x_1,u) \in \R \times \R^n \times \R^n \times \R^n \times \R^m$ mit
$$t_0 \leq t\leq t_1, \quad \|x-x(t)\| < \gamma, \quad \|x_0-x(t_0)\| < \gamma, \quad \|x_1-x(t_1)\| < \gamma, \quad u \in \R^m$$
stetig in der Gesamtheit der Variablen und stetig differenzierbar bez"uglich $x$ sind. \\[2mm]
Das Tripel $\big(x(\cdot),u(\cdot),\mathscr{A}\big) \in W^1_\infty([t_0,t_1],\R^n) \times L_\infty([t_0,t_1],U) \times \Zt$
hei"st ein zul"assiger Multiprozess in der Aufgabe (\ref{HA1})--(\ref{HA5}),
falls $\big(x(\cdot),u(\cdot),\mathscr{A}\big)$ dem System (\ref{HA2}) gen"ugt, die Randbedingungen (\ref{HA3})
und die Zustandsbeschr"ankungen (\ref{HA5}) erf"ullt.
Die Menge $\mathscr{B}^{\,\mathcal{M}}_{\rm adm}$ bezeichnet die Menge der zul"assigen Multiprozesse. \\[2mm]
Ein zul"assiger Multiprozess $\big(x_*(\cdot),u_*(\cdot),\mathscr{A}_*\big)$ ist eine starke lokale
Minimalstelle\index{Minimum, starkes lokales!Multi@-- Multiprozesse}
der Aufgabe (\ref{HA1})--(\ref{HA5}),
falls eine Zahl $\varepsilon > 0$ derart existiert, dass die Ungleichung 
$$J\big(x(\cdot),u(\cdot),\mathscr{A}\big) \geq J\big(x_*(\cdot),u_*(\cdot),\mathscr{A}_*\big)$$
f"ur alle $\big(x(\cdot),u(\cdot),\mathscr{A}\big) \in \mathscr{B}^{\,\mathcal{M}}_{\rm adm}$ mit 
$\|x(\cdot)-x_*(\cdot) \|_\infty < \varepsilon$ gilt. \\[2mm]
In der Aufgabe (\ref{HA1})--(\ref{HA5})
definieren wir f"ur die einzelnen Steuerungssysteme die partiellen Pontrjaginschen Funktionen
$H_s:\R \times \R^n \times \R^{m_s} \times \R^n \times \R \to \R$
gem"a"s
$$H_s(t,x,u_s,p,\lambda_0) = \langle p , \varphi_s(t,x,u_s) \rangle - \lambda_0 f_s(t,x,u_s).$$
Wir fassen die Pontrjaginschen Funktionen $H_s$ zu der Vektorfunktion $H^{\mathcal{M}}=(H_1,...,H_k)$ zusammen.
Dann hei"st f"ur die Aufgabe optimaler Multiprozesse die Funktion
\begin{eqnarray*}
\mathscr{A} \circ H^{\mathcal{M}}&:&\R \times \R^n \times \R^{m_1 + ... + m_k} \times \R^n \times \R \to \R, \\
\mathscr{A} \circ H^{\mathcal{M}}(t,x,u,p,\lambda_0\big)
&=& \chi_{\mathscr{A}}(t) \circ H^{\mathcal{M}}(t,x,u,p,\lambda_0\big) = \sum_{s=1}^k \chi_{\mathscr{A}_s}(t) \cdot H_s(t,x,u_s,p,\lambda_0) \\
&=& \sum_{s=1}^k \chi_{\mathscr{A}_s}(t) \cdot \big( \langle p , \varphi_s(t,x,u_s) \rangle - \lambda_0 f_s(t,x,u_s) \big)
\end{eqnarray*}
die Pontrjagin-Funktion der Aufgabe optimaler Multiprozesse.

%% file: 3-32-PMPeinfach.tex
\subsubsection{Die elementare Aufgabe mit freiem Endpunkt} \label{AbschnittPMPeinfachMP}
Für die Aufgabe (\ref{HA1})--(\ref{HA5}) untersuchen wir zunächst den Spezialfall im Rahmen der stückweise stetigen Steuerungen und
stückweise stetig differenzierbaren Zustände.
Dafür bezeichne $\mathscr{Z}^k_0([t_0,t_1]) \subset \Zt$ die Menge aller $k$-fachen Zerlegungen,
in denen sich die Mengen $\mathscr{A}_1,...,\mathscr{A}_k$ durch endliche Vereingungen von rechtsseitig halboffenenen oder abgeschlossenen
Teilintervallen von $[t_0,t_1]$ ergeben.
Ganz entsprechend ergibt sich für die charakteristischen Vektorfunktionen $\chi_{\mathscr{A}}(\cdot)$ die Menge
$$\mathscr{Y}^k(I) = \bigg\{ y(\cdot) \in PC([t_0,t_1],\R)\,\bigg|\, y(t)=\big(y_1(t),...,y_k(t)\big), y_s(t) \in \{0,1\},
                             \sum_{s=1}^k y_s(t)=1\bigg\}.$$
Die Aufgabe mit freiem Endpunkt ist das Steuerungsproblem 
\begin{eqnarray}
&&\label{PMPeinfach1MP} J\big(x(\cdot),u(\cdot),\mathscr{A}\big)
               = \int_{t_0}^{t_1} \chi_{\mathscr{A}}(t) \circ f\big(t,x(t),u(t)\big) dt \to \inf, \\
&&\label{PMPeinfach2MP} \dot{x}(t) = \chi_{\mathscr{A}}(t) \circ \varphi\big(t,x(t),u(t) \big), \quad x(t_0)=x_0, \\
&&\label{PMPeinfach3MP} u(t) \in U= U_1 \times ... \times U_k, \quad U_s \not= \emptyset, \quad \mathscr{A} \in \mathscr{Z}^k_0([t_0,t_1])
\end{eqnarray}
auf dem gegebenen Intervall $[t_0,t_1]$ und mit festem Punkt $x_0 \in \R^n$. \\[2mm]
Die Aufgabe (\ref{PMPeinfach1MP})--(\ref{PMPeinfach3MP}) untersuchen wir bez"uglich der Tripel
$$\big(x(\cdot),u(\cdot),\mathscr{A}\big) \in PC_1([t_0,t_1],\R^n) \times PC([t_0,t_1],U) \times \mathscr{Z}^k_0([t_0,t_1]).$$
Mit $\mathscr{D}^{\,\mathcal{M}}_{\rm Lip}$ bezeichnen wir die Menge aller Tripel $\big(x(\cdot),u(\cdot),\mathscr{A}\big)$,
für die es ein $\gamma>0$ derart gibt,
dass die Abbildungen $f_s(t,x,u_s)$, $\varphi_s(t,x,u_s)$ für $s=1,...,k$ auf der Menge aller Punkte $(t,x,u) \in \R \times \R^n \times \R^m$ mit 
$$t \in [t_0,t_1], \qquad \|x-x(t)\| < \gamma, \qquad u \in \R^m$$
stetig in der Gesamtheit aller Variablen und stetig differenzierbar bezüglich $x$ sind. \\[2mm]
Ein Tripel $\big(x(\cdot),u(\cdot),\mathscr{A}\big) \in PC_1([t_0,t_1],\R^n) \times PC([t_0,t_1],U) \times \mathscr{Z}^k_0([t_0,t_1])$
ist ein zul"assiger Steuerungsprozess der Aufgabe (\ref{PMPeinfach1MP})--(\ref{PMPeinfach3MP}),
falls $\big(x(\cdot),u(\cdot),\mathscr{A}\big)$ dem System (\ref{PMPeinfach2MP}) zu $x(t_0)=x_0$ gen"ugt.
Mit $\mathscr{D}^{\,\mathcal{M}}_{\rm adm}$ bezeichnen wir die Menge der zul"assigen Steuerungsprozesse. \\[2mm]
Ein zul"assiger Steuerungsprozess $\big(x_*(\cdot),u_*(\cdot),\mathscr{A}_*\big)$ ist eine
starke lokale Minimalstelle\index{Minimum, starkes lokales!elementar@-- elementare Aufgabe}
der Aufgabe (\ref{PMPeinfach1MP})--(\ref{PMPeinfach3MP}),
falls eine Zahl $\varepsilon > 0$ derart existiert, dass die Ungleichung 
$$J\big(x(\cdot),u(\cdot),\mathscr{A}\big) \geq J\big(x_*(\cdot),u_*(\cdot),\mathscr{A}_*\big)$$
f"ur alle $\big(x(\cdot),u(\cdot),\mathscr{A}\big) \in \mathscr{D}^{\,\mathcal{M}}_{\rm adm}$ mit $\|x(\cdot)-x_*(\cdot)\|_\infty < \varepsilon$ gilt. \\[2mm]
Weiterhin bezeichnet $\mathscr{A} \circ H^{\mathcal{M}}$ die Pontrjagin-Funktion
$$\mathscr{A} \circ H^{\mathcal{M}}(t,x,u,p,\lambda_0\big) =  \chi_{\mathscr{A}}(t) \circ H^{\mathcal{M}}(t,x,u,p,\lambda_0\big)
  = \sum_{s=1}^k \chi_{\mathscr{A}_s}(t) \cdot H_s(t,x,u_s,p,\lambda_0).$$

\begin{theorem}[Pontrjaginsches Maximumprinzip] \label{SatzPMPeinfachMP}
\index{Pontrjaginsches Maximumprinzip!Multi@-- Multiprozesse} 
\index{Pontrjaginsches Maximumprinzip!elementar@-- elementare Aufgabe} 
Sei $\big(x_*(\cdot),u_*(\cdot),\mathscr{A}_*\big) \in \mathscr{D}^{\,\mathcal{M}}_{\rm adm} \cap \mathscr{D}^{\,\mathcal{M}}_{\rm Lip}$. 
Ist $\big(x_*(\cdot),u_*(\cdot),\mathscr{A}_*\big)$ ein starkes lokales Minimum der Aufgabe (\ref{PMPeinfach1MP})--(\ref{PMPeinfach3MP}),
dann existiert eine Vektorfunktion $p(\cdot) \in PC_1([t_0,t_1],\R^n)$ derart, dass
\begin{enumerate}
\item[(a)] die adjungierte Gleichung
           \index{adjungierte Gleichung!elementar@-- elementare Aufgabe}
           \index{adjungierte Gleichung!Multi@-- Multiprozesse}
           \begin{equation}\label{PMPeinfach4MP} 
           \dot{p}(t) = - \chi_{\mathscr{A}_*}(t) \circ H^{\mathcal{M}}_x\big(t,x_*(t),u_*(t),p(t),1\big),
           \end{equation}
\item[(b)] die Transversalitätsbedingung
           \index{Transversalitätsbedingungen!elementar@-- elementare Aufgabe}
           \index{Transversalitätsbedingungen!Multi@-- Multiprozesse}
           \begin{equation}\label{PMPeinfach5MP} p(t_1)=0 \end{equation}   
\item[(c)] und in fast allen Punkten $t \in [t_0,t_1]$ die Maximumbedingung           
           \index{Maximumbedingung!elementar@-- elementare Aufgabe}
           \index{Maximumbedingung!Multi@-- Multiprozesse}
           \begin{equation}\label{PMPeinfach6MP} 
           \chi_{\mathscr{A}_*}(t) \circ H^{\mathcal{M}}\big(t,x_*(t),u_*(t),p(t),1\big)
           = \max_{\substack{u_s \in U_s\\ 1\leq s\leq k}} H^{\mathcal{M}}_s\big(t,x_*(t),u_s,p(t),1\big)
           \end{equation}
\end{enumerate}
erfüllt sind.
\end{theorem}

{\bf Beweis} 
Genauso wie im Abschnitt \ref{AbschnittPMPBeweiseinfach} ergibt sich wieder,
dass die Lösung $p(\cdot)$ der adjungierten Gleichung (\ref{PMPeinfach4MP}) zur Transversalitätsbedingung (\ref{PMPeinfach5MP})
im Fall $\lambda_0=0$ im Widerspruch zur Nichttrivialität der Multiplikatoren über $[t_0,t_1]$ identisch verschwinden würde.
Also gilt $\lambda_0 >0$ und Theorem \ref{SatzPMPeinfachMP} gilt in der normalen Form mit $\lambda_0=1$. \\[1mm]
Wir passen den Beweis für die elementare Standardaufgabe in Abschnitt \ref{AbschnittPMPBeweiseinfach}
an die zusätzliche Steuerung $\mathscr{A} \in \mathscr{Z}^k_0([t_0,t_1])$ an:
Es sei $\tau \in (t_0,t_1)$ ein Stetigkeitspunkt der Steuerungen $u_*(\cdot)$ und $\chi_{\mathcal{A}_*}(\cdot)$.
Weiterhin seien $v=(v_1,...,v_k) \in U=U_1 \times ... \times U_k$ und $s \in \{1,...,k\}$.
Damit legen wir die einfachen Nadelvariationen $u_\lambda(\cdot)$, $\chi_{\mathcal{A}_\lambda}(\cdot)$ durch \index{Nadelvariation, einfache}
$$u_{\lambda}(t) = 
  \left\{ \begin{array}{ll}
          u_*(t) & \mbox{ f"ur } t \not\in [\tau-\lambda,\tau) \\
          v      & \mbox{ f"ur } t     \in [\tau-\lambda,\tau)
          \end{array} \right. \mbox{ und }
   \left\{ \begin{array}{ll}
           \chi_{\mathcal{A}_\lambda}(t) = \chi_{\mathcal{A}_*}(t) & \mbox{ f"ur } t \not\in [\tau-\lambda,\tau) \\
           \chi_{\mathcal{A}_{s,\lambda}}(t)=1     & \mbox{ f"ur } t     \in [\tau-\lambda,\tau)
          \end{array} \right.$$ 
fest.
Dabei wählt $\chi_{\mathcal{A}_\lambda}(\cdot)$ über $[\tau-\lambda,\tau)$ von den gegebenen Steuerungssystemen
dasjenige zum Index $s \in \{1,...,k\}$ aus.
Weiter sei $x_\lambda(\cdot)$ die eindeutige L"osung der Gleichung
$$\dot{x}_\lambda(t) = \chi_{\mathcal{A}_\lambda}(t) \circ \varphi\big(t,x(t),u_\lambda(t)\big), \quad x(t_0)=x_0.$$
Dann ist $x_{\lambda}(t) = x_*(t)$ f"ur $t_0 \leq t \leq  \tau - \lambda$. \\[2mm]
F"ur $t \geq \tau$ untersuchen wir den Grenzwert
$\displaystyle y(t)=\lim_{\lambda \to 0^+}\frac{x_{\lambda}(t) - x_*(t)}{\lambda}$:
Nach Voraussetzung sind $u_*(\cdot)$ und $\chi_{\mathcal{A}_*}(\cdot)$ f"ur hinreichend kleine $\lambda>0$
auf $[\tau-\lambda,\tau]$ stetig.
Daher gelten
\begin{eqnarray*}
x_*(\tau) &=& x_*(\tau - \lambda)
              + \lambda \cdot \chi_{\mathcal{A}_*}(\tau - \lambda) \circ \varphi\big(\tau - \lambda,x_*(\tau - \lambda),u_*(\tau - \lambda)\big) + o(\lambda), \\
x_{\lambda}(\tau) &=& x_*(\tau - \lambda)
             + \lambda \cdot \chi_{\mathcal{A}_\lambda}(\tau - \lambda) \circ \varphi\big(\tau - \lambda,x_*(\tau - \lambda),u_\lambda(\tau - \lambda)\big) + o(\lambda).
\end{eqnarray*}
Deswegen existiert der Grenzwert
$\displaystyle y(\tau)=\lim_{\lambda \rightarrow 0+0}\frac{x_{\lambda}(\tau) - x_*(\tau)}{\lambda}$
und ist gleich
\begin{equation} \label{BeweisPMPeinfachMP2}
y(\tau)=\chi_{\mathcal{A}_\lambda}(\tau) \circ \varphi\big(\tau,x_*(\tau),u_\lambda(\tau)\big) - \chi_{\mathcal{A}_*}(\tau) \circ \varphi\big(\tau,x_*(\tau),u_*(\tau)\big).
\end{equation}
Über $[\tau,t_1]$ ergibt sich aus den Sätzen \ref{SatzEEglobal} und \ref{SatzDGLDifferenzierbarkeit} "uber die Stetigkeit und Differenzierbarkeit der L"osung
eines Differentialgleichungssystems in Abh"angigkeit von den Anfangsdaten,
dass $y(t)$ existiert und der Integralgleichung
\begin{equation} \label{BeweisPMPeinfachMP1}
y(t)=y(\tau) + \int_{\tau}^{t} \big[\chi_{\mathcal{A}_*}(s) \circ \varphi_x\big(s,x_*(s),u_*(s)\big) \big] \, y(s) \, ds
\end{equation}
genügt.
Mit (\ref{PMPeinfach4MP}),  (\ref{PMPeinfach5MP}) und (\ref{BeweisPMPeinfachMP1}) erhalten wir f"ur $t \geq \tau$ die Beziehung
\begin{equation} \label{BeweisPMPeinfachMP3}
\langle p(\tau) , y(\tau) \rangle
 = - \int_{\tau}^{t_1} \big\langle \chi_{\mathcal{A}_*}(t) \circ f_x\big(t,x_*(t),u_*(t)\big), y(t) \big\rangle dt.
\end{equation}
Da $\big(x_*(\cdot),u_*(\cdot)\big)$ ein starkes lokales Minimum ist, gilt f"ur hinreichend kleine $\lambda >0$
$$\frac{J\big(x_\lambda(\cdot),u_\lambda(\cdot),\mathcal{A}_\lambda\big) - J\big(x_*(\cdot),u_*(\cdot),\mathcal{A}_*\big)}{\lambda} \geq 0.$$
Im Grenzwert $\lambda \to 0^+$ ergibt sich in dieser Ungleichung
\begin{eqnarray*}
0 &\leq& \lim_{\lambda \to 0^+} \frac{J\big(x_\lambda(\cdot),u_\lambda(\cdot),\mathcal{A}_\lambda\big)- J\big(x_*(\cdot),u_*(\cdot),\mathcal{A}_*\big)}{\lambda} \\
  &=& \chi_{\mathcal{A}_\lambda}(\tau) \circ f\big(\tau,x_*(\tau),u_\lambda(\tau)\big) - \chi_{\mathcal{A}_*}(\tau) \circ f\big(\tau,x_*(\tau),u_*(\tau)\big) \\
  & & + \int_{\tau}^{t_1} \big\langle \chi_{\mathcal{A}_*}(t) \circ f_x\big(t,x_*(t),u_*(t)\big) , y(t) \big\rangle \, dt.
\end{eqnarray*}
Bei Anwendung von Gleichung (\ref{BeweisPMPeinfachMP3}) ergibt sich die Ungleichung
\begin{eqnarray*}
\lefteqn{ \big\langle p(\tau) , \chi_{\mathcal{A}_*}(\tau) \circ \varphi\big(\tau,x_*(\tau),u_*(\tau)\big) \big\rangle
                               - \chi_{\mathcal{A}_*}(\tau) \circ f\big(\tau,x_*(\tau),u_*(\tau)\big)} \\[1mm]
&\geq& \big\langle p(\tau) , \chi_{\mathcal{A}_\lambda}(\tau) \circ \varphi\big(\tau,x_*(\tau),u_\lambda(\tau)\big) \big\rangle
                             - \chi_{\mathcal{A}_\lambda}(\tau) \circ f\big(\tau,x_*(\tau),u_\lambda(\tau)\big).
\end{eqnarray*}
Darin ist $u_\lambda(\tau)=v$ und $\mathscr{A}_\lambda$ wählt das Steuerungssystem zum Index $s \in \{1,...,k\}$ aus.
Mit der Pontrjagin-Funktion erhält die Ungleichung die Form 
\begin{equation*}
\chi_{\mathscr{A}_*}(\tau) \circ H^{\mathcal{M}}\big(\tau,x_*(\tau),u_*(\tau),p(\tau),1\big) 
\geq H_s\big(\tau,x_*(\tau),v_s,p(\tau),1\big).
\end{equation*}
Nun ist $\tau$ ein beliebiger Stetigkeitspunkt von $u_*(\cdot)$, $v$ ein beliebiger Punkt der Menge $U$ und
$s$ ein beliebiger Index aus der Menge $\{1,...,k\}$.
Daher ist die Beziehung (\ref{PMPeinfach6MP}) in allen Stetigkeitspunkten von $u_*(\cdot)$, $\chi_{\mathcal{A}_*}(\cdot)$
wahr und das Maximumprinzip bewiesen. \hfill $\blacksquare$

%% file: 3-33-Maximumprinzip.tex
\subsubsection{Das Pontrjaginsche Maximumprinzip}
\begin{theorem}[Pontrjaginsches Maximumprinzip] \label{SatzPMPhybridFest}
\index{Pontrjaginsches Maximumprinzip!Multi@-- Multiprozesse} 
Sei $\big(x_*(\cdot),u_*(\cdot),\mathscr{A}_*\big) \in \mathscr{B}^{\,\mathcal{M}}_{\rm adm} \cap \mathscr{B}^{\,\mathcal{M}}_{\rm Lip}$.
Ist $\big(x_*(\cdot),u_*(\cdot),\mathscr{A}_*\big)$ ein starkes lokales Minimum der Aufgabe (\ref{HA1})--(\ref{HA5}),
dann existieren eine Zahl $\lambda_0 \geq 0$,
Vektoren $l_0 \in \R^{s_0}$, $l_1 \in \R^{s_1}$,
eine Vektorfunktion $p(\cdot):[t_0,t_1] \to \R^n$
und auf den Mengen
$$T_j=\big\{t \in [t_0,t_1] \,\big|\, g_j\big(t,x_*(t)\big)=0\big\}, \quad j=1,...,l,$$
konzentrierte nichtnegative regul"are Borelsche Ma"se $\mu_j$ endlicher Totalvariation
(wobei s"amtliche Gr"o"sen nicht gleichzeitig verschwinden) derart,
dass die Vektorfunktion $p(\cdot)$ von beschr"ankter Variation und rechtsseitig stetig ist, und
\begin{enumerate}
\item[(a)] die adjungierte Gleichung
           \index{adjungierte Gleichung!Multi@-- Multiprozesse}
           \begin{eqnarray} 
           p(t) &=& - {h_1'}^T\big(x_*(t_1)\big) \, l_1
               + \int_t^{t_1} \chi_{\mathscr{A}_*}(s) \circ H^{\mathcal{M}}_x\big(s,x_*(s),u_*(s),p(s),\lambda_0\big) \, ds \nonumber \\
                & & \label{SatzHAPMP1}  - \sum_{j=1}^l \int_t^{t_1} g_{j,x}\big(s,x_*(s)\big) \, d\mu_j(s),
           \end{eqnarray}
\item[(b)] die Transversalit"atsbedingungen
           \index{Transversalitätsbedingungen!Multi@-- Multiprozesse}
           \begin{equation} \label{SatzHAPMP2} 
            \left. \begin{array}{lcl} p(t_0) &=& {h_0'}^T\big(x_*(t_0)\big) l_0, \\[2mm]
                  p(t_1^-)-p(t_1) &=& \displaystyle - \sum_{j=1}^l \mu_j(\{t_1\}) \, g_{j,x}\big(t_1,x_*(t_1)\big)
                             \end{array} \right\}
           \end{equation}
\item[(c)] und in fast allen Punkten $t \in [t_0,t_1]$ die Maximumbedingung
           \index{Maximumbedingung!Multi@-- Multiprozesse}
           \begin{equation} \label{SatzHAPMP3} 
           \chi_{\mathscr{A}_*}(t) \circ H^{\mathcal{M}}\big(t,x_*(t),u_*(t),p(t),\lambda_0\big)
           = \max_{\substack{u_s \in U_s\\ 1\leq s\leq k}} H_s\big(t,x_*(t),u_s,p(t),\lambda_0\big)
           \end{equation}
\end{enumerate}
erfüllt sind.
\end{theorem}


\begin{bemerkung} {\rm
Mit den Hamilton-Funktionen
\begin{eqnarray*}
\mathscr{H}^{\mathcal{M}}(t,x,p,\lambda_0) &=& \chi_{\mathscr{A}_*}(t) \circ H^{\mathcal{M}}\big(t,x,u_*(t),p,\lambda_0\big), \\
\mathscr{H}^{\mathcal{M}}_s(t,x,p,\lambda_0) &=& H_s\big(t,x,u_{*s}(t),p,\lambda_0\big)
\end{eqnarray*}
liefert die Maximumbedingung (\ref{SatzHAPMP3}) bez"uglich der optimalen Wechselstrategie die Beziehung
$$\mathscr{H}^{\mathcal{M}}\big(t,x_*(t),p(t),\lambda_0\big) = \max_{1\leq s\leq k} \mathscr{H}^{\mathcal{M}}_s\big(t,x_*(t),p(t),\lambda_0\big),$$
die man auch als Prinzip der konkurrierenden Hamilton-Funktionen bezeichnet (Bock \& Longman \cite{Bock}). \hfill $\square$}
\end{bemerkung}

%% file: 3-34-Beweis.tex
\subsubsection{Der Nachweis der notwendigen Optimalit\"atsbedingungen}
Wir betrachten f"ur die Tripel
$$\big(x(\cdot),u(\cdot),\mathscr{A}\big) \in C([t_0,t_1],\R^n) \times L_\infty([t_0,t_1],\R^m) \times \Zt$$
die Abbildungen
\begin{eqnarray*}
J\big(x(\cdot),u(\cdot),\mathscr{A}\big) &=& \int_{t_0}^{t_1}  \chi_{\mathscr{A}}(t) \circ f\big(t,x(t),u(t)\big) \, dt, \\
F\big(x(\cdot),u(\cdot),\mathscr{A}\big)(t) &=&
   x(t) -x(t_0) -\int_{t_0}^t \chi_{\mathscr{A}}(s) \circ \varphi\big(s,x(s),u(s)\big) \, ds, \quad t \in [t_0,t_1],\\
H_i\big(x(\cdot)\big) &=& h_i\big(x(t_i)\big), \quad i=0,1, \\
G_j\big(x(\cdot)\big)(t) &=& \max_{t \in [t_0,t_1]} g_j\big(t,x(t)\big), \quad j=1,...,l.
\end{eqnarray*}
Da $x(\cdot)$ zu $C([t_0,t_1],\R^n)$ geh"ort, gilt f"ur diese Abbildungen
\begin{eqnarray*}
J &:& C([t_0,t_1],\R^n) \times L_\infty([t_0,t_1],\R^m) \times \Zt \to \R, \\
F &:& C([t_0,t_1],\R^n) \times L_\infty([t_0,t_1],\R^m) \times \Zt \to C_0([t_0,t_1],\R^n), \\
H_i &:& C([t_0,t_1],\R^n) \to \R^{s_i}, \quad i=0,1, \\
G_j &:& C([t_0,t_1],\R^n) \to \R, \quad j=1,...,l.
\end{eqnarray*}
Wir setzen $\mathscr{F}=(F,H_0,H_1)$ und pr"ufen f"ur die Extremalaufgabe
\begin{equation} \label{ExtremalaufgabePMPHA}
J\big(x(\cdot),u(\cdot),\mathscr{A}\big) \to \inf, \; \mathscr{F}\big(x(\cdot),u(\cdot),\mathscr{A}\big)=0, \;
G_j\big(x(\cdot)\big) \leq 0, \; u(\cdot) \in L_\infty([t_0,t_1],U)
\end{equation}
im Punkt $\big(x_*(\cdot),u_*(\cdot),\mathscr{A}_*\big) \in \mathscr{B}^{\,\mathcal{M}}_{\rm Lip}$,
wobei wir $x_*(\cdot)$ als Element des Raumes $C([t_0,t_1],\R^n)$ auffassen,
die Voraussetzungen von Theorem \ref{SatzExtremalprinzipStark}:

\begin{enumerate}
\item[(A)] Zum Nachweis der Voraussetzungen (A$_1$)--(A$_2$) merken wir an,
           dass jedes einzelne Steuerungssystem die Eigenschaften der Aufgabe (\ref{PMP1})--(\ref{PMP5}) im letzten Kapitel
           besitzt.
           Beachten wir nun ferner die Ungleichung
           $$\|\chi_{\mathscr{A}}(t) \circ h(t,x,u)\| \leq \sum_{s=1}^k \|h^{s}(t,x,u^{s})\|,$$
           so lassen sich die Voraussetzungen (A$_1$)--(A$_2$) wie im Kapitel \ref{KapitelStark} pr"ufen.
\item[(B)] Die endliche Kodimension des Operators $\mathscr{F}_x\big(x_*(\cdot),u_*(\cdot),\mathscr{A}_*\big)$ folgt wie in Abschnitt \ref{AbschnittBeweisSOP}.
\item[(C)] Diese Eigenschaften sind im Anhang \ref{AnhangNV},
           speziell in Lemma \ref{LemmaEigenschaftNadelvariationHybrid1} und Lemma \ref{LemmaEigenschaftNadelvariationHybrid2}
           im Abschnitt \ref{AnhangNVMP}, dargestellt.
\end{enumerate}

Zur Extremalaufgabe (\ref{ExtremalaufgabePMPHA}) definieren wir auf
$$C([t_0,t_1],\R^n) \times L_\infty([t_0,t_1],\R^m) \times \Zt \times \R \times C_0^*([t_0,t_1],\R^n) \times \R^{s_0} \times \R^{s_1} \times \R^l$$
die Lagrange-Funktion $\mathscr{L}=\mathscr{L}\big(x(\cdot),u(\cdot),\mathscr{A},\lambda_0,y^*,l_0,l_1,\lambda\big)$,
\begin{eqnarray*}
\mathscr{L}
&=& \lambda_0 J\big(x(\cdot),u(\cdot),\mathscr{A}\big)+ \big\langle y^*,F\big(x(\cdot),u(\cdot),\mathscr{A}\big)\big\rangle \\
& & \hspace*{15mm}+ l_0^T H_0\big(x(\cdot)\big)+l_1^T H_1\big(x(\cdot)\big)  + \sum_{j=1}^l \lambda_j G_j\big(x(\cdot)\big).
\end{eqnarray*}
Ist $\big(x_*(\cdot),u_*(\cdot),\mathscr{A}_*\big)$ eine starke lokale Minimalstelle der Extremalaufgabe (\ref{ExtremalaufgabePMPHA}),
dann existieren nach Theorem \ref{SatzExtremalprinzipStark}
nicht gleichzeitig verschwindende Lagrangesche Multiplikatoren $\lambda_0 \geq 0$, $y^* \in C_0^*([t_0,t_1],\R^n)$, $l_i \in \R^{s_i}$
und $\lambda_1 \geq 0,...,\lambda_l \geq 0$ derart,
dass gelten:
\begin{enumerate}
\item[(a)] Die Lagrange-Funktion besitzt bez"uglich $x(\cdot)$ in $x_*(\cdot)$ einen station"aren Punkt, d.\,h.
           \begin{equation}\label{SatzPMPHALMR1}
           0 \in \partial_x \mathscr{L}\big(x_*(\cdot),u_*(\cdot),\mathscr{A}_*,\lambda_0,y^*,l_0,l_1,\lambda\big);
           \end{equation}         
\item[(b)] Die Lagrange-Funktion erf"ullt bez"uglich $\big(u(\cdot),\mathscr{A}\big)$ in $\big(u_*(\cdot),\mathscr{A}_*\big)$ die Bedingung
           \begin{equation}\label{SatzPMPHALMR2}
           \mathscr{L}\big(x_*(\cdot),u_*(\cdot),\mathscr{A}_*,\lambda_0,...,\lambda\big)
           = \min_{\substack{u(\cdot) \in L_\infty ([t_0,t_1],U) \\ \mathscr{A} \in \Zt}}
             \mathscr{L}\big(x_*(\cdot),u(\cdot),\mathscr{A},\lambda_0,...,\lambda\big);
           \end{equation}
\item[(c)] Die komplement"aren Schlupfbedingungen gelten, d.\,h.
           \begin{equation}\label{SatzPMPHALMR3}
           0 = \lambda_j G_j\big(x(\cdot)\big), \qquad i=1,...,l.
           \end{equation}
\end{enumerate}
Ebenso wie in Kapitel \ref{KapitelStark} liefert (\ref{SatzPMPHALMR3}),
dass nur diejenigen Multiplikatoren $\lambda_j \geq 0$ von Null verschieden sein k"onnen,
f"ur die die zugeh"origen Ma"se $\tilde{\mu}_j$ die Totalvariation $\|\tilde{\mu}_j\|=1$ besitzen und die auf den Mengen
$$T_j=\big\{t \in [t_0,t_1] \,\big|\, g_j\big(t,x_*(t)\big)=0\big\}$$
konzentriert sind.
Deswegen k"onnen wir ohne Einschr"ankung annehmen,
dass alle Ma"se $\mu_j=\lambda_j \tilde{\mu}_j$ auf den Mengen $T_j$ konzentriert sind. \\[2mm]
Aufgrund (\ref{SatzPMPHALMR1}) ist folgende Variationsgleichung f"ur alle $x(\cdot) \in C([t_0,t_1],\R^n)$ erf"ullt: 
\begin{eqnarray}
0 &=& \lambda_0 \int_{t_0}^{t_1} \big\langle \chi_{\mathscr{A}_*}(t) \circ f_x\big(t,x_*(t),u_*(t)\big),x(t) \big\rangle\, dt \nonumber \\
  & & + \big\langle l_0, h_0'\big(x_*(t_0)\big) x(t_0) \big\rangle + \big\langle l_1, h_1'\big(x_*(t_1)\big) x(t_1) \big\rangle
      + \sum_{j=1}^l \int_{t_0}^{t_1}\big\langle g_{j,x}\big(t,x_*(t)\big),x(t) \big\rangle \,d\mu_j(t) \nonumber \\
  & & \label{BeweisschlussPMPHA1}
      + \int_{t_0}^{t_1} \bigg[ x(t)-x(0) - \int_{t_0}^t \chi_{\mathscr{A}_*}(s) \circ \varphi_x\big(s,x_*(s),u_*(s)\big) x(s) \,ds \bigg]^T d\mu(t).
\end{eqnarray}
In der Gleichung (\ref{BeweisschlussPMPHA1}) "andern wir die Integrationsreihenfolge im letzten Summanden und bringen sie in
die Form
\begin{eqnarray}
0 &=& \int_{t_0}^{t_1} \big\langle \lambda_0 \big[\chi_{\mathscr{A}_*}(t) \circ f_x\big(t,x_*(t),u_*(t)\big)\big] \nonumber \\
  & & \hspace*{20mm} - \big[\chi_{\mathscr{A}_*}(t) \circ \varphi_x^T\big(t,x_*(t),u_*(t)\big)\big] \int_{t}^{t_1} d\mu(s) ,
                       x(t) \big\rangle \, dt \nonumber \\
  & & + \int_{t_0}^{t_1} [x(t)]^T \, d\mu(t) + \Big\langle {h_0'}^T\big(x_*(t_0)\big)l_0 - \int_{t_0}^{t_1} d\mu(t) , x(t_0) \Big\rangle
      + \langle {h_1'}^T\big(x_*(t_1)\big)l_1 , x(t_1) \rangle \nonumber \\
  & & \label{BeweisschlussPMPHA2}
      + \sum_{j=1}^l \int_{t_0}^{t_1} \big\langle g_{j,x}\big(t,x_*(t)\big),x(t) \big\rangle \,d\mu_j(t).
\end{eqnarray}
Wir setzen $p(t)=\displaystyle \int_t^{t_1} \, d\mu(s)$,
so ist $p(\cdot)$ eine Funktion von beschränkter Variation und gemäß den Eigenschaften einer Verteilungsfunktion rechtsseitig stetig. \\[1mm]
Die rechte Seite in (\ref{BeweisschlussPMPHA2}) definiert ein stetiges lineares Funktional im Raum $C([t_0,t_1],\R^n)$.
Wenden wir den Darstellungssatz von Riesz an,
so folgen aus der eindeutigen Darstellung eines stetigen linearen Funktionals im Raum $C([t_0,t_1],\R^n)$ die adjungierte Gleichung
\begin{eqnarray*}
p(t) &=& -{h_1'}^T\big(x_*(t_1)\big)l_1
         + \int_t^{t_1} \Big[ \big[\chi_{\mathscr{A}_*}(s) \circ \varphi_x^T\big(s,x_*(s),u_*(s)\big)\big]p(s) \\
     & & \hspace*{50mm} -\lambda_0 \big[\chi_{\mathscr{A}_*}(s) \circ f_x\big(s,x_*(s),u_*(s)\big)\big] \Big]\, ds \\
     &=& -{h_1'}^T\big(x_*(t_1)\big)l_1
         + \int_t^{t_1} \chi_{\mathscr{A}_*}(s) \circ H^{\mathcal{M}}_x\big(s,x_*(s),u_*(s),p(s),\lambda_0\big) ds \\
     & & - \sum_{j=1}^l \int_t^{t_1} g_{j,x}\big(s,x_*(s)\big) \,d\mu_j(s)
\end{eqnarray*}
und die Transveralitätsbedingung
$$p(t_0) = {h_0'}^T\big(x_*(t_0)\big)l_0.$$
Damit sind (\ref{SatzHAPMP1}) und (\ref{SatzHAPMP2}) gezeigt. \\[2mm]
Genauso wie im letzen Kapitel erhalten wir aus (\ref{SatzPMPHALMR2}) die Beziehung
$$\int_{t_0}^{t_1} \chi_{\mathscr{A}_*}(t) \circ H^{\mathcal{M}}\big(t,x_*(t),u_*(t),p(t),\lambda_0\big) \, dt
  \geq \int_{t_0}^{t_1} \chi_{\mathscr{A}}(t) \circ H^{\mathcal{M}}\big(t,x_*(t),u(t),p(t),\lambda_0\big) \, dt$$
f"ur alle $u(\cdot) \in L_\infty([t_0,t_1],U)$ und alle $\mathscr{A} \in \Zt$. \\[2mm]
Wir bemerken, dass die Funktion $p(\cdot)$ von beschr"ankter Variation ist.
Daher kann man $p(\cdot)$ als Differenz zweier monotoner Funktionen schreiben und
$p(\cdot)$ besitzt h"ochstens abz"ahlbar viele Unstetigkeiten.
Mit Hilfe dieser Anmerkung folgt abschlie"send die Maximumbedingung (\ref{SatzHAPMP3}) via Standardtechniken
f"ur Lebesguesche Punkte. \hfill $\blacksquare$

%% file: 3-35-FreieZeit.tex
\subsubsection{Freier Anfangs- und Endzeitpunkt}
Wir wenden uns der Aufgabe mit freiem Anfangs- und Endzeitpunkt zu:
\begin{eqnarray}
&& \label{HAFrei1} J\big(x(\cdot),u(\cdot),\mathscr{A} \big)
           = \int_{t_0}^{t_1} \chi_{\mathscr{A}}(t) \circ f\big(t,x(t),u(t)\big) dt \to \inf, \\
&& \label{HAFrei2} \dot{x}(t) = \chi_{\mathscr{A}}(t) \circ \varphi \big(t,x(t),u(t) \big), \\
&& \label{HAFrei3} h_0\big( t_0,x(t_0) \big) = 0, \quad h_1\big( t_1,x(t_1) \big) = 0, \\
&& \label{HAFrei4} u(t) \in U= U^1 \times ... \times U^k, \quad U^s \not= \emptyset, \quad \mathscr{A} \in \Zt, \\
&& \label{HAFrei5} g_j\big(t,x(t)\big) \leq 0, \quad t \in [t_0,t_1], \quad j=1,...,l.
\end{eqnarray}
Die Aufgabe (\ref{HAFrei1})--(\ref{HAFrei5}) betrachten wir bez"uglich
$$\big(x(\cdot),u(\cdot),\mathscr{A}\big) \in W^1_\infty([t_0,t_1],\R^n) \times L_\infty([t_0,t_1],U) \times \Zt,
  \qquad [t_0,t_1] \subset \R.$$
Zur Menge $\mathscr{B}^{\,\mathcal{F}}_{\rm Lip}$ geh"oren diejenigen
$\big(x(\cdot),u(\cdot),\mathscr{A}\big)$,
f"ur die es eine Zahl $\gamma>0$ derart gibt,
dass die Abbildungen $f_s(t,x,u_s)$, $\varphi_s(t,x,u_s)$, $g_j(t,x)$ und $h_i(\tau_i,x_i)$
auf der Menge aller Punkte $(t,\tau_0,\tau_1,x,x_0,x_1,u) \in \R \times \R \times \R \times \R^n \times \R^n \times \R^n \times \R^m$ mit
\begin{eqnarray*}
&& t_0-\gamma < t < t_1+\gamma,\quad t_0-\gamma < \tau_0< t_0+\gamma,\quad t_1-\gamma < \tau_1 < t_1+\gamma, \\
&& \|x-x(t)\| < \gamma, \quad \|x_0-x(t_0)\| < \gamma, \quad \|x_1-x(t_1)\| < \gamma, \quad u \in \R^m
\end{eqnarray*}
stetig in der Gesamtheit aller Variablen und
stetig differenzierbar bez"uglich der Variablen $t, t_0, t_1, x, x_0, x_1$ sind.
(In $h_0,h_1$ treten wieder die Zeitvariablen $\tau_0,\tau_1$ auf.) \\[2mm]
In der Aufgabe mit freiem Anfangs- und Endzeitpunkt nennen wir $\big([t_0,t_1],x(\cdot),u(\cdot),\mathscr{A}\big)$ mit
$[t_0,t_1] \subset \R$, $x(\cdot) \in W^1_\infty\big([t_0,t_1],\R^n\big)$, $u(\cdot) \in L^\infty\big([t_0,t_1],U\big)$ und $\mathscr{A} \in \Zt$
einen Multiprozess.
Ein Multiprozess hei"st zul"assig in der Aufgabe (\ref{HAFrei1})--(\ref{HAFrei5}),
wenn auf dem Intervall $[t_0,t_1]$ die Funktion $x(\cdot)$ fast "uberall der Gleichung (\ref{HAFrei2}) gen"ugt,
die Randbedingungen (\ref{HAFrei3})
und die Zustandsbeschr"ankungen (\ref{HAFrei5}) erf"ullt.
Die Menge $\mathscr{B}^{\,\mathcal{F}}_{\rm adm}$ bezeichnet die Menge der zul"assigen Multiprozesse. \\[2mm]
Ein zul"assiger Multiprozess $\big([t_{0*},t_{1*}],x_*(\cdot),u_*(\cdot),\mathscr{A}_*\big)$ ist ein starkes lokales
Minimum\index{Minimum, starkes lokales!Multi@-- Multiprozesse},
wenn eine Zahl $\varepsilon > 0$ derart existiert,
dass f"ur jeden anderen zul"assigen Multiprozess $\big([t_0,t_1],x(\cdot),u(\cdot),\mathscr{A}\big)$ mit den Eigenschaften
$$|t_0 - t_{0*}| < \varepsilon, \quad |t_1 - t_{1*}| < \varepsilon, \qquad
  \| x(t)-x_*(t) \| < \varepsilon \quad\mbox{ f"ur jedes } t \in [t_{0*},t_{1*}] \cap [t_0,t_1]$$
die Ungleichung $J\big(x(\cdot),u(\cdot),\mathscr{A}\big) \geq J\big(x_*(\cdot),u_*(\cdot),\mathscr{A}_*\big)$ gilt. \\[2mm]
In der Aufgabe (\ref{HAFrei1})--(\ref{HAFrei5}) bezeichnet
$$\mathscr{H}^{\mathcal{M}}(t,x,p,\lambda_0) = \max_{\substack{u_s \in U_s\\ 1\leq s\leq k}} H_s(t,x,u_s,p,\lambda_0)$$
die Hamilton-Funktion. \\[1mm]
Zum Beweis notwendiger Optimalitätsbedingungen l"asst sich wieder die Methode der Substitution der Zeit aus Abschnitt \ref{AbschnittFreieZeitSOP} anwenden
und wir erhalten das Pontrjaginschen Maximumprinzip für die Aufgabe mit freier Zeit in der Form:

\begin{theorem}[Pontrjaginsches Maximumprinzip] \label{SatzFreieZeitPMP}
\index{Pontrjaginsches Maximumprinzip!Multi@-- Multiprozesse} 
In der Aufgabe (\ref{HAFrei1})--(\ref{HAFrei5}) sei der Multiprozess
$\big([t_{0*},t_{1*}],x_*(\cdot),u_*(\cdot),\mathscr{A}_*\big) \in \mathscr{B}^{\,\mathcal{F}}_{\rm adm} \cap \mathscr{B}^{\,\mathcal{F}}_{\rm Lip}$.
Ist $\big([t_{0*},t_{1*}],x_*(\cdot),u_*(\cdot),\mathscr{A}_*\big)$ ein starkes lokales Minimum der Aufgabe
(\ref{HAFrei1})--(\ref{HAFrei5}),
dann existieren $\lambda_0 \geq 0$,
Vektoren $l_0 \in \R^{s_0}$, $l_1 \in \R^{s_1}$,
eine Vektorfunktion $p(\cdot):[t_{0*},t_{1*}] \to \R^n$
und auf den Mengen
$$T_j=\big\{t \in [t_{0*},t_{1*}] \,\big|\, g_j\big(t,x_*(t)\big)=0\big\}, \quad j=1,...,l,$$
konzentrierte nichtnegative regul"are Borelsche Ma"se $\mu_j$ endlicher Totalvariation
(wobei s"amtliche Gr"o"sen nicht gleichzeitig verschwinden) derart,
dass die Vektorfunktion $p(\cdot)$ von beschr"ankter Variation und rechtsseitig stetig ist, und
\begin{enumerate}
\item[(a)] die adjungierte Gleichung
           \index{adjungierte Gleichung!Multi@-- Multiprozesse}
           \begin{eqnarray} 
           p(t) &=& -h^T_{1,x_1}\big(t_{1*},x_*(t_{1*})\big) l_1 +
                    \int_t^{t_{1*}} \chi_{\mathscr{A}_*}(s) \circ H^{\mathcal{M}}_x\big(s,x_*(s),u_*(s),p(s),\lambda_0\big) \, ds \nonumber \\
                & & \label{SatzFreieZeitPMP1} - \sum_{j=1}^l \int_t^{t_{1*}} g_{j,x}\big( s,x_*(s)\big) \, d\mu_j(s),
           \end{eqnarray}
\item[(b)] die Transversalit"atsbedingungen
           \index{Transversalitätsbedingungen!Multi@-- Multiprozesse}
           \begin{equation} \label{SatzFreieZeitPMP2}
            \left. \begin{array}{lcl} p(t_{0*}) &=& h^T_{0,x_0} \big(t_{0*},x_*(t_{0*})\big) l_0, \\[2mm]
                  p(t_{1*}^-)-p(t_{1*}) &=& \displaystyle - \sum_{j=1}^l \mu_j(\{t_{1*}\}) \, g_{j,x}\big(t_{1*},x_*(t_{1*})\big)
                             \end{array} \right\}
           \end{equation} 
\item[(c)] und in fast allen Punkten $t\in [t_{0*},t_{1*}]$ die Maximumbedingung 
           \index{Maximumbedingung!Multi@-- Multiprozesse}
           \begin{equation} \label{SatzFreieZeitPMP3} 
           \chi_{\mathscr{A}_*}(t) \circ H^{\mathcal{M}}\big(t,x_*(t),u_*(t),p(t),\lambda_0\big)
           = \max_{\substack{u_s \in U_s\\ 1\leq s\leq k}} H_s\big(t,x_*(t),u_s,p(t),\lambda_0\big)
           \end{equation}
\end{enumerate}
erfüllt sind und weiterhin
\begin{enumerate}
\item[(d)] die Beziehungen 
           \begin{eqnarray}
           \mathscr{H}^{\mathcal{M}}\big(t,x_*(t),p(t),\lambda_0\big)
                       &=& \big\langle h_{1,t_1}\big(t_{1*},x_*(t_{1*})\big) , l_1 \big\rangle \nonumber \\
           && -\int_t^{t_{1*}} \chi_{\mathscr{A}_*}(s) \circ H^{\mathcal{M}}_t\big(s,x_*(s),u_*(s),p(s),\lambda_0\big) \, ds  \nonumber \\
           && \label{SatzFreieZeitPMP4} 
               + \sum_{j=1}^l \int_t^{t_{1*}} g_{j,t}\big( s,x_*(s)\big) \, d\mu_j(s), \\
           \mathscr{H}^{\mathcal{M}}\big(t_{0*},x_*(t_{0*}),p(t_{0*}),\lambda_0\big)
           &=& \label{SatzFreieZeitPMP5} - \big\langle h_{0,t_0}\big(t_{0*},x_*(t_{0*})\big) , l_0 \big\rangle
           \end{eqnarray}
\end{enumerate}
gelten.
\end{theorem}

%% file: 3-36-Investitionsmodell.tex
\subsubsection{Ein Investitionsmodell} \label{AbschnittBeispielhybrid}
Wir\index{Kapitalakkumulation} untersuchen ein Investitionsmodell,
dass sich aus dem linearen und dem konkaven Beispielen \ref{BeispielLinInv}, \ref{BeispielKonInv} zusammensetzt:
\begin{eqnarray*}
&& J\big(x(\cdot),u(\cdot),\mathscr{A} \big)
  = \int_0^T \chi_{\mathscr{A}}(t) \circ f\big(t,x(t),u(t)\big) dt \to \sup, \\
&& \dot{x}(t) = \chi_{\mathscr{A}}(t) \circ \varphi\big(t,x(t),u(t)\big), \quad x(0)=x_0>0, \\
&& u(t) \in [0,1] \times [0,1], \quad \mathscr{A}=\{ \mathscr{A}_1,\mathscr{A}_2 \} \in \mathscr{Z}^2([0,T]),
\end{eqnarray*}
mit den einzelnen Steuerungssystemen
$$\begin{array}{ll}
     f_1(t,x,u_1) = (1-u_1)x,        & \varphi_1(t,x,u_1) = u_1 x, \\[1mm]
     f_2(t,x,u_2) = (1-u_2)x^\alpha, & \varphi_2(t,x,u_2) = u_2 x^\alpha,
  \end{array} \quad U^1=U^2=[0,1]$$
und mit den Modellparametern
$$\alpha \in (0,1) \mbox{ konstant}, \qquad T \mbox{ fest mit } T > \max \left\{1, \frac{x_0^{1-\alpha}}{\alpha} \right\}.$$
Da $U_1=U_2$ ist, unterscheiden wir nicht zwischen den Steuervariablen $u_1,u_2$ und schreiben $u$.
Wir wenden Theorem \ref{SatzPMPhybridFest} mit $\lambda_0=1$ an.
Es gilt nach (\ref{SatzHAPMP3}):
\begin{eqnarray*}
\lefteqn{\chi_{\mathscr{A}_*}(t) \circ H^{\mathcal{M}}\big(t,x_*(t),u_*(t),p(t),1\big)
         =\max_{s \in \{1,2\}} H_s\big(t,x_*(t),u_*(t),p(t),1\big)} \\
&=& \max_{s \in \{1,2\}} \Big\{ \max_{u \in [0,1]} \big[ \big( p(t) -1 \big) u + 1 \big] \cdot x_*(t),
                \max_{u \in [0,1]} \big[ \big( p(t) -1 \big) u + 1 \big] \cdot x_*^\alpha(t) \Big\}.
\end{eqnarray*}
Dies k"onnen wir weiterhin in die Form
$$\max_{u \in [0,1]} \big[ \big( p(t) -1 \big) u \big] \cdot \max_{s \in \{1,2\}}\{ x_*(t) , x_*^\alpha(t) \}$$
bringen.
Daraus erhalten wir f"ur die optimale Investitionsrate und Wechselstrategie
$$u_*(t) = \left\{ \begin{array}{ll} 1, & p(t) < 1, \\ 0, & p(t) > 1, \end{array} \right. \qquad
  \chi_{\mathscr{A}_*}(t) = \left\{ \begin{array}{ll} ( 1,0 ), & x_*(t) > 1, \\ ( 0,1 ), & 0 < x_*(t) < 1. \end{array} \right.$$
Die Funktion $p(\cdot)$ ist die L"osung der adjungierten Gleichung (\ref{SatzHAPMP1}):
$$\dot{p}(t) = \left\{ \begin{array}{ll} - \big[ \big( p(t) -1 \big) u_*(t) + 1 \big],    & t \in \mathscr{A}_1, \\[1mm]
          - \big[ \big( p(t) -1 \big) u_*(t) + 1 \big] \cdot \alpha x_*^{\alpha-1}(t), & t \in \mathscr{A}_2, \end{array} \right.
   \quad p(T)=0.$$
Betrachten wir die einzelnen Steuerungssysteme,
so sind nach den Beispielen \ref{BeispielLinInv}, \ref{BeispielKonInv} durch $t_1=T-1$ bzw. $t_2=\alpha T -x_0^{1-\alpha}$
die Zeitpunkte f"ur den optimalen Wechsel von vollst"andiger Investition in komplette Kosumption gegeben.
Au"serdem ist f"ur die L"osung $x(\cdot)$ der Differentiagleichung $\dot{x}(t)=x^\alpha(t)$ mit $x_0 \in (0,1)$ der
Zeitpunkt $\sigma$, in dem $x(\sigma)=1$ gilt, durch $\sigma=\frac{1-x_0^{1-\alpha}}{1-\alpha}$ bestimmt.
Im Weiteren seien also
$$t_1=T-1,  \qquad t_2=\alpha T -x_0^{1-\alpha}, \qquad \sigma=\frac{1-x_0^{1-\alpha}}{1-\alpha}.$$
Wir diskutieren die einzelnen m"oglichen Szenarien:
\begin{enumerate}
\item[(a)] Sei $x_0 \geq 1$: In diesem Fall lautet der Kandidat
           \begin{eqnarray*}
           && x_*(t) = \left\{ \begin{array}{ll} x_0 \cdot e^t, & t \in [0,t_1), \\
                                                 x_0 \cdot e^{t_1}, & t \in [t_1,T], \end{array} \right.
              \qquad u_*(t) = \left\{ \begin{array}{ll} 1, & t \in [0,t_1), \\
                                                         0, & t \in [t_1,T], \end{array} \right. \\
           && \chi_{\mathscr{A}_*}(t) = (1,0), \qquad J\big(x_*(\cdot),u_*(\cdot),\mathscr{A}_*\big) = x_0 \cdot e^{t_1}.          
           \end{eqnarray*}
\item[(b)] Seien $x_0<1$ und $T < \frac{1 - \alpha x_0^{1-\alpha}}{\alpha (1-\alpha)}$:
           Dann lautet der Kandidat
           \begin{eqnarray*}
           && y_*(t) = \left\{ \begin{array}{ll}
                       \big[ (1-\alpha)t +  x_0^{1-\alpha} \big]^\frac{1}{1-\alpha}, & t \in [0,t_2), \\
                       \big[ \alpha(1-\alpha)T + \alpha x_0^{1-\alpha} \big]^\frac{1}{1-\alpha}, & t \in [t_2,T],
                       \end{array} \right.
              \quad v_*(t) = \left\{ \begin{array}{ll} 1, & t \in [0,t_2), \\
                                                         0, & t \in [t_2,T], \end{array} \right. \\
           && \chi_{\mathscr{B}_*}(t) = (0,1), \qquad 
              J\big(y_*(\cdot),v_*(\cdot),\mathscr{B}_*\big) = 
              \alpha^\frac{\alpha}{1-\alpha} \cdot \big[ (1-\alpha)T + x_0^{1-\alpha} \big]^\frac{1}{1-\alpha}.          
           \end{eqnarray*}
\item[(c)] Seien $x_0<1$ und $T > \frac{1 - \alpha x_0^{1-\alpha}}{\alpha (1-\alpha)}$:
           Nun erhalten wir den Kandidaten
           \begin{eqnarray*}
           && z_*(t) = \left\{ \begin{array}{ll}
                       \big[ (1-\alpha)t +  x_0^{1-\alpha} \big]^\frac{1}{1-\alpha}, & t \in [0,\sigma), \\
                       e^{t-\sigma}, & t \in [\sigma,t_1), \\
                       e^{t_1-\sigma}, & t \in [t_1,T], \end{array} \right.                      
              \quad w_*(t)= \left\{ \begin{array}{ll} 1, & t \in [0,t_1), \\
                                                      0, & t \in [t_1,T], \end{array} \right. \\
           && \chi_{\mathscr{C}_*}(t) = \left\{ \begin{array}{ll} (0,1), & t \in [0,\sigma), \\
                                                         (1,0), & t \in [\sigma,T], \end{array} \right.
              \qquad J\big(z_*(\cdot),w_*(\cdot),\mathscr{C}_*\big) = e^{T-1-\sigma}.        
           \end{eqnarray*}
\item[(d)] Seien $x_0<1$ und $T = \frac{1 - \alpha x_0^{1-\alpha}}{\alpha (1-\alpha)}$:
           Wegen $\sigma=t_2$ sind alle Bedingungen von Theorem \ref{SatzPMPhybridFest} f"ur die Multiprozesse
           $\big(y_*(\cdot),v_*(\cdot),\mathscr{B}_*\big)$ und $\big(z_*(\cdot),w_*(\cdot),\mathscr{C}_*\big)$ erf"ullt. \\
           Sei $n_0 \in \N$ mit $n_0 >\alpha/(1-\alpha)$.
           F"ur $n \geq n_0$ betrachten wir die Steuerungen
           $$v_n(t)= v_*(t) + \chi_{[\sigma,\sigma+\frac{1}{n})}(t) \cdot \big(w_*(t)-v_*(t)\big),
             \quad \chi_{\mathscr{B}_n}(t)= \chi_{\mathscr{C}_*}(t),$$
           die wie im Fall (c) einen Wechsel des Steuerungssystems zum Zeitpunkt $t=\sigma=t_2$ und eine verl"angerte Investitionsphase
           vorgeben.
           F"ur den zugeh"origen Kapitalbestand $y_n(\cdot)$ gilt $\|y_n(\cdot)-y_*(\cdot)\|_\infty=e^{\frac{1}{n}}-1$
           und f"ur den Wert des Zielfunktionals
           \begin{eqnarray*}
           \lefteqn{J\big(y_n(\cdot),v_n(\cdot),\mathscr{B}_n\big)-J\big(y_*(\cdot),v_*(\cdot),\mathscr{B}_*\big)
              =e^{\frac{1}{n}}\bigg(T-\sigma-\frac{1}{n}\bigg)- (T-\sigma)} \\
           && >\bigg(1+\frac{1}{n}\bigg)\bigg(\frac{1}{\alpha}-\frac{1}{n}\bigg)-\frac{1}{\alpha}
              =\frac{1}{n}\bigg(\frac{1}{\alpha}-1-\frac{1}{n}\bigg)>0 \quad\mbox{ f"ur alle } n \geq n_0.
           \end{eqnarray*}
           Daher stellt $\big(y_*(\cdot),v_*(\cdot),\mathscr{B}_*\big)$ in diesem Fall kein starkes lokales Maximum dar.
           Der Multiprozess $\big(z_*(\cdot),w_*(\cdot),\mathscr{C}_*\big)$ ist der einzige der Kandidat. \hfill $\square$
\end{enumerate}

\newpage
Die Aufgabe eines optimalen Multiprozesses enth"alt die Menge $\Zt$,
welche die Wechselstrategien charakterisiert und keine konvexe Menge darstellt.
Im vorliegenden Beispiel sind deswegen die hinreichenden Arrow-Bedingungen 
\index{hinreichende Bedingungen, Arrow!Multi@-- Multiprozesse}
nur eingeschränkt anwendbar.
Der Beweis der nachstehenden Bedingungen erfolgt wie im Abschnitt \ref{AbschnittHBPMP}.

\begin{theorem} \label{SatzHBMP}
In der Aufgabe (\ref{PMPeinfach1MP})--(\ref{PMPeinfach3MP}) sei
$\big(x_*(\cdot),u_*(\cdot),\mathscr{A}_*\big) \in \mathscr{B}^{\,\mathcal{M}}_{\rm adm} \cap \mathscr{B}^{\,\mathcal{M}}_{\rm Lip}$
und es sei $p(\cdot) \in PC_1([t_0,t_1],\R^n)$. Ferner gelte:
\begin{enumerate}
\item[(a)] Das Quadrupel $\big(x_*(\cdot),u_*(\cdot),\mathscr{A}_*,p(\cdot)\big)$
           erf"ullt (\ref{PMPeinfach4MP})--(\ref{PMPeinfach6MP}) in Theorem \ref{SatzPMPeinfachMP}.        
\item[(b)] F"ur jedes $t \in [t_0,t_1]$ ist die Funktion $\mathscr{H}^{\mathcal{M}}\big(t,x,p(t),1\big)$ konkav in $x$ auf $V^{\mathcal{S}}_\gamma(t)$.
\end{enumerate}
Dann ist $\big(x_*(\cdot),u_*(\cdot),\mathscr{A}_*\big)$ ein starkes lokales Minimum der Aufgabe (\ref{PMPeinfach1MP})--(\ref{PMPeinfach3MP}).
\end{theorem}

In unserem Investitionsmodell lauten die Pontrjagin-Funktionen $H_1,H_2$
zu den beiden Steuerungssystemen
$$H_1(t,x,u,p,1)=[(p -1) u + 1 ] \cdot x, \quad H_2(t,x,u,p,1)=[(p -1) u + 1 ] \cdot x^\alpha.$$
Dies führt zu der Hamilton-Funktion
$$\mathscr{H}^{\mathcal{M}}(t,x,p,1) = \max_{s \in \{1,2\}} \big\{\max_{u_s \in [0,1]} H_s(t,x,u_s,p,1)\big\}= \max\{p,1\} \cdot \max\{ x , x^\alpha \},$$
welche die Funktion
$$f(x)= \max\{p,1\} \cdot \left\{\begin{array}{ll} x^\alpha, & x \in (0,1), \\ x, & x\geq 1, \end{array} \right.$$
enthält.
Die Funktion $f(x)$ ist nicht konkav an der Stelle $x=1$. \\[2mm]
Damit ist die Hamilton-Funktion $\mathscr{H}^{\mathcal{M}}\big(t,x,p(t),1\big)$ bezüglich der Trajektorie $x_*(\cdot)$
genau dann f"ur jedes $t \in [t_0,t_1]$ konkav in der Variable $x$ auf $V^{\mathcal{S}}_\gamma(t)$,
wenn $x_*(t) \not=1$ f"ur alle $t \in [t_0,t_1]$ gilt.
F"ur die einzelnen F"alle ergeben sich damit:
\begin{enumerate}
\item[(a)] Theorem \ref{SatzHBMP} ist nur dann anwendbar, wenn $x_0>1$ gilt.
           In diesem Fall ist der Kandidat $\big(x_*(\cdot),u_*(\cdot),\mathscr{A}_*\big)$ optimal.
\item[(b)] Da in diesem Fall $x_0\leq y_*(t)<1$ auf $[0,T]$ gilt, ist Theorem \ref{SatzHBMP} anwendbar und
           der Kandidat $\big(y_*(\cdot),v_*(\cdot),\mathscr{B}_*\big)$ optimal.
\item[(c)] Wegen $z_*(\sigma)=1$ ist Theorem \ref{SatzHBMP} nicht anwendbar.
\item[(d)] Wegen $y_*(\sigma)=1$ und $z_*(\sigma)=1$ ist Theorem \ref{SatzHBMP} auf keinen der ermittelten Kandidaten
           $\big(y_*(\cdot),v_*(\cdot),\mathscr{B}_*\big)$, $\big(z_*(\cdot),w_*(\cdot),\mathscr{C}_*\big)$ anwendbar.
           Die lokale Optimalität des Kandidaten $\big(y_*(\cdot),v_*(\cdot),\mathscr{B}_*\big)$ wurde bereits
           ausgeschlossen. \hfill $\square$
\end{enumerate}

%% file: 3-37-Wassercontainer.tex
\subsubsection{Zeitoptimale Steuerung gekoppelter Kompartimente}
\begin{figure}[h]
	\centering
	\fbox{\includegraphics[width=6cm]{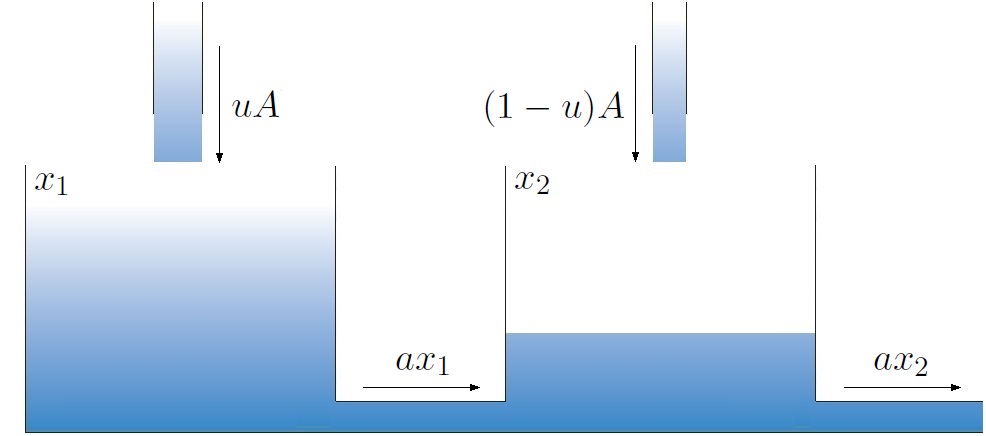}} \hspace*{1cm} \fbox{\includegraphics[width=6cm]{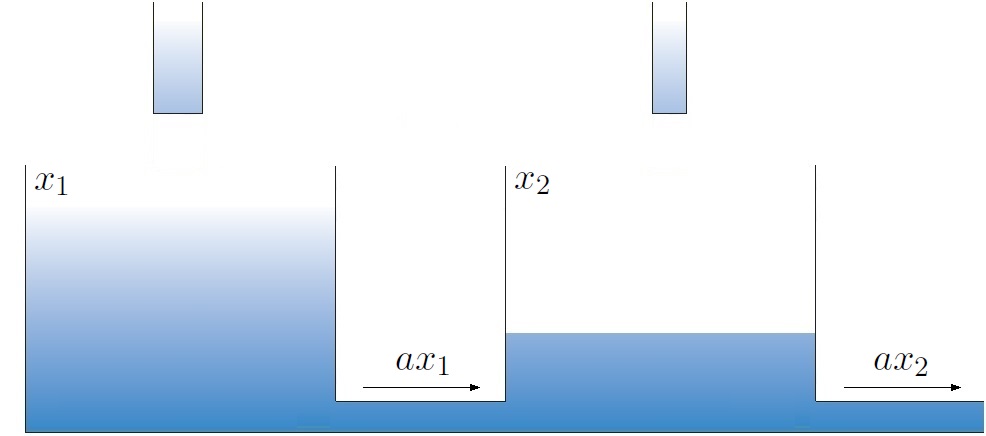}}
	\caption[Gekoppelte Kompartimente]{Die Beh"alter in den verschiedenen Steuerungssystemen.}
	\label{AbbildungContainer}
\end{figure}

Kompartimentmodelle\index{Kompartimentmodelle} sind vereinfachte Darstellungen f"ur komplexe Abl"aufe und Interaktionen.
Anwendung finden sie z.\,B. bei der Abbildung von Mischungsprozessen oder der chemischen Reaktionskinetik oder bei einer
medikament"osen Dauertherapie (Heuser \cite{HeuserGD}). \\[2mm]
Das System in Abbildung \ref{AbbildungContainer} zeigt zwei gekoppelte Beh"alter,
die mit einer Fl"ussigkeit gef"ullt sind.
In diesem Modell besteht die Kopplung der Beh"alter durch den zustandsabh"angigen Zufluss $ax_1$ vom ersten
in den zweiten Beh"alter.
Der zweite Beh"alter gibt seinen Inhalt mit der Rate $ax_2$ an die Umwelt ab.
Die beiden Steuerungssysteme bestehen einerseits aus dem Auff"ullen der Beh"alter mit konstanter Rate $A$,
die zwischen beiden Beh"altern durch die Steuerung $u$ aufgeteilt wird.
Andererseits kann die Bef"ullung komplett gestoppt werden.

Dies liefert folgende Dynamiken in den jeweligen Steuerungssystemen:
$$\begin{array}{l} \dot{x}_1(t) = -ax_1(t)+u(t)A, \\ \dot{x}_2(t)=ax_1(t)-ax_2(t)+(1-u(t))A, \end{array}
  \qquad\qquad \begin{array}{l} \dot{x}_1(t)=-ax_1(t),  \\ \dot{x}_2(t)=ax_1(t)-ax_2(t). \end{array}$$
Dabei beachte man, dass die ungesteuerte zweite Dynamik kein Spezialfall der ersten Dynamik ist.
In diesem Modell besteht das Ziel darin,
die Anfangszust"ande $(x_1^0,x_2^0)$ in k"urzester Zeit in einen Sollzustand $(x_1^f,x_2^f)$ zu "uberf"uhren.
Die Aufgabe des zeitoptimalen Multiprozesses lautet damit:
\begin{eqnarray}
&& \label{Container1} J\big(x(\cdot),u(\cdot),\mathscr{A} \big) = \int_0^T 1 \, dt \to \inf, \\
&& \label{Container2} \dot{x}(t) = \chi_{\mathscr{A}}(t) \circ \varphi \big(t,x(t),u(t) \big), \\
&& \label{Container3} \big(x_1(0),x_2(0)\big)=(x_1^0,x_2^0), \qquad \big(x_1(T),x_2(T)\big)=(x_1^f,x_2^f), \\
&& \label{Container4} u(t) \in [0,1], \quad \mathscr{A} \in \mathscr{Z}^2([0,T]), \quad a,A>0.
\end{eqnarray}
Wir wenden Theorem \ref{SatzFreieZeitPMP} an:
Die Pontrjagin-Funktion der Aufgabe (\ref{Container1})--(\ref{Container4}) lautet
$$\mathscr{A} \circ H^{\mathcal{M}}(t,x,u,p,q,\lambda_0\big)=
  \left\{ \begin{array}{ll}
          p[-ax_1+uA]+q[ax_1-ax_2+(1-u)A]-\lambda_0, & t \in \mathscr{A}_1, \\[1mm]
          p[-ax_1]+q[ax_1-ax_2]-\lambda_0, & t \in \mathscr{A}_2. 
         \end{array}\right.$$
Unabh"angig von $\mathscr{A}$ erhalten wir daraus f"ur die Adjungierten $p(\cdot),q(\cdot)$:
\begin{equation} \label{Container5}
\bigg(\begin{array}{l} \dot{p}(t) \\ \dot{q}(t) \end{array}\bigg) =
  a \cdot \bigg(\begin{array}{r} p(t)-q(t) \\ q(t) \end{array}\bigg) \qquad\Rightarrow\qquad
  \bigg(\begin{array}{l} p(t) \\ q(t) \end{array}\bigg) =
  e^{at} \cdot \bigg(\begin{array}{l} p_0-atq_0 \\ q_0 \end{array}\bigg).
\end{equation}
Die Maximumbedingung (\ref{SatzFreieZeitPMP3}) ist "aquivalent zur Maximierungsaufgabe
\begin{equation} \label{Container6} \max_{s \in \{1,2\}} \Big\{ \max_{u\in [0,1]} [up(t)+(1-u)q(t)],\; 0 \Big\}. \end{equation}
Au"serdem folgt, da die Aufgabe autonom ist,
$\mathscr{H}^{\mathcal{M}}\big(t,x_*(t),p(t),\lambda_0\big) \equiv 0$. \\[2mm]
Zur Auswertung der Maximierungsaufgabe (\ref{Container6}) sind die Relationen zwischen $p(t),q(t)$ und $r(t)\equiv 0$ zu kl"aren.
Zun"achst gilt $(p_0,q_0)\not=(0,0)$ in (\ref{Container5}), denn sonst w"urde f"ur
\begin{enumerate}
\item[(i)] $\lambda_0=0$ der Fall trivialer Multiplikatoren in Theorem \ref{SatzFreieZeitPMP} bzw.
\item[(ii)] $\lambda_0\not=0$ ein Widerspruch zu $\mathscr{H}^{\mathcal{M}}\big(t,x_*(t),p(t),\lambda_0\big) \equiv 0$
\end{enumerate}
auftreten.
Dann ergeben sich in Abh"angigkeit von $q_0 \in \R$ folgende Situation f"ur die Adjungierten,
wobei $\sigma_S$ den m"oglichen Schnittpunkt der Graphen von $p(\cdot)$ und $q(\cdot)$ bezeichnet:
\begin{enumerate}
\item[(1)] $q_0=0$: $p(t)>q(t)\equiv 0$ f"ur $p_0>0$ oder $p(t)<q(t)\equiv 0$ f"ur $p_0<0$.
\item[(2)] $q_0>0$: $p(t)<q(t)$ f"ur $p_0\leq q_0$ und $t>0$ oder es existiert $\sigma_S>0$ f"ur $p_0>q_0$.
\item[(3)] $q_0<0$: $p(t)>q(t)$ f"ur $p_0\geq q_0$ und $t>0$ oder es existiert $\sigma_S>0$ f"ur $p_0<q_0$.
\end{enumerate}
\begin{figure}[h]
	\centering
	\fbox{\includegraphics[height=4cm]{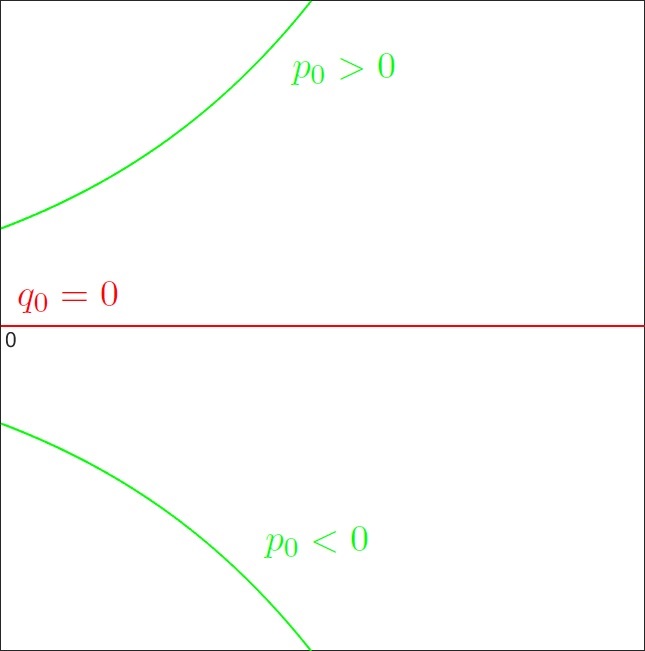} \hspace*{5mm}
          \includegraphics[height=4cm]{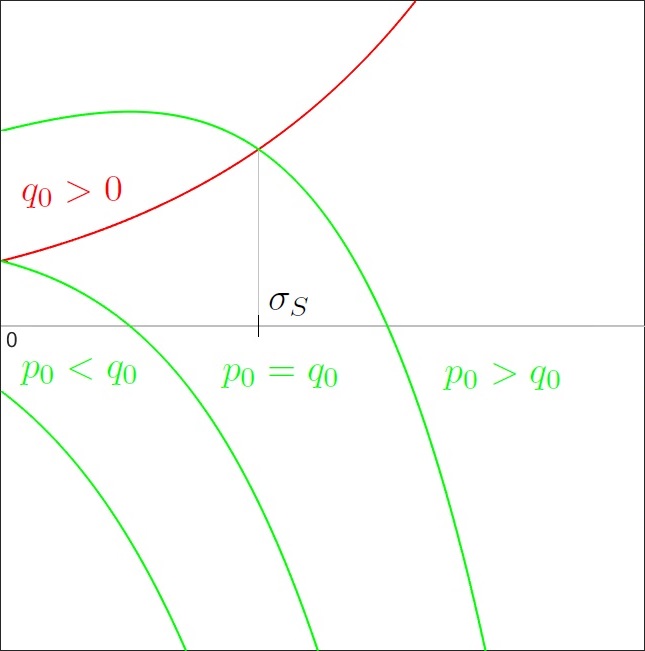} \hspace*{5mm}
          \includegraphics[height=4cm]{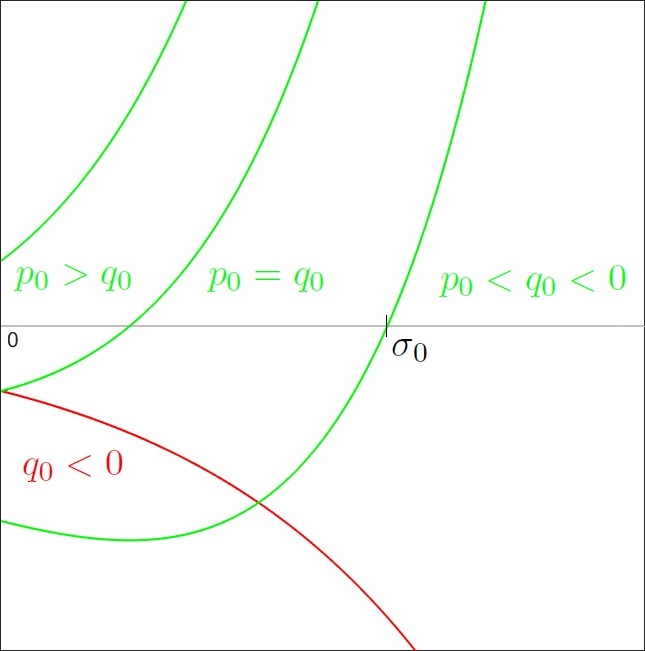}}
	\caption[Gekoppelte Kompartimente - Darstellung der Adjungierten]{Lagebeziehungen (1)--(3) der Adjungierten $p(\cdot)$ und $q(\cdot)$.}
	\label{AbbildungKompAdj}
\end{figure}
Es sei $\sigma_0$ die Nullstelle von $p(\cdot)$.
Unter Beachtung der optimalen Stoppzeit $T_*$ lassen sich der Abbildung \ref{AbbildungKompAdj}
folgende Wechselstrategien $\mathscr{A}_*$ und Steuerungen $u_*(t)$ entnehmen:
\begin{enumerate}
\item[(a)] Es ist nur ein Steuerungssystem aktiv auf $[0,T_*]$, wenn
          \begin{enumerate}
          \item[(i)] $p(t)>\max\{q(t),0\}$ f"ur alle $t \in (0,T_*)$ gilt und es ist
                     $$\big(u_*(t),\chi_{\mathscr{A}_{1*}}(t),\chi_{\mathscr{A}_{2*}}(t)\big) \equiv (1,1,0);$$
          \item[(ii)] $q(t)>\max\{p(t),0\}$ f"ur alle $t \in (0,T_*)$ gilt und es ist
                      $$\big(u_*(t),\chi_{\mathscr{A}_{1*}}(t),\chi_{\mathscr{A}_{2*}}(t)\big) \equiv (0,1,0);$$
          \item[(iii)] $0>\max\{p(t),q(t)\}$ f"ur alle $t \in (0,T_*)$ gilt und es ist
                       $$\big(\chi_{\mathscr{A}_{1*}}(t),\chi_{\mathscr{A}_{2*}}(t)\big) \equiv (0,1).$$
          \end{enumerate}
\item[(b)] Es gibt einen Strategiewechsel, und zwar wenn
          \begin{enumerate}
          \item[(i)] $p_0>q_0>0$ und $\sigma_S \in (0,T_*)$:
                    $$\big(u_*(t),\chi_{\mathscr{A}_{1*}}(t),\chi_{\mathscr{A}_{2*}}(t)\big) = (1,1,0)
                     \quad\longrightarrow \quad \big(u_*(t),\chi_{\mathscr{A}_{1*}}(t),\chi_{\mathscr{A}_{2*}}(t)\big) = (0,1,0);$$
          \item[(ii)] $p_0,q_0<0$ und $\sigma_0 \in (0,T_*)$:
                     $$\big(\chi_{\mathscr{A}_{1*}}(t),\chi_{\mathscr{A}_{2*}}(t)\big) = (0,1)
                      \quad\longrightarrow \quad\big(u_*(t),\chi_{\mathscr{A}_{1*}}(t),\chi_{\mathscr{A}_{2*}}(t)\big) = (1,1,0).$$
          \end{enumerate}
\item[(c)] Der Strategiewechsel ist nicht eindeutig bestimmbar, wenn $p_0<q_0=0$:
          $$\big(u_*(t),\chi_{\mathscr{A}_{1*}}(t),\chi_{\mathscr{A}_{2*}}(t)\big) = (0,1,0)
            \quad\longleftrightarrow \quad  \big(\chi_{\mathscr{A}_{1*}}(t),\chi_{\mathscr{A}_{2*}}(t)\big) = (0,1).$$
\end{enumerate}

\begin{figure}[h]
	\centering
	\fbox{\includegraphics[width=3.8cm]{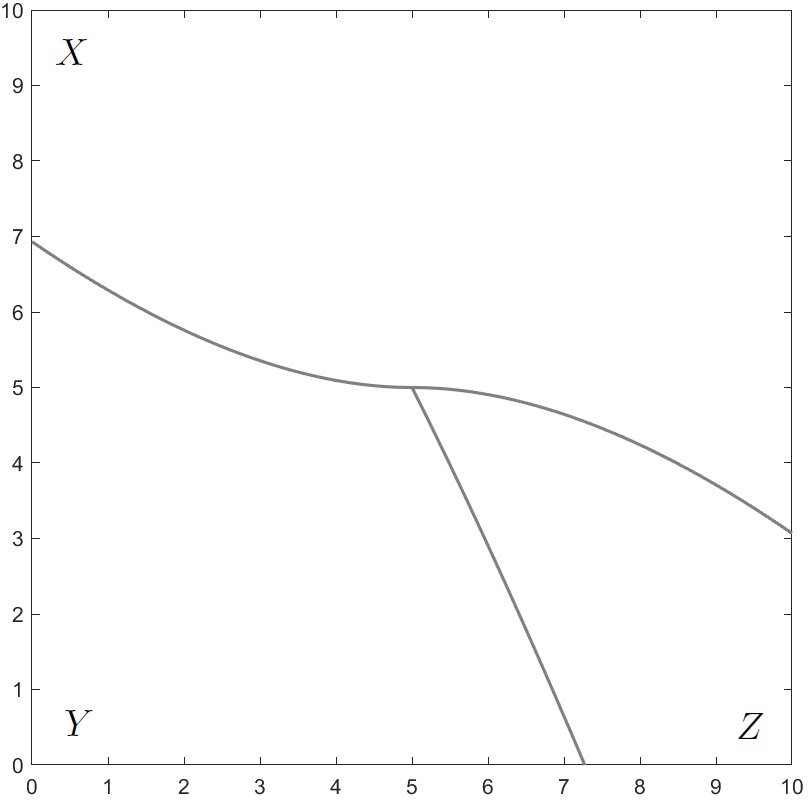}} \hspace*{1cm} \fbox{\includegraphics[width=3.8cm]{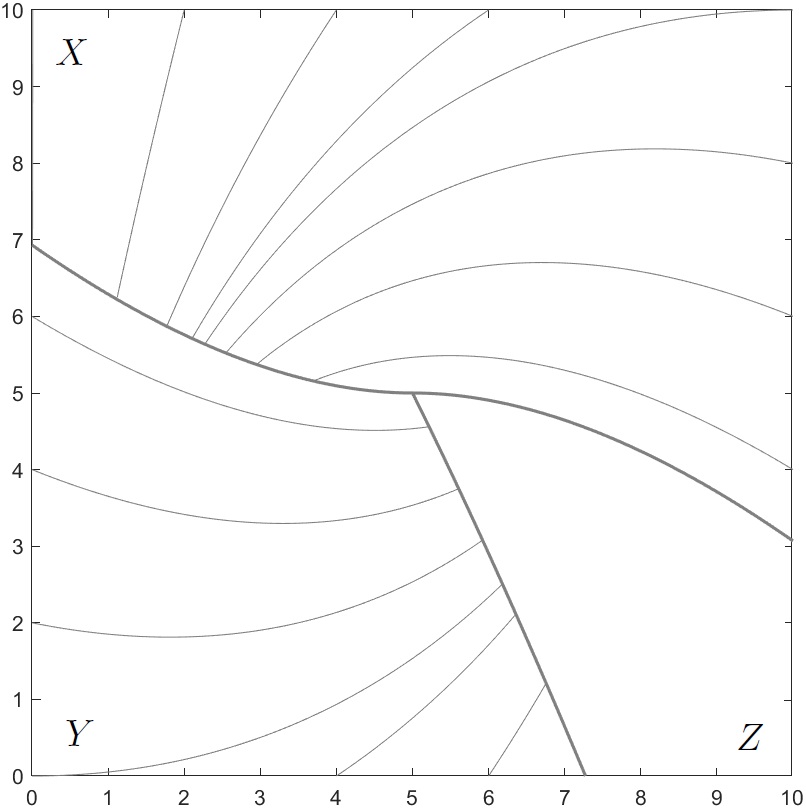}}
	\caption[Gekoppelte Kompartimente - Phasendiagramm]{Die Einteilung des Phasendiagramms und die optimalen Trajektorien.}
	\label{AbbildungBereiche}
\end{figure}

Zur Illustration sei $a=1$, $A=10$, $x_1^0,x_2^0 \in [0,10]$ und $(x_1^f,x_2^f)=(5,5)$.
Die Zust"ande $\big(x_{1*}(t),x_{2*}(t)\big)$,
die zu den Steuerungen und dem aktiven System im Fall (a) geh"oren,
teilen das Phasendiagramm $(x_1,x_2) \in [0,10] \times [0,10]$ in die Bereiche $X,Y,Z$ ein (Abbildung \ref{AbbildungBereiche}).

Zur Erl"auterung der Vorg"ange innerhalb der einzelnen Bereiche:
\begin{enumerate}
\item[($X$)] Die optimale Steuerung und Wechselstrategie wird durch (ii) im Fall (b) gegeben.
             Zur Zeit $t=\sigma_0$ erreicht die Trajektorie $\big(x_{1*}(t),x_{2*}(t)\big)$ den Rand des Bereiches $X$.
\item[($Y$)] Die optimale Steuerung und Wechselstrategie wird durch (i) im Fall (b) gegeben.
             Zur Zeit $t=\sigma_S$ erreicht die Trajektorie $\big(x_{1*}(t),x_{2*}(t)\big)$ den Rand des Bereiches $Y$.
\item[($Z$)] Die optimale Strategie ist nicht eindeutig, denn es tritt der Fall (c) ein.
             In diesem Bereich bestimmt sich die optimale Stoppzeit $T_*$ durch die schnellste Reduktion des Zustandes $x_1^0$ zu $x^f_1$.
             Im zweiten Beh"alter muss zum Zeitpunkt $t=T_*$ lediglich die Bedingung $x_2(T_*)=x_2^f$ sichergestellt werden.
\end{enumerate}

Wir f"ugen dem Beh"altermodell (\ref{Container1})--(\ref{Container4}) eine Zustandsbeschr"ankung hinzu:
$$x_2(t) \geq S \mbox{ f"ur alle } t \in [0,T], \qquad x_2^0,x_2^f,A>S.$$
Sei nun $\big(x_{1*}(t),x_{2*}(t)\big)$ eine optimale Trajektorie in der Aufgabe (\ref{Container1})--(\ref{Container4})
ohne Zustandsbeschr"ankung mit $x_{2*}(t)<S$ f"ur ein $t \in (0,T_*)$.
Dann existieren Zeitpunkte $t_1,t_2$ mit $0<t_1<t_2<T_*$ und $x_2(t_1)=x_2(t_2)=S$. \\
Aus Theorem \ref{SatzFreieZeitPMP} k"onnen wir folgenden Kandidaten ableiten:
Auf dem Teilst"uck $[t_1,t_2]$,
auf dem die Zustandsbeschr"ankung aktiv wird, erhalten wir f"ur die Adjungierten
$$p(t)=q(t)>0, \qquad \dot{p}(t)=\dot{q}(t)=0.$$
F"ur das Borelsche Ma"s liefert der Ansatz $d\mu(t)=\lambda(t)\,dt$ die Darstellung
$$\lambda(t)=a\big(2q(t)-p(t)\big)>0.$$
Eingeschr"ankt auf das Teilst"uck $[t_1,t_2]$ erhalten wir f"ur den Multiprozess
\begin{eqnarray*}
&& x_{1*}(t)=x_{1*}(t_1)+(t-t_1)(A-S), \quad x_{2*}(t)=S, \\
&& u_*(t)= \frac{x_{1*}(t)-S+A}{A}, \quad \big(\chi_{\mathscr{A}_{1*}}(t),\chi_{\mathscr{A}_{2*}}(t)\big)= (1,0).
\end{eqnarray*}
F"ur die bereits verwendeten Parameter $a=1$, $A=10$, $x_1^0,x_2^0 \in [0,10]$, $(x_1^f,x_2^f)=(5,5)$
und f"ur die untere Schranke $S=3,5$ ergibt sich das folgende Phasenportrait:
\begin{figure}[h]
	\centering
	\fbox{\includegraphics[width=4.3cm]{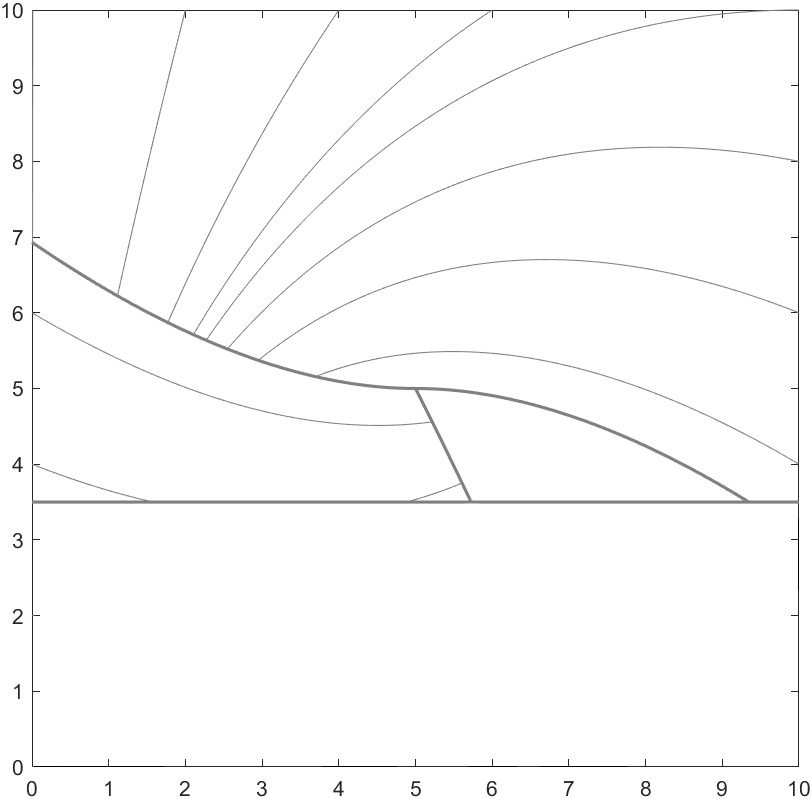}}
	\caption[Gekoppelte Kompartimente - Phasendiagramm bei Zustandsbeschränkung]{Phasenportrait unter der Zustandsbeschr"ankung.}
	\label{AbbildungBereiche2}
\end{figure}

Mit der Betrachtung der Zustandsbeschränkung ist unsere Diskussion der gekoppelten Kompartimente abgeschlossen. \hfill $\square$

%% file: 3-40-RetardierteSysteme.tex
\subsection{Zeitverz\"ogerte dynamische Systeme} \label{KapitelRetardierteSysteme}\index{Zeitverzögerte Systeme}
Dynamische Systeme mit Zeitverzögerungen sind wesentliche Elemente in der Modellierung und Auswertung existenzieller
und lebensnaher Phänomene in Anwendungsfeldern wie der Biologie oder der Biomedizin. 
So zählt zu den gängigen Beispielen einer retardierten Arznei die Form,
bei der der Wirkstoff zeitlich verzögert freigesetzt wird. \\[2mm]
Das verzögerte System entsteht dadurch,
dass die Zustandsänderung $\dot{x}(t)$ nicht ausschließlich auf den vorliegenden Zustand $x(t)$ und
die simultanen Ma"snahmen $u(t)$ interagiert,
sondern mit zeitlichen Verzögerungen $\delta_x,\,\delta_u >0$ auf den vorherigen Status $x(t-\delta_x)$ und auf die getroffenen
Maßnahmen $u(t-\delta_u)$ reagiert:
$$\dot{x}(t) = \varphi\big(t,x(t),x(t-\delta_x),u(t),u(t-\delta_u)\big).$$
Für die Erweiterung der Standardaufgabe um zeitverzögerte Systeme werden wir eine rekursive Herangehensweise vorstellen.
Dabei greifen wir die Annahme kommensurabler Gr"o"sen $\delta_x$ und $\delta_u$ auf,
die z.\,B. bei Göllmann \& Maurer \cite{GoMa,GoMa2} Anwendung findet.
Die erzielten Ergebnisse in Form des Pontrjaginschen Maximumprinzips spiegeln dabei die Komplexität der herausfordernden
Aufgabenstellung wider.
So ist die Auswertung verzögerter Steuerungsprobleme und die Anwendung von Optimalitätsbedingungen im Wesentlichen nur mit
Computerunterstützung möglich. \\[2mm]
Als Anwendung für ein verzögertes Steuerungsproblem geben wir das Modell einer Chemoimmuntherapie nach
Rihan\,et.\,al. \cite{Rihan} und nach Göllmann \& Maurer in \cite{GoMa} an.
Die numerischen Ergebnisse sind \cite{GoMa} entnommen.

%% file: 3-41-PMPeinfach.tex
\subsubsection{Die elementare Aufgabe mit freiem Endpunkt} \label{AbschnittPMPeinfachRet}
In der elementaren Aufgabe mit freiem rechten Endpunkt liegt nur eine Zeitverzögerung $\delta >0$ bezüglich
dem Zustand $x(t)$ und der Steuerung $u(t)$ vor. 
Damit gelangen wir zu dem zeitverzögerten Steuerungsproblem
\begin{eqnarray}
&& \label{Ret1} \hspace*{-5mm} J\big(x(\cdot),u(\cdot)\big) = \int_{t_0}^{t_1} f\big(t,x(t),x(t-\delta),u(t),u(t-\delta)\big) \, dt\to \inf, \\
&& \label{Ret2} \hspace*{-5mm} \dot{x}(t) = \varphi\big(t,x(t),x(t-\delta),u(t),u(t-\delta)\big),\quad x(t)=x_0 \mbox{ für } t \in [t_0-\delta,t_0], \\
&& \label{Ret3} \hspace*{-5mm} u(t) \in U \subseteq \R^m, \quad u(t)=u_0 \in U \mbox{ für } t \in [t_0-\delta,t_0),\quad U\not= \emptyset.
\end{eqnarray} 
Im Vergleich zur Aufgabe (\ref{PMPeinfach1})--(\ref{PMPeinfach3}) ergänzen wir in den Abbildungen $f$ und $\varphi$
die Variablen $y$ und $v$ bez"uglich des verzögerten Zustandes bzw. der verzögerten Steuerung:
\begin{eqnarray*}
&& f=f(t,x,y,u,v):\R \times \R^n \times \R^n \times \R^m \times \R^m \to \R, \\
&& \varphi=\varphi(t,x,y,u,v):\R \times \R^n \times \R^n \times \R^m \times \R^m \to \R^n.
\end{eqnarray*}
\newpage

Bei der Behandlung von Steuerungsproblemen mit einer Zeitverzögerung stellt sich die Frage,
wie das dynamische System (\ref{Ret2})
zu behandeln ist.
Für den Teilabschnitt $[t_0,t_0+\delta]$ reduziert es sich auf die Gleichung
$$\dot{x}(t) = \varphi\big(t,x(t),x_0,u(t),u_0\big), \quad x(t_0)=x_0,$$
welche sich für eine stückweise stetige rechte Seite $\varphi$ mit den Werkzeugen im Anhang über Differentialgleichungen handhaben lässt.
Dieser Gedankengang kann anschließend schrittweise über den Teilabschnitten $[t_0+\delta,t_0+2\delta]$, $[t_0+2\delta,t_0+3\delta]$, ... wiederholt werden.
Daher kann die Dynamik mit Zeitverzögerung abschnittsweise wie eine Differentialgleichung mit stückweise stetigen rechten Seiten
aufgefasst werden. \\[2mm]
Die Aufgabe (\ref{Ret1})--(\ref{Ret3}) betrachten wir bez"uglich der Paare
$$\big(x(\cdot),u(\cdot)\big) \in PC_1([t_0,t_1],\R^n) \times PC([t_0,t_1],U).$$
Mit $\mathscr{D}^{\,\mathcal{Z}}_{\rm Lip}$ bezeichnen wir die Menge aller Paare $\big(x(\cdot),u(\cdot)\big)$,
für die es ein $\gamma>0$ derart gibt,
dass die Abbildungen $f(t,x,y,u,v)$, $\varphi(t,x,y,u,v)$ auf der Menge aller Punkte
$(t,x,y,u,v) \in \R \times \R^n \times \R^n \times \R^m \times \R^m$ mit
$$t \in [t_0,t_1], \qquad \|x-x(t)\| < \gamma, \qquad \|y-x(t-\delta)\| < \gamma, \qquad u,v \in \R^m$$
stetig in der Gesamtheit aller Variablen und stetig differenzierbar bezüglich $x$ und $y$ sind. \\[2mm]
Das Paar $\big(x(\cdot),u(\cdot)\big) \in PC_1([t_0,t_1],\R^n) \times PC([t_0,t_1],U)$
heißt ein zul"assiger Steuerungsprozeß der Aufgabe (\ref{Ret1})--(\ref{Ret3}),
falls $\big(x(\cdot),u(\cdot)\big)$ dem System (\ref{Ret2}) gen"ugt und die Steuerbeschränkungen (\ref{Ret3}) erfüllt.
Mit $\mathscr{D}^{\,\mathcal{Z}}_{\rm adm}$ bezeichnen wir die Menge der zul"assigen Steuerungsprozesse. \\[2mm]
Ein zul"assiger Steuerungsprozess $\big(x_*(\cdot),u_*(\cdot)\big)$ ist eine
starke lokale Minimalstelle\index{Minimum, starkes lokales!ZDelay@-- zeitverzögerte Systeme}
der Aufgabe (\ref{Ret1})--(\ref{Ret3}),
falls eine Zahl $\varepsilon > 0$ derart existiert, dass die Ungleichung 
$$J\big(x(\cdot),u(\cdot)\big) \geq J\big(x_*(\cdot),u_*(\cdot)\big)$$
f"ur alle $\big(x(\cdot),u(\cdot)\big) \in \mathscr{D}^{\,\mathcal{Z}}_{\rm adm}$ mit $\|x(\cdot)-x_*(\cdot)\|_\infty < \varepsilon$ gilt. \\[2mm]
Es bezeichnet $H^{\mathcal{Z}}: \R \times \R^n \times \R^n \times \R^m \times \R^m \times \R^n \times \R \to \R$ die Pontrjagin-Funktion
$$H^{\mathcal{Z}}(t,x,y,u,v,p,\lambda_0) = \langle p,\varphi(t,x,y,u,v)\rangle - \lambda_0 f(t,x,y,u,v).$$
Weiterhin führen wir die folgenden abkürzenden Bezeichnungen ein:
\begin{eqnarray*}
H^{\mathcal{Z}}[t] &=& H^{\mathcal{Z}}\big(t,x_*(t),x_*(t-\delta), u_*(t),u_*(t-\delta),p(t),1\big), \\
H^{\mathcal{Z}}[t,u,v] &=& H^{\mathcal{Z}}\big(t,x_*(t),x_*(t-\delta), u,v,p(t),1\big).
\end{eqnarray*}
Im Gegensatz zu später folgenden Abschnitten beziehen sich diese Abkürzungen auf den normalen Fall $\lambda_0=1$. \\[2mm]
Damit formulieren wir das Pontrjaginsche Maximumprinzip der Aufgabe (\ref{Ret1})--(\ref{Ret3}): 
\newpage
\begin{theorem}[Pontrjaginsches Maximumprinzip] \label{SatzRetPMPeinfach}
\index{Pontrjaginsches Maximumprinzip!ZDelay@-- zeitverzögerte Systeme} 
Es sei $\big(x_*(\cdot),u_*(\cdot)\big) \in \mathscr{D}^{\,\mathcal{Z}}_{\rm adm} \cap \mathscr{D}^{\,\mathcal{Z}}_{\rm Lip}$.
Ist $\big(x_*(\cdot),u_*(\cdot)\big)$ ein starkes lokales Minimum der Aufgabe (\ref{Ret1})--(\ref{Ret3}),
dann existiert eine Vektorfunktion $p(\cdot) \in PC_1([t_0,t_1],\R^n)$ derart, dass
\begin{enumerate}
\item[(a)] fast überall in $[t_0,t_1]$ die adjungierte Gleichung
           \index{adjungierte Gleichung!ZDelay@-- zeitverzögerte Systeme}
           \begin{equation}\label{SatzRetPMPeinfach1}
           \dot{p}(t) = - H^{\mathcal{Z}}_x[t] - \chi_{[t_0,t_1-\delta]}(t) \cdot H^{\mathcal{Z}}_y[t+\delta],
           \end{equation}
\item[(b)] in $t =t_1$ die Transversalitätsbedingung
           \index{Transversalitätsbedingungen!ZDelay@-- zeitverzögerte Systeme}
           \begin{equation}\label{SatzRetPMPeinfach2} 
           p(t_1)=0
           \end{equation}
\item[(c)] und in fast allen Punkten $t \in [t_0,t_1]$ die Maximumbedingung
           \index{Maximumbedingung!ZDelay@-- zeitverzögerte Systeme}
           \begin{eqnarray}
           && \hspace*{-10mm}
           H^{\mathcal{Z}}[t,u_*(t),u_*(t-\delta)]+\chi_{[t_0,t_1-\delta]}(t) \cdot H^{\mathcal{Z}}[t+\delta,u_*(t+\delta),u_*(t)] \nonumber \\
           \label{SatzRetPMPeinfach3}
           && \hspace*{-10mm}
              = \max_{u \in U} \Big\{H^{\mathcal{Z}}[t,u,u_*(t-\delta)]+\chi_{[t_0,t_1-\delta]}(t) \cdot H^{\mathcal{Z}}[t+\delta,u_*(t+\delta),u] \Big\}
           \end{eqnarray}
\end{enumerate}
erfüllt sind.
\end{theorem}

In Theorem \ref{SatzRetPMPeinfach} sind die Bedingungen kompakt formuliert.
Die adjungierte Gleichung (\ref{SatzRetPMPeinfach1}) besitzt ausführlich aufgeschrieben die Gestalt
\begin{eqnarray*}
\dot{p}(t) &=& - \varphi_x^T\big(t,x_*(t),x_*(t-\delta), u_*(t),u_*(t-\delta)\big)p(t) \\
           & & \hspace*{2cm} +f_x\big(t,x_*(t),x_*(t-\delta), u_*(t),u_*(t-\delta)\big) \\
           & & + \chi_{[t_0,t_1-\delta]}(t) \cdot
                 \Big[ - \varphi_y^T\big(t+\delta,x_*(t+\delta),x_*(t), u_*(t+\delta),u_*(t)\big)p(t+\delta) \\
           & &\hspace*{4.5cm}  +f_y\big(t+\delta,x_*(t+\delta),x_*(t), u_*(t+\delta),u_*(t)\big)\Big].
\end{eqnarray*}
Ähnlich wie das System (\ref{Ret2}) besitzt die adjungierte Gleichung (\ref{SatzRetPMPeinfach1}) über $[t_1-\delta,t_1]$ die
``gewöhnliche'' Form $\dot{p}(t) = - H^{\mathcal{Z}}_x[t]$ zur Transversalitätsbedingung (\ref{SatzRetPMPeinfach2}) 
und lässt sich dann abschnittsweise über $[t_1-2\delta,t_1-\delta]$, $[t_1-3\delta,t_1-2\delta]$, ... behandeln.
Ferner hat die Maximumbedingung (\ref{SatzRetPMPeinfach3}) ausführlich die Form
\begin{eqnarray*}
& & \big\langle \varphi\big(t,x_*(t),x_*(t-\delta), u_*(t),u_*(t-\delta)\big),p(t) \big\rangle \\
& & \hspace*{2cm} f\big(t,x_*(t),x_*(t-\delta), u_*(t),u_*(t-\delta)\big) \\
& & + \chi_{[t_0,t_1-\delta]}(t) \cdot
      \Big[ \big\langle \varphi\big(t+\delta,x_*(t+\delta),x_*(t), u_*(t+\delta),u_*(t)\big) , p(t+\delta)\big\rangle \\
& &\hspace*{4.5cm}  -f\big(t+\delta,x_*(t+\delta),x_*(t), u_*(t+\delta),u_*(t)\big)\Big] \\
&=& \max_{u \in U}\Big\{ \big\langle \varphi\big(t,x_*(t),x_*(t-\delta), u,u_*(t-\delta)\big),p(t) \big\rangle \\
& & \hspace*{3cm} f\big(t,x_*(t),x_*(t-\delta), u,u_*(t-\delta)\big) \\
& & \hspace*{1cm}+ \chi_{[t_0,t_1-\delta]}(t) \cdot
      \Big[ \big\langle \varphi\big(t+\delta,x_*(t+\delta),x_*(t), u_*(t+\delta),u\big) , p(t+\delta)\big\rangle \\
& &\hspace*{5.5cm}  -f\big(t+\delta,x_*(t+\delta),x_*(t), u_*(t+\delta),u\big)\Big]\Big\}.
\end{eqnarray*}

{\bf Beweis}
Im Zuge der Herleitung des Pontrjaginschen Maximumprinzips für die elementare Aufgabe
gelangen wir zu Ausdrücken, welche Verzögerungen in verschiedene Zeitrichtungen aufweisen.
Mit der Substitution $t = s+ \delta$ gilt $f(t)g(t-\delta)=f(s+\delta)g(s)$ und es ergibt sich (nach Umbenennung von $s$ zu $t$) die Beziehung
\begin{equation}\label{Zeitshift}
\int_a^b \chi_{[c+\delta,d]}(t) \cdot f(t)g(t-\delta) \, dt = \int_a^b \chi_{[c,d-\delta]}(t) \cdot f(t+\delta)g(t) \, dt.
\end{equation}

Beachten wir die oben getroffenen Bemerkungen zur Lösung von Differentialgleichungen mit Zeitverzögerungen,
so gibt es nach Lemma \ref{LemmaDGL3} und Lemma \ref{LemmaDGL5} eine eindeutige L"osung $p(\cdot) \in PC_1([t_0,t_1],\R^n)$
der Gleichung (\ref{SatzRetPMPeinfach1}) zur Randbedingung (\ref{SatzRetPMPeinfach2}). \\[1mm]
Für $\lambda < \delta$ betrachten wir wieder die einfache Nadelvariation\index{Nadelvariation, einfache}
$$u(t;v,\tau,\lambda) = u_{\lambda}(t) = 
  \left\{ \begin{array}{ll}
          u_*(t) & \mbox{ f"ur } t \not\in [\tau-\lambda,\tau), \\
          v      & \mbox{ f"ur } t     \in [\tau-\lambda,\tau), 
          \end{array} \right.$$
und untersuchen f"ur $t \geq \tau$ den Grenzwert (vgl. Abschnitt \ref{AbschnittPMPBeweiseinfach})
\begin{equation} \label{BeweisRet1}
y(t)=\lim_{\lambda \to 0^+}\frac{x_{\lambda}(t) - x_*(t)}{\lambda}.
\end{equation}
Im weiteren Vorgehen müssen wir beachten,
dass in der Aufgabe (\ref{Ret1})--(\ref{Ret3}) die Nadelvariation $u_\lambda(\cdot)$ über den beiden Intervall
$[\tau-\lambda,\tau)$ und $[\tau+\delta-\lambda,\tau+\delta)$ in die Abbildungen $f$ und $\varphi$ einfließt.
Für den Quotienten auf der rechten Seite in (\ref{BeweisRet1}) ergeben sich für $t \geq \tau$ abschnittsweise
über dem Intervall $[\tau-\lambda,t_1]$ die folgenden Ausdrücke:
\begin{eqnarray*}
&& \frac{1}{\lambda} \int_{\tau - \lambda}^{\tau} \chi_{[t_0,t]}(s) \cdot
   \Big[ \varphi\big(s,x_\lambda(s),x_*(s-\delta),v,u_*(s-\delta)\big) \nonumber\\
&& \hspace*{4cm} - \varphi\big(s,x_*(s),x_*(s-\delta),u_*(s),u_*(s-\delta)\big) \Big] \, ds \nonumber\\
&& + \frac{1}{\lambda} \int_{\tau}^{\tau+\delta - \lambda} \chi_{[t_0,t]}(s) \cdot
   \Big[ \varphi\big(s,x_\lambda(s),x_*(s-\delta),u_*(s),u_*(s-\delta)\big) \nonumber\\
&& \hspace*{4cm} - \varphi\big(s,x_*(s),x_*(s-\delta),u_*(s),u_*(s-\delta)\big) \Big] \, ds \\
&& + \frac{1}{\lambda} \int_{\tau+\delta - \lambda}^{\tau+\delta} \chi_{[t_0,t]}(s) \cdot
   \Big[ \varphi\big(s,x_\lambda(s),x_*(s-\delta),u_*(s),v\big) \nonumber\\
&& \hspace*{4cm} - \varphi\big(s,x_*(s),x_*(s-\delta),u_*(s),u_*(s-\delta)\big) \Big] \, ds \\
&& + \frac{1}{\lambda} \int_{\tau+\delta}^{t_1} \chi_{[t_0,t]}(s) \cdot
   \Big[ \varphi\big(s,x_\lambda(s),x_\lambda(s-\delta),u_*(s),u_*(s-\delta)\big) \nonumber\\
&& \hspace*{4cm} - \varphi\big(s,x_*(s),x_*(s-\delta),u_*(s),u_*(s-\delta)\big) \Big] \, ds. \nonumber
\end{eqnarray*}
Wir erweitern diese Summe um die Terme
$$\pm \frac{1}{\lambda} \int_{\tau+\delta - \lambda}^{\tau+\delta} \chi_{[t_0,t]}(s) \cdot
  \varphi\big(s,x_\lambda(s),x_*(s-\delta),u_*(s),u_*(s-\delta)\big)\, ds,$$
formen damit den zweiten und dritten Summanden um und bringen die gesamte Summe in die Gestalt
\begin{eqnarray*}
&& \frac{1}{\lambda} \int_{\tau - \lambda}^{\tau} \chi_{[t_0,t]}(s) \cdot
   \Big[ \varphi\big(s,x_\lambda(s),x_*(s-\delta),v,u_*(s-\delta)\big) \nonumber\\
&& \hspace*{4cm} - \varphi\big(s,x_*(s),x_*(s-\delta),u_*(s),u_*(s-\delta)\big) \Big] \, ds \nonumber\\
&& + \frac{1}{\lambda} \int_{\tau}^{\tau+\delta} \chi_{[t_0,t]}(s) \cdot
   \Big[ \varphi\big(s,x_\lambda(s),x_*(s-\delta),u_*(s),u_*(s-\delta)\big) \nonumber\\
&& \hspace*{4cm} - \varphi\big(s,x_*(s),x_*(s-\delta),u_*(s),u_*(s-\delta)\big) \Big] \, ds \\
&& + \frac{1}{\lambda} \int_{\tau+\delta - \lambda}^{\tau+\delta} \chi_{[t_0,t]}(s) \cdot
   \Big[ \varphi\big(s,x_\lambda(s),x_*(s-\delta),u_*(s),v\big) \nonumber\\
&& \hspace*{4cm} - \varphi\big(s,x_\lambda(s),x_*(s-\delta),u_*(s),u_*(s-\delta)\big) \Big] \, ds \nonumber\\
&& + \frac{1}{\lambda} \int_{\tau+\delta}^{t_1} \chi_{[t_0,t]}(s) \cdot
   \Big[ \varphi\big(s,x_\lambda(s),x_\lambda(s-\delta),u_*(s),u_*(s-\delta)\big) \nonumber\\
&& \hspace*{4cm} - \varphi\big(s,x_*(s),x_*(s-\delta),u_*(s),u_*(s-\delta)\big) \Big] \, ds.
\end{eqnarray*}
Nun ergänzen wir in dieser Summe die Terme
$$\pm \frac{1}{\lambda} \int_{\tau+\delta}^{t_1} \chi_{[t_0,t]}(s) \cdot
  \varphi\big(s,x_\lambda(s),x_*(s-\delta),u_*(s),u_*(s-\delta)\big)\, ds,$$
formen damit den zweiten und vierten Summanden um und bringen die Summe in diejenige Form,
für die wir den Grenzübergang $\lambda \to 0^+$ bilden wollen:
\begin{eqnarray}
&& \frac{1}{\lambda} \int_{\tau - \lambda}^{\tau}\chi_{[t_0,t]}(s) \cdot
   \Big[ \varphi\big(s,x_\lambda(s),x_*(s-\delta),v,u_*(s-\delta)\big) \nonumber\\
&& \hspace*{4cm} - \varphi\big(s,x_*(s),x_*(s-\delta),u_*(s),u_*(s-\delta)\big) \Big] \, ds  \nonumber\\
&& + \frac{1}{\lambda} \int_{\tau}^{t_1} \chi_{[t_0,t]}(s) \cdot
   \Big[ \varphi\big(s,x_\lambda(s),x_*(s-\delta),u_*(s),u_*(s-\delta)\big) \nonumber\\
&& \hspace*{4cm} - \varphi\big(s,x_*(s),x_*(s-\delta),u_*(s),u_*(s-\delta)\big) \Big] \, ds  \nonumber\\
&& + \frac{1}{\lambda} \int_{\tau+\delta - \lambda}^{\tau+\delta} \chi_{[t_0,t]}(s) \cdot
   \Big[ \varphi\big(s,x_\lambda(s),x_*(s-\delta),u_*(s),v\big) \nonumber\\
&& \hspace*{4cm} - \varphi\big(s,x_\lambda(s),x_*(s-\delta),u_*(s),u_*(s-\delta)\big) \Big] \, ds  \nonumber\\
&& + \frac{1}{\lambda} \int_{\tau+\delta}^{t_1} \chi_{[t_0,t]}(s) \cdot
   \Big[ \varphi\big(s,x_\lambda(s),x_\lambda(s-\delta),u_*(s),u_*(s-\delta)\big) \nonumber\\
&& \label{BeweisRet2} \hspace*{4cm} - \varphi\big(s,x_\lambda(s),x_*(s-\delta),u_*(s),u_*(s-\delta)\big) \Big] \, ds.
\end{eqnarray}
Um uns von den klobigen Ausdrücken der Art (\ref{BeweisRet2}) zu lösen,
verwenden wir im Weiteren die Kurzschreibweisen
\begin{eqnarray*}
\varphi[t] &=& \varphi\big(t,x_*(t),x_*(t-\delta), u_*(t),u_*(t-\delta)\big), \\
\varphi[t,u,v] &=& \varphi\big(t,x_*(t),x_*(t-\delta), u,v\big), \\
f[t] &=& f\big(t,x_*(t),x_*(t-\delta), u_*(t),u_*(t-\delta)\big), \\
f[t,u,v] &=& f\big(t,x_*(t),x_*(t-\delta), u,v\big).
\end{eqnarray*}
Mit diesen Schreibweisen ergibt der Grenzübergang $\lambda \to 0^+$ in (\ref{BeweisRet2}):
\begin{eqnarray*}
y(t)&=& \underbrace{\varphi[\tau,v,u_*(\tau-\delta)] - \varphi[\tau,u_*(\tau),u_*(\tau-\delta)]}_{=y(\tau)} 
        + \int_{\tau}^{t_1} \chi_{[t_0,t]}(s) \cdot \varphi_x[s] y(s) \,ds \\
    & & + \chi_{[t_0,t]}(\tau+\delta) \cdot
        \Big[ \underbrace{\varphi[\tau+\delta,u_*(\tau+\delta),v] - \varphi[\tau+\delta,u_*(\tau+\delta),u_*(\tau)]}_{=\eta(\tau+\delta)} \Big] \\
    & & + \int_{\tau+\delta}^{t_1}\chi_{[t_0,t]}(s) \cdot \varphi_y[s] y(s-\delta) \,ds.
\end{eqnarray*}
Das erste Integral beschränken wir auf das Intervall $[\tau,t]$ und entfernen dafür die charakteristische Funktion im Integranden.
Im zweiten Integral stimmt der Geltungsbereich $[\tau+\delta,t_1] \cap [t_0,t]$ mit $[\tau,t] \cap [\tau+\delta,t_1]$ überein.
Ferner dürfen wir $\chi_{[t_0,t- \delta]}(\tau)$ anstelle von $\chi_{[t_0,t]}(\tau+\delta)$ schreiben.
Es ergibt sich
\begin{eqnarray}
y(t) &=& y(\tau) +  \chi_{[t_0,t-\delta]}(\tau) \cdot \eta(\tau+\delta) \nonumber\\
     & & \label{BeweisRet21} + \int_{\tau}^t \varphi_x[s] y(s) \,ds + \int_{\tau}^t \chi_{[\tau+\delta,t_1]}(s) \cdot \varphi_y[s] y(s-\delta) \,ds. 
\end{eqnarray}
Es kann $y(t)$ in $t=\tau+\delta \in (t_0,t_1)$ unstetig sein und dann gilt
$$\lim_{\varepsilon \to 0^+} [y(\tau+\delta+\varepsilon) -y(\tau+\delta-\varepsilon)] =\eta(\tau+\delta).$$
Wir werten nun die Ableitung
$\displaystyle \frac{d}{dt} \langle p(t), y(t) \rangle = \langle \dot{p}(t), y(t) \rangle + \langle p(t), \dot{y}(t) \rangle$
aus.
Es gelten für $y(\cdot)$ und nach (\ref{SatzRetPMPeinfach1}) für die Adjungierte $p(\cdot)$ in $[\tau,t_1]$:
\begin{eqnarray*}
\dot{y}(t) &=& \varphi_x[t] y(t) + \chi_{[\tau+\delta,t_1]}(t) \cdot \varphi_y[t] y(t-\delta), \\
\dot{p}(t) &=& - \varphi^T_x[t] p(t)+f_x[t] + \chi_{[\tau,t_1-\delta]}(t) \cdot \big( - \varphi^T_y[t+\delta] p(t+\delta)+f_y[t+\delta] \big).
\end{eqnarray*}
Für die $\varphi$-Terme mit Zeitverschiebung $\pm\delta$ ergibt sich mit Hilfe von (\ref{Zeitshift})
$$\int_\tau^{t_1} \chi_{[\tau+\delta,t_1]}(t) \cdot \big\langle p(t) , \varphi_y[t] y(t-\delta) \big\rangle \, dt
   =\int_\tau^{t_1} \chi_{[\tau,t_1-\delta]}(t) \cdot \big\langle \varphi^T_y[t+\delta] p(t+\delta) , y(t) \big\rangle \, dt.$$
Da außerdem $\langle - \varphi^T_x[t] p(t), y(t) \rangle = - \langle p(t), \varphi_x[t] y(t) \rangle$ gilt, erhalten wir
\begin{eqnarray*}
    \langle p(t_1), y(t_1) \rangle
&=& \langle p(\tau), y(\tau) \rangle + \chi_{[t_0,t_1-\delta]}(\tau) \cdot \langle p(\tau+\delta), \eta(\tau+\delta) \rangle
     + \int_\tau^{t_1} \frac{d}{dt} \langle p(t), y(t) \rangle \, dt \\
&=& \langle p(\tau), y(\tau) \rangle + \chi_{[t_0,t_1-\delta]}(\tau) \cdot \langle p(\tau+\delta), \eta(\tau+\delta) \rangle \\
& & \hspace*{1cm} + \int_{\tau}^{t_1} \langle f_x[t] , y(t) \rangle \, dt 
    + \int_{\tau}^{t_1} \chi_{[\tau,t_1-\delta]}(t) \cdot \langle f_y[t+\delta], y(t) \rangle \, dt.
\end{eqnarray*}
Im zweiten Integral ebenso die Umformung (\ref{Zeitshift}) angewendet liefert dann
\begin{eqnarray}
    \hspace*{-5mm} \langle p(t_1), y(t_1) \rangle
&=& \langle p(\tau), y(\tau) \rangle + \chi_{[t_0,t_1-\delta]}(\tau) \cdot \langle p(\tau+\delta), \eta(\tau+\delta) \rangle \nonumber \\
\hspace*{-5mm} & & \label{BeweisRet3} + \int_{\tau}^{t_1} \langle f_x[t] , y(t) \rangle \, dt 
    + \int_{\tau}^{t_1} \chi_{[\tau+\delta,t_1]}(t) \cdot \langle f_y[t], y(t-\delta) \rangle \, dt, \\
\hspace*{-5mm} y(\tau) &=& \label{BeweisRet4} \varphi[\tau,v,u_*(\tau-\delta)] - \varphi[\tau,u_*(\tau),u_*(\tau-\delta)], \\
\hspace*{-5mm} \eta(\tau+\delta) &=& \label{BeweisRet5} \varphi[\tau+\delta,u_*(\tau+\delta),v] - \varphi[\tau+\delta,u_*(\tau+\delta),u_*(\tau)].
\end{eqnarray}
Die Auswertung der dynamischen Wirkung der Nadelvariation ist (endlich) abgeschlossen und wir kommen zum Abschluss des Beweises:
Da $\big(x_*(\cdot),u_*(\cdot)\big)$ ein starkes lokales Minimum ist, gilt für alle hinreichend kleine $\lambda >0$ die Ungleichung
$$0 \leq \frac{J\big(x_\lambda(\cdot),u_\lambda(\cdot)\big)- J\big(x_*(\cdot),u_*(\cdot)\big)}{\lambda}.$$
Den Quotienten auf der rechten Seite der Ungleichung formen wir auf die gleiche Weise wie zur Herleitung von (\ref{BeweisRet2}) bzw. (\ref{BeweisRet21}) um
und erhalten im Grenzübergang $\lambda \to 0^+$:
\begin{eqnarray}
\hspace*{-10mm} 0 &\leq& f[\tau,v,u_*(\tau-\delta)] -f[\tau,u_*(\tau),u_*(\tau-\delta)]  \nonumber\\
\hspace*{-10mm} & & + \chi_{[t_0,t_1-\delta]}(\tau) \cdot \Big[ f[\tau+\delta,u_*(\tau+\delta),v] - f[\tau+\delta,u_*(\tau+\delta),u_*(\tau)] \Big]  \nonumber\\
\hspace*{-10mm} & & \label{BeweisRet6}
                   + \int_{\tau}^{t_1} \langle f_x[t] , y(t) \rangle \,dt
                   + \int_{\tau}^{t_1} \chi_{[\tau+\delta,t_1]}(t) \cdot \langle f_y[t] , y(t-\delta) \rangle \,dt. 
\end{eqnarray}
Die Bedingungen (\ref{BeweisRet3})--(\ref{BeweisRet6}) und die Transversalitätsbedingung (\ref{SatzRetPMPeinfach1}) ergeben 
\begin{eqnarray}
0 &\leq& f[\tau,v,u_*(\tau-\delta)] -f[\tau,u_*(\tau),u_*(\tau-\delta)]  \nonumber\\
    & & + \chi_{[t_0,t_1-\delta]}(\tau) \cdot\Big[ f[\tau+\delta,u_*(\tau+\delta),v] - f[\tau+\delta,u_*(\tau+\delta),u_*(\tau)] \Big]  \nonumber\\
    & & \label{BeweisRet7} - \langle p(\tau), y(\tau) \rangle - \chi_{[t_0,t_1-\delta]}(\tau)\cdot\langle p(\tau+\delta), \eta(\tau+\delta) \rangle. 
\end{eqnarray}
Da $\tau \in (t_0,t_1)$ als ein Stetigkeitspunkt der Steuerung $u_*(\cdot)$ und zudem $v \in U$ beliebig gewählt waren,
liefert (\ref{BeweisRet6}) die Gültigkeit der Maximumbedingung (\ref{SatzRetPMPeinfach3}).
Der Beweis von Theorem \ref{SatzRetPMPeinfach} ist damit abgeschlosssen. \hfill $\blacksquare$

\begin{beispiel} \label{BeispielInvDelay}
{\rm In Anlehnung an Beispiel \ref{BeispielLinInv} betrachten wir eine Aufgabe,
in welcher sich die Investitionen mit der Verzögerung von einer Zeiteinheit $\delta = 1$ auf den Kapitalbestand auswirken.
Es ergibt sich das Investitionsmodell
\begin{eqnarray}
&&\label{BeispielRet1} J\big(K(\cdot),u(\cdot)\big) = \int_0^T \big( 1-u(t)\big) \cdot K(t) \, dt \to \sup, \\
&&\label{BeispielRet2} \dot{K}(t) = u(t-1) \cdot K(t-1), \quad K(0)=K_0 >0,\\
&&\label{BeispielRet3} u(t) \in [0,1] \mbox{ für } t \in [0,T], \quad u(t)=0 \mbox{ für } t \in [-1,0), \quad T > 2 \mbox{ fest}.
\end{eqnarray}
In der Aufgabe (\ref{BeispielRet1})--(\ref{BeispielRet3}) verwenden wir für $K(t),K(t-1),u(t),u(t-1)$ die Variablenbezeichnungen $K,L,u,v$.
Damit lauten die Abbildungen der Aufgabe
$$f(t,K,L,u,v) = -( 1-u) \cdot K, \quad \varphi(t,K,L,u,v) = v \cdot L$$
und die Pontrjaginsche Funktion
$$H^{\mathcal{Z}}(t,K,L,u,v,p) = p vL +(1-u)K.$$
Die Maximumbedingung (\ref{SatzRetPMPeinfach3}) besitzt die Form
\begin{eqnarray*}
\lefteqn{\max_{u \in U} \Big\{H^{\mathcal{Z}}[t,u,u_*(t-1)]+\chi_{[0,T-1]}(t) \cdot H^{\mathcal{Z}}[t+1,u_*(t+1),u] \Big\}} \\
&=& \max_{u \in [0,1]}\Big\{p(t)u_*(t-1)K_*(t-1) +(1-u)K_*(t) \\
& & \hspace*{1.5cm} + \chi_{[0,T-1]}(t)\cdot\big[p(t+1)u K_*(t) +\big(1-u_*(t+1)\big)K_*(t+1)\Big\}
\end{eqnarray*}
und ist äquivalent zu folgenden Maximierungen über Teilabschnitte:
\begin{equation} \label{BeispielRet4}
\max_{u \in [0,1]} \left\{\begin{array}{ll}
u \big(p(t+1)-1\big) K_*(t), & t \in [0,T-1], \\[1mm] 
-u K_*(t),& t \in [T-1,T].
\end{array} \right.
\end{equation}
Die adjungierte Gleichung (\ref{SatzRetPMPeinfach1}) erhält die Gestalt
$$\dot{p}(t) = - H^{\mathcal{Z}}_K[t] - \chi_{[0,T-1]}(t) \cdot H^{\mathcal{Z}}_L[t+1] = -\big(1-u_*(t)\big) - \chi_{[0,T-1]}(t)\cdot p(t+1)u_*(t)$$
und lautet abschnittsweise
\begin{equation} \label{BeispielRet5}
\dot{p}(t)=\left\{\begin{array}{ll} 
\displaystyle -\big(1-u_*(t)\big) - p(t+1)u_*(t), & t \in (0,T-1), \\[1mm]
\displaystyle -\big(1-u_*(t)\big),& t \in (T-1,T),
\end{array} \right.
\end{equation}
zur Transversalitätsbedingung $p(T)=0$. \\[2mm]
Da $\dot{K}(t) \geq 0$ und $\dot{p}(t)<0$ in $(0,T)$ gelten, führen (\ref{BeispielRet4}) und (\ref{BeispielRet5}) zu:
\begin{enumerate}
\item[$\cdot$] Über $[T-1,T]$ lauten (\ref{BeispielRet4}) und (\ref{BeispielRet5})
               $$\max_{u \in [0,1]} -u K_*(t), \quad \dot{p}(t)=-\big(1-u_*(t)\big), \quad p(T)=0.$$
               Daraus ist unmittelbar $u_*(t)=0$ und $p(t)=T-t$ über $[T-1,T]$ einzusehen.
\item[$\cdot$] Damit folgt für (\ref{BeispielRet4}) über $[T-2,T-1)$:
               $\displaystyle \max_{u \in [0,1]} u \big(p(t+1)-1\big) K_*(t)$.
               Da $p(t)<1$ für $t \in (T-1,T]$ ausfällt, ist $u_*(t)=0$ für $t \in [T-2,T-1)$.
               In (\ref{BeispielRet5}) ergibt sich ferner $\dot{p}(t)=-1$ und es gilt weiterhin $p(t)=T-t$ für $t \in [T-2,T-1)$.
\item[$\cdot$] Da stets $\dot{p}<0$ ausfällt, ist $p(t+1)>2$ über $[0,T-2)$.
               Dies hat anhand der Maximumbedingung (\ref{BeispielRet4}) die Investitionsrate
               $u_*(t)=1$ für $t \in [0,T-2)$ zur Folge.
               Weiterhin ergibt sich $\dot{p}(t)=-p(t+1)$ in (\ref{BeispielRet5}).              
\end{enumerate}
In der Aufgabe (\ref{BeispielRet1})--(\ref{BeispielRet3}) mit verzögerten Investitionen konnten wir die Steuerung
$$u_*(t)=\left\{\begin{array}{rl} 0, & t \in [-1,0), \\ 1, & t \in [0,T-2), \\ 0, & t \in [T-2,T], \end{array} \right. \qquad
  u_*(t-1)=\left\{\begin{array}{rl} 0, & t \in [0,1), \\ 1, & t \in [1,T-1), \\ 0, & t \in [T-1,T], \end{array} \right.$$
ermitteln.
Mit dieser ergeben sich für den Kaptialbestand und für die Adjungierte die dynamischen Entwicklungen
$$\dot{K}_*(t)=\left\{\begin{array}{rl} 0, & t \in (0,1), \\ K_*(t-1), & t \in (1,T-1), \\ 0, & t \in (T-1,T), \end{array} \right. \qquad
  \dot{p}(t)=\left\{\begin{array}{rl} -p(t+1), & t \in (0,T-2), \\ -1 , & t \in (T-2,T). \end{array} \right.$$
Zusammen mit $K(0)=K_0$ und $p(T)=0$ führt dies abschnittsweise auf die Darstellungen
\begin{eqnarray*}
K_*(t) &=& \left\{\begin{array}{rl} K_0, & t \in [0,1), \\ K_0+K_0(t-1), & t \in [1,2), \\ 
                                    2K_0 + K_0(t-2)+ \frac{1}{2}K_0\big((t-1)^2-1\big), & t \in [2,3), \\ 
                                    \vdots \\ K_*(T-1), & t \in [T-1,T], \end{array} \right. \\[2mm]
p(t) &=& \left\{\begin{array}{rl} \vdots \\ 2+\frac{1}{2}\big[\big(T-(t+1)\big)^2-1\big], & t \in [T-3,T-2), \\ T-t , & t \in [T-2,T]. \end{array} \right.
\end{eqnarray*}
Aufgrund der abschnittsweise rekursiven Bildung des optimalen Kapitalbestandes $K_*(\cdot)$ verzichten wir an dieser Stelle auf die Berechnung
des Optimalwertes. \hfill $\square$}
\end{beispiel}

%% file: 3-42-Aufgabenstellung.tex
\subsubsection{Die Aufgabenstellung}
Unsere Betrachtungen schränken wir jeweils auf eine Verzögerung $\delta_x>0$ bez"uglich dem Zustand
und eine Verzögerung $\delta_u > 0$ bez"uglich der Steuerung ein.
Wir betrachten im Folgenden die allgemeine Aufgabenstellung eines verzögerten Steuerungsproblems der Form
\begin{eqnarray}
&& \label{RetPMP1} J\big(x(\cdot),u(\cdot)\big) = \int_{t_0}^{t_1} f\big(t,x(t),x(t-\delta_x),u(t),u(t-\delta_u)\big) \, dt \to \inf, \\
&& \label{RetPMP2} \dot{x}(t) = \varphi\big(t,x(t),x(t-\delta_x),u(t),u(t-\delta_u)\big), \\
&& \label{RetPMP3} x(t_0)=x_0, \quad x(t_1)=x_1,\qquad x(t)=x_0 \mbox{ für } t \in [t_0-\delta_x,t_0], \\
&& \label{RetPMP4} u(t) \in U \subseteq \R^m,\quad u(t)=u_0 \in U \mbox{ für } t \in [t_0-\delta_u,t_0), \quad U\not= \emptyset, \\
&& \label{RetPMP5} g_j\big(t,x(t),x(t-\delta_x)\big) \leq 0 \quad \mbox{f"ur alle } t \in [t_0,t_1] \mbox{ und alle } j=1,...,l.
\end{eqnarray} 
Im Vergleich zur Aufgabe (\ref{PMP1})--(\ref{PMP5}) ergänzen wir in den Abbildungen $f,\varphi$ und $g_j$
die Variablen $y$ und $v$ bez"uglich
des verzögerten Zustandes bzw. der verzögerten Steuerung:
\begin{eqnarray*}
&& f=f(t,x,y,u,v):\R \times \R^n \times \R^n \times \R^m \times \R^m \to \R, \\
&& \varphi=\varphi(t,x,y,u,v):\R \times \R^n \times \R^n \times \R^m \times \R^m \to \R^n, \\
&& g_j=g_j(t,x,y):\R \times \R^n \times \R^n \to \R^n.
\end{eqnarray*}
Wir treffen die Annahme kommensurabler Zeitgrößen,
d.\,h. es gibt ein $\Delta>0$ und natürliche Zahlen $k,\kappa,N \in \N$ mit 
\begin{equation} \label{RetPMPKommens}
k \cdot \Delta = \delta_x, \qquad \kappa \cdot \Delta = \delta_u, \qquad (N+1) \cdot \Delta = t_1-t_0.
\end{equation}
Wir nennen die Trajektorie $x(\cdot)$ eine L"osung der Gleichung (\ref{RetPMP2}) zur Steuerung $u(\cdot)$,
falls $x(\cdot)$ auf $[t_0,t_1]$ definiert ist und die Dynamik im Sinn von Carath\'eodory l"ost. \\[2mm]
Mit $\mathscr{B}^{\,\mathcal{Z}}_{\rm Lip}$ bezeichnen wir die Menge aller Paare $\big(x(\cdot),u(\cdot)\big)$,
für die es ein $\gamma>0$ derart gibt, dass die Abbildungen
$f(t,x,y,u,v)$, $\varphi(t,x,y,u,v)$ und $g(t,x,y)$ auf der Menge aller Punkte
$(t,x,y,u,v) \in \R \times \R^n \times \R^n \times \R^m \times \R^m$ mit
$$t \in [t_0,t_1], \qquad \|x-x(t)\| < \gamma, \qquad \|y-x(t-\delta)\| < \gamma, \qquad u,v \in \R^m$$
stetig in der Gesamtheit aller Variablen und stetig differenzierbar bezüglich $x$ und $y$ sind. \\[2mm]
Das Paar $\big(x(\cdot),u(\cdot)\big) \in W^1_\infty([t_0,t_1],\R^n) \times L_\infty([t_0,t_1],U)$
hei"st ein zul"assiger Steuerungsprozess in der Aufgabe (\ref{RetPMP1})--(\ref{RetPMP5}),
falls $\big(x(\cdot),u(\cdot)\big)$ dem System (\ref{RetPMP2}) gen"ugt,
sowie die Randbedingungen (\ref{RetPMP3}) und Zustandsbeschr"ankung (\ref{RetPMP5}) erf"ullt.
Die Menge $\mathscr{B}^{\,\mathcal{Z}}_{\rm adm}$ bezeichnet die Menge der zul"assigen Steuerungsprozesse. \\
Ein zul"assiger Steuerungsprozess $\big(x_*(\cdot),u_*(\cdot)\big)$ ist ein starkes lokales
Minimum\index{Minimum, starkes lokales!ZDelay@-- zeitverzögerte Systeme}
der Aufgabe (\ref{RetPMP1})--(\ref{RetPMP5}),
falls eine Zahl $\varepsilon > 0$ derart existiert, dass die Ungleichung 
$$J\big(x(\cdot),u(\cdot)\big) \geq J\big(x_*(\cdot),u_*(\cdot)\big)$$
f"ur alle $\big(x(\cdot),u(\cdot)\big) \in \mathscr{B}^{\,\mathcal{Z}}_{\rm adm}$ mit 
$\|x(\cdot)-x_*(\cdot) \|_\infty < \varepsilon$ gilt.

%% file: 3-43-Maximumprinzip.tex
\subsubsection{Das Pontrjaginsche Maximumprinzip}
Wir definieren zur Aufgabe (\ref{RetPMP1})--(\ref{RetPMP5}) die Pontrjagin-Funktion $H^{\mathcal{Z}}$ gemäß
$$H^{\mathcal{Z}}(t,x,y,u,v,p,\lambda_0) = \langle p,\varphi(t,x,y,u,v)\rangle - \lambda_0 f(t,x,y,u,v).$$
Weiterhin führen wir die folgenden abkürzenden Bezeichnungen ein:
\begin{eqnarray*}
\label{RetPontr} H^{\mathcal{Z}}[t] &=& H\big(t,x_*(t),x_*(t-\delta_x), u_*(t),u_*(t-\delta_u),p(t),\lambda_0 \big), \\
H^{\mathcal{Z}}[t,u,v] &=& H\big(t,x_*(t),x_*(t-\delta_x), u,v,p(t),\lambda_0 \big).
\end{eqnarray*}

Damit formulieren wir das Pontrjaginsche Maximumprinzip der Aufgabe (\ref{RetPMP1})--(\ref{RetPMP4}): 
\begin{theorem}[Pontrjaginsches Maximumprinzip] \label{SatzRetPMP}
\index{Pontrjaginsches Maximumprinzip!ZDelay@-- zeitverzögerte Systeme} 
Es sei $\big(x_*(\cdot),u_*(\cdot)\big) \in \mathscr{B}^{\,\mathcal{Z}}_{\rm adm} \cap \mathscr{B}^{\,\mathcal{Z}}_{\rm Lip}$
und es sei (\ref{RetPMPKommens}) erfüllt.
Ist $\big(x_*(\cdot),u_*(\cdot)\big)$ ein starkes lokales Minimum der Aufgabe (\ref{RetPMP1})--(\ref{RetPMP4}),
dann existieren nicht gleichzeitig verschwindende Multiplikatoren $\lambda_0 \geq 0$,
$p(\cdot) \in W^1_\infty([t_0,t_1],\R^n)$ und $l_0,\,l_1 \in \R^n$ derart, dass
\begin{enumerate}
\item[(a)] die adjungierte Gleichung
           \index{adjungierte Gleichung!ZDelay@-- zeitverzögerte Systeme}
           \begin{equation}\label{SatzRetPMP1}
           \dot{p}(t) = - H^{\mathcal{Z}}_x[t] - \chi_{[t_0,t_1-\delta_x]}(t) \cdot H^{\mathcal{Z}}_y[t+\delta_x]
           \end{equation}
\item[(b)] die Transversalit"atsbedingungen
           \index{Transversalitätsbedingungen!ZDelay@-- zeitverzögerte Systeme}
           \begin{equation}\label{SatzRetPMP2} 
           p(t_0)=l_0,\qquad p(t_1)=l_1
           \end{equation}
\item[(c)] und in fast allen Punkten $t \in [t_0,t_1]$ die Maximumbedingung
           \index{Maximumbedingung!ZDelay@-- zeitverzögerte Systeme}
           \begin{eqnarray}
           && \hspace*{-20mm} H^{\mathcal{Z}}[t,u_*(t),u_*(t-\delta_u)]+\chi_{[t_0,t_1-\delta_u]}(t) \cdot H^{\mathcal{Z}}[t+\delta_u,u_*(t+\delta_u),u_*(t)] \nonumber \\
           \label{SatzRetPMP3} && \hspace*{-10mm} = \max_{u \in U} \Big\{
                H^{\mathcal{Z}}[t,u,u_*(t-\delta_u)]+\chi_{[t_0,t_1-\delta_u]}(t) \cdot H^{\mathcal{Z}}[t+\delta_u,u_*(t+\delta_u),u] \Big\}
           \end{eqnarray}
\end{enumerate}
erfüllt sind.
\end{theorem}

Die komplexen Strukturen der Bedingungen (\ref{SatzRetPMP1}) und (\ref{SatzRetPMP3}) wurden bereits im Anschluss
an Theorem \ref{SatzRetPMPeinfach} für den normalen Fall mit $\lambda_0=1$ angegeben. \\[2mm]
Wir wenden uns den notwendigen Bedingungen der Aufgabe (\ref{RetPMP1})--(\ref{RetPMP5}) mit Zustandsbeschränkungen zu.
Bezüglich der Abbildungen in den Zustandsbeschränkungen (\ref{RetPMP5}) setzen wir der Kürze halber: 
$$g_j[t] = g_j\big(t,x_*(t),x_*(t-\delta_x)\big), \quad j=1,...,l.$$
Damit lautet das Pontrjaginsche Maximumprinzip der Aufgabe (\ref{RetPMP1})--(\ref{RetPMP5}): 

\begin{theorem}[Pontrjaginsches Maximumprinzip] \label{SatzRetPMPZB}
\index{Pontrjaginsches Maximumprinzip!ZDelay@-- zeitverzögerte Systeme} 
Es sei $\big(x_*(\cdot),u_*(\cdot)\big) \in \mathscr{B}^{\,\mathcal{Z}}_{\rm adm} \cap \mathscr{B}^{\,\mathcal{Z}}_{\rm Lip}$
und es sei (\ref{RetPMPKommens}) erfüllt.
Ist $\big(x_*(\cdot),u_*(\cdot)\big)$ ein starkes lokales Minimum der Aufgabe (\ref{RetPMP1})--(\ref{RetPMP5}),
dann existieren nicht gleichzeitig verschwindende Multiplikatoren $\lambda_0 \geq 0$,
$p(\cdot) \in W^1_\infty([0,T],\R^n)$, $l_0,\,l_1 \in \R^n$ 
und auf den Mengen
$$T_j=\big\{t \in [t_0,t_1] \,\big|\, g_j\big(t,x_*(t),x_*(t-\delta_x)\big)=0\big\}, \quad j=1,...,l,$$
konzentrierte nichtnegative regul"are Borelsche Ma"s $\mu_j$ endlicher Totalvariation
derart, dass die Funktion $p(\cdot)$ von beschr"ankter Variation und rechtsseitig stetig ist, und
\begin{enumerate}
\item[(a)] die adjungierte Gleichung
           \index{adjungierte Gleichung!ZDelay@-- zeitverzögerte Systeme}
           \begin{eqnarray} 
           \hspace*{-50mm} p(t) &=& l_1 + \int_t^{t_1} H^{\mathcal{Z}}_x[s] \, ds - \sum_{j=1}^l \int_t^{t_1} g_{j,x}[s]\, d\mu_j(s) \nonumber \\                           
                                & & + \chi_{[t_0,t_1-\delta_x]}(t) \cdot \bigg(\int_t^{t_1-\delta_x} H^{\mathcal{Z}}_y[s+\delta_x] \, ds \nonumber \\
          \label{SatzRetPMPZB1} & & \hspace*{3.5cm} - \sum_{j=1}^l \int_t^{t_1-\delta_x} g_{j,y}[s+\delta_x]\, d\mu_j(s+\delta_x)\bigg),
           \end{eqnarray}
\item[(b)] die Transversalitätsbedingungen
           \index{Transversalitätsbedingungen!ZDelay@-- zeitverzögerte Systeme}
           \begin{equation} \label{SatzRetPMPZB2}
            \left. \begin{array}{lcl} p(t_0) &=& l_0, \\[2mm]
                  p(t_1^-)-p(t_1) &=& \displaystyle - \sum_{j=1}^l \mu_j(\{t_1\}) \, g_{j,x}\big(t_1,x_*(t_1),x_*(t_1-\delta_x)\big)
                             \end{array} \right\}
           \end{equation}
\item[(c)] und in fast allen Punkten $t \in [t_0,t_1]$ die Maximumbedingung
           \index{Maximumbedingung!ZDelay@-- zeitverzögerte Systeme}
           \begin{eqnarray}
           && \hspace*{-15mm} H^{\mathcal{Z}}[t,u_*(t),u_*(t-\delta_u)]+\chi_{[t_0,t_1-\delta_u]}(t) \cdot H^{\mathcal{Z}}[t+\delta_u,u_*(t+\delta_u),u_*(t)] \nonumber \\
           \label{SatzRetPMPZB3} && \hspace*{-12mm} = \max_{u \in U} \Big\{ 
                 H^{\mathcal{Z}}[t,u,u_*(t-\delta_u)]+\chi_{[t_0,t_1-\delta_u]}(t) \cdot H^{\mathcal{Z}}[t+\delta_u,u_*(t+\delta_u),u] \Big\}
           \end{eqnarray}
\end{enumerate}
erfüllt sind.
\end{theorem}

Die Optimalitätsbedingungen für die Aufgabe mit Zustandsbeschränkungen in Theorem \ref{SatzRetPMPZB} erweisen sich durch das gleichzeitige Auftreten der
regulären Borelschen Maße zu den Zeitpunkten $t$ und $t+\delta_x$ als sehr komplex.
Dadurch gestaltet sich die Auswertung der Bedingungen in Theorem  \ref{SatzRetPMPZB} mehr als ``nur'' schwierig. \\[2mm]
Die Herleitung von Theorem \ref{SatzRetPMPZB} bzw. Theorem \ref{SatzRetPMP} baut wesentlich auf der Annahme (\ref{RetPMPKommens}) der
kommensurablen Zeitgrößen mit $k \cdot \Delta = \delta_x$, $\kappa \cdot \Delta = \delta_u$, $(N+1) \cdot \Delta = t_1-t_0$,
der Zahl $\Delta>0$ und den natürlichen Zahlen $k,\kappa,N$ auf.

%% file: 3-44-Beweis.tex
\subsubsection{Der Beweis des Maximumprinzips} \label{AbschnittBeweisPMPRet}
Wir folgen im Weiteren \cite{GoMa2} und führen die Aufgabe stückweise auf das ``Elementarintervall'' $[t_0,t_0+\Delta]$ zurück.
Es gelte $(N+1)\Delta=t_1-t_0$ mit der Zahl $\Delta>0$ aus (\ref{RetPMPKommens}).
Für $t \in [t_0,t_0+\Delta]$ setzen wir 
zu $\big(x(\cdot),u(\cdot)\big) \in W^1_\infty([t_0,t_1],\R^n) \times L_\infty([t_0,t_1],U)$:
\begin{equation} \label{RetBeweis1}
\left.\begin{array}{lll}
\xi_i(t) = x(t+i\Delta), & i=0,...,N, & \xi_i(\cdot) \in W^1_\infty([t_0,t_0+\Delta],\R^n), \\[1mm]
\omega_i(t) = u(t+i\Delta), & i=0,...,N, & \omega_i(\cdot) \in L_\infty([t_0,t_0+\Delta],U).
\end{array}\right\}
\end{equation}

Nach (\ref{RetBeweis1}) wird jeder Einschränkung von $x(\cdot)$ bzw. $u(\cdot)$ auf eines der Teilintervalle \\
$I_i=[t_0+i\Delta,t_0+(i+1)\Delta]$, $i=0,...,N$,
eine neue Größe über $[t_0,t_0+\Delta]$ zugeordnet:
\begin{figure}[h]
	\centering
	\fbox{\includegraphics[width=12cm]{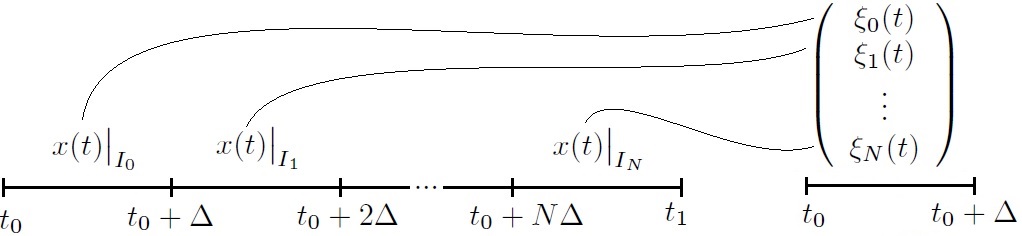}}
\caption[Zeitverzögerte Systeme - Zerlegung der Trajektorien]{Zerlegung der Trajektorien.}
	\label{AbbAufteilung1}
\end{figure}

Durch die Zeitverzögerungen $\delta_x$ und $\delta_u$ entstehen außerdem Variablen zu negativen Indizes.
Diese Variablen stellen keine eigentlichen Optimierungsvariablen dar,
denn mit den gegebenen Größen $x_0 \in \R^n$ und $u_0 \in  U$ gelten
\begin{equation} \label{RetBeweisNegIndex}
\xi_i(t) = x_0 \mbox{ für } i=-k,...,-1 \mbox{ und } \omega_i(t) = u_0 \mbox{ für } i=-\kappa,...,-1.
\end{equation}
Die Stetigkeit von $x(\cdot)$ über $[t_0,t_1]$ liefert für die $\xi_i(\cdot)$ die gekoppelten Randbedingungen
$$\xi_i(t_0+\Delta)=\xi_{i+1}(t_0) \quad\Leftrightarrow\quad \xi_{i+1}(t_0)-\xi_i(t_0+\Delta)=0
  \quad\mbox{für } i=0,...,N-1.$$
Die Umkehrung von (\ref{RetBeweis1}) ergibt die Zusammenhänge
$$x(t)= \xi_i(t-i\Delta), \; u(t)=\omega_i(t-i\Delta) \mbox{ für } t \in [t_0+i\Delta,t_0+(i+1)\Delta) \mbox{ und } i=0,...,N,$$
sowie die ``Wiederherstellung'' der Eingangsfunktionen über $[t_0,t_1]$ durch
\begin{equation} \label{RetBeweis1Einschub}
\left.\begin{array}{l}
\displaystyle x(t) = \sum_{i=0}^N \chi_{[t_0+i\Delta,t_0+(i+1)\Delta]}(t) \cdot \xi_i(t-i\Delta), \\[1mm]
\displaystyle u(t) = \sum_{i=0}^N \chi_{[t_0+i\Delta,t_0+(i+1)\Delta]}(t) \cdot \omega_i(t-i\Delta).
\end{array}\right\}
\end{equation}
Trotz der abgeschlossenen Teilintervalle wird $x(\cdot)$ in den inneren Nahtstellen $t=t_0+i\Delta$ durch die gekoppelten
Randbedingungen eindeutig und stetig festgelegt.
Für die messbare Steuerung $u(\cdot)$ ist die Setzung in den Nahtstellen ohne Bedeutung.

\begin{figure}[h]
	\centering
	\fbox{\includegraphics[width=14cm]{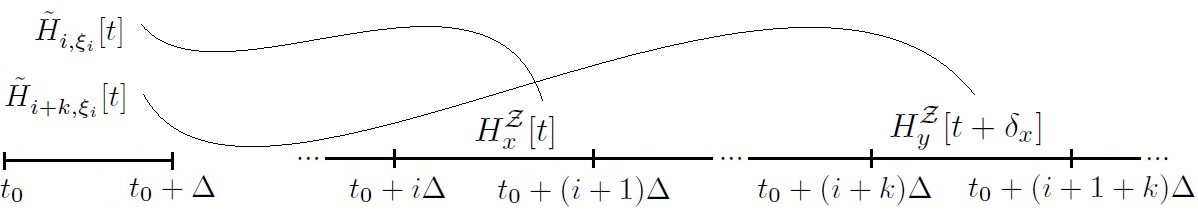}}
\caption[Zeitverzögerte Systeme - Pontrjagin-Funktionen]{Zuordnungen der Pontrjagin-Funktionen.}
	\label{AbbAufteilung2}
\end{figure}

\newpage
Wir fassen $\xi_0(\cdot),...,\xi_N(\cdot)$ und $\omega_0(\cdot),...,\omega_N(\cdot)$ zu den Vektorfunktionen
$$\big(\xi(\cdot),\omega(\cdot)\big) \in W^1_\infty([t_0,t_0+\Delta],\R^{n \cdot (N+1)}) \times L_\infty([t_0,t_0+\Delta],U^{N+1}), \quad
  U^{N+1}=\underbrace{U \times ... \times U}_{(N+1)-\mbox{mal}},$$
zusammen.
Auf diese Weise erhalten wir anstelle von (\ref{RetPMP1})--(\ref{RetPMP5}) die Aufgabe
\begin{equation} \label{RetBeweis2}
\left. \hspace*{-5mm}\begin{array}{l}
\displaystyle \tilde{J}\big(\xi(\cdot),\omega(\cdot)\big) = \int_{t_0}^{t_0+\Delta} \bigg(
                      \sum_{i=0}^N f\big(t+i\Delta,\xi_i(t),\xi_{i-k}(t),\omega_i(t),\omega_{i-\kappa}(t)\big) \bigg) \, dt \to \inf, \\[1mm]
\displaystyle \dot{\xi}_i(t) = \varphi\big(t+i\Delta,\xi_i(t),\xi_{i-k}(t),\omega_i(t),\omega_{i-\kappa}(t)\big), \; i=0,...,N, \\[1mm]
\displaystyle \xi_0(t_0)=x_0, \quad \xi_N(t_0+\Delta)=x_1, \quad \xi_{i+1}(t_0)-\xi_i(t_0+\Delta)=0,\; i=0,...,N-1,\\[1mm]
\displaystyle \omega(t)=\big(\omega_0(t),...,\omega_N(t)\big) \in U^{N+1}, \quad U \subseteq \R^m, \quad U\not= \emptyset,  \\[1mm]
\displaystyle g_j\big(t+i\Delta,\xi_i(t),\xi_{i-k}(t)\big) \leq 0 \quad \mbox{f"ur alle } t \in [t_0,t_0+\Delta],\, j=1,...,l,\, i=0,...,N.
\end{array}\right\}
\end{equation}
In (\ref{RetBeweis2}) führen wir folgende Abbildungen ein:
\begin{eqnarray*}
\tilde{f}\big(t,\xi(t),\omega(t)\big) &=& \sum_{i=0}^N f\big(t+i\Delta,\xi_i(t),\xi_{i-k}(t),\omega_i(t),\omega_{i-\kappa}(t)\big), \\
\tilde{\varphi}\big(t,\xi(t),\omega(t)\big) &=&
     \left( \begin{array}{l}
            \varphi\big(t,\xi_0(t),\xi_{-k}(t),\omega_0(t),\omega_{-\kappa}(t)\big) \\
            \varphi\big(t+\Delta,\xi_1(t),\xi_{1-k}(t),\omega_1(t),\omega_{1-\kappa}(t)\big) \\
            \hspace*{20mm} \vdots \\
            \varphi\big(t+N\Delta,\xi_N(t),\xi_{N-k}(t),\omega_N(t),\omega_{N-\kappa}(t)\big)
            \end{array}\right), \\
\tilde{h}\big(\xi(t_0),\xi(t_1)\big) &=& 
     \left( \begin{array}{l}
            \xi_0(t_0)-x_0 \\
            \xi_1(t_0) - \xi_0(t_0+\Delta) \\
            \hspace*{20mm} \vdots \\
            \xi_N(t_0)-\xi_{N-1}(t_0+\Delta) \\
            x_1-\xi_N(t_0+\Delta)
            \end{array}\right), \\
\tilde{g}_{ij}\big(t,\xi(t)\big) &=& g_j\big(t+i\Delta,\xi_i(t),\xi_{i-k}(t)\big), \quad i=0,...,N,\, j=1,...,l.
\end{eqnarray*}
So erhält (\ref{RetBeweis2}) die Gestalt der Aufgabe (\ref{GZF1})--(\ref{GZF5}):
\begin{eqnarray}
&& \label{RetBeweis3} \tilde{J}\big(\xi(\cdot),\omega(\cdot)\big) = \int_{t_0}^{t_0+\Delta} \tilde{f}\big(t,\xi(t),\omega(t)\big) \, dt  \to \inf, \\
&& \label{RetBeweis4} \dot{\xi}(t) = \tilde{\varphi}\big(t,\xi(t),\omega(t)\big), \\
&& \label{RetBeweis5} \tilde{h}\big(\xi(t_0),\xi(t_0+\Delta)\big)=0, \\
&& \label{RetBeweis6} w(t) \in U^{N+1}, \quad U \subseteq \R^m, \quad U\not= \emptyset, \\
&& \label{RetBeweis7} \tilde{g}_{ij}\big(t,\xi(t)\big) \leq 0 \quad \mbox{f"ur alle } t \in [t_0,t_0+\Delta], \quad i=0,...,N,\, j=1,...,l.
\end{eqnarray}
In dieser Aufgabenstellung sind $\xi_{-k}(\cdot),...,\xi_{-1}(\cdot),\omega_{-\kappa}(\cdot),...,\omega_{-1}(\cdot)$ enthalten,
die weiterhin keine echten Optimierungsgrößen darstellen, sondern durch (\ref{RetBeweisNegIndex}) bereits festgelegt sind. \\
Die Pontrjagin-Funktion $\tilde{H}(t,\xi,w,\tilde{p},\lambda_0)$, $\tilde{p}=(\tilde{p}_0,...,\tilde{p}_N)^T$, besitzt die Gestalt
\begin{eqnarray}
\lefteqn{\tilde{H}(t,\xi,w,\tilde{p},\lambda_0) = \langle \tilde{p},\tilde{\varphi}(t,\xi,w)\rangle - \lambda_0 \tilde{f}(t,\xi,w)} \nonumber \\
\label{RetBeweisPontr} && \hspace*{-5mm} 
  = \sum_{i=0}^N \Big(\underbrace{\langle \tilde{p}_i,\varphi(t+i\Delta,\xi_i,\xi_{i-k},\omega_i,\omega_{i-\kappa})\rangle 
             - \lambda_0 f(t+i\Delta,\xi_i,\xi_{i-k},\omega_i,\omega_{i-\kappa})}_{=\tilde{H}_i(t,\xi_i,\xi_{i-k},\omega_i,\omega_{i-\kappa},\tilde{p}_i,\lambda_0)} \Big).
\end{eqnarray}
Ist $\big(\xi_*(\cdot),\omega_*(\cdot)\big)$ ein starkes lokales Minimum in (\ref{RetBeweis3})--(\ref{RetBeweis6}),
dann existieren nach Theorem \ref{SatzPMPGZF} nicht gleichzeitig verschwindende Multiplikatoren
$\lambda_0 \geq 0$, 
die adjungierte Funktion $\tilde{p}(\cdot) = \big(\tilde{p}_0(\cdot),...,\tilde{p}_N(\cdot)\big)^T$ mit
$\tilde{p}_0(\cdot),...,\tilde{p}_N(\cdot) \in W^1_\infty([t_0,t_0+\Delta],\R^n)$
und der Vektor $\tilde{l} = (l_0,\tilde{l}_0,...,\tilde{l}_{N-1},l_1)^T$ mit $l_0, ...,l_1 \in \R^n$
derart, dass die adjungierte Gleichung
\begin{equation} \label{RetBeweisAdjGl}
\dot{\tilde{p}}(t)=-\tilde{H}_\xi\big(t,\xi_*(t),w_*(t),\tilde{p}(t),\lambda_0\big),
\end{equation}
mit der Abbildung $\tilde{h}=\tilde{h}(\xi^0,\xi^1)$ in (\ref{RetBeweis5}),
deren Ableitungen $\tilde{h}_{\xi^0}$ und $\tilde{h}_{\xi^1}$ Diagonalmatrizen mit Eintragungen aus dem $\R^n$ sind,
die Transversalit"atsbedingungen
\begin{equation} \label{RetBeweisTrBed}
\left. \begin{array}{lclcl}
\tilde{p}(t_0)        &=& \tilde{h}_{\xi^0}^T\big(\xi_*(t_0),\xi_*(t_0+\Delta)\big) \tilde{l}
                      &=& (l_0,\tilde{l}_0,...,\tilde{l}_{N-1},0)^T, \\[2mm]
\tilde{p}(t_0+\Delta) &=& - \tilde{h}_{\xi^1}^T\big(\xi_*(t_0),\xi_*(t_0+\Delta)\big) \tilde{l}
                      &=& (0,\tilde{l}_0,...,\tilde{l}_{N-1},l_1)^T
                             \end{array} \right\}
\end{equation}
und f"ur fast alle $t\in [t_0,t_1]$ die Maximumbedingung
\begin{equation} \label{RetBeweisMaxBed}
\tilde{H}\big(t,\xi_*(t),w_*(t),\tilde{p}(t),\lambda_0\big)
   = \max_{\omega \in U^{N+1}} \tilde{H}\big(t,\xi_*(t),w,\tilde{p}(t),\lambda_0\big)
\end{equation}
erfüllt sind. \\[2mm]
Zum weiteren Vorgehen geben wir die Bedingungen (\ref{RetBeweisAdjGl})-- (\ref{RetBeweisMaxBed})
für alle Komponenten zu den Indizes $i=0,...,N$ explizit an.
Dazu führen wir mit Hilfe von (\ref{RetBeweisPontr}) folgende abkürzenden Bezeichnungen ein:
\begin{eqnarray*}
\tilde{H}_i[t] &=& \tilde{H}_i\big(t,\xi_i^*(t),\xi^*_{i-k}(t),\omega_i^*(t),\omega_{i-\kappa}^*(t),\tilde{p}_i(t),\lambda_0\big), \\
\tilde{H}_i[t,u,v] &=& \tilde{H}_i\big(t,\xi_i^*(t),\xi^*_{i-k}(t),u,v,\tilde{p}_i(t),\lambda_0\big).
\end{eqnarray*} 
Bezüglich der einzelnen Zustandsvariablen $\xi_0,...,\xi_N$ führt die adjungierte Gleichung (\ref{RetBeweisAdjGl}) zu dem System
\begin{equation} \label{RetBeweis8}
\dot{\tilde{p}}_i (t) = -\tilde{H}_{i,\xi_i}[t] - \chi_{\{0,...,N-k\}}(i+k) \cdot \tilde{H}_{i+k,\xi_i}[t], \quad i=0,...,N.
\end{equation}
Darin ist zu beachten,
dass sich die Indizes $i$ bzw. $i+k$ auf verschiedene Teilintervalle über dem Intervall $[t_0,t_1]$
und dass sich die Ableitungen $\tilde{H}_{i,\xi_i}[t]$ und $\tilde{H}_{i+k,\xi_i}[t]$ nach der Variable $\xi_i$
auf verschiedene Variablenpositionen beziehen. \\[2mm]
Die Transversalitätsbedingungen (\ref{RetBeweisTrBed}) ergeben die Kopplungen
\begin{equation} \label{RetBeweis9}
\tilde{p}_0(t_0+\Delta)= \tilde{p}_1(t_0)=\tilde{l}_0, \; ...\;, \;
  \tilde{p}_{N-1}(t_0+\Delta)= \tilde{p}_{N}(t_0)=\tilde{l}_{N-1},
\end{equation}
sowie die zusätzlichen Bedingungen
\begin{equation} \label{RetBeweis10}
\tilde{p}_0(t_0)=l_0, \qquad \tilde{p}_N(t_0+\Delta)=l_1.
\end{equation}
Wir wählen einen Index $i \in \{0,...,N\}$ und betrachten $\omega(i) \in U^{N+1}$ mit $\omega_j(i)= \omega_j^*$ für $j \not= i$.
So führt die Maximumbedingung (\ref{RetBeweisMaxBed}) auf die Beziehung
\begin{eqnarray}
\lefteqn{\tilde{H}_i[t,\omega^*,\omega_{i-\kappa}^*(t)]
        + \chi_{\{0,...,N-\kappa\}}(i+\kappa) \cdot \tilde{H}_{i+\kappa}[t,\omega_{i+\kappa}^*(t),\omega^*]} \nonumber \\
&=& \label{RetBeweis12} \max_{\omega \in U} \Big\{ \tilde{H}_i[t,\omega,\omega_{i-\kappa}^*(t)]
        + \chi_{\{0,...,N-\kappa\}}(i+\kappa) \cdot \tilde{H}_{i+\kappa}[t,\omega_{i+\kappa}^*(t),\omega]\Big\},
\end{eqnarray}
die für jeden Index $i \in \{0,...,N\}$ und fast alle $t \in [t_0,t_0+\Delta]$ gültig ist. \\[2mm]
Zum Beweis von Theorem \ref{SatzRetPMP} müssen die Ergebnisse der Aufgabe (\ref{RetBeweis3})--(\ref{RetBeweis6})
auf die Aufgabe (\ref{RetPMP1})--(\ref{RetPMP4}) zurückgeführt werden.
Durch die Umkehrung (\ref{RetBeweis1Einschub})
lassen sich die Steuerungsprozesse $\big(x(\cdot),u(\cdot)\big)$ der Aufgabe (\ref{RetPMP1})--(\ref{RetPMP4})
und $\big(\xi(\cdot),\omega(\cdot)\big)$ der Aufgabe (\ref{RetBeweis3})--(\ref{RetBeweis6}) einander zuordnen.
In gleicher Weise verfahren wir mit den Adjungierten und den Pontrjagin-Funktionen:
\begin{eqnarray*}
p(t) &=& \sum_{i=0}^N \chi_{[t_0+i\Delta,t_0+(i+1)\Delta]}(t) \cdot \tilde{p}_i(t-i\Delta), \\
H^{\mathcal{Z}}[t] &=& \sum_{i=0}^N \chi_{[t_0+i\Delta,t_0+(i+1)\Delta]}(t) \cdot \tilde{H}_i[t-i\Delta].
\end{eqnarray*}
In den inneren Nahtstellen $t_0+i\Delta$ ist $p(\cdot)$ durch die Kopplungsbedingungen (\ref{RetBeweis9}) eindeutig festgelegt.
Weiterhin ist zu beachten,
dass der Steuerungsprozess $\big(\xi(\cdot),\omega(\cdot)\big)$ in $\tilde{H}_i[t,u,v]$ einfließt,
während $\big(x(\cdot),u(\cdot)\big)$ in $H[t,u,v]$ berücksichtigt wird. \\[2mm]
Wir verwenden nun diese Zuordnungen zwischen den Steuerungsprozessen $\big(x_*(\cdot),u_*(\cdot)\big)$ und $\big(\xi_*(\cdot),\omega_*(\cdot)\big)$.
Ferner beachten wir (\ref{RetPMPKommens}), d.\,h. $k \cdot \Delta = \delta_x$ und $\kappa \cdot \Delta = \delta_u$.
Au"serdem erinnern wir an die Bezeichnungen $f=f(t,x,y,u,v)$ und $\varphi=\varphi(t,x,y,u,v)$ der Variablen der Abbildungen $f$ und $\varphi$
in der Aufgabe (\ref{RetPMP1})--(\ref{RetPMP4}),
die zu den Ableitungen $H^{\mathcal{Z}}_x$ und $H^{\mathcal{Z}}_y$ der Pontrjagin-Funktion $H^{\mathcal{Z}}$
anstelle der Ableitungen $\tilde{H}_{i,\xi_i}$ und $\tilde{H}_{i+k,\xi_i}$
der Pontrjagin-Funktion $\tilde{H}$ der Aufgabe (\ref{RetPMP1})--(\ref{RetPMP4}) führen.
Zu alledem beziehen sich die Indizes $i$ bzw. $i+k$ auf verschiedene Teilintervalle über dem Intervall $[t_0,t_1]$.
Diese Sachverhalte sind in der Abbildung \ref{AbbAufteilung2} dargestellt. \\[2mm]
Es ergibt sich also,
dass die Funktion $p(\cdot)$ aufgrund von (\ref{RetBeweis8}) in jedem Teilintervall $(t_0+i\Delta,t_0+(i+1)\Delta)$ die verallgemeinerte Ableitung
$$\dot{p}(t) = - H^{\mathcal{Z}}_x[t] - \chi_{[0,T-\delta_x]}(t) \cdot H^{\mathcal{Z}}_y[t+\delta_x]$$
besitzt, aufgrund der Bedingungen (\ref{RetBeweis9}) in den Stellen $t=t_0+i\Delta$ für $i=1,...,N$ stetig ist
(also dem Raum $W^1_\infty([t_0,t_1],\R^n)$ angehört) und
nach (\ref{RetBeweis10}) die Bedingungen $p(t_0)=l_0$ und $p(t_1)=l_1$ erfüllt.
Da ferner die Maximumbedingung (\ref{RetBeweis12}) für jedes $i \in \{0,...,N\}$ erfüllt ist,
ergibt sich in jedem Teilintervall $(t_0+i\Delta,t_0+(i+1)\Delta)$ fast überall die Gültigkeit der Beziehung

\begin{eqnarray*}
\lefteqn{H^{\mathcal{Z}}[t,u_*(t),u_*(t-\delta_u)]+\chi_{[0,T-\delta_u]}(t) \cdot H^{\mathcal{Z}}[t+\delta_u,u_*(t+\delta_u),u_*(t)]} \nonumber \\
&& = \max_{u \in U} \Big\{ H^{\mathcal{Z}}[t,u,u_*(t-\delta_u)]+\chi_{[0,T-\delta_u]}(t) \cdot H^{\mathcal{Z}}[t+\delta_u,u_*(t+\delta_u),u] \Big\}.
\end{eqnarray*}
Die Verschiebungen um die Verzögerung $\delta_u$ ergeben sich hierin auf die gleiche Weise wie sie in Abbildung \ref{AbbAufteilung2}
für die Verzögerung $\delta_x$ dargestellt sind.
Zusammenfassend sind über jedem Teilintervall $(t_0+i\Delta,t_0+(i+1)\Delta)$ die adjungierte Gleichung und die Maximumbedingung erfüllt;
somit gelten diese Beziehungen fast überall über $[t_0,t_1]$.
Ferner gehört die Funktion $p(\cdot)$ dem Raum $W^1_\infty([t_0,t_1],\R^n)$ an und erfüllt die Randbedingungen $p(t_0)=l_0$, $p(t_1)=l_1$.
Damit ist Theorem \ref{SatzRetPMP} bewiesen. \\[2mm]
Wir betrachten nun die Zustandsbeschränkungen (\ref{RetBeweis7}).
Es bezeichne:
$$\tilde{g}_{ij}[t] = g_j\big(t+i\Delta,\xi_i^*(t),\xi^*_{i-k}(t)\big), \quad i=0,...,N, \; j=1,...,l.$$
Ist $\big(\xi_*(\cdot),\omega_*(\cdot)\big)$ ein starkes lokales Minimum der Aufgabe (\ref{RetBeweis3})--(\ref{RetBeweis7}),
dann gibt es nach Theorem \ref{SatzPMPGZFZA} auf den Mengen
$$\tilde{T}_{ij}=\big\{t \in [t_0,t_0+\Delta] \,\big|\, g_{ij}\big(t,\xi_*(t)\big)=0\big\}, \quad i=0,...,N,\; j=1,...,l,$$
konzentrierte nichtnegative regul"are Borelsche Ma"se $\tilde{\mu}_{ij}$ endlicher Totalvariation derart,
dass zu $i \in \{0,...,N\}$ die adjungierte Gleichung die Gestalt
\begin{eqnarray}
\tilde{p}_i(t) &=& \tilde{p}_i(t_0+\Delta)
     + \int_t^{t_0+\Delta} \tilde{H}_{i,\xi_i}[s] \, ds -\sum_{j=1}^l \int_t^{t_0+\Delta} \tilde{g}_{ij,\xi_i}[s]\, d\tilde{\mu}_{ij}(s) \nonumber \\
\label{RetBeweis11}  & & \hspace*{-25mm} + \chi_{\{0,...,N-k\}}(i+k) \cdot \bigg(
       \int_t^{t_0+\Delta} \tilde{H}_{i+k,\xi_i}[s] \, ds  -\sum_{j=1}^l \int_t^{t_0+\Delta} \tilde{g}_{(i+k)j,\xi_i}[s]\, d\tilde{\mu}_{(i+k)j}(s)\bigg)
\end{eqnarray}
besitzt und die Sprung-Transversalitätsbedingungen
\begin{eqnarray*}
\tilde{p}_i(t_0+\Delta^-) - \tilde{p}_i(t_0+\Delta)
&=& - \sum_{j=1}^l \tilde{\mu}_{ij}(\{t_0+\Delta\})\,\tilde{g}_{ij,\xi_i}[t_0+\Delta]  \\
& & \hspace*{-30mm} -  \chi_{\{0,...,N-k\}}(i+k) \cdot\sum_{j=1}^l \tilde{\mu}_{(i+k)j}(\{t_0+\Delta\})\, \tilde{g}_{(i+k)j,\xi_i}[t_0+\Delta], \quad i<N, \\
\tilde{p}_N(t_0+\Delta^-) - \tilde{p}_N(t_0+\Delta)
&=& - \sum_{j=1}^l \tilde{\mu}_{Nj}(\{t_0+\Delta\})\, \tilde{g}_{Nj,\xi_N}[t_0+\Delta]
\end{eqnarray*}
gelten.
Zur Rückführung über $[t_0,t_1]$ verfahren wir wie eben und übernehmen die obigen Bezeichnungen, z.\,B. für $p[t]$ und $H^{\mathcal{Z}}[t]$.
Die Sprung-Transversalitätsbedingungen beziehen sich für $i<N$ auf die inneren Nahtstellen $t=t_0+(i+1)\Delta$ des Intervalls $[t_0,t_1]$,
in denen durch diese Bedingungen Sprungstellen entstehen können.
Für $i=N$ ergibt sich die Transversalitätsbedingung in $t=t_1$. \\[2mm]
Wir setzen nun bezüglich der Zustandsbeschränkungen für $j=1,...,l$:
\begin{eqnarray*}
g_j[t] &=& \sum_{i=0}^N \chi_{[t_0+i\Delta,t_0+(i+1)\Delta)}(t) \cdot \tilde{g}_{ij}[t-i\Delta], \\
\mu_j(t) &=& \sum_{i=0}^N \int_{t_0}^t \chi_{[t_0+i\Delta,t_0+(i+1)\Delta]}(s) \, d\tilde{\mu}_{ij}(s-i\Delta).
\end{eqnarray*}
In den Abbildungen $\tilde{g}_{ij}$ ersetzen wir $\xi^*_i(t)$, $\xi_{i-k}^*(t)$ in den Teilintervallen $[t_0+i\Delta,t_0+(i+1)\Delta]$
durch $x_*(t)$ bzw. $x_*(t-\delta_x)$.
Die Abbildungen $g_j[\cdot]$ sind dann über $[t_0,t_1]$ stetig,
da sämtliche einfließende Abbildungen in den Anschlusspunkten $t=t_0+i\Delta$ stetig sind. \\
Bezüglich der Maße muss die charakteristische Funktion über dem Abschluss der Teilintervalle betrachtet werden,
damit keine atomaren Anteile in den Anschlusspunkten verloren gehen.
Die Maße $\mu_j$ sind damit nichtnegativ und von beschränkter Totalvariation.
Ferner verschwindet für $j \in \{1,...,l\}$ das Maß $\mu_j$ genau dann über $[t_0,t_1]$, 
wenn sämtliche Maße $\tilde{\mu}_{ij}$ für $i=0,...,N$ über $[t_0,t_0+\Delta]$ verschwinden. \\
Die Mengen $\tilde{T}_{ij}$ erhalten bei der Rückführung die Gestalt
$$T_{ij}=\big\{t \in [t_0+i\Delta,t_0+(i+1)\Delta] \,\big|\, g_j\big(t,x_*(t),x_*(t-\delta_x)\big) =0\big\}$$
und wir setzen lediglich noch $T_j= T_{0j} \cup ... \cup T_{Nj}$ für $j=1,...,l$. \\
Zusammenfassend gibt es auf den Mengen $T_j=\big\{t \in [t_0,t_1] \,\big|\, g_j\big(t,x_*(t),x_*(t-\delta_x)\big)=0\big\}$, $j=1,...,l$,
konzentrierte nichtnegative regul"are Borelsche Ma"s $\mu_j$ endlicher Totalvariation
derart, dass die Funktion $p(\cdot)$ von beschr"ankter Variation und rechtsseitig stetig ist,
der adjungierten Gleichung (\ref{SatzRetPMPZB1}) genügt und die Transversalitätsbedingungen (\ref{SatzRetPMPZB2}) erfüllt.
Schließlich ist somit auch Theorem \ref{SatzRetPMPZB} nachgewiesen. \hfill $\blacksquare$

%% file: 3-45-Arrow.tex
\subsubsection{Hinreichende Bedingungen nach Arrow}
Die Ableitung von hinreichenden Optimalitästbedingungen nach Arrow \index{hinreichende Bedingungen, Arrow!Standard@-- Zeitverzögerte Systeme}
gestaltet sich durch das Auftreten der zeitlichen Verschiebung um $\delta_x>0$ schwieriger.
Dennoch führt uns das Vorgehen in Abschnitt \ref{AbschnittHBPMP} zum Ziel.
Auf die Betrachtung der Zustandsbeschränkungen (\ref{RetPMP5}) verzichten wir und konzentrieren uns auf die Behandlung
der Terme zu verschiedenen Zeitpunkten. \\[2mm]
In der Untersuchung der Aufgabe (\ref{RetPMP1})--(\ref{RetPMP4}) verwenden wir wieder die Menge
$$V^{\mathcal{S}}_\gamma(t)=\{ x \in \R^n \,|\, \|x-x_*(t)\| < \gamma\}$$
aus Abschnitt \ref{AbschnittHBPMP},
sowie die Menge $V^{\mathcal{Z}}_\gamma(t) =V^{\mathcal{S}}_\gamma(t) \times V^{\mathcal{S}}_\gamma(t-\delta_x)$, d.\,h.
$$V^{\mathcal{Z}}_\gamma(t) = \{ (x,y) \in \R^n \times \R^n \,|\, \|x-x_*(t)\| < \gamma, \|y-x_*(t-\delta_x)\| < \gamma\}.$$
Ferner bezeichnet $\mathscr{H}^{\mathcal{Z}}$ die Hamilton-Funktion
$$\mathscr{H}^{\mathcal{Z}}(t,x,y,p) = \sup_{u,v \in U} H^{\mathcal{Z}}\big(t,x,y,u,v,p,1).$$

\begin{theorem} 
In der Aufgabe (\ref{RetPMP1})--(\ref{RetPMP4}) sei
$\big(x_*(\cdot),u_*(\cdot)\big) \in \mathscr{B}^{\,\mathcal{Z}}_{\rm adm} \cap \mathscr{B}^{\,\mathcal{Z}}_{\rm Lip}$
und es sei $p(\cdot) \in W^1_\infty([t_0,t_1],\R^n)$. Ferner gelte:
\begin{enumerate}
\item[(a)] Das Tripel $\big(x_*(\cdot),u_*(\cdot),p(\cdot)\big)$
           erf"ullt (\ref{SatzRetPMP1})--(\ref{SatzRetPMP3}) in Theorem \ref{SatzRetPMP}.        
\item[(b)] F"ur jedes $t \in [t_0,t_1]$ ist die Hamilton-Funktion $\mathscr{H}^{\mathcal{Z}}\big(t,x,y,p(t)\big)$ konkav auf $V^{\mathcal{Z}}_\gamma(t)$.
\end{enumerate}
Dann ist $\big(x_*(\cdot),u_*(\cdot)\big)$ ein starkes lokales Minimum der Aufgabe (\ref{RetPMP1})--(\ref{RetPMP4}).
\end{theorem}

{\bf Beweis} Wie im Abschnitt \ref{AbschnittHBPMP} folgt,
dass $\alpha_*= -\mathscr{H}^{\mathcal{Z}}\big(t,x_*(t),x_*(t-\delta_x),p(t)\big)$ ein Randpunkt der konvexen Menge 
$$Z=\big\{ (\alpha,x,y) \in \R \times \R^n \times \R^n \,\big|\, 
           (x,y) \in V^{\mathcal{Z}}_\gamma(t), \alpha \geq -\mathscr{H}^{\mathcal{Z}}\big(t,x,y,p(t)\big) \big\}$$
ist.
Daher existiert ein nichttrivialer Vektor $\big(a_0(t),a_1(t),a_2(t)\big) \in \R \times \R^n \times \R^n$ mit
$$a_0(t) \alpha + \langle a_1(t),x\rangle + \langle a_2(t),y\rangle
  \geq a_0(t) \alpha_* + \langle a_1(t),x_*(t)\rangle+ \langle a_2(t),x_*(t-\delta_x)\rangle$$
f"ur alle $(\alpha,x,y) \in Z$.
Weiterhin ergibt sich für $y=x_*(t-\delta_x)$ mit der gleichen Argumentation wie in Abschnitt \ref{AbschnittHBPMP},
dass wir ohne Einschr"ankung $a_0(t)=1$ annehmen können und erhalten
\begin{eqnarray}
\lefteqn{\langle a_1(t),x-x_*(t)\rangle + \langle a_2(t),y-x_*(t-\delta_x)\rangle} \nonumber\\
&\geq& \label{BeweisHBPMPZ1} \mathscr{H}^{\mathcal{Z}}\big(t,x,y,p(t)\big) - \mathscr{H}^{\mathcal{Z}}\big(t,x_*(t),x_*(t-\delta_x),p(t)\big).
\end{eqnarray}
Aus (\ref{BeweisHBPMPZ1}) ergeben sich für alle $x \in V^{\mathcal{S}}_\gamma(t)$ und alle
$y \in V^{\mathcal{S}}_\gamma(t-\delta_x)$ die Relationen
\begin{eqnarray*}
       \langle a_1(t),x-x_*(t)\rangle
&\geq& \mathscr{H}^{\mathcal{Z}}\big(t,x,x_*(t-\delta_x),p(t)\big) - \mathscr{H}^{\mathcal{Z}}\big(t,x_*(t),x_*(t-\delta_x),p(t)\big), \\
       \langle a_2(t),y-x_*(t-\delta_x)\rangle
&\geq& \mathscr{H}^{\mathcal{Z}}\big(t,x_*(t),y,p(t)\big) - \mathscr{H}^{\mathcal{Z}}\big(t,x_*(t),x_*(t-\delta_x),p(t)\big).
\end{eqnarray*}
In der zweiten Ungleichung nehmen wir eine Zeitverschiebung um $\delta_x$ vor und überführen den Zeitpunkt $t$ in $t + \delta_x$.
Damit gehört $y$ der Menge $V^{\mathcal{S}}_\gamma(t)$ an und wir dürfen nach der Verschiebung $y=x \in V^{\mathcal{S}}_\gamma(t)$ wählen.
So ergibt sich für alle $x \in V^{\mathcal{S}}_\gamma(t)$ die Ungleichung
\begin{eqnarray}
      \lefteqn{\langle a_1(t)+\chi_{[t_0,t_1-\delta_x]}(t) \cdot a_2(t+\delta_x),x-x_*(t)\rangle} \nonumber \\
&\geq& \mathscr{H}^{\mathcal{Z}}\big(t,x,x_*(t-\delta_x),p(t)\big) - \mathscr{H}^{\mathcal{Z}}\big(t,x_*(t),x_*(t-\delta_x),p(t)\big) \nonumber \\
&& + \chi_{[t_0,t_1-\delta_x]}(t) \cdot \Big(\mathscr{H}^{\mathcal{Z}}\big(t+\delta_x,x_*(t+\delta_x),x,p(t+\delta_x)\big) \nonumber \\
&& \label{BeweisHBPMPZ2} \hspace*{3,5cm} - \mathscr{H}^{\mathcal{Z}}\big(t+\delta_x,x_*(t+\delta_x),x_*(t),p(t+\delta_x)\big)\Big).
\end{eqnarray} 
Es sei in $t \in [t_0,t_1]$ die Maximumbedingung (\ref{SatzRetPMP3}) erf"ullt.
Mit Hilfe der Pontrjagin-Funktion können wir die Ungleichung (\ref{BeweisHBPMPZ2}) weiterführen und erhalten
\begin{eqnarray}
&\geq& H^{\mathcal{Z}}\big(t,x,x_*(t-\delta_x),u_*(t),u_*(t-\delta_x),p(t),1\big) - \mathscr{H}^{\mathcal{Z}}\big(t,x_*(t),x_*(t-\delta_x),p(t)\big) \nonumber \\
&& + \chi_{[t_0,t_1-\delta_x]}(t) \cdot \Big(H^{\mathcal{Z}}\big(t+\delta_x,x_*(t+\delta_x),x,u_*(t+\delta_x),u_*(t),p(t+\delta_x),1\big) \nonumber \\
&& \label{BeweisHBPMPZ3} \hspace*{3,5cm} - \mathscr{H}^{\mathcal{Z}}\big(t+\delta_x,x_*(t+\delta_x),x_*(t),p(t+\delta_x)\big)\Big)=\Psi(x)
\end{eqnarray} 
f"ur alle $x \in V^{\mathcal{S}}_\gamma(t)$.
Mit Hilfe der rechten Seite $\Psi(x)$ von (\ref{BeweisHBPMPZ3}) bilden wir die Funktion
$$\Phi(x) = \Psi(x) - \langle a_1(t)+\chi_{[t_0,t_1-\delta_x]}(t) \cdot a_2(t+\delta_x),x-x_*(t)\rangle,$$
welche in dem inneren Punkt $x_*(t)$ der Menge $V^{\mathcal{S}}_\gamma(t)$ ihr globales Maximum annimmt.
Also gilt $0=\Phi'(x_*(t))$, d.\,h.
\begin{eqnarray}
&& \hspace*{-15mm} -a_1(t)-\chi_{[t_0,t_1-\delta_x]}(t) \cdot a_2(t+\delta_x) \nonumber \\
&& \hspace*{-10mm} = -H_x^{\mathcal{Z}}\big(t,x_*(t),x_*(t-\delta_x),u_*(t),u_*(t-\delta_x),p(t),1\big) \nonumber \\
&& \label{BeweisHBPMPZ4} \hspace*{-10mm} - \chi_{[t_0,t_1-\delta_x]}(t) \cdot H_y^{\mathcal{Z}}\big(t+\delta_x,x_*(t+\delta_x),x_*(t),u_*(t+\delta_x),u_*(t),p(t+\delta_x),1\big).
\end{eqnarray}
Die adjungierte Gleichung (\ref{SatzRetPMP1}) zeigt nun
\begin{equation} \label{BeweisHBPMPZ5}
\dot{p}(t) = -a_1(t)-\chi_{[t_0,t_1-\delta_x]}(t) \cdot a_2(t+\delta_x)\quad \mbox{ für fast alle } t \in [t_0,t_1].
\end{equation}
Es sei $\big(x(\cdot),u(\cdot)\big) \in \mathscr{B}^{\,\mathcal{Z}}_{\rm adm}$ mit $\|x(\cdot)-x_*(\cdot)\|_\infty < \gamma$.
Wir erhalten mit (\ref{BeweisHBPMPZ1}) und (\ref{BeweisHBPMPZ5}):
\begin{eqnarray*}
\lefteqn{\int_{t_0}^{t_1} \big[H^{\mathcal{Z}}\big(t,x_*(t),x_*(t-\delta_x),u_*(t),u_*(t-\delta_x),p(t),1\big)} \\
& & \hspace*{3cm} -H^{\mathcal{Z}}\big(t,x(t),x(t-\delta_x),u(t),u(t-\delta_x),p(t),1\big)\big] \, dt \\
&=& \int_{t_0}^{t_1} \big[\mathscr{H}^{\mathcal{Z}}\big(t,x_*(t),x_*(t-\delta_x),p(t)\big)
                           -H^{\mathcal{Z}}\big(t,x(t),x(t-\delta_x),u(t),u(t-\delta_x),p(t),1\big)\big] \, dt \\
&\geq& \int_{t_0}^{t_1} \big[\mathscr{H}^{\mathcal{Z}}\big(t,x_*(t),x_*(t-\delta_x),p(t)\big)
              -\mathscr{H}^{\mathcal{Z}}\big(t,x(t),x(t-\delta_x),p(t)\big)\big] \, dt \\
&\geq& \int_{t_0}^{t_1} \big(\langle -a_1(t),x(t)-x_*(t)\rangle + \langle -a_2(t),x(t-\delta_x)-x_*(t-\delta_x)\rangle \big) \, dt \\
&=& \int_{t_0}^{t_1} \langle \underbrace{-a_1(t)-\chi_{[t_0,t_1-\delta_x]}(t) \cdot a_2(t+\delta_x)}_{=\dot{p}(t)},x(t)-x_*(t)\rangle \, dt.
\end{eqnarray*}
Abschließend können wir festhalten:
\begin{eqnarray*}
\lefteqn{J\big(x(\cdot),u(\cdot)\big)-J\big(x_*(\cdot),u_*(\cdot)\big)=\int_{t_0}^{t_1} \big[ f\big(t,x(t),x(t-\delta_x),u(t),u(t-\delta_x)\big)} \\
& & \hspace*{50mm}  -f\big(t,x_*(t),x_*(t-\delta_x),u_*(t),u_*(t-\delta_x)\big)\big] \, dt \\
&=& \int_{t_0}^{t_1} \big[H^{\mathcal{Z}}\big(t,x_*(t),x_*(t-\delta_x),u_*(t),u_*(t-\delta_x),p(t),1\big) \\
& & \hspace*{3cm} -H^{\mathcal{Z}}\big(t,x(t),x(t-\delta_x),u(t),u(t-\delta_x),p(t),1\big)\big] \, dt  \\
& & + \int_{t_0}^{t_1} \langle p(t), \dot{x}(t)-\dot{x}_*(t) \rangle dt \\
&\geq& \int_{t_0}^{t_1}  \big[\langle \dot{p}(t),x(t)-x_*(t)\rangle + \langle p(t), \dot{x}(t)-\dot{x}_*(t) \rangle \big] \, dt \\
&=& \langle p(t_1),x(t_1)-x_*(t_1)\rangle-\langle p(t_0),x(t_0)-x_*(t_0)\rangle \geq 0
\end{eqnarray*}
für alle $\big(x(\cdot),u(\cdot)\big) \in \mathscr{B}^{\mathcal{Z}}_{\rm adm}$ mit $\|x(\cdot)-x_*(\cdot)\|_\infty < \gamma$. \hfill $\blacksquare$

\begin{beispiel}
{\rm Die Hamilton-Funktion im Beispiel \ref{BeispielInvDelay},
$$\mathscr{H}^{\mathcal{Z}}(t,K,L,u,v,p) = \sup_{u,v \in [0,1]} \{ p vL +(1-u)K\} =pL+K,$$
ist für $p,K,L \geq 0$ konkav und die notwendigen Optimalitätsbedingungen sind in diesem Beispiel gleichzeitig hinreichend. \hfill $\square$}
\end{beispiel}

%% file: 3-46-Chemoimmuntherapie.tex
\subsubsection{Optimale Steuerung einer Chemoimmuntherapie}\index{Chemoimmuntherapie}
Eine Chemoimmuntherapie ist die Kombination von Chemotherapie und Immuntherapie. \\
Bei einer Chemotherapie verwendet man verschiedene Medikamente (Zytostatika),
um Krebszellen abzutöten oder das Wachstum der Krebszellen zu verlangsamen.
Bei der Behandlung bösartiger Tumorerkrankungen nutzen die meisten Zytostatika die schnelle Teilungsfähigkeit der Tumorzellen,
da diese empfindlicher als gesunde Zellen auf Störungen der Zellteilung reagieren.
Auf gesunde Zellen mit ähnlich guter Teilungsfähigkeit üben sie allerdings eine ähnliche Wirkung aus,
wodurch sich Nebenwirkungen wie Haarausfall oder Durchfall einstellen können.
Immuntherapien beinhalten Behandlungen,
die das Immunsystem im Kampf gegen Krebs stimulieren oder stärken.
Effektorzellen wie B-Lymphozyten (kurz B-Zellen) gehören zu den Leukozyten (weiße Blutkörperchen),
die Plasmazellen  bilden, die wiederum Antikörper bilden. \\[2mm]
Das folgende Model wurde von Rihan et.\,al. \cite{Rihan} vorgeschlagen.
Die folgenden numerischen Untersuchungen wurden von Göllmann \& Maurer in \cite{GoMa} unternommen.
Im Rahmen des Modells tretendie folgenden Größen auf:
\begin{center}{\small \begin{tabular}{|ll|}
\hline
Variable & Beschreibung \\
\hline
$E$ & Konzentration der Effektorzellen (B-Zellen) \\
$K$ & Konzentration der Krebszellen \\
$G$ & Konzentration der gesunden Zellen \\
$Z$ & Konzentration der Zytostatika der Chemotherapie \\
$u_1$ & Dosierung der Chemotherapie \\
$u_2$ & Dosierung der Immuntherapie zur Bildung von Effektorzellen \\
\hline
\end{tabular} }\end{center}
Die Dynamik der Therapie besteht aus dem nachstehenden nicht-linearen Differentiagleichungssystem,
in dem die Entwicklungen der verschiedenen Zellpopulationen miteinander verknüpft sind:
\begin{eqnarray*}
\dot{E}(t) &=& \sigma + \bigg( \frac{\varrho}{\eta+ K(t-\Delta)}- \mu_E \bigg) E(t-\Delta) K(t-\Delta) \\
           & & -\big(\vartheta + a_1(1-e^{-Z(t)})\big)E(t)+r_1u_2(t-\delta), \\
\dot{K}(t) &=& \Big(r_2(1-\beta_K K(t))- n_K E(t) - c_1 G(t)-a_2 (1-e^{-Z(t)})\Big) K(t), \\
\dot{G}(t) &=& \Big(r_3(1-\beta_G G(t))-c_2K(t)-a_3 (1-e^{-Z(t)})\Big) G(t), \\
\dot{Z}(t) &=& u_1(t)-d Z(t).
\end{eqnarray*}
Im Vergleich zu \cite{Rihan} beinhaltet die Gleichung für $E(\cdot)$ nach \cite{GoMa} die verzögerte Steuerung $u_2(t-\delta)$,
da das menschliche Immunsystem auf die Immuntherapie mit einer gewissen Verzögerung reagiert.
Ferner ist in dem verzögerten Differentialgleichungssystem die realistische Annahme getroffen,
dass vor Beginn der Immuntherapie keine extern erwirkte Bildung von Effektorzellen voliegt,
d.\,h. $u_2(t)=0$ für $t<0$ ist. \\[1mm]
Den Erwartungen entsprechend besteht das Ziel einer Chemoimmuntherapie in der Minimierung der Krebszellenpopulation.
Gleichzeitig muss aber der Einsatz der Chemo- und Immuntherapie einbezogen werden,
sowie die Bildung der Effektorzellen. \\[2mm]
In \cite{GoMa} wird die Minimierung des Zielfunktionals   
$$J_1=\int_0^T \big( K(t)-E(t)+B_1u_1(t)+B_2 u_2(t)\big) \, dt \to \inf$$
vorgeschlagen und anschlie"send die Ergebnisse beim Einsatz des Funktionals 
$$J_2=\int_0^T \big( K(t)-E(t)+B_1u^2_1(t)+B_2 u^2_2(t)\big) \, dt \to \inf,$$
welches in \cite{Rihan} Anwendung findet, verglichen. \\[1mm]
Die numerische Auswertung in \cite{GoMa} basiert auf folgenden Paramaterwerten:
\begin{center}{\small \begin{tabular}{|lll|}
\hline 
Parameter    & Beschreibung & Parameterwerte \\
\hline
$T$ & Dauer der Therapie & $30$ Tage \\
$(E_0,K_0,G_0,Z_0)$ & Anfangsbedingungen & $(0,\!3;\, 300;\, 0,\!9;\,0,\!0)$ \\
$\Delta$ & Verzögerung der Zustände $E$, $K$ & $1,\!5$ Tage \\
$\delta$ & Verzögerung der Steuerung $u_2$ & $3,\!0$ Tage \\
$U$ & Steuerbereich für $u_1$ und $u_2$ & $[0,1]$ \\
$(a_1,a_2,a_3)$ & Reaktion auf die Zellsterberate & $(0,\!2;\, 0,\!4; \, 0,\!1)$ \\
$(\beta_K,\beta_G)$ & gegenseitige Tragfähigkeit von Krebszellen und & \\
                  & gesunden Zellen & $(0,\!002;\,1,\!0)$ \\
$(c_1,c_2)$ & Skalierungsparameter & $(3 \cdot 10^{-5}, 3 \cdot 10^{-8})$ \\
$d$ & Medikamentenabbaurate & $0,\!01$ \\
$\vartheta$ & Absterberate der Immunzellen & $0,\!2$ \\
$\eta$ & Anstieg der Immunantwort & $0,\!3$ \\
$\mu_E$ & Reduktionsrate der uninfizierten Effektorzellen & $0,\!003611$ \\
$(\sigma,\varrho)$ & Zuwachs- und Abnahmerate der Immunzellen & $(0,\!2;\, 0,\!2)$ \\
$(r_1,r_2,r_3)$ & Zellwachstumsfaktoren & $(0,\!3;\, 1,\!03;\, 1,\!0)$ \\
$n_K$ & Abnahmerate der immunen Effektorzellen & $1,\!0$ \\
$(B_1,B_2)$ & Bewertungsfaktoren & $(5;\,10)$ \\
\hline
\end{tabular} }\end{center}

\begin{figure}[h]
	\centering
	\fbox{\includegraphics[width=10cm]{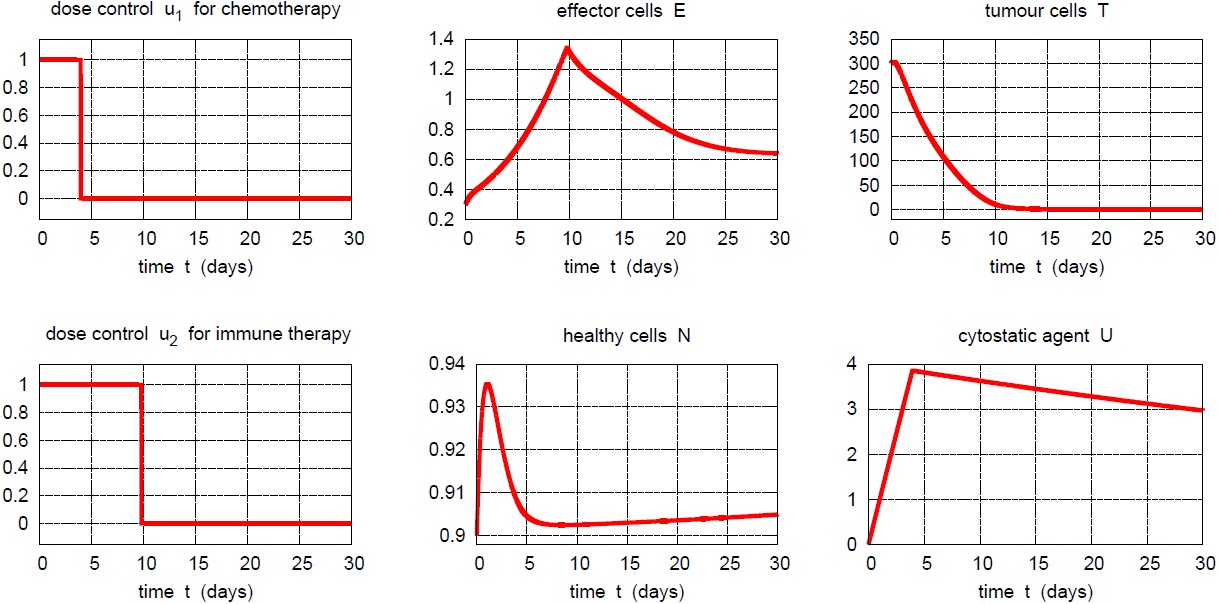}}
\caption[Chemoimmuntherapie - unverzögertes Problem]{Darstellung der optimalen Lösung in \cite{GoMa} im unverzögerten Problem.}
	\label{AbbChemo1}
\end{figure}

\begin{figure}[h]
	\centering
	\fbox{\includegraphics[width=10cm]{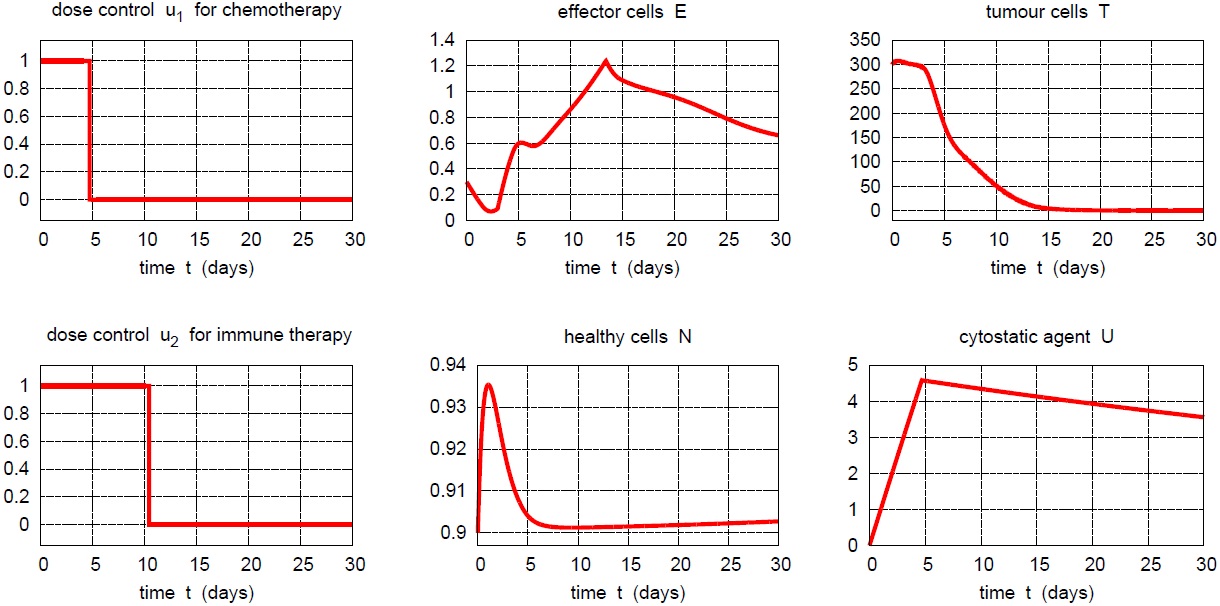}}
	\caption[Chemoimmuntherapie - verzögertes Problem]{Darstellung der optimalen Lösung in \cite{GoMa} im verzögerten Problem.}
	\label{AbbChemo2}
\end{figure}

Während die Optimierung bezüglich $J_2$ stetige Steuerungen $u_1(\cdot)$, $u_2(\cdot)$ liefert,
führt die Optimierung bezüglich dem Kriterium $J_1$ zu ``bang-bang'' Steuerungen.
Die numerischen Lösungen in \cite{GoMa} im unverzögerten Problem ($\Delta=\delta=0$\,Tage)
bzw. im verzögerten Problem ($\Delta=1,\!5$\,Tage, $\delta=3,\!0$\,Tage)
sind in den Abbildungen \ref{AbbChemo1} und \ref{AbbChemo2} dargestellt. \\[2mm]
Die Steuerungen $u_1(\cdot)$ und $u_2(\cdot)$ in den Abbildungen \ref{AbbChemo1} und \ref{AbbChemo2} implizieren
eine blockweise Verabreichung der Medikamente bzw. eine Impulstherapie
(therapeutischer ``Impuls'' durch kurze Gabe eines Medikamentes);
die Dosierungen $u_1(\cdot)$, $u_2(\cdot)$ spiegeln damit die gängige medizinische Praxis wider.

%% file: 3-50-Differentialspiele.tex
\subsection{Differentialspiele} \label{KapitelDifferentialspiele} \index{Differentialspiel}\index{Nash-Gleichgewicht}
In vielen Problemstellungen sind verschiedene Entscheidungsträger involviert,
deren Interessen nicht im Einklang stehen müssen.
Eine derartige Situation bezeichnet man als mathematisches Spiel. 
Unterliegen die Zielkriterien der Konfliktgruppen dynamischen Nebenbedingungen,
so spricht man von einem Differentialspiel. \\
Als Wegbereiter der Differentialspiele gilt Rufus Isaacs (1914--1981),
dessen Untersuchungen von Differentialspielen bei der Rand Corporation begannen und
erste Ergebnisse in der Arbeit \cite{Isaacs0} im Jahr 1951 angegeben sind. \\
Die Entwicklung der Differentialpiele erfolgte unabhängig vom Pontrjaginschen Maximumprinzip.
Vielfach werden Differentialspiele, in denen anstatt einem mehrere Entscheidungsträger agieren,
als Verallgemeinerung von Steuerungsproblemen angesehen (vgl. \cite{Feichtinger}).
Erst nach der Veröffentlichung der Monographie \cite{Isaacs} von Isaacs (1963) wurde man sich der Verbindung zwischen
dynamischen Spielen und der Steuerungstheorie bewusst. \\[2mm]
Wir betrachten eine dynamisches Spiel,
in dem die Steuerungen $u_1(\cdot)$ und $u_2(\cdot)$ die Einflussnahmen durch den ersten bzw. zweiten Spieler widergeben.
Wir gelangen zur Spielsituation durch das Hinzufügen eines weiteren Zielfunktionals in der Standardaufgabe:
\begin{eqnarray*}
&& J_1\big(x(\cdot),u_1(\cdot),u_2(\cdot)\big) = \int_{t_0}^{t_1} f_1\big(t,x(t),u_1(t),u_2(t)\big) \, dt \to \inf, \\
&& J_2\big(x(\cdot),u_1(\cdot),u_2(\cdot)\big) = \int_{t_0}^{t_1} f_2\big(t,x(t),u_1(t),u_2(t)\big) \, dt \to \inf, \\
&& \dot{x}(t) = \varphi\big(t,x(t),u_1(t),u_2(t)\big), \quad x(t_0)=x_0, \; x(t_1)=x_1, \quad
    u_1(t) \in U_1, \; u_2(t) \in U_2.
\end{eqnarray*}
Ein L"osungskonzept zur Behandlung eines solchen Spiels ist das
Nash-Gleichgewicht\index{Nash-Gleichgewicht} f"ur nichtkooperative Spiele.
Neben diesem existieren weitere Gleichgewichtskonzepte. \\
Beispielsweise besitzt bei einem Stackelberg-Gleichgewicht ein Spieler eine dominante Verhandlungsposition
(``Stackelbergführer''),
während der andere Spieler (``Stackelbergfolger'') die bestmögliche Antwort aus der unterlegenen Position heraus finden muss.
Der Stackelbergführer muss bei seiner Entscheidungsfindung die Antwort nach Nachfolgers berücksichtigen,
denn die einseitige Optimierung kann im Extremalfall zu einer destruktiven Antwort führen. \\[2mm]
In unserem Differentialspiel ist ein Nash-Gleichgewicht ein Tripel
$\big(x^*(\cdot),u_1^*(\cdot),u_2^*(\cdot)\big)$ mit
\begin{eqnarray*}
    J_1\big(x^*(\cdot),u_1^*(\cdot),u_2^*(\cdot)\big)
&=& \min_{(x(\cdot),u_1(\cdot))} J_1\big(x(\cdot),u_1(\cdot),u_2^*(\cdot)\big), \\
    J_2\big(x^*(\cdot),u_1^*(\cdot),u_2^*(\cdot)\big)
&=& \min_{(x(\cdot),u_2(\cdot))} J_2\big(x(\cdot),u_1^*(t),u_2(\cdot)\big).
\end{eqnarray*}
Befindet sich also das Spiel im Nash-Gleichgewicht $\big(x^*(\cdot),u_1^*(\cdot),u_2^*(\cdot)\big)$,
so hat kein Spieler einen Anreiz seine Strategie zu ändern,
denn das bestmögliche Ergebnis zur Strategie des Konkurrenten wurde gefunden. \\
Ausf"uhrlich geschrieben ist das Nash-Gleichgewicht ein System von zwei gekoppelten Steuerungsproblemen.
Für den ersten Spieler ergibt sich das Steuerungsproblem
\begin{eqnarray*}
&& J_1\big(x(\cdot),u_1(\cdot),u^*_2(\cdot)\big) = \int_{t_0}^{t_1} f_1\big(t,x(t),u_1(t),u^*_2(t)\big) \, dt \to \inf, \\
&& \dot{x}(t) = \varphi\big(t,x(t),u_1(t),u^*_2(t)\big), \quad x(t_0)=x_0, \; x(t_1)=x_1, \quad
   u_1(t) \in U_1
\end{eqnarray*}
bei gegebener Strategie $u_2^*(\cdot)$ des zweiten Spielers; entsprechend für den zweiten Spieler:
\begin{eqnarray*}
&& J_2\big(x(\cdot),u_1^*(\cdot),u_2(\cdot)\big) = \int_{t_0}^{t_1} f_2\big(t,x(t),u^*_1(t),u_2(t)\big) \, dt \to \inf, \\
&& \dot{x}(t) = \varphi\big(t,x(t),u^*_1(t),u_2(t)\big), \quad x(t_0)=x_0, \; x(t_1)=x_1, \quad
   u_2(t) \in U_2.
\end{eqnarray*}
Formal sind die Abbildungen $f_1\big(t,x,u_1,u^*_2(t)\big)$, $f_2\big(t,x,u^*_1(t),u_2\big)$ und
$\varphi\big(t,x,u_1,u^*_2(t)\big)$, $\varphi\big(t,x,u^*_1(t),u_2\big)$ nicht stetig bezüglich der Variable $t$
und verletzen die Anforderungen an die Standardaufgabe in den ersten Kapiteln.
Sind die Steuerungen $u_1^*(\cdot)$ und $u_2^*(\cdot)$ zumindest stückweise stetig,
so bleiben die Optimalitätsbedingungen in Form des Schwachen Optimalitätsprinzips und des Pontrjaginschen Maximumprinzips
gültig.

%% file: 3-51-Kapitalismusspiel.tex
\subsubsection{Kapitalismusspiel} \label{AbschnittKapitalismusspiel}\index{Kapitalismusspiel}
Wir betrachten ein Differentialspiel zwischen Arbeitern und Unternehmern.
Beide Parteien sind bestrebt ihr Konsumverlangen zu befriedigen.
Zudem zeichnen sich die Unternehmer für Investitionen verantwortlich.
F"ur die Arbeiter stellt sich dabei das Problem,
inwieweit sie den produzierten Ertrag konsumieren oder den Unternehmern "uberlassen sollen,
damit eine k"unftige hohe G"uterproduktion erm"oglicht wird, von der auch die Arbeiter wieder profitieren.
Das Dilemma, das sich dem Arbeiter stellt, ist,
dass sie keine Garantien "uber ausreichende Neuinvestitionen seitens der Unternehmer haben. \\
Die Unternehmer stehen ihrerseits vor der Frage,
wie sie mit dem verbliebenen Teil des Ertrages, der nicht dem Arbeiter zugesprochen wird, umgehen:
Sollen sie diesen investieren oder konsumieren? \\
Die Spielsituation ensteht dadurch,
dass der Nutzen f"ur den Arbeiter und den Unternehmern durch den aufgeteilten Ertrag miteinander gekoppelt ist.
D.\,h. der Gewinn eines jeden Spielers h"angt von der Entscheidung des anderen Spielers ab. \\[2mm]
In diesem Modell werden die Geld- und G"uterwerte im Kapitalstock $K(\cdot)$ zusammengefasst.
Weiter bezeichne $u(t)$ denjenigen relativen Anteil an der Produktion, der dem Arbeiter zum Zeitpunkt $t$ zugesprochen wird.
Von dem verbliebenen Teil kann der Unternehmer mit der Rate $v(t)$ zur Zeit $t$ investieren.
Mit den Zielfunktionalen $W$ bzw. $C$ f"ur den Arbeiter bzw. Unternehmer entsteht damit die folgende Aufgabe
\begin{eqnarray}
&&\label{KapSpiel1} J^W\big(K(\cdot),u(\cdot),v(\cdot)\big) = \int_0^T u(t)K(t) \, dt \to \sup,\\
&&\label{KapSpiel2} J^C\big(K(\cdot),u(\cdot),v(\cdot)\big) = \int_0^T \big(1-v(t)\big)\big(1-u(t)\big)K(t) \, dt \to \sup, \\
&&\label{KapSpiel3} \dot{K}(t) = v(t)\big(1-u(t)\big)K(t), \qquad K(0)=K_0>0, \quad K(T) \mbox{ frei}, \\
&&\label{KapSpiel4} \big(u(t),v(t)\big) \in [a,b] \times [0,1], \quad 0<a<b<1, \quad b > \frac{1}{2},\quad T>\frac{1}{1-b}.\hspace*{5mm}
\end{eqnarray}
Das Modell geht auf Lancaster \cite{Lancaster} zur"uck.
Wir folgen der Darstellung und untersuchen das Nash-Gleichgewicht und die Kollusionsl"osung,
die wir abschlie"send miteinander vergleichen.
Wir verweisen au"serdem auf die umfassenden Ausf"uhrungen in
Feichtinger \& Hartl \cite{Feichtinger} und Seierstad \& Syds\ae ter \cite{Seierstad}. \\[2mm]
Wir platzieren die zwei Spieler Arbeiter und Unternehmer in einer Spielsituation und benutzen als L"osungskonzept das
Nash-Gleichgewicht\index{Nash-Gleichgewicht} f"ur nichtkooperative Spiele.
Im vorliegenden Differentialspiel ist ein Nash-Gleichgewicht ein Tripel
$$\big(K_*(\cdot),u_*(\cdot),v_*(\cdot)\big),$$
f"ur das die folgenden Gleichungen erf"ullt sind:
\begin{eqnarray*}
J^W\big(K_*(t),u_*(t),v_*(t)\big) &=& \max_{(K(\cdot),u(\cdot))} J^W\big(K(\cdot),u(\cdot),v_*(t)\big), \\
J^C\big(K_*(t),u_*(t),v_*(t)\big) &=& \max_{(K(\cdot),v(\cdot))} J^C\big(K(\cdot),u_*(t),v(\cdot)\big).
\end{eqnarray*}
Ausf"uhrlich geschrieben ist das Nash-Gleichgewicht das folgende System von zwei gekoppelten Steuerungsproblemen
f"ur den Arbeiter,
\begin{equation} \label{KapSpiel5} \left. \begin{array}{l}
J^W\big(K(\cdot),u(\cdot),v_*(\cdot)\big) = \displaystyle\int_0^T u(t)K(t) \, dt  \to \sup, \\[1mm]
\dot{K} =   v_*(t)\big(1-u(t)\big)K(t), \quad K(0)=K_0>0, \quad u \in [a,b], \\[1mm]
\end{array} \right\}
\end{equation}
und für den Unternehmer,
\begin{equation} \label{KapSpiel6} \left. \begin{array}{l}
J^C\big(K(\cdot),u_*(\cdot),v(\cdot)\big)
 = \displaystyle\int_0^T \big(1-v(t)\big)\big(1-u_*(t)\big)K(t) \, dt \to \sup, \\[1mm]
\dot{K} = v(t)\big(1-u_*(t)\big)K(t), \quad K(0)=K_0>0, \quad v \in [0,1].
\end{array} \right\}
\end{equation} 
Auf beide Steuerungsprobleme wenn wir die Optimalitätsbedingungen im Pontrjaginsche Maximumprinzip
unter Beachtung der Strategie des Gegenspielers an: \\
In der Aufgabe der Arbeiter (\ref{KapSpiel5}) erhalten wir die Maximumbedingung
$$H^W\big(t,K_*(t),u_*(t),p(t),1\big) = \max_{u \in [a,b]} \Big\{ \big(1-p(t)v_*(t)\big)uK_*(t) \Big\} + p(t)v_*(t)K_*(t).$$
Da jede zul"assige Trajektorie $K(\cdot)$ stets positiv ist, gilt
\begin{equation} \label{KapSpiel7} \left. 
   \begin{array}{ll} 
      u_*(t)=a,         & \mbox{ wenn } p(t)v_*(t)>1, \\
      u_*(t)=b,         & \mbox{ wenn } p(t)v_*(t)<1, \\
      u_*(t) \in [a,b], & \mbox{ wenn } p(t)v_*(t)=1.    
    \end{array} \right\}
\end{equation}
Weiterhin gen"ugt $p(\cdot)$ der adjungierten Gleichung (\ref{PMPeinfach4}) und der Bedingung (\ref{PMPeinfach5}):
\begin{equation} \label{KapSpiel8}
\dot{p}(t) = -v_*(t)\big(1-u_*(t)\big)p(t) - u_*(t), \qquad p(T)=0.   
\end{equation}
In Gleichung (\ref{KapSpiel8}) erkennt man unmittelbar $\dot{p}(t)\leq  - u_*(t)\leq -a$ f"ur $t \in (0,T)$.
Daher gilt $p(t)>0$ f"ur $t \in [0,T)$. \\[2mm]
Entsprechend ergibt sich in der Aufgabe des Unternehmers (\ref{KapSpiel6}) die Maximumbedingung
$$H^C\big(t,K_*(t),v_*(t),q(t),1\big)=
  \max_{v \in [0,1]} \Big\{ \big(q(t)-1\big)v\big(1-u_*(t)\big)K_*(t) \Big\} + \big(1-u_*(t)\big)K_*(t)$$
und es gelten daher
\begin{equation} \label{KapSpiel9} \left. 
   \begin{array}{ll} 
      v_*(t)=0         & \mbox{ wenn } q(t)<1, \\
      v_*(t)=1         & \mbox{ wenn } q(t)>1, \\
      v_*(t) \in [0,1] & \mbox{ wenn } q(t)=1.    
    \end{array} \right\}
\end{equation}
Ferner gen"ugt $q(\cdot)$ der adjungierten Gleichung zur Transversalit"atsbedingung:
\begin{equation} \label{KapSpiel10}
\dot{q}(t) = -v_*(t)\big(1-u_*(t)\big)q(t) - \big(1-v_*(t)\big)\big(1-u_*(t)\big),\qquad q(T)=0.
\end{equation}
Zusammen erhalten wir aus (\ref{KapSpiel9}) und (\ref{KapSpiel10}) 
$$\dot{q}(t) = -\big(1-u_*(t)\big) \cdot \max\{1,\,q(t)\}<0.$$
Aufgrund der strengen Monotonie existiert eine eindeutig bestimmte L"osung von
$$\tau= \min \{ t \in [0,T] \,|\, q(t) \leq 1 \}.$$
Wegen (\ref{KapSpiel9}) gilt $v_*(t) = 0$ f"ur alle $t \in [\tau,T]$ und aus (\ref{KapSpiel7}) folgt nun $u_*(t) = b$ auf $[\tau,T]$.
Auf $[\tau,T]$ erhalten wir daraus
$$q(t)=(1-b)(T-t), \quad q(\tau)=1, \quad p(t)=b(T-t),\quad p(\tau)=\frac{b}{1-b}.$$
Weiterhin folgt nach (\ref{KapSpiel4}) und aus der Bedingung $q(\tau)=1$
\begin{equation} \label{KapSpiel11}
\tau=T-\frac{1}{1-b}>0.
\end{equation}
Aufgrund $b \geq 1/2$ in (\ref{KapSpiel4}) ist $p(\tau)\geq 1$.
Da $p(\cdot)$ streng monoton fallend ist,
erhalten wir $u_*(t)=a$ und $v_*(t)=1$ f"ur alle $t \in [0,\tau)$. \\[2mm]
Zusammenfassend lauten im Nash-Gleichgewicht die jeweiligen Strategien
$$u_*(t)=\left\{\begin{array}{l} a,\; t \in [0,\tau), \\ b,\; t \in [\tau,T], \end{array}\right. \quad
  v_*(t)=\left\{\begin{array}{l} 1,\; t \in [0,\tau), \\ 0,\; t \in [\tau,T], \end{array}\right.$$
zum Kapitalstock 
$$K_*(t)=\left\{\begin{array}{l} K_0e^{(1-a)t},\; t \in [0,\tau), \\ K_0e^{(1-a)\tau},\; t \in [\tau,T]. \end{array}\right.$$
F"ur die Zielfunktionale ergeben sich
\begin{equation} \label{KapSpiel12} \left. 
\begin{array}{l}
J^W\big(K_*(t),u_*(t),v_*(t)\big)
   = \displaystyle\left(\frac{a}{1-a}+\frac{b}{1-b}\right)K_0 e^{(1-a)\tau} - \frac{a}{1-a}K_0, \\[5mm]
J^C\big(K_*(t),u_*(t),v_*(t)\big) = \displaystyle K_0 e^{(1-a)\tau}. 
\end{array} \right\}
\end{equation}
Die Untersuchung des Nash-Gleichgewichtes ist damit abgeschlossen. \\[2mm]
Im Gegensatz zur Spielsituation besprechen sich die Arbeiter und Unternehmer,
wie sie gemeinsam am besten agieren k"onnen.
F"ur das gemeinsame Zielfunktional ergibt sich
$$\int_0^T u(t)K(t) \, dt + \int_0^T \big(1-v(t)\big)\big(1-u(t)\big)K(t) \, dt =
     \int_0^T \big[1- v(t)\big(1-u(t)\big)\big]K(t) \, dt.$$
F"uhren wir die neue Steuerungsvariable $w(\cdot)$ durch
$$w(t)=v(t)\big(1-u(t)\big), \qquad w(t) \in [0,1-a],$$
ein, so erhalten wir die Aufgabe 
\begin{eqnarray*}
&& J\big(K(\cdot),w(\cdot)\big) = \int_0^T \big(1-w(t)\big)K(t) \, dt \to \sup, \\
&& \dot{K}(t) = w(t)K(t), \qquad K(0)=K_0>0, \quad K(T) \mbox{ frei}, \\
&& w(t) \in [0,1-a], \quad 0<a<b<1, \quad b > \frac{1}{2},\quad T>\frac{1}{1-b}. \hspace*{5mm}
\end{eqnarray*}
Die kooperative Lösung ist
$$w_*(t)=\left\{\begin{array}{ll} 1-a,& t \in [0,T-1), \\ b,& 0 \in [T-1,T], \end{array}\right. \qquad
  K_*(t)=\left\{\begin{array}{ll} K_0e^{(1-a)t},& t \in [0,T-1), \\ K_0e^{(1-a)(T-1)},& t \in [T-1,T], \end{array}\right.$$
mit folgendem Optimalwert f"ur das Zielfunktional
\begin{equation} \label{KapSpiel13}
J\big(K_*(t),w_*(t)\big) = \frac{1}{1-a}K_0 e^{(1-a)(T-1)} - \frac{a}{1-a}K_0.
\end{equation}
Der Gesamtnutzen im Nash-Gleichgewicht nach (\ref{KapSpiel12}) ist gleich
\begin{eqnarray}
\lefteqn{J^W\big(K_*(t),u_*(t),v_*(t)\big) + J^C\big(K_*(t),u_*(t),v_*(t)\big)} \nonumber \\
&&\label{KapSpiel14} = \left(\frac{1}{1-a}+\frac{b}{1-b}\right)K_0 e^{(1-a)\tau} - \frac{a}{1-a}K_0.
\end{eqnarray}
Zum Vergleich der Werte (\ref{KapSpiel13}) und (\ref{KapSpiel14}) verwenden wir die Beziehungen
$$T-1=\tau+\frac{b}{1-b}.$$
Damit erhalten wir
\begin{eqnarray*}
J\big(K_*(t),w_*(t)\big) &>& J^W\big(K_*(t),u_*(t),v_*(t)\big) + J^C\big(K_*(t),u_*(t),v_*(t)\big) \\
\Leftrightarrow \qquad \frac{1}{1-a}e^{(1-a)\frac{b}{1-b}} &>& \frac{1}{1-a}+\frac{b}{1-b}.
\end{eqnarray*}
Benutzen wir die elementare Ungleichung $e^x>1+x$ f"ur alle $x>0$ folgt abschlie"send
$$\frac{1}{1-a}e^{(1-a)\frac{b}{1-b}}>\frac{1}{1-a}\left(1+(1-a)\frac{b}{1-b}\right)= \frac{1}{1-a}+\frac{b}{1-b}.$$
D.\,h. der Gesamtnutzen ist im sozialen Optimum gr"o"ser als im Wettbewerb.
Aber es bleibt die spannende Frage nach der ``gerechten'' Aufteilung des Nutzens auf Unternehmer und Arbeiter
an dieser Stelle offen. 

%% file: 3-52-Fischerei.tex
\subsubsection{Ein Fischerei-Differentialspiel}
\label{BeispielDockner} {\rm Wir betrachten nach Dockner et\,al. \cite{Dockner} das Differentialspiel
\index{Fischereimodell}
\begin{eqnarray*}
&& \tilde{J}_i\big(x(\cdot),u_1(\cdot),u_2(\cdot)\big) =\int_0^T e^{-\varrho t}\big(p x(t)-c_i\big)u_i(t) \, dt \to \sup, \\
&& \dot{x}(t)=x(t) \big(\alpha-r\ln x(t) \big) -u_1(t)x(t)-u_2(t)x(t), \quad x(0)=x_0>0, \\
&& u_i > 0, \quad \alpha,c_i,p,r,\varrho >0, \quad \alpha> \frac{1}{c_1+c_2}, \quad i=1,2.
\end{eqnarray*}
In dieser Aufgabe sei der Preis $p$ nicht konstant, sondern umgekehrt proportional zum Angebot:
$$p=p(u_1x+u_2x)= \frac{1}{u_1x+u_2x}.$$
Dieser Ansatz spiegelt die "Okonomie einer ``Eskimo''-Gesellschaft wider,
in der der Fischbestand die wichtigste Nahrungsgrundlage darstellt und kein echtes Ersatzprodukt existiert.
Unter diesen Umst"anden f"uhrt eine prozentuale Preissteigerung zu einem Umsatzr"uckgang in gleicher Relation. \\[2mm]
Wenden wir die Transformation $z=\ln x$ an, so erhalten wir das Spielproblem
\begin{eqnarray*}
&& J_i\big(z(\cdot),u_1(\cdot),u_2(\cdot)\big) =\int_0^T e^{-\varrho t}\bigg(\frac{1}{u_1(t)+u_2(t)}-c_i\bigg)u_i(t) \, dt \to \sup, \\
&& \dot{z}(t)=-r z(t)+\alpha-u_1(t)-u_2(t), \quad z(0)=\ln x_0>0, \\
&& u_i>0, \quad \alpha,c_i,p,r,\varrho >0, \quad \alpha> \frac{1}{c_1+c_2}, \quad i=1,2.
\end{eqnarray*}
Ein zul"assiger Steuerungsprozess $\big(z^*(\cdot),u^*_1(\cdot),u^*_2(\cdot)\big)$ ist ein Nash-Gleichgewicht des Spiels,
falls f"ur alle anderen zul"assigen Steuerungsprozesse $\big(z(\cdot),u_1(\cdot),u^*_2(\cdot)\big)$,
$\big(z(\cdot),u^*_1(\cdot),u_2(\cdot)\big)$ die Ungleichungen
$$J_1\big(z^*(\cdot),u^*_1(\cdot),u^*_2(\cdot)\big) \geq J_1\big(z(\cdot),u_1(\cdot),u^*_2(\cdot)\big), \;
J_2\big(z^*(\cdot),u^*_1(\cdot),u^*_2(\cdot)\big) \geq J_2\big(z(\cdot),u^*_1(\cdot),u_2(\cdot)\big)$$
gelten.
Halten wir die optimale Strategie des Gegenspielers fest,
dann ergeben sich f"ur $i,j \in \{1,2\}$, $i \not= j$, folgende miteinander gekoppelte Steuerungsprobleme:
\begin{eqnarray*}
&& J_i\big(z(\cdot),u_i(\cdot),u^*_j(\cdot)\big) =
      \int_0^T e^{-\varrho t}\bigg(\frac{1}{u_i(t)+u^*_j(t)}-c_i\bigg)u_i(t) \, dt \to \sup, \\
&& \dot{z}(t)=-r z(t)+\alpha-u_i(t)-u^*_j(t), \quad z(0)=\ln x_0>0, \quad u_i>0.
\end{eqnarray*}
Die zul"assigen Steuerungen geh"oren dem Raum $L_\infty([0,T],U)$ an.
Deshalb sind im Ansatz des Nash-Gleichgewichtes die Abbildungen
$$f_i(t,z,u_i)=\bigg(\frac{1}{u_i+u^*_j(t)}-c_i\bigg)u_i, \quad \varphi_i(t,z,u_i)=-r z+\alpha-u_i-u^*_j(t)$$
nicht stetig bez"uglich der Variable $t$.
Die Voraussetzungen an die Aufgabe sind nicht erf"ullt.
Da aber $\varphi_i(t,z,u_i)$ linear in $z$ und $u_i$ sind und die Variable $z$ in $f_i(t,z,u_i)$ nicht einflie"st,
sind die gleichm"a"sigen Eigenschaften der Voraussetzungen f"ur fast alle $t \in [0,T]$ erf"ullt.
Unter diesen Rahmenbedingungen beh"alt Theorem \ref{SatzPMP} seine G"ultigkeit. \\[2mm]
Die Pontrjagin-Funktionen lauten
$$H_i(t,z,u_i,p_i,\lambda_0) = p_i\big(-r z+\alpha-u_i-u^*_j(t)\big)+\lambda_0 e^{-\varrho t}\bigg(\frac{1}{u_i+u^*_j(t)}-c_i\bigg)u_i.$$
Es liegt ein freier Endpunkt vor.
Deswegen d"urfen wir $\lambda_0=1$ annehmen und es gelten die Transversalit"atsbedingungen $p_i(T)=0$.  
Damit erhalten wir aus der adjungierten Gleichung (\ref{SatzPMP1})
$$\dot{p}_i(t)=r p_i(t), \quad p_i(T)=0,\qquad\Rightarrow\qquad p_i(t) \equiv 0.$$
Aus der Maximumbedingung (\ref{SatzPMP3}) folgt schlie"slich
$$u_1^*(t)\equiv \frac{c_2}{(c_1+c_2)^2}, \qquad u_2^*(t)\equiv \frac{c_1}{(c_1+c_2)^2}$$
und f"ur die optimale Trajektorie
$$z_*(t)=(z_0-c_0)e^{-rt} + c_0, \qquad c_0=  \frac{1}{r}\bigg(\alpha-\frac{1}{c_1+c_2}\bigg).$$
Die Funktion $z_*(\cdot)$ ist streng monoton und nimmt nur Werte des Segments $[z_0,c_0]$ an.
Da $c_0$ und $z_0$ positiv sind, ist $x_*(t)=\exp\big(z_*(t)\big)$ "uber $[0,T]$ wohldefiniert und
$\big(x_*(\cdot),u^*_1(\cdot),u^*_2(\cdot)\big)$ liefert einen Kandidaten f"ur das urspr"ungliche Differentialspiel.

%% file: 4-0-StarkeVariationen.tex
\section{Starkes lokales Minimum \"uber unendlichem Zeithorizont} \label{KapitelStrong}
Wir untersuchen in diesem Kapitel Aufgaben der Optimalen Steuerung "uber unendlichem Zeithorizont unter Zustandsbeschr"ankungen und mit
Randwerten im Unendlichen.
Die erzielten Resultate beinhalten das Pontrjaginsche Maximumprinzip,
d.\,h. einen vollst"andigen Satz notwendiger Optimalit"atsbedingungen f"ur ein starkes lokales Minimum,
und die G"ultigkeit verschiedener Transversalit"atsbedingungen wie der ``nat"urlichen'' Transversalit"atsbedingungen und der Bedingung von Michel \cite{Michel}.
Au"serdem nehmen wir Bezug zu den hinreichenden Optimalit"atsbedingungen nach Arrow f"ur Aufgaben mit freiem und festem Randwert im Unendlichen und
f"ur die Aufgabe mit Zustandsrestriktionen. \\[2mm]
Die Herleitung notwendiger Optimalit"atsbedingungen f"ur starke lokale Minimalstellen in der Klassischen Variationsrechnung und
in der Optimalen Steuerung basiert vorallem auf der Nadelvariationsmethode.
F"ur Aufgaben mit einem unbeschr"ankten Zeitintervall sieht man sich dabei mit verschiedenen H"urden konfrontiert.
Zur Verdeutlichung betrachten wir diesbez"uglich zu Beginn dieses Kapitels die einfachen Nadelvariationen.
Dabei zeigen sich bereits die Schwierigkeiten bei der Variierbarkeit einer Trajektorie,
bei der Beschr"anktheit des Zielfunktionals und beim Nachweis von Transversalit"atsbedingungen. \\[2mm]
In der Folge weisen wir das Pontrjaginsche Maximumprinzip f"ur Steuerungsprobleme mit unendlichem Zeithorizont nach.
Dabei erweitern wir die mehrfachen Nadelvariationen nach Ioffe \& Tichomirov \cite{Ioffe} f"ur das unbeschr"ankte Zeitintervall.
Weiterhin gehen wir auf Zustandsbeschr"ankungen "uber dem unendlichen Zeithorizont ein.
Zu dieser Aufgabenklasse existieren bisher nur sehr wenige Beitr"age.
Wir pr"asentieren einen Zugang im Rahmen der stetigen Funktionen, die im Unendlichen einen Grenzwert besitzen.
Auf dieser Grundlage besitzen die erzielten Ergebnisse eine enge Verwandtschaft mit dem Maximumprinzip
f"ur Aufgaben mit endlichem Zeitintervall.
Die Anforderung, die wir dabei an die Aufgabe stellen,
sind "ahnlich zu denen bei Brodskii \cite{Brodskii}. \\[2mm]
In \cite{Brodskii} wird die Aufgabe mit unendlichem Zeithorizont mittels Richtungsvariationen auf schwache lokale
Minimalstellen im Rahmen beschränkter Funktionen untersucht.
Im Gegensatz dazu behandeln wir die Aufgabe mit unendlichem Zeithorizont mittels der Nadelvariationsmethode
und zeigen Optimialitätsbedingungen für starke lokale Minimalstellen.
Dazu benötigen wir einen geeigneten Rahmen und schlagen die stetigen Funktionen, die im Unendlichen konvergieren, vor.
Dementsprechend sind zusätzliche Anforderungen an die Standardaufgabe zu stellen,
die die Konvergenz zulässiger Trajektorien gewährleisten. \\[2mm]
Die Form der notwendigen Optimalitätsbedingungen in Gestalt der Pontrjaginschen Maximumprinzips geben Anlass die Resultate
über endlich und unendlichem Zeithorizont miteinander in Bezug zu setzen.
Das ist Gegenstand eines eigenständigen Abschnitts zur Einordnung der Aufgabenklassen.

%% file: 4-1-EinfacheNV.tex
\subsection{Die elementare Aufgabe mit freiem Endpunkt} \label{AbschnittPMPeinfachUH}
In der Literatur werden die Steuerungsprobleme mit unendlichem Zeithorizont häufig als eine Folgerung
der Standardaufgabe dargestellt.
Dabei sind die Hürden, die das unbeschränkte Zeitintervall mit sich führt,
von grundlegender Natur und mit den bisher zur Verfügung gestellten Methoden nicht erfassbar.
Anhand der elementaren Aufgabe mit freiem Endpunkt lassen sich diese Schwierigkeiten bereits sehr gut verdeutlichen. \\
Unsere Betrachtungen beschränken wir auf Problemstellungen,
in denen der optimale Wert des Zielfunktionals endlich ist.
Darin liegt unsere grundsätzliche Position,
dass sich bei Vorhandensein des Zielwertes ``$+\infty$'' oder ``$-\infty$'' die Suche nach einem ``Maximum'' bzw. ``Minimum'' bereits von selbst ausschließt. \\[2mm]
Unsere Untersuchungen der elementaren Aufgaben auf starke lokale Optimalstellen basieren auf der einfachen Nadelvariation\index{Nadelvariation, einfache}
$$u(t;v,\tau,\lambda) = u_{\lambda}(t) = 
  \left\{ \begin{array}{ll}
          u_*(t) & \mbox{ f"ur } t \not\in [\tau-\lambda,\tau), \\
          v      & \mbox{ f"ur } t     \in [\tau-\lambda,\tau), 
          \end{array} \right. \qquad v \in U.$$
Im Vergleich zum endlichen Zeitintervall ist die Variation mittels $u_\lambda(\cdot)$ nicht in jedem Fall möglich.
Zur Veranschaulichung betrachten wir die Dynamik
$$\dot{x}(t)=x(t)+x^2(t)+u(t), \qquad x(0)=0, \qquad u(t)\in U = [0,\infty).$$
Die Trajektorie $x_0(t) \equiv 0$ zur Steuerung $u_0(t)\equiv 0$ ist nicht variierbar.
Denn jede Steuerung $u(t;v,\tau,\lambda)$ mit $\lambda >0$ und $v>0$ liefert eine Trajektorie $x_\lambda(\cdot)$,
die als Lösung der Differentialgleichung nicht über der gesamten Halbachse $[0,\infty)$ existiert. \\
Im Beweis des Pontrjaginschen Maximumprinzips für die elementare Aufgabe sind die Abhängigkeitssätze im Anhang \ref{AbschnittDGLAbhaengigkeit} zentraler Bestandteil.
Auf das unbeschränkte Zeitintervall sind diese Resultate nicht unmittelbar übertragbar. \\[2mm]
Die einfache Nadelvariation $u_{\lambda}(\cdot)$ kann zu einer Trajektorie $x_\lambda(\cdot)$ f"uhren,
die "uber $[0,\infty)$ wohldefiniert ist,
aber das Zielfunktional keinen endlichen Wert annimmt.
In der folgenden Variante des Beispiels von Halkin \cite{Halkin},
\begin{eqnarray*}
&& J\big(x(\cdot),u(\cdot)\big) = \int_0^\infty e^{-\varrho t}\big(u(t)-x(t)\big) \, dt \to \sup, \\
&& \dot{x}(t) = u^2(t)+x(t), \quad x(0)=0, \quad  u(t) \in [0,1], \quad \varrho \in (0,1),
\end{eqnarray*}
ist der Steuerungsprozess $\big(x_*(t),u_*(t)\big) \equiv (0,0)$ das globale Maximum,
da f"ur jeden anderen Steuerungsprozess zur Steuerung $u(t) \not\equiv 0$ das Zielfunktional den Wert $-\infty$ besitzt. \\[2mm]
In unseren weiteren Untersuchungen legen wir uns auf beschränkte Steuerungsprozesse fest,
für die das Zielfunktional endlich ausfällt.
Im Rahmen der elementaren Aufgabe beschränken wir uns außerdem auf Variationen $x_\lambda(\cdot)$,
die für $\lambda \to 0^+$ gleichmäßig gegen $x_*(\cdot)$ über $[0,\infty)$ konvergieren.
Die Existenz derartiger Variationen bleibt an dieser Stelle offen,
da wir auf adäquate Abhängigkeitssätze nicht zugreifen können.
Schließlich ist zudem die Existenz der adjungierten Funktion über $[0,\infty)$ sicherzustellen.

\newpage      
Die Aufgabe mit unendlichem Zeithorizont besitzt die Gestalt
\begin{eqnarray}
&& \label{UHA1} J\big(x(\cdot),u(\cdot)\big) = \int_0^\infty e^{-\varrho t} f\big(t,x(t),u(t)\big) \, dt \to \inf, \\
&& \label{UHA2} \dot{x}(t) = \varphi\big(t,x(t),u(t)\big), \quad x(0)=x_0, \\
&& \label{UHA3} u(t) \in U \subseteq \R^m, \quad U \not= \emptyset, \quad \varrho > 0.
\end{eqnarray}
Die Aufgabe (\ref{UHA1})--(\ref{UHA3}) untersuchen wir bez"uglich der Steuerungsprozesse
$$\big(x(\cdot),u(\cdot)\big) \in PC_1([0,\infty),\R^n) \times PC([0,\infty),U).$$
Dabei sind $PC([0,\infty),\R)$ bzw. $PC_1([0,\infty),\R)$ die Räume derjenigen Funktionen,
die über $[0,\infty)$ beschränkt und die über jedem endlichen Intervall $[0,T]$ stückweise stetig bzw. stückweise stetig differenzierbar sind. \\[2mm]
Mit $\mathscr{D}^{\,\mathcal{U}}_{\rm Lip}$ bezeichnen wir die Menge aller Paare $\big(x(\cdot),u(\cdot)\big)$,
für die es ein $\gamma>0$ derart gibt,
dass die Abbildungen $f(t,x,u)$, $\varphi(t,x,u)$ für jede kompakte Teilmenge $U_1 \subset \R^m$
auf der Menge aller Punkte $(t,x,u) \in \R \times \R^n \times \R^m$ mit
$0 \leq t < \infty$, $\|x-x(t)\| < \gamma$ und $u \in U_1$
beschränkt, stetig in der Gesamtheit aller Variablen und stetig differenzierbar bezüglich $x$
mit beschränkten Ableitungen $f_x(t,x,u)$, $\varphi_x(t,x,u)$ sind. \\[2mm]
Das Paar $\big(x(\cdot),u(\cdot)\big) \in PC_1([0,\infty),\R^n) \times PC([0,\infty),U)$
hei"st ein zul"assiger Steuerungsprozess in der Aufgabe (\ref{UHA1})--(\ref{UHA3}),
falls $\big(x(\cdot),u(\cdot)\big)$ über jedem endlichen Intervall der Dynamik (\ref{UHA2}) zu $x(t_0)=x_0$ gen"ugt,
die Steuerbeschränkungen (\ref{UHA3}) erfüllt und das Zielfunktional (\ref{UHA1}) endlich ist.
Mit $\mathscr{D}^{\,\mathcal{U}}_{\rm adm}$ bezeichnen wir die Menge der zul"assigen Steuerungsprozesse. \\[2mm]
Ein zul"assiger Steuerungsprozess $\big(x_*(\cdot),u_*(\cdot)\big)$ ist eine
starke lokale Minimalstelle\index{Minimum, starkes lokales!elementar@-- unendlicher Zeithorizont}
der Aufgabe (\ref{UHA1})--(\ref{UHA3}),
falls eine Zahl $\varepsilon > 0$ derart existiert, dass die Ungleichung 
$$J\big(x(\cdot),u(\cdot)\big) \geq J\big(x_*(\cdot),u_*(\cdot)\big)$$
f"ur alle $\big(x(\cdot),u(\cdot)\big) \in \mathscr{D}^{\,\mathcal{U}}_{\rm adm}$ mit $\|x(\cdot)-x_*(\cdot)\|_\infty < \varepsilon$ gilt. \\[2mm]
Eine wesentliche Herausforderung bei der Behandlung von Steuerungsproblemen mit unendlichem Zeithorizont
ist der Nachweis einer adjungierten Funktion,
die sowohl die adjungierte Gleichung als auch eine Transversalitätsbedingung im Unendlichen erfüllt.
Die Existenz der Adjungierten wird im Folgenden mit Hilfe der Bedingung
\begin{equation} \label{PMPBedingungeinfach}
\int_0^\infty \big\|\varphi_x\big(t,x(t),u(t)\big)\big\| \, dt < \infty
\end{equation}
nach Lemma \ref{LemmaDGL4} und \ref{LemmaDGL6} gesichert.
Allerdings ist diese Bedingung nicht unkritisch,
da sie zum Beispiel lineare Dynamiken mit konstanten Koeffizienten ausschließt.
Ungeachtet der umfassenden Einschränkungen können und wollen wir auf die Bedingung (\ref{PMPBedingungeinfach}) nicht verzichten
und bezeichnen mit $\mathscr{D}^{\,\mathcal{U}}_{\lim}$ die Menge aller Paare $\big(x(\cdot),u(\cdot)\big)$,
welche der Bedingung (\ref{PMPBedingungeinfach}) genügen. \\[2mm] 
Im Folgenden bezeichnet $H^{\mathcal{U}}: \R \times \R^n \times \R^m \times \R^n \times \R \to \R$ die Pontrjagin-Funktion
$$H^{\mathcal{U}}(t,x,u,p,\lambda_0) = \langle p, \varphi(t,x,u) \rangle - \lambda_0 e^{-\varrho t} f(t,x,u).$$

\begin{theorem}[Pontrjaginsches Maximumprinzip] \label{SatzPMPUHA}
\index{Pontrjaginsches Maximumprinzip!elementar@-- elementare Aufgabe}
\index{Pontrjaginsches Maximumprinzip!unendlich@-- unendlicher Zeithorizont} 
Sei $\big(x_*(\cdot),u_*(\cdot)\big) \in \mathscr{D}^{\,\mathcal{U}}_{\rm adm} \cap \mathscr{D}^{\,\mathcal{U}}_{\rm Lip} \cap \mathscr{D}^{\,\mathcal{U}}_{\lim}$.
Ist $\big(x_*(\cdot),u_*(\cdot)\big)$ ein starkes lokales Minimum der Aufgabe (\ref{UHA1})--(\ref{UHA3}),
dann existiert eine Vektorfunktion $p(\cdot) \in PC_1([0,\infty),\R^n)$ derart, dass
\begin{enumerate}
\item[(a)] die adjungierte Gleichung
           \index{adjungierte Gleichung!elementar@-- elementare Aufgabe}
           \index{adjungierte Gleichung!unendlich@-- unendlicher Zeithorizont}
           \begin{equation}\label{PMPUHA1} 
           \dot{p}(t) = -H_x^{\mathcal{U}}\big(t,x_*(t),u_*(t),p(t),1\big),
           \end{equation} 
\item[(b)] für $t \to \infty$ die Transversalitätsbedingung
           \index{Transversalitätsbedingungen!elementar@-- elementare Aufgabe}
           \index{Transversalitätsbedingungen!unendlich@-- unendlicher Zeithorizont}
           \begin{equation}\label{PMPUHA2} 
           \lim_{t \to \infty} p(t)=0
           \end{equation} 
\item[(c)] und in fast allen Punkten $t \in [0,\infty)$ die Maximumbedingung 
           \index{Maximumbedingung!elementar@-- elementare Aufgabe}
           \index{Maximumbedingung!unendlich@-- unendlicher Zeithorizont}
           \begin{equation}\label{PMPUHA3} 
           H^{\mathcal{U}}\big(t,x_*(t),u_*(t),p(t),1\big) = \max_{u \in U} H^{\mathcal{U}}\big(t,x_*(t),u,p(t),1\big).
           \end{equation}
\end{enumerate}
erfüllt sind.
\end{theorem}

{\bf Beweis} Da der Steuerungsprozess $\big(x_*(\cdot),u_*(\cdot)\big)$ der Menge
$\mathscr{D}^{\,\mathcal{U}}_{\rm adm} \cap \mathscr{D}^{\,\mathcal{U}}_{\rm Lip} \cap \mathscr{D}^{\,\mathcal{U}}_{\lim}$ angehört,
ist er beschränkt.
Ferner sind $t \to e^{-\varrho t} f_x\big(t,x_*(t),u_*(t)\big)$ und
$t \to \varphi_x\big(t,x_*(t),u_*(t)\big)$ über $[0,\infty)$ integrierbar, beschränkt und stückweise stetig.
Damit sind die Voraussetzungen von Lemma \ref{LemmaDGL4} und von Lemma \ref{LemmaDGL6} erfüllt
und es existiert eine eindeutige stetige L"osung $p(\cdot)$ der adjungierten Gleichung (\ref{PMPUHA1})
zur Transversalitätsbedingung  (\ref{PMPUHA2}).
Wegen der stückweisen Stetigkeit der einfließenden Abbildungen ist $\dot{p}(\cdot)$ ebenfalls stückweise stetig
und die Adjungierte $p(\cdot)$ gehört dem Raum $PC_1([0,\infty),\R^n)$ an.  \\
Außerdem ergibt sich mit der gleichen Argumentation wie im Abschnitt \ref{AbschnittPMPBeweiseinfach},
dass man ohne Einschränkung $\lambda_0=1$ annehmen darf. \\[2mm]
Da $t \to \varphi_x\big(t,x_*(t),u_*(t)\big)$ integrierbar ist, besitzt die Integralgleichung
\begin{equation} \label{BeweisPMPUH1}
y(t)=y(\tau) + \int_{\tau}^{t} \varphi_x\big(s,x_*(s),u_*(s)\big)\,y(s) \, ds, \quad t \in [0,\infty),
\end{equation}
nach Lemma \ref{LemmaDGL4} für jedes $\tau \in [0,\infty)$ und jedes $y(\tau) \in \R^n$ eine eindeutige Lösung. \\[2mm]
Wir wenden wieder die einfache Nadelvariation\index{Nadelvariation, einfache}
$u(t;v,\tau,\lambda) = u_{\lambda}(t)$ des Beweises in Abschnitt \ref{AbschnittPMPBeweiseinfach} an.
Für $T > \tau$ folgt mit den gleichen Argumenten wie ebenda,
dass für $t \in [\tau,T]$ der Grenzwert
$$y(t;T)=\lim_{\lambda \to 0^+}\frac{x_{\lambda}(t) - x_*(t)}{\lambda}$$
existiert, die Funktion $y(\cdot;T)$ der Integralgleichung
\begin{eqnarray*}
y(t;T) &=& y(\tau) + \int_{\tau}^{t} \varphi_x\big(s,x_*(s),u_*(s)\big)\,y(s;T) \, ds, \\
y(\tau) &=& \varphi\big(\tau,x_*(\tau),v\big) - \varphi\big(\tau,x_*(\tau),u_*(\tau)\big)
\end{eqnarray*}
gen"ugt.
Da die Integralgleichung (\ref{BeweisPMPUH1}) eine eindeutige Lösung $y(\cdot)$ über $[0,\infty)$ besitzt,
gilt $y(t;T)= y(t)$ für alle $[\tau,T]$, und dies für beliebiges $T > \tau$.
Ferner ist die Beziehung
$$\langle p(\tau) , y(\tau) \rangle-\langle p(T) , y(T) \rangle = - \int_{\tau}^T e^{-\varrho t} \big\langle f_x\big(t,x_*(t),u_*(t)\big) , y(t) \big\rangle \, dt$$
erf"ullt.
Da $y(\cdot)$ nach Lemma \ref{LemmaDGL4} einen Grenzwert im Unendlichen besitzt, ergibt sich
\begin{equation} \label{BeweisPMPUH2}
\langle p(\tau) , y(\tau) \rangle = - \int_{\tau}^\infty e^{-\varrho t} \big\langle f_x\big(t,x_*(t),u_*(t)\big) , y(t) \big\rangle \, dt.
\end{equation}
Es sei $x_\lambda(\cdot)$ die Variation von $x_*(\cdot)$ zur Steuerung $u_\lambda(\cdot)$,
die auf $[\tau,\infty)$ gleichmäßig konvergent gegen $x_*(\cdot)$ sei.
Es zeigt sich abschließend wie in Abschnitt \ref{AbschnittPMPBeweiseinfach}
\begin{eqnarray*}
0 &\leq& \lim_{\lambda \to 0^+} \frac{J\big(x_\lambda(\cdot),u_\lambda(\cdot)\big)- J\big(x_*(\cdot),u_*(\cdot)\big)}{\lambda} \\
  &=   & f\big(\tau,x_*(\tau),v\big) - f\big(\tau,x_*(\tau),u_*(\tau)\big)
         + \int_{\tau}^\infty e^{-\varrho t} \big\langle f_x\big(t,x_*(t),u_*(t)\big),y(t) \big\rangle \, dt,
\end{eqnarray*}
welche zusammen mit (\ref{BeweisPMPUH2}) der Gültigkeit der Maximumbedingung (\ref{PMPUHA3}) entspricht. \hfill $\blacksquare$

\begin{bemerkung} {\rm
Im Beweis sind die Wohldefiniertheit von $x_\lambda(\cdot)$ auf $[0,\infty)$ und die gleichmäßige Konvergenz von $x_\lambda(\cdot)$ gegen $x_*(\cdot)$
nicht gesichert! \hfill $\square$}
\end{bemerkung}

\begin{beispiel}
{\rm Wir betrachten eine Variante des Differentialspiels \ref{BeispielDockner} nach Dockner et\,al. \cite{Dockner},
welches bereits im Abschnitt \ref{BeispielDockner} über einem endlichen Zeitraum diskutiert wurde:
\index{Fischereimodell}
\begin{eqnarray*}
&& \tilde{J}_i\big(x(\cdot),u_1(\cdot),u_2(\cdot)\big) =\int_0^\infty e^{-\varrho t}\big(p x(t)-c_i\big)u_i(t) \, dt \to \sup, \\
&&  \dot{x}(t)=e^{-\delta t} \big[x(t) \big(\alpha-r\ln x(t) \big) -u_1(t)x(t)-u_2(t)x(t)\big], \quad x(0)=x_0>0, \\
&& u_i > 0, \quad \alpha,c_i,p,r,\varrho >0, \quad \alpha> \frac{1}{c_1+c_2}, \quad i=1,2.
\end{eqnarray*}
Im Vergleich zu \cite{Dockner} haben wir den Faktor $e^{-\delta t}$ in der Dynamik hinzugefügt,
um die Stabilität des dynamischen Systems gemäß der Bedingungen (\ref{PMPBedingungeinfach}) zu sichern. \\[2mm]
Wir wenden wieder die Transformation $z=\ln x$ an und erhalten für ein Nash-Gleichgewicht $\big(z^*(\cdot),u^*_1(\cdot),u^*_2(\cdot)\big)$
f"ur $i,j \in \{1,2\}$, $i \not= j$, die gekoppelten Steuerungsprobleme:
\begin{eqnarray*}
&& J_i\big(z(\cdot),u_i(\cdot),u^*_j(\cdot)\big) =
      \int_0^\infty e^{-\varrho t}\bigg(\frac{1}{u_i(t)+u^*_j(t)}-c_i\bigg)u_i(t) \, dt \to \sup, \\
&& \dot{z}(t)=e^{-\delta t} \big[ -r z(t)+\alpha-u_i(t)-u^*_j(t)\big], \quad z(0)=\ln x_0>0, \quad u_i>0.
\end{eqnarray*}

Die zul"assigen Steuerungen geh"oren dem Raum $PC([0,\infty),U)$ an.
Deshalb sind im Ansatz des Nash-Gleichgewichtes die Abbildungen
$$f_i(t,z,u_i)=e^{-\varrho t}\bigg(\frac{1}{u_i+u^*_j(t)}-c_i\bigg)u_i, \quad \varphi_i(t,z,u_i)=e^{-\delta t}[-r z+\alpha-u_i-u^*_j(t)]$$
nicht stetig bez"uglich der Variable $t$.
An dieser Stelle verweisen auf die Bemerkungen im Rahmen der Differentialspiele im Abschnitt \ref{KapitelDifferentialspiele},
dass unter diesen Rahmenbedingungen Theorem \ref{SatzPMPUHA} seine G"ultigkeit beh"alt. \\[2mm]
Die Pontrjagin-Funktionen lauten
\begin{eqnarray*}
H^{\mathcal{U}}_1(t,z,u_1,p_1) &=& p_1 e^{-\delta t} \big[-r z+\alpha-u_1-u^*_2(t)\big]+ e^{-\varrho t}\bigg(\frac{1}{u_1+u^*_2(t)}-c_1\bigg)u_1, \\
H^{\mathcal{U}}_2(t,z,u_2,p_2) &=& p_2 e^{-\delta t} \big[-r z+\alpha-u^*_1(t)-u_2\big]+ e^{-\varrho t}\bigg(\frac{1}{u^*_1(t)+u_2}-c_2\bigg)u_2.
\end{eqnarray*}  
Mit Lemma \ref{LemmaDGL6} ergibt sich die eindeutige Lösung der adjungierten Gleichung (\ref{PMPUHA1}):
$$\dot{p}_i(t)=r e^{-\delta t} p_i(t), \quad \lim_{t \to \infty} p_i(t)=0,\qquad\Rightarrow\qquad p_i(t) \equiv 0.$$
Die Maximumbedingung (\ref{PMPUHA2}) liefert das Gleichungssystem
$$\frac{\partial}{\partial u_i} H^{\mathcal{U}}_i\big(t,z_*(t),u_i^*(t),p_i(t)\big)=0 \quad\Rightarrow\quad
  \frac{1}{u_1^*(t)+u_2^*(t)}- \frac{u_i^*(t)}{\big(u_1^*(t)+u_2^*(t)\big)^2}=c_i, \; i=1,2.$$
Aus der Summe beider Gleichungen ergibt sich
$$\frac{2}{u_1^*(t)+u_2^*(t)}- \frac{u_1^*(t)+u_2^*(t)}{\big(u_1^*(t)+u_2^*(t)\big)^2}=\frac{1}{u_1^*(t)+u_2^*(t)}=c_1+c_2$$
und wir erhalten die optimalen Steuerungen
$$u_1^*(t)\equiv \frac{c_2}{(c_1+c_2)^2}, \qquad u_2^*(t)\equiv \frac{c_1}{(c_1+c_2)^2}.$$
Die optimale Trajektorie besitzt die Gestalt
$$z_*(t)=\bigg(z_0-\frac{c_0}{r}\bigg)\exp\bigg\{\frac{r}{\delta}( e^{-\delta t}-1)\bigg\} + \frac{c_0}{r}, \qquad c_0=  \alpha-\frac{1}{c_1+c_2}.$$
Für diese gilt im Unendlichen
$$z_\infty=\lim_{t \to \infty}z_*(t)=\bigg(z_0-\frac{c_0}{r}\bigg)e^{-\frac{r}{\delta}} + \frac{c_0}{r}
                           = z_0 e^{-\frac{r}{\delta}}+\bigg(\alpha-\frac{1}{c_1+c_2}\bigg)(1-e^{-\frac{r}{\delta}})>0.$$
Die Funktion $z_*(\cdot)$ ist streng monoton und nimmt nur Werte zwischen $z_0$ und $z_\infty$ an.
Somit ist $x_*(t)=\exp\big(z_*(t)\big)$ "uber $[0,\infty)$ wohldefiniert, besitzt eine untere positive Schranke und
$\big(x_*(\cdot),u^*_1(\cdot),u^*_2(\cdot)\big)$ liefert einen Kandidaten f"ur das urspr"ungliche Differentialspiel. \hfill $\square$}
\end{beispiel}

%% file: 4-2-Aufgabenstellung.tex
\subsection{Die Aufgabenstellung}
Wir betrachten starke lokale Minimalstellen in der folgenden Aufgabe:
\begin{eqnarray}
&& \label{PAUH1} J\big(x(\cdot),u(\cdot)\big) = \int_0^\infty \omega(t)f\big(t,x(t),u(t)\big) \, dt \to \inf, \\
&& \label{PAUH2} \dot{x}(t) = \varphi\big(t,x(t),u(t)\big), \\
&& \label{PAUH3} h_0\big(x(0)\big)=0, \qquad \lim_{t \to \infty} h_1\big(t,x(t)\big)=0, \\
&& \label{PAUH4} u(t) \in U \subseteq \R^m, \quad U \not= \emptyset, \\
&& \label{PAUH5} g_j\big(t,x(t)\big) \leq 0 \quad \mbox{f"ur alle } t \in \R_+, \quad j=1,...,l.
\end{eqnarray}
Dabei ist $\omega(\cdot) \in L_1(\R_+,\R_+)$ und es gelten f"ur die eingehenden Abbildungen
$$f:\R \times \R^n \times \R^m \to \R, \qquad \varphi:\R \times \R^n \times \R^m \to \R^n, \qquad g_j: \R \times \R^n \to \R,$$
sowie f"ur die Randbedingungen $h_0:\R^n \to \R^{s_0}$ und $h_1:\R \times \R^n \to \R^{s_1}$. \\[2mm]
Wir nennen die Trajektorie $x(\cdot)$ eine L"osung des dynamischen Systems (\ref{PAUH2}),
falls $x(\cdot)$ auf $\R_+$ definiert ist und auf jedem endlichen Intervall die Dynamik mit Steuerung $u(\cdot)$
im Sinn von Carath\'eodory l"ost. \\[2mm]
Mit $\mathscr{B}^{\,\mathcal{U}}_{\rm Lip}$ bezeichnen wir die Menge aller Paare $\big(x(\cdot),u(\cdot)\big)$,
für die es ein $\gamma>0$ derart gibt,
dass die Abbildungen $f(t,x,u)$, $\varphi(t,x,u)$, $h_0(x_0)$, $h_1(t,x)$ und $g_j(t,x)$ auf der Menge aller
$(t,x,x_0,u) \in \R \times \R^n \times \R^n \times \R^m$ mit
$0 \leq t < \infty$, $\|x-x(t)\| < \gamma$, $\|x_0-x(0)\| < \gamma$ und $u \in U_1$
beschränkt und gleichmäßig stetig in allen Variablen, sowie gleichmäßig stetig differenzierbar bezüglich $x,x_0$
mit beschränkten Ableitungen $f_x(t,x,u)$, $\varphi_x(t,x,u)$, $h_0'(x_0)$, $h_{1,x}(t,x)$ und $g_{j,x}(t,x)$ sind.
Außerdem seien diese Eigenschaften für $g_j(t,x)$ und $h_1(t,x)$ in $t=\infty$ fortsetzbar,
da sich dieser Zeitpunkt auf die Aufgabe maßgeblich auswirkt. \\[2mm]
Der Steuerungsprozess $\big(x(\cdot),u(\cdot)\big) \in W^1_\infty(\R_+,\R^n) \times L_\infty(\R_+,U)$
hei"st zul"assig in der Aufgabe (\ref{PAUH1})--(\ref{PAUH5}),
falls $\big(x(\cdot),u(\cdot)\big)$ dem System (\ref{PAUH2}) gen"ugt,
die Restriktionen (\ref{PAUH3}) und (\ref{PAUH4}) erf"ullt,
sowie das Lebesgue-Integral in (\ref{PAUH1}) endlich ist.
Die Menge $\mathscr{B}^{\,\mathcal{U}}_{\rm adm}$ bezeichnet in der Aufgabe (\ref{PAUH1})--(\ref{PAUH5}) die Menge der zul"assigen Steuerungsprozesse. \\[2mm]
Der zul"assige Steuerungprozess $\big(x_*(\cdot),u_*(\cdot)\big)$ ist eine starke lokale
Minimalstelle\index{Minimum, starkes lokales!unendlich@-- unendlicher Zeithorizont}
in der Aufgabe (\ref{PAUH1})--(\ref{PAUH5}),
falls eine Zahl $\varepsilon > 0$ derart existiert, dass die Ungleichung
$$J\big(x(\cdot),u(\cdot)\big) \geq J\big(x_*(\cdot),u_*(\cdot)\big)$$
f"ur alle $\big(x(\cdot),u(\cdot)\big) \in \mathscr{B}^{\,\mathcal{U}}_{\rm adm}$ mit 
$\|x(\cdot)-x_*(\cdot)\|_\infty < \varepsilon$ gilt. \\[2mm]
Bez"uglich des speziellen Rahmens, in die wir die Aufgabe (\ref{PAUH1})--(\ref{PAUH5}) einbetten,
sind zus"atzliche Einschr"ankungen erforderlich:
Zu $\big(x(\cdot),u(\cdot)\big) \in W^1_\infty(\R_+,\R^n) \times L_\infty(\R_+,U)$ seien die folgenden Integralterme endlich
\begin{equation} \label{PMPBedingung}
\int_0^\infty \big\|\varphi\big(t,x(t),u(t)\big)\big\| \, dt < \infty, \qquad 
\int_0^\infty \big\|\varphi_x\big(t,x(t),u(t)\big)\big\| \, dt < \infty.
\end{equation}
Weiterhin nehmen wir an, dass es zu jedem $\delta>0$ ein $T>0$ existiert mit
\begin{eqnarray}
&& \int_T^\infty \big\| \varphi\big(t,\xi(t),u(t)\big)-\varphi\big(t,\xi'(t),u(t)\big)
                         - \varphi_x\big(t,x(t),u(t)\big)\big(\xi(t)-\xi'(t)\big) \big\| \, dt \nonumber \\
&& \label{PMPBedingung2} \hspace*{20mm} \leq \delta \|\xi(\cdot)-\xi'(\cdot)\|_\infty
\end{eqnarray}
f"ur alle $\xi(\cdot), \xi'(\cdot) \in W^1_\infty(\R_+,\R^n)$ mit
$\|\xi(\cdot)-x(\cdot)\|_\infty < \gamma$, $\|\xi'(\cdot)-x(\cdot)\|_\infty < \gamma$.
Die Menge $\mathscr{B}^{\,\mathcal{U}}_{\lim}$ bezeichnet die Menge aller $\big(x(\cdot),u(\cdot)\big) \in W^1_\infty(\R_+,\R^n) \times L_\infty(\R_+,U)$,
die die Eigenschaften (\ref{PMPBedingung}) und (\ref{PMPBedingung2}) besitzen.

\begin{bemerkung}
{\rm Aus der ersten Bedingung in (\ref{PMPBedingung}) ergibt sich,
dass die Trajektorie $x(\cdot)$ im Unendlichen einen Grenzwert besitzt.
Damit lässt sich die Dynamik als eine Abbildung in den Raum $C_{\lim}(\R_+,\R^n)$ auffassen.
Die zweite Bedingung in (\ref{PMPBedingung}) und die Bedingung (\ref{PMPBedingung2}) sind in dem vorliegenden Rahmen entscheidend für die
stetige Differenzierbarkeit der Dynamik im Sinne einer Abbildung in Banachräumen.
Andererseits sind die strengen Einschr"ankungen (\ref{PMPBedingung}) und (\ref{PMPBedingung2}) zum Beispiel nicht
f"ur lineare Systeme mit konstanten Koeffizienten erf"ullt.
An dieser Stelle sind hinreichende Bedingungen nach Arrow besser geeignet,
da die Annahmen (\ref{PMPBedingung}) und (\ref{PMPBedingung2}) fallen gelassen werden k"onnen. \hfill $\square$}
\end{bemerkung}

Abschlie"send bezeichnen wir mit $H^{\mathcal{U}}: \R \times \R^n \times \R^m \times \R^n \times \R \to \R$ die Pontrjagin-Funktion
$$H^{\mathcal{U}}(t,x,u,p,\lambda_0) = \langle p, \varphi(t,x,u) \rangle-\lambda_0 \omega(t)f(t,x,u)$$
der Aufgabe (\ref{PAUH1})--(\ref{PAUH5}).

%% file: 4-3-PontrjaginAufgabe.tex
\subsection{Das Pontrjaginsche Maximumprinzip} \label{AbschnittPMPUH}
\subsubsection{Notwendige Optimalit\"atsbedingungen}
\begin{theorem}[Pontrjaginsches Maximumprinzip] \label{SatzPAUH}
\index{Pontrjaginsches Maximumprinzip!unendlich@-- unendlicher Zeithorizont} 
Sei $\big(x_*(\cdot),u_*(\cdot)\big) \in \mathscr{B}^{\,\mathcal{U}}_{\rm adm} \cap \mathscr{B}^{\,\mathcal{U}}_{\rm Lip} \cap \mathscr{B}^{\,\mathcal{U}}_{\lim}$.
Ist $\big(x_*(\cdot),u_*(\cdot)\big)$ ein starkes lokales Minimum der Aufgabe (\ref{PAUH1})--(\ref{PAUH4}),
dann existieren nicht gleichzeitig verschwindende Multiplikatoren $\lambda_0 \geq 0$,
eine absolutstetige Funktion $p(\cdot):\R_+ \to \R^n$ und $l_i \in \R^{s_i}$, $i=0,1$, derart, dass
\begin{enumerate}
\item[(a)] die adjungierte Gleichung
           \index{adjungierte Gleichung!unendlich@-- unendlicher Zeithorizont}
           \begin{equation}\label{SatzPAUH1}
           \dot{p}(t) = -H_x^{\mathcal{U}}\big(t,x_*(t),u_*(t),p(t),\lambda_0\big),
           \end{equation}
\item[(b)] die Transversalit"atsbedingungen
           \index{Transversalitätsbedingungen!unendlich@-- unendlicher Zeithorizont}
           \begin{equation}\label{SatzPAUH2}
           p(0) = {h_0'}^T\big(x_*(0)\big)l_0, \qquad \lim_{t \to \infty} p(t)= - \lim_{t \to \infty}h_{1,x}^T\big(t,x_*(t)\big) l_1
           \end{equation}
\item[(c)] und in fast allen Punkten $t \in \R_+$ die Maximumbedingung
           \index{Maximumbedingung!unendlich@-- unendlicher Zeithorizont}
           \begin{equation}\label{SatzPAUH3}
           H^{\mathcal{U}}\big(t,x_*(t),u_*(t),p(t),\lambda_0\big) = \max_{u \in U} H^{\mathcal{U}}\big(t,x_*(t),u,p(t),\lambda_0\big)
           \end{equation}
\end{enumerate}
erfüllt sind.
Anhand (\ref{SatzPAUH1}) kann für die verallgemeinerte Ableitung $\dot{p}(\cdot) \in L_1(\R_+,\R^n)$
im Fall einer unbeschränkten Verteilungsfunktion $\omega(\cdot) \in L_1(\R_+,\R_+)$ gelten.
\end{theorem}

\newpage
\begin{beispiel} \label{BeispielRWUnendlich2}
{\rm Wir betrachten die Aufgabe
\begin{eqnarray*}
&& J\big(x(\cdot),u(\cdot)\big) = \int_0^\infty e^{-\varrho t} \big(1-u(t)\big) x(t) \, dt \to \sup,\\
&& \dot{x}(t)=u(t)x(t), \quad x(0)=1, \quad \lim_{t \to \infty} x(t)=x_1>1, \quad u(t) \in [0,1], \quad \varrho \in (0,1).
\end{eqnarray*}
Wir stellen zun"achst die Pontrjagin-Funktion auf:
$$H^{\mathcal{U}}(t,x,z,u,p,q,\lambda_0) = pux+\lambda_0 e^{-\varrho t}(1-u)x.$$
Mit Hilfe der Bedingungen (\ref{SatzPAUH1})--(\ref{SatzPAUH3}) k"onnen wir den Fall $\lambda_0=0$ ausschlie"sen.
Weiterhin ergeben sich der Steuerungsprozess
$$x_*(t) = \left\{ \begin{array}{ll} e^t, & t \in [0,\tau), \\ x_1, & t \in [\tau, \infty), \end{array}\right. \quad
  u_*(t)= \left\{ \begin{array}{ll} 1, & t \in [0,\tau), \\ 0, & t \in [\tau, \infty), \end{array}\right. \quad \tau=\ln x_1$$
und die Adjungierte
$$p(t) = \left\{ \begin{array}{ll}
          e^{(1-\varrho)\tau} e^{-t}, & t \in [0,\tau), \\[1mm]
          \displaystyle\frac{\varrho-1}{\varrho}e^{-\varrho \tau} + \frac{1}{\varrho} e^{-\varrho t}, & t \in [\tau, \infty). \end{array}\right.$$
F"ur die Adjungierte gilt dabei $\displaystyle\lim_{t \to \infty} p(t)= \frac{\varrho-1}{\varrho}e^{-\varrho \tau} \not=0$ im Unendlichen.
Damit konnten wir aus den notwendigen Bedingungen (\ref{SatzPAUH1})--(\ref{SatzPAUH3}) des Maximumprinzips einen eindeutigen
Kandidaten bestimmen. \hfill $\square$}
\end{beispiel}

\begin{beispiel} \label{BeipielRWUnendlich}
{\rm Wir betrachten die Aufgabe
\begin{eqnarray*}
&& J\big(x(\cdot),u(\cdot)\big)=\int_0^\infty e^{-\varrho t}\big(1-u(t)\big)x(t) \, dt \to \sup, \\
&& \dot{x}(t)=u(t)x(t), \qquad x(0)=1, \qquad u \in [0,1], \qquad \varrho \in (0,1)
\end{eqnarray*}
mit der Budgetbeschr"ankung
$$\int_0^\infty e^{-\varrho t}x(t) \, dt = Z, \qquad Z> \frac{1}{\varrho}.$$
Bez"uglich der Budgetbeschr"ankung f"uhren wir die folgende Zustandsgleichung mit Randwert im Unendlichen ein:
$$\dot{z}(t)=e^{-\varrho t}x(t), \qquad z(0)=0, \qquad \lim_{t \to \infty} z(t) = Z.$$
Offensichtlich ist $\dot{z}(t) >0$ auf $\R_+$ und damit $z(t)$ streng monoton wachsend.
Demzufolge kann die Beschr"ankung $z(t) \leq Z$ erst im Unendlichen aktiv werden und greift nur durch das Verhalten in
$t=\infty$ in die gestellte Aufgabe ein. \\
Da stets $\dot{x}(t)\geq 0$ ist,
muss f"ur jede zul"assige Trajektorie $e^{-\varrho t}x(t) \to 0$ f"ur $t \to \infty$ gelten,
denn nur dann ist $z(t) \leq Z$ erf"ullt.
Damit erhalten wir f"ur zul"assige Steuerungsprozesse im Zielfunktional zun"achst
\begin{eqnarray*}
    J\big(x(\cdot),z(\cdot),u(\cdot)\big)
&=& \int_0^\infty e^{-\varrho t}\big(1-u(t)\big)x(t) \, dt
     = \int_0^\infty \dot{z}(t) \, dt - \int_0^\infty e^{-\varrho t} \dot{x}(t) \, dt \\
&=& \int_0^\infty \dot{z}(t) \, dt + 1 - \varrho \int_0^\infty e^{-\varrho t} x(t) \, dt \leq 1 + (1-\varrho) Z.
\end{eqnarray*}
Es ergibt sich also f"ur das Zielfunktional die obere Schranke $1 + (1-\varrho) Z$.
D.\,h., dass jeder Steuerungsprozess $\big(x(\cdot),z(\cdot),u(\cdot)\big)$,
f"ur den die Zustandsbeschr"ankung $z(t) \leq Z$ im Unendlichen aktiv wird,
global optimal ist und $J\big(x(\cdot),z(\cdot),u(\cdot)\big) = 1 + (1-\varrho) Z$ gilt. \\
Die Voraussetzungen des Pontrjaginschen Maximumprinzips an einen zul"assigen Steuerungsprozess sind in dem vorliegenden
Beispiel genau dann erf"ullt,
wenn die Steuerung $u(\cdot)$ dem Raum $L_1(\R_+,[0,1])$ angeh"ort.
Denn in diesem Fall gelten  f"ur die korrespondierende Trajektorie $x(\cdot)$ die Bedingungen
(\ref{PMPBedingung}), (\ref{PMPBedingung2}).
Wir diskutieren zwei F"alle:
\begin{enumerate}
\item[(A)] Es liefert der Steuerungsprozess $\big(x_*(\cdot),z_*(\cdot),u_*(\cdot)\big)$ mit
           $$x_*(t) = e^{\alpha t},\qquad z_*(t)= \frac{1}{\alpha - \varrho} (e^{(\alpha-\varrho)t}-1),\qquad
             u_*(t)=\alpha,\qquad \alpha=\varrho-\frac{1}{Z} \in (0, \varrho)$$ ein globales Maximum.
           Da die vorgeschlagene Steuerung $u_*(\cdot)$ "uber $\R_+$ nicht integrierbar ist,
           sind f"ur $x_*(\cdot)$ die Bedingungen (\ref{PMPBedingung}) nicht erf"ullt.
           Das Maximumprinzip ist in diesem Fall nicht anwendbar.
\item[(B)] Es stellt der Steuerungsprozess $\big(y_*(\cdot),\zeta_*(\cdot),w_*(\cdot)\big)$ mit
           $$y_*(t) = \left\{ \begin{array}{ll} e^t,& t \in [0,\tau), \\ e^\tau, & t \in [\tau,\infty), \end{array} \right.
                      \quad
             w_*(t) = \left\{ \begin{array}{ll} 1,& t \in [0,\tau), \\ 0, & t \in [\tau,\infty) \end{array} \right.$$
           mit dem Umschaltzeitpunkt $\tau >0$, der der Bedingung
           $$e^{(1-\varrho)\tau}\bigg(\frac{1}{\varrho}+\frac{1}{1-\varrho}\bigg) = Z+\frac{1}{1-\varrho},
             \qquad Z> \frac{1}{\varrho},$$
           gen"ugt,
           ebenfalls ein globales Maximum dar.
           Die zugeh"orige Trajektorie $\zeta_*(\cdot)$ lautet
           $$\zeta_*(t) = \left\{ \begin{array}{ll}
             \displaystyle \frac{1}{1-\varrho}\big(e^{(1-\varrho)t} - 1 \big) ,& t \in [0,\tau), \\
       \displaystyle \zeta(\tau) + \frac{1}{\varrho}\big(e^{(1-\varrho)\tau} - e^{\tau-\varrho t} \big), & t \in [\tau,\infty).
             \end{array} \right.$$
           Da die Steuerung $w_*(\cdot)$ dem Raum $L_1(\R_+,[0,1])$ angeh"ort,
           gelten s"amtliche Voraussetzungen von Theorem \ref{SatzPAUH}.
           Mit den Multiplikatoren
           $$\lambda_0=1, \qquad p(t)= e^{-\varrho t}, \qquad q(t)= \varrho-1$$
           sind dann die notwendigen Bedingungen (\ref{SatzPAUH1})--(\ref{SatzPAUH3}) erf"ullt. \hfill $\square$
\end{enumerate}}
\end{beispiel}

%% file: 4-31-Beweis.tex
\subsubsection{Der Nachweis der notwendigen Optimalit\"atsbedingungen} \label{AbschnittBeweisPMPUH}
Wir betrachten f"ur $\big(x(\cdot),u(\cdot)\big) \in C_{\lim}(\R_+,\R^n) \times L_\infty(\R_+,\R^m)$ die Abbildungen
\begin{eqnarray*}
J\big(x(\cdot),u(\cdot)\big) &=& \int_0^\infty \omega(t)f\big(t,x(t),u(t)\big) \, dt, \\
F\big(x(\cdot),u(\cdot)\big)(t) &=& x(t) -x(t_0) -\int_0^t \varphi\big(s,x(s),u(s)\big) \, ds, \quad t \in \R_+,\\
H_0\big(x(\cdot)\big) &=& h_0\big(x(0)\big), \qquad H_1\big(x(\cdot)\big) = \lim_{t \to \infty} h_1\big(t,x(t)\big).
\end{eqnarray*}
Die Abbildungen fassen wir in folgenden Funktionenr"aumen auf:
\begin{eqnarray*}
J &:& C_{\lim}(\R_+,\R^n) \times L_\infty(\R_+,\R^m) \to \R, \\
F &:& C_{\lim}(\R_+,\R^n) \times L_\infty(\R_+,\R^m) \to C_{\lim}(\R_+,\R^n), \\
H_i &:& C_{\lim}(\R_+,\R^n) \to \R^{s_i}, \quad i=0,1.
\end{eqnarray*}

Wir setzen $\mathscr{F}=(F,H_0,H_1)$, sowie die Menge $\mathscr{U}$ gem"a"s
\begin{eqnarray*}
\mathscr{U}=\big\{ u(\cdot)  \in L_\infty(\R_+,U) \!\!\!&\big|&\!\!\! u(t)=u_*(t) + \chi_M(t)\big(w(t)-u_*(t)\big), \; w(\cdot) \in L_\infty(\R_+,U), \\
                                                  \!\!\!&     &\!\!\! M \subset \R_+ \mbox{ me"sbar und beschr"ankt} \big\}
\end{eqnarray*}
und pr"ufen f"ur die Extremalaufgabe
\begin{equation} \label{ExtremalaufgabePMPUH}
J\big(x(\cdot),u(\cdot)\big) \to \inf, \qquad \mathscr{F}\big(x(\cdot),u(\cdot)\big)=0, \qquad u(\cdot) \in \mathscr{U}
\end{equation}
im Punkt $\big(x_*(\cdot),u_*(\cdot)\big) \in \mathscr{B}^{\,\mathcal{U}}_{\rm Lip} \cap \mathscr{B}^{\,\mathcal{U}}_{\lim}$
die Voraussetzungen von Theorem \ref{SatzExtremalprinzipStark}:

\begin{enumerate}
\item[(A$_1$)] F"ur jedes $u(\cdot) \in \mathscr{U}$ ist die Abbildung $x(\cdot) \to J\big(x(\cdot),u(\cdot)\big)$
               nach Beispiel \ref{DiffZielfunktionalS} im Punkt $x_*(\cdot)$ Fr\'echet-differenzierbar.
\item[(A$_2$)] Die Abbildung $F$ ist die Summe der Abbildung $x(\cdot) \to x(t)$ und der Abbildung
               $$\big(x(\cdot),u(\cdot)\big) \to -\int_0^t \varphi\big(s,x(s),u(s)\big) \, ds.$$
               Im Beispiel \ref{DiffDynamikS} ist die Fr\'echet-Differenzierbarkeit der zweiten Abbildung im Punkt $x_*(\cdot)$
               f"ur jedes $u(\cdot) \in \mathscr{U}$ nachgewiesen.
               Da $x_*(\cdot)$ dem Raum $C_{\lim}(\R_+,\R^n)$ angeh"ort und $h_1(t,x)$ in $t=\infty$ stetig und stetig differenzierbar bez"uglich $x$ ist,
               sind die Abbildungen $H_i$ stetig differenzierbar.
\item[(B)] Nach Lemma \ref{LemmaDGL4} besitzt der Operator $\mathscr{F}_x\big(x_*(\cdot),u_*(\cdot)\big)$ eine endliche Kodimension.
\item[(C)] Der Nachweis dieser Voraussetzungen sind in Lemma \ref{LemmaNVUH1} und Lemma \ref{LemmaNVUH2} "uber mehrfache Nadelvariationen
           "uber dem unendlichen Zeithorizont dargestellt.
\end{enumerate}

Zur Extremalaufgabe (\ref{ExtremalaufgabePMPUH}) definieren wir auf
$$C_{\lim}(\R_+,\R^n) \times L_\infty(\R_+,\R^m) \times \R \times C_{\lim}^*(\R_+,\R^n) \times \R^{s_0} \times \R^{s_1}$$
die Lagrange-Funktion $\mathscr{L}=\mathscr{L}\big(x(\cdot),u(\cdot),\lambda_0,y^*,l_0,l_1\big)$,
$$\mathscr{L}= \lambda_0 J\big(x(\cdot),u(\cdot)\big)+ \big\langle y^*, F\big(x(\cdot),u(\cdot)\big) \big\rangle
                         +l_0^T H_0\big(x(\cdot)\big)+l_1^T H_1\big(x(\cdot)\big).$$
Ist $\big(x_*(\cdot),u_*(\cdot)\big)$ eine starke lokale Minimalstelle der Aufgabe (\ref{ExtremalaufgabePMPUH}),
dann existieren nach Theorem \ref{SatzExtremalprinzipStark}
nicht gleichzeitig verschwindende Lagrangesche Multiplikatoren 
$\lambda_0 \geq 0$, das Funktional $y^* \in C_{\lim}^*(\R_+,\R^n)$ und $l_i \in \R^{s_i}$ derart,
dass folgende Bedingungen gelten:
\begin{enumerate}
\item[(a)] Die Lagrange-Funktion besitzt bez"uglich $x(\cdot)$ in $x_*(\cdot)$ einen station"aren Punkt, d.\,h.
          \begin{equation}\label{SatzPMPLMRUH1}
          \mathscr{L}_x\big(x_*(\cdot),u_*(\cdot),\lambda_0,y^*,l_0,l_1\big)=0;
          \end{equation}         
\item[(b)] Die Lagrange-Funktion erf"ullt bez"uglich $u(\cdot)$ in $u_*(\cdot)$ die Minimumbedingung
           \begin{equation}\label{SatzPMPLMRUH2}
           \mathscr{L}\big(x_*(\cdot),u_*(\cdot),\lambda_0,y^*,l_0,l_1\big)
           = \min_{u(\cdot) \in \mathscr{U}} \mathscr{L}\big(x_*(\cdot),u(\cdot),\lambda_0,y^*,l_0,l_1\big).
           \end{equation}
\end{enumerate}
Aufgrund (\ref{SatzPMPLMRUH1}) ist folgende Variationsgleichung f"ur alle $x(\cdot) \in C_{\lim}(\R_+,\R^n)$ erf"ullt: 
\begin{eqnarray*}
0 &=& \lambda_0 \cdot \int_0^\infty \omega(t) \big\langle f_x\big(t,x_*(t),u_*(t)\big),x(t) \big\rangle\, dt \\
  & & + \int_0^\infty \bigg[ x(t)-x(0) - \int_0^t \varphi_x\big(s,x_*(s),u_*(s)\big) x(s) \,ds \bigg]^T  d\mu(t) \\
  & & + \big\langle l_0, h_0'\big(x_*(0)\big) x(0) \big\rangle + \big\langle l_1, \lim_{t \to \infty} h_{1,x}\big(t,x_*(t)\big) x(t) \big\rangle.
\end{eqnarray*}
Dabei ist $\mu$ nach Folgerung \ref{FolgerungRieszClim} ein signiertes reguläres Borelsches Ma"s über $\overline{\R}_+$,
welche im Anhang \ref{AnhangMassR} in Definition \ref{DefinitionMassR} eingeführt sind.
Wir "andern die Integrationsreihenfolge und setzen $p(t)= \displaystyle \int_t^\infty \, d\mu(s)$.
Damit ist $p(\cdot)$ von beschränkter Variation.
Wir erhalten wegen der eindeutigen Darstellung eines stetigen linearen Funktionals im Raum $C_{\lim}(\R_+,\R^n)$
\begin{eqnarray*}
p(t) &=& -\lim_{t \to \infty} h_{1,x}\big(t,x_*(t)\big)l_1 \\
     & & + \int_{t}^{\infty} \Big( \varphi_x^T\big(s,x_*(s),u_*(s)\big) p(s) - \lambda_0 \omega(s)f_x\big(s,x_*(s),u_*(s)\big) \Big) \, ds, \\
p(t_0) &=& {h_0'}^T\big(x_*(t_0)\big)l_0, \qquad \lim_{t\to \infty} p(t)= -\lim_{t \to \infty} h_{1,x}\big(t,x_*(t)\big)l_1.
\end{eqnarray*}
Die erste Gleichung liefert die absolute Stetigkeit von $p(\cdot)$.
Damit sind die Bedingungen (\ref{SatzPAUH1}) und (\ref{SatzPAUH2}) von Theorem \ref{SatzPAUH} gezeigt. \\[2mm]
Gem"a"s (\ref{SatzPMPLMRUH2}) gilt f"ur alle $u(\cdot) \in \mathscr{U}$ die Ungleichung
$$\int_0^\infty H^{\mathcal{U}}\big(t,x_*(t),u_*(t),p(t),\lambda_0\big) \, dt \geq
  \int_0^\infty H^{\mathcal{U}}\big(t,x_*(t),u(t),p(t),\lambda_0\big) \, dt.$$
Daraus folgt abschlie"send durch Standardtechniken f"ur Lebesguesche Punkte die Maximumbedingung (\ref{SatzPAUH3}).
Der Beweis von Theorem \ref{SatzPAUH} ist abgeschlossen. \hfill $\blacksquare$

%% file: 4-32-Normalform.tex
\subsubsection{Zur normalen Form und zu Transversalit\"atsbedingungen} \label{AbschnittNormalformUH}
F"ur die Aufgabe (\ref{PAUH1})--(\ref{PAUH4}) lassen sich Aussagen "uber die Normalform des Pontrjaginschen Maximumprinzips und zu diversen
Transversalit"atsbedingungen ableiten. \\[2mm]
Wir betrachten zun"achst die Aufgabe (\ref{PAUH1})--(\ref{PAUH4}) mit freiem Endpunkt im Unendlichen.
Dann gen"ugt die Adjungierte $p(\cdot)$ nach Theorem \ref{SatzPAUH} dem Randwertproblem
$$\dot{p}(t) = -\varphi_x^T\big(t,x_*(t),u_*(t)\big) p(t) + \lambda_0 \omega(t)f_x\big(t,x_*(t),u_*(t)\big), \qquad
  \lim_{t \to \infty} p(t)=0.$$
  
\begin{folgerung} \label{FolgerungPAUH1}
In der Aufgabe (\ref{PAUH1})--(\ref{PAUH4}) mit freiem Endpunkt im Unendlichen
ergeben sich aus dem Maximumprinzip die ``nat"urlichen''
Transversalit"atsbedingungen\index{Transversalitätsbedingungen!nat@-- natürliche}:
$$\lim_{t \to \infty} p(t) =0, \qquad \lim_{t \to \infty} \langle p(t),x(t) \rangle = 0 \mbox{ f"ur alle } x(\cdot) \in W^1_\infty(\R_+,\R^n).$$
\end{folgerung}

Nach Voraussetzung des Maximumprinzips ist
$$\int_0^\infty \big\|\varphi_x\big(t,x_*(t),u_*(t)\big)\big\| \, dt < \infty.$$
Unter der Annahme $\lambda_0=0$ w"urde die Adjungierte als L"osung der Gleichung
$$\dot{p}(t) = -\varphi_x^T\big(t,x_*(t),u_*(t)\big) p(t), \qquad \lim_{t \to \infty} p(t)=0,$$
nach Lemma \ref{LemmaDGL6} im Widerspruch zu Theorem \ref{SatzPAUH} identisch verschwinden.

\begin{folgerung} \label{FolgerungPAUH2}
In der Aufgabe (\ref{PAUH1})--(\ref{PAUH4}) mit freiem Endpunkt im Unendlichen gilt $\lambda_0 \not= 0$ und wir k"onnen ohne Einschr"ankung
$\lambda_0=1$ annehmen.
\end{folgerung}

Wegen der Integrierbarkeit der Abbildung $t \to \varphi_x\big(t,x_*(t),u_*(t)\big)$ "uber $\R_+$ k"onnen wir die 
die in $t=0$ normalisierten Fundamentalmatrizen $Y_*(t)$ bzw. $Z_*(t)$ der homogenen Systeme
$$\dot{y}(t)=\varphi_x\big(t,x_*(t),u_*(t)\big) y(t), \qquad \dot{z}(t)=-\varphi_x^T\big(t,x_*(t),u_*(t)\big) z(t)$$
im Rahmen des Raumes $C_{\lim}(\R_+,\R^n)$ betrachten.
Die Adjungierte $p(\cdot)$ erf"ullt nach Theorem \ref{SatzPAUH} die Gleichung (\ref{SatzPAUH1}) und es gilt $\lambda_0=1$
nach Folgerung \ref{FolgerungPAUH2}.
Damit besitzt die Adjungierte $p(\cdot)$ f"ur $t,T \in \R_+$ die Darstellung (vgl. Aseev \& Kryazhimskii \cite{AseKry})
$$p(t)=Z_*(t)\bigg(Z^{-1}_*(T) p(T) + \int_T^t \omega(s) Z_*^{-1}(s)f_x\big(s,x_*(s),u_*(s)\big) ds \bigg).$$
Es bezeichne ferner $y_\xi(\cdot) \in C_{\lim}(\R_+,\R^n)$ die Funktion $y_\xi(t)=Y_*(t) Y^{-1}_*(T) \xi$ mit $\xi \in \R^n$ 
und $\|\xi\|=1$.
Verwenden wir nun $Z^{-1}_*(t)=Y^T_*(t)$, so ergibt sich
$$\langle p(t), y_\xi(t) \rangle = 
  \Big\langle  p(T) + Z_*(T)\int_T^t \omega(s) Z_*^{-1}(s)f_x\big(s,x_*(s),u_*(s)\big) ds , \xi \Big\rangle.$$
Unter Verwendung der ``nat"urlichen'' Transversalit"atsbedingungen in Folgerung \ref{FolgerungPAUH1} und wegen der 
Willk"urlichkeit von $\xi$ erhalten wir daraus die Integraldarstellung
der Arbeiten von Aseev \& Kryazhimskii und Aseev \& Veliov \cite{AseKry,AseVel,AseVel2,AseVel3}:

\begin{folgerung} \label{FolgerungPAUH3}
Es gen"uge $\big(x_*(\cdot),u_*(\cdot)\big)$ den Voraussetzungen des Pontrjaginschen Maximumprinzips \ref{SatzPAUH}.
Ist $\big(x_*(\cdot),u_*(\cdot)\big)$ ein starkes lokales Minimum der Aufgabe (\ref{PAUH1})--(\ref{PAUH4}) mit freiem 
Endpunkt im Unendlichen,
dann besitzt die Adjungierte $p(\cdot)$ die Darstellung
$$p(t)= -Z_*(t) \int_t^\infty \omega(s) Z^{-1}_*(s) f_x\big(s,x_*(s),u_*(s)\big) \, ds.$$
Dabei ist $Z_*(t)$ die in $t=0$ normalisierte Fundamentalmatrix des linearen Systems
$$\dot{z}(t)=-\varphi^T_x\big(t,x_*(t),u_*(t)\big) z(t).$$
\end{folgerung}

In der Aufgabe (\ref{PAUH1})--(\ref{PAUH4}) seien nun gewisse Randwerte im Unendlichen explizit gegeben, d.\,h. $h_1\big(t,x(t)\big)=x(t)-x_1$.
Wir schlie"sen dabei nicht aus,
dass dabei gewisse Komponenten von $x_1$ nicht fest vorgegeben, sondern ohne Einschr"ankung sind.
Wir sprechen dabei von expliziten Randbedingungen wenn f"ur $x(t)=\big(x_1(t),...,x_n(t)\big)$ gilt:
$$\lim_{t \to \infty} x_i(t)=x_i \in \R, \qquad \lim_{t \to \infty} x_j(t) \mbox{ frei}, \qquad i=1,...,l,\; j=l+1,...,n.$$
Daher verschwinden einerseits bei expliziten Randwerten im Unendlichen die entsprechenden Komponenten $x_i(t)-x_{*i}(t)$ f"ur $t \to \infty$, $i=1,...,l$.
F"ur diejenigen Komponenten, f"ur die die Randwerte im Unendlichen frei sind,
verschwinden die entsprechenden Komponenten $p_j(t)$ der Adjungierten f"ur $t \to \infty$, $j=l+1,...,n$.
Damit k"onnen wir festhalten:

\begin{folgerung} \label{FolgerungPAUH4}
In der Aufgabe (\ref{PAUH1})--(\ref{PAUH4}) mit expliziten Randbedingungen im Unendlichen
ergibt sich aus dem Maximumprinzip die ``nat"urliche'' Transversalit"atsbedingung
\index{Transversalitätsbedingungen!nat@-- natürliche}
$$\lim_{t \to \infty} \langle p(t),x(t)-x_*(t) \rangle = 0 \quad \mbox{ f"ur alle zul"assige } x(\cdot) \in W^1_\infty(\R_+,\R^n).$$
\end{folgerung}

In der Aufgabe (\ref{PAUH1})--(\ref{PAUH4}) mit freiem Endpunkt nehmen wir nun an,
dass die Verteilungsfunktion $\omega(\cdot) \in L_1(\R_+,\R_+)$ im Unendlichen verschwindet.
Dann folgen aus der gleichm"a"sigen Stetigkeit der Abbildungen $f,\,\varphi$ und aus der Beschr"anktheit des Steuerungsprozesses $\big(x_*(\cdot),u_*(\cdot)\big)$ 
die Grenzwerte
$$\lim_{t \to \infty} \omega(t) f\big(t,x_*(t),u_*(t)\big)=0, \qquad
  \lim_{t \to \infty} \big\langle p(t)\,,\, \varphi\big(t,x_*(t),u_*(t)\big)\big\rangle=0.$$

\begin{folgerung} \label{FolgerungPAUH5}
In der Aufgabe (\ref{PAUH1})--(\ref{PAUH4}) mit freiem Endpunkt im Unendlichen besitze die Verteilungsfunktion $\omega(\cdot)$
einen Grenzwert im
Unendlichen, d.\,h.
$$\lim_{t \to \infty} \omega(t)=0.$$
Dann ergibt sich die Bedingung von Michel \cite{Michel} \index{Transversalitätsbedingungen!von@-- von Michel}:
$$\lim_{t \to \infty} H\big(t,x_*(t),u_*(t),p(t),\lambda_0\big)=0.$$
\end{folgerung}

%% file: 4-33-Hinreichend.tex
\subsubsection{Hinreichende Bedingungen nach Arrow} \label{AbschnittArrowPMP}
Unser\index{hinreichende Bedingungen, Arrow!unendlich@-- unendlicher Zeithorizont}
Vorgehen zur Herleitung der hinreichenden Bedingungen basiert wieder auf der Darstellung in 
Seierstad \& Syds\ae ter \cite{Seierstad} (vgl. auch Abschnitt \ref{AbschnittHBPMP}).
Wir betrachten das Steuerungsproblem
\begin{eqnarray}
&& \label{HBPAUH1} J\big(x(\cdot),u(\cdot)\big) = \int_0^\infty \omega(t)f\big(t,x(t),u(t)\big) \, dt \to \inf, \\
&& \label{HBPAUH2} \dot{x}(t) = \varphi\big(t,x(t),u(t)\big), \\
&& \label{HBPAUH3} x(0)=x_0, \qquad \lim_{t \to \infty} x(t)=x_1, \\
&& \label{HBPAUH4} u(t) \in U \subseteq \R^m, \quad U \not= \emptyset.
\end{eqnarray}
Ebenso wie im Abschnitt \ref{AbschnittHBPMP} liegen in der Aufgabenstellung die Bedingungen (\ref{HBPAUH3}) in gemischter Form mit freien und festen
Randwerten vor. \\[2mm]
Wir betrachten die Menge $V^{\mathcal{S}}_\gamma(t)=\{ x \in \R^n \,|\, \|x-x_*(t)\| < \gamma\}$.
Au"serdem bezeichnet
$$\mathscr{H}^{\,\mathcal{U}}(t,x,p) = \sup\limits_{u \in U} H^{\mathcal{U}}(t,x,u,p,1)$$
die Hamilton-Funktion $\mathscr{H}^{\,\mathcal{U}}$ im normalen Fall.

\begin{theorem} \label{SatzHBPMPUH}
In der Aufgabe (\ref{HBPAUH1})--(\ref{HBPAUH4}) sei $\big(x_*(\cdot),u_*(\cdot)\big) \in \mathscr{B}^{\,\mathcal{U}}_{\rm Lip} \cap \mathscr{B}^{\,\mathcal{U}}_{\rm adm}$ 
und es sei $p(\cdot):\R_+ \to \R^n$. Ferner gelte:
\begin{enumerate}
\item[(a)] Das Tripel $\big(x_*(\cdot),u_*(\cdot),p(\cdot)\big)$
           erf"ullt (\ref{SatzPAUH1})--(\ref{SatzPAUH3}) in Theorem \ref{SatzPAUH} mit $\lambda_0=1$.        
\item[(b)] F"ur jedes $t \in \R_+$ ist die Funktion $\mathscr{H}^{\,\mathcal{U}}\big(t,x,p(t)\big)$ konkav in $x$ auf $V^{\mathcal{S}}_\gamma(t)$.
\end{enumerate}
Dann ist $\big(x_*(\cdot),u_*(\cdot)\big)$ ein starkes lokales Minimum der Aufgabe (\ref{HBPAUH1})--(\ref{HBPAUH4}).
\end{theorem}

{\bf Beweis} Auf die gleiche Weise wie im Beweis von Theorem \ref{SatzHBPMP} in Abschnitt \ref{AbschnittHBPMP}
ergibt sich zu $T \in \R_+$ die Beziehung
\begin{eqnarray*}
    \Delta(T)
&=& \int_0^T \omega(t)\big[f\big(t,x(t),u(t)\big)-f\big(t,x_*(t),u_*(t)\big)\big] \, dt \\
&\geq & \int_0^T\big[\mathscr{H}^{\,\mathcal{U}}\big(t,x_*(t),p(t)\big)-\mathscr{H}^{\,\mathcal{U}}\big(t,x(t),p(t)\big)\big] \, dt 
        + \int_0^T \langle p(t), \dot{x}(t)-\dot{x}_*(t) \rangle \, dt \\
&\geq& \langle p(T),x(T)-x_*(T)\rangle-\langle p(0),x(0)-x_*(0)\rangle.
\end{eqnarray*}

Im Fall fester Anfangs- und Endbedingungen verschwinden die Differenzen $x(0)-x_*(0)$ und $x(T)-x_*(T)$ f"ur $T \to \infty$.
Sind jedoch gewissen Komponenten im Anfangs- oder Endpunkt $x_0$ bzw. $x_1$ frei,
dann liefern die Transversalit"atsbedingungen,
dass die entsprechenden Komponenten der Adjungierten $p(\cdot)$ zum Zeitpunkt $t=0$ bzw. im Unendlichen verschwinden.
Daher folgt die Beziehung
$$\lim_{T \to \infty} \Delta(T) \geq \lim_{T \to \infty} \langle p(T),x(T)-x_*(T)\rangle-\langle p(0),x(0)-x_*(0)\rangle=0$$
f"ur alle zul"assigen $\big(x(\cdot),u(\cdot)\big)$ mit $\|x(\cdot)-x_*(\cdot)\|_\infty < \gamma$. \hfill $\blacksquare$

\begin{bemerkung}{\rm
Die Herleitung der hinreichenden Bedingungen in Theorem \ref{SatzHBPMPUH} basiert im Wesentlichen auf der Konkavit"at der 
Hamilton-Funktion $\mathscr{H}^{\,\mathcal{U}}$.
Daher m"ussen bei der Anwendung von Theorem \ref{SatzHBPMPUH} die Einschr"ankungen (\ref{PMPBedingung})
und (\ref{PMPBedingung2}) nicht gelten.
Dies hat weiterhin zur Folge, dass im Fall des freien Endpunktes die Trajektorie $x_*(\cdot)$ lediglich beschr"ankt sein muss.
\hfill $\square$}
\end{bemerkung}

\begin{beispiel} {\rm Im Beispiel \ref{BeispielRWUnendlich2} lieferten
die notwendigen Bedingungen (\ref{SatzPAUH1})--(\ref{SatzPAUH3})
den Schaltpunkt $\tau=\ln x_1$ und das Tripel
\begin{eqnarray*}
x_*(t) &=& \left\{ \begin{array}{ll} e^t, & t \in [0,\tau), \\ x_1, & t \in [\tau, \infty), \end{array}\right. \quad
  u_*(t)= \left\{ \begin{array}{ll} 1, & t \in [0,\tau), \\ 0, & t \in [\tau, \infty), \end{array}\right. \\
p(t) &=& \left\{ \begin{array}{ll}
        e^{(1-\varrho)\tau} e^{-t}, & t \in [0,\tau), \\
        \displaystyle\frac{\varrho-1}{\varrho}e^{-\varrho \tau} + \frac{1}{\varrho} e^{-\varrho t}, & t \in [\tau, \infty).
        \end{array}\right.
\end{eqnarray*}
Es ist die Hamilton-Funktion $\mathscr{H}^{\,\mathcal{U}}$ konkav in $x$ und damit $\big(x_*(\cdot),u_*(\cdot)\big)$
ein starkes lokales Maximum der Aufgabe.  \hfill $\square$}
\end{beispiel}

\begin{beispiel} 
{\rm In der Aufgabe des linear-quadratischen Reglers 
\begin{eqnarray*}
&& J\big(x(\cdot),u(\cdot)\big) = \int_0^\infty e^{-2t} \cdot \frac{1}{2}\big( x^2(t)+u^2(t)\big) \, dt \to \inf, \\
&& \dot{x}(t) = 2 x(t)+u(t), \qquad x(0)=2, \qquad u(t) \in \R
\end{eqnarray*}
sind die restriktiven Bedingungen (\ref{PMPBedingung}), (\ref{PMPBedingung2}) nicht erf"ullt und es ist das
Pontrjaginsche Maximumprinzip \ref{SatzPAUH} nicht anwendbar. \\
Es liefern die Bedingungen (\ref{SatzPAUH1})--(\ref{SatzPAUH3}) den Steuerungsprozess und die Adjungierte
$$x_*(t)=2e^{(1-\sqrt{2})t}, \quad u_*(t)=-2(1+\sqrt{2})e^{(1-\sqrt{2})t}, \quad p(t)=e^{-2t}u_*(t).$$
Da die Hamilton-Funktion $\mathscr{H}^{\,\mathcal{U}}$ konkav bez"uglich $x$ ist,
ist $\big(x_*(\cdot),u_*(\cdot)\big)$ nach Theorem \ref{SatzHBPMPUH} ein starkes lokales Minimum. \hfill $\square$}
\end{beispiel}

\begin{beispiel}
{\rm Im Beispiel \ref{BeipielRWUnendlich} mit Budgetbeschr"ankung,
\begin{eqnarray*}
&& J\big(x(\cdot),z(\cdot),u(\cdot)\big)=\int_0^\infty e^{-\varrho t}\big(1-u(t)\big)x(t) \, dt \to \sup, \\
&& \dot{x}(t)=u(t)x(t), \; x(0)=1,\qquad \dot{z}(t)=e^{-\varrho t}x(t), \; z(0)=0, \; \lim_{t \to \infty} z(t) = Z, \\
&& u \in [0,1], \qquad \varrho \in (0,1),
\end{eqnarray*}
ist jeder zul"assige Steuerungsprozess $\big(x(\cdot),z(\cdot),u(\cdot)\big)$ global optimal.
Die Voraussetzung $\big(x(\cdot),z(\cdot),u(\cdot)\big) \in \mathscr{B}^{\,\mathcal{U}}_{\rm adm}$ in Theorem \ref{SatzHBPMPUH} ist genau
dann erf"ullt,
wenn die Steuerung $u(\cdot)$ dem Raum $L_1(\R_+,[0,1])$ angeh"ort.
F"ur jeden zul"assigen Steuerungprozess sind mit den Multiplikatoren
$$\lambda_0=1, \qquad p(t)= e^{-\varrho t}, \qquad q(t)= \varrho-1$$
die notwendigen Bedingungen (\ref{SatzPAUH1})--(\ref{SatzPAUH3}) erf"ullt.
Weiterhin ist die Hamilton-Funktion $\mathscr{H}^{\,\mathcal{U}}$ offenbar bez"uglich $(x,z)$ konkav und es gilt Theorem
\ref{SatzHBPMPUH}.} \hfill $\square$
\end{beispiel}

%% file: 4-4-Zustandsaufgabe.tex
\subsection{Aufgaben mit Zustandsbeschr\"ankungen}
\subsubsection{Notwendige Optimalit\"atsbedingungen}
\begin{theorem}[Pontrjaginsches Maximumprinzip] \label{SatzPAUHZA}
\index{Pontrjaginsches Maximumprinzip!unendlich@-- unendlicher Zeithorizont} 
In der Aufgabe (\ref{PAUH1})--(\ref{PAUH5}) sei
$\big(x_*(\cdot),u_*(\cdot)\big) \in \mathscr{B}^{\,\mathcal{U}}_{\rm adm} \cap \mathscr{B}^{\,\mathcal{U}}_{\rm Lip} \cap \mathscr{B}^{\,\mathcal{U}}_{\lim}$.
Ist $\big(x_*(\cdot),u_*(\cdot)\big)$ ein starkes lokales Minimum  der Aufgabe (\ref{PAUH1})--(\ref{PAUH5}),
dann existieren eine Zahl $\lambda_0 \geq 0$,
Vektoren $l_0 \in \R^{s_0}$ und $l_1 \in \R^{s_1}$,
eine Vektorfunktion $p(\cdot):\R_+ \to \R^n$
und auf den Mengen
$$T_j=\big\{t \in \overline{\R}_+ \,\big|\, g_j\big(t,x_*(t)\big)=0\big\}, \quad j=1,...,l,$$
konzentrierte nichtnegative reguläre Borelsche Ma"se $\mu_j$ endlicher Totalvariation
(wobei s"amtliche Gr"o"sen nicht gleichzeitig verschwinden) derart, dass
die Vektorfunktion $p(\cdot)$ von beschr"ankter Variation und rechtsseitig stetig ist, und
\begin{enumerate}
\item[(a)] die adjungierte Gleichung
           \index{adjungierte Gleichung!unendlich@-- unendlicher Zeithorizont}
           \begin{eqnarray}
           p(t)&=&- \lim_{t \to \infty} h_{1,x}^T\big(t,x_*(t)\big) l_1 
                    + \int_t^\infty H^{\mathcal{U}}_x\big(s,x_*(s),u_*(s),p(s),\lambda_0\big) \, ds \nonumber \\
           \label{SatzPAUHZA1} & & -\sum_{j=1}^l \int_t^\infty g_{j,x}\big(s,x_*(s)\big)\, d\mu_j(s),
           \end{eqnarray}
\item[(b)] die Transversalit"atsbedingungen
           \index{Transversalitätsbedingungen!unendlich@-- unendlicher Zeithorizont}
           \begin{eqnarray}
          \label{SatzPAUHZA2} p(0) &=& {h_0'}^T\big(x_*(0)\big)l_0, \\
          \label{SatzPAUHZA4} \lim_{t \to \infty} p(t) 
                                   &=& \lim_{t \to \infty} \bigg[-h_{1,x}^T\big(t,x_*(t)\big) l_1 -\sum_{j=1}^l \mu_j(\{\infty\})\,g_{j,x}\big(t,x_*(t)\big) \bigg]
           \end{eqnarray}           
\item[(c)] und f"ur fast alle $t\in \R_+$ die Maximumbedingung
           \index{Maximumbedingung!unendlich@-- unendlicher Zeithorizont}
           \begin{equation}\label{SatzPAUHZA3}
           H^{\mathcal{U}}\big(t,x_*(t),u_*(t),p(t),\lambda_0\big) = \max_{u \in U} H^{\mathcal{U}}\big(t,x_*(t),u,p(t),\lambda_0\big)
           \end{equation}
\end{enumerate}
erfüllt sind.
\end{theorem}

Gem"a"s Definition \ref{DefinitionMassR} besitzen die regulären Borelschen Maße $\mu_j$
in Theorem \ref{SatzPAUHZA} die eindeutige Darstellung $\mu_j=\mu_{0j}+\mu_j(\{\infty\})$ mit
einem regul"aren Borelschen Ma"s $\mu_{0j}$ "uber $\R_+$ und mit einem in $t=\infty$ konzentrierten Maß $\mu_j(\{\infty\})$.
Ferner gilt nach Anhang \ref{AnhangMassR}: 
$$\int_t^\infty g_{j,x}\big(s,x_*(s)\big)\, d\mu_j(s) = \int_t^\infty g_{j,x}\big(s,x_*(s)\big)\, d\mu_{0j}(s)
  +\lim_{s \to \infty} \mu_j(\{\infty\})\, g_{j,x}\big(s,x_*(s)\big).$$
Da die Totalvariation der regulären Borelschen Maße $\mu_{0j}$ im Unendlichen verschwindet,
ergibt sich im Unendlichen die Transversalit"atsbedingung (\ref{SatzPAUHZA4}).

\begin{beispiel}[Abbau einer nicht erneuerbaren Ressource] \label{ExampleRessource}\index{Ressourcenabbau}
{\rm Wir betrachten die Aufgabe
\begin{eqnarray}
&& \label{Ressource1} J\big(x(\cdot),y(\cdot),u(\cdot)\big)
   =\int_0^\infty e^{-\varrho t}\big[pf\big(u(t)\big)-ry(t)-qu(t)\big] \, dt \to \sup, \\
&& \label{Ressource2} \dot{x}(t) = -u(t),\quad \dot{y}(t)=cf\big(u(t)\big), \quad x(0)=x_0>0,\quad y(0)=y_0\geq 0, \\
&& \label{Ressource3} x(t) \geq 0, \qquad u(t) \geq 0, \qquad b,c,\varrho, q,r >0, \qquad \varrho - rc>0.
\end{eqnarray}
Die Funktion $f$ sei zweimal stetig differenzierbar, $f'>0$, $f'(0)<\infty$, $f''<0$ und es sei
$f'(u) \to 0$ f"ur $u \to \infty$.
In der vorliegenden Formulierung der Aufgabe wurde im Vergleich zu Seierstad \& Syds\ae ter \cite{Seierstad} die
Restriktion $\liminf\limits_{t \to \infty} x(t) \geq 0$ durch die Zustandsbeschr"ankung $x(t) \geq 0$ in (\ref{Ressource3})
ersetzt. \\[2mm]
"Okonomische Interpretation:
$x(t)$ bezeichnet die Menge einer nat"urlichen Ressource und $u(t)$ ist die industrielle Abbaurate dieser Ressource.
Auf Basis der Ressource werden G"uter mit der Produktionsrate $f\big(u(t)\big)$ hergestellt.
Die Kosten der Herstellung einer Produktionseinheit ist $q$ und der Preis einer G"utereinheit am Markt betr"agt $p$.
Bei der Herstellung der G"uter entstehen proportional zur Produktion Abf"alle,
deren Gesamtmenge durch $y(t)$ beschrieben wird.
Die Kosten der Beseitigung der negativen Auswirkungen der Abfallmenge sind $ry(t)$.
Im Weiteren gehen wir von dem Preis $p=1$ aus. \\[2mm]
Wegen der Zustandsbeschr"ankung sind f"ur jeden zul"assigen Steuerungprozess die Beschr"ankungen an die Dynamik in
Theorem \ref{SatzPAUHZA} erf"ullt.
Wir stellen die Optimalit"atsbedingungen von Theorem \ref{SatzPAUHZA} mit $\lambda_0=1$ auf:
\begin{enumerate}
\item[(a)] Die Pontrjagin-Funktion der Aufgabe (\ref{Ressource1})--(\ref{Ressource3}) lautet
           $$H^{\mathcal{U}}(t,x,y,u,p_1,p_2,1) = p_1 (-u)+p_2 cf(u) + e^{-\varrho t}[f(u)-ry-qu].$$
\item[(b)] Die Adjungierten gen"ugen den Gleichungen
           $$p_1(t)=\int_t^\infty  \, d\mu(s), \qquad \dot{p}_2(t)=r e^{-\varrho t} \Rightarrow p_2(t)=-\frac{r}{\varrho}e^{-\varrho t} + K.$$
           Das auf der Menge $T=\{t \in \overline{\R}_+ \,|\, x_*(t)=0\}$ konzentrierte Ma"s $\mu$ ist nichtnegativ.
           Daher ist $p_1(t) \geq 0$ "uber $\R_+$ und monoton fallend.
           Ferner erhalten wir $K=0$ aus der Transversalit"atsbedingung (\ref{SatzPAUHZA4}) bez"uglich dem Zustand $y$.
\item[(c)] Die Maximumbedingung k"onnen wir auf folgende Aufgabe reduzieren:
           $$\max_{u \geq 0} \Big[ -p_1(t) u +c p_2(t) f(u) + e^{-\varrho t}[f(u)-qu]\Big].$$
           Das Einsetzen der Darstellung f"ur $p_2(t)$ liefert mit $d=(\varrho -rc)/\varrho$ weiterhin:
           $$\max_{u \geq 0} \Big( d f(u)e^{-\varrho t}-u\big(p_1(t)+ qe^{-\varrho t}\big)\Big).$$
\end{enumerate}
Die Reduktion der Maximumbedingung f"uhrt f"ur festes $t$ zu der Funktion
$$g(u)= d f(u)e^{-\varrho t}-u\big(p_1(t)+ qe^{-\varrho t}\big).$$
Diese Funktion ist zweimal stetig differenzierbar und es gilt
$$g'(u)= \big(d f'(u)-q\big)e^{-\varrho t}-p_1(t), \quad g''(u)=df''(u)e^{-\varrho t}, \quad d=(\varrho -rc)/\varrho>0.$$
Daher ist $g$ streng konkav und besitzt auf der Menge $U=\{u\geq 0\}$ ein Maximum,
da $f'(u)>0$ und $f'(u) \to 0$ f"ur $u \to \infty$ gelten.
Wir diskutieren drei F"alle:
\begin{enumerate}
\item[(A)] $df'(0)\leq q$: In diesem Fall ist $g'(0) \leq 0$ und man erh"alt
           $$u_*(t) \equiv 0, \quad x_*(t) \equiv x_0, \quad
             y_*(t)=y_0 + cf(0)t, \quad p_1(t) \equiv 0, \quad p_2(t)=-\frac{r}{\varrho}e^{-\varrho t}.$$
\item[(B)] $df'(0)> q$ und $p_1(0)=0$:
           Aus $g'(u)=\big(d f'(u)-q\big)e^{-\varrho t}=0$ erhalten wir die optimale Strategie $u_*(t)=u_0>0$ f"ur alle $t \in \R_+$.
           Also gilt $x_*(t)=x_0-u_0t$ auf $\R_+$, was $x_*(t) \geq 0$ widerspricht.
           In diesem Fall lässt sich kein optimaler Kandidat ableiten.
\item[(C)] $df'(0)> q$ und $p_1(0)>0$:
           Wegen $p_1(0)>0$ wird die Ressource vollst"andig abgebaut.
           Andernfalls w"are $p_1(t)\equiv p_1(0)>0$; im Widerspruch zu (\ref{SatzPAUHZA4}). \\
           Da die Ressource vollst"andig abgebaut wird,
           gibt es ein $t'>0$ mit $x_*(t)>0$ f"ur $t \in [0,t')$ und $x_*(t)=0$ f"ur $t\geq t'$.
           Demnach folgt unmittelbar $u_*(t)=0$ f"ur $t\geq t'$. \\
           F"ur $t\geq t'$ ist $p_1(\cdot)$ monoton fallend.
           Ferner erhalten wir f"ur $t \in \R_+$ die Beziehung
           $$g'(u)=0 \qquad\Rightarrow\qquad f'\big(u_*(t)\big)=\frac{1}{d}(q+p_1(t)e^{\varrho t}).$$
           W"urde demnach die Adjungierte $p_1(\cdot)$ f"ur $t \geq t'$ eine Unstetigkeitstelle besitzen,
           dann folgt aus der Monotonie von $p_1(\cdot)$, dass die Abbaurate sich wieder sprunghaft vergr"o"sert.
           Diese Steuerung f"uhrt zu einem erneuten Abbau der Ressource, obwohl diese bereits vollst"andig aufgebraucht ist.
           Daher ist $p_1(\cdot)$ in $t=t'$ stetig. \\
           F"ur die Adjungierte $p_1(\cdot)$ erhalten wir damit
           $$p_1(t) = \big(df'(0)-q\big)e^{-\varrho t'} \mbox{ f"ur } t \leq t', \qquad
             p_1(t) = \big(df'(0)-q\big)e^{-\varrho t} \mbox{ f"ur } t \geq t'.$$
           Wir zeigen noch, dass der Zeitpunkt $t' \in (0,\infty)$ existiert und eindeutig ist:
           Durch
           $$f'\big(u_\tau(t)\big)=\frac{1}{d}(q+p_1(0)e^{\varrho t})=\frac{1}{d}(q+[df'(0)-q]e^{\varrho(t-\tau)}), \quad t \in [0,\tau],$$
           und $u_\tau(t)=0$ f"ur $t \geq \tau$ wird wegen $f'\big(u_\tau(\tau)\big)=f'(0)$ eine Familie $u_\tau(\cdot)$ stetiger Funktionen definiert.
           Dabei gilt $f'\big(u_\tau(t)\big) < f'\big(u_s(t)\big)$, d.\,h. $u_\tau(t) > u_s(t)$, f"ur alle $t \in [0,\tau)$ und $\tau>s$.
           Ferner besitzen die Funktionen $u_\tau(\cdot)$ die obere Schranke $\hat{u}$ mit $f'(\hat{u})=q/d$.
           Außerdem gelten $u_1(0)=u_1>0$ und $u_{\tau+\Delta}(t+\Delta)=u_\tau(t)$ für $\tau >0$.
           Damit ist die Funktion $U(\tau)$,
           $$U(\tau):= \int_0^\infty u_\tau(t) \, dt, \quad U(\tau)=\int_0^\tau u_\tau(t) \, dt > \int_0^{\tau-1} u_1(0) \, dt=(\tau -1)u_1 \mbox{ für } \tau >1,$$
           streng monoton wachsend und es gelten $U(0)=0$, $U(\tau) \to \infty$ für $\tau \to \infty$.
           Der Parameter $t'$ ergibt sich dann aus der Bedingung $U(t')=x_0$. \hfill $\square$
\end{enumerate}}
\end{beispiel}

%% file: 4-41-Beweis.tex
\subsubsection{Der Nachweis der notwendigen Optimalit\"atsbedingungen}
Wir betrachten f"ur $\big(x(\cdot),u(\cdot)\big) \in C_{\lim}(\R_+,\R^n) \times L_\infty(\R_+,\R^m)$ die Abbildungen
\begin{eqnarray*}
J\big(x(\cdot),u(\cdot)\big) &=& \int_0^\infty \omega(t)f\big(t,x(t),u(t)\big) \, dt, \\
F\big(x(\cdot),u(\cdot)\big)(t) &=& x(t) -x(t_0) -\int_0^t \varphi\big(s,x(s),u(s)\big) \, ds, \quad t \in \R_+,\\
H_0\big(x(\cdot)\big) &=& h_0\big(x(0)\big), \qquad H_1\big(x(\cdot)\big) = \lim_{t \to \infty} h_1\big(t,x(t)\big), \\
G_j\big(x(\cdot)\big) &=& \max_{t \in \overline{\R}_+} g_j\big(t,x(t)\big), \quad j=1,...,l.
\end{eqnarray*}
Die Abbildungen fassen wir in folgenden Funktionenr"aumen auf:
\begin{eqnarray*}
J &:& C_{\lim}(\R_+,\R^n) \times L_\infty(\R_+,\R^m) \to \R, \\
F &:& C_{\lim}(\R_+,\R^n) \times L_\infty(\R_+,\R^m) \to C_{\lim}(\R_+,\R^n), \\
H_i &:& C_{\lim}(\R_+,\R^n) \to \R^{s_i}, \quad i=0,1, \\
G_j &:& C_{\lim}(\R_+,\R^n) \to \R, \quad j=1,...,l.
\end{eqnarray*}

Wir setzen $\mathscr{F}=(F,H_0,H_1)$, sowie die Menge $\mathscr{U}$ gem"a"s
\begin{eqnarray*}
\mathscr{U}=\big\{ u(\cdot)  \in L_\infty(\R_+,U) \!\!\!&\big|&\!\!\! u(t)=u_*(t) + \chi_M(t)\big(w(t)-u_*(t)\big), \; w(\cdot) \in L_\infty(\R_+,U), \\
                                                  \!\!\!&     &\!\!\! M \subset \R_+ \mbox{ me"sbar und beschr"ankt} \big\}
\end{eqnarray*}
und pr"ufen f"ur die Extremalaufgabe
\begin{equation} \label{ExtremalaufgabePMPUHZA}
J\big(x(\cdot),u(\cdot)\big) \to \inf, \quad \mathscr{F}\big(x(\cdot),u(\cdot)\big)=0, \quad G_j\big(x(\cdot)\big) \leq 0,
\quad u(\cdot) \in \mathscr{U}
\end{equation}
im Punkt $\big(x_*(\cdot),u_*(\cdot)\big) \in \mathscr{B}^{\,\mathcal{U}}_{\rm Lip} \cap \mathscr{B}^{\,\mathcal{U}}_{\lim}$
die Voraussetzungen von Theorem \ref{SatzExtremalprinzipStark}:

\begin{enumerate}
\item[(A$_2$)] Mit Verweis auf Abschnitt \ref{AbschnittBeweisPMPUH} sind nur noch die Abbildungen $G_j$ zu diskutieren.
               In den Beispielen \ref{SubdifferentialMaximum3} und \ref{SubdifferentialMaximum4} wird gezeigt,
               dass die Funktionen $G_j$ Hintereinanderausf"uhrungen einer stetigen, konvexen, eigentlichen Funktion und
               einer Fr\'echet-differenzierbaren Abbildung sind.
               Daher sind nach Lemma \ref{LemmaRichtungsableitung} die Funktionen $G_j$ in $x_*(\cdot)$ lokalkonvex und bez"uglich jeder Richtung
               gleichm"a"sig differenzierbar.              
\end{enumerate}

Zur Extremalaufgabe (\ref{ExtremalaufgabePMPUHZA}) definieren wir auf
$$C_{\lim}(\R_+,\R^n) \times L_\infty(\R_+,\R^m) \times \R \times C_{\lim}^*(\R_+,\R^n) \times \R^n \times \R^{s_0} \times  \R^{s_1}$$
die Lagrange-Funktion $\mathscr{L}=\mathscr{L}\big(x(\cdot),u(\cdot),\lambda_0,y^*,l_0,l_1,\lambda\big)$,
$$\mathscr{L}= \lambda_0 J\big(x(\cdot),u(\cdot)\big)+ \big\langle y^*, F\big(x(\cdot),u(\cdot)\big) \big\rangle
                         +l_0^T H_0\big(x(\cdot)\big)+l_1^T H_1\big(x(\cdot)\big) + \sum_{j=1}^l \lambda_j G_j\big(x(\cdot)\big).$$
Ist $\big(x_*(\cdot),u_*(\cdot)\big)$ eine starke lokale Minimalstelle der Extremalaufgabe (\ref{ExtremalaufgabePMPUHZA}),
dann existieren nach Theorem \ref{SatzExtremalprinzipStark}
nicht gleichzeitig verschwindende Lagrangesche Multiplikatoren
$$\lambda_0 \geq 0, \qquad y^* \in C_{\lim}^*(\R_+,\R^n), \qquad l_0 \in  \R^{s_0}, \qquad l_1 \in  \R^{s_1}$$
und $\lambda_1 \geq 0,...,\lambda_l \geq 0$ derart,
dass gelten:
\begin{enumerate}
\item[(a)] Die Lagrange-Funktion besitzt bez"uglich $x(\cdot)$ in $x_*(\cdot)$ einen station"aren Punkt, d.\,h.
           \begin{equation}\label{SatzPMPUHZALMR1}
           0 \in \partial_x \mathscr{L}\big(x_*(\cdot),u_*(\cdot),\lambda_0,y^*,l_0,l_1,\lambda\big);
           \end{equation}         
\item[(b)] Die Lagrange-Funktion erf"ullt bez"uglich $u(\cdot)$ in $u_*(\cdot)$ die Minimumbedingung
           \begin{equation}\label{SatzPMPUHZALMR2}
           \hspace*{-3mm} \mathscr{L}\big(x_*(\cdot),u_*(\cdot),\lambda_0,y^*,l_0,l_1,\lambda\big)
           = \min_{u(\cdot) \in \mathscr{U}} \mathscr{L}\big(x_*(\cdot),u(\cdot),\lambda_0,y^*,l_0,l_1,\lambda\big);
           \end{equation}
\item[(c)] Die komplement"aren Schlupfbedingungen gelten, d.\,h.
           \begin{equation}\label{SatzPMPUHZALMR3}
           0 = \lambda_j G_j\big(x(\cdot)\big), \qquad i=1,...,l.
           \end{equation}
\end{enumerate}
Aufgrund (\ref{SatzPMPUHZALMR1}) ist folgende Variationsgleichung f"ur alle $x(\cdot) \in C_{\lim}(\R_+,\R^n)$ erf"ullt: 
\begin{eqnarray*}
0 &=& \lambda_0 \cdot \int_0^\infty \omega(t) \big\langle f_x\big(t,x_*(t),u_*(t)\big),x(t) \big\rangle\, dt \\
  & & + \int_0^\infty \bigg[ x(t)-x(0) - \int_0^t \varphi_x\big(s,x_*(s),u_*(s)\big) x(s) \,ds \bigg]^T  d\mu(t) \\
  & & + \big\langle l_0, h_0'\big(x_*(0)\big) x(0) \big\rangle + \big\langle l_1, \lim_{t \to \infty} h_{1,x}\big(t,x_*(t)\big) x(t) \big\rangle \\
  & & + \sum_{j=1}^l \lambda_j \int_0^\infty \big\langle g_{j,x}\big(t,x_*(t)\big),x(t) \big\rangle \,d\tilde{\mu}_j(t).
\end{eqnarray*}
Ebenso wie im Abschnitt \ref{AbschnittBeweisPMPZA} ergibt sich,
dass alle Ma"se $\mu_j=\lambda_j \tilde{\mu}_j$ auf den Mengen
$$T_j=\big\{t \in \overline{\R}_+ \,\big|\, g_j\big(t,x_*(t)\big)=0\big\}, \quad j=1,...,l,$$
konzentriert sind.
Wir "andern die Integrationsreihenfolge und setzen $p(t)= \displaystyle \int_t^\infty \, d\mu(s)$.
Die Funktion $p(\cdot)$ ist von beschränkter Variation und gemäß den Eigenschaften einer Verteilungsfunktion rechtsseitig stetig über $\R_+$.
Dann ergeben sich mit den gleichen Argumenten wie in den Abschnitten \ref{AbschnittBeweisPMPZA} und \ref{AbschnittBeweisPMPUH}
die Bedingungen (\ref{SatzPAUHZA1}) und (\ref{SatzPAUHZA2}).
Gem"a"s (\ref{SatzPMPUHZALMR2}) gilt f"ur alle $u(\cdot) \in \mathscr{U}$ die Ungleichung
$$\int_0^\infty H^{\mathcal{U}}\big(t,x_*(t),u_*(t),p(t),\lambda_0\big) \, dt \geq
  \int_0^\infty H^{\mathcal{U}}\big(t,x_*(t),u(t),p(t),\lambda_0\big) \, dt.$$
Daraus folgt abschlie"send durch Standardtechniken f"ur Lebesguesche Punkte die Maximumbedingung (\ref{SatzPAUHZA3}).

%% file: 4-42-Hinreichend.tex
\subsubsection{Hinreichende Bedingungen nach Arrow} \label{AbschnittArrowPMPZB}
Die folgende Herleitung von hinreichenden Bedingungen
\index{hinreichende Bedingungen, Arrow!unendlich@-- unendlicher Zeithorizont}
in einem Steuerungsproblem mit gemischten freien und festen Randwerten im Unendlichen
basiert wieder auf der Darstellung in Seierstad \& Syds\ae ter \cite{Seierstad}.
Es wird das folgende Steuerungsproblem untersucht:
\begin{eqnarray}
&& \label{HBPAUHZB1} J\big(x(\cdot),u(\cdot)\big) = \int_0^\infty \omega(t)f\big(t,x(t),u(t)\big) \, dt \to \inf, \\
&& \label{HBPAUHZB2} \dot{x}(t) = \varphi\big(t,x(t),u(t)\big), \\
&& \label{HBPAUHZB3} x(0)=x_0, \qquad \lim_{t \to \infty} x(t)=x_1, \\
&& \label{HBPAUHZB4} u(t) \in U \subseteq \R^m, \quad U \not= \emptyset, \\
&& \label{HBPAUHZB5} g_j\big(t,x(t)\big) \leq 0 \quad \mbox{f"ur alle } t \in \R_+, \quad j=1,...,l.
\end{eqnarray}
Gemäß den vorhergehenden Betrachtungen und Beweisführungen besitzen in dieser Aufgabe die einfließenden Elemente sinnvolle Grenzwerte im Unendlichen.
Deswegen beruht der Beweis der hinreichenden Arrow-Bedingungen auf den bekannten Argumenten über dem endlichen Intervall $[0,T]$ und dem
anschließenden Grenzübergang $T \to \infty$. \\[2mm]
Wir betrachten die Menge $V^{\mathcal{S}}_\gamma(t)=\{ x \in \R^n \,|\, \|x-x_*(t)\| < \gamma\}$.
Au"serdem bezeichnet
$$\mathscr{H}^{\,\mathcal{U}}(t,x,p) = \sup\limits_{u \in U} H^{\mathcal{U}}(t,x,u,p,1)$$
die Hamilton-Funktion $\mathscr{H}^{\,\mathcal{U}}$ im normalen Fall.

\begin{theorem} \label{SatzHBPMPUHZB}
In der Aufgabe (\ref{HBPAUHZB1})--(\ref{HBPAUHZB5}) sei $\big(x_*(\cdot),u_*(\cdot)\big) \in \mathscr{B}^{\,\mathcal{U}}_{\rm Lip} \cap \mathscr{B}^{\,\mathcal{U}}_{\rm adm}$.
Au"serdem sei die Vektorfunktion $p(\cdot):\R_+ \to \R^n$ st"uckweise stetig,
besitze h"ochstens abz"ahlbar viele Sprungstellen $s_k \in (0,\infty)$,
die sich nirgends im Endlichen h"aufen,
und $p(\cdot)$ sei zwischen diesen Spr"ungen stetig differenzierbar. 
Ferne gelte:
\begin{enumerate}
\item[(a)] Das Tripel $\big(x_*(\cdot),u_*(\cdot),p(\cdot)\big)$
           erf"ullt (\ref{SatzPAUH1})--(\ref{SatzPAUH3}) in Theorem \ref{SatzPAUH} mit $\lambda_0=1$.        
\item[(b)] F"ur jedes $t \in \R_+$ ist die Funktion $\mathscr{H}^{\,\mathcal{U}}\big(t,x,p(t)\big)$ konkav 
           und es sind die Funktionen $g_j(t,x)$, $j=1,...,l$, konvex bez"uglich $x$ auf $V^{\mathcal{S}}_\gamma(t)$.
\end{enumerate}
Dann ist $\big(x_*(\cdot),u_*(\cdot)\big)$ ein starkes lokales Minimum der Aufgabe (\ref{HBPAUHZB1})--(\ref{HBPAUHZB5}).
\end{theorem}

{\bf Beweis} Wie im Abschnitt \ref{AbschnittHBPMP} ergibt sich die Ungleichung
\begin{eqnarray*}
    \Delta(T)
&=& \int_0^T \omega(t)\big[f\big(t,x(t),u(t)\big)-f\big(t,x_*(t),u_*(t)\big)\big] \, dt \\
&\geq & \int_0^T\big[\mathscr{H}^{\,\mathcal{U}}\big(t,x_*(t),p(t)\big)-\mathscr{H}^{\,\mathcal{U}}\big(t,x(t),p(t)\big)\big] \, dt 
        + \int_0^T \langle p(t), \dot{x}(t)-\dot{x}_*(t) \rangle \, dt.
\end{eqnarray*}
F"ur $T \not=s_k$ ergibt sich ferner mit den Argumenten aus Abschnitt \ref{AbschnittArrowZB}:
\begin{eqnarray*}
    \Delta(T)
&\geq& \int_0^T \langle \dot{p}(t),x(t)-x_*(t)\rangle + \langle p(t), \dot{x}(t)-\dot{x}_*(t) \rangle \, dt \\
&&    - \int_0^T \sum_{j=1}^l \lambda_j(t) \big\langle g_{j,x}\big(t,x_*(t)\big) , x(t)-x_*(t) \big\rangle \, dt \\
&&    \hspace*{20mm} - \sum_{s_k <T} \langle p(s_k)-p(s_k^-),x(s_k)-x_*(s_k)\rangle \\
&\geq& \int_0^T \langle \dot{p}(t),x(t)-x_*(t)\rangle + \langle p(t), \dot{x}(t)-\dot{x}_*(t) \rangle \, dt \\
&=& \langle p(T),x(T)-x_*(T)\rangle-\langle p(0),x(0)-x_*(0)\rangle.
\end{eqnarray*}
Da $p(\cdot)$ eine Funktion beschr"ankter Variation ist, gilt au"serdem
$$\lim_{k \to \infty} \|p(s_k)-p(s_k^-)\|=0.$$
Damt ist der Grenz"ubergang $T \to \infty$ gerechtfertigt und wir erhalten aus der Transversalit"atsbedingung
$$\lim_{t \to \infty} p(t) = \lim_{t \to \infty} \bigg[l_1-\sum_{j=1}^l \mu_j(\{\infty\})\, g_{j,x}\big(t,x_*(t)\big)\bigg]$$
die Ungleichung
$$\lim_{T \to \infty} \Delta(T) = \lim_{T \to \infty} \int_0^T \omega(t)\big[f\big(t,x(t),u(t)\big)-f\big(t,x_*(t),u_*(t)\big)\big] \, dt \geq 0$$
f"ur alle zul"assigen $\big(x(\cdot),u(\cdot)\big)$ mit $\|x(\cdot)-x_*(\cdot)\|_\infty < \gamma$. \hfill $\blacksquare$

\begin{beispiel}
{\rm Im Beispiel \ref{ExampleRessource}, \index{Ressourcenabbau}
\begin{eqnarray*}
&& J\big(x(\cdot),y(\cdot),u(\cdot)\big) =\int_0^\infty e^{-\varrho t}\big[pf\big(u(t)\big)-ry(t)-qu(t)\big] \, dt \to \sup, \\
&& \dot{x}(t) = -u(t),\quad \dot{y}(t)=cf\big(u(t)\big), \quad x(0)=x_0>0,\quad y(0)=y_0\geq 0, \\
&& x(t) \geq 0, \qquad u(t) \geq 0, \qquad b,c,\varrho, q,r >0, \qquad \varrho - rc>0,
\end{eqnarray*}
enth"alt die Pontrjagin-Funktion
$$H^{\,\mathcal{U}}(t,x,y,u,p_1,p_2,1) = p_1 (-u)+p_2 cf(u) + e^{-\varrho t}[f(u)-ry-qu]$$
bez"uglich der Zustandsvariablen $x$ und $y$ nur den Term $e^{-\varrho t}ry$.
Damit ist die Hamilton-Funktion $\mathscr{H}^{\mathcal{U}}$ konkav bez"uglich $(x,y)$.
Au"serdem ist die Zustandsbeschr"ankung linear in $x$.
Damit sind die gefundenen Kandidaten in den F"allen (A) und (C) optimal.}
\end{beispiel}

%% file: 4-5-Einordnung.tex
\subsection{Steuerungsprobleme mit endlichem und unendlichem Zeithorizont}
Der unendliche Zeithorizont wird oft als sehr langes, aber endliches Zeitintervall aufgefasst.
Daraus ergibt sich die Frage,
ob ein Steuerungsproblem mit unendlichem Zeithorizont lediglich ein Standardproblem mit endlichem Zeithorizont darstellt.
Damit ist zu klären in welcher Beziehung die Standardaufgaben über endlichem und unendlichem Zeithorizont zueinander stehen;
es fehlt bisher eine Einordnung der Aufgabe mit unendlichem Zeithorizont in die Theorie der Optimalen Steuerung.
Dieser Unternehmung wollen wir uns in diesem Abschnitt stellen. \\[2mm]
Ein häufig angewandter Zugang zu Steuerungsproblemen mit unendlichem Zeithorizont ist
die Approximation durch eine Aufgabe mit endlichem Zeithorizont.
Allerdings zeigt sich bei diesem Ansatz,
dass der vollst"andige Satz an notwendigen Optimalit"atsbedingungen verloren geht und
das resultierende Pontrjaginsche Maximumprinzip im Hinblick auf die Angabe korrekter Transversalitätsbedingungen
unvollständig ist. \\[2mm]
Weiterhin basieren Zugänge zu Aufgaben mit dem unbeschränkten Zeitintervall auf der R"uckf"uhrung auf ein endliches
Zeitintervall mittels der Substitution der Zeit.
Hierbei offenbart sich das Wesen des unendlichen Zeithorizonts als eine Singularit"at in der transformierten Aufgabe
über endlichem Zeitintervall.
Auf diese Weise entsteht also keine Standardaufgabe,
sondern durch Auftreten der Singularität eine kaum untersuchte Aufgabenklasse. \\[2mm]
Damit sind wir wieder an unserer zentralen Fragestellung:
Welche Beziehung besteht zwischen den Aufgaben mit endlichem und unendlichem Zeithorizont?
Wir werden zum Abschluss dieses Abschnitts zeigen,
dass die entwickelten Methoden und erzielten Ergebnisse für Steuerungsprobleme mit unendlichem Zeithorizont
eine direkte Verallgemeinerung für die entsprechenden Argumente und Resultate über endlichem Zeithorizont darstellen.
Dies unterstreicht die besonderen Herausforderungen,
mit denen wir bei der Behandlung der Aufgabe mit unendlichem Zeithorizont konfrontiert wurden.

%% file: 4-51-Approximation.tex
\subsubsection{Zur Approximation durch eine Aufgabe \"uber endlichem Horizont}
Es bezeichnet $\mathscr{B}_{\rm adm}$ die Menge aller zul"assigen Steuerungsprozesse in der Aufgabe
\begin{eqnarray}
&& \label{EUH1} J\big(x(\cdot),u(\cdot)\big) = \int_0^\infty \omega(t)f\big(t,x(t),u(t)\big) \, dt \to \inf, \\
&& \label{EUH2} \dot{x}(t) = \varphi\big(t,x(t),u(t)\big), \\
&& \label{EUH3} h_0\big(x(0)\big)=0, \qquad \lim_{t \to \infty} h_1\big(t,x(t)\big)=0, \\
&& \label{EUH4} u(t) \in U \subseteq \R^m, \quad U \not= \emptyset, \\
&& \label{EUH5} g_j\big(t,x(t)\big) \leq 0 \quad \mbox{f"ur alle } t \in \R_+, \quad j=1,...,l.
\end{eqnarray}
D.\,h. die Menge aller $\big(x(\cdot),u(\cdot)\big)$,
die der Dynamik (\ref{EUH2}) zur Anfangs- und Randbedingung (\ref{EUH3}),
den Steuerrestriktionen (\ref{EUH4}) und den Zustandsbeschr"ankungen (\ref{EUH5}) gen"ugen,
und f"ur die das Zielfunktional (\ref{EUH1}) endlich ist.
Zur Gegen"uberstellung von verschiedenen Optimalit"atsbegriffen betrachten wir globale Optimalit"atskriterien.

\begin{enumerate}
\item[(a)] Der Steuerungsprozess $\big(x_*(\cdot),u_*(\cdot)\big) \in \mathscr{B}_{\rm adm}$ hei"st
           global optimal (GO)\index{optimal, catching up!global@--, global}, falls
           $$J\big(x_*(\cdot),u_*(\cdot)\big) \leq J\big(x(\cdot),u(\cdot)\big)$$
           f"ur alle $\big(x(\cdot),u(\cdot)\big) \in \mathscr{B}_{\rm adm}$ gilt.
\item[(b)] Der Steuerungsprozess $\big(x_*(\cdot),u_*(\cdot)\big) \in \mathscr{B}_{\rm adm}$ hei"st
           streng global optimal (SGO)\index{optimal, catching up!streng global@--, streng global}, falls
           $$J\big(x_*(\cdot),u_*(\cdot)\big) < J\big(x(\cdot),u(\cdot)\big)$$
           f"ur alle $\big(x(\cdot),u(\cdot)\big) \in \mathscr{B}_{\rm adm}$ mit
           $\big(x(\cdot),u(\cdot)\big) \not= \big(x_*(\cdot),u_*(\cdot)\big)$ gilt.
\end{enumerate}
Zur Untersuchung des Steuerungsproblems mit unendlichem Zeithorizont wird in den meisten Ans"atzen der Aufgabe
(\ref{EUH1})--(\ref{EUH5})
ein Problem "uber einem endlichen Zeitintervall $[0,T]$ zugeordnet und der Grenz"ubergang $T \to \infty$ betrachtet.
Wir nennen dies die ``Approximation durch endlichen Horizont''. \index{Approximation durch endlichen Horizont} \\[2mm]
Das zur Aufgabe (\ref{EUH1})--(\ref{EUH5}) geh"orende Problem mit endlichem Zeithorizont lautet
\begin{eqnarray}
&& \label{EUH1a} J_T\big(x(\cdot),u(\cdot)\big) = \int_0^T \omega(t) f\big(t,x(t),u(t)\big)\, dt \to \inf, \\
&& \label{EUH2a} \dot{x}(t) = \varphi\big(t,x(t),u(t)\big), \\
&& \label{EUH3a} h_0\big(x(0)\big)=0, \qquad h_1\big(T,x(T)\big)=0, \\
&& \label{EUH4a} u(t) \in U \subseteq \R^m, \quad U \not= \emptyset, \\
&& \label{EUH5a} g_j\big(t,x(t)\big) \leq 0 \quad \mbox{f"ur alle } t \in [0,T], \quad j=1,...,l.
\end{eqnarray}
Wegen des "Uberganges zu einer Aufgabe "uber endlichem Zeithorizont und der anschlie"senden Betrachtung des Grenzwertes
$T \to \infty$ werden die Optimalit"atsbegriffe angepasst. \\
Zur Aufgabe (\ref{EUH1a})--(\ref{EUH5a}) definieren wir den Defekt $\Delta(T)$:
$$\Delta(T)= J_T\big(x(\cdot),u(\cdot)\big) - J_T\big(x_*(\cdot),u_*(\cdot)\big).$$
\begin{enumerate}
\item[(c)] Der Steuerungsprozess $\big(x_*(\cdot),u_*(\cdot)\big) \in \mathscr{B}_{\rm adm}$ hei"st
           overtaking optimal (OT)\index{optimal, catching up!overtaking@--, overtaking},
           falls es zu jedem $\big(x(\cdot),u(\cdot)\big) \in \mathscr{B}_{\rm adm}$ eine Zahl $T_0$ gibt mit
           $$\Delta(T)= J_T\big(x(\cdot),u(\cdot)\big) - J_T\big(x_*(\cdot),u_*(\cdot)\big)\geq 0$$
           f"ur alle $T \geq T_0$ (von Weizs"acker \cite{VonWeiz}).
\item[(d)] Der Steuerungsprozess $\big(x_*(\cdot),u_*(\cdot)\big) \in \mathscr{B}_{\rm adm}$ hei"st
           catching up optimal (CU)\index{optimal, catching up}, falls
           $$\liminf_{T \to \infty} \Delta(T)
             = \liminf_{T \to \infty} J_T\big(x(\cdot),u(\cdot)\big) - J_T\big(x_*(\cdot),u_*(\cdot)\big) \geq 0$$
           f"ur jedes $\big(x(\cdot),u(\cdot)\big) \in \mathscr{B}_{\rm adm}$ gilt (Gale \cite{Gale}).
\item[(e)] Das Paar $\big(x_*(\cdot),u_*(\cdot)\big) \in \mathscr{B}_{\rm adm}$ hei"st
           sporadically catching up optimal (SCU)\index{optimal, catching up!sporadically@--, sporadically catching up}, falls
           $$\limsup_{T \to \infty} \Delta(T)
              = \limsup_{T \to \infty} J_T\big(x(\cdot),u(\cdot)\big) - J_T\big(x_*(\cdot),u_*(\cdot)\big)\geq 0$$
           f"ur jedes $\big(x(\cdot),u(\cdot)\big) \in \mathscr{B}_{\rm adm}$ gilt (Halkin \cite{Halkin}).
\end{enumerate}
Durch die Betrachtung auf der Menge $\mathscr{B}_{\rm adm}$ haben wir in den Definitionen bereits eingef"ugt,
dass f"ur einen zul"assigen Steuerungsprozess das Zielfunktional endlich ist.
In diversen Arbeiten zu Steuerungsproblemen mit unendlichem Zeithorizont,
z.\,B. bei Aseev \& Veliov \cite{AseVel,AseVel2,AseVel3}, darf das Integral in (\ref{EUH1}) divergieren. \\
Konvergiert das Zielfunktional f"ur jeden zul"assigen Steuerungsprozess,
so sind die Optimalit"atsbegriffe (CU), (SCU) und (GO) "aquivalent, denn es gilt
$$\liminf_{T \to \infty} \Delta(T) \geq 0 \quad\Leftrightarrow\quad
  \limsup_{T \to \infty} \Delta(T) \geq 0 \quad\Leftrightarrow\quad \lim_{T \to \infty} \Delta(T)\geq 0.$$
Zur Einordnung der Begriffe (SGO) und (OT) f"uhren wir folgendes Beispiel an:

\begin{beispiel}
{\rm Wir betrachten die Aufgabe
$$\int_0^\infty -u(t)x(t) \, dt \to \inf, \qquad \dot{x}(t) = -u(t)x(t), \quad x(0)=1,
  \quad u(t) \in [0,1], \quad \varrho >0.$$
Offenbar ist jede zul"assige Trajektorie $x(\cdot)$ monoton fallend und wegen $1 \geq x(t) > 0$ beschr"ankt.
Also besitzt sie deswegen einen Grenzwert f"ur $t \to \infty$.
Damit gilt
$$J\big(x(\cdot),u(\cdot)\big)=\int_0^\infty -u(t)x(t) \, dt = \int_0^\infty \dot{x}(t) \, dt = \lim_{t \to \infty} x(t)-1 \geq -1.$$
Daher ist jedes Paar $\big(x(\cdot),u(\cdot)\big) \in \mathscr{B}_{\rm adm}$ mit $x(t) \to 0$ global optimal.
Dabei gilt:
$$\lim_{t \to \infty} x(t) = 0 \quad \Leftrightarrow \quad \int_0^\infty u(t) \, dt = \infty.$$
In diesem Beispiel ist jeder Steuerungsprozess $\big(x(\cdot),u(\cdot)\big) \in \mathscr{B}_{\rm adm}$
mit $x(t) \to 0$ f"ur $t \to \infty$
optimal im Sinn von (GO), (CU) und (SCU).
Wegen $J_T\big(x(\cdot),u(\cdot)\big)\geq x(T)-1$ ist unter diesen Steuerungsprozessen nur derjenige optimal im Sinn von (OT),
f"ur den $u(t) \equiv 1$ auf $\R_+$ gilt.
Eine (SGO)-L"osung existiert nicht. \hfill $\square$}
\end{beispiel}
Dieses Beispiel verdeutlicht,
dass unter der Annahme eines endlichen Zielfunktionals zu den bereits erw"ahnten "Aquivalenzen folgende weitere Relationen
zwischen den verschiedenen Optimalit"atsbegriffen bestehen:
$$\mbox{(SGO) } \Longrightarrow \mbox{ (OT) } \Longrightarrow \mbox{ (CU) } \Longleftrightarrow \mbox{ (SCU) }
  \Longleftrightarrow \mbox{ (GO)}.$$
Ist demnach die globale Optimalstelle eindeutig, so fallen bei endlichem Zielfunktional
die aufgef"uhrten Optimalit"atsbegriffe zusammen. \\[2mm]
Die Betrachtungen dieses Abschnitts ordnen die Begriffe der globalen Optimalit"at in den Kontext 
g"angiger Optimalit"atsbegriffe ein,
die bei der Approximation mit einem endlichen Horizont angewendet werden.
Unbehandelt bleibt dabei allerdings die wesentliche Frage,
wann eine Familie $\big\{\big(x_*^T(\cdot),u_*^T(\cdot)\big)\big\}_{T \in \R_+}$ von optimalen Steuerungsprozessen,
die sich in der Approximation ergibt,
gegen ein globales Optimum konvergiert.
Dazu muss man sich vorab im Klaren sein,
dass f"ur die Approximation die fundamentale Beziehung
$$\lim_{T \to \infty} \int_0^T f(t) \, dt = \int_0^\infty f(t) \, dt$$
im Allgemeinen nur bei vorliegender Existenz des Lebesgue-Integrals im Zielfunktional gilt.
Ein bekanntes Gegenbeispiel ist der Integralsinus:
$$\lim_{T \to \infty} \int_0^T \frac{\sin t}{t} \, dt = \frac{\pi}{2}, \qquad
  \int_0^\infty \frac{\sin t}{t} \, dt \mbox{ existiert nicht.}$$
Noch weniger "uberschaubar wird die Situation in der Aufgabe (\ref{EUH1})--(\ref{EUH5}),
in der die Dynamik, sowie Zustands- und Steuerungsbeschr"ankungen vorkommen.
Dann k"onnen noch weitere unerw"unschte Situationen eintreten.
Dazu geben wir das folgende Beispiel an:

\begin{beispiel} \label{BeispielInvestmentUH}
{\rm Wir betrachten die Aufgabe
$$\int_0^T e^{-\varrho t} \big(1-u(t)\big)x(t) \, dt \to \sup, \quad
  \dot{x}(t) = u(t)x(t), \quad x(0)=1, \quad u(t) \in [0,1], \quad \varrho \in (0,1).$$
F"ur jedes feste $T>\tau$ liefert der Steuerungsprozess
$$x^T_*(t)= \left\{ \begin{array}{ll} e^t,& t \in [0,\tau), \\[1mm] e^\tau,& t \in [\tau,T], \end{array} \right. \quad
  u^T_*(t)= \left\{ \begin{array}{ll} 1,& t \in [0,\tau), \\[1mm] 0,& t \in [\tau,T], \end{array} \right. \quad
  \tau = T + \frac{\ln (1-\varrho)}{\varrho}$$
das globale Maximum.
Betrachten wir den Grenz"ubergang $T \to \infty$,
dann konvergiert die Familie $\big\{\big(x_*^T(\cdot),u_*^T(\cdot)\big)\big\}_{T \in \R_+}$ punktweise gegen den
Steuerungsprozess
$$x_*(t) = e^t, \qquad u_*(t)=1, \qquad t \in \R_+.$$
Dieses Paar ist das globale Minimum der Aufgabe mit unendlichem Zeithorizont. \hfill $\square$}
\end{beispiel}

Ein wesentlicher Augenmerk bei der Herleitung notwendiger Optimalit"atsbedingungen liegt auf der G"ultigkeit der ``nat"urlichen'' Transversalit"atsbedingungen
$$\lim_{t \to \infty} \|p(t)\| = 0, \qquad \lim_{t \to \infty} \langle p(t),x(t) \rangle = 0 \quad
  \mbox{f"ur alle zul"assigen } x(\cdot)$$
in der Aufgabe  (\ref{EUH1})--(\ref{EUH4}) mit freiem Endpunkt.
In diesem Zusammenhang gibt Halkin \cite{Halkin} folgendes Gegenbeispiel an:
\begin{beispiel} {\rm In der Aufgabe
$$\int_0^\infty u(t)\big(1-x(t)\big)\, dt \to \sup, \quad \dot{x}(t)=u(t)\big(1-x(t)\big), \quad x(0)=0,\quad u \in [0,1],$$
liefert der Steuerungsprozess
$$u_*(t)=\left\{\begin{array}{ll}
  0,& t \in [0,1], \\ 1,& t \in [1,\infty),\end{array} \right. \quad
  x_*(t)=\left\{\begin{array}{ll}
  0,& t \in [0,1], \\ 1-e^{1-t},& t \in [1,\infty),\end{array} \right.$$
ein globales Optimum.
F"ur die zugeh"origen Multiplikatoren $\big(\lambda_0,p(\cdot)\big)\not=(0,0)$,
mit denen die Bedingungen des Maximumprinzips erf"ullt sind,
gilt dann $p(t) \equiv -\lambda_0$; ein Widerspruch zur Gültigkeit der ``nat"urlichen'' Transversalit"atsbedingung. \hfill $\square$}
\end{beispiel}

\begin{beispiel} {\rm In der folgenden Variante eines Beispieles nach Halkin \cite{Halkin},
\begin{eqnarray*}
&& J\big(x(\cdot),u(\cdot)\big) = \int_0^\infty e^{-\varrho t}\big(u(t)-x(t)\big) \, dt \to \sup, \\
&& \dot{x}(t) = x(t)+u(t), \quad x(0)=0, \quad  u(t) \in [0,1], \quad \varrho \in (0,1),
\end{eqnarray*}
erhalten wir für einen zulässigen Steuerungsprozess $\big(x(\cdot),u(\cdot)\big)$ für den Zustand
$$x(t)=e^t \cdot \int_0^t u(s)e^{-s} \, ds, \quad t \geq 0.$$
Damit liefert der Steuerungsprozess $\big(x_*(t),u_*(t)\big) \equiv (0,0)$ das globale Maximum der Aufgabe.
Denn f"ur jeden anderen zulässigen Steuerungsprozess mit einer stückweise stetigen Steuerung,
die nicht identisch verschwindet, besitzt das Zielfunktional den Wert $-\infty$.  \\
Wir werten die Bedingungen des Pontrjaginschen Maximumprinzips im allgemeinen Fall aus:
Die Pontrjagin-Funktion lautet $H^{\mathcal{U}}(t,x,u,p,\lambda_0)=p(x+u)+\lambda_0 e^{-\varrho t}(u-x)$.
Für die Lösung $p(\cdot)$ der adjungierten Gleichung ergibt sich
$$\dot{p}(t)=-p(t)+\lambda_0 e^{-\varrho t} \quad\Rightarrow\quad
  p(t)=\bigg(p(0)-\frac{\lambda_0}{1-\varrho}\bigg) e^{-t}+\frac{\lambda_0}{1-\varrho} e^{-\varrho t},\; p(0) \in \R,$$
und es gilt $p(t) \to 0$ für $t \to \infty$.
Die Maximumbedingung lautet
$$H^{\mathcal{U}}\big(t,x_*(t),u_*(t),p(t),\lambda_0\big) =
  \max_{u \in [0,1]} \Big[p(t)\big(x_*(t)+u\big)+\lambda_0 e^{-\varrho t}\big(u-x_*(t)\big)\Big]$$
und ist äquivalent zu der Maximierungsaufgabe
$$\max_{u \in [0,1]} u\big[p(t)+\lambda_0 e^{-\varrho t}\big]
  =\max_{u \in [0,1]} u \cdot \bigg[ \bigg(p(0)-\frac{\lambda_0}{1-\varrho}\bigg) e^{-t}+\frac{2-\varrho}{1-\varrho} \lambda_0e^{-\varrho t}\bigg].$$
Wäre nun $\lambda_0 >0$,
so fällt wegen $\varrho \in (0,1)$ der Ausdruck in der eckigen Klammer für alle hinreichend große $t$ positiv aus,
und $u_*(t) \equiv 0$ kann nicht der Maximumbedingung für fast alle $t \geq 0$ genügen.
Deswegen muss der anormale Fall der notwendigen Optimalitätsbedingungen mit $\lambda_0=0$ eintreten. \hfill $\square$}
\end{beispiel}

Die letzten beiden Beispiele zeigen,
dass ohne strenge Vorkehrungen die Transversalitätsbedingungen in Abschnitt \ref{AbschnittNormalformUH}
oder für die Aufgabe mit freiem Endpunkt der normale Fall mit $\lambda_0=1$
nicht gültig sein müssen.

%% file: 4-52-Approximationssatz.tex
\subsubsection{Ein Resultat zur Approximation mit endlichem Horizont}
Im letzten Abschnitt wurden die verschiedenen Optimaliätätsbegriffe (GO), (SGO) und (OT), (CU), (SCU) eingeführt
und die Beziehungen bei endlichem Zielfunktional diskutiert.
Dabei kam der Konvergenz des Zielfunktionals die tragende Rolle zu. \\
Dies erweitern wir nun durch folgenden Approximationssatz:

\begin{satz} \label{LemmaKonvergenzFiniteApproximation}
In der Aufgabe (\ref{EUH1a})--(\ref{EUH5a}) sei $\big\{\big(x_*^T(\cdot),u_*^T(\cdot)\big)\big\}_{T \in \R_+}$ eine Familie
von globalen Minimalstellen, f"ur die folgendes erfüllt sei:
\begin{enumerate}
\item[(A)] Es existiert der endliche Grenzwert $\displaystyle \lim_{T \to \infty} J_T\big(x_*^T(\cdot),u_*^T(\cdot)\big)$.
\item[(B)] Es gibt einen zul"assigen Steuerungsprozess $\big(x_*(\cdot),u_*(\cdot)\big)$ 
           der Aufgabe (\ref{EUH1})--(\ref{EUH5}) mit 
           $$\lim_{T \to \infty} J_T\big(x_*^T(\cdot),u_*^T(\cdot)\big) = J\big(x_*(\cdot),u_*(\cdot)\big).$$
\end{enumerate}
Dann ist $\big(x_*(\cdot),u_*(\cdot)\big)$ ein globales Minimum der Aufgabe (\ref{EUH1})--(\ref{EUH5}).
\end{satz}

{\bf Beweis} Sei $\big(x(\cdot),u(\cdot)\big)$ zul"assig in der Aufgabe (\ref{EUH1})--(\ref{EUH5}).
Insbesondere ist dann das Zielfunktional in (\ref{EUH1}) endlich.
Dann l"asst sich zu jedem $\varepsilon>0$ eine Zahl $T'>0$ derart angeben,
dass die Einschr"ankung des Zielfunktionals $J\big(x(\cdot),u(\cdot)\big)$
auf das Intervall $[T,\infty)$ vom Betrag kleiner oder gleich $\varepsilon$ f"ur alle $T\geq T'$ ausf"allt.
Ferner kann die Zahl $T'$ so gew"ahlt werden, dass
$\big|J_T\big(x_*^T(\cdot),u_*^T(\cdot)\big) - J\big(x_*(\cdot),u_*(\cdot)\big)\big| \leq \varepsilon$ f"ur alle $T\geq T'$ gilt.
Mit den suggestiven Bezeichnungen, dass $J_T$ bzw. $J$ die Integration des Zielfunktionals "uber $[0,T]$ bzw. "uber $\R_+$ angeben,
k"onnen wir die Differenz
$J\big(x(\cdot),u(\cdot)\big) - J\big(x_*(\cdot),u_*(\cdot)\big)$
in die Form
\begin{eqnarray*}
&& \Big[J\big(x(\cdot),u(\cdot)\big) - J_T\big(x(\cdot),u(\cdot)\big)\Big] 
+\Big[J_T\big(x(\cdot),u(\cdot)\big) - J_T\big(x_*^T(\cdot),u_*^T(\cdot)\big)\Big] \\
 &&   + \Big[J_T\big(x_*^T(\cdot),u_*^T(\cdot)\big) - J\big(x_*(\cdot),u_*(\cdot)\big)\Big]
\end{eqnarray*}
bringen.
Aufgrund der Wahl von $T'$ fallen der erste und dritte Summand gr"o"ser oder gleich $-\varepsilon$ f"ur alle $T\geq T'$ aus.
Der zweite Ausdruck ist nichtnegativ,
da $\big(x_*^T(\cdot),u_*^T(\cdot)\big)$ ein globales Minimum "uber $[0,T]$ darstellt.
Daher gilt $J\big(x(\cdot),u(\cdot)\big) - J\big(x_*(\cdot),u_*(\cdot)\big) \geq -2\varepsilon$
mit beliebigem $\varepsilon>0$. \hfill $\blacksquare$ \\[2mm]
Der Satz \ref{LemmaKonvergenzFiniteApproximation} ist zwar eine hinreichende Bedingung, aber:
\begin{enumerate}
\item[1.)] Der Satz trifft keine Annahmen an eine Konvergenz der Familie  
           $\big\{\big(x_*^T(\cdot),u_*^T(\cdot)\big)\big\}$.
           Im Allgemeinen bleibt damit die Frage nach der Gestalt von $\big(x_*(\cdot),u_*(\cdot)\big)$ offen.
\item[2.)] Selbst bei gewisser Konvergenz der Familie $\big\{\big(x_*^T(\cdot),u_*^T(\cdot)\big)\big\}$ muss der
           resultierende Steuerungsprozess keine Lösung liefern.
\item[3.)] Der Satz enthält keinen Existenznachweis für den Steuerungsprozesses $\big(x_*(\cdot),u_*(\cdot)\big)$.
\end{enumerate}
Bei Aseev \& Kryazhimskii \cite{AseKry}, Theorem 2.1, sind Voraussetzungen angegeben,
unter denen in der Aufgabe (\ref{EUH1})--(\ref{EUH4}) mit freiem Endpunkt die Annahmen von
Satz \ref{LemmaKonvergenzFiniteApproximation} erf"ullt sind.
Ferner liefert dabei in \cite{AseKry} die Approximation durch die Aufgabe (\ref{EUH1a})--(\ref{EUH4a}) mit endlichem Horizont
eine Familie von globalen Minimalstellen,
die auf jedem endlichen Intervall gleichm"a"sig gegen das globale Minimum der Aufgabe (\ref{EUH1})--(\ref{EUH4})
konvergiert. \\[2mm]
Wir geben einige Beispiele an, die verschiedene Aspekte von Satz \ref{LemmaKonvergenzFiniteApproximation} dokumentieren.
In den ersten beiden Beispielen wird die Bedeutung der Voraussetzungen (A) und (B) hervorgehoben;
in den weiteren Beispielen wird die Anwendbarkeit verdeutlicht:

\begin{beispiel}
{\rm Im Beispiel \ref{BeispielInvestmentUH} lieferten
$$x^T_*(t)= \left\{ \begin{array}{ll} e^t,& t \in [0,\tau), \\[1mm] e^\tau,& t \in [\tau,T], \end{array} \right. \quad
  u^T_*(t)= \left\{ \begin{array}{ll} 1,& t \in [0,\tau), \\[1mm] 0,& t \in [\tau,T], \end{array} \right. \quad
  \tau = T + \frac{\ln (1-\varrho)}{\varrho}$$
eine Familie globaler Maxima.
Da der Grenzwert von
$$J_T\big(x_*^T(\cdot),u_*^T(\cdot)\big) = \int_0^T e^{-\varrho t} \big(1-u_*(t)\big)x_*(t) \, dt
  =\int_\tau^T e^{-\varrho t} e^\tau\, dt =(1-\varrho)^{\frac{1-\varrho}{\varrho}} \cdot e^{(1-\varrho)T}$$
für $T \to \infty$ nicht endlich ist,
ist die Voraussetzung (A) nicht erfüllt und
Satz \ref{LemmaKonvergenzFiniteApproximation} darf nicht angewendet werden. \hfill $\square$}
\end{beispiel}

\begin{beispiel} 
{\rm Wesentlicher Punkt beim Nachweis von Satz \ref{LemmaKonvergenzFiniteApproximation} war die Endlichkeit des
Zielfunktionals.
Das Auftreten einer Gewichtsfunktion $\omega(\cdot)$ wurde nicht ausdrücklich gefordert.
Diesbezüglich betrachten wir eine Aufgabe ohne gewichtetes Zielfunktional:
$$J\big(x(\cdot)\big) = \int_0^\infty \big[\big(x(t)-1\big)^2+\dot{x}^2(t)\big] \, dt \to \inf,
   \qquad x(0)=0,\quad \lim_{t \to \infty} x(t)=2.$$
Für $T>0$ ergibt sich die Grundaufgabe im Beispiel \ref{BeispielELGRegler},
$$J_T\big(x(\cdot)\big) = \int_0^T \big[\big(x(t)-1\big)^2+\dot{x}^2(t)\big] \, dt \to \inf, \qquad x(0)=0,\quad x(T)=2,$$
mit der Lösung
$$x^T_*(t)=1+\frac{1+e^{-T}}{e^T-e^{-T}}e^t- \frac{1+e^T}{e^T-e^{-T}}e^{-t}.$$
Im Grenzwert $T \to \infty$ ergibt sich die Trajektorie $x_*(t)=1-e^{-t}$,
die in der Aufgabe mit unendlichem Zeithorizont nicht zulässig ist,
da sie die Randbedingung $\displaystyle \lim_{t \to \infty} x(t)=2$ verletzt.
Obwohl der Grenzwert $\displaystyle \lim_{T \to \infty} J_T\big(x_*^T(\cdot),u_*^T(\cdot)\big)$ existiert und endlich ist,
gibt es keine zulässige Trajektorie $\big(x_*(\cdot),u_*(\cdot)\big)$,
die die Voraussetzung (B) in Satz \ref{LemmaKonvergenzFiniteApproximation} erfüllt.
Die Aufgabe mit unendlichem Zeithorizont besitzt keine Lösung. \hfill $\square$}
\end{beispiel}

\begin{beispiel} 
{\rm Die Voraussetzungen an die Aufgabe (\ref{PAUH1})--(\ref{PAUH5}) schlie"st Probleme mit zyklischen Daten bzw.
Lösungen aus,
da gewisse Voraussetzungen über Konvergenz im Unendlichen erfüllt sein müssen.
Eine einfache Aufgabe mit periodischen Daten ist
\begin{eqnarray*}
&& J\big(x(\cdot),u(\cdot)\big)
   = \frac{1}{2}\int_0^\infty \big[\big(x(t)-\cos t\big)^2+\big(u(t)+\sin t\big)^2\big] \, dt \to \inf, \\
&& \dot{x}(t)=u(t), \qquad x(0)=0, \qquad u(t) \in \R.
\end{eqnarray*}
Über jedem endlichen Horizont $[0,T]$ liefern das Pontrjaginsche Maximumprinzip und die Arrow-Bedingungen
in Kapitel \ref{KapitelStark} die globale Lösung
$$x_*^T(t)=c_1(T)e^t-c_2(T)e^{-t}+\cos t, \qquad u_*^T(t)=c_1(T)e^t+c_2(T)e^{-t}-\sin t$$
zu den Koeffizienten
$$c_1(T)=-\frac{e^{-T}}{e^T+e^{-T}}, \qquad c_2(T)=\frac{e^{T}}{e^T+e^{-T}}.$$
In der nachstehenden Abbildung sind die optimalen Trajektorien $x_*^1(t)$, $x_*^2(t)$, $x_*^3(t)$ zu $T=1,2,3,$ und
$x_*(t)$ für $T= \infty$ dargestellt:
\begin{figure}[h]
	\centering
	\fbox{\includegraphics[width=10cm]{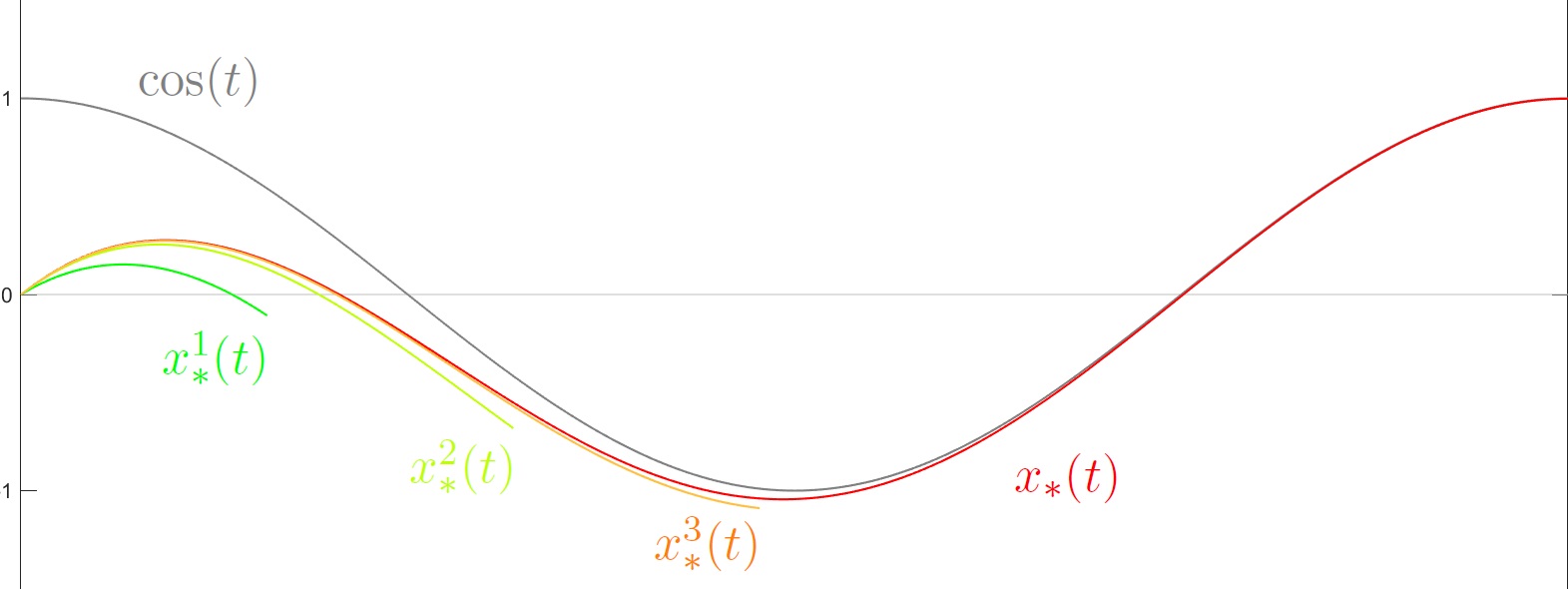}}
\caption[Linear-quadratischer Regler mit periodischen Daten]{Optimale Trajektorien für $T=1,\,2,\,3$ und $T=\infty$.}
\end{figure}

Im Grenzwert $T \to \infty$ ergibt sich der Steuerungsprozess
$$x_*(t)=-e^{-t}+\cos t, \qquad u_*(t)=e^{-t}-\sin t$$
und für den Wert des Zielfunktionals
$$\lim_{T \to \infty} J_T\big(x_*^T(\cdot),u_*^T(\cdot)\big)
  = \lim_{T \to \infty} \frac{1}{2}\cdot\frac{1-e^{-4T}}{(1+e^{-2T})^2}
  = \frac{1}{2} = J\big(x_*(\cdot),u_*(\cdot)\big).$$
Damit ist der im Unendlichen divergente Prozess $\big(x_*(\cdot),u_*(\cdot)\big)$ die globale Lösung. \hfill $\square$}
\end{beispiel}

%% file: 4-53-Zeittransformation.tex
\subsubsection{Transformation der Zeit}
Ein Zugang zu Aufgaben mit unendlichem Zeithorizont ist das Zur"uckf"uhren auf eine Aufgabe
"uber endlichem Zeitintervall mit Hilfe der Substitution der Zeit.
Dabei muss eine Zeitransformation stets eine eineindeutige Zuordnung zwischen den Intervallen $[0,1]$ und $[0,\infty]$
liefern; also der Stelle $1$ den Wert $\infty$ zuordnen.
In der Literatur gebräuchliche Vertreter derartiger Transformationen $s:[0,\infty] \to [0,1]$ bzw. $t:[0,1] \to [0,\infty]$
sind
\begin{enumerate}
\item[$\cdot$] die kanonische Kompaktifizierung
               $$s(t)=\frac{t}{1+t}, \qquad t(s)=\frac{s}{1-s},$$
\item[$\cdot$] die Abbildung
               $$s(t)=1-e^{-t}, \qquad t(s)=-\ln(1-s).$$
\end{enumerate}
Betrachten wir z.\,B. die Transformation $t(s)= -\ln(1-s)$,
so gilt
$$t'(s)=v(s)=\frac{1}{1-s}>0 \quad \mbox{ für } s \in [0,1)$$
und es wird die Aufgabe (\ref{EUH1})--(\ref{EUH5}) eineindeutig
in die Aufgabe (\ref{EUH1a})--(\ref{EUH5a}) "uber dem endlichen Intervall $[0,1]$ "uberf"uhrt.
Im Abschnitt \ref{AbschnittFreieZeitSOP} wurde gezeigt,
dass damit über dem Intervall $[0,1]$ die Beziehungen (\ref{FZSOP6}) und (\ref{FZSOP7}),
d.\,h.
$$t'(s) = v(s), \quad y'(s) = v(s) \cdot \varphi\big(t(s),y(s),w(s)\big), \quad
  \int_0^1 v(s) \cdot f\big(t(s),y(s),w(s)\big) \, ds,$$
folgen.
Beachten wir au"serdem die Zustandsbeschränkungen,
so besitzt die transformierte Aufgabe nach Abschnitt \ref{AbschnittFreieZeit} die Form
\begin{eqnarray*}
&& \int_0^1 v(s) \cdot \omega\big(t(s)\big) f\big(t(s),y(s),w(s)\big) \, ds \to \inf, \\
&& y'(s) = v(s) \cdot \varphi\big(t(s),y(s),w(s)\big), \qquad t'(s) = v(s), \\
&& y(0)=x_0,\quad h\big(t(1),y(1)\big)=0, \qquad t(0)=0, \quad t(1)=\infty, \\
&& w(s) \in U, \qquad v(s) > 0 \mbox{ für } s \in [0,1), \\
&& g_j\big(t(s),y(s)\big) \leq 0, \quad s \in [0,1], \quad j=1,...,l.
\end{eqnarray*}
Die Aufgabe enth"alt für den Zustand $t(\cdot)$ die Singularit"at $t(1)=\infty$,
die in s"amtlichen Elementen dieser Aufgabe einflie"st,
und au"serdem die "uber $[0,1]$ nicht integrable und unbeschr"ankte Funktion $v(\cdot)$.
Die auf $[0,1]$ "uberf"uhrte Aufgabenstellung geh"ort damit nicht in den Rahmen der klassischen Steuerungsprobleme.
Sondern die Substitution der Zeit verdeutlicht viel mehr,
dass der unendliche Horizont in der Aufgabe (\ref{EUH1})--(\ref{EUH5}) selbst eine Singularit"at ist,
die man nicht aus der Betrachtung argumentieren kann.
Also kann die Aufgabe (\ref{EUH1})--(\ref{EUH5}) kein klassisches Problem der Optimalen Steuerung sein. \\[2mm]
Das Wesen des unendlichen Zeithorizontes als eine Singularität in der Aufgabe bei der Zurückführung auf ein endliches Zeitintervall
mittels der Substitution der Zeit ist nicht unbekannt (vgl. \cite{AseKry}),
findet aber in der Literatur weitestgehend keine Beachtung.

%% file: 4-54-Einordnung.tex
\subsubsection{Die Einordnung der Aufgabenklassen mit verschiedenen Horizonten}
In diesem Abschnitt gehen wir der grundlegenden Frage nach,
wie sich die Aufgabenklassen mit endlichem und mit unendlichem Zeithorizont zueinander einordnen. \\
In den letzten Abschnitten zeigte sich,
dass weder die Approximation der Aufgabe mit einem endlichen Zeithorizont noch die Transformation der Zeit im Allgemeinen
geeignete Methoden sind,
um Steuerungsprobleme mit unendlichem Zeithorizont zu behandeln.
Ferner wurde durch das Auftreten der Singularität nach einer Zeittransformation verdeutlicht,
in welch umfassend Ma"se sich die Steuerungsprobleme mit unendlichem Zeithorizont von Aufgaben über endlichem Zeitintervall
unterscheiden. \\[1mm]
Damit stellt sich umgekehrt die Frage,
ob sich die Aufgabe mit endlichem Zeitintervall der Aufgabe mit unendlichem Horizont unterordnet:
Die Aufgabe (\ref{PMP1})--(\ref{PMP5}) "uber endlichem Zeithorizont in Kapitel \ref{KapitelStark} lautete
\begin{equation} \label{EinordnungEH} \left. \begin{array}{l}
  J\big(x(\cdot),u(\cdot)\big) = \displaystyle \int_{t_0}^{t_1} f\big(t,x(t),u(t)\big) \, dt \to \inf, \\[1mm]
  \dot{x}(t) = \varphi\big(t,x(t),u(t)\big), \\[1mm]
  h_0\big(x(t_0)\big)=0, \qquad h_1\big(x(t_1)\big)=0, \\[1mm]
  u(t) \in U \subseteq \R^m, \quad U\not= \emptyset, \\[1mm]
  g_j\big(t,x(t)\big) \leq 0 \quad \mbox{f"ur alle } t \in [t_0,t_1], \quad j=1,...,l.
  \end{array} \right\}
\end{equation}
Bez"uglich des Integranden $f$ und der rechten Seite $\varphi$ der Dynamik setzen wir:
$$\tilde{f}(t,x,u)= \left\{ \begin{array}{ll} f(t,x,u), & t \in [t_0,t_1], \\ 0, & t \not \in [t_0,t_1], \end{array} \right. \qquad
  \tilde{\varphi}(t,x,u)= \left\{ \begin{array}{ll} \varphi(t,x,u), & t \in [t_0,t_1], \\ 0, & t \not \in [t_0,t_1]. \end{array} \right.$$
Mit den Randbedingungen verkn"upfen wir die Abbildungen
$$\tilde{h}_0(x)= h_0(x), \qquad \tilde{h}_1(t,x)=h_1(x).$$
Au"serdem setzen wir bez"uglich der Zustandsbeschr"ankungen
$$\tilde{g}_j(t,x) = \left\{ \begin{array}{ll} g_j(t_0,x) - (1-e^{(t-t_0)^2}), & t < t_0, \\
                                               g_j(t,x), & t \in [t_0,t_1], \\
                                               g_j(t_1,x) - (1-e^{(t-t_1)^2}), & t > t_1. \end{array}\right.$$
Mit der Verteilungsfunktion $\omega(t) = \chi_{[t_0,t_1]}(t)$ ergibt sich auf diese Weise "uber dem unendlichen Zeithorizont die Aufgabe
\begin{equation} \label{EinordnungUH} \left. \begin{array}{l}
  J\big(x(\cdot),u(\cdot)\big) = \displaystyle \int_0^\infty \omega(t)\tilde{f}\big(t,x(t),u(t)\big) \, dt \to \inf, \\[1mm]
  \dot{x}(t) = \tilde{\varphi}\big(t,x(t),u(t)\big), \\[1mm]
  \tilde{h}_0\big(x(0)\big)=0, \qquad \displaystyle\lim_{t \to \infty} \tilde{h}_1\big(t,x(t)\big)=0, \\[1mm]
  u(t) \in U \subseteq \R^m, \quad U \not= \emptyset, \\[1mm]
  \tilde{g}_j\big(t,x(t)\big) \leq 0 \quad \mbox{f"ur alle } t \in \R_+, \quad j=1,...,l.
  \end{array} \right\}
\end{equation}
Die Aufgabe (\ref{EinordnungUH}) besitzt damit die Gestalt der Aufgabe (\ref{EUH1})--(\ref{EUH5}).
Bez"uglich der Verweise auf die Theoreme \ref{SatzPAUH} und \ref{SatzPAUHZA} ist anzumerken,
dass die unstetigen Anschl"usse der Abbildungen $\tilde{f}$ und $\tilde{\varphi}$ in den Stellen $t=t_0$ und $t=t_1$ sich nicht
nachteilig auf die Beweise der Maximumprinzipien auswirken. \\
Aus den Festlegungen ist ersichtlich,
dass ein Steuerungsprozess $\big(\tilde{x}(\cdot),\tilde{u}(\cdot)\big)$ genau dann in der Aufgabe (\ref{EinordnungUH}) zul"assig bzw. 
ein starkes lokales Minimum ist,
wenn die Einschr"ankung $\big(x(\cdot),u(\cdot)\big)$ mit 
$\big(x(t),u(t)\big)= \big(\tilde{x}(t),\tilde{u}(t)\big)$ f"ur $t \in [t_0,t_1]$
zul"assig bzw. ein starkes lokales Minimum in der Aufgabe (\ref{EinordnungEH}) darstellt.
Insbesondere stimmen dann bez"uglich der Zustandsbeschr"ankungen die Mengen
$$T_j = \big\{t \in [t_0,t_1] \,\big|\, g_j\big(t,x(t)\big)=0\big\}, \qquad 
  \tilde{T}_j = \big\{t \in \overline{\R}_+=[0,\infty] \,\big|\, \tilde{g}_j\big(t,\tilde{x}(t)\big)=0\big\}$$
f"ur $j=1,...,l$ "uberein.
Bei der Anwendung der Theoreme \ref{SatzPAUH} und \ref{SatzPAUHZA} auf die Aufgabe (\ref{EinordnungUH}) sind
die restriktiven Annahmen (\ref{PMPBedingung}) und (\ref{PMPBedingung2}) zu beachten,
die der Setzung nach f"ur $\tilde{\varphi}(t,x,u)$ offensichtlich erf"ullt sind.
Daher existieren nach dem Pontrjaginschen Maximumprinzip (Theorem \ref{SatzPAUHZA}) 
f"ur ein starkes lokales Minimum $\big(x_*(\cdot),u_*(\cdot)\big)$ der Aufgabe (\ref{EinordnungUH})
eine Zahl $\lambda_0 \geq 0$,
Vektoren $l_0 \in \R^{s_0}$, $l_1 \in \R^{s_1}$,
eine Vektorfunktion $p(\cdot):\R_+ \to \R^n$
und auf den Mengen
$$\tilde{T}_j=\big\{t \in \overline{\R}_+ \,\big|\, \tilde{g}_j\big(t,x_*(t)\big)=0\big\}, \quad j=1,...,l,$$
konzentrierte nichtnegative Borelsche Ma"se $\mu_j$ endlicher Totalvariation
(wobei s"amtliche Gr"o"sen nicht gleichzeitig verschwinden) derart, dass
die Vektorfunktion $p(\cdot)$ von beschr"ankter Variation und rechtsseitig stetig ist, und
\begin{enumerate}
\item[(a)] die adjungierte Gleichung
           \begin{eqnarray*}
           p(t)&=& - \lim_{t \to \infty}\tilde{h}_{1,x}^T\big(t,x_*(t)\big) l_1
                   + \int_t^\infty H^{\mathcal{U}}_x\big(s,x_*(s),u_*(s),p(s),\lambda_0\big) \, ds \\
               & & -\sum_{j=1}^l \int_t^\infty \tilde{g}_{j,x}\big(s,x_*(s)\big)\, d\mu_j(s),
           \end{eqnarray*}
\item[(b)] die Transversalit"atsbedingungen
           \begin{eqnarray*}
           p(0) &=& \tilde{h}_0{'}^T\big(x_*(0)\big)l_0, \\
           \lim_{t \to \infty} p(t) &=& \lim_{t \to \infty} \bigg[-\tilde{h}_{1,x}^T\big(x_*(t)\big)l_1
                                         -\sum_{j=1}^l \mu_j(\{\infty\})\, \tilde{g}_{j,x}\big(t,x_*(t)\big)\bigg]
           \end{eqnarray*}
\item[(c)] und f"ur fast alle $t\in \R_+$ die Maximumbedingung
           $$H^{\mathcal{U}}\big(t,x_*(t),u_*(t),p(t),\lambda_0\big) = \max_{u \in U} H^{\mathcal{U}}\big(t,x_*(t),u,p(t),\lambda_0\big)$$
\end{enumerate}
erfüllt sind. 

\newpage
In der Aufgabe ohne Zustandsbeschr"ankungen gehört $p(\cdot)$ dem Raum $W^1_\infty(\R_+,\R^n)$ an,
denn die Verteilungsfunktion $\omega(t) = \chi_{[t_0,t_1]}(t)$ ist messbar und beschränkt.
Weiterhin ergeben sich wegen $\dot{p}(t)=0$ f"ur $t \not \in [t_0,t_1]$ aus den Transversalit"atsbedingungen unmittelbar
$$p(t_0) = {h_0'}^T\big(x_*(t_0)\big)l_0, \qquad p(t_1)= -\tilde{h}^T_{1,x}\big(t_1,x_*(t_1)\big)l_1 =-{h_1'}^T\big(x_*(t_1)\big)l_1.$$
Das Pontrjaginsche Maximumprinzip \ref{SatzPMP} "uber dem endlichen Horizont in Kapitel \ref{KapitelStark} l"asst sich also unmittelbar aus
Theorem \ref{SatzPAUH} ableiten. \\[2mm]
In der Aufgabe mit Zustandsbeschr"ankungen gelten $g_j\big(t,x_*(t)\big) <0$ f"ur alle $t \not \in [t_0,t_1]$
und es sind die Ma"se $\mu_j$ auf den Mengen $T_j = \big\{t \in [t_0,t_1] \,\big|\, g_j\big(t,x(t)\big)=0\big\}$ konzentriert.
Daher gelten für $j=1,...,l$ die Beziehungen
$$ \lim_{t \to t_1^+} \int_t^\infty \tilde{g}_{j,x}\big(s,x_*(s)\big)\, d\mu_j(s) = 0$$
und es ist wieder $\dot{p}(t)=0$ f"ur $t \not \in [t_0,t_1]$.
Da die Adjungierte $p(\cdot)$ rechtsseitig stetig ist und die einseitigen Grenzwerte in $t=t_1$ existieren,
ergeben sich die Transversalitätsbedingung
$$p(t_1) = -{h_1'}^T\big(x_*(t_1)\big)l_1$$
und die Sprung-Transversalitätsbedingung
$$p(t_1) - p(t_1^-) = \sum_{j=1}^l \mu_j(\{t_1\})\, g_{j,x}\big(t_1,x_*(t_1)\big).$$
Zusammenfassend ergibt sich 
also das Pontrjaginsche Maximumprinzip \ref{SatzPMPZA} "uber dem endlichen Horizont f"ur die Aufgabe mit Zustandsbeschr"ankungen. \\[2mm]
Zusammenfassend k"onnen wir das Pontrjaginsche Maximumprinzip f"ur die Aufgaben "uber endlichem Zeithorizont aus den
Theoremen \ref{SatzPAUH} und \ref{SatzPAUHZA} vollst"andig ableiten.
Die Aufgabe mit unendlichem Zeithorizont und die vorgestellten Ergebnisse im Kapitel \ref{KapitelStrong} stellen somit
im Vergleich zum Kapitel \ref{KapitelStark} echte Verallgemeinerungen dar. \\
Insbesondere darf im Zielfunktional (\ref{PAUH1}) die Funktion $\omega(\cdot) \in L_1(\R_+,\R_+)$ lokal unbeschränkt sein.
Ein Vertreter einer lokal unbeschränkten Dichtefunktion ist die für die Finanz- und Ingenieurmathematik wichtige
Weibull-Verteilung mit einer Dichte der Gestalt $\omega(t)=t^{k-1}e^{-t^k}$ und einem Formparameter $k \in (0,1)$. 
Die Einbindung von lokal unbeschränkten Dichtefunktionen ist ein Detail,
das in der Literatur für die Standardaufgabe (\ref{PMP1})--(\ref{PMP5}) über endlichem Zeithorizont
bisher wenig Beachtung gefunden hat.

%% file: 5-0-Integralgleichungen.tex
\section{Steuerung Volterrascher Integralgleichungen} \label{KapitelIGL}
Die nach Vito Volterra (1860--1940) benannten Integralgleichungen treten in natürlicher Weise in dynamischen Problemen auf,
in denen das System eine Form von ``Erinnerungsvermögen'' besitzt.
Da sich im Rahmen der Standardaufgabe und deren Erweiterungen der Einfluss einer Steuervariable stets
unmittelbar auf den Zustand auswirkt,
kann ein Effekt, der sich im Laufe der Zeit entwickelt, oft nicht modelliert werden.
Hier ist die Beschreibung des dynamischen Systems mit Hilfe einer Integralgleichung ein probates Mittel. \\[2mm]
Ein Effekt über Zeit kann bei der Zusammenstellung der aggregierten Produktionskapazitäten vorliegen.
Die Kapazitäten ergeben sich als die gesamten Investitionen in Produktionsanlagen in den Jahrgängen $s \leq t$.
Durch Verschlei"s, Wartung oder technologischen Fortschritt wird die Produktionsfähigkeit zur Zeit $t$ der Anlagen des Jahrgangs $s \leq t$ durch
eine Funktion $\pi(t,s)$ beschrieben.
Zum Zeitpunkt $t$ ergeben sich demnach die gesamten vorliegenden Produktionskapazitäten durch 
$\displaystyle P(t)=\int_{t_0}^t \pi(t,s) \, ds$. \\
Im Rahmen der Werbeindustrie wird durch den strategischen Auf- und Ausbau der Produktbekanntheit dem Konsument ein gewisses
Produktimage und eine Verbundenheit zu dem Produkt suggeriert.
Dabei ist das Management mit der Problemstellung konfrontiert,
dass der Bekanntheitsgrad des Produktes durch Werbema"snahmen über eine gewisse Zeit aufgebaut werden muss.
Andererseits lässt die Produktbekanntheit (in Folge dessen auch die Nachfrage) umso stärker nach,
je länger der Werbeimpuls in der Vergangenheit liegt.
Die aggregierte Auswirkung der Werbekampagne zum einem Zeitpunkt ist demanch das Resultat der gesamten Werbeanstrengungen,
die im Vorfeld unternommen wurden. \\[2mm]
Optimalitätsbedingungen für die Steuerung von Integralgleichungen werden z.\,B. bei Bakke \cite{Bakke}, Bonnans \& De la Vega \cite{Bonnans},
Carlson \cite{Carlson}, Dmitruk \& Osmolovskii \cite{DmitrukOsmo1,DmitrukOsmo2} und Vinokurov \cite{Vinokurov} behandelt.
Im Gegensatz zu den aufgeführten Referenzen erweitern wir im vorliegenden Kapitel die Nadelvariationsmethode für die Standardaufgabe.
Der Ausbau der Methodik ist notwendig,
da bei der Steuerung von Integralgleichungen die Dynamik nicht mehr über einem Zeitintervall $\{t \in \R \,|\, t_0 \leq t \leq t_1\}$,
sondern über einem zweidimensionalen Zeitbereich $\{(s,t) \in \R^2 \,|\, t_0 \leq s, t \leq t_1\}$ agiert.
Auf diesem Merkmal liegt der Fokus bei der Untersuchung der Aufgabenklasse. \\[2mm]
Unsere Darstellungen beginnen wir mit der Herleitung des Maximumprinzips für die elementare Aufgabe mit freiem Endpunkt
auf Basis der einfachen Nadelvariation.
Anschlie"send untersuchen wir die allgemeinen Aufgabenklassen mit festen Randwerten und mit Zustandsbeschränkungen.
Au"serdem geben wir hinreichende Bedingungen nach Arrow an.
Als ein Anwendungsbeispiel für die Steuerung von Integralgleichungen betrachten wir in Anlehnung an Feichtinger \& Hartl \cite{Feichtinger}
und Kamien \& Schartz \cite{Kamien} die Bestimmung von optimalen Werbestrategien. 
Weitere ökonomische Problemstellungen sind bei Hritonenko \& Yatsenko \cite{Hritonenko} zu finden.
Abschließend behandeln wir die Aufgabe mit freiem Anfangs- und Endzeitpunkt.
Dabei halten wir uns eng an Dmitruk \& Osmolovskii \cite{DmitrukOsmo2}. 

%% file: 5-1-PMPeinfach.tex
\subsection{Die elementare Aufgabe mit freiem Endpunkt} \label{AbschnittPMPeinfachIGL}
Die Steuerung einer Volterraschen Integralgleichung bezeichnet die Aufgabe 
\begin{eqnarray}
&&\label{PMPeinfachIGL1} J\big(x(\cdot),u(\cdot)\big) = \int_{t_0}^{t_1} f\big(t,x(t),u(t)\big) \, dt \to \inf, \\
&&\label{PMPeinfachIGL2} x(t) = x_0 + \int_{t_0}^t \varphi\big(t,s,x(s),u(s)\big) \, ds, \quad t \in [t_0,t_1], \\
&&\label{PMPeinfachIGL3} u(t) \in U \subseteq \R^m, \quad U \not= \emptyset,
\end{eqnarray}
die wir bez"uglich $\big(x(\cdot),u(\cdot)\big) \in PC_1([t_0,t_1],\R^n) \times PC([t_0,t_1],U)$ untersuchen. \\[2mm]
Mit $\mathscr{D}^{\mathcal{I}}_{\rm Lip}$ bezeichnen wir die Menge $\big(x(\cdot),u(\cdot)\big)$,
für die es ein $\gamma>0$ derart gibt,
dass die Abbildungen $f(s,x,u)$, $\varphi(t,s,x,u)$ auf der Menge aller $(t,s,x,u) \in \R \times \R \times \R^n \times \R^m$ mit
$$t_0 \leq s, t \leq t_1, \qquad \|x-x(s)\| < \gamma, \qquad u \in \R^m$$
stetig in der Gesamtheit aller Variablen und stetig differenzierbar bezüglich $x$ sind. \\[2mm]
Das Paar $\big(x(\cdot),u(\cdot)\big) \in PC_1([t_0,t_1],\R^n) \times PC([t_0,t_1],U)$
hei"st ein zul"assiger Steuerungsprozess in der Aufgabe (\ref{PMPeinfachIGL1})--(\ref{PMPeinfachIGL3}),
falls $\big(x(\cdot),u(\cdot)\big)$ dem System (\ref{PMPeinfachIGL2}) zu $x(t_0)=x_0$ gen"ugt.
Mit $\mathscr{D}^{\mathcal{I}}_{\rm adm}$ bezeichnen wir die Menge der zul"assigen Steuerungsprozesse. \\[2mm]
Ein zul"assiger Steuerungsprozess $\big(x_*(\cdot),u_*(\cdot)\big)$ ist eine
starke lokale Minimalstelle\index{Minimum, starkes lokales!Integral@-- Integralgleichungen}
der Aufgabe (\ref{PMPeinfachIGL1})--(\ref{PMPeinfachIGL3}),
falls eine Zahl $\varepsilon > 0$ derart existiert, dass die Ungleichung 
$$J\big(x(\cdot),u(\cdot)\big) \geq J\big(x_*(\cdot),u_*(\cdot)\big)$$
f"ur alle $\big(x(\cdot),u(\cdot)\big) \in \mathscr{D}^{\mathcal{I}}_{\rm adm}$ mit $\|x(\cdot)-x_*(\cdot)\|_\infty < \varepsilon$ gilt. \\[2mm]
Es bezeichnet $H_0^{\mathcal{I}}: \R \times \R^n \times \R^m \times \R^n \times \R \to \R$ die ``einfache'' Pontrjagin-Funktion
$$H_0^{\mathcal{I}}(t,x,u,p(\cdot),\lambda_0) = -\int_t^{t_1} \langle \varphi(\tau,t,x,u), \dot{p}(\tau) \rangle \, dt -\lambda_0 f(t,x,u).$$

\begin{theorem}[Pontrjaginsches Maximumprinzip] \label{SatzPMPeinfachIGL}
\index{Pontrjaginsches Maximumprinzip!Integral@-- Integralgleichungen} 
Es sei $\big(x_*(\cdot),u_*(\cdot)\big) \in \mathscr{D}^{\mathcal{I}}_{\rm adm} \cap \mathscr{D}^{\mathcal{I}}_{\rm Lip}$. 
Ist $\big(x_*(\cdot),u_*(\cdot)\big)$ ein starkes lokales Minimum der Aufgabe (\ref{PMPeinfachIGL1})--(\ref{PMPeinfachIGL3}),
dann existiert eine Vektorfunktion $p(\cdot) \in PC_1([t_0,t_1],\R^n)$ derart, dass
\begin{enumerate}
\item[(a)] die adjungierte Gleichung
           \index{adjungierte Gleichung!Integral@-- Integralgleichungen}
           \begin{equation}\label{PMPeinfachIGL4} 
           \dot{p}(t) = -H_{0,x}^{\mathcal{I}}\big(t,x_*(t),u_*(t),p(\cdot),1\big),
           \end{equation}
\item[(b)] in $t=t_1$ die Transversalitätsbedingung
           \index{Transversalitätsbedingungen!Integral@-- Integralgleichungen}
           \begin{equation}\label{PMPeinfachIGL5} p(t_1)=0 \end{equation}
\item[(c)] und in fast allen Punkten $t \in [t_0,t_1]$ die Maximumbedingung
           \index{Maximumbedingung!Integral@-- Integralgleichungen}
           \begin{equation}\label{PMPeinfachIGL6} 
           H_0^{\mathcal{I}}\big(t,x_*(t),u_*(t),p(\cdot),1\big) = \max_{u \in U} H_0^{\mathcal{I}}\big(t,x_*(t),u,p(\cdot),1\big)
           \end{equation}
\end{enumerate}
erfüllt sind.
\end{theorem}



\newpage
{\bf Beweis} Da s"amtliche Abbildungen stetig und stetig differenzierbar bezüglich $x$ sind, und $u_*(\cdot)$ stückweise stetig ist,
existiert nach Lemma \ref{LemmaDGL2} eine eindeutige stückweise stetige L"osung $\dot{p}(\cdot)$ der adjungierten Gleichung (\ref{PMPeinfachIGL4}).
Angenommen, es ist $\lambda_0=0$.
Dann ergibt sich $\dot{p}(t_1)=0$ in (\ref{PMPeinfachIGL4}) und damit die eindeutige Lösung $\dot{p}(\cdot)=0$.
Wegen der Transversalitätsbedingung $p(t_1)=0$ würde au"serdem $p(\cdot)=0$ im Widerspruch zur
Nichttrivialität der Multiplikatoren $\big(\lambda_0,p(\cdot)\big)$ gelten.
Daher ist $\lambda_0>0$ und man darf $\lambda_0=1$ wählen. \\[1mm]
Es sei $\tau \in (t_0,t_1)$ ein Stetigkeitspunkt der Steuerung $u_*(\cdot)$.
Dann ist $u_*(\cdot)$ auch in einer gewissen hinreichend kleinen Umgebung von $\tau$ stetig und wir w"ahlen ein festes
$\lambda$ positiv und hinreichend klein, so dass sich $\tau-\lambda$ in dieser Umgebung befindet.
Weiter sei nun $v$ ein beliebiger Punkt aus $U$.
Wir setzen \index{Nadelvariation, einfache}
$$u(t;v,\tau,\lambda) = u_{\lambda}(t) = 
  \left\{ \begin{array}{ll}
          u_*(t) & \mbox{ f"ur } t \not\in [\tau-\lambda,\tau), \\
          v      & \mbox{ f"ur } t     \in [\tau-\lambda,\tau), 
          \end{array} \right.$$
und es bezeichne $x_\lambda(\cdot)$, $x_\lambda(t)=x(t;v,\tau,\lambda)$, die eindeutige L"osung der Gleichung
$$x(t) = x_0 + \int_{t_0}^t \varphi\big(t,s,x(s),u_\lambda(s)\big), \qquad t \in [t_0,t_1].$$
Dann ist $x_{\lambda}(t) = x_*(t)$ f"ur $t_0 \leq t \leq  \tau - \lambda$.
F"ur $t \geq \tau$ betrachten wir den Grenzwert
$$y(t)=\lim_{\lambda \to 0^+}\frac{x_{\lambda}(t) - x_*(t)}{\lambda}.$$
Nach den Sätzen \ref{SatzEElokal}--\ref{SatzDGLDifferenzierbarkeit} "uber die Stetigkeit und Differenzierbarkeit
der L"osung einer Integralgleichung in Abh"angigkeit von den Anfangsdaten ergibt sich für $y(t)$
im Grenzwert $\lambda \to 0^+$ die lineare Integralgleichung
\begin{eqnarray}
y(t) &=&  \varphi\big(t,\tau,x_*(\tau),v\big) - \varphi\big(t,\tau,x_*(\tau),u_*(\tau)\big)
          + \int_{\tau}^{t} \varphi_x\big(t,s,x_*(s),u_*(s)\big) y(s) \, ds \nonumber \\
     &=& \label{BeweisPMPeinfachIGL1} y(\tau;t) + \int_{\tau}^{t} \varphi_x\big(t,s,x_*(s),u_*(s)\big) y(s) \, ds.
\end{eqnarray}
Im Weiteren schreiben wir der Kürze halber
$$f_x[t]=f_x\big(t,x_*(t),u_*(t)\big), \qquad \varphi_x[t,s]=\varphi_x\big(t,s,x_*(s),u_*(s)\big).$$
Mit der adjungierten Gleichung (\ref{PMPeinfachIGL4}) erhalten wir f"ur $t \geq \tau$:
\begin{eqnarray}
    \int_{\tau}^{t_1} \langle f_x[t] , y(t) \rangle \, dt
&=& \int_{\tau}^{t_1} \Big\langle \dot{p}(t)- \int_t^{t_1} \varphi^T_x[s,t] \dot{p}(s) \, ds, y(t) \Big\rangle \, dt \nonumber \\
&=& \label{BeweisPMPeinfachIGL2}
    \int_{\tau}^{t_1} \langle \dot{p}(t) , y(t) \rangle \, dt - \int_{\tau}^{t_1} \Big\langle \int_t^{t_1} \varphi^T_x[s,t] \dot{p}(s) \, ds, y(t) \Big\rangle \, dt.
\end{eqnarray}
In der Differenz erhalten wir mit (\ref{BeweisPMPeinfachIGL1}) für den ersten Term 
\begin{eqnarray*}
    \int_{\tau}^{t_1} \langle \dot{p}(t) , y(t) \rangle \, dt
&=& \int_{\tau}^{t_1} \Big\langle \dot{p}(t), y(\tau;t) + \int_{\tau}^{t} \varphi_x[t,s]y(s) \, ds \Big\rangle \, dt \\
&=& \int_{\tau}^{t_1} \langle \dot{p}(t) , y(\tau;t) \rangle \, dt
    + \int_{\tau}^{t_1} \bigg( \int_{\tau}^{t} \langle \dot{p}(t), \varphi_x[t,s]y(s) \rangle \, ds \bigg) \, dt.
\end{eqnarray*}
Durch Vertauschen der Integrationsreihenfolge und der Variablenbezeichnung von $t$ und $s$ ergibt sich im zweiten Term
\begin{eqnarray*}
    \int_{\tau}^{t_1} \Big\langle \int_t^{t_1} \varphi^T_x[s,t] \dot{p}(s) \, ds, y(t) \Big\rangle \, dt
&=& \int_{\tau}^{t_1} \bigg( \int_{\tau}^s \langle \varphi^T_x[s,t] \dot{p}(s), y(t) \rangle \, dt \bigg) \, ds \\
&=& \int_{\tau}^{t_1} \bigg( \int_{\tau}^t \langle \dot{p}(t) , \varphi_x[t,s]  y(s) \rangle \, ds \bigg) \, dt.
\end{eqnarray*}
Damit lässt sich die rechte Seite in (\ref{BeweisPMPeinfachIGL2}) zusammenfassen und führt zu
\begin{equation} \label{BeweisPMPeinfachIGL3}
\int_{\tau}^{t_1} \langle f_x[t] , y(t) \rangle \, dt= \int_{\tau}^{t_1} \langle \dot{p}(t) , y(\tau;t) \rangle \, dt.
\end{equation}
Da $\big(x_*(\cdot),u_*(\cdot)\big)$ ein starkes lokales Minimum ist, gilt
\begin{eqnarray*}
0 &\leq& \lim_{\lambda \to 0^+} \frac{J\big(x_\lambda(\cdot),u_\lambda(\cdot)\big)- J\big(x_*(\cdot),u_*(\cdot)\big)}{\lambda} \\
  &=   & f\big(\tau,x_*(\tau),v\big) - f\big(\tau,x_*(\tau),u_*(\tau)\big)
         +\int_{\tau}^{t_1} \big\langle f_x\big(t,x_*(t),u_*(t)\big),y(t) \big\rangle \, dt.
\end{eqnarray*}
Bei Anwendung von Gleichung (\ref{BeweisPMPeinfachIGL3}) ergibt sich hieraus die Ungleichung
\begin{eqnarray*}
\lefteqn{-\int_{\tau}^{t_1} \big\langle \varphi\big(t,\tau,x_*(\tau),u_*(\tau)\big),\dot{p}(t) \big\rangle \, dt
                          - f\big(\tau,x_*(\tau),u_*(\tau)\big)}\\
&\geq& \bigg[-\int_{\tau}^{t_1} \big\langle \varphi\big(t,\tau,x_*(\tau),v\big),\dot{p}(t) \big\rangle \, dt
                          - f\big(\tau,x_*(\tau),v\big)\bigg].
\end{eqnarray*}                          
Damit folgt aus der Definition der Pontrjagin-Funktion mit $\lambda_0=1$:
\begin{equation*}
H_0^{\mathcal{I}}\big(\tau,x_*(\tau),u_*(\tau),p(\tau),1\big) \geq H_0^{\mathcal{I}}\big(\tau,x_*(\tau),v,p(\tau),1\big).
\end{equation*}
Nun ist $\tau$ ein beliebiger Stetigkeitspunkt von $u_*(\cdot)$ und $v$ ein beliebiger Punkt der Menge $U$. 
Demzufolge ist die Beziehung (\ref{PMPeinfachIGL6}) in allen Stetigkeitspunkten von $u_*(\cdot)$ wahr und damit ist das
Maximumprinzip bewiesen. \hfill $\blacksquare$

\begin{beispiel} \index{Werbestrategien} \label{AbschnittWerbung1}
{\rm Werbekampagnen zeichnen sich dadurch aus,
dass sie den Absatz nicht unmittelbar,
sondern mit einer zeitlichen Verzögerung beeinflussen.
In Anlehnung an Feichtinger \& Hartl \cite{Feichtinger} betrachten wir das Modell:
\begin{eqnarray}
&& \label{Werbung1} J\big(x(\cdot),u(\cdot)\big) = \int_0^T e^{-\varrho t} \big[\pi x(t)-u(t)\big] \, dt \to \sup, \\
&& \label{Werbung2} x(t)=x_0+\int_0^t e^{-\delta(t-s)} \big[f\big(u(s)\big)-\alpha x(s)\big] \, ds, \quad t \in [0,T], \\
&& \label{Werbung3} x_0 \in [0,1],\quad u(t) \geq 0, \quad \varrho >0, \quad \alpha, \delta \geq 0.
\end{eqnarray}
Dabei bezeichnet $x(\cdot)$ den Bekanntheitsgrad eines Produktes, der auf $x(t) \in [0,1]$ skaliert sei.
Der Einsatz von Werbung wird durch die Steuerung $u(\cdot)$ widergegeben.
Die Auswirkungen der Werbestrategie auf den Bekanntheitsgrad wird durch (\ref{Werbung2}) festgelegt.
Dabei sei $f$ über $[0,\infty)$ stetig, differenzierbar und streng monoton wachsend mit
$$f(0)=0, \quad \lim_{u \to \infty} f(u)=\gamma \leq \alpha, \quad \lim_{u \to 0+} f'(u)=\infty, \quad \lim_{u \to \infty} f'(u)=0.$$
Für $0< \sigma < 1$ genügen die Funktionen
$\displaystyle f_{\sigma}(u)=\gamma\frac{u^\sigma}{1+ u^\sigma}$
diesen Bedingungen. \\[2mm]
Die Differentation der Integralgleichung (\ref{Werbung2}) nach $t$ ergibt
\begin{eqnarray*}
\dot{x}(t) &=& -\delta \int_0^t e^{-\delta(t-s)} \big[f\big(u(s)\big)-\alpha x(s)\big] \, ds
               + \big[f\big(u(t)\big)-\alpha x(t)\big] \nonumber \\
           &=&  -\delta [x(t)-x_0]+ \big[f\big(u(t)\big)-\alpha x(t)\big]
               = -(\alpha+\delta) x(t) + \delta x_0 + f\big(u(t)\big)
\end{eqnarray*}
und die Aufgabe (\ref{Werbung1})--(\ref{Werbung3}) erhält die Gestalt der Standardaufgabe (\ref{PMP1})--(\ref{PMP4}):
\begin{eqnarray}
&& \label{Werbung4} J\big(x(\cdot),u(\cdot)\big) = \int_0^T e^{-\varrho t} \big[\pi x(t)-u(t)\big] \, dt \to \sup, \\
&& \label{Werbung5} \dot{x}(t) = -(\alpha+\delta) x(t) + \delta x_0 + f\big(u(t)\big), \quad x(0)=x_0, \\
&& \label{Werbung6} x_0 \in [0,1],\quad u(t) \geq 0, \quad \varrho >0, \quad \alpha, \delta \geq 0.
\end{eqnarray}
Im Fall $\delta =0$ bekommt die Differentialgleichung (\ref{Werbung5}) die Gestalt
$$\dot{x}(t)= f\big(u(t)\big)-\alpha x(t), \quad x(0)=x_0,$$
die die Form eines Ebbinghausenschen Vergessensmodell trägt.
Hermann Ebbinghausen (1850--1909) gilt als Begründer der experimentellen Gedächtnisforschung. \\[2mm]
Der Fall $u(s) \equiv 0$ überführt (\ref{Werbung5}) in die Gestalt $\dot{x}(t)=-(\alpha+\delta)x(t)+ \delta x_0$
zum Anfangswert $x(0)=x_0$ und besitzt die Lösung
\begin{equation} \label{Werbung7}
x(t)=\frac{\alpha}{\alpha+\delta} x_0e^{-(\alpha+\delta)t}+\frac{\delta}{\alpha+\delta} x_0.
\end{equation}
Die Popularität $x(t)$ ergibt sich anhand Gleichung (\ref{Werbung7}) durch das Zusammenspiel der Parameter
$\alpha$ und $\delta$.
Dabei spiegelt der Parameter $\alpha$ einem der Natur der Sache innewohnenden Attraktivitätsverlust für das Produkt wider,
während $\delta$ ein Stabilitätsmaß,
d.\,h. einen eigendynamischen Erhalt der Popularität innerhalb der Interessengruppe, darstellt.
Speziell ergeben sich  $x(t) \equiv x_0$ für $\alpha=0$ und es liegt kein Attraktivitätsverlust für das Produkt vor.
Im Fall $\delta =0$ stellt sich der Popularitätsverlust $x(t)=x_0e^{-\alpha t}$ zum Parameter $\alpha >0$ ein.
Ferner liegt für $\delta \to \infty$ eine stabile Faszination für das Produkt vor und es bleibt die Beliebtheit
$x(t) \equiv x_0$ unabhängig vom Parameter $\alpha \geq 0$ erhalten. \\[2mm]
Im Fall der Aufgabe (\ref{Werbung4})--(\ref{Werbung6}) liefert das Pontrjaginsche Maximumprinzip in Form von
Theorem \ref{SatzPMP} die adjungierte Gleichung
\begin{equation} \label{Werbung8}
\dot{p}(t)=(\alpha+\delta)p(t)-\pi e^{-\varrho t}, \quad p(T)=0,
\end{equation}
die die Lösung
\begin{equation} \label{Werbung9}
p(t)=-\frac{\pi \cdot e^{-(\alpha + \delta + \varrho)T}}{\alpha + \delta + \varrho} \cdot e^{(\alpha+\delta)t}
       +\frac{\pi}{\alpha + \delta + \varrho} \cdot e^{-\varrho t}
\end{equation}
besitzt.
Die Maximumbedingung führt nun weiter zur der Beziehung
\begin{equation} \label{Werbung10}
\max_{u \geq 0} \big[p(t)\cdot f(u)- e^{-\varrho t} \cdot u\big],
\end{equation}
aus der die Gültigkeit der Gleichung $p(t)\cdot f'\big(u_*(t)\big)-e^{-\varrho t} =0$ bzw.
\begin{equation} \label{Werbung1u}
f'\big(u_*(t)\big) = \frac{1}{p(t) \cdot e^{\varrho t}} =\frac{\alpha + \delta + \varrho}{\pi ( 1-e^{-(\alpha + \delta + \varrho)(T-t)})}=\psi(t)
\end{equation}
resultiert.
Die Funktion $\psi(t)$ ist über $[0,T)$ positiv und streng monton wachsend.
Im Grenzwert $t \to T^-$ ergibt sich $\psi(t) \to \infty$.
Aus den Eigenschaften von $f$ folgt,
dass $u_*(\cdot)$ über $[0,T]$ stetig, streng monoton fallend mit $u_*(t) \in [0,1)$ und $u_*(T)=0$ ist. \\[2mm]
Wir wenden uns nun der Aufgabe (\ref{Werbung1})--(\ref{Werbung3}) zu:
Das Pontrjaginsche Maximumprinzip in Form von Theorem \ref{SatzPMPeinfachIGL} liefert die adjungierte Gleichung
\begin{equation} \label{Werbung11}
\dot{q}(t)=\int_t^T -\alpha e^{-\delta(\tau-t)} \dot{q}(\tau) \, d\tau -\pi e^{-\varrho t}, \quad q(T)=0.
\end{equation}
Ferner führt die Maximumbedingung zu der Beziehung
\begin{equation} \label{Werbung12}
\max_{u \geq 0} \bigg[-\int_t^T e^{-\delta(\tau-t)} \dot{q}(\tau) \, d\tau \cdot f(u) - e^{-\varrho t} u\bigg].
\end{equation}
Darin besitzt $u_*(\cdot)$ wieder die implizite Darstellung $f'\big(u_*(t)\big) = 1/ [p(t) \cdot e^{\varrho t}]$. 
Der Vergleich von (\ref{Werbung10}) und (\ref{Werbung12}) ergibt den Zusammenhang
$\displaystyle p(t)=-\int_t^T e^{-\delta(\tau-t)} \dot{q}(\tau) \, d\tau$.
Durch Differentation dieser Gleichung nach $t$ erhalten wir weiter
$$\dot{p}(t) = -\delta \cdot \int_t^T e^{-\delta(\tau-t)} \dot{q}(\tau) \, d\tau + \dot{q}(t) = \delta p(t) + \dot{q}(t).$$
Es gilt damit $\dot{q}(t)=\dot{p}(t)-\delta p(t)=\alpha p(t)-\pi e^{-\varrho t}$.
Daher ist
$$\dot{q}(t)=\alpha p(t)-\pi e^{-\varrho t}=\int_t^T -\alpha e^{-\delta(\tau-t)} \dot{q}(\tau) \, d\tau -\pi e^{-\varrho t}$$
und $\dot{q}(\cdot)$ genügt der adjungierten Gleichung (\ref{Werbung11}).
Abschlie"send erhalten wir mit Hilfe von (\ref{Werbung9}) für $q(\cdot)$ die explizite Darstellung
$$\dot{q}(t) = -\frac{\alpha \pi \cdot e^{-(\alpha + \delta + \varrho)T}}{\alpha + \delta + \varrho} \cdot e^{(\alpha+\delta)t}
       -\frac{(\delta + \varrho)\pi}{\alpha + \delta + \varrho} \cdot e^{-\varrho t},
  \quad q(t)=\int_T^t \dot{q}(s)\,ds.$$
In dem Modell (\ref{Werbung4})--(\ref{Werbung6}) ist es vernünftig anzunehmen,
dass zum Abschluss der Werbekampagne das Produkt einen weitergehenden Absatz mit Gesamtumsatz $S\big(x(T)\big)$ in
Abhängigkeit vom Bekanntheitsgrad $x(T)$ besitzt.
Folglich ergibt sich das Zielfunktional
$$\tilde{J}\big(x(\cdot),u(\cdot)\big) = \int_0^T e^{-\varrho t} \big[\pi x(t)-u(t)\big] \, dt
                     + e^{-\varrho T} S\big(x(T)\big) \to \sup.$$
Die Abbildung $S$ sei stetig und stetig differenzierbar mit $S'(x)>0$ für alle $x\geq 0$. \\
Das Pontrjaginsche Maximumprinzip für die Aufgabe mit gemischtem Zielfunktional 
führt wieder
auf die adjungierte Gleichung (\ref{Werbung8}), jedoch zur Randbedingung
$$p(T)=e^{-\varrho T} S'\big(x_*(T)\big).$$
Daraus ergibt sich die Adjungierte
$$p(t)=\bigg(e^{-\varrho T} S'\big(x_*(T)\big)-\frac{\pi}{\alpha + \delta + \varrho}\bigg)
       \cdot e^{-(\alpha + \delta + \varrho)T}\cdot e^{(\alpha+\delta)t}
       +\frac{\pi}{\alpha + \delta + \varrho} \cdot e^{-\varrho t}.$$
Die Maximumbedingung führt erneut auf die Beziehung $f'\big(u_*(t)\big) = 1/ [p(t) \cdot e^{\varrho t}]$, d.\,h.
$$f'\big(u_*(t)\big)
  =\frac{1}{ \displaystyle \frac{\pi}{\alpha + \delta + \varrho}( 1-e^{-(\alpha + \delta + \varrho)(T-t)})
             +e^{-\varrho T} S'\big(x_*(T)\big)\cdot e^{-(\alpha + \delta + \varrho)(T-t)}}.$$
Die Werbestrategie $u_*(\cdot)$ ist erneut stetig, aber es gilt am Ende der Planungsperiode
$$f'\big(u_*(T)\big)
  =\frac{1}{e^{-\varrho T} S'\big(x_*(T)\big)}\in (0,\infty),$$
d.\,h. $u_*(T)>0$.
Demnach werden bis zum Abschluss der Periode Werbemaßnahmen getroffen. \hfill $\square$}
\end{beispiel}

%% file: 5-2-Aufgabenstellung.tex
\subsection{Die Aufgabenstellung} \label{AbschnittAufgabeIGL}
Wir betrachten als Steuerung einer Volterraschen Integralgleichung die Aufgabe
\begin{eqnarray}
&& \label{IG1} J\big(x(\cdot),u(\cdot)\big) = \int_{t_0}^{t_1} f\big(t,x(t),u(t)\big) \, dt \to \inf, \\
&& \label{IG2} x(t) = x(t_0)+ \int_{t_0}^t\varphi\big(t,s,x(s),u(s)\big) \, ds, \quad t \in [t_0,t_1], \\
&& \label{IG3} h_0\big(x(t_0)\big)=0, \qquad h_1\big(x(t_1)\big)=0, \\
&& \label{IG4} u(t) \in U \subseteq \R^m, \quad U\not= \emptyset, \\
&& \label{IG5} g_j\big(t,x(t)\big) \leq 0 \quad \mbox{f"ur alle } t \in [t_0,t_1], \quad j=1,...,l.
\end{eqnarray}
Darin gelten f"ur die Abbildungen
$f:\R \times \R^n \times \R^m \to \R$ und $\varphi:\R \times \R \times \R^n \times \R^m \to \R^n$,
sowie $h_i:\R^n \to \R^{s_i}$ für $i=0,1$ und $g_j:\R \times \R^n \to \R$ für $j=1,...,l$.
Wir betrachten die Aufgabe (\ref{IG1})--(\ref{IG5}) bez"uglich der Paare
$\big(x(\cdot),u(\cdot)\big) \in C([t_0,t_1],\R^n) \times L_\infty([t_0,t_1],U).$ \\[2mm]
Mit $\mathscr{B}^{\,\mathcal{I}}_{\rm Lip}$ bezeichnen wir die Menge aller Paare $\big(x(\cdot),u(\cdot)\big)$,
für die es ein $\gamma>0$ derart gibt, dass die Abbildungen
$f(s,x,u)$, $\varphi(t,s,x,u)$, $h_i(x_i)$ und $g_j(s,x)$ auf der Menge aller Punkte 
$(t,s,x,x_0,x_1) \in \R \times \R \times \R^n \times \R^n \times \R^n$ mit
$$t_0 \leq s, t\leq t_1, \quad \|x-x(s)\| < \gamma, \quad \|x_0-x(t_0)\| < \gamma, \quad \|x_1-x(t_1)\| < \gamma$$
stetig in der Gesamtheit der Variablen und stetig differenzierbar bez"uglich $x$ bzw. $x_i$ sind.
Daher sind für jede kompakte Menge $U_1 \subset \R^m$ die Abbildungen $\varphi$ und $\varphi_x$ gleichmäßig stetig auf
$V^{\mathcal{I}}_\gamma \times U_1$:
Zu jedem $\varepsilon >0$ gibt es ein $\delta >0$ mit
$$\|\varphi(t,s,x,u) - \varphi(t',s',x',u')\| < \varepsilon, \quad
  \|\varphi_x(t,s,x,u) - \varphi_x(t',s',x',u')\| < \varepsilon$$
für alle $(t,s,x,u),\, (t',s',x',u') \in V^{\mathcal{I}}_\gamma \times U_1$ mit
$\|(t,s,x,u) -(t',s',x',u')\| < \delta$. \\[2mm]
Das Paar $\big(x(\cdot),u(\cdot)\big) \in C([t_0,t_1],\R^n) \times L_\infty([t_0,t_1],U)$
hei"st ein zul"assiger Steuerungsprozess in der Aufgabe (\ref{IG1})--(\ref{IG5}),
falls $\big(x(\cdot),u(\cdot)\big)$ der Dynamik (\ref{IG2}) gen"ugt, sowie die Randbedingungen (\ref{IG3}) und
Zustandsbeschr"ankungen (\ref{IG5}) erf"ullt.
Die Menge $\mathscr{B}^{\,\mathcal{I}}_{\rm adm}$ bezeichnet die Menge der zul"assigen Steuerungsprozesse. \\[2mm]
Ein zul"assiger Steuerungsprozess $\big(x_*(\cdot),u_*(\cdot)\big)$ ist eine starke lokale
Minimalstelle\index{Minimum, starkes lokales!Integral@-- Integralgleichungen}
der Aufgabe (\ref{IG1})--(\ref{IG5}),
falls eine Zahl $\varepsilon > 0$ derart existiert, dass die Ungleichung 
$$J\big(x(\cdot),u(\cdot)\big) \geq J\big(x_*(\cdot),u_*(\cdot)\big)$$
f"ur alle $\big(x(\cdot),u(\cdot)\big) \in \mathscr{B}^{\,\mathcal{I}}_{\rm adm}$ mit 
$\|x(\cdot)-x_*(\cdot) \|_\infty < \varepsilon$ gilt. \\[2mm]
Zur Aufgabe (\ref{IG1})--(\ref{IG5}) definieren wir die Pontrjagin-Funktion wir folgt:
\begin{equation} \label{PontrjaginFunktionIG}
H^{\mathcal{I}}(t,x,u,p(\cdot),\lambda_0) = -\int_t^{t_1} \langle \varphi(\tau,t,x,u),dp(\tau) \rangle -\lambda_0 f(t,x,u).
\end{equation}
Dabei sind $(t,x,u,\lambda_0) \in [t_0,t_1] \times \R^n \times \R^m \times \R$
und $p(\cdot)$ eine rechtsseitig stetige Vektorfunktion von beschränkter Variation.

%% file: 5-3-PontrjaginAufgabe.tex
\subsection{Das Pontrjaginsche Maximumprinzip} \label{AbschnittPMPIGL}
\subsubsection{Notwendige Optimalit\"atsbedingungen}
\begin{theorem}[Pontrjaginsches Maximumprinzip] \label{SatzPMPIGS}
\index{Pontrjaginsches Maximumprinzip!Integral@-- Integralgleichungen} 
Es sei $\big(x_*(\cdot),u_*(\cdot)\big) \in \mathscr{B}^{\,\mathcal{I}}_{\rm adm} \cap \mathscr{B}^{\,\mathcal{I}}_{\rm Lip}$.
Ist $\big(x_*(\cdot),u_*(\cdot)\big)$ ein starkes lokales Minimum der Aufgabe (\ref{IG1})--(\ref{IG4}),
dann existieren nicht gleichzeitig verschwindende Multiplikatoren $\lambda_0 \geq 0$,
$p(\cdot) \in W^1_\infty([t_0,t_1),\R^n)$ und $l_i \in \R^{s_i}$, $i=0,1$, derart, dass
\begin{enumerate}
\item[(a)] die adjungierte Gleichung
           \index{adjungierte Gleichung!Integral@-- Integralgleichungen}
           \begin{equation} \label{SatzPMPIG1}
           \dot{p}(t) = -H^{\mathcal{I}}_x\big(t,x_*(t),u_*(t),p(\cdot),\lambda_0\big),
           \end{equation}
\item[(b)] die Transversalit"atsbedingungen
           \index{Transversalitätsbedingungen!Integral@-- Integralgleichungen}
           \begin{equation}\label{SatzPMPIG2}
           p(t_0)= {h_0'}^T\big(x_*(t_0)\big)l_0, \quad p(t_1^-)=-{h_1'}^T\big(x_*(t_1)\big)l_1, \quad p(t_1)=0
           \end{equation}
\item[(c)] und in fast allen Punkten $t \in [t_0,t_1]$ die Maximumbedingung           
           \index{Maximumbedingung!Integral@-- Integralgleichungen}
           \begin{equation}\label{SatzPMPIG3}
           H^{\mathcal{I}}\big(t,x_*(t),u_*(t),p(\cdot),\lambda_0\big) = \max_{u \in U} H^{\mathcal{I}}\big(t,x_*(t),u,p(\cdot),\lambda_0\big)
           \end{equation}
\end{enumerate}
erfüllt sind.
\end{theorem}

\begin{bemerkung} \label{BemerkungDeutungIGL}
{\rm \index{Pontrjaginsches Maximumprinzip!oekonomische@--, ökonomische Interpretation}
In Vergleich zu Abschnitt \ref{AbschnittDeutung},
wo $v \to H^{\mathcal{S}}\big(\tau,x_*(\tau),v,p(\tau),1\big)$ die Profitrate aus direktem $f\big(\tau,x_*(\tau),v\big)$
und indirektem Gewinn $\big\langle p(\tau) , \varphi\big(\tau,x_*(\tau),v\big) \big\rangle$ bemisst,
enthält $H^{\mathcal{I}}$ den Erwartungswert
$\displaystyle -\int_\tau^{t_1} \big\langle \varphi\big(t,\tau,x_*(\tau),v\big), dp(t) \big\rangle$
zum Gewicht $p(\cdot)$ über die gesamte künftige Entwicklung des indirekten Gewinnes.
\hfill $\square$}
\end{bemerkung}

\begin{bemerkung} \label{BemerkungIGLPontrjaginFunktion}
{\rm Die durch (\ref{SatzPMPIG2}) formulierte Sprungbedingung in $t=t_1$ ist wesentlicher Bestandteil der Transversalitätsbedingungen.
Aufgrund der absoluten Stetigkeit der Adjungierten $p(\cdot)$ kann die adjungierte Gleichung (\ref{SatzPMPIG1}) in die Form
\begin{eqnarray*}
\dot{p}(t) &=& \int_t^{t_1} \varphi_x^T\big(\tau,t,x_*(t),u_*(t)\big) \dot{p}(\tau) \, d\tau
                          + \varphi_x^T\big(t_1,t,x_*(t),u_*(t)\big) \big(p(t_1)-p(t_1^-) \big) \nonumber \\
           & & + \lambda_0 f_x\big(t,x_*(t),u_*(t)\big)
\end{eqnarray*}
zu den Transversalitätsbedingungen $p(t_1^-)=-{h_1'}^T\big(x_*(t_1)\big)l_1$ und $p(t_1)=0$ gebracht werden.
Trifft man an dieser Stelle die Vereinbarung, dass die Adjungierte $p(\cdot)$ in $t=t_1$ linksseitig stetig ist,
so entsteht die adjungierte Gleichung
\begin{eqnarray}
\dot{p}(t) &=& \int_t^{t_1} \varphi_x^T\big(\tau,t,x_*(t),u_*(t)\big) \dot{p}(\tau) \, d\tau
                          - \varphi_x^T\big(t_1,t,x_*(t),u_*(t)\big)p(t_1) \nonumber \\
\label{SatzPMPIG4} & & + \lambda_0 f_x\big(t,x_*(t),u_*(t)\big)
\end{eqnarray}
zur Transversalitätsbedingung $p(t_1)=-{h_1'}^T\big(x_*(t_1)\big)l_1$
als abgewandelte Darstellung der adjungierten Gleichung (\ref{SatzPMPIG1}).
\hfill $\square$}
\end{bemerkung}

\begin{bemerkung}
{\rm In der elementaren Aufgabe (\ref{IG1})--(\ref{IG4}) mit freiem rechten Endpunkt ist $p(t_1^-)=p(t_1)=0$.
Au"serdem muss in diesem Fall $\lambda_0>0$ gelten und es kann $\lambda_0=1$ gewählt werden.
Ferner erhält in diesem Fall die adjungierte Gleichung $(\ref{SatzPMPIG4})$ die Form
$$\dot{p}(t) = \int_t^{t_1} \varphi^T_x\big(\tau,t,x_*(t),u_*(t)\big) \dot{p}(\tau) \, d\tau + f_x\big(t,x_*(t),u_*(t)\big)$$
und weiterhin die Maximumbedingung $(\ref{SatzPMPIG3})$ die Gestalt
\begin{eqnarray*}
\lefteqn{-\int_t^{t_1} \big\langle \varphi\big(\tau,t,x_*(t),u_*(t)\big),\dot{p}(\tau) \big\rangle \, d\tau
                          - f\big(t,x_*(t),u_*(t)\big)} \\
&=& \max_{u \in U} \bigg[-\int_t^{t_1} \big\langle \varphi\big(\tau,t,x_*(t),u\big),\dot{p}(\tau) \big\rangle \, d\tau
                          - f\big(t,x_*(t),u\big)\bigg].
\end{eqnarray*} 
Zusammen ergeben sich mit der ``einfachen'' Pontrjagin-Funktion $H^{\mathcal{I}}_0$ die Bedingungen des Pontrjaginschen Maximumprinzips
in Form von Theorem \ref{SatzPMPeinfachIGL} in Abschnitt \ref{AbschnittPMPeinfachIGL}.
\hfill $\square$}
\end{bemerkung}

\begin{bemerkung}
{\rm In der Literatur, z.\,B. bei Feichtinger \& Hartl \cite{Feichtinger} oder Kamien \& Schwartz \cite{Kamien},
wird die adjungierte Gleichung (\ref{SatzPMPIG1}) häufig in der Form
$$\lambda(t) = \int_t^{t_1} \varphi_x^T\big(\tau,t,x_*(t),u_*(t)\big) \lambda(\tau) \, d\tau + f_x\big(t,x_*(t),u_*(t)\big)$$
angegeben.
In unseren Darstellungen korrespondiert dabei die Funktion $\lambda(\cdot)$ mit der Ableitung $\dot{p}(\cdot)$ der Adjungierten.
Bei Kamien \& Schwartz \cite{Kamien} wird nun keine Transversalitätsbedingung für $\lambda(\cdot)$ angegeben.
Jedoch stellen Feichtinger \& Hartl \cite{Feichtinger} bei einer Aufgabe mit dem gemischten Zielfunktional
$$\tilde{J}\big(x(\cdot),u(\cdot)\big) = \int_{t_0}^{t_1} f\big(t,x(t),u(t)\big) \, dt + S_1\big(x(t_1)\big) \to \inf$$
die Transversalitätsbedingung $\lambda(t_1)= - S_1'\big(x_*(t_1)\big)$ auf.
Man beachte, dass weiterhin $\lambda(\cdot)$ mit $\dot{p}(\cdot)$ korrespondiert.
Insbesondere bei fehlendem Terminalfunktional $S_1(x)$ würde diese Transversalitätsbedingung $\lambda(t_1)=\dot{p}(t_1)=0$ bedeuten.
Allerdings ergibt sich im Beispiel optimaler Werbestrategien im Abschnitt \ref{AbschnittWerbung1} die Gleichung (\ref{Werbung11}), d.\,h.
$$\dot{q}(t)=\int_t^T -\alpha e^{-\delta(\tau-t)} \dot{q}(\tau) \, d\tau -\pi e^{-\varrho t},$$
die mit $\lambda(\cdot)=\dot{q}(\cdot)$ die Gestalt
$$\lambda(t)=\int_t^T -\alpha e^{-\delta(\tau-t)} \lambda(\tau) \, d\tau -\pi e^{-\varrho t}$$
erhält und zu der Transversalitätsbedingung $\lambda(T)= -\pi e^{-\varrho T} \not=0$ führt.
Unsere Betrachtungen liefern jedoch neben der adjungierte Gleichung (\ref{SatzPMPIG1}) bezüglich $\dot{p}(\cdot)$
außerdem in der Aufgabe mit gemischtem Zielfunktional im Abschnitt \ref{AbschnittGemischtesZF} die Transversalitätsbedingung
$p(t_1) = -{h_1'}^T\big(x_*(t_1)\big)l_1 - \lambda_0 S_1'\big(x_*(t_1)\big)$ zum Terminalfunktional $S_1$.
\hfill $\square$}
\end{bemerkung}

%% file: 5-31-Beweis.tex
\subsubsection{Der Nachweis der notwendigen Optimalit\"atsbedingungen} \label{AbschnittBeweisPMPIG}
Bei der Herleitung der notwendigen Optimalitätsbedingungen ist bei der Konstruktion der mehrfachen Nadelvariation zur Dynamik
(\ref{IG1}) deren besondere Charakteristik zu beachten:
Im Vergleich zur Dynamik (\ref{PMP2}) in der Standardaufgabe (\ref{PMP1})--(\ref{PMP5}), die in Integraldarstellung die Gestalt
$$x(t)=x(t_0)+\int_{t_0}^t \varphi\big(s,x(s),u(s)\big) \, ds$$
besitzt,
liegt nun die Dynamik in der Form
$$x(t)=x(t_0)+\int_{t_0}^t \varphi\big(t,s,x(s),u(s)\big) \, ds$$
vor.
Die zusätzliche Variable $t$ hat dabei den Einfluss,
dass anstelle einer rechten Seite $s \to\varphi\big(s,x(s),u(s)\big)$ nun
im Rahmen der Steuerung einer Volterraschen Integralgleichung für jedes $t \in [t_0,t_1]$
die Abbildungsfamilie $s \to \big\{\varphi\big(\tau,s,x(s),u(s)\big) \,|\, \tau \in [t_0,t]\big\}$ vorliegt.
Deswegen sind bezüglich der Aufgabe (\ref{IG1})--(\ref{IG5}) mehrfache Nadelvariationen zu diesen Abbildungsfamilien einzuführen.
Daher müssen nun Nadelvariationen über dem Zeitbereich $\{(t,s) \in \R^2 \,|\, t_0 \leq s, t\leq t_1\}$ betrachtet werden.
Der Nachweis, dass erneut mehrfache Nadelvariationen mit den entsprechenden Eigenschaften existieren,
ist bei der Steuerung von Integralgleichungen die Herausforderung.

\begin{bemerkung}
{\rm Bei genauerem Hinsehen agiert die Volterrasche Integralgleichung über dem dreieckigen Zeitbereich $\{(t,s) \in \R^2 \,|\, t_0 \leq s \leq t\leq t_1\}$.
Die (bisher stillschweigende) Erweiterung auf den quadratischen Zeitbereich $\{(t,s) \in \R^2 \,|\, t_0 \leq s, t\leq t_1\}$ ist eine formale Sache,
die die Darstellungen im Anhang \ref{AnhangNVIG} erleichtern
und die Einführung der mehrfachen Nadelvariationen besser illustrieren lässt.
\hfill $\square$}
\end{bemerkung}

Wir betrachten f"ur
$\big(x(\cdot),u(\cdot)\big) \in C([t_0,t_1],\R^n) \times L_\infty([t_0,t_1],\R^m)$
die Abbildungen
\begin{eqnarray*}
J\big(x(\cdot),u(\cdot)\big) &=& \int_{t_0}^{t_1} f\big(s,x(s),u(s)\big) \, ds, \\
F\big(x(\cdot),u(\cdot)\big)(t) &=& x(t) -x(t_0) -\int_{t_0}^t \varphi\big(t,s,x(s),u(s)\big) \, ds, \quad t \in [t_0,t_1],\\
H_i\big(x(\cdot)\big) &=& h_i\big(x(t_i)\big), \quad i=0,1.
\end{eqnarray*}
F"ur diese Abbildungen gelten
\begin{eqnarray*}
J &:& C([t_0,t_1],\R^n) \times L_\infty([t_0,t_1],\R^m) \to \R, \\
F &:& C([t_0,t_1],\R^n) \times L_\infty([t_0,t_1],\R^m) \to C_0([t_0,t_1],\R^n), \\
H_i &:& C([t_0,t_1],\R^n) \to \R^{s_i}, \quad i=0,1.
\end{eqnarray*}
Wir setzen $\mathscr{F}=(F,H_0,H_1)$ und pr"ufen f"ur die Extremalaufgabe
\begin{equation} \label{ExtremalaufgabePMPIG}
J\big(x(\cdot),u(\cdot)\big) \to \inf, \quad \mathscr{F}\big(x(\cdot),u(\cdot)\big)=0,\quad u(\cdot) \in L_\infty([t_0,t_1],U)
\end{equation}
im Punkt $\big(x_*(\cdot),u_*(\cdot)\big) \in \mathscr{B}^{\,\mathcal{I}}_{\rm adm} \cap \mathscr{B}^{\,\mathcal{I}}_{\rm Lip}$
die Voraussetzungen von Theorem \ref{SatzExtremalprinzipStark}:

\begin{enumerate}
\item[(A)] Im Vergleich zur Aufgabe (\ref{PMP1})--(\ref{PMP5}) ist nur noch der Operator $F$ zu diskutieren.
           Da $x_*(\cdot) \in \mathscr{B}^{\,\mathcal{I}}_{\rm Lip}$ gilt,
           ist die Abbildung $\varphi(t,s,x,u)$ für jedes $u(\cdot) \in L_\infty([t_0,t_1],U)$ über einem Umgebungsstreifen von $x_*(\cdot)$
           gleichmäßig beschränkt.
           Daher bildet der Operator $F$ eine Umgebung des Punktes $x_*(\cdot)$ in den Raum $C_0([t_0,t_1],\R^n)$ ab. \\
           Beachten wir ferner die stetige Differenzierbarkeit von $\varphi$ bezüglich $x$,
           so gilt
           $$\frac{\varphi(t,s,x+\lambda y,u)-\varphi(t,s,x,u)}{\lambda}
              =\int_0^1 \varphi_x(t,s,x+\tau\lambda y,u) y \, d\tau.$$
           Die Differenzierbarkeit von $F$ bezüglich $x(\cdot)$ ergibt sich nun wie im Beispiel \ref{DiffDynamik2}.
\item[(B)] Für die Abbildung $A(t,s)=\varphi_x\big(t,s,x_*(s),u_*(s)\big)$
           wählen wir in Lemma \ref{LemmaDGL1} die integrierbare Funktion $c$ als das wesentliche Supremum von $A(t,s)$
           über der Menge $t_0 \leq s, t \leq t_1$.
           Die endliche Kodimension des Operators $\mathscr{F}_x\big(x_*(\cdot),u_*(\cdot)\big)$ folgt daher unmittelbar mit
           Lemma \ref{LemmaDGL1} und Bemerkung \ref{BemDGL}.
\item[(C)] Diese Voraussetzungen von Theorem \ref{SatzExtremalprinzipStark} beziehen sich auf die mehrfachen Nadelvariationen
           bei Integralgleichungen.
           Dass die entsprechenden Eigenschaften erfüllt sind, sind im Anhang \ref{AnhangNV}
           in Lemma \ref{LemmaEigenschaftNadelvariationIG1} und Lemma \ref{LemmaEigenschaftNadelvariationIG2} dargestellt.
\end{enumerate}

Zur Extremalaufgabe (\ref{ExtremalaufgabePMPIG}) definieren wir auf
$$C([t_0,t_1],\R^n)\times L_\infty([t_0,t_1],\R^m)\times \R\times C_0^*([t_0,t_1],\R^n)\times\R^{s_0}\times\R^{s_1}$$
die Lagrange-Funktion $\mathscr{L}=\mathscr{L}\big(x(\cdot),u(\cdot),\lambda_0,y^*,l_0,l_1\big)$ gemä"s
$$\mathscr{L} = \lambda_0 J\big(x(\cdot),u(\cdot)\big)+ \big\langle y^*, F\big(x(\cdot),u(\cdot)\big) \big\rangle
                +l_0^T H_0\big(x(\cdot)\big)+l_1^T H_1\big(x(\cdot)\big).$$
Ist $\big(x_*(\cdot),u_*(\cdot)\big)$ eine starke lokale Minimalstelle der Extremalaufgabe (\ref{ExtremalaufgabePMPIG}),
dann existieren nach Theorem \ref{SatzExtremalprinzipStark}
nicht gleichzeitig verschwindende Lagrangesche Multiplikatoren
$\lambda_0 \geq 0$, $y^* \in C_0^*([t_0,t_1],\R^n)$ und $l_i \in \R^{s_i}$ derart,
dass die folgenden notwendigen Bedingungen gelten:
\begin{enumerate}
\item[(a)] Die Lagrange-Funktion besitzt bez"uglich $x(\cdot)$ in $x_*(\cdot)$ einen station"aren Punkt, d.\,h.
           \begin{equation}\label{SatzPMPIGLMR1}
           0= \mathscr{L}_x\big(x_*(\cdot),u_*(\cdot),\lambda_0,y^*,l_0,l_1\big);
           \end{equation}         
\item[(b)] Die Lagrange-Funktion erf"ullt bez"uglich $u(\cdot)$ in $u_*(\cdot)$ die Bedingung
           \begin{equation}\label{SatzPMPIGLMR2}
           \mathscr{L}\big(x_*(\cdot),u_*(\cdot),\lambda_0,y^*,l_0,l_1\big)
           = \min_{u(\cdot) \in L_\infty ([t_0,t_1],U)}
             \mathscr{L}\big(x_*(\cdot),u(\cdot),\lambda_0,y^*,l_0,l_1\big).
           \end{equation}
\end{enumerate}
Aufgrund (\ref{SatzPMPIGLMR1}) ist folgende Variationsgleichung f"ur alle $x(\cdot) \in C([t_0,t_1],\R^n)$ erf"ullt: 
\begin{eqnarray}
0 &=& \lambda_0 \int_{t_0}^{t_1} \big\langle f_x\big(t,x_*(t),u_*(t)\big),x(t) \big\rangle\, dt 
      + \big\langle l_0, h_0'\big(x_*(t_0)\big) x(t_0) \big\rangle + \big\langle l_1, h_1'\big(x_*(t_1)\big) x(t_1) \big\rangle
      \nonumber \\
  & & \label{BeweisschlussPMPIG1}
      + \int_{t_0}^{t_1} \bigg[ x(t)-x(0) - \int_{t_0}^t \varphi_x\big(t,s,x_*(s),u_*(s)\big) x(s) \,ds \bigg]^T d\mu(t).
\end{eqnarray}
Hierin besitzt $\mu$ nach Folgerung \ref{FolgerungRiesz3} in $t=t_0$ kein Atom.
Wir "andern die Integrationsreihenfolge im letzten Summanden und bringen Gleichung (\ref{BeweisschlussPMPIG1}) in die Form
\begin{eqnarray}
0 = \int_{t_0}^{t_1} \Big\langle \lambda_0 f_x\big(t,x_*(t),u_*(t)\big)
              -  \int_{t}^{t_1} \varphi_x^T\big(\tau,t,x_*(t),u_*(t)\big) d\mu(\tau) , x(t) \Big\rangle \, dt \hspace*{15mm} \nonumber \\
 \label{BeweisschlussPMPIG2} + \int_{t_0}^{t_1} \langle x(t) , d\mu(t) \rangle
      + \Big\langle {h_0'}^T\big(x_*(t_0)\big)l_0 - \int_{t_0}^{t_1} d\mu(t) , x(t_0) \Big\rangle
      + \langle {h_1'}^T\big(x_*(t_1)\big)l_1 , x(t_1) \rangle.
\end{eqnarray}
Hierin setzen wir $\displaystyle p(t) = \int_t^{t_1} d\mu(\tau)$.
Ohne Einschränkung können wir annehmen,
dass das Ma"s $\mu$ über $\R$ regulär fortgesetzt, von beschränkter Totalvariation und auf $[t_0,t_1]$ konzentriert ist.
Andernfalls, wenn $\mu$ nicht auf $[t_0,t_1]$ konzentriert wäre,
so würde sich eine Funktion $\tilde{p}(\cdot)=p(\cdot)+c$ ergeben,
wobei die Verschiebung um die Konstante $c \in \R$ keinen Einfluss auf die nachstehende Argumentation besitzt.
Daher ist die Funktion $p(\cdot)$ von beschränkter Variation, rechtsseitig stetig und es gilt $p(t_1)=0$, da
$\mu$ auf $[t_0,t_1]$ konzentriert ist.
Damit kann $p(\cdot)$ in $t=t_1$ eine Sprungstelle besitzen. \\[2mm]
Die rechte Seite in (\ref{BeweisschlussPMPIG2}) definiert ein stetiges lineares Funktional im Raum $C([t_0,t_1],\R^n)$.
Nach dem Satz von Riesz ist die Darstellung eines stetigen linearen Funktionals eindeutig und 
es ergeben sich aus (\ref{BeweisschlussPMPIG2}) die Beziehungen $p(t_0) = {h_0'}^T\big(x_*(t_0)\big)l_0$ und
\begin{eqnarray}
p(t) \!\!&=&\!\!  -{h_1'}^T\big(x_*(t_1)\big)l_1
          - \int_t^{t_1} \!\!\bigg[\int_s^{t_1} \varphi^T_x\big(\tau,s,x_*(s),u_*(s)\big) dp(\tau) + \lambda_0 f_x\big(s,x_*(s),u_*(s)\big) \!\bigg] \, ds \nonumber \\
     \!\!&=&\!\! \label{BeweisschlussPMPIG3}  
         -{h_1'}^T\big(x_*(t_1)\big)l_1 + \int_t^{t_1} H^{\mathcal{I}}_x\big(s,x_*(s),u_*(s),p(\cdot),\lambda_0\big) \, ds.
\end{eqnarray}
Aufgrund (\ref{BeweisschlussPMPIG3}) ist $p(\cdot)$ über $[t_0,t_1)$ absolutstetig,
$\dot{p}(t) = -H^{\mathcal{I}}_x\big(t,x_*(t),u_*(t),p(\cdot),\lambda_0\big)$,
$p(t_1^-)=-{h_1'}^T\big(x_*(t_1)\big)l_1$ und $p(t_1)=0$.
Damit sind (\ref{SatzPMPIG1}) und (\ref{SatzPMPIG2}) gezeigt. \\[2mm]
Die Beziehung (\ref{SatzPMPIGLMR2}) führt nach Vertauschen der Integrationsreihenfolge und durch die Einführung der
Funktion $p(\cdot)$ zu der Ungleichung
\begin{eqnarray*}
\lefteqn{\int_{t_0}^{t_1} \bigg[\lambda_0 f\big(t,x_*(t),u_*(t)\big)
              +  \int_{t}^{t_1} \big\langle\varphi\big(\tau,t,x_*(t),u_*(t)\big) , dp(\tau) \big\rangle \bigg] \, dt} \\
  &\leq& \int_{t_0}^{t_1} \bigg[\lambda_0 f\big(s,x_*(s),u(s)\big)
              +  \int_{t}^{t_1} \big\langle\varphi\big(\tau,t,x_*(t),u(t)\big) , dp(\tau) \big\rangle \bigg] \, dt,
\end{eqnarray*}
die f"ur alle $u(\cdot) \in L_\infty([t_0,t_1],U)$ gültig ist.
In dieser Ungleichung betrachten wir nun Steuerungen der Form
$$u_\lambda(t)=u_*(t)+ \chi_{[s,s+\lambda]}(t) \cdot \big(u-u_*(t)\big)$$
mit $s \in (t_0,t_1)$, $0<\lambda <t_1-s$ und $u \in U$.
Für $u_\lambda(\cdot)$ erhält die Ungleichung die Form
\begin{eqnarray}
\lefteqn{\int_{t_0}^{t_1} \bigg[ \chi_{[s,s+\lambda]}(t) \cdot
                          \lambda_0 \big[ f\big(t,x_*(t),u_*(t)\big)-f\big(t,x_*(t),u\big)\big]} \nonumber \\
&& \hspace*{10mm} + \int_{t}^{t_1} \chi_{[s,s+\lambda]}(t) \cdot
    \big\langle\varphi\big(\tau,t,x_*(t),u_*(t)\big)-\varphi\big(\tau,t,x_*(t),u\big), dp(\tau) \big\rangle \bigg] \, dt
   \nonumber \\
&=& \int_{t_0}^{t_1} \chi_{[s,s+\lambda]}(t) \cdot \bigg[ 
                          \lambda_0 \big[ f\big(t,x_*(t),u_*(t)\big)-f\big(t,x_*(t),u\big)\big] \nonumber \\
&& \label{BeweisschlussPMPIG5} \hspace*{15mm}
   + \int_{t}^{t_1} \big\langle\varphi\big(\tau,t,x_*(t),u_*(t)\big)-\varphi\big(\tau,t,x_*(t),u\big), dp(\tau) \big\rangle
   \bigg] \, dt \leq 0.
\end{eqnarray}
Zu $u \in U$ bzw. zu $u_*(\cdot)$ führen wir die Abbildungen
$$[\Phi(u)](t)=\int_{t}^{t_1} \big\langle \varphi\big(\tau,t,x_*(t),u\big), dp(\tau) \big\rangle, \quad
  [\Phi(u_*)](t)=\int_{t}^{t_1} \big\langle \varphi\big(\tau,t,x_*(t),u_*(t)\big), dp(\tau) \big\rangle$$
ein.
Es sind $(\tau,t) \to \varphi\big(\tau,t,x_*(t),u\big)$ für jedes $u \in U$ stetig, $x_*(\cdot)$ stetig und $p(\cdot)$ über $[t_0,t_1)$ absolutstetig.
Daher ist die Abbildung $t \to [\Phi(u)](t)$ für jedes $u \in U$ stetig über $(t_0,t_1)$.
Ganz entsprechend ist die Abbildung $t \to f\big(t,x_*(t),u\big)$ für jedes $u \in U$ über $(t_0,t_1)$ stetig.
Deswegen besitzen die Abbildungen  $t \to [\Phi(u)](t)$ und $t \to f\big(t,x_*(t),u\big)$  in jeder Stelle $t \in (t_0,t_1)$ einen Lebesgueschen Punkt. \\
Es sei $t \in (t_0,t_1)$ ein Lebesguescher Punkt der messbar und beschränkten Abbildungen
$$t \to f\big(t,x_*(t),u_*(t)\big), \qquad t \to [\Phi(u_*)](t).$$
So ergibt sich aus (\ref{BeweisschlussPMPIG5}) die Beziehung
$$\frac{1}{\lambda} \int_{t}^{t+\lambda} \Big[ 
                          \lambda_0 \big( f\big(\tau,x_*(\tau),u_*(\tau)\big)-f\big(\tau,x_*(\tau),u\big)\big)
  + \big([\Phi(u_*)](t) - [\Phi(u)](t)\big) \Big] \, dt \leq 0.$$
Aus der Eigenschaft eines Lebesgueschen Punktes folgt die Gültigkeit von
$$\lambda_0 \big( f\big(t,x_*(t),u_*(t)\big)-f\big(t,x_*(t),u\big)\big) + \big([\Phi(u_*)](t) - [\Phi(u)](t)\big) \leq 0$$
bzw. in äquivalenter Form mit Hilfe der Pontrjagin-Funktion (\ref{PontrjaginFunktionIG}) die Gültigkeit von
$$H^{\mathcal{I}}\big(t,x_*(t),u_*(t),p(\cdot),\lambda_0\big) \geq H^{\mathcal{I}}\big(t,x_*(t),u,p(\cdot),\lambda_0\big)$$
f"ur alle $u \in U$.
Da unabhängig von $u \in U$ fast überall in $(t_0,t_1)$ ein Lebesguescher Punkt vorliegt,
gilt fast überall über $[t_0,t_1]$ die Maximumbedingung (\ref{SatzPMPIG3}). \hfill $\blacksquare$

%% file: 5-32-Hinreichend.tex
\subsubsection{Hinreichende Bedingungen nach Arrow} \label{AbschnittArrowIGL}
Wir wenden uns nun hinreichenden Bedingungen der Aufgabe (\ref{IG1})--(\ref{IG4}) mit
freien bzw. festen Randwerten $x(t_0)=x_0$ und $x(t_1)=x_1$
nach Arrow\index{hinreichende Bedingungen, Arrow!Integral@-- Integralgleichungen} zu: 
Die Hamilton-Funktion $\mathscr{H}^{\mathcal{I}}$ besitzt für diese Aufgabe in der normalen Form die Gestalt
$$\mathscr{H}^{\mathcal{I}}\big(t,x,p(\cdot)\big)
  = \sup_{u \in U} \bigg[-\int_t^{t_1} \langle \varphi(\tau,t,x,u),dp(\tau) \rangle - f(t,x,u)\bigg].$$
Außerdem führen wir die Menge $V^{\mathcal{S}}_\gamma(t)=\{ x \in \R^n \,|\, \|x-x_*(t)\| < \gamma\}$ ein.

\begin{theorem} \label{SatzHBPMPIG}
In der Aufgabe (\ref{IG1})--(\ref{IG4})
sei $\big(x_*(\cdot),u_*(\cdot)\big) \in  \mathscr{B}^{\,\mathcal{I}}_{\rm adm} \cap \mathscr{B}^{\,\mathcal{I}}_{\rm Lip}$.
Es gelte:
\begin{enumerate}
\item[(a)] Das Tripel $\big(x_*(\cdot),u_*(\cdot),p(\cdot)\big)$
           erf"ullt (\ref{SatzPMPIG1})--(\ref{SatzPMPIG3}) in Theorem \ref{SatzPMPIGS} mit $\lambda_0=1$.        
\item[(b)] F"ur jedes $t \in [t_0,t_1]$ ist die Funktion $\mathscr{H}^{\mathcal{I}}\big(t,x,p(\cdot)\big)$ konkav in $x$ auf $V^{\mathcal{S}}_\gamma(t)$.
\end{enumerate}
Dann ist $\big(x_*(\cdot),u_*(\cdot)\big)$ ein starkes lokales Minimum der Aufgabe (\ref{IG1})--(\ref{IG4}).
\end{theorem}

{\bf Beweis}  Es sei $t \in [t_0,t_1]$ so gew"ahlt,
dass die Maximumbedingung (\ref{SatzPMPIG3}) zu diesem Zeitpunkt erf"ullt ist.
Mit den gleichen Argumenten wie im Abschnitt \ref{AbschnittHBPMP} folgt aus der Konkavität der
Hamilton-Funktion $\mathscr{H}^{\mathcal{I}}$ die Gültigkeit der Ungleichung
\begin{eqnarray*}
      -\langle a(t),x-x_*(t)\rangle
&\geq& \mathscr{H}^{\mathcal{I}}\big(t,x,p(\cdot)\big) - \mathscr{H}^{\mathcal{I}}\big(t,x_*(t),p(\cdot)\big) \\
&\geq& \int_t^{t_1} \langle \varphi\big(\tau,t,x_*(t),u_*(t)\big),dp(\tau) \rangle + f\big(t,x_*(t),u_*(t)\big) \\
&    & - \int_t^{t_1} \langle \varphi\big(\tau,t,x,u_*(t)\big),dp(\tau) \rangle - f\big(t,x,u_*(t)\big)
\end{eqnarray*}
f"ur alle $x \in V^{\mathcal{S}}_\gamma(t)$,
und es führt weiterhin die Funktion
\begin{eqnarray*}
\Phi(x) &=& \int_t^{t_1} \langle \varphi\big(\tau,t,x_*(t),u_*(t)\big),dp(\tau) \rangle + f\big(t,x_*(t),u_*(t)\big) \\
        & & -\int_t^{t_1} \langle \varphi\big(\tau,t,x,u_*(t)\big),dp(\tau) \rangle - f\big(t,x,u_*(t)\big) + \langle a(t),x-x_*(t)\rangle
\end{eqnarray*}
zu der Beziehung
$$0=\Phi'(x_*(t)) = - \int_t^{t_1} \langle \varphi_x\big(\tau,t,x_*(t),u_*(t)\big),dp(\tau) \rangle - f_x\big(t,x_*(t),u_*(t)\big) + a(t).$$
Daher gilt f"ur fast alle $t \in [t_0,t_1]$ über der Menge $V^{\mathcal{S}}_\gamma(t)$ nach (\ref{SatzPMPIG1}) die Ungleichung
\begin{equation} \label{HBPMPIG1}
\langle \dot{p}(t),x-x_*(t)\rangle \leq \mathscr{H}^{\mathcal{I}}\big(t,x_*(t),p(\cdot)\big)- \mathscr{H}^{\mathcal{I}}\big(t,x,p(\cdot)\big).
\end{equation}
Ferner ist für einen zulässigen Steuerungsprozess $\big(x(\cdot),u(\cdot)\big)$ nach (\ref{PMPeinfachIGL2}):
\begin{equation} \label{HBPMPIG2}
x(t)-x_*(t) = x(t_0)-x_*(t_0) +\int_{t_0}^t \Big[\varphi\big(t,s,x(s),u(s)\big)-\varphi\big(t,s,x_*(s),u_*(s)\big)\Big] \, ds.
\end{equation}
Es sei $\big(x(\cdot),u(\cdot)\big) \in  \mathscr{B}^{\,\mathcal{I}}_{\rm adm}$ mit $\|x(\cdot)-x_*(\cdot)\|_\infty < \gamma$.
Dann gilt
\begin{eqnarray*}
\Delta &=& J\big(x(\cdot),u(\cdot)\big)- J\big(x_*(\cdot),u_*(\cdot)\big)
\geq \int_{t_0}^{t_1} \big[\mathscr{H}^{\mathcal{I}}\big(t,x_*(t),p(\cdot)\big)-\mathscr{H}^{\mathcal{I}}\big(t,x(t),p(\cdot)\big)\big] \, dt \\
& & \hspace*{15mm} + \int_{t_0}^{t_1} \bigg[\int_t^{t_1}
    \big\langle \varphi\big(\tau,t,x_*(t),u_*(t)\big) - \varphi\big(\tau,t,x(t),u(t)\big),dp(\tau) \big\rangle \bigg]\, dt.
\end{eqnarray*}
Wir ändern im letzten Ausdruck die Integrationsreihenfolge,
\begin{eqnarray*}
\lefteqn{\int_{t_0}^{t_1} \bigg[\int_t^{t_1}
     \big\langle \varphi\big(\tau,t,x_*(t),u_*(t)\big) - \varphi\big(\tau,t,x(t),u(t)\big),dp(\tau) \big\rangle
     \bigg]\, dt} \\
&=& \int_{t_0}^{t_1} \bigg[\int_{t_0}^t \big[ \varphi\big(t,s,x_*(s),u_*(s)\big) - \varphi\big(t,s,x(s),u(s)\big) \big] \, ds \bigg]^T \, dp(t),
\end{eqnarray*}
und erhalten unter Verwendung von (\ref{HBPMPIG1}) und (\ref{HBPMPIG2}):
\begin{equation} \label{HBPMPIG3}
\Delta \geq \int_{t_0}^{t_1} \langle \dot{p}(t),x(t)-x_*(t)\rangle \, dt + \int_{t_0}^{t_1} \langle x_*(t)-x(t)-x_*(t_0)+x(t_0), dp(t) \rangle.
\end{equation}
Nach Theorem \ref{SatzPMPIGS} ist $p(\cdot)$ absolutstetig und erfüllt die Transversalitätsbedingungen (\ref{SatzPMPIG2});
insbesondere eine Sprungbedingung in $t=t_1$.
Demnach ergeben sich
\begin{eqnarray*}
    \int_{t_0}^{t_1} \langle x_*(t)-x(t), dp(t) \rangle
&=& \int_{t_0}^{t_1} \langle \dot{p}(t),x_*(t)-x(t)\rangle \, dt + \langle p(t_1)- p(t_1^-),x_*(t_1)-x(t_1) \rangle, \\
\int_{t_0}^{t_1} \langle x(t_0) - x_*(t_0), dp(t) \rangle
&=& \langle [p(t_1^-)-p(t_0)]+[p(t_1)-p(t_1^-)],x(t_0)-x_*(t_0) \rangle.
\end{eqnarray*}
Unter Beachtung dieser Beziehungen und von $p(t_1)=0$ erhalten wir abschlie"send
$$\Delta \geq \langle p(t_1^-),x(t_1) - x_*(t_1)\rangle - \langle p(t_0),x(t_0)-x_*(t_0) \rangle = 0.$$
Somit ist $\big(x_*(\cdot),u_*(\cdot)\big)$ ein starkes lokales Minimum der Aufgabe (\ref{IG1})--(\ref{IG4}). \hfill $\blacksquare$

\begin{beispiel} {\rm  \index{Werbestrategien}
In der Aufgabe (\ref{Werbung1})--(\ref{Werbung3}) ist der rechte Endpunkt frei und es flie"st die Zustandsvariable $x$ linear ein.
Somit ist die Hamilton-Funktion konkav in $x$ und die Bedingungen des Pontrjaginschen Maximumprinzips \ref{SatzPMPeinfachIGL}
sind nach den Arrow-Bedingungen in Theorem \ref{SatzHBPMPIG} hinreichend.
Also ist die in der Aufgabe (\ref{Werbung1})--(\ref{Werbung3}) ermittelte Werbestrategie, die durch die implizite Beziehung
$$f'\big(u_*(t)\big) = \frac{\alpha + \delta + \varrho}{\pi ( 1-e^{-(\alpha + \delta + \varrho)(T-t)})}$$
bestimmt ist, tatsächlich optimal.  \hfill $\square$}
\end{beispiel}

%% file: 5-33-Werbestrategien.tex
\subsubsection{Optimale Werbestrategien bei festem Werbebudget} \index{Werbestrategien} \label{AbschnittWerbung2}
Wir betrachten erneut das Beispiel (\ref{Werbung1})--(\ref{Werbung3}) aus Abschnitt \ref{AbschnittWerbung1}:
\begin{eqnarray}
&& \label{Werbung21} J\big(x(\cdot),u(\cdot)\big) = \int_0^T e^{-\varrho t} \big[\pi x(t)-u(t)\big] \, dt \to \sup, \\
&& \label{Werbung22} x(t)=x_0+\int_0^t e^{-\delta(t-s)} \big[f\big(u(s)\big)-\alpha x(s)\big] \, ds, \quad t \in [0,T], \\
&& \label{Werbung23} x_0 \in [0,1],\quad u(t) \geq 0, \quad \varrho >0, \quad \alpha, \delta \geq 0.
\end{eqnarray}
In dieser Aufgabe ergaben sich mit den in  Abschnitt \ref{AbschnittWerbung1} verwendeten Bezeichnungen:
\begin{enumerate}
\item[1.)] Die Adjungierte $q(\cdot)$ ist Lösung der adjungierten Gleichung
           $$\dot{q}(t)=\int_t^T -\alpha e^{-\delta(\tau-t)} \dot{q}(\tau) \, d\tau -\pi e^{-\varrho t}, \quad q(T^-)=q(T)=0,$$
           wobei $\dot{q}(\cdot)$ der Beziehung $\dot{q}(t)=\dot{p}(t)-\delta p(t)=\alpha p(t)-\pi e^{-\varrho t}$ genügt.
           Dabei ist
           $$p(t)=-\frac{\pi \cdot e^{-(\alpha + \delta + \varrho)T}}{\alpha + \delta + \varrho} \cdot e^{(\alpha+\delta)t}
                  +\frac{\pi}{\alpha + \delta + \varrho} \cdot e^{-\varrho t}$$
           die Adjungierte des Problems (\ref{Werbung4})--(\ref{Werbung6}) und erfüllt die adjungierte Gleichung 
           $$\dot{p}(t)=(\alpha+\delta)p(t)-\pi e^{-\varrho t}$$
           zur Transversalitätsbedingung $p(T)=0$.         
\item[2.)] Die implizit definierte, optimale Werbestrategie genügt nach (\ref{Werbung1u}) der Beziehung
           \begin{equation} \label{Werbung2u}
           f'\big(u_*(t)\big) = \frac{1}{p(t) \cdot e^{\varrho t}} =\frac{\alpha + \delta + \varrho}{\pi ( 1-e^{-(\alpha + \delta + \varrho)(T-t)})}.
           \end{equation}
\end{enumerate}
Wir erweitern die Aufgabe (\ref{Werbung21})--(\ref{Werbung23}) um die Beschränkung
\begin{equation} \label{Werbung24} 
z(t)=\int_0^t e^{-\varrho s} u(s) \, ds, \quad t \in [0,T], \quad z(T)=Z>0,
\end{equation}
die den vollständigen Einsatz eines fest vorgegebenen Werbebudgets $Z$ fordert. 
Ferner bezeichne $Z_0$ das Werbebudget zu der in (\ref{Werbung2u}) definierten Steuerung,
die in der Aufgabe (\ref{Werbung1})--(\ref{Werbung3}) optimal ist.  \\[2mm]
Bezüglich dem Zustand $z(\cdot)$ gehen wir zur Differentialgleichung $\dot{z}(t)=e^{-\varrho t}u(t)$ über und erhalten zur
Aufgabe (\ref{Werbung21})--(\ref{Werbung24}) die Pontrjagin-Funktion
$$H^{\mathcal{I}}(t,x,z,u,q(\cdot),r,\lambda_0)
  = -\int_t^T e^{-\delta(\tau-t)} \big[f\big(u(t)\big)-\alpha x(t)\big] dq(\tau) + r[ e^{-\varrho t} u] + \lambda_0 e^{-\varrho t} [\pi x-u].$$
\newpage
Die notwendigen Optimalitätsbedingungen in Theorem \ref{SatzPMPIGS} ergeben mit $\lambda_0=1$:
\begin{enumerate}
\item[$\cdot$] die adjungierten Gleichungen und Transversalitätsbedingungen
           \begin{eqnarray}
            \label{Werbung25} -H^{\mathcal{I}}_x\big(t,...,1\big) &=& \dot{q}(t)
                              = \int_t^{t_1} -\alpha e^{-\delta(\tau-t)}  dq(\tau) - \pi e^{-\varrho t}, \quad q(T)=0, \quad \\            
            \label{Werbung26} -H^{\mathcal{I}}_z\big(t,...,1\big) &=& \dot{r}(t) =0, \quad r(T)=-l_1;
           \end{eqnarray}           
\item[$\cdot$] in fast allen Punkten $t \in [t_0,t_1]$ die Maximumbedingung
           \begin{equation} \label{Werbung27} 
           H^{\mathcal{I}}\big(t,x_*(t),z_*(t),u_*(t),q(\cdot),r(\cdot),1\big)
           = \max_{u \in U} H^{\mathcal{I}}\big(t,x_*(t),z_*(t),u,q(\cdot),r(\cdot),1\big).
           \end{equation}
\end{enumerate}
Nach Theorem \ref{SatzHBPMPIG} sind die Bedingungen (\ref{Werbung25})--(\ref{Werbung27}) hinreichend. \\[2mm]
Die Lösung $q(t)$ der ersten Gleichung ist unter dem Punkt 1.) bereits mit Hilfe der Funktion $p(\cdot)$ angegeben,
genauer ist
$$\dot{q}(t) = \alpha p(t)-\pi e^{-\varrho t}
  = -\frac{\alpha \pi \cdot e^{-(\alpha + \delta + \varrho)T}}{\alpha + \delta + \varrho} \cdot e^{(\alpha+\delta)t}
       -\frac{(\delta + \varrho)\pi}{\alpha + \delta + \varrho} \cdot e^{-\varrho t}.$$
Die Lösung der zweiten Gleichung ist $r(t) \equiv -l_1$. \\[1mm]
Die Maximumbedingung (\ref{Werbung27}) ist äquivalent zu der Beziehung
$$\max_{u \geq 0} \bigg[-\int_t^T e^{-\delta(\tau-t)} \dot{q}(\tau) \, d\tau \cdot f(u) - e^{-\varrho t}l_1 \cdot u- e^{-\varrho t} \cdot u\bigg].$$
Anstelle von (\ref{Werbung1u}) ergibt sich daraus mit Hilfe der Funktion $p(\cdot)$ in Punkt 1.) nun
$$f'\big(u_*(t)\big) = (1-r_0) \cdot \frac{1}{p(t) \cdot e^{\varrho t}}
  = (1+l_1) \cdot \frac{\alpha + \delta + \varrho}{\pi ( 1-e^{-(\alpha + \delta + \varrho)(T-t)})}.$$
Es lassen sich daraus folgende Aussagen ableiten:
\begin{enumerate}
\item[$\cdot$] Aus den Eigenschaften der Funktion $f$ folgt $l_1>-1$.
\item[$\cdot$] Der Wert für den Werbeeinsatz $u_*(t)>0$ fällt umso kleiner aus, desto größer der Parameter $l_1 >-1$ ist.
\item[$\cdot$] Für $l_1=0$ stimmt $u_*(t)$ mit der optimalen Werbestrategie der Aufgabe (\ref{Werbung1})--(\ref{Werbung3}) überein.
               In diesem Fall bezeichnet $Z_0$ die gesamten Werbeausgaben.
\item[$\cdot$] Bei einem reduzierten Werbebudget $0 < Z < Z_0$ führt die Budgetbeschränkung auf $l_1 >0$.
               Bei einem überdimensioniertem Werbebudget $Z > Z_0$ ist $-1< l_1 <0$. \hfill $\square$
\end{enumerate}

%% file: 5-4-Zustandsaufgabe.tex
\subsection{Aufgaben mit Zustandsbeschr\"ankungen}
\subsubsection{Notwendige Optimalit\"atsbedingungen} 
\begin{theorem}[Pontrjaginsches Maximumprinzip] \label{SatzPMPIGZB}
\index{Pontrjaginsches Maximumprinzip!Integral@-- Integralgleichungen} 
Es sei $\big(x_*(\cdot),u_*(\cdot)\big) \in \mathscr{B}^{\,\mathcal{I}}_{\rm adm} \cap \mathscr{B}^{\,\mathcal{I}}_{\rm Lip}$.
Ist $\big(x_*(\cdot),u_*(\cdot)\big)$ ein starkes lokales Minimum der Aufgabe (\ref{IG1})--(\ref{IG5}),
dann existieren eine Zahl $\lambda_0 \geq 0$,
Vektoren $l_0 \in \R^{s_0}$ und $l_1 \in \R^{s_1}$,
eine Vektorfunktion $p(\cdot):[t_0,t_1] \to \R^n$
und auf den Mengen
$$T_j=\big\{t \in [t_0,t_1] \,\big|\, g_j\big(t,x_*(t)\big)=0\big\}, \quad j=1,...,l,$$
konzentrierte nichtnegative regul"are Borelsche Ma"se $\mu_j$ endlicher Totalvariation
(wobei s"amtliche Gr"o"sen nicht gleichzeitig verschwinden) derart,
dass die Vektorfunktion $p(\cdot)$ von beschr"ankter Variation und rechtsseitig stetig ist, und
\begin{enumerate}
\item[(a)] die adjungierte Gleichung
           \index{adjungierte Gleichung!Integral@-- Integralgleichungen}
           \begin{eqnarray}
           p(t) &=& -{h_1'}^T\big(x_*(t_1)\big)l_1 + \int_t^{t_1} H^{\mathcal{I}}_x\big(s,x_*(s),u_*(s),p(\cdot),\lambda_0\big) \, ds \nonumber \\
           \label{SatzPMPIGZB1} & & -\sum_{j=1}^l \int_t^{t_1} g_{j,x}\big(s,x_*(s)\big)\, d\mu_j(s),
           \end{eqnarray}
\item[(b)] die Transversalit"atsbedingungen
           \index{Transversalitätsbedingungen!Integral@-- Integralgleichungen}
           \begin{equation}\label{SatzPMPIGZB2}
           p(t_0) = {h_0'}^T\big(x_*(t_0)\big)l_0, \qquad p(t_1)=0
           \end{equation}
\item[(c)] und f"ur fast alle $t \in [t_0,t_1]$ die Maximumbedingung
           \index{Maximumbedingung!Integral@-- Integralgleichungen}
           \begin{equation}\label{SatzPMPIGZB3}
           H^{\mathcal{I}}\big(t,x_*(t),u_*(t),p(\cdot),\lambda_0\big) = \max_{u \in U} H^{\mathcal{I}}\big(t,x_*(t),u,p(\cdot),\lambda_0\big)
           \end{equation}
\end{enumerate}
erfüllt sind.
\end{theorem}

Die adjungierte Gleichung (\ref{SatzPMPIGZB1}) lautet ausführlich geschrieben
\begin{eqnarray*}
p(t) &=& -{h_1'}^T\big(x_*(t_1)\big)l_1 - \int_t^{t_1} \bigg[\int_s^{t_1} \varphi_x^T\big(\tau,s,x_*(s),u_*(s)\big) \, dp(\tau) \bigg] \,ds \\
     & & - \int_t^{t_1} \lambda_0 f_x\big(s,x_*(s),u_*(s)\big) \, ds -\sum_{j=1}^l \int_t^{t_1} g_{j,x}\big(s,x_*(s)\big)\, d\mu_j(s).
\end{eqnarray*}
An der Stelle $t=t_1$ kann $p(\cdot)$ unstetig sein und eine Sprungstelle besitzen.
In diesem Fall gilt in $t_1$ die Transversalit"atsbedingung
\begin{equation} \label{SatzPMPIGZB4}
p(t_1^-) = -{h_1'}^T\big(x_*(t_1)\big)l_1-\sum_{j=1}^l \mu_j(\{t_1\}) \, g_{j,x}\big(t_1,x_*(t_1)\big).
\end{equation}

%% file: 5-41-Beweis.tex
\subsubsection{Der Nachweis der notwendigen Optimalit"atsbedingungen}
Wir erweitern im Beweis von Theorem \ref{SatzPMPIGS} in Abschnitt \ref{AbschnittBeweisPMPIG} die 
Extremalaufgabe (\ref{ExtremalaufgabePMPIG}) um die Abbildungen
$$G_j\big(x(\cdot)\big) = \max_{t \in [t_0,t_1]} g_j\big(t,x(t)\big), \quad G_j: C([t_0,t_1],\R^n) \to \R, \quad j=1,...,l$$
Damit entsteht zur Aufgabe (\ref{IG1})--(\ref{IG5}) die Extremalaufgabe
\begin{equation} \label{ExtremalaufgabePMPIGZA}
J\big(x(\cdot),u(\cdot)\big) \to \inf, \quad \mathscr{F}\big(x(\cdot),u(\cdot)\big)=0, \quad G_j\big(x(\cdot)\big) \leq 0,
\quad u(\cdot) \in L_\infty([t_0,t_1],U),
\end{equation}
für die die Voraussetzungen in den Punkten (A)--(C) des Extremalprinzips \ref{SatzExtremalprinzipStark} nach den Ausführungen
in den Abschnitten \ref{AbschnittBeweisPMPIG} und \ref{AbschnittBeweisPMPZA} erfüllt sind.
Zur Extremalaufgabe (\ref{ExtremalaufgabePMPIGZA}) hat die Lagrange-Funktion $\mathscr{L}=\mathscr{L}\big(x(\cdot),u(\cdot),\lambda_0,y^*,l_0,l_1,\lambda\big)$,
die über
$$C([t_0,t_1],\R^n)\times L_\infty([t_0,t_1],\R^m)\times \R\times C_0^*([t_0,t_1],\R^n)\times\R^{s_0}\times\R^{s_1} \times \R^l$$
definiert ist, die Form
$$\mathscr{L}= \lambda_0 J\big(x(\cdot),u(\cdot)\big)+ \big\langle y^*, F\big(x(\cdot),u(\cdot)\big) \big\rangle
                         +l_0^T H_0\big(x(\cdot)\big)+l_1^T H_1\big(x(\cdot)\big) + \sum_{j=1}^l \lambda_j G_j\big(x(\cdot)\big).$$
Ist $\big(x_*(\cdot),u_*(\cdot)\big)$ eine starke lokale Minimalstelle der Extremalaufgabe (\ref{ExtremalaufgabePMPIGZA}),
dann existieren nach Theorem \ref{SatzExtremalprinzipStark}
nicht gleichzeitig verschwindende Lagrangesche Multiplikatoren $\lambda_0 \geq 0$, $y^* \in C_0^*([t_0,t_1],\R^n)$, $l_i \in \R^{s_i}$
und $\lambda_1 \geq 0,...,\lambda_l \geq 0$ derart,
dass gelten:
\begin{enumerate}
\item[(a)] Die Lagrange-Funktion besitzt bez"uglich $x(\cdot)$ in $x_*(\cdot)$ einen station"aren Punkt, d.\,h.
           \begin{equation}\label{SatzPMPIGZALMR1}
           0 \in \partial_x \mathscr{L}\big(x_*(\cdot),u_*(\cdot),\lambda_0,y^*,l_0,l_1,\lambda\big);
           \end{equation}         
\item[(b)] Die Lagrange-Funktion erf"ullt bez"uglich $u(\cdot)$ in $u_*(\cdot)$ die Minimumbedingung
           \begin{equation}\label{SatzPMPIGZALMR2}
           \hspace*{-3mm} \mathscr{L}\big(x_*(\cdot),u_*(\cdot),\lambda_0,y^*,l_0,l_1,\lambda\big)
           = \min_{u(\cdot) \in L_\infty([t_0,t_1],U)} \mathscr{L}\big(x_*(\cdot),u(\cdot),\lambda_0,y^*,l_0,l_1,\lambda\big);
           \end{equation}
\item[(c)] Die komplement"aren Schlupfbedingungen gelten, d.\,h.
           \begin{equation}\label{SatzPMPIGZALMR3}
           0 = \lambda_j G_j\big(x(\cdot)\big), \qquad i=1,...,l.
           \end{equation}
\end{enumerate}

Ebenso wie im Beweis der Standardaufgabe unter Zustandsbeschränkungen im Abschnitt \ref{AbschnittBeweisPMPZA} liefert (\ref{SatzPMPIGZALMR3}),
dass nur diejenigen Multiplikatoren $\lambda_j \geq 0$ von Null verschieden sein können,
für die die zugehörigen nichtnegativen regul"aren Borelschen Ma"se $\tilde{\mu}_j$ die Totalvariation $\|\tilde{\mu}_j\|=1$ besitzen
und die auf den Mengen $T_j = \big\{ t \in [t_0,t_1] \big| g_j\big(t,x_*(t)\big) = 0 \big\}$ konzentriert sind. 
Deswegen können wir annehmen,
dass alle Ma"se $\mu_j=\lambda_j \tilde{\mu}_j$ auf den Mengen $T_j$ konzentriert sind.
\newpage
Aufgrund (\ref{SatzPMPIGZALMR1}) ist folgende Variationsgleichung f"ur alle $x(\cdot) \in C([t_0,t_1],\R^n)$ erf"ullt: 
\begin{eqnarray}
0 &=& \lambda_0 \int_{t_0}^{t_1} \big\langle f_x\big(t,x_*(t),u_*(t)\big),x(t) \big\rangle\, dt \nonumber \\
  & & + \big\langle l_0, h_0'\big(x_*(t_0)\big) x(t_0) \big\rangle + \big\langle l_1, h_1'\big(x_*(t_1)\big) x(t_1) \big\rangle
      + \sum_{j=1}^l \int_{t_0}^{t_1} \big\langle g_{j,x}\big(t,x_*(t)\big),x(t) \big\rangle \,d\mu_j(t)
      \nonumber \\
  & & \label{BeweisschlussPMPIGZA1}
      + \int_{t_0}^{t_1} \bigg[ x(t)-x(0) - \int_{t_0}^t \varphi_x\big(t,s,x_*(s),u_*(s)\big) x(s) \,ds \bigg]^T d\mu(t).
\end{eqnarray}
In der Gleichung (\ref{BeweisschlussPMPIGZA1}) "andern wir die Integrationsreihenfolge im letzten Summanden und bringen
diese Gleichung in die Form
\begin{eqnarray}
0 &=& \int_{t_0}^{t_1} \Big\langle \lambda_0 f_x\big(t,x_*(t),u_*(t)\big)
              -  \int_{t}^{t_1} \varphi_x^T\big(\tau,t,x_*(t),u_*(t)\big) d\mu(\tau) , x(t) \Big\rangle \, dt\nonumber \\
  & & + \int_{t_0}^{t_1} \langle x(t) , d\mu(t) \rangle
      + \Big\langle {h_0'}^T\big(x_*(t_0)\big)l_0 - \int_{t_0}^{t_1} d\mu(t) , x(t_0) \Big\rangle
      + \langle {h_1'}^T\big(x_*(t_1)\big)l_1 , x(t_1) \rangle \nonumber \\
  & & \label{BeweisschlussPMPIGZA2} + \sum_{j=1}^l \int_{t_0}^{t_1} \big\langle g_{j,x}\big(t,x_*(t)\big),x(t) \big\rangle \,d\mu_j(t).
\end{eqnarray}
Hierin setzen wir $\displaystyle p(t) = \int_t^{t_1} d\mu(\tau)$.
Dann ist $p(\cdot)$ eine rechtsseitig stetige Funktion beschränkter Variation und es gilt $p(t_1)=0$. \\[2mm]
Die rechte Seite in (\ref{BeweisschlussPMPIGZA2}) definiert ein stetiges lineares Funktional im Raum $C([t_0,t_1],\R^n)$.
Nach dem Satz von Riesz ist die Darstellung eines stetigen linearen Funktionals eindeutig und 
es ergeben sich aus (\ref{BeweisschlussPMPIGZA2}) die Beziehungen 
\begin{eqnarray*}
p(t) &=& -{h_1'}^T\big(x_*(t_1)\big)l_1 + \int_t^{t_1} H^{\mathcal{I}}_x\big(s,x_*(s),u_*(s),p(\cdot),\lambda_0\big) \, ds \\
& & - \sum_{j=1}^l \int_{t_0}^{t_1} \big\langle g_{j,x}\big(t,x_*(t)\big),x(t) \big\rangle \,d\mu_j(t), \\
p(t_0) &=&  {h_0'}^T\big(x_*(t_0)\big)l_0, \quad p(t_1)=0.
\end{eqnarray*}
Die Adjungierte $p(\cdot)$ kann als rechtsseitig stetige Funktion in $t=t_1$ unstetig sein.
Für den linksseitigen Grenzwert ergibt sich
$$p(t_1^-)=-{h_1'}^T\big(x_*(t_1)\big)l_1-\sum_{j=1}^l \mu_j(\{t_1\}) \, g_{j,x}\big(t_1,x_*(t_1)\big).$$
Damit sind (\ref{SatzPMPIGZB1}) und (\ref{SatzPMPIGZB2}), sowie die Sprungbedingung (\ref{SatzPMPIGZB4}) gezeigt. \\[2mm]
Die Beziehung (\ref{SatzPMPIGZALMR2}) führt nach Vertauschen der Integrationsreihenfolge zur Ungleichung
\begin{eqnarray*}
\lefteqn{\int_{t_0}^{t_1} \bigg[\lambda_0 f\big(t,x_*(t),u_*(t)\big)
              +  \int_{t}^{t_1} \big\langle\varphi\big(\tau,t,x_*(t),u_*(t)\big) , dp(\tau) \big\rangle \bigg] \, dt} \\
  &\leq& \int_{t_0}^{t_1} \bigg[\lambda_0 f\big(t,x_*(t),u(t)\big)
              +  \int_{t}^{t_1} \big\langle\varphi\big(\tau,t,x_*(t),u(t)\big) , dp(\tau) \big\rangle \bigg] \, dt,
\end{eqnarray*}
die f"ur alle $u(\cdot) \in L_\infty([t_0,t_1],U)$ erfüllt ist.
Wir passen die Argumentation im Abschnitt \ref{AbschnittBeweisPMPIG} an das Vorhandensein aktiver Zustandsbeschränkungen an:
Für die Steuerungen $u_\lambda(\cdot)$ ergibt sich die Ungleichung
\begin{eqnarray*}
\lefteqn{\int_{t_0}^{t_1} \chi_{[s,s+\lambda]}(t) \cdot \bigg[ 
                          \lambda_0 \big[ f\big(t,x_*(t),u_*(t)\big)-f\big(t,x_*(t),u\big)\big]} \\
&& 
   + \int_{t}^{t_1} \big\langle\varphi\big(\tau,t,x_*(t),u_*(t)\big)-\varphi\big(\tau,t,x_*(t),u\big), dp(\tau) \big\rangle
   \bigg] \, dt \leq 0.
\end{eqnarray*}
Zu $u \in U$ bzw. $u_*(\cdot)$ führen wir die Abbildungen
$$[\Phi(u)](t)=\int_{t}^{t_1} \big\langle \varphi\big(\tau,t,x_*(t),u\big), dp(\tau) \big\rangle, \quad
  [\Phi(u_*)](t)=\int_{t}^{t_1} \big\langle \varphi\big(\tau,t,x_*(t),u_*(t)\big), dp(\tau) \big\rangle$$
ein.
Aus der Setzung folgt, dass die Abbildung $t \to [\Phi(u)](t)$ über $[t_0,t_1]$ messbar und beschränkt ist.
Au"serdem ist $(\tau,t) \to \varphi\big(\tau,t,x_*(t),u\big)$ für jedes $u \in U$ stetig.
Deswegen besitzt $t \to [\Phi(u)](t)$ nur an den Stellen $t \in [t_0,t_1]$ keinen Lebesgueschen Punkt,
in denen die Funktion $p(t)=\displaystyle \int_t^{t_1} d\mu(\tau)$ unstetig ist.
Da $\mu$ ein signiertes reguläres Borelsches Vektorma"s ist,
ist $p(\cdot)$ über $[t_0,t_1]$ von beschränkter Variation und besitzt nur abzählbar viele Unstetigkeitsstellen. \\
Aus diesem Grund besitzt die Abbildung $t \to [\Phi(u)](t)$ über $[t_0,t_1]$ nur an abzählbar vielen Stellen keine Lebesgueschen Punkte.
Da diese Stellen sich nur in den Unstetigkeitsstellen der Funktion $p(\cdot)$ ergeben können,
befinden sich die Lebesgueschen Punkte unabhängig von der Wahl von $u \in U$ stets an den gleichen Stellen. \\[1mm]
Es sei $t \in (t_0,t_1)$ ein Lebesguescher Punkt der messbar und beschränkten Abbildungen
$$t \to f\big(t,x_*(t),u_*(t)\big), \quad t \to [\Phi(u_*)](t), \quad t \to [\Phi(u)](t).$$
Wie im Abschnitt \ref{AbschnittBeweisPMPIG} ergibt sich die Gültigkeit von
$$\lambda_0 \big( f\big(t,x_*(t),u_*(t)\big)-f\big(t,x_*(t),u\big)\big) + \big([\Phi(u_*)](t) - [\Phi(u)](t)\big) \leq 0$$
bzw. in äquivalenter Form mit Hilfe der Pontrjagin-Funktion (\ref{PontrjaginFunktionIG}) die Gültigkeit von
$$H^{\mathcal{I}}\big(t,x_*(t),u_*(t),p(\cdot),\lambda_0\big) \geq H^{\mathcal{I}}\big(t,x_*(t),u,p(\cdot),\lambda_0\big)$$
f"ur alle $u \in U$.
Da unabhängig von $u \in U$ fast überall in $(t_0,t_1)$ ein Lebesguescher Punkt vorliegt,
gilt fast überall über $[t_0,t_1]$ die Maximumbedingung (\ref{SatzPMPIGZB3}). \hfill $\blacksquare$

%% file: 5-42-HinreichendZB.tex
\subsubsection{Hinreichende Bedingungen nach Arrow} \label{AbschnittArrowIGLZB}
Bei \index{hinreichende Bedingungen, Arrow!Integral@-- Integralgleichungen}
der Herleitung hinreichender Bedingungen nach Arrow für die Steuerung von Integralgleichungen
unter Zustandsbeschränkungen greifen wir auf die Argumente der Abschnitte \ref{AbschnittHBPMP}, \ref{AbschnittArrowZB} und \ref{AbschnittArrowIGL} zurück. \\[2mm]
Wir betrachten das Steuerungsproblem mit festen bzw. freien Randwerten
\begin{eqnarray}
&& \label{HBZAIGL1} J\big(x(\cdot),u(\cdot)\big) = \int_{t_0}^{t_1} f\big(t,x(t),u(t)\big) \, dt \to \inf, \\
&& \label{HBZAIGL2} x(t) = x(t_0)+ \int_{t_0}^t\varphi\big(t,s,x(s),u(s)\big) \, ds, \quad t \in [t_0,t_1], \\
&& \label{HBZAIGL3} x(t_0)=x_0, \qquad x(t_1)=x_1, \\
&& \label{HBZAIGL4} u(t) \in U \subseteq \R^m, \quad U\not= \emptyset, \\
&& \label{HBZAIGL5} g_j\big(t,x(t)\big) \leq 0 \quad \mbox{f"ur alle } t \in [t_0,t_1], \quad j=1,...,l.
\end{eqnarray}

Im Weiteren bezeichnet $V^{\mathcal{S}}_\gamma(t)$ die Menge $V^{\mathcal{S}}_\gamma(t)=\{ x \in \R^n \,|\, \|x-x_*(t)\| < \gamma\}$.
Die Hamilton-Funktion $\mathscr{H}^{\mathcal{I}}$ besitzt für die Aufgabe (\ref{HBZAIGL1})--(\ref{HBZAIGL5}) die Gestalt
$$\mathscr{H}^{\mathcal{I}}\big(t,x,p(\cdot)\big)
  = \sup_{u \in U} \bigg[-\int_t^{t_1} \langle \varphi(\tau,t,x,u),dp(\tau) \rangle - f(t,x,u)\bigg].$$

\begin{theorem} \label{SatzHBIGLZA}
In der Aufgabe (\ref{HBZAIGL1})--(\ref{HBZAIGL5})
sei $\big(x_*(\cdot),u_*(\cdot)\big) \in \mathscr{B}^{\,\mathcal{I}}_{\rm Lip} \cap \mathscr{B}^{\,\mathcal{I}}_{\rm adm}$.
Au"serdem sei die Adjungierte $p(\cdot)$ über $[t_0,t_1]$ st"uckweise absolutstetig
und besitze in $[t_0,t_1]$ h"ochstens endlich viele Sprungstellen.
Ferner gelte:
\begin{enumerate}
\item[(a)] Das Tripel $\big(x_*(\cdot),u_*(\cdot),p(\cdot)\big)$
           erf"ullt (\ref{SatzPMPIGZB1})--(\ref{SatzPMPIGZB3}) in Theorem \ref{SatzPMPIGZB} mit $\lambda_0=1$.        
\item[(b)] F"ur jedes $t \in [t_0,t_1]$ ist die Funktion $\mathscr{H}^{\mathcal{I}}\big(t,x,p(\cdot)\big)$ konkav 
           und es sind die Funktionen $g_j(t,x)$, $j=1,...,l$, konvex bez"uglich $x$ auf $V^{\mathcal{S}}_\gamma(t)$.
\end{enumerate}
Dann ist $\big(x_*(\cdot),u_*(\cdot)\big)$ ein starkes lokales Minimum der Aufgabe (\ref{HBZAIGL1})--(\ref{HBZAIGL5}).
\end{theorem}

{\bf Beweis} Wir greifen auf Bemerkung \ref{BemHBSOPZB} zurück:
Es bezeichnen $t_0<s_1<...<s_d<t_1$ die Stellen im Intervall $(t_0,t_1)$,
in denen die Adjungierte $p(\cdot)$ unstetig ist.
Dann folgen aus der adjungierten Gleichung \ref{SatzPMPIGZB1} und der Nichtnegativität der Ma"se $\mu_j$ die Sprungbedingungen
\begin{equation} \label{HBPMPIGZB0}
p(s_k)-p(s_k^-)= \sum_{j=1}^l \beta_j^k g_{j,x}\big(s_k,x_*(s_k)\big),\qquad \beta_j^k = \mu_j(\{s_k\}) \geq 0, \quad k=1,...,d.
\end{equation}
Weiterhin gibt es eine st"uckweise absolutstetige Funktion $\lambda(\cdot):[t_0,t_1] \to \R^l$ derart, dass
$$\dot{p}(t)=- H^{\mathcal{I}}_x\big(t,x_*(t),u_*(t),p(t),1\big) + \sum_{j=1}^l \lambda_j(t) g_{j,x}\big(t,x_*(t)\big)$$
st"uckweise auf $(s_{k-1},s_k)$ für $k=1,...,d+1$ gilt. Dabei sind $s_0=t_0$ und $s_{d+1}=t_1$. 
Ferner gelten aufgrund der komplementären Schlupfbedingungen
$$\beta_j^k g_j\big(s_k,x_*(s_k)\big) =0, \quad \lambda_j(t)g_j\big(t,x_*(t)\big)=0, \quad \beta_j^k \geq 0, \quad \lambda_j(t) \geq 0$$
f"ur $j=1,...,l$, $k=1,...,d$ und $t \in (s_{k-1},s_k)$ für $k=1,...,d+1$. \\[2mm]
Mit den Argumenten aus dem Abschnitt  \ref{AbschnittHBPMP} folgt die Gültigkeit der Ungleichung
$$\langle a(t),x-x_*(t)\rangle \leq \mathscr{H}^{\mathcal{I}}\big(t,x_*(t),p(t)\big) - \mathscr{H}^{\mathcal{I}}\big(t,x,p(t)\big)$$
für fast alle $t$ und alle $x \in V^{\mathcal{S}}_\gamma(t)$.
Au"serdem gilt $a(t)=-H^{\mathcal{I}}_x\big(t,x_*(t),u_*(t),p(t),1\big)$ f"ur fast alle $t \in [t_0,t_1]$.
Mit der adjungierten Gleichung \ref{SatzPMPIGZB1} ergibt sich daraus ferner
\begin{equation} \label{HBPMPIGZB1}
\mathscr{H}^{\mathcal{I}}\big(t,x_*(t),p(t)\big) - \mathscr{H}^{\mathcal{I}}\big(t,x,p(t)\big) 
  \geq \Big\langle \dot{p}(t) - \sum_{j=1}^l \lambda_j(t) g_{j,x}\big(t,x_*(t)\big) , x-x_*(t) \Big\rangle.
\end{equation}
Die Konvexit"at der Abbildungen $g_j(t,x)$ bez"uglich $x$ auf $V^{\mathcal{S}}_\gamma(t)$ liefert die Ungleichungen
\begin{equation} \label{HBPMPIGZB2}
-\lambda_j(t) \big\langle g_{j,x}\big(t,x_*(t)\big) , x(t)-x_*(t) \big\rangle \geq - \lambda_j(t) g_j\big(t,x(t)\big)\geq 0
\end{equation}
f"ur $j=1,...,l$, $t \in (s_{k-1},s_k)$ und $k=1,...,d+1$;
au"serdem
\begin{equation} \label{HBPMPIGZB3}
-\sum_{j=1}^l \big\langle \beta_j^k g_{j,x}\big(s_k,x_*(s_k)\big), x(t)-x_*(t) \big\rangle \geq 0
\end{equation}
in den Unstetigkeitsstellen $s_k \in (t_0,t_1),\, k=1,...,d$ der Adjungierten. \\[2mm]
Im Endpunkt $t=t_1$ gelten nach (\ref{SatzPMPIGZB2}) und (\ref{SatzPMPIGZB4}) die Transversalitätsbedingungen
\begin{equation} \label{HBPMPIGZB4}
p(t_1^-)=-l_1-\sum_{j=1}^l  \mu_j(\{t_1\}) \, g_{j,x}\big(t_1,x_*(t_1)\big), \quad p(t_1)=0.
\end{equation}
Wiederum aufgrund der Konvexit"at der Abbildungen $g_j(t,x)$ und der komplementären Schlupfbedingungen gilt in $t=t_1$ die Ungleichung
\begin{equation} \label{HBPMPIGZB5}
-\sum_{j=1}^l \big\langle \mu_j(\{t_1\}) \, g_{j,x}\big(t_1,x_*(t_1)\big), x(t_1)-x_*(t_1) \big\rangle \geq 0
\end{equation}
für alle $x(t_1) \in V_\gamma(t_1)$ mit $g_j\big(t_1,x(t_1)\big) \leq 0$.

\newpage
Wir betrachten die Differenz
$$\Delta = J\big(x(\cdot),u(\cdot)\big)- J\big(x_*(\cdot),u_*(\cdot)\big)
  = \int_{t_0}^{t_1} \big[f\big(t,x(t),u(t)\big)-f\big(t,x_*(t),u_*(t)\big)\big] \, dt.$$
Es sei $\big(x(\cdot),u(\cdot)\big) \in  \mathscr{B}^{\,\mathcal{I}}_{\rm adm}$ mit $\|x(\cdot)-x_*(\cdot)\|_\infty < \gamma$.
Dann gilt
\begin{eqnarray*}
\Delta &=& J\big(x(\cdot),u(\cdot)\big)- J\big(x_*(\cdot),u_*(\cdot)\big)
\geq \int_{t_0}^{t_1} \big[\mathscr{H}^{\mathcal{I}}\big(t,x_*(t),p(\cdot)\big)-\mathscr{H}^{\mathcal{I}}\big(t,x(t),p(\cdot)\big)\big] \, dt \\
& & \hspace*{15mm} + \int_{t_0}^{t_1} \bigg[\int_t^{t_1}
    \big\langle \varphi\big(\tau,t,x_*(t),u_*(t)\big) - \varphi\big(\tau,t,x(t),u(t)\big),dp(\tau) \big\rangle \bigg]\, dt.
\end{eqnarray*}
Wir ändern im letzten Ausdruck die Integrationsreihenfolge,
\begin{eqnarray*}
\lefteqn{\int_{t_0}^{t_1} \bigg[\int_t^{t_1}
     \big\langle \varphi\big(\tau,t,x_*(t),u_*(t)\big) - \varphi\big(\tau,t,x(t),u(t)\big),dp(\tau) \big\rangle
     \bigg]\, dt} \\
&=& \int_{t_0}^{t_1} \bigg[\int_{t_0}^t \big[ \varphi\big(t,s,x_*(s),u_*(s)\big) - \varphi\big(t,s,x(s),u(s)\big) \big] \, ds \bigg]^T \, dp(t),
\end{eqnarray*}  
und erhalten unter Verwendung von (\ref{HBPMPIG2}), sowie (\ref{HBPMPIGZB1}) und (\ref{HBPMPIGZB2})
\begin{eqnarray*}
\Delta &\geq& \int_{t_0}^{t_1} \big[\mathscr{H}^{\mathcal{I}}\big(t,x_*(t),\mu\big)-\mathscr{H}^{\mathcal{I}}\big(t,x(t),\mu\big)\big] \, dt  \\
       &    & \hspace*{2cm} + \int_{t_0}^{t_1} \langle x_*(t)-x(t)-x_*(t_0)+x(t_0), dp(t) \rangle \\
       &\geq&  -\int_{t_0}^{t_1} \langle x(t)-x_*(t), \dot{p}(t)\rangle \, dt + \int_{t_0}^{t_1} \langle x(t)-x_*(t), dp(t)\rangle \\
       &    & \hspace*{2cm} - \sum_{j=1}^l  \int_{t_0}^{t_1} \big\langle \lambda_j(t) g_{j,x}\big(t,x_*(t)\big) , x-x_*(t) \big\rangle \,dt \\
       &\geq&  -\int_{t_0}^{t_1} \langle x(t)-x_*(t), \dot{\mu}(t)\rangle \, dt + \int_{t_0}^{t_1} \langle x(t)-x_*(t), d\mu(t)\rangle.
\end{eqnarray*}
Die Adjungierte $p(\cdot)$ ist stückweise absolutstetig.
Daher gilt
\begin{eqnarray*}
\lefteqn{\int_{t_0}^{t_1} \langle x(t)-x_*(t), dp(t)\rangle = \int_{t_0}^{t_1} \langle x(t)-x_*(t), \dot{p}(t)\rangle \, dt} \\
& & +\sum_{k=1}^d \langle x(s_k)-x_*(s_k), p(s_k)-p(s_k^-)\rangle + \langle x(t_1)-x_*(t_1), p(t_1)-p(t_1^-)\rangle.
\end{eqnarray*}
In den Unstetigkeitsstellen $s_k$ gelten (\ref{HBPMPIGZB0}) und (\ref{HBPMPIGZB3}),
im Endpunkt $t=t_1$ die Ungleichungen (\ref{HBPMPIGZB4}) und (\ref{HBPMPIGZB5}).
Zusammen erhalten wir abschlie"send
$$\Delta \geq -\langle l_1, x(t_1)-x_*(t_1)\rangle - \langle l_0 , x(t_0)-x_*(t_0) dt=0.$$
Daher gilt $\Delta \geq 0$ f"ur alle zul"assigen $\big(x(\cdot),u(\cdot)\big)$ mit 
$\|x(\cdot)-x_*(\cdot)\|_\infty < \gamma$. \hfill $\blacksquare$

%% file: 5-43-Werbestrategien.tex
\subsubsection{Optimale Werbestrategien bei beschränktem Werbebudget}
Wir betrachten wieder die Aufgabe optimaler Werbestrategien:
\begin{eqnarray}
&& \label{Werbung31} J\big(x(\cdot),u(\cdot)\big) = \int_0^T e^{-\varrho t} \big[\pi x(t)-u(t)\big] \, dt \to \sup, \\
&& \label{Werbung32} x(t)=x_0+\int_0^t e^{-\delta(t-s)} \big[f\big(u(s)\big)-\alpha x(s)\big] \, ds, \quad t \in [0,T], \\
&& \label{Werbung33} x_0 \in [0,1],\quad u(t) \geq 0, \quad \varrho >0, \quad \alpha, \delta \geq 0.
\end{eqnarray}
Die Aufgabe (\ref{Werbung31})--(\ref{Werbung33}) erweitern wir nun um die Zustandsbeschränkung
\begin{equation} \label{Werbung34} 
z(t)=\int_0^t e^{-\varrho s} u(s) \, ds, \quad t \in [0,T], \quad z(t) \leq Z , \quad Z >0.
\end{equation}
Weiterhin führen wir die Bezeichnungen aus den Abschnitten \ref{AbschnittWerbung1} und \ref{AbschnittWerbung2} fort. \\[1mm]
Da die Funktion $z(\cdot)$ über $[0,T]$ monoton wachsend,
ist es vernünftig die Diskussion der Aufgabe (\ref{Werbung31})--(\ref{Werbung34}) auf die folgenden Fälle zu reduzieren:

\begin{enumerate}
\item[(A)] Die Zustandsbeschränkung ist über $[0,T]$ nicht aktiv.
           Dann liefern die Optimalitätsbedingungen die optimale Werbestrategie aus Abschnitt \ref{AbschnittWerbung1}:
           $$f'\big(u_*(t)\big) = \frac{\alpha + \delta + \varrho}{\pi ( 1-e^{-(\alpha + \delta + \varrho)(T-t)})}.$$
           Dabei ist $u_*(t)$ stetig, streng monoton fallend und $u_*(T)=0$. 
\item[(B)] Die Zustandsbeschränkung ist nur in $t=T$ aktiv. 
           Es bezeichne $Z_0$ die Werbeausgaben im Fall (A) ohne Einschränkungen an das Werbebudget.
           Dann greift die Beschränkung (\ref{Werbung34}) nur dann, wenn $0<Z\leq Z_0$ erfüllt ist. \\[1mm]
           Im Fall $Z=Z_0$ ergeben sich bezüglich dem Zustand $z(\cdot)$ die Adjungierte $r(t) \equiv 0$,
           sowie bezüglich $x(\cdot)$ die Adjungierte $q(\cdot)$ und die Werbestrategie $u_*(\cdot)$ wie unter (A). \\[1mm]
           Im Fall $Z<Z_0$ ist nach Abschnitt \ref{AbschnittWerbung2} bezüglich dem Zustand $z(\cdot)$ die Adjungierte $r(t) \equiv r_0 < 0$.
           Für die Zustandsbeschränkung $z(t)-Z \leq 0$ ergibt sich wegen $r(T)=- \mu(\{T\}) < 0$, dass das in $t=T$ konzentrierte Ma"s positiv sein muss.
           Bezüglich $x(\cdot)$ erhalten wir die Adjungierte $q(\cdot)$ wie im Fall (A). 
           Au"serdem genügt die optimale Werbestrategie $u_*(\cdot)$ der Beziehung
           $$f'\big(u_*(t)\big) = (1-r_0) \cdot \frac{\alpha + \delta + \varrho}{\pi ( 1-e^{-(\alpha + \delta + \varrho)(T-t)})}.$$
           Die Steuerung $u_*(t)$ ist wieder stetig, streng monoton fallend und $u_*(T)=0$.
           Im Vergleich zum Fall (A) bewirkt der Faktor $(1-r_0)$ eine reduzierte Werberate $u_*(t)$.           
\item[(C)] Die Zustandsbeschränkung ist über dem Intervall $[\tau,T]$ mit $0<\tau <T$ aktiv.
           Dies würde $u_*(t)=0$ und damit insbesondere $f'\big(u_*(t)\big)=\infty$ über $[\tau,T]$ bedeuten.
           Das ist ein Widerspruch zur Gültigkeit der Maximumbedingung.
           Deswegen kann dieser Fall nicht eintreten. \hfill $\square$
\end{enumerate} 

%% file: 5-5-FreieZeit.tex
\subsection{Freier Anfangs- und Endzeitpunkt} \label{AbschnittFreieZeitIGL}
Die nachstehenden Betrachtungen sind durch die Arbeit von Dmitruk \& Osmolovskii \cite{DmitrukOsmo2} angeregt.
Im Wesentlichen folgen wir deren Bearbeitung.
Allerdings fügen wir der Aufgabe im Zielfunktional die äußeren Integrationsgrenzen hinzu:
\begin{eqnarray}
&& \label{IGLFZ1} J\big(x(\cdot),u(\cdot)\big) = \int_{t_0}^{t_1} f\big(t,t_0,t_1,x(t),u(t)\big) \, dt \to \inf, \\
&& \label{IGLFZ2} x(\zeta) = x(t_0) + \int_{t_0}^\zeta \varphi\big(\zeta,t,x(t),u(t)\big), \quad \zeta \in [t_0,t_1], \\
&& \label{IGLFZ3} h_0\big(t_0,x(t_0)\big)=0, \qquad h_1\big(t_1,x(t_1)\big)=0, \\
&& \label{IGLFZ4} u(t) \in U \subseteq \R^m, \quad U \not= \emptyset.
\end{eqnarray}
Der Umgang mit den Bezeichnungen bezüglich der ``inneren'' und ``äußeren'' Zeitvariablen gestaltet sich durch die Methode der Substitution
der Zeit und dem Auftreten von Zeitableitungen unübersichtlicher.
Deswegen bezeichnet im Folgenden $t$ die ``innere'' bzw. die in der Standardaufgabe ``übliche'' Zeit,
während $\zeta$ für die ``äußere'' Zeitvariable steht.
Damit liegen die Variablenbezeichnungen $f(t,t_0,t_1,x,u)$, $\varphi(\zeta,t,x,u)$ und $h_i(t_i,x_i)$ vor. \\[2mm]
Die Aufgabe (\ref{IGLFZ1})--(\ref{IGLFZ4}) untersuchen wir f"ur Tripel $\big([t_0,t_1],x(\cdot),u(\cdot)\big)$ mit
$$[t_0,t_1] \subset \R, \qquad x(\cdot) \in C([t_0,t_1],\R^n), \qquad u(\cdot) \in L_\infty([t_0,t_1],U).$$
Zur Menge $\mathscr{B}^{\,\mathcal{F}}_{\rm Lip}$ geh"oren diejenigen Tripel
$\big([t_0,t_1],x(\cdot),u(\cdot)\big)$,
f"ur die es eine Zahl $\gamma>0$ derart gibt,
dass die Abbildungen $f(t,\tau_0,\tau_1,x,u)$, $\varphi(\zeta,t,x,u)$ und $h_i(\tau_i,x_i)$
auf der Menge aller Punkte $(\zeta,t,\tau_0,\tau_1,x,x_0,x_1,u) \in \R \times \R \times \R \times \R \times \R^n \times \R^n \times \R^n \times \R^m$ mit
\begin{eqnarray*}
&& t_0-\gamma < t,\zeta < t_1+\gamma,\quad t_0-\gamma < \tau_0< t_0+\gamma,\quad t_1-\gamma < \tau_1 < t_1+\gamma, \\
&& \|x-x(t)\| < \gamma, \quad \|x_0-x(t_0)\| < \gamma, \quad \|x_1-x(t_1)\| < \gamma, \quad u \in \R^m
\end{eqnarray*}
stetig in der Gesamtheit aller Variablen und
stetig differenzierbar bez"uglich der Zeit- und Zustandsvariablen $\zeta, t, t_0, t_1, x, x_0, x_1$ sind.
(Zur unmissverständlichen Angabe der Punktemenge treten $\tau_0,\tau_1$ in $f$ und $h_i$ anstelle der Zeitvariablen $t_0,t_1$ auf.) \\[2mm]
In der Aufgabe (\ref{IGLFZ1})--(\ref{IGLFZ4}) mit freiem Anfangs- und Endzeitpunkt nennen wir ein Tripel
$\big([t_0,t_1],x(\cdot),u(\cdot)\big)$ mit
$[t_0,t_1] \subset \R$, $x(\cdot) \in C([t_0,t_1],\R^n)$ und $u(\cdot) \in L_\infty([t_0,t_1],U)$ einen Steuerungsprozess.
Ein Steuerungsprozess hei"st zul"assig in dem Steuerungsproblem (\ref{IGLFZ1})--(\ref{IGLFZ4}),
wenn auf dem Intervall $[t_0,t_1]$ die Funktion $x(\cdot)$ fast "uberall der Gleichung (\ref{IGLFZ2}) gen"ugt und
die Randbedingungen (\ref{IGLFZ3}) erf"ullt.
Die Menge $\mathscr{B}^{\,\mathcal{F}}_{\rm adm}$ bezeichnet die Menge der zul"assigen Steuerungsprozesse. \\[2mm]
Einen zul"assigen Steuerungsprozess $\big([t_{0*},t_{1*}],x_*(\cdot),u_*(\cdot)\big)$ nennen wir ein starkes lokales
Minimum\index{Minimum, starkes lokales!Integral@-- Integralgleichungen},
wenn eine Zahl $\varepsilon > 0$ derart existiert,
dass f"ur jeden zul"assigen Steuerungsprozess $\big([t_0,t_1],x(\cdot),u(\cdot)\big)$ mit den Eigenschaften
$$|t_0 - t_{0*}| < \varepsilon, \qquad |t_1 - t_{1*}| < \varepsilon, \qquad
  \| x(t)-x_*(t) \| < \varepsilon \quad \mbox{ f"ur alle } t \in [t_{0*},t_{1*}] \cap [t_0,t_1]$$
die Ungleichung $J\big(x(\cdot),u(\cdot)\big) \geq J\big(x_*(\cdot),u_*(\cdot)\big)$ gilt.

%% file: 5-51-Substitution.tex
\subsubsection{Die Methode der Substitution der Zeit}
Wir wenden die Methode der Substitution der Zeit aus Abschnitt \ref{AbschnittFreieZeitSOP} an,
$$t(s) = t(0) + v \cdot s, \quad y(s) = x\big(t(s)\big), \quad w(s) = u\big(t(s)\big), \quad v >0,$$
und erhalten folgende Aufgabe über dem festen Zeitintervall $[0,1]$:
\begin{equation} \label{IGLFZ6}
\left.\begin{array}{l}
\displaystyle \tilde{J}\big(t(\cdot),y(\cdot),w(\cdot),v\big)  = \int_0^1 v \cdot f\big(t(\sigma),t(0),t(1),y(\sigma),w(\sigma)\big) \, d\sigma \to \inf, \\[3mm]
\displaystyle y(s) = y(0) + \int_0^s v \cdot \varphi\big(t(s),t(\sigma),y(\sigma),w(\sigma)\big) \, d\sigma, \quad s \in [0,1], \\[3mm]
\displaystyle t(s) = t(0) + \int_0^s v \, d\sigma, \quad s \in [0,1], \\[3mm]
h_0\big(t(0),y(0)\big)=0, \qquad h_1\big(t(1),y(1)\big)=0, \\[1mm]
w(s) \in U \subseteq \R^m, \quad U \not= \emptyset, \quad v >0. \\[1mm]
\end{array}\right\}
\end{equation}
In der Aufgabe (\ref{IGLFZ6}) enthalten
das Zielfunktional bzw. die Dynamik die Größen $t(0)$ und $t(1)$, sowie $t(s)$,
die wir als Zustände der äußeren Zeit bezeichnen. \index{Zustände der äußeren Zeit}
Das Auftreten von Zuständen der äußeren Zeit wird durch die Methode der Substitution der Zeit hervorgerufen.
Auf dieses Detail sind wir durch die bemerkenswerte Arbeit \cite{DmitrukOsmo2} aufmerksam geworden.
Unserer Kenntnis nach wird die Steuerung einer Volterraschen Integralgleichung
mit freiem Anfangs- und Endzeitpunkt erstmals in \cite{DmitrukOsmo2} umfassend untersucht. \\[2mm]
Durch das Auftreten der Zustände der äußeren Zeit gliedert sich die Aufgabe (\ref{IGLFZ6}) nicht in die bisher betrachtete Klasse (\ref{IG1})--(\ref{IG4}) ein.
Wir sehen uns also mit einem allgemeineren Steuerungsproblem für Volterrasche Integralgleichungen konfrontiert. 
Deswegen ist die ausführliche Auswertung dieser Aufgabe erforderlich. \\[2mm]
In der Aufgabe (\ref{IGLFZ6}) nennen wir ein Quadrupel
$$\big(t(\cdot),y(\cdot),w(\cdot),v\big) \in C([0,1],\R) \times C([0,1],\R^n) \times L_\infty([0,1],U)\times \R$$
einen zul"assigen Steuerungspozess,
wenn die Dynamiken und die Randbedingungen dieser Aufgabe erf"ullt sind. \\[2mm]
Einen zul"assigen Steuerungsprozess $\big(t_*(\cdot),y_*(\cdot),w_*(\cdot),v_*\big)$ nennen wir ein starkes lokales
Minimum der Aufgabe (\ref{IGLFZ6}),
wenn eine Zahl $\varepsilon > 0$ derart existiert,
dass die Ungleichung
$$\tilde{J}\big(t(\cdot),y(\cdot),w(\cdot),v\big) \geq \tilde{J}\big(t_*(\cdot),y_*(\cdot),w_*(\cdot),v_*\big)$$
für jeden zul"assigen Steuerungsprozess $\big(t(\cdot),y(\cdot),w(\cdot),v\big)$
mit $\|t(\cdot)-t_*(\cdot) \|_\infty < \varepsilon$ und mit $\|y(\cdot)-y_*(\cdot) \|_\infty < \varepsilon$ gilt.

%% file: 5-52-Beweis.tex
\subsubsection{Die Auswertung der Aufgabe mit Zuständen der äußeren Zeit}
Wir wollen mit Hilfe des Extremalprinzips in Form von Theorem \ref{SatzExtremalprinzipStark}
notwendige Optimalitätsbedingungen für die Aufgabe (\ref{IGLFZ6}) ableiten.
Dazu betrachten wir für 
$$\big(t(\cdot),y(\cdot),w(\cdot),v\big) \in C([0,1],\R) \times C([0,1],\R^n) \times L_\infty([0,1],U)\times \R$$
die Abbildungen
\begin{eqnarray*}
&& \hspace*{-1cm} \tilde{J}\big(t(\cdot),y(\cdot),w(\cdot),v\big)  = \int_0^1 v \cdot f\big(t(s),t(0),t(1),y(s),w(s)\big) \, ds, \\
&& \hspace*{-1cm} \tilde{F}_1\big(t(\cdot),y(\cdot),w(\cdot),v\big)(s)
   = y(s) - y(0) - \int_0^s v \cdot \varphi\big(t(s),t(\sigma),y(\sigma),w(\sigma)\big) \, d\sigma, \\
&& \hspace*{-1cm} \tilde{F}_2\big(t(\cdot),v\big)(s) = t(s) - t(0) - \int_0^s v \, d\sigma, \quad s \in [0,1], \\
&& \hspace*{-1cm} \tilde{H}_0\big(t(\cdot),y(\cdot)\big) = h_0\big(t(0),y(0)\big)=0,
   \qquad \tilde{H}_1\big(t(\cdot),y(\cdot)\big) = h_1\big(t(1),y(1)\big)=0.
\end{eqnarray*}
Da die Zustände $t(\cdot), x(\cdot)$ stetige Funktionen sind, gilt f"ur diese Abbildungen
\begin{eqnarray*}
\tilde{J} &:& C([0,1],\R) \times C([0,1],\R^n) \times L_\infty([0,1],U) \times \R \to \R, \\
\tilde{F}_1 &:& C([0,1],\R) \times C([0,1],\R^n) \times L_\infty([0,1],U)\times \R  \to C_0([0,1],\R^n), \\
\tilde{F}_2 &:& C([0,1],\R) \times \R  \to C_0([0,1],\R), \\
\tilde{H}_i &:& C([0,1],\R) \times C([0,1],\R^n) \to \R^{s_i}, \quad i=0,1.
\end{eqnarray*}

Wir setzen $\mathscr{F}=(\tilde{F}_1,\tilde{F}_2,\tilde{H}_0,\tilde{H}_1)$ und pr"ufen f"ur die Extremalaufgabe
\begin{equation} \label{ExtremalaufgabePMPIGLFZ}
\tilde{J}\big(t(\cdot),y(\cdot),w(\cdot),v\big) \to \inf, \; \mathscr{F}\big(t(\cdot),y(\cdot),w(\cdot),v\big)=0, \;
\big(w(\cdot),v\big) \in L_\infty([0,1],U) \times \R
\end{equation}
im starken lokalen Minimum $\big(t_*(\cdot),y_*(\cdot),w_*(\cdot),v_*\big)$
die Voraussetzungen von Theorem \ref{SatzExtremalprinzipStark}:

\begin{enumerate}
\item[(A$_1$)] Wie im Beispiel \ref{DiffZielfunktional2} folgt,
               dass für jedes Paar $\big(w(\cdot),v\big) \in L_\infty([0,1],U) \times \R$
               die Abbildung $\big(t(\cdot),y(\cdot)\big) \to \tilde{J}\big(t(\cdot),y(\cdot),w(\cdot),v\big)$
               im Punkt $\big(t_*(\cdot),y_*(\cdot)\big)$ Fr\'echet-differenzierbar ist und für die Ableitungen die Darstellungen
               \begin{eqnarray*}
               \tilde{J}_{t(\cdot)} \big(t_*(\cdot),y_*(\cdot),w(\cdot),v\big)t(\cdot)
               &=& \int_0^1 v \cdot \Big[ f_t\big(t_*(s),t_*(0),t_*(1),y_*(s),w(s)\big) t(s) \\
               &&  \hspace*{12mm} + f_{t_0}\big(t_*(s),t_*(0),t_*(1),y_*(s),w(s)\big) t(0) \\
               &&  \hspace*{12mm} + f_{t_1}\big(t_*(s),t_*(0),t_*(1),y_*(s),w(s)\big) t(1) \Big] \, ds, \\
               \tilde{J}_{y(\cdot)} \big(t_*(\cdot),y_*(\cdot),w(\cdot),v\big) y(\cdot) 
               &=& \int_0^1 v \cdot \big\langle f_x\big(t_*(s),t_*(0),t_*(1),y_*(s),w(s)\big), y(s) \big\rangle \, ds
               \end{eqnarray*}
               gelten. Man beachte dabei weiterhin $f=f(t,t_0,t_1,x,u)$.
\item[(A$_2$)] Die Differenzierbarkeit der Abbildungen $\tilde{F}_1$ und $\tilde{F}_2$ bezüglich
               $\big(t(\cdot),y(\cdot)\big)$ folgt wie im Beispiel \ref{DiffDynamik2}.
               F"ur die Abbildungen $\tilde{H}_i$ ist die stetige Differenzierbarkeit offensichtlich.
               Es ergeben sich für die Dynamiken die Ableitungen
               \begin{eqnarray*}
               && \big[\tilde{F}_{1,t(\cdot)}\big(t(\cdot),y(\cdot),w(\cdot),v\big) \tau(\cdot)\big](s) \\
               && \quad = - \int_0^s v \cdot \big[\varphi_\zeta\big(t(s),t(\sigma),y(\sigma),w(\sigma)\big) \tau(s)
                          + \varphi_t\big(t(s),t(\sigma),y(\sigma),w(\sigma)\big) \tau(\sigma) \big] \, d\sigma, \\
               && \big[\tilde{F}_{1,y(\cdot)}\big(t(\cdot),y(\cdot),w(\cdot),v\big) \eta(\cdot)\big](s) \\
               && \quad = \eta(s) - \eta(0) - \int_0^s v \cdot \varphi_x\big(t(s),t(\sigma),y(\sigma),w(\sigma)\big) \eta(\sigma) \, d\sigma, \\
               && \big[\tilde{F}_{2,t(\cdot)}\big(t(\cdot),v\big)\tau(\cdot)\big](s) = \tau(s) - \tau(0), \quad s \in [0,1],
               \end{eqnarray*}
               und für die Ableitungen der Randbegingungen die Darstellungen
               \begin{eqnarray*}
               && \tilde{H}_{0,t(\cdot)}\big(t(\cdot),y(\cdot)\big) = h_{0,t_0}\big(t(0),y(0)\big),
                  \quad \tilde{H}_{0,y(\cdot)}\big(t(\cdot),y(\cdot)\big) = h_{0,x_0}\big(t(0),y(0)\big) \\
               && \tilde{H}_{1,t(\cdot)}\big(t(\cdot),y(\cdot)\big) = h_{1,t_1}\big(t(1),y(1)\big),
                   \quad \tilde{H}_{1,y(\cdot)}\big(t(\cdot),y(\cdot)\big) = h_{1,x_1}\big(t(1),y(1)\big).
               \end{eqnarray*}
               Auch hier bestehen noch die Zuweisungen $\varphi=\varphi(\zeta,t,x,u)$ und $h_i=h_i(t_i,x_i)$.
\item[(B)] Für den Nachweis der endlichen Kodimension müssen wir zeigen,
           dass es zu jedem $\big(z_1(\cdot),z_2(\cdot)\big) \in C_0([0,1],\R) \times C_0([0,1],\R^n)$
           ein $\big(\overline{t}(\cdot),\overline{y}(\cdot)\big) \in C([0,1],\R) \times C([0,1],\R^n)$ gibt mit
           \begin{eqnarray*}
           z_1(\cdot) &=& \tilde{F}_{1,t(\cdot)}\big(t_*(\cdot),y_*(\cdot),w_*(\cdot),v_*\big) \overline{t}(\cdot)
                          + \tilde{F}_{1,y(\cdot)}\big(t_*(\cdot),y_*(\cdot),w_*(\cdot),v_*\big) \overline{y}(\cdot), \\
           z_2(\cdot) &=& \tilde{F}_{2,t(\cdot)}\big(t_*(\cdot),v_*\big) \overline{t}(\cdot).
           \end{eqnarray*}
           Die zweite Gleichung besitzt die einfache Form $z_2(s)= \overline{t}(s)-\overline{t}(0)$, $s \in [0,1]$,
           und wir können unmittelbar durch $\overline{t}(s)=z_2(s)$ eine Lösung angeben.
           Setzen wir nun $\overline{t}(\cdot)$ in die erste Gleichung ein,
           so geht diese mit den abkürzenden Schreibweisen
           \begin{eqnarray*}
           a(s) &=& - \int_0^s v_* \cdot \big[ \varphi_\zeta\big(t_*(s),t_*(\sigma),y_*(\sigma),w_*(\sigma)\big) \overline{t}(s) \\
                & & \hspace*{2cm} + \varphi_t\big(t_*(s),t_*(\sigma),y_*(\sigma),w_*(\sigma)\big) \overline{t}(\sigma) \big] \, d\sigma, \\
           A(s,\sigma) &=& v_* \cdot \varphi_x\big(t_*(s),t_*(\sigma),y_*(\sigma),w_*(\sigma)\big)
           \end{eqnarray*}
           in die Form
           $$z_1(s) = a(s) + \overline{y}(s) - \overline{y}(0) - \int_0^s A(s,\sigma)\overline{y}(\sigma)\, d\sigma$$
           über.
           Setzen wir noch $\tilde{z}_1(s) = z_1(s)-a(s)$,
           erhält sie die Gestalt
           $$\tilde{z}_1(s) = \overline{y}(s) - \overline{y}(0) - \int_0^s A(s,\sigma)\overline{y}(\sigma)\, d\sigma$$
           einer Volterraschen Integralgleichung,
           welche nach Lemma \ref{LemmaDGL1} und Bemerkung \ref{BemDGL} eine Lösung $\overline{y}(\cdot)$ besitzt.
\item[(C)] Wir fassen die Zustände zu $t(\cdot),\, y(\cdot)$ zu $x(\cdot)$,
           die Steuerungen $w(\cdot),\, v$ zu $u(\cdot)$ und die Operatoren
           $\tilde{F}_1,\, \tilde{F}_2$ zu $F$ zusammen,
           setzen $\psi = v \cdot \varphi$
           und führen die allgemeinere Integralgleichung
           $$\tilde{F}\big(x(\cdot),u(\cdot)\big)(s) = x(s) - x(0)
                            + \int_0^s \psi\big(s,\sigma,x(s),x(\sigma),u(\sigma)\big) \, d\sigma, \quad s \in [0,1],$$
           mit dem Zustand $x(s)$ der oberen Integrationsgrenze im Integranden ein.
           Der Nachweis der Voraussetzungen (C) von Theorem \ref{SatzExtremalprinzipStark} sind im Anhang \ref{AnhangNV}
           in Lemma \ref{LemmaEigenschaftNadelvariationIG1} und in Lemma \ref{LemmaEigenschaftNadelvariationIG2} dargestellt.
\end{enumerate}

Zur Extremalaufgabe (\ref{ExtremalaufgabePMPIGLFZ}) definieren wir auf
$$C([0,1],\R) \times C([0,1],\R^n) \times L_\infty([0,1],\R^m) \times \R \times \R \times C_0^*([0,1],\R^n) \times C_0^*([0,1],\R) \times \R^{s_0} \times \R^{s_1}$$
die Lagrange-Funktion $\mathscr{L}=\mathscr{L}\big(t(\cdot),y(\cdot),w(\cdot),v,\lambda_0,y_1^*,y_2^*,l_0,l_1\big)$,
\begin{eqnarray*}
\mathscr{L} &=& \lambda_0 \tilde{J}\big(t(\cdot),y(\cdot),w(\cdot),v\big)+ \big\langle y_1^*, \tilde{F}_1\big(t(\cdot),y(\cdot),w(\cdot),v\big) \big\rangle
                + \big\langle y_2^*, \tilde{F}_2\big(t(\cdot),v\big) \big\rangle \\
            & & + l_0^T \tilde{H}_0\big(t(\cdot),y(\cdot)\big)+l_1^T \tilde{H}_1\big(t(\cdot),y(\cdot)\big).
\end{eqnarray*}
Ist $\big(t_*(\cdot),y_*(\cdot),w_*(\cdot),v_*\big)$ eine starke lokale Minimalstelle der Aufgabe (\ref{ExtremalaufgabePMPIGLFZ}),
dann existieren nach Theorem \ref{SatzExtremalprinzipStark}
nicht gleichzeitig verschwindende Lagrangesche Multiplikatoren 
$$\lambda_0 \geq 0, \quad y_1^* \in C_0^*([0,1],\R^n), \quad y_2^* \in C_0^*([0,1],\R), \quad l_i \in \R^{s_i},\; i=0,1,$$
derart,
dass folgende Bedingungen gelten:
\begin{enumerate}
\item[(a)] Die Lagrange-Funktion besitzt in $\big(t_*(\cdot),y_*(\cdot)\big)$ einen station"aren Punkt, d.\,h.
          \begin{equation}\label{SatzPMPLMRIGLFZ1}
          \left. \begin{array}{l}
          \mathscr{L}_{t(\cdot)}\big(t_*(\cdot),y_*(\cdot),w_*(\cdot),v_*,\lambda_0,y_1^*,y_2^*,l_0,l_1\big)=0, \\[2mm]
          \mathscr{L}_{y(\cdot)}\big(t_*(\cdot),y_*(\cdot),w_*(\cdot),v_*,\lambda_0,y_1^*,y_2^*,l_0,l_1\big)=0.
          \end{array} \right\}
          \end{equation}         
\item[(b)] Die Lagrange-Funktion erf"ullt bez"uglich $w(\cdot) \in L_\infty([0,1],U)$ und $v>0$ die Minimumbedingung
           \begin{eqnarray}
           \lefteqn{\mathscr{L}\big(t_*(\cdot),y_*(\cdot),w_*(\cdot),v_*,\lambda_0,y_1^*,y_2^*,l_0,l_1\big)} \nonumber \\
           \label{SatzPMPLMRIGLFZ2}
           &&= \min_{w(\cdot),\, v >0}
             \mathscr{L}\big(t_*(\cdot),y_*(\cdot),w(\cdot),v,\lambda_0,y_1^*,y_2^*,l_0,l_1\big).
           \end{eqnarray}
\end{enumerate}
Zu $y_1^* \in C_0^*([0,1],\R)$ und $y_2^* \in C_0^*([0,1],\R)$ gibt es ein eindeutig bestimmtes regul"ares Borelsches Vektorma"s $\mu_1$
bzw. Maß $\mu_2$ derart,
dass nach (\ref{SatzPMPLMRIGLFZ1}) f"ur alle $y(\cdot) \in C([0,1],\R^n)$ die Gleichung
\begin{eqnarray*}
0 &=& \lambda_0 \int_0^1 v_* \cdot \big\langle f_x\big(t_*(s),t_*(0),t_*(1),y_*(s),w_*(s)\big),y(s)\big\rangle \, ds \\
  & & + \int_0^1 \Big\langle y(s)-y(0) -\int_0^s v_* \cdot \varphi_x\big(t_*(s),t_*(\sigma),y_*(\sigma),w_*(\sigma)\big) y(\sigma) \, d\sigma, d\mu_1(s) \Big\rangle \\
  & & + \big\langle l_0, h_{0,x_0}\big(t_*(0),y_*(0)\big) y(0) \big\rangle + \big\langle l_1, h_{1,x_1}\big(t_*(1),y_*(1)\big) y(1) \big\rangle
\end{eqnarray*}
und für alle $t(\cdot) \in C([0,1],\R)$ die Gleichung
\begin{eqnarray*}
0 &=& \lambda_0 \int_0^1 v_* \cdot \Big[ f_t\big(t_*(s),t_*(0),t_*(1),y_*(s),w_*(s)\big) t(s) \\
  & & \hspace*{20mm} + f_{t_0} \big(t_*(s),t_*(0),t_*(1),y_*(s),w_*(s)\big) t(0) \\
  & & \hspace*{20mm} + f_{t_1} \big(t_*(s),t_*(0),t_*(1),y_*(s),w_*(s)\big) t(1) \Big]\, ds \\
  & & - \int_0^1 \Big\langle \int_0^s v_* \cdot \big[ \varphi_\zeta\big(t_*(s),t_*(\sigma),y_*(\sigma),w_*(\sigma)\big) t(s) \\
  & &  \hspace*{3.5cm} +  \varphi_t\big(t_*(s),t_*(\sigma),y_*(\sigma),w_*(\sigma)\big) t(\sigma) \big] \, d\sigma , d\mu_1(s) \Big\rangle \\
  & & +\int_0^1 \big(t(s)-t(0)\big) \, d\mu_2(s) 
      + \big\langle l_0, h_{0,t_0}\big(t_*(0),y_*(0)\big)\big\rangle t(0) + \big\langle l_1, h_{1,t_1}\big(t_*(1),y_*(1)\big) \big\rangle t(1)
\end{eqnarray*}
erf"ullt sind.
In diesen Gleichungen "andern wir die Integrationsreihenfolge in den Doppelintegralen.
Dann ergeben sich die Beziehungen
\begin{eqnarray*}
0 &=& \int_0^1 v_* \cdot \lambda_0 \Big\langle f_x\big(t_*(s),t_*(0),t_*(1),y_*(s),w_*(s)\big) \\
  & & \hspace*{3cm} -\int_s^1 \varphi_x^T\big(t_*(\sigma),t_*(s),y_*(s),w_*(s)\big) \, d\mu_1(\sigma), y(s) \Big\rangle \, ds \\
  & & \hspace*{-5mm} +\int_0^1 \langle y(s)-y(0), d\mu_1(s) \rangle
      + \big\langle l_0 , h_{0,x_0}\big(t_*(0),y_*(0)\big) y(0) \big\rangle + \big\langle l_1 , h_{1,x_1}\big(t_*(1),y_*(1)\big) y(1) \big\rangle
\end{eqnarray*}
für alle $y(\cdot) \in C([0,1],\R^n)$ und
\begin{eqnarray*}
0 &=& \int_0^1 v_* \cdot \Big[ \lambda_0 f_t\big(t_*(s),t_*(0),t_*(1),y_*(s),w_*(s)\big) \\
  & & \hspace*{3cm} -\int_s^1 \big\langle \varphi_t\big(t_*(\sigma),t_*(s),y_*(s),w_*(s)\big), d\mu_1(\sigma) \big\rangle \Big] t(s) \, ds \\
  & & - \int_0^1 \int_s^1 v_* \cdot \big\langle \varphi_\zeta\big(t_*(\sigma),t_*(s),y_*(s),w_*(s)\big) t(\sigma), d\mu_1(\sigma) \big\rangle  \, ds \\
  & & + \bigg[\int_0^1 v_* \cdot \lambda_0 f_{t_0}\big(t_*(s),t_*(0),t_*(1),y_*(s),w_*(s)\big) \, ds \bigg] \, t(0) \\ 
  & & + \bigg[\int_0^1 v_* \cdot \lambda_0 f_{t_1}\big(t_*(s),t_*(0),t_*(1),y_*(s),w_*(s)\big) \, ds \bigg] \, t(1) \\ 
  & & +\int_0^1 \big(t(s)-t(0)\big) \, d\mu_2(s) 
      + \big\langle l_0, h_{0,t_0}\big(t_*(0),y_*(0)\big)\big\rangle t(0) + \big\langle l_1, h_{1,t_1}\big(t_*(1),y_*(1)\big)\big\rangle t(1)
\end{eqnarray*}
für alle $t(\cdot) \in C([0,1],\R^n)$.
Wir setzen $\tilde{p}(s)=\displaystyle \int_s^1 d\mu_1(\sigma)$.
Damit folgt wegen der eindeutigen Darstellung eines stetigen linearen Funktionals im Raum $C([0,1],\R^n)$
\begin{equation} \label{IGLFZ7}
\left. \begin{array}{lll}
\tilde{p}'(s) &=& \displaystyle - v_* \cdot \bigg[ \lambda_0 f_x\big(t_*(s),t_*(0),t_*(1),y_*(s),w_*(s)\big) \\[3mm]
              & & \hspace*{2cm} \displaystyle  + \int_{s}^1 \varphi_x^T\big(t_*(\tau),t_*(s),y_*(s),w_*(s)\big)\, d\tilde{p}(\tau) \bigg], \\[3mm]
\tilde{p}(0) &=& h_{0,x_0}^T\big(t_*(0),y_*(0)\big)l_0, \quad \tilde{p}(1^-) \,=\, -h_{1,x_1}^T\big(t_*(1),y_*(1)\big)l_1, \quad \tilde{p}(1)=0.
\end{array} \right\}
\end{equation}
Weiterhin setzen wir $\tilde{q}(s)= \displaystyle \int_s^1 d\mu_2(\sigma)$.
So ergibt sich ferner aus der eindeutigen Darstellung eines stetigen linearen Funktionals im Raum $C([0,1],\R)$
\begin{equation} \label{IGLFZ8}
\left. \begin{array}{lll}
\tilde{q}'(s) &=& \displaystyle -v_* \bigg[ \lambda_0 f_t\big(t_*(s),t_*(0),t_*(1),y_*(s),w_*(s)\big) \\[3mm]
              & & \hspace*{2cm} \displaystyle
                  + \int_s^1 \big\langle \varphi_t\big(t_*(\tau),t_*(s),y_*(s),w_*(s)\big), d\tilde{p}(\tau)\big\rangle \bigg] \\[3mm]
              & & \displaystyle
                  -v_* \int_s^1 \big\langle \varphi_\zeta\big(t_*(\tau),t_*(s),y_*(s),w_*(s)\big), d\tilde{p}(\tau)\big\rangle, \\[3mm]
\tilde{q}(0) &=& \big\langle h_{0,t_0}\big(t_*(0),y_*(0)\big), l_0 \big\rangle \\[3mm]
             & &\displaystyle + \lambda_0 \int_0^1 v_* f_{t_0}\big(t_*(s),t_*(0),t_*(1),y_*(s),w_*(s)\big) \, ds \\[3mm]
\tilde{q}(1) &=& -\big\langle h_{1,t_1}\big(t_*(1),y_*(1)\big), l_1 \big\rangle \\[3mm]
             & & \displaystyle - \lambda_0 \int_0^1 v_* f_{t_1}\big(t_*(s),t_*(0),t_*(1),y_*(s),w_*(s)\big) \, ds.
\end{array} \right\}
\end{equation}
Die Beziehung (\ref{SatzPMPLMRIGLFZ2}) führt nach Vertauschen der Integrationsreihenfolge zur Ungleichung
\begin{eqnarray}
\lefteqn{v_* \cdot \int_0^1 \bigg[\lambda_0 f\big(t_*(s),t_*(0),t_*(1),y_*(s),w_*(s)\big)} \nonumber \\
&& \hspace*{2cm} + \int_s^1 \big\langle\varphi\big(t_*(\sigma),t_*(s),y_*(s),w_*(s)\big) , d\tilde{p}(\sigma) \big\rangle - \tilde{q}(s) \bigg] \, ds \nonumber \\
\lefteqn{ \leq v \cdot \int_0^1 \bigg[\lambda_0 f\big(t_*(s),t_*(0),t_*(1),y_*(s),w(s)\big)} \nonumber \\
&& \label{IGLFZ9} \hspace*{2cm} + \int_s^1 \big\langle\varphi\big(t_*(\sigma),t_*(s),y_*(s),w(s)\big) , d\tilde{p}(\sigma) \big\rangle -\tilde{q}(s) \bigg] \, ds,
\end{eqnarray}
die f"ur alle $w(\cdot) \in L_\infty([0,1],U)$ und $v>0$ gültig ist.

%% file: 5-53-Maximumprinzip.tex
\subsubsection{Notwendige Optimalitätsbedingungen}
Auf Basis der Beziehungen (\ref{IGLFZ7})--(\ref{IGLFZ9}) leiten wir das Pontrjaginsche Maximumprinzip der Aufgabe (\ref{IGLFZ1})--(\ref{IGLFZ4}) ab.
Dazu betrachten wir zu $t_*(\cdot)$ die inverse Funktion $s_*(\cdot)$, 
$$t_*(s)=t_{0*}+(t_{1*}-t_{0*})\cdot s=t_{0*}+v_*\cdot s, \; s \in [0,1], \qquad
  s_*(t) = \frac{t-t_{0*}}{t_{1*}-t_{0*}}, \; t \in [t_{0*},t_{1*}],$$
und folgen der Rückführung der Aufgabe des Abschnittes \ref{AbschnittFreieZeitSOP}: \\[2mm]
Wir f"uhren die Bezeichnungen $p(t) = \tilde{p}\big(s_*(t)\big)$ und $q(t) = \tilde{q}\big(s_*(t)\big)$ ein.
So lassen sich die Beziehungen (\ref{IGLFZ7}) in die Form
\begin{eqnarray*} 
\dot{p}(t) &=& \int_{t}^{t_{1*}} \varphi_x^T\big(\tau,t,x_*(t),u_*(t)\big)\, dp(\tau)
               + \lambda_0 f_x\big(t,t_{0*},t_{1*},x_*(t),u_*(t)\big), \\
p(t_{0*}) &=& h^T_{0,x_0} \big(t_{0*},x_*(t_{0*})\big) l_0, \quad
              p(t_{1*}^-) \;=\; -h^T_{1,x_1}\big(t_{1*},x_*(t_{1*})\big) l_1, \quad p(t_{1*})=0
\end{eqnarray*}
bringen und die Beziehungen (\ref{IGLFZ8}) erhalten die Gestalt
\begin{eqnarray*}
\dot{q}(t) &=& \int_{t}^{t_{1*}} \big\langle \varphi_t^T\big(\tau,t,x_*(t),u_*(t)\big), dp(\tau) \big\rangle
               + \lambda_0 f_t\big(t,t_{0*},t_{1*},x_*(t),u_*(t)\big), \\
           & & + \int_t^{t_{1*}} \big\langle \varphi_\zeta\big(\tau,t,x_*(t),u_*(t)\big), dp(\tau)\big\rangle, \\
q(t_{0*}) &=& \big\langle h_{0,t_0}\big(t_{0*},x_*(t_{0*})\big) , l_0 \big\rangle
              + \lambda_0 \int_{t_{0*}}^{t_{1*}} f_{t_0}\big(t,t_{0*},t_{1*},x_*(t),u_*(t)\big) \, dt, \\
q(t_{1*}) &=& - \big\langle h_{1,t_1}\big(t_{1*},x_*(t_{1*})\big) , l_1 \big\rangle 
              - \lambda_0 \int_{t_{0*}}^{t_{1*}} f_{t_1}\big(t,t_{0*},t_{1*},x_*(t),u_*(t)\big) \, dt.
\end{eqnarray*}
In der Adjungierten $q(\cdot)$ manifestiert sich der Einfluss des freien Anfangs- und Endzeitpunktes
in Form der Ableitungen $\varphi_\zeta$ und $f_{t_0},f_{t_1}$ nach den äußeren Zeitvariablen. \\[2mm] 
Die Bedingung (\ref{IGLFZ9}) führt bezüglich $v>0$ zu dem Zusammenhang
$$q(t)=\int_{t}^{t_{1*}} \big\langle \varphi\big(\tau,t,x_*(t),u_*(t)\big), dp(\tau) \big\rangle + \lambda_0 f\big(t,t_{0*},t_{1*},x_*(t),u_*(t)\big)$$
und liefert zudem die Maximumbedingung
\begin{eqnarray*}
\lefteqn{-\int_{t}^{t_{1*}} \big\langle \varphi\big(\tau,t,x_*(t),u_*(t)\big), dp(\tau) \big\rangle 
                               - \lambda_0 f\big(t,t_{0*},t_{1*},x_*(t),u_*(t)\big)} \\
&=& \max_{u \in U} \bigg(-\int_{t}^{t_{1*}} \big\langle \varphi\big(\tau,t,x_*(t),u\big), dp(\tau) \big\rangle
                           - \lambda_0 f\big(t,t_{0*},t_{1*},x_*(t),u\big)\bigg),
\end{eqnarray*}
die für fast alle $t \in [t_{0*},t_{1*}]$ gültig ist.

\newpage
Zur Formulierung der notwendigen Bedingungen der Aufgabe (\ref{IGLFZ1})--(\ref{IGLFZ4}) mit freiem Anfangs- und Endzeitpunkt
führen wir die erweiterte Pontrjagin-Funktion $H^{\mathcal{I}}$ ein:
$$H^{\mathcal{I}}(t,t_0,t_1,x,u,p(\cdot),\lambda_0) = -\int_t^{t_1} \big\langle \varphi(\tau,t,x,u), dp(\tau) \big\rangle - \lambda_0 f(t,t_0,t_1,x,u).$$

Zusammenfassend erhalten wir für die Aufgabe mit freiem Anfangs- und Endzeitpunkt:

\begin{theorem}[Pontrjaginsches Maximumprinzip] \label{SatzIGLFZ}
\index{Pontrjaginsches Maximumprinzip!Integral@-- Integralgleichungen} 
In der Aufgabe (\ref{IGLFZ1})--(\ref{IGLFZ4}) sei der Steuerungsprozess
$\big([t_{0*},t_{1*}],x_*(\cdot),u_*(\cdot)\big) \in \mathscr{B}^{\,\mathcal{F}}_{\rm adm} \cap \mathscr{B}^{\,\mathcal{F}}_{\rm Lip}$.
Ist $\big([t_{0*},t_{1*}],x_*(\cdot),u_*(\cdot)\big)$ ein starkes lokales Minimum der Aufgabe (\ref{IGLFZ1})--(\ref{IGLFZ4}),
dann existieren nicht gleichzeitig verschwindende Multiplikatoren $\lambda_0 \geq 0$,
$p(\cdot) \in W^1_\infty([t_{0*},t_{1*}),\R^n)$ und $l_i \in \R^{s_i}$, $i=0,1$, 
derart, dass 
\begin{enumerate}
\item[(a)] die adjungierte Gleichung
           \index{adjungierte Gleichung!Integral@-- Integralgleichungen}
           \begin{equation} \label{SatzPMPIGFZ1}
           \dot{p}(t) = -H^{\mathcal{I}}_x\big(t,t_{0*},t_{1*},x_*(t),u_*(t),p(\cdot),\lambda_0\big),
           \end{equation}
\item[(b)] die Transversalit"atsbedingungen
           \index{Transversalitätsbedingungen!Integral@-- Integralgleichungen}
           \begin{equation} \label{SatzPMPIGFZ2}
           p(t_{0*}) = h^T_{0,x_0} \big(t_{0*},x_*(t_{0*})\big) l_0, \quad
           p(t_{1*}^-) = -h^T_{1,x_1}\big(t_{1*},x_*(t_{1*})\big) l_1, \quad p(t_{1*})=0
           \end{equation} 
\item[(c)] und in fast allen Punkten $t\in [t_{0*},t_{1*}]$ die Maximumbedingung 
           \index{Maximumbedingung!Integral@-- Integralgleichungen}
           \begin{equation}\label{SatzPMPIGFZ3}
           \hspace*{-5mm} H^{\mathcal{I}}\big(t,t_{0*},t_{1*},x_*(t),u_*(t),p(\cdot),\lambda_0\big) 
            = \max_{u \in U} H^{\mathcal{I}}\big(t,t_{0*},t_{1*},x_*(t),u,p(\cdot),\lambda_0\big)
           \end{equation}
\end{enumerate}
erfüllt sind.
Außerdem definieren wir die erweiterte Hamilton-Funktion gemäß
$$\mathscr{H}^{\mathcal{I}}(t,x,p(\cdot),\lambda_0) = \max_{u \in U} H^{\mathcal{I}}(t,t_{0*},t_{1*},x,u,p(\cdot),\lambda_0).$$
\begin{enumerate}
\item[(d)] Es gelten außerdem für die Abbildung $t \to \mathscr{H}^{\mathcal{I}}\big(t,x_*(t),p(t),\lambda_0\big)$ die Beziehungen 
           \begin{eqnarray}
           \label{SatzPMPIGFZ4} \hspace*{-2cm} \frac{d}{dt} \mathscr{H}^{\mathcal{I}}\big(t,x_*(t),p(\cdot),\lambda_0\big) 
           &=& \int_t^{t_{1*}} \big\langle \varphi_t\big(\tau,t,x_*(t),u_*(t)\big), dp(\tau)\big\rangle \nonumber \\
           && \hspace*{-5cm} + \int_t^{t_{1*}} \big\langle \varphi_\zeta\big(\tau,t,x_*(t),u_*(t)\big), dp(\tau)\big\rangle
                               + \lambda_0 f_t\big(t,t_{0*},t_{1*},x_*(t),u_*(t)\big), \\
           \label{SatzPMPIGFZ5} \hspace*{-2cm} \mathscr{H}^{\mathcal{I}}\big(t_{0*},x_*(t_{0*}),p(\cdot),\lambda_0\big)
           &=& - \big\langle h_{0,t_0}\big(t_{0*},x_*(t_{0*})\big) , l_0 \big\rangle \nonumber \\
           && \hspace*{-0.5cm} - \lambda_0 \int_{t_{0*}}^{t_{1*}} f_{t_0}\big(t,t_{0*},t_{1*},x_*(t),u_*(t)\big) \, dt, \\
           \label{SatzPMPIGFZ6} \hspace*{-2cm} \mathscr{H}^{\mathcal{I}}\big(t_{1*},x_*(t_{1*}),p(\cdot),\lambda_0\big)
           &=& \big\langle h_{1,t_1}\big(t_{1*},x_*(t_{1*})\big) , l_1 \big\rangle \nonumber \\
           && \hspace*{-0.5cm} + \lambda_0 \int_{t_{0*}}^{t_{1*}} f_{t_1}\big(t,t_{0*},t_{1*},x_*(t),u_*(t)\big) \, dt.
           \end{eqnarray}
\end{enumerate}
\end{theorem}

%% file: 5-54-Werbestrategien.tex
\subsubsection{Optimaler Planungszeitraum einer Werbestrategie} \index{Werbestrategien} \label{AbschnittWerbung4}
In dem Beispiel (\ref{Werbung21})--(\ref{Werbung24}) einer optimalen Werbestrategie mit vorgegebenem Werbebudget im Abschnitt \ref{AbschnittWerbung2}
blieb die Frage nach dem Planungszeitraum offen.
Der Grund ist die Form der Aufgabe,
in der diese Frage nicht vernünftig beantwortet werden kann. \\[2mm]
Erfahrungsgemäß schwindet aus verschiedensten Gründen die Marktrelevanz im Laufe der Zeit.
Deswegen betrachten wir einen gewissen Mittelwert des Gewinnes.
Da der Bedeutungsschwund in der Zukunft in immer schnelleren Maße zunimmt,
machen wir den Ansatz eines gewichteten Mittels mit Hilfe des Faktors $e^{-\sigma T}$.
Zusätzlich nehmen wir an,
dass dieser Faktor das verfügbare Werbebudget (\ref{Werbung24}) beeinflusst,
da zugewiesene Finanzierungsmittel zeitlich nicht unbegrenzt in voller Höhe zur Verfügung stehen. \\[2mm]
Zur Unterscheidung der Bezeichnungen,
die wir im Beispiel (\ref{Werbung21})--(\ref{Werbung24}) verwendeten, versehen wir
in der vorliegenden Aufgabe das Zielfunktional, die Zustände, die Steuerung und die Adjungierten mit einem Symbol.
Die Aufgabe optimaler Werbestrategien mit gewichtetem Mittel des Gewinnes und mit freiem Endzeitpunkt besitzt die Form
\begin{eqnarray}
&& \label{Werbung41} \tilde{J}\big(\tilde{x}(\cdot),\tilde{u}(\cdot)\big)
   = e^{-\sigma T} \int_0^T e^{-\varrho t} \big[\pi \tilde{x}(t)-\tilde{u}(t)\big] \, dt \to \sup, \\
&& \label{Werbung42} \tilde{x}(t)=x_0+\int_0^t e^{-\delta(t-s)}
   \big[f\big(\tilde{u}(s)\big)-\alpha \tilde{x}(s)\big] \, ds, \quad t \in [0,T], \\
&& \label{Werbung43} \tilde{z}(t)=\int_0^t e^{-\varrho s} \tilde{u}(s) \, ds, \quad t \in [0,T], \quad e^{-\sigma T} \tilde{z}(T)=Z>0, \\
&& \label{Werbung44} x_0 \in [0,1],\quad \tilde{u}(t) \geq 0, \quad \varrho, \sigma >0, \quad \alpha, \delta \geq 0, \quad T >0 \mbox{ frei}.
\end{eqnarray}
In den Integranden fließt die äußere Zeit $T$ ein.
Deswegen besitzt der Integrand $f$ nach Übergang zu einem Minimierungsproblem die Gestalt
$$f(t,T,\tilde{x},\tilde{u})= -e^{-\sigma T} e^{-\varrho t} (\pi \tilde{x}-\tilde{u}).$$
Bezüglich dem Zustand $z(\cdot)$ gehen wir wieder zur Differentialgleichung $\dot{z}(t)=e^{-\varrho t}u(t)$ über.
Damit erhalten wir für die Aufgabe (\ref{Werbung41})--(\ref{Werbung44}) die erweiterte Pontrjagin-Funktion
\begin{eqnarray*}
    H^{\mathcal{I}}(t,T,\tilde{x},\tilde{z},\tilde{u},\tilde{q}(\cdot),\tilde{r},\lambda_0)
&=& -\int_t^T e^{-\delta(\tau-t)} \big[f(\tilde{u})-\alpha \tilde{x}\big] \, d\tilde{q}(\tau) \\
& & \hspace*{10mm} + \tilde{r} [e^{-\varrho t} \tilde{u}] + \lambda_0 e^{-\sigma T}e^{-\varrho t} [\pi \tilde{x}-\tilde{u}].
\end{eqnarray*}
Die notwendigen Optimalitätsbedingungen in Theorem \ref{SatzIGLFZ} ergeben mit $\lambda_0=1$: \\[2mm]
Nach (\ref{SatzPMPIGFZ1}) und (\ref{SatzPMPIGFZ2}) gelten die adjungierten Gleichungen
\begin{eqnarray*}
-H^{\mathcal{I}}_{\tilde{x}}\big(t,...,1\big)
&=& \dot{\tilde{q}}(t)=\int_t^{T_*} -\alpha e^{-\delta(\tau-t)} \, d\tilde{q}(\tau) -\pi e^{-\sigma T_*} e^{-\varrho t}, \\
-H^{\mathcal{I}}_{\tilde{z}}\big(t,...,1\big)
&=& \dot{\tilde{r}}(t) =0
\end{eqnarray*}
zu den Transversalitätsbedingungen
$$\tilde{q}(T_*^-)=\tilde{q}(T_*)=0, \qquad \tilde{r}(T_*)=-l_1 e^{-\sigma T_*},$$
wobei die Randbedingung die Form $h_1\big(T,\tilde{z}(T)\big)=e^{-\sigma T}\tilde{z}(T)-Z=0$ besitzt. \\
Mit dem Verweis auf die Adjungierten $q(\cdot)$ im Abschnitt \ref{AbschnittWerbung2} ergibt sich für $\tilde{q}(\cdot)$
die Darstellung $\tilde{q}(t)=e^{-\sigma T_*}q(t)$, denn es gilt offenbar
\begin{eqnarray*}
\dot{\tilde{q}}(t) = e^{-\sigma T_*}\cdot \dot{q}(t)
&=& e^{-\sigma T_*} \cdot \bigg[\int_t^{T_*} -\alpha e^{-\delta(\tau-t)} \dot{q}(\tau) \, d\tau -\pi e^{-\varrho t} \bigg] \\
&=& \int_t^{T_*} -\alpha e^{-\delta(\tau-t)} \dot{\tilde{q}}(\tau) \, d\tau - e^{-\sigma T_*} \cdot \pi e^{-\varrho t}.
\end{eqnarray*}
Aus der Maximumbedingung (\ref{SatzPMPIGFZ3}) folgt für fast alle $t \in [0,T_*]$ die Gültigkeit der Aufgabe
$$\max_{u \geq 0} \bigg[-\int_t^{T_*} e^{-\delta(\tau-t)} \dot{\tilde{q}}(\tau) \, d\tau \cdot f(\tilde{u})
  + \tilde{r}(t) e^{-\varrho t} \cdot \tilde{u}- e^{-\sigma T_*}e^{-\varrho t} \cdot \tilde{u}\bigg].$$
Unter Verwendung der Darstellung von $\dot{\tilde{q}}(\cdot)$ und dem Zusammenhang $\tilde{q}(t)=e^{-\sigma T_*}q(t)$,
sowie $\tilde{r}(t) \equiv -l_1 e^{-\sigma T_*}$
erhalten wir mit den Beziehungen im Abschnitt \ref{AbschnittWerbung2} für die optimale Steuerung $\tilde{u}(\cdot)$ die implizite Darstellung
\begin{eqnarray*}
f'\big(\tilde{u}_*(t)\big)
&=& e^{-\sigma T_*} (1+l_1) e^{-\varrho t} \cdot \frac{\alpha}{\dot{\tilde{q}}(t)+e^{-\sigma T_*} \pi e^{-\varrho t}}
    = (1+l_1) e^{-\varrho t} \cdot \frac{\alpha}{\dot{q}(t)+ \pi e^{-\varrho t}} \\
&=& (1+l_1) \cdot \frac{\alpha + \delta + \varrho}{\pi ( 1-e^{-(\alpha + \delta + \varrho)(T_*-t)})}.
\end{eqnarray*}
Abschließend gehen wir auf den Teil (d) von Theorem \ref{SatzIGLFZ} ein:
Die Aufgabe beinhaltet den Integranden $f(t,T,\tilde{x},\tilde{u}) = - e^{-\sigma T}e^{-\varrho t} [\pi \tilde{x}-\tilde{u}]$,
sowie bezüglich der Zustände $\tilde{x}, \tilde{z}$ die rechten Seiten
$\varphi_1(\zeta,t,\tilde{x},\tilde{u}) = e^{-\delta(\zeta-t)} [f(\tilde{u})-\alpha \tilde{x}]$
 und $\varphi_2(t,\tilde{u}) = e^{-\varrho t} \tilde{u}$. \\
Diesen Abbildungen entnimmt man unmittelbar die verschiedenen Zeitableitungen
\begin{eqnarray*}
f_t(t,T,\tilde{x},\tilde{u}) &=& -\varrho f(t,T,\tilde{x},\tilde{u}), \quad
f_T(t,T,\tilde{x},\tilde{u}) = -\sigma f(t,T,\tilde{x},\tilde{u}), \\
\varphi_{1,t}(\zeta,t,\tilde{x},\tilde{u}) &=& -\varphi_{1,\zeta}(\zeta,t,\tilde{x},\tilde{u})
                                             = \delta \varphi_1(\zeta,t,\tilde{x},\tilde{u}), \quad
\varphi_{2,t}(t,\tilde{u}) = -\varrho \varphi_2(t,\tilde{u}).
\end{eqnarray*}
Daraus ergibt sich nach (\ref{SatzPMPIGFZ4}) für die erweiterte Hamilton-Funktion
\begin{eqnarray*}
\frac{d}{dt} \mathscr{H}^{\mathcal{I}}\big(t,\tilde{x}_*(t),\tilde{z}_*(t),\tilde{q}(\cdot),\tilde{r},1\big)
&=& - \varphi_{2,t}\big(t,\tilde{u}_*(t)\big) \tilde{r}(t) + f_t\big(t,T_*,\tilde{x}_*(t),\tilde{u}_*(t)\big) \\
&& \hspace*{-2cm} = -\varrho e^{-\varrho t}\tilde{u}_*(t)\cdot l_1e^{-\sigma T_*} - \varrho f\big(t,T_*,\tilde{x}_*(t),\tilde{u}_*(t)\big).
\end{eqnarray*}
Im Weiteren verwenden wir zur Abkürzung das Zielfunktional $\tilde{J}$, d.\,h.
$$-\tilde{J}\big(x_*(\cdot),u_*(\cdot)\big) = \int_0^{T_*} f\big(t,T_*,\tilde{x}_*(t),\tilde{u}_*(t)\big) \, dt.$$
Die Randbedingungen (\ref{SatzPMPIGFZ5}) und (\ref{SatzPMPIGFZ6}) liefern die Beziehungen
\begin{eqnarray*}
\mathscr{H}^{\mathcal{I}}\big(0,\tilde{x}_*(0),\tilde{z}_*(0),\tilde{q}(\cdot),\tilde{r},1\big)
&=& 0, \\
\mathscr{H}^{\mathcal{I}}\big(T_*,\tilde{x}_*(T_*),\tilde{z}_*(T_*),\tilde{q}(\cdot),\tilde{r},1\big)
&=& h_{1,T}\big(T_*,\tilde{z}_*(T_*)\big) l_1 + \sigma \tilde{J}\big(x_*(\cdot),u_*(\cdot)\big) \\
&=& - \sigma l_1 e^{-\sigma T_*} \tilde{z}_*(T_*) + \sigma \tilde{J}\big(x_*(\cdot),u_*(\cdot)\big).
\end{eqnarray*}
Außerdem gilt für den Randwert in $t=T_*$:
\begin{eqnarray*}
\mathscr{H}^{\mathcal{I}}\big(T_*,\tilde{x}_*(T_*),\tilde{z}_*(T_*),\tilde{q}(\cdot),\tilde{r},1\big)
&=& \int_0^{T_*} \frac{d}{dt} \mathscr{H}^{\mathcal{I}}\big(t,\tilde{x}_*(t),\tilde{z}_*(t),\tilde{q}(\cdot),\tilde{r},1\big) \, dt \\
&=& -\varrho l_1 e^{-\sigma T_*} \tilde{z}_*(T_*) +  \varrho \tilde{J}\big(x_*(\cdot),u_*(\cdot)\big).
\end{eqnarray*} 
Damit ergibt sich aus der erweiterten Hamilton-Funktion $\mathscr{H}^{\mathcal{I}}$ die Bedingung
$$\tilde{J}\big(x_*(\cdot),u_*(\cdot)\big)=l_1 Z.$$
Zusammenfassend liefern die Optimalitätsbedingungen von Theorem \ref{SatzIGLFZ}
eine Familie von Steuerungen $\tilde{u}(\cdot \,;T,l_1)$ der Parameter $T$ und $l_1$,
welche gemäß
$$f'\big(\tilde{u}(t;T,l_1)\big) = (1+l_1) \cdot \frac{\alpha + \delta + \varrho}{\pi ( 1-e^{-(\alpha + \delta + \varrho)(T-t)})}$$
über $[0,T]$ streng monoton fallend sind und für die im Endpunkt $\tilde{u}(T;T,l_1)=0$ gilt. \\
Innerhalb dieser Familie selektiert der Parameter $T$ diejenigen Elemente, für die
$$e^{-\sigma T} \int_0^T e^{-\varrho t} \tilde{u}(t;T,l_1) \, dt =Z$$
erfüllt ist.
Der Multiplikator $l_1$ wiederum bedeutet anhand der Beziehung
$$l_1 = \frac{\tilde{J}\big(x_*(\cdot),u_*(\cdot)\big)}{Z}$$
die Rentabilität des gewichteten mittleren Gewinns zu den eingesetzten Mitteln.
Damit ist $l_1$ ein Bewertungsmaß für die Wirtschaftlichkeit der Unternehmung.
Also wird außerdem innerhalb der Familie $\tilde{u}(\cdot \,;T,l_1)$ diejenige Steuerung $\tilde{u}_*(\cdot)$ gesucht,
mit der sich sinnvollerweise die maximale Rentabilität realisieren lässt. \\[2mm]
Die Rentabilität $l_1=\tilde{J}/Z$ und die Werbeintensität $\tilde{u}_*(\cdot)$ können nicht beliebig groß ausfallen,
denn es lassen sich unmittelbar obere Schranken angeben:
Das Zielfunktional besitzt für $x(t) \equiv 1$ und $u(t) \equiv 0$ die obere Beschränkung
$$\tilde{J}\big(x(\cdot),u(\cdot)\big) < \max_{T>0} \bigg[ e^{-\sigma T} \int_0^T e^{-\varrho t} \pi \, dt \bigg] = \tilde{J}_{\max}.$$
Unter der Annahme einer nichtnegativen Rentabilität $l_1$ ergibt sich
$$f'(\tilde{u}) > \frac{\alpha + \delta + \varrho}{\pi} =f'(\tilde{u}_{\max}).$$
Damit beschließen wir unsere Diskussionen optimaler Werbestrategien.  \hfill $\square$

%% file: 5-55-Standardaufgabe.tex
\subsubsection{Zurück zur Standardaufgabe mit freiem Anfangs- und Endzeitpunkt}
Durch das Weglassen der ``äußeren'' Zeitvariablen $\zeta$ in der Volterraschen Integralgleichung (\ref{IGLFZ2})
gelangen wir wieder zur Standardaufgabe.
Allerdings liegt durch das Auftreten der Zeitvariablen $t_0,t_1$ im Integranden eine gewisse Verallgemeinerung
zu Abschnitt \ref{AbschnittFreieZeit} vor:
\begin{equation} \label{IGLFZS}
\left. \begin{array}{lll}
&& \displaystyle J\big(x(\cdot),u(\cdot)\big) = \int_{t_0}^{t_1} f\big(t,t_0,t_1,x(t),u(t)\big) \, dt \to \inf, \\[3mm]
&& \dot{x}(t) = \varphi\big(t,x(t),u(t)\big), \\[2mm]
&& h_0\big(t_0,x(t_0)\big)=0, \qquad h_1\big(t_1,x(t_1)\big)=0, \\[2mm]
&& u(t) \in U \subseteq \R^m, \quad U \not= \emptyset.
\end{array} \right\}
\end{equation}

Als Spezialfall der Aufgabe (\ref{IGLFZ1})--(\ref{IGLFZ4}) übernehmen wir sämtliche bereits getroffene Voraussetzungen und Definitionen.
Die erweiterte Pontrjagin-Funktion $H^{\mathcal{S}}$ der Aufgabe (\ref{IGLFZS}) besitzt die Gestalt:
$$H^{\mathcal{S}}(t,t_0,t_1,x,u,p,\lambda_0) = \langle p \,,\, \varphi(t,x,u) \rangle - \lambda_0 f(t,t_0,t_1,x,u).$$
Bei Weglassen der Ableitungen $\varphi_\zeta$ erhalten wir unmittelbar:

\begin{theorem}[Pontrjaginsches Maximumprinzip] \label{SatzIGLFZS}
\index{Pontrjaginsches Maximumprinzip!Integral@-- Integralgleichungen} 
In der Aufgabe (\ref{IGLFZS}) sei der Steuerungsprozess
$\big([t_{0*},t_{1*}],x_*(\cdot),u_*(\cdot)\big) \in \mathscr{B}^{\,\mathcal{F}}_{\rm adm} \cap \mathscr{B}^{\,\mathcal{F}}_{\rm Lip}$.
Ist $\big([t_{0*},t_{1*}],x_*(\cdot),u_*(\cdot)\big)$ ein starkes lokales Minimum der Aufgabe (\ref{IGLFZS}),
dann existieren nicht gleichzeitig verschwindende Multiplikatoren $\lambda_0 \geq 0$,
$p(\cdot) \in W^1_\infty([t_{0*},t_{1*}],\R^n)$ und $l_i \in \R^{s_i}$, $i=0,1$, 
derart, dass 
\begin{enumerate}
\item[(a)] die adjungierte Gleichung
           \index{adjungierte Gleichung!Integral@-- Integralgleichungen}
           \begin{equation} \label{SatzPMPIGFZS1}
           \dot{p}(t) = -H^{\mathcal{S}}_x\big(t,t_{0*},t_{1*},x_*(t),u_*(t),p(t),\lambda_0\big),
           \end{equation}
\item[(b)] die Transversalit"atsbedingungen
           \index{Transversalitätsbedingungen!Integral@-- Integralgleichungen}
           \begin{equation} \label{SatzPMPIGFZS2}
           p(t_{0*}) = h^T_{0,x_0} \big(t_{0*},x_*(t_{0*})\big) l_0, \quad
           p(t_{1*}) = -h^T_{1,x_1}\big(t_{1*},x_*(t_{1*})\big) l_1
           \end{equation} 
\item[(c)] und in fast allen Punkten $t\in [t_{0*},t_{1*}]$ die Maximumbedingung 
           \index{Maximumbedingung!Integral@-- Integralgleichungen}
           \begin{equation}\label{SatzPMPIGFZS3}
           \hspace*{-5mm} H^{\mathcal{S}}\big(t,t_{0*},t_{1*},x_*(t),u_*(t),p(t),\lambda_0\big) 
            = \max_{u \in U} H^{\mathcal{S}}\big(t,t_{0*},t_{1*},x_*(t),u,p(t),\lambda_0\big)
           \end{equation}
\end{enumerate}
erfüllt sind.
Außerdem definieren wir die erweiterte Hamilton-Funktion gemäß
$$\mathscr{H}^{\mathcal{S}}(t,x,p,\lambda_0) = \max_{u \in U} H^{\mathcal{S}}(t,t_{0*},t_{1*},x,u,p,\lambda_0).$$
\begin{enumerate}
\item[(d)] Es gelten außerdem für die Abbildung $t \to \mathscr{H}^{\mathcal{S}}\big(t,x_*(t),p(t),\lambda_0\big)$ die Beziehungen 
           \begin{eqnarray}
           \label{SatzPMPIGFZS4} \hspace*{-2cm} \frac{d}{dt} \mathscr{H}^{\mathcal{S}}\big(t,x_*(t),p(t),\lambda_0\big) 
           &=& -H_t^{\mathcal{S}}\big(t,x_*(t),p(t),\lambda_0\big) , \\
           \label{SatzPMPIGFZS5} \hspace*{-2cm} \mathscr{H}^{\mathcal{S}}\big(t_{0*},x_*(t_{0*}),p(t_{0*}),\lambda_0\big)
           &=& - \big\langle h_{0,t_0}\big(t_{0*},x_*(t_{0*})\big) , l_0 \big\rangle \nonumber \\
           && \hspace*{-0.5cm} - \lambda_0 \int_{t_{0*}}^{t_{1*}} f_{t_0}\big(t,t_{0*},t_{1*},x_*(t),u_*(t)\big) \, dt, \\
           \label{SatzPMPIGFZS6} \hspace*{-2cm} \mathscr{H}^{\mathcal{S}}\big(t_{1*},x_*(t_{1*}),p(t_{1*}),\lambda_0\big)
           &=& \big\langle h_{1,t_1}\big(t_{1*},x_*(t_{1*})\big) , l_1 \big\rangle \nonumber \\
           && \hspace*{-0.5cm} + \lambda_0 \int_{t_{0*}}^{t_{1*}} f_{t_1}\big(t,t_{0*},t_{1*},x_*(t),u_*(t)\big) \, dt.
           \end{eqnarray}
\end{enumerate}
\end{theorem}

\begin{beispiel} \index{Werbestrategien}
{\rm Mit der Umformung der Volterraschen Integralgleichung in eine Differentialgleichung in Beispiel \ref{AbschnittWerbung1}
ergibt sich für die Aufgabe (\ref{Werbung41})--(\ref{Werbung44}) die Form
\begin{equation*} 
\left. \begin{array}{lll}
&& \displaystyle \tilde{J}\big(\tilde{x}(\cdot),\tilde{u}(\cdot)\big)
   = e^{-\sigma T} \int_0^T e^{-\varrho t} \big[\pi \tilde{x}(t)-\tilde{u}(t)\big] \, dt \to \sup, \\[3mm]
&& \dot{\tilde{x}}(t) = -(\alpha+\delta) \tilde{x}(t) + \delta x_0 + f\big(\tilde{u}(t)\big), \quad \tilde{x}(0)=x_0, \\[2mm]
&& \dot{\tilde{z}}(t)=e^{-\varrho t} \tilde{u}(t), \quad \tilde{z}(0)=0, \quad e^{-\sigma T} \tilde{z}(T)=Z>0, \\[2mm]
&& x_0 \in [0,1],\quad \tilde{u}(t) \geq 0, \quad \varrho >0, \quad \alpha, \delta \geq 0.
\end{array} \right\}
\end{equation*}
Wir setzen $\tilde{p}(t)=e^{-\sigma T_*} \cdot p(t)$ mit der Adjungierten $p(\cdot)$,
die den Gleichungen (\ref{Werbung8}) und (\ref{Werbung9}) genügt.
Dann gelten für die Adjungierte $\tilde{p}(\cdot)$ bezüglich dem Zustand $\tilde{x}(\cdot)$:
\begin{eqnarray*}
\dot{\tilde{p}}(t) &=& e^{-\sigma T_*} \cdot \dot{\tilde{p}}(t) = (\alpha+\delta)\tilde{p}(t)-\pi e^{-\varrho t}e^{-\sigma T_*},
                       \quad \tilde{p}(T_*)=0, \\
\tilde{p}(t) &=& e^{-\sigma T_*} \bigg[
                 -\frac{\pi \cdot e^{-(\alpha + \delta + \varrho)T}}{\alpha + \delta + \varrho} \cdot e^{(\alpha+\delta)t}
                 +\frac{\pi}{\alpha + \delta + \varrho} \cdot e^{-\varrho t} \bigg].
\end{eqnarray*}
Für die Adjungierte $\tilde{r}(\cdot)$ bezüglich dem Zustand $\tilde{z}(\cdot)$ ergibt sich unmittelbar
$$\dot{\tilde{r}}(t) \equiv 0, \qquad \tilde{r}(t) \equiv -l_1 e^{-\sigma T_*}.$$
Die Maximumbedingung führt für $\tilde{u}_*(\cdot)$ zu der Beziehung
$$f'\big(\tilde{u}_*(t)\big)
  = \frac{e^{-\sigma T_*} e^{-\varrho t} - \tilde{r}(t) e^{-\varrho t}}{\tilde{p}(t)}
  = \frac{e^{-\sigma T_*} e^{-\varrho t}(1+ l_1)}{e^{-\sigma T_*} \cdot p(t)}
  = \frac{(1+l_1)(\alpha + \delta + \varrho)}{\pi ( 1-e^{-(\alpha + \delta + \varrho)(T_*-t)})},$$
die wir bereits in Abschnitt \ref{AbschnittWerbung4} für die äquivalente Aufgabe mit einer Volterraschen Integralgleichung
identifizieren konnten. \\[2mm]
Die notwendigen Optimalitätsbedingungen im Teil (d) von Theorem \ref{SatzIGLFZS} führen
auf die gleiche Weise wie im Abschnitt \ref{AbschnittWerbung4} auf
$$\tilde{J}\big(x_*(\cdot),u_*(\cdot)\big)=l_1 Z,$$
wonach der Multiplikator $l_1$ die maximal mögliche Rentabilität widerspiegelt. \hfill $\square$}
\end{beispiel}

%% file: I-Masstheorie.tex
\section{Grundz\"uge der Maß- und Integrationstheorie} \label{AnhangMasstheorie}
Unsere Darstellungen zur Ma"s- und Integrationstheorie basieren auf 
Bauer \cite{Bauer}, Elstrodt \cite{Elstrodt}, Heuser \cite{HeuserII} und Natanson \cite{Natanson}.
Die Eigenschaften der Funktionen von beschr"ankter Variation und des Stieltjes-Integral sind Heuser \cite{HeuserI}
und Natanson \cite{Natanson} entnommen.

\subsection{Maßtheoretische Grundlagen}
Sei $M$ eine Menge, und sei $\Sigma$ ein System von Teilmengen von $M$ mit den Eigenschaften
$$\emptyset \in \Sigma, \qquad A \in \Sigma \Rightarrow M \setminus A \in \Sigma, \qquad
  A_1,A_2, ... \in \Sigma \Rightarrow \bigcup\limits_{i=1}^\infty A_i \in \Sigma,$$
dann hei"st $\Sigma$ eine $\sigma$-Algebra\index{Sigma-Algebra@$\sigma$-Algebra} und jede Menge $A \in \Sigma$ messbar. \\
Die offenen Teilmengen der Menge $M$ erzeugen die $\sigma$-Algebra der Borelmengen\index{Borelmengen}
$\mathscr{B}(M)$.
D.\,h. $\mathscr{B}(M)$ ist die kleinste $\sigma$-Algebra, die $M$ umfasst. \\
Eine Abbildung $\mu: \Sigma \to [0,\infty]$ hei"st ein Ma"s\index{Maß}, wenn $\mu(\emptyset) = 0$ und
$\mu$ $\sigma$-additiv\index{Maß!Sigmaadditiv@--, $\sigma$-additives} ist, d.\,h.,
wenn $ \displaystyle \mu\bigg( \bigcup\limits_{i=1}^\infty A_i \bigg) = \sum_{i=1}^\infty \mu(A_i)$ f"ur paarweise disjunkte $A_1,A_2, ... \in \Sigma$ gilt.
Das Ma"s $\mu$ nennen wir $\sigma$-endlich\index{Maß!Sigmaendlich@--, $\sigma$-endliches},
wenn es in $M$ eine Folge $A_1,A_2, ...$ von Teilmengen gibt mit
$\mu(A_i) < \infty$ f"ur alle $i$ und $A_1 \cup A_2 \cup ...= M$. 
Eine Abbildung $\mu: \Sigma \to \R $ hei"st ein signiertes Ma"s\index{Maß!signiertes@--, signiertes}, falls $\mu$ $\sigma$-additiv ist.
Ein Ma"s ist demnach im Unterschied zu signierten Ma"sen eine Abbildung mit ausschlie"slich nichtnegativen Werten.
Sei $\mu$ ein signiertes Ma"s auf $\Sigma$ und $A \in \Sigma$.
Wir setzen $\displaystyle |\mu|(A) = \sup \sum_{k=1}^n |\mu(A_k)|,$
wobei das Supremum "uber alle endlichen Zerlegungen $A = A_1 \cup ... \cup A_n$ mit paarweise disjunkten $A_k \in \Sigma$ zu bilden ist.
$|\mu|(A)$ hei"st totale Variation\index{Maß!totale Variation@-- totale Variation eines}
\index{totale Variation, einer Funktion!eines@-- eines Maßes} von $\mu$ auf $A$.
Ist $\mu$ ein signiertes Ma"s, so definiert $A \to |\mu|(A)$ ein endliches Ma"s.
Ferner nennen wir $A \in \Sigma$ mit $\mu(A)=0$ eine $\mu$-Nullmenge\index{muNullmenge@$\mu$-Nullmenge}. \\
Ein Ma"s $\mu$ hei"st von innen bzw. von au"sen regul"ar, falls f"ur alle $A \in \Sigma$
$$\mu(A) = \sup\{\mu(C) \, | \, C\subset A, C \mbox{ kompakt} \}, \quad\mbox{bzw.}\quad
  \mu(A) = \inf\{\mu(O) \, | \, A \subset O, O \mbox{ offen} \}$$
gilt.
Ein von innen und au"sen regul"ares Ma"s hei"st regul"ar\index{Maß!regul"ares@--, regul"ares}.
Ein signiertes Ma"s hei"st regul"ar, wenn seine totale Variation $|\mu|$ regul"ar ist.
Ein (signiertes) Ma"s $\mu$ auf der $\sigma$-Algebra $\mathscr{B}(M)$ der Menge $M$ hei"st (signiertes)
Borelsches Ma"s\index{Borelsches Maß}\index{Maß!Borelsches@--, Borelsches} "uber $M$.

\begin{satz}[Hahn-Jordan-Zerlegung] 
Sei $\mu$ ein signiertes Ma"s über der $\sigma$-Algebra $\Sigma$.
\begin{enumerate}
\item[(a)] Es existiert ein bis auf eine $|\mu|$-Nullmenge eindeutig bestimmtes $A^+ \in \Sigma$ mit
           $$0 \leq \mu(A) \leq \mu(A^+) \quad\mbox{ und }\quad 0 \geq \mu(B) \geq \mu(A^-)$$
           f"ur alle $A, B \in \Sigma$ mit $A \subseteq A^+$ und $B \subseteq A^-:=M \setminus A^+$.
\item[(b)] $\mu^+ = \chi_{A^+} \mu$ und  $\mu^- = \chi_{A^-} \mu$ sind endliche Ma"se mit
           $\mu(A) = \mu^+(A) - \mu^-(A)$ und $|\mu|(A) = \mu^+(A) + \mu^-(A)$ f"ur alle $A \in \Sigma$.
\end{enumerate}
Die Zerlegungen $A=A^+ \cup A^-$ und $\mu = \mu^+ - \mu^-$ hei"sen Hahn- bzw. Jordan-Zerlegung.
\end{satz}

\subsection{Integration bez\"uglich eines Borelschen Maßes} \label{AbschnittIntegrationMass}
Es sei $I \subseteq \R$ ein Intervall, 
d.\,h. eine beschr"ankte oder unbeschr"ankte konvexe Teilmengen von $\R$,
und $\mathscr{B}$ die Borelsche $\sigma$-Algebra auf $I$.
Die Funktion $f:I \to \R$ ist eine Treppenfunktion 
wenn $f$ nur endliche viele Werte annimmt.
D.\,h. gibt es paarweise disjunkte messbare Mengen $A_1,...,A_k \subseteq I$ und Punkte $y_1,...,y_k \in \R$ mit
$A_1 \cup ... \cup A_k = I$ und 
$$f(t)=\sum_{i=1}^k y_i \chi_{A_i}(t), \quad t \in I.$$
Es sei $\mu$ ein Ma"s "uber der $\sigma$-Algebra $\mathscr{B}$.
Eine nichtnegative reelle Funktion $f$ ist genau dann me"sbar,
wenn es eine monoton wachsende Folge $\{f_k\}$ von nichtnegativen Treppenfunktionen gibt,
die punktweise gegen $f$ konvergiert.
Jede beschr"ankte me"sbare Funktion ist gleichm"a"siger Grenzwert von Treppenfunktionen.
Eine nichtnegative reelle Funktion $f$ hei"st $\mu$-integrierbar\index{muintegrierbar@$\mu$-integrierbar},
falls eine monoton wachsende Folge $\{ f_k \}$ nichtnegativer Treppenfunktionen existiert,
die punktweise gegen $f$ konvergiert. Die Zahl 
$$\int_I f \, d\mu = \int_I f(t) \, d\mu(t) = \lim_{k \to \infty} \int_I f_k(t) \, d\mu(t) = \lim_{k \to \infty} \sum_{i=1}^k y_{i,k} \cdot \mu(A_{i,k}),$$
welche unabh"angig von der Darstellung mittels $\{f_k\}$ ist, hei"st das $\mu$-Integral\index{muintegral@$\mu$-Integral} von $f$. \\[1mm]
Eine me"sbare reelle Funktion $f$ hei"st $\mu$-integrierbar,
falls ihr Positivteil $f^+$ und ihr Negativteil $f^-$ $\mu$-integrierbar sind und
die $\mu$-Integrale $\displaystyle \int_I f^+ \, d\mu, \, \int_I f^- \, d\mu$ reelle Zahlen sind.
Dann hei"st $\displaystyle \int_I f \, d\mu = \int_I f^+ \, d\mu - \int_I f^- \, d\mu$ das $\mu$-Integral von $f$.
Es sei $\mu$ ein signiertes Ma"s, $\mu = \mu^+ - \mu^-$ seine Jordan-Zerlegung und die reelle Funktion $f$ me"sbar.
Dann hei"st $f$ $\mu$-integrierbar,
falls $f$ $|\mu|$-integrierbar ist und wir setzen $\displaystyle\int_I f \, d\mu = \int_I f \, d\mu^+ - \int_I f \, d\mu^-$.
Für das Lebesgue-Ma"s sagen wir integrierbar und schreiben
$\displaystyle\int_I f \, dt$.\index{Lebesgue-Integral} \\[2mm]
Ein innerer Punkt $t_0 \in I$ hei"st Lebesguescher Punkt\index{Lebesguescher Punkt} der Funktion $f$, wenn
$$\lim_{h \to 0^+} \frac{1}{h} \int_0^h |f(t_0 \pm t)-f(t_0)| \, dt =0.$$

\begin{satz}
Ist $f$ "uber $I$ integrierbar, so ist fast jeder Punkt ein Lebesguescher Punkt.
\end{satz}

Es sei $\mu$ eine nichtnegatives reguläres Borelsches Maß über $\R$ mit $\mu(\R)=1$.
Dann definiert $F(t)=\displaystyle \int_{-\infty}^t d\mu(s)$ eine Verteilungsfunktion, d.\,h.
\index{Verteilungsfunktion} \index{Funktion, absolutstetige!Verteilung@--, Verteilungsfunktion}
$F(t)$ ist monoton wachsend, rechtsseitig stetig  und $\lim\limits_{t \to -\infty} F(t)=0$, $\lim\limits_{t \to \infty} F(t)=1$.
In der Tradition von Andrei Nikolajewitsch Kolmogorow (1903--1987) wird abweichend von dieser Definition
in der osteuropäischen Literatur häufig die linksseitige Stetigkeit gefordert.

\begin{satz}[Satz von Lusin] \label{SatzLusin} \index{Satz, Darstellungssatz von Riesz!von Lusin@-- von Lusin}
Für $f:I \to \R$ und  das Lebesgue-Ma"s $\lambda$ sind "aquivalent:
\begin{enumerate}
\item[(a)] Es gibt eine me"sbare Funktion $g:I \to \R$ mit $f=g$ $\lambda$-f."u.
\item[(b)] Zu jedem offenen $U \subseteq I$ mit $\lambda(U)<\infty$ und jedem $\delta>0$ gibt es ein Kompaktum $K \subset U$ mit
           $\lambda(U \setminus K) < \delta$, so dass $f |_K$ stetig ist.
\item[(c)] Zu jedem Kompaktum $T \subset I$ und jedem $\delta>0$ gibt es ein Kompaktum $K \subset T$ mit
           $\lambda(T \setminus K) < \delta$, so dass $f |_K$ stetig ist.
\end{enumerate}
\end{satz}

\subsection{Funktionen beschr\"ankter Variation und Stieltjes-Integral}
Eine Funktion $g:[a,b] \to \R$ hei"st
von beschr"ankter Variation\index{Funktion, absolutstetige!von beschränkter@-- von beschränkter Variation} auf $[a,b]$,
wenn es eine Konstante $M>0$ derart gibt, dass f"ur jede Zerlegung $Z:=\{x_0,x_1,...,x_n\}$ von $[a,b]$
mit $a=x_0<x_1<...<x_n=b$ stets $\displaystyle \Var(g,Z):=\sum_{k=1}^n|g(x_k)-g(x_{k-1})| \leq M$ bleibt.
In diesem Fall wird die reelle Zahl $\displaystyle\Var_a^b(g) = \sup_{Z} \Var(g,Z)$ die totale Variation\index{totale Variation, einer Funktion}
\index{Funktion, absolutstetige!totale Variation@--, totale Variation einer} von $g$ auf $[a,b]$ genannt.
Eine Funktion $g(\cdot)$ von beschr"ankter Variation hei"st normalisiert,
wenn $g(a)=0$ gilt und $g(\cdot)$ f"ur alle $a<t<b$ rechtsseitig stetig ist. \\[2mm]
F"ur unbeschr"ankte Intervalle definieren wir die totale Variation durch
$$\Var_{-\infty}^\infty(g)=\lim_{n \to \infty} \Var_{-n}^n(g), \qquad \Var_a^\infty(g)=\lim_{n \to \infty} \Var_a^n(g).$$
Sind diese Variationen endlich,
dann gelten $\displaystyle\lim_{a \to -\infty} \Var_{-\infty}^a(g)=0$ und $\displaystyle \lim_{a \to \infty} \Var_a^{+\infty}(g)=0$.

\begin{satz}
Eine Funktion ist genau dann von beschr"ankter Variation, wenn sie als Differenz zweier monotoner Funktionen darstellbar ist.
\end{satz}

Auf $[a,b]$ seien zwei endliche Funktionen $f$ und $g$ definiert.
Weiterhin sei $Z=\{x_0,...,x_n\}$ eine Zerlegung von $[a,b]$.
In jedem Teilintervall $[x_{k-1},x_k]$ w"ahlen wir einen Punkt $\xi_k$ und bilden die Summe
$$\sigma=\sum_{k=1}^n f(\xi_k) [g(x_k)-g(x_{k-1})].$$
Konvergiert $\sigma$ bei Verfeinerung der Zelegung gegen einen endlichen Grenzwert $I$,
der weder von der Zerlegung noch von der Auswahl der Punkte $\xi_k$ abh"angig ist,
so hei"st der Grenzwert $I$ das Stieltjes-Integral\index{Stieltjes-Integral} von $f$ nach der Funktion $g$.
Wir schreiben $\displaystyle \int_a^b f(x) \, dg(x)$.

\begin{satz}
Existiert eines der Integrale $\displaystyle \int_a^b f(x) \, dg(x)$ oder $\displaystyle \int_a^b g(x) \, df(x)$,
so existiert auch das andere, und es gilt die Regel der partiellen Integration:
$$\int_a^b f(x) \, dg(x) + \int_a^b g(x) \, df(x) =f(b)g(b)-f(a)g(a).$$
\end{satz}

\begin{satz}
Ist $f$ auf $[a,b]$ stetig und $g$ auf $[a,b]$ von beschr"ankter Variation,
dann existiert das Stieltjes-Integral von $f$ nach der Funktion $g$.
\end{satz}

Eine Funktion $g:[a,b] \to \R$ hei"st absolutstetig\index{Funktion, absolutstetige}, 
wenn zu jedem $\varepsilon >0$ ein $\delta>0$ derart existiert,
dass f"ur jedes endliche System von disjunkten Intervallen $\{(a_k,b_k)\}$ mit Gesamtl"ange
$\displaystyle\sum (b_k-a_k) \leq \delta$ die Ungleichung $\sum |g(b_k)-g(a_k)| \leq \varepsilon$ gilt.

\begin{satz}
Ist $f$ auf $[a,b]$ stetig und $g$ absolutstetig auf $[a,b]$,
so f"uhrt die Berechnung des Stieltjes-Integral $\displaystyle \int_a^b f(x) \, dg(x)$ auf das Lebesgue-Integral $\displaystyle \int_a^b f(x) g'(x)\, dx$:
$$\int_a^b f(x) \, dg(x) = \int_a^b f(x) g'(x)\, dx.$$
\end{satz}

\subsection{Borelsche Maße über der erweiterten reellen Zahlengeraden} \label{AnhangMassR}
\index{Borelsches Maß}\index{Maß!Borelsches@--, Borelsches}\index{Maß!BorelschesM@-- Borelsches Ma"s über $\overline{\R}$}
Die erweiterte reelle Zahlengerade $\overline{\R}$ entsteht durch Hinzufügen der uneigentlichen Punkte $-\infty$ und $+\infty$;
also $\overline{\R}=\R \cup \{-\infty,+\infty\}$.
Über $\overline{\R}$ stellt die Abbildung $f:\overline{\R} \to [-1,1]$ mit
$$f(x)=\frac{x}{1+|x|} \mbox{ für } x \in \R, \quad f(-\infty)=-1,\quad f(+\infty)=1,$$
einen Homöomorphismus dar.
Topologisch ist daher $\overline{\R}$ völlig gleichwertig zu einem abgeschlossenen Intervall.
Demzufolge stellen für jedes $a \in \R$ die Mengen
$$[-\infty,a)=\{x \in \R\,|\, x<a\}\cup \{-\infty\} \quad\mbox{bzw.}\quad (a,+\infty]=\{x \in \R\,|\, x>a\}\cup \{\infty\}$$
offene Umgebungen des uneigentlichen Punktes $-\infty$ bzw. $+\infty$ dar.
Ferner ergibt sich damit die Kompaktheit der Mengen $[-\infty,a]$ mit $-\infty<a\leq \infty$ und $[b,\infty]$ mit
$-\infty \leq b < \infty$. \\[2mm]
Es bezeichnet $\mathscr{B}(\R)$ die Borelsche $\sigma$-Algebra\index{Borelmengen},
die von den offenen Teilmengen von $\R$ erzeugt wird.
Die Borelsche $\sigma$-Algebra $\mathscr{B}(\overline{\R})$ besitzt die Darstellung
$$\mathscr{B}(\overline{\R})
  =\{ A \subseteq \overline{\R} \,|\, A=B\cup E, \; B \in \mathscr{B}(\R),\; E \subseteq \{-\infty,+\infty\}\}.$$
Für messbare Mengen $A \in \mathscr{B}(\overline{\R})$ kommt im Weiteren der Zerlegung $A=B\cup E$ mit $B \in \mathscr{B}(\R)$
und $E \subseteq \{-\infty,+\infty\}$ eine gesonderte Betrachtung zu,
da die uneigentlichen Punkte $-\infty$ und $+\infty$ eigenständig behandelt werden müssen. \\[2mm]
Ein Borelsches Ma"s ist ein Maß auf einer Borelschen $\sigma$-Algebra.
Für ein signiertes regul"ares Borelsches Ma"s $\mu_0$ "uber $\R$ ist die Totalvariation $|\mu_0|$ ein endliches Ma"s "uber $\R$
und es folgt aus $\sigma$-Additivit"at und der absoluten Konvergenz der Reihen
$$\sum_{n =0}^\infty \mu_0([n,n+1)), \qquad \sum_{n =0}^\infty \mu_0([-n,-n-1)),$$
dass $\mu_0([\alpha,\infty)) \to 0$ und $\mu_0((-\infty,-\alpha]) \to 0$ f"ur $\alpha \to \infty$ gelten
(vgl. \cite{Rudin}, S.\,116).

\begin{definition}
Für $A \in \mathscr{B}(\overline{\R})$ legen wir $\mu(A)=\mu(B\cup E):=\mu(B)+\mu(E)$
mit der Zerlegung $A=B\cup E$, $B \in \mathscr{B}(\R)$ und $E \subseteq \{-\infty,+\infty\}$ fest.
Genauer hei"st dies, dass die uneigentlichen Punkte stets gesondert betrachtet werden müssen.
\end{definition}

Die Separierung der uneigentlichen Punkte ist erforderlich;
andernfalls ergeben sich grundlegende Widersprüche.
Zur Veranschaulichung betrachten wir das Ma"s $\mu$ mit
$$\mu(\{-\infty\})=0, \; \mu((\alpha,0))=0,\; \mu([0,\beta])=1-e^{-\beta},\; \mu(\{\infty\})=1,$$
wobei $-\infty<\alpha<0<\beta<\infty$ erfüllt sind.
Hier entsteht der Widerspruch
$$\mu(\overline{\R}_+)=2
 \not= \mu(\{0\})+\sum_{n =0}^\infty \mu((n,n+1])=\lim_{k \to \infty}\sum_{n =0}^k \mu((n,n+1])=1.$$
Separiert man andererseits den uneigentlichen Punkt $+\infty$, d.\,h. $\mu(\overline{\R}_+)=\mu(\R_+)+\mu(\{\infty\})$,
so darf $\R_+$ bzw. $\R$ beliebig zerlegt werden und die $\sigma$-Additivit"at bleibt erhalten. \\
Es sei $\mu$ ein signiertes regul"ares Borelsches Ma"s "uber $\overline{\R}$.
Wie oben ergeben sich aus der $\sigma$-Additivit"at und der absoluten Konvergenz der Reihen
$$\lim_{\alpha \to \infty} \mu((-\infty,-\alpha]) \to 0, \qquad
  \lim_{\alpha \to \infty} \mu([\alpha,\infty)) \to 0.$$
Mit der Separierung der uneigentlichen Punkte $-\infty$ und $+\infty$ folgen weiterhin
\begin{eqnarray*}
\lim_{\alpha \to \infty} \mu([-\infty,-\alpha]) 
&=& \mu(\{-\infty\}) + \lim_{\alpha \to \infty}  \mu((-\infty,-\alpha]) = \mu(\{-\infty\}), \\
\lim_{\alpha \to \infty} \mu([\alpha,\infty]) 
&=& \lim_{\alpha \to \infty} \mu([\alpha,\infty)) + \mu(\{\infty\}) = \mu(\{\infty\}).
\end{eqnarray*}

\begin{definition} \label{DefinitionMassR} \index{Maß!regul"ares@--, regul"ares}
Ein signiertes Borelsches Ma"s $\mu$ ist "uber $\overline{\R}$ regulär,
wenn $\mu$ die eindeutige Zerlegung $\mu=\mu(\{-\infty\})+\mu_0+\mu(\{\infty\})$ mit einem signierten regul"aren Borelschen
Ma"s $\mu_0$ "uber $\R$ besitzt.
Dabei ist die Endlichkeit aller Ausdrücke Teil der Anforderungen.
\end{definition}

Die Separation der uneigentlichen Punkte $-\infty$ und $+\infty$ ist ebenfalls Bestandteil bei der Integration über
$\overline{\R}$ bezüglich eines Borelschen Ma"ses $\mu=\mu(\{-\infty\})+\mu_0+\mu(\{\infty\})$.
Ohne weitere Anforderung an die Funktion $f$ setzen wir zunächst formal
$$\int_0^\infty f(t) \, d\mu(t):=
  f(-\infty) \cdot \mu(\{-\infty\}) + \int_\R f(t) \, d\mu_0(t) + f(\infty) \cdot \mu(\{\infty\}).$$
Dabei ist der Integralterm über $\R$ entsprechend des Zuganges in Abschnitt \ref{AbschnittIntegrationMass} aufzufassen.
Die Ausdrücke $f(-\infty) \cdot \mu(\{-\infty\})$ und  $f(\infty) \cdot \mu(\{\infty\})$ sind allerdings nur sinnvoll,
wenn die Funktion $f$ in einem gewissen Sinn die Grenzwerte $f(-\infty)$ und $f(\infty)$ in den uneigentlichen Punkten 
$-\infty$ bzw. $+\infty$ besitzt.
Beispiele derartiger Funktionenklassen sind stetige Funktionen, die im Unendlichen konvergieren, oder
messbare Funktionen, die im Unendlichen dem Maße nach konvergieren.
Diese Funktionenklassen werden im nachfolgenden Abschnitt zu den funktionalanalytischen Hilfsmitteln eingeführt und untersucht.

%% file: II-Hilfsmittel.tex
\section{Funktionalanalytische Hilfsmittel} \label{AnhangHilfsmittel}
Die Grundbegriffe und die Haupts"atze zur Funktionalanalysis, sowie die Darstellung stetiger linearer Funktionale sind
Heuser \cite{HeuserFA}, Ioffe \& Tichomirov \cite{Ioffe} und Werner \cite{Werner}
entnommen. Den Fixpunktsatz von Weissinger findet man in Heuser \cite{HeuserGD},
den Darstellungssatz f"ur die stetigen linearen Funktionale auf dem Raum der im Unendlichen verschwindenden stetigen Funktionen
in Rudin \cite{Rudin}.
Bei der Darstellung des Satzes von Ljusternik folgen wir Ioffe \& Tichomirov \cite{Ioffe}.
Dieser Satz geht auf Ljusternik \cite{Ljusternik} zur"uck.

\subsection{Stetige lineare Operatoren auf normierten R\"aumen} \label{AbschnittNormRaum}
Es seien $X,Y$ normierte lineare R"aume.
Der normierte lineare Raum $X$ mit Norm $\|\cdot\|$ hei"st vollst"andig oder ein Banachraum\index{Banachraum},
wenn jede Cauchyfolge in $X$ bez"uglich $\|\cdot\|$ konvergiert.
Der normierte linearer Raum $X$ mit Skalarprodukt $\langle \cdot,\cdot \rangle$ hei"st ein
Hilbertraum\index{Hilbertraum},
wenn $X$ bez"uglich der induzierten Norm $\|x\|=\sqrt{\langle x,x \rangle}$ vollst"andig ist. \\[2mm]
Auf $X$ wird eine Norm $\|\cdot\|$ genau dann durch ein Skalarprodukt $\langle \cdot,\cdot \rangle$ induziert,
wenn die Parallelogrammgleichung $\|x+y\|^2 + \|x-y\|^2 = 2 (\|x\|^2+\|y\|^2)$ f"ur alle $x,y \in X$ gilt.
In diesem Fall ist das Skalarprodukt eindeutig bestimmt und wird in einem reellen Vektorraum $X$ durch folgende Formel gegeben:
$$\langle x,y \rangle = \frac{1}{4}(\|x+y\|^2 - \|x-y\|^2).$$
Wir geben einige Beispiele f"ur normierte R"aume an (vgl. \cite{Werner}): Sei $I\subseteq \R$.
\begin{enumerate}
\item[(a)] Der Raum $C_b(I,\R^n)$\label{Cb} der beschr"ankten stetigen Funktionen ist bez"uglich der Supremumsnorm $\|\cdot\|_\infty$ vollst"andig.
\item[(b)] $C(I,\R^n)$\label{C}\index{Raum!stetiger@-- stetiger Funktionen} ist der Raum der stetigen Funktionen auf $I\subseteq \R$.
           Bez"uglich der Supremumsnorm $\|\cdot\|_\infty$ ist dieser Raum f"ur eine kompakte Menge $I$ vollst"andig.
\item[(c)] $C_0(\R_+,\R^n)$\label{C0}
           bezeichnet den Raum der stetigen Funktionen $x(\cdot)$, die im Unendlichen verschwinden.
           D.\,h., zu jedem $\varepsilon >0$ ist die Menge $M = \{t \in \R_+ \,|\, \|x(t)|\ \geq \varepsilon\}$ kompakt.
           Als abgeschlossener Unterraum des Raumes $C_b(\R_+,\R^n)$ ist $C_0(\R_+,\R^n)$ vollst"andig.
\item[(d)] Der\index{Raum!stetig@-- stetig differenzierbarer Funktionen}\label{C1}
           Raum der stetig differenzierbaren Funktionen $C_1(I,\R^n)$ ist "uber einer kompakten Menge $I$ bez"uglich der Norm
           $\|x(\cdot)\|_{C_1} = \|x(\cdot)\|_\infty + \|\dot{x}(\cdot)\|_\infty$ vollst"andig.
\item[(e)] Die\index{Raum!Lebesgue@--, Lebesgue-Raum}\label{Lp}
           Lebesgue-R"aume $L_p(I,\R^n)$, $1\leq p \leq \infty$, sind bez"uglich der $L_p$-Norm $\|\cdot\|_{L_p}$ vollst"andig.
           F"ur $p=2$ liegt ein Hilbertraum vor.
\item[(f)] Die\index{Raum!Sobolev@--, Sobolev-Raum}\label{W1p}
           Sobolev-R"aume $W^1_p(I,\R^n)$ sind bez"uglich 
           $\|x(\cdot)\|_{W^1_p} = \|x(\cdot)\|_{L_p} + \|\dot{x}(\cdot)\|_{L_p}$ f"ur $1\leq p \leq \infty$ vollst"andig.
           Der Raum $W^1_2(I,\R^n)$ ist ein Hilbertraum.
\item[(g)] Die\label{Ellp} $\ell_p$-Folgenr"aume, $1\leq p \leq \infty$, sind bez"uglich der $\ell_p$-Norm $\|\cdot\|_{\ell_p}$ vollst"andig.
           F"ur $p=2$ liegt ein Hilbertraum vor.
\end{enumerate}

Wir erweitern unsere Aufz"ahlung f"ur spezielle normierte R"aume von Funktionen und Folgen,
die (in einem gewissen Sinn) im Unendlichen Konvergenz aufweisen.

\begin{definition}
\begin{enumerate}
\item[(a)] $x(\cdot):\R_+ \to \R^n$ konvergiert im Unendlichen gegen $a \in \R^n$,
           wenn es zu jedem $\varepsilon >0$ eine Zahl $N>0$ gibt mit
           $\|x(t)-a\|\leq \varepsilon$ für alle $t \geq N$.
\item[(b)] $x(\cdot):\R_+ \to \R^n$ konvergiert im Lebesgue-Maß $\lambda$ im Unendlichen gegen $a \in \R^n$,
           wenn $\lim\limits_{N \to \infty} \lambda\big(\{t\geq N \,|\, \|x(t)-a\|\geq \varepsilon\}\big)=0$
           für jedes $\varepsilon >0$ gilt.
\end{enumerate}
\end{definition}

Wir betrachten nun Räume mit folgender Struktur:
Es sei $\Omega$ eine nichtleere Menge und es sei $X_0=\{x_0(\cdot): \Omega \to \R^n\}$ ein normierter Raum mit Norm $\|\cdot\|_{X_0}$.
Dann sei $X$ der Raum derjenigen Elemente $x(\cdot)$,
die die eindeutige Darstellung $x(\cdot)=x_0(\cdot)+a$,
d.\,h. $x(\omega)=x_0(\omega)+a$ für alle $\omega \in \Omega$,
mit $x_0(\cdot) \in X_0$ und $a \in \R^n$ besitzen.
Weiterhin versehen wir $X$ mit einer $p$-Norm, d.\,h.
$$\|x(\cdot)\|^p_X = \|x_0(\cdot)\|^p_{X_0} + \|a\|^p \mbox{ f"ur } 1 \leq p < \infty, \qquad
  \|x(\cdot)\|_X = \|x_0(\cdot)\|_{X_0} + \|a\| \mbox{ f"ur } p=\infty.$$

\begin{lemma} \label{LemmaBanachraum}
\begin{enumerate}
\item[(a)] Ist $(X_0,\|\cdot\|_{X_0})$ ein Banachraum, dann auch $(X,\|\cdot\|_X)$.
\item[(b)] Ist $(X_0,\|\cdot\|_{X_0})$ ein Hilbertraum, dann auch $(X,\|\cdot\|_X)$ bezüglich der $2$-Norm.
\end{enumerate}
\end{lemma}

In den nachstehenden Beispielen sind die Elemente der Räume gerade durch die Struktur
$x(\cdot)=x_0(\cdot)+a$ mit $x_0(\cdot) \in X_0$ und $a \in \R^n$ charakterisiert.\index{Raum!im@-- im Unendlichen konvergenter Funktionen}
\begin{enumerate}
\item[(h)] $C_{\lim}(\R_+,\R^n)$ \label{Clim}
           bezeichnet den Raum der stetigen Funktionen $x(\cdot)$, die im Unendlichen einen Grenzwert besitzen:
           $$C_{\lim}(\R_+,\R^n)
             =\{ x(\cdot) \,|\, x(\cdot)=  x_0(\cdot)+x_0,\; x_0(\cdot) \in C_0(\R_+,\R^n), \; x_0 \in \R^n \}.$$
           Als abgeschlossener Unterraum des Raumes $C_b(\R_+,\R^n)$ ist $C_{\lim}(\R_+,\R^n)$ bezüglich der
           Supremumsnorm $\|\cdot\|_\infty$ vollst"andig.
           Au"serdem lassen sich mit Hilfe der erweiterten reellen Zahlengeraden die normierten Räume
           $\big(C_{\lim}(\R_+,\R^n),\|\cdot\|_\infty\big)$ und $\big(C(\overline{\R}_+,\R^n),\|\cdot\|_\infty\big)$
           miteinander identifizieren.
\item[(i)] F"ur $1 \leq p< \infty$ f"uhren wir die Lebesgue-R"aume $L_{p,\lim}(\R_+,\R^n)$ ein,
           deren Elemente im Unendlichen dem Maße nach konvergieren: \label{Lplim}
           $$L_{p,\lim}(\R_+,\R^n) =\{ x(\cdot) \,|\,
               x(\cdot)= x_0(\cdot)+ x_0,\; x_0(\cdot) \in L_p(\R_+,\R^n), \; x_0 \in \R^n \}.$$
           Bez"uglich \label{LplimNorm}
           $\displaystyle \|x(\cdot)\|_{L_{p,\lim}} = \big(\|x_0(\cdot)\|^p_{L_p} + \|x_0\|^p \big)^{1/p}$
           ist $L_{p,\lim}(\R_+,\R^n)$ vollst"andig. \\
           Im Fall $p=2$ gilt die Parallelogrammgleichung und $L_{2,\lim}(\R_+,\R^n)$ stellt mit dem Skalarprodukt
           $\displaystyle \langle x(\cdot), y(\cdot) \rangle
             = \int_0^\infty \langle x_0(t), y_0(t) \rangle \, dt + \langle x_0, y_0 \rangle$
           einen Hilbertraum dar.
\item[(j)] F"ur $1 \leq p< \infty$ definieren wir die Sobolev-R"aume $W^1_{p,\lim}(\R_+,\R^n)$ wie folgt: \label{W1plim}
           $$W^1_{p,\lim}(\R_+,\R^n) =\{ x(\cdot) \,|\,
                   x(\cdot)=x_0(\cdot) + x_0,\; x_0(\cdot) \in W^1_p(\R_+,\R^n), \; x_0 \in \R^n \}.$$
           Bez"uglich \label{W1plimNorm}
           $\displaystyle\|x(\cdot)\|_{W^1_{p,\lim}} = \big(\|x_0(\cdot)\|^p_{W^1_p} + \|x_0\|^p \big)^{1/p}$
           ist $W^1_{p,\lim}(\R_+,\R^n)$ vollst"andig. \\
           Speziell ist $W^1_{2,\lim}(\R_+,\R^n)$ ein Hilbertraum.
\item[(k)] F"ur $1 \leq p< \infty$ definieren wir die R"aume $\ell_{p,\lim}$ konvergenter Folgen: \label{Ellplim}
           $$\ell_{p,\lim} =\{ x=(x_n)_{n \in \N} \,|\, x_n= x^0_n+ x_0,\; x^0=(x^0_n)_{n \in \N} \in \ell_p, \; x_0 \in \R \}.$$
           Bez"uglich der Norm \label{EllplimNorm} $\displaystyle\|x\|_{\ell_{p,\lim}} = \big(\|x^0\|^p_{\ell_p} + \|x_0\|^p \big)^{1/p}$
           ist $\ell_{p,\lim}$ vollst"andig.
           Ferner ist $\ell_{2,\lim}$ ein Hilbertraum.
\end{enumerate}

Es bezeichnet $L(X,Y)$ die Menge der stetigen linearen Abbildungen $T:X \to Y$.
Ist $Y=\R$, dann nennen wir $T:X \to \R$ ein Funktional.
Dabei hei"st eine (nicht notwendig lineare) Abbildung $T$ stetig\index{Abbildung, beschränkte!stetige@--, stetige},
wenn das Urbild $T^{-1}(O)=\{x \in X \,|\, Tx \in O\}$ offener Mengen $O \subseteq Y$ offen in $X$ ist.
Dabei verwenden wir die "ubliche Schreibweise $Tx$ statt $T(x)$.
Zu dieser Definition der Stetigkeit sind "aquivalent:
\begin{enumerate}
\item[(i)] Aus $x_n \to x$ folgt $Tx_n \to Tx$.
\item[(ii)] Zu jedem $x_0 \in X$ und jedem $\varepsilon>0$ existiert ein $\delta=\delta(x_0,\varepsilon) >0$ mit
           $$\|Tx-Tx_0\| \leq \varepsilon \quad \mbox{ f"ur alle } \|x-x_0\| \leq \delta.$$
\end{enumerate}
Insbesondere sind f"ur lineare Abbildungen "aquivalent:
\begin{enumerate}
\item[(a)] $T$ ist stetig, d.\,h. $T$ ist in jedem Punkt $x \in X$ stetig.
\item[(b)] $T$ ist in $x=0$ stetig.
\item[(c)] $T$ ist gleichm"a"sig stetig.
\item[(d)] $T$ ist beschr"ankt\index{Abbildung, beschränkte}, d.\,h. es existiert eine Zahl $K\geq 0$ mit
           $\|Tx\| \leq K \|x\|$ f"ur alle $x \in X$.
\end{enumerate}

Die kleinste Zahl $K \geq 0$, mit der (d) gilt, hei"st die Operatornorm\index{Operatornorm} der linearen Abbildung $T$ und wird mit $\|T\|$ bezeichnet.
F"ur $\|T\|$ gilt
$$\|T\|= \sup_{x \not=0} \frac{\|Tx\|}{\|x\|}=  \sup_{\|x\|=1} \|Tx\| = \sup_{\|x\| \leq 1} \|Tx\|.$$
Der Raum $L(X,\R)$ der stetigen linearen Funktional $x^*:X \to \R$ auf dem normierten Raum $X$ hei"st der Dualraum\index{Dualraum}
von $X$ und wird mit $X^*$ bezeichnet.
Ferner schreiben wir $\langle x^*,x \rangle$ statt $x^*(x)$.

\begin{theorem}[Satz von Fr\'echet-Riesz]  \index{Satz, Darstellungssatz von Riesz!von FrechetRiesz@-- von Fr\'echet-Riesz}
Sei $H$ ein Hilbertraum und $x^* \in H^*$ ein stetiges lineares Funktional.
Dann existiert genau ein $y \in H$ mit $\langle x^*,x \rangle = \langle x, y \rangle$ f"ur alle $x \in H$.
\end{theorem}


\subsection{Differentialrechnung und Differenzierbarkeit konkreter Abbildungen}
Es seien $X,Y$ normierte R"aume, $U \subseteq X$ offen und $F:U \to Y$ eine Abbildung.
Existiert in $x_0 \in U$ f"ur jedes $x \in X$ der Grenzwert
$$\lim_{\lambda \to 0^+} \frac{F(x_0+\lambda x)-F(x_0)}{\lambda}=\delta F(x_0)x,$$
dann hei"st die Abbildung $x \to \delta F(x_0)x$ die erste Variation der Abbildung $F$ in $x_0$.
Gibt es ein $T \in L(X,Y)$ mit $Tx=\delta F(x_0)x$,
dann hei"st $F$ in $x_0$ G\^ateaux-differenzierbar\index{Abbildung, beschränkte!Gateaux@--, G\^ateaux-differenzierbare} und $T$
die G\^ateaux-Ableitung\index{Ableitung, Fr\'echet-!Gateaux@--, G\^ateaux-} von $F$
in $x_0$.
Ferner hei"st $F$ in $x_0$ Fr\'echet-differenzierbar\index{Abbildung, beschränkte!Fr\'echet@--, Fr\'echet-differenzierbare}
und $F'(x_0)=T$ die Fr\'echet-Ableitung\index{Ableitung, Fr\'echet-} von $F$ in $x_0$,
falls der Grenzwert der G\^ateaux-Ableitung bez"uglich $\|x\|\leq 1$ gleichm"a"sig konvergiert.
D.\,h. zu jedem $\varepsilon >0$ gibt es ein $\lambda_0 >0$ mit
$$\bigg\|\frac{F(x_0+\lambda x)-F(x_0)}{\lambda}-F'(x_0)x\bigg\| \leq \varepsilon \quad\mbox{ f"ur alle }
  0<\lambda<\lambda_0, \; \|\xi\|\leq 1.$$
Eine Abbildung $F$ ist in $x_0$ genau dann Fr\'echet-differenzierbar,
falls ein $T \in L(X,Y)$ und eine Abbildung $r:X \to Y$ mit
$$F(x_0+x)=F(x_0)+Tx+r(x), \qquad \lim_{\|x\| \to 0} \frac{\|r(x)\|}{\|x\|} =0.$$
Die Abbildung $F$ nennen wir im Punkt $x$ regul"ar\index{Abbildung, beschränkte!reguläre@--, reguläre},
wenn sie in diesem Punkt Fr\'echet-differenzierbar ist und ${\rm Im\,}F'(x)=Y$ gilt.
Die Abbildung $F$ ist in $x_0$ bzw. auf $U$ stetig
differenzierbar\index{Abbildung, beschränkte!stetig differenzierbare@--, stetig differenzierbare},
wenn f"ur alle Punkte der Menge $U$ die Ableitung $F'(x)$ existiert und die Abbildung $x \to F'(x)$ bez"uglich der Operatornorm des Raum $L(X,Y)$
in $x_0$ bzw. auf $U$ stetig ist. \\[2mm]
Es sei $h:\R^n \to \R^m$ eine differenzierbare Abbildung.
Dann gilt f"ur $x,y\in \R^n$ und $\lambda \not=0$ folgende Gleichung, die wir im Weiteren verwenden werden:
$$\frac{h(x+\lambda y) - h(x)}{\lambda}= \int_0^1 h'(x+ s\lambda y) y  \, ds.$$

\begin{beispiel} \label{DiffAbbildung}
{\rm Sei $I \subseteq \overline{\R}$ ein abgeschlossenes Intervall, d.\,h. eine abgeschlossene konvexe Teilmenge von $\overline{\R}$,
und $y_*(\cdot)$ eine stetige Vektorfunktion auf $I$.
Zu $y_*(\cdot)$ definieren wir die Menge $V_\gamma = \{ (t,y) \in I \times \R^n \,|\, \|y-y_*(t)\| \leq \gamma\}$. \\
Wir nehmen an,
dass die Abbildung $g(t,y) : \R \times \R^n \to \R^m$ auf der Menge $V_\gamma$
gleichm"a"sig stetig und bez"uglich $y$ gleichm"a"sig stetig differenzierbar ist. \\
Dann ist die Abbildung $G: C(I,\R^n) \to C(I,\R^m)$,
$G\big(y(\cdot)\big)(t) = g\big(t,y(t)\big)$,
im Punkt $y_*(\cdot)$ stetig Fr\'echet-differenzierbar und es gilt
$$\big[G'\big(y_*(\cdot)\big) y(\cdot)\big](t) = g_y\big(t,y_*(t)\big) y(t), \quad t \in I.$$
Denn: 
F"ur $t \in I$, $\|y(\cdot)-y_*(\cdot)\|_\infty \leq \gamma$ und $0 < \lambda \leq \lambda_0$ gilt
$$\bigg[\frac{G\big(y(\cdot) + \lambda x(\cdot)\big) - G\big(y(\cdot)\big)}{\lambda} - G'\big(y(\cdot)\big)x(\cdot)\bigg](t)
  = \hspace*{-1mm}\int_0^1  \hspace*{-1mm} \big[g_y\big(t,y(t)+s\lambda x(t)\big) - g_y\big(t,y(t)\big) \big] x(t) ds.$$
Da die Abbildung $g_y(t,y)$ auf der Menge $V_\gamma$ gleichm"a"sig stetig ist,
gibt es ein $C>0$ mit 
$$\bigg\| \frac{G\big(y(\cdot) + \lambda x(\cdot)\big) - G\big(y(\cdot)\big)}{\lambda}
   - G'\big(y(\cdot)\big)x(\cdot) \bigg\|_\infty
  \leq \sup_{t \in I} \int_0^1 C \|\lambda x(t)\| ds \cdot \|x(\cdot)\|_\infty = \lambda C,$$
d.\,h. der Grenzwert $\lambda \to 0^+$ konvergiert gleichm"a"sig bez"uglich $\|x(\cdot)\|_\infty = 1$.
Also ist die Abbildung $G$ auf einer Umgebung von $y_*(\cdot)$ Fr\'echet-differenzierbar. \\
Weiterhin ergibt sich wegen der gleichm"a"sigen Stetigkeit von $g_y(t,y)$ auf der Menge $V_\gamma$
f"ur die Abbildung $y(\cdot) \to G'\big(y(\cdot)\big)$ bez"uglich der Operatornorm in $y_*(\cdot)$:
\begin{eqnarray*}
    \big\|G'\big(y(\cdot)\big) - G'\big(y_*(\cdot)\big) \big\|
&=& \sup_{\|x(\cdot)\|_\infty=1} \big\| \big[G'\big(y(\cdot)\big) - G'\big(y_*(\cdot)\big)\big] x(\cdot) \big\|_\infty \\
&\leq& \sup_{t \in I} \big\|  g_y\big(t,y(t)\big) - g_y\big(t,y_*(t)\big) \big\|
       \leq C \|y(\cdot)-y_*(\cdot) \|_\infty.
\end{eqnarray*}
Somit ist die stetige Fr\'echet-Differenzierbarkeit nachgewiesen. \hfill $\square$}
\end{beispiel}

\begin{beispiel} \label{DiffZielfunktional1}
{\rm Es sei $y_*(\cdot) \in L_\infty([t_0,t_1],\R^n)$.
Dann k"onnen wir annehmen, dass $y_*(\cdot)$ auf $[t_0,t_1]$ beschr"ankt ist.
Ferner sei die Funktion $f(t,y): \R \times \R^n \to \R$ auf dem Abschluss der Menge
$V_\gamma = \{ (t,y) \in [t_0,t_1] \times \R^n \,|\, \|y-y_*(t)\| \leq \gamma\}$
stetig und stetig differenzierbar bez"uglich $y$.
Dann ist $J: L_\infty([t_0,t_1],\R^n) \to \R$, $\displaystyle J\big(y(\cdot)\big) = \int_{t_0}^{t_1} f\big(t,y(t)\big) \, dt$,
im Punkt $y_*(\cdot)$ Fr\'echet-differenzierbar und es gilt
$$J'\big(y_*(\cdot)\big) y(\cdot) = \int_{t_0}^{t_1} \big\langle f_y\big(t,y_*(t)\big), y(t) \big\rangle \, dt.$$
Denn: Nach Voraussetzung ist die Abbildung $f(t,y)$ auf der Menge $\overline{V}_\gamma$ stetig und stetig differenzierbar bez"uglich $y$.
Daher ist $J$ auf einer Umgebung von $y_*(\cdot)$ wohldefiniert.
Ferner ist der lineare Operator $J'\big(y_*(\cdot)\big)$ stetig.
Sei $\varepsilon >0$ gegeben.
Aufgrund der gleichm"a"sigen Stetigkeit von $f_y(t,y)$ auf der Menge $\overline{V}_\gamma$ gibt es eine Zahl $\lambda_0>0$ mit
$$\big\|f_y\big(t,y_*(t)+\lambda y\big)-f_y\big(t,y_*(t)\big)\big\| \leq \frac{\varepsilon}{t_1-t_0}$$
f"ur fast alle $t \in [t_0,t_1]$, f"ur alle $\|y\| \leq 1$ und alle $0 < \lambda \leq \lambda_0$.
Damit erhalten wir
\begin{eqnarray*}
\lefteqn{\bigg|\frac{J\big(y_*(\cdot)+\lambda y(\cdot)\big) - J\big(y_*(\cdot)\big)}{\lambda} - J'\big(y_*(\cdot)\big)y(\cdot)\bigg|} \\
&=& \bigg|\int_{t_0}^{t_1} \bigg[\int_0^1 \langle f_y\big(t,y_*(t)+\lambda s y(t)\big) - f_y\big(t,y_*(t)\big), y(t) \rangle \, ds \bigg] \, dt \bigg|
    \leq \varepsilon
\end{eqnarray*}
f"ur alle $\|y(\cdot) \|_{L_\infty} \leq 1$ und alle $0 < \lambda \leq \lambda_0$. \hfill $\square$}
\end{beispiel}

\begin{beispiel} \label{DiffZielfunktional2}
{\rm Es seien $[t_0,t_1] \subset \R$, $U \subseteq \R^m$, $x_*(\cdot) \in C([t_0,t_1],\R^n)$ und
$$V_\gamma = \{ (t,x) \in [t_0,t_1] \times \R^n \,|\, \|x-x_*(t)\| \leq \gamma\}.$$
Ferner sei die Funktion $f(t,x,u): [t_0,t_1] \times \R^n \times \R^m \to \R$ f"ur jede kompakte Menge $U_1 \subseteq \R^m$
auf der Menge $V_\gamma \times U_1$ stetig und bez"uglich $x$ stetig differenzierbar.
Dann ist f"ur jedes $u(\cdot) \in L_\infty([t_0,t_1],U)$ die Abbildung
$$x(\cdot) \to J\big(x(\cdot),u(\cdot)\big) = \int_{t_0}^{t_1} f\big(t,x(t),u(t)\big) \, dt, \quad J: C([t_0,t_1],\R^n) \times L_\infty([t_0,t_1],U) \to \R,$$
im Punkt $x_*(\cdot)$ Fr\'echet-differenzierbar und es gilt
$$J_x\big(x_*(\cdot),u(\cdot)\big) x(\cdot) = \int_{t_0}^{t_1} \big\langle f_x\big(t,x_*(t),u(t)\big), x(t) \big\rangle \, dt.$$
Denn: Wegen $u(\cdot) \in L_\infty([t_0,t_1],U)$ l"asst sich eine kompakte Menge $U_1 \subseteq \R^m$ mit $u(t) \in U_1$ f"ur fast alle $t$ angeben.
Nach Voraussetzung ist die Abbildung $f(t,x,u)$ auf der Menge $V_\gamma \times U_1$ stetig und stetig differenzierbar bez"uglich $x$.
Daher ist $J$ wohldefiniert.
Weiterhin ist f"ur jedes $u(\cdot) \in L_\infty([t_0,t_1],U)$ der lineare Operator $J_x\big(x_*(\cdot),u(\cdot)\big)$ stetig. \\
Sei $\varepsilon >0$ gegeben.
Aufgrund der gleichm"a"sigen Stetigkeit von $f_x(t,x,u)$ auf der Menge $V_\gamma \times U_1$ gibt es eine Zahl $\lambda_0>0$ mit
$$\big\|f_x\big(t,x_*(t)+\lambda x,u\big)-f_x\big(t,x_*(t),u\big)\big\| \leq \frac{\varepsilon}{t_1-t_0}$$
f"ur fast alle $t \in [t_0,t_1]$ und f"ur alle $\|x\| \leq 1$, $u \in U_1$, $0 < \lambda \leq \lambda_0$.
Damit erhalten wir
\begin{eqnarray*}
\lefteqn{\bigg|\frac{J\big(x_*(\cdot)+\lambda x(\cdot),u(\cdot)\big) - J\big(x_*(\cdot),u(\cdot)\big)}{\lambda}
               - J_x\big(x_*(\cdot),u(\cdot)\big)x(\cdot)\bigg|} \\
&=& \bigg|\int_{t_0}^{t_1} \bigg[\int_0^1 \langle f_x\big(t,x_*(t)+\lambda s x(t),u(t)\big) - f_x\big(t,x_*(t),u(t)\big), x(t) \rangle \, ds \bigg] \, dt \bigg|
    \leq \varepsilon
\end{eqnarray*}
f"ur alle $\|x(\cdot) \|_\infty \leq 1$ und alle $0 < \lambda \leq \lambda_0$. \hfill $\square$}
\end{beispiel}

\begin{beispiel} \label{DiffZielfunktionalS}
{\rm Es seien $U \subseteq \R^m$, $\omega(\cdot)$ eine Gewichtsfunktion aus dem Raum $L_1(\R_+,\R_+)$,
weiterhin $x_*(\cdot) \in C_{\lim}(\R_+,\R^n)$ und $V_\gamma= \{ (t,x) \in \overline{\R}_+ \times \R^n \,|\, \|x-x(t)\| \leq \gamma\}$. \\
Ferner sei die Funktion $f(t,x,u): \R \times \R^n \times \R^m \to \R$ f"ur jede kompakte Menge $U_1 \subseteq \R^m$
auf der Menge $V_\gamma \times U_1$ gleichm"a"sig stetig und bez"uglich $x$ gleichm"a"sig stetig differenzierbar.
Dann ist f"ur jedes $u(\cdot) \in L_\infty(\R_+,U)$ die Abbildung 
$$x(\cdot) \to J\big(x(\cdot),u(\cdot)\big) = \int_0^\infty \omega(t) f\big(t,x(t),u(t)\big) \, dt,$$
wobei $J:C_{\lim}(\R_+,\R^n) \times L_\infty(\R_+,U) \to \R$ gilt,
im Punkt $x_*(\cdot)$ Fr\'echet-differenzierbar und
$$J_x\big(x_*(\cdot),u(\cdot)\big) x(\cdot) = \int_0^\infty \omega(t) \big\langle f_x\big(t,x_*(t),u(t)\big), x(t) \big\rangle \, dt.$$
Denn: Wegen $u(\cdot) \in L_\infty(\R_+,U)$, $\omega(\cdot) \in L_1(\R_+,\R_+)$ und wegen der gleichm"a"sigen Differenzierbarkeit lassen sich
zu jedem $\varepsilon>0$ Zahlen $\lambda_0,T>0$ und nach dem Satz von Lusin eine kompakte Menge $K \subseteq [0,T]$ derart angeben,
dass $\omega(\cdot)$ auf $K$ stetig ist, die Ungleichung
$$\int_{\R_+ \setminus K} \omega(t) \,dt \leq \frac{\varepsilon}{2}$$
gilt und die Beziehung
$$\big\|f_x\big(t,x_*(t)+\lambda x,u(t)\big)-f_x\big(t,x_*(t),u(t)\big)\big\| \leq 1$$
f"ur fast alle $t\not \in K$ und f"ur alle $\|x\| \leq 1$, $0 < \lambda \leq \lambda_0$ erf"ullt ist.
Benutzen wir "uber der kompakten Menge $K$ die gleiche Argumentation wie in Beispiel \ref{DiffZielfunktional2},
so erhalten wir
\begin{eqnarray*}
\lefteqn{\bigg|\frac{J\big(x_*(\cdot)+\lambda x(\cdot),u(\cdot)\big) - J\big(x_*(\cdot),u(\cdot)\big)}{\lambda}
               - J_x\big(x_*(\cdot),u(\cdot)\big)x(\cdot)\bigg|} \\
&=& \bigg|\int_0^\infty \omega(t)
          \bigg[\int_0^1 \langle f_x\big(t,x_*(t)+\lambda s x(t),u(t)\big) - f_x\big(t,x_*(t),u(t)\big), x(t) \rangle \, ds \bigg] \, dt \bigg|
    \leq \varepsilon
\end{eqnarray*}
f"ur alle $\|x(\cdot) \|_\infty \leq 1$ und alle $0 < \lambda \leq \lambda_0$. \hfill $\square$}
\end{beispiel}

\begin{beispiel} \label{DiffDynamik1}
{\rm  Es sei $y_*(\cdot) \in L_\infty([t_0,t_1],\R^n)$.
Dann k"onnen wir annehmen, dass $y_*(\cdot)$ f"ur $t \in [t_0,t_1]$ beschr"ankt ist.
Ferner sei die Abbildung $\varphi(t,y): \R \times \R^n \to \R^n$ auf dem Abschluss der Menge
$V_\gamma = \{ (t,y) \in [t_0,t_1] \times \R^n \,|\, \|y-y_*(t)\| \leq \gamma\}$
stetig und stetig differenzierbar bez"uglich $y$.
Dann ist $F: L_\infty([t_0,t_1],\R^n) \to C([t_0,t_1],\R^n)$,
$$F\big(y(\cdot)\big)(t) = \int_{t_0}^t \varphi\big(s,y(s)\big) \, ds, \quad t \in [t_0,t_1],$$
im Punkt $y_*(\cdot)$ stetig Fr\'echet-differenzierbar und es gilt
$$\big[F'\big(y_*(\cdot)\big) y(\cdot)\big](t) = \int_{t_0}^t  \varphi_y\big(s,y_*(s)\big) y(s) \, ds, \quad t \in [t_0,t_1].$$
Denn: Nach Voraussetzung ist die Abbildung $F$ auf einer Umgebung von $y_*(\cdot)$ wohldefiniert und der lineare Operator $F'\big(y_*(\cdot)\big)$ stetig. \\
Es seien $y(\cdot)$ mit $\|y(\cdot)-y_*(\cdot)\|_{L_\infty} \leq \gamma$ und $\varepsilon >0$ gegeben.
Da die Abbildung $\varphi_y(t,y)$ auf $\overline{V}_\gamma$ gleichm"a"sig stetig ist,
gibt es eine Zahl $\lambda_0>0$ mit 
$$\big\|\varphi_y\big(t,y(t)+\lambda x\big)-\varphi_y\big(t,y(t)\big)\big\| \leq \frac{\varepsilon}{t_1-t_0}$$
f"ur fast alle $[t_0,t_1]$, f"ur alle $\|x\| \leq 1$ und alle $0 < \lambda \leq \lambda_0$.
Damit folgt
\begin{eqnarray*}
\lefteqn{\bigg\|
     \frac{F\big(y(\cdot)+\lambda x(\cdot)\big) - F\big(y(\cdot)\big)}{\lambda} - F'\big(y(\cdot)\big) x(\cdot)\bigg\|_\infty } \\
&\leq& \int_{t_0}^{t_1} \bigg[ \int_0^1
       \big\|\varphi_y\big(t,y(t)+\lambda s x(t)\big)-\varphi_y\big(t,y(t)\big)\big\|\, ds \bigg]  \, dt \leq \varepsilon
\end{eqnarray*}
f"ur alle $\|x(\cdot) \|_{L_\infty} \leq 1$ und alle $0 < \lambda \leq \lambda_0$. \\
Wegen der gleichm"a"sigen Stetigkeit von $f_y(t,y)$ auf $\overline{V}_\gamma$ erhalten wir
f"ur die Abbildung $y(\cdot) \to F'\big(y(\cdot)\big)$ bez"uglich der Operatornorm in $y_*(\cdot)$
\begin{eqnarray*}
    \big\|F'\big(y(\cdot)\big) - F'\big(y_*(\cdot)\big) \big\|
&=& \sup_{\|x(\cdot)\|_{L_\infty}=1} \big\| \big[F'\big(y(\cdot)\big) - F'\big(y_*(\cdot)\big)\big] x(\cdot) \big\|_\infty \\
&\leq& \esssup_{t \in I} \big\|  \varphi_y\big(t,y(t)\big) - \varphi_y\big(t,y_*(t)\big) \big\|
       \leq C \|y(\cdot)-y_*(\cdot) \|_{L_\infty}
\end{eqnarray*}
f"ur alle $y(\cdot)$ mit $\|y(\cdot)-y_*(\cdot)\|_{L_\infty} \leq \gamma$. \hfill $\square$}
\end{beispiel}

\begin{beispiel} \label{DiffDynamik2}
{\rm Es seien $[t_0,t_1] \subset \R$, $U \subseteq \R^m$, $x_*(\cdot) \in C([t_0,t_1],\R^n)$ und
$$V_\gamma = \{ (t,x) \in [t_0,t_1] \times \R^n \,|\, \|x-x_*(t)\| \leq \gamma\}.$$
Ferner sei die Abbildung $\varphi(t,x,u): [t_0,t_1] \times \R^n \times \R^m \to \R^n$ f"ur jede kompakte Menge $U_1 \subseteq \R^m$
auf der Menge $V_\gamma \times U_1$ stetig und bez"uglich $x$ stetig differenzierbar.
Dann ist f"ur jedes $u(\cdot) \in L_\infty([t_0,t_1],U)$ die Abbildung 
$F: C([t_0,t_1],\R^n) \times L_\infty([t_0,t_1],\R^n) \to C([t_0,t_1],\R^n)$,
$$x(\cdot) \to F\big(x(\cdot),u(\cdot)\big)(t) = \int_{t_0}^t \varphi\big(s,x(s),u(s)\big) \, ds, \quad t \in [t_0,t_1],$$
im Punkt $x_*(\cdot)$ Fr\'echet-differenzierbar und es gilt
$$\big[F_x\big(x_*(\cdot),u(\cdot)\big) x(\cdot)\big](t) = \int_{t_0}^t  \varphi_x\big(s,x_*(s),u(s)\big) x(s) \, ds, \quad t \in [t_0,t_1].$$
Denn wie im Beispiel \ref{DiffZielfunktional2} l"asst sich zu $u(\cdot) \in L_\infty([t_0,t_1],U)$ eine kompakte Menge $U_1 \subseteq \R^m$ mit $u(t) \in U_1$ 
f"ur fast alle $t$ angeben,
so dass $\varphi_x(t,x,u)$ auf der Menge $V_\gamma \times U_1$ gleichm"a"sig stetig ist.
Wie im Beispiel \ref{DiffDynamik1} zeigt man damit die Fr\'echet-Differenzierbarkeit. \hfill $\square$}
\end{beispiel}

\begin{beispiel} \label{DiffDynamikS}
{\rm Es seien $U \subseteq \R^m$, $x_*(\cdot) \in C_{\lim}(\R_+,\R^n)$, $u_*(\cdot) \in L_\infty(\R_+,U)$, und
$$V_\gamma= \{ (t,x) \in \overline{\R}_+ \times \R^n \,|\, \|x-x(t)\| \leq \gamma\}.$$
Ferner sei die Abbildung $\varphi(t,x,u): \R \times \R^n \times \R^m \to \R$ f"ur jede kompakte Menge $U_1 \subseteq \R^m$
auf der Menge $V_\gamma \times U_1$ gleichm"a"sig stetig und bez"uglich $x$ gleichm"a"sig stetig differenzierbar.
Weiterhin nehmen wir an,
dass 
$$\int_0^\infty \big\|\varphi\big(t,x_*(t),u_*(t)\big)\big\| \, dt < \infty, \qquad \int_0^\infty \big\|\varphi_x\big(t,x_*(t),u_*(t)\big)\big\| \, dt < \infty$$
gelten.
Au"serdem m"oge zu jedem $\delta>0$ ein $T>0$ existieren, dass
\begin{eqnarray*}
&& \int_T^\infty \big\| \varphi\big(t,x(t),u_*(t)\big)-\varphi\big(t,x'(t),u_*(t)\big) \nonumber \\
&& \hspace*{20mm} - \varphi_x\big(t,x_*(t),u_*(t)\big)\big(x(t)-x'(t)\big) \big\| \, dt \leq \delta \|x(\cdot)-x'(\cdot)\|_\infty
\end{eqnarray*}
f"ur alle $x(\cdot), x'(\cdot) \in C_{\lim}(\R_+,\R^n)$ mit
$\|x(\cdot)-x_*(\cdot)\|_\infty \leq \gamma$, $\|x'(\cdot)-x_*(\cdot)\|_\infty \leq \gamma$ erf"ullt ist.
Au"serdem bezeichnen wir mit $\mathscr{U}$ die Menge:
\begin{eqnarray*}
\mathscr{U}=\big\{ u(\cdot)  \in L_\infty(\R_+,U) &\big|& u(t)=u_*(t) + \chi_M(t)\big(w(t)-u_*(t)\big), \; w(\cdot) \in L_\infty(\R_+,U), \\
                                                  &     & M \subset \R_+ \mbox{ me"sbar und beschr"ankt} \big\}.
\end{eqnarray*}
Die Menge ist verbunden mit Nadelvariationen des Steuerungsprozesses $\big(x_*(\cdot),u_*(\cdot)\big)$,
die nur "uber beschr"ankten Mengen erfolgen. \\[2mm]
Unter diesen Voraussetzungen ist f"ur jedes $u(\cdot) \in \mathscr{U}$ die Abbildung 
$$x(\cdot) \to F\big(x(\cdot),u(\cdot)\big) = \int_0^t \varphi\big(s,x(s),u(s)\big) \, ds, \quad t \in \R_+,$$
$F:C_{\lim}(\R_+,\R^n) \times \mathscr{U} \to C_{\lim}(\R_+,\R^n)$,
in $x_*(\cdot)$ Fr\'echet-differenzierbar und
$$\big[F_x\big(x_*(\cdot),u(\cdot)\big) x(\cdot)\big](t) = \int_0^t  \varphi_x\big(s,x_*(s),u(s)\big) x(s) \, ds, \quad t \in \R_+.$$
Denn: Wegen $u(\cdot) \in \mathscr{U}$ k"onnen wir eine Zahl $T>0$ angeben mit $u(t)=u_*(t)$ f"ur alle $t\geq T$.
Daher ergibt sich f"ur alle $x(\cdot) \in C_{\lim}(\R_+,\R^n)$ mit $\|x(\cdot)-x_*(\cdot)\|_\infty \leq \gamma$
$$\int_0^\infty \big\|\varphi\big(t,x(t),u(t)\big)\big\| \, dt
   \leq \int_0^T \big\|\varphi\big(t,x(t),u(t)\big)\big\| \, dt + \int_T^\infty \big\|\varphi\big(t,x(t),u_*(t)\big)\big\| \, dt.$$
Darin folgt nach Voraussetzung au"serdem
$$\int_T^\infty \big\|\varphi\big(t,x(t),u_*(t)\big)\big\| \, dt
   \leq \int_T^\infty \Big[\big\|\varphi\big(t,x_*(t),u_*(t)\big)\big\| + \gamma\big\|\varphi_x\big(t,x_*(t),u_*(t)\big)\big\|\Big] \, dt
   + \delta\gamma.$$
Damit bildet der Operator $F$ in den Raum $C_{\lim}(\R_+,\R^n)$ ab.
Ferner k"onnen wir zu jedem $\varepsilon >0$ Zahlen $\lambda_0, T>0$ derart angeben,
dass $u(t)=u_*(t)$ f"ur alle $t\geq T$ gilt und 
$$\int_T^\infty \big\| \varphi\big(t,x_*(t)+\lambda x(t),u_*(t)\big)-\varphi\big(t,x_*(t),u_*(t)\big)
                       - \varphi_x\big(t,x_*(t),u_*(t)\big)\big(\lambda x(t)\big) \big\| \, dt \leq \frac{\varepsilon}{2}$$
f"ur alle $0<\lambda \leq \lambda_0$ und alle $\|x(\cdot)\|_\infty \leq 1$ ausf"allt. \\
Nutzen wir "uber dem Intervall $[0,T]$ die Argumentation aus Beispiel \ref{DiffDynamik2},
so erhalten wir
$$\bigg\| \frac{F\big(x_*(\cdot)+\lambda x(\cdot),u(\cdot)\big) - F\big(x_*(\cdot),u(\cdot)\big)}{\lambda}
          - F_x\big(x_*(\cdot),u(\cdot)\big) x(\cdot)\bigg\|_\infty \leq \varepsilon$$
f"ur alle $\|x(\cdot) \|_\infty \leq 1$ und alle $0 < \lambda \leq \lambda_0$. \hfill $\square$}
\end{beispiel}


\subsection{Grundprinzipien der Funktionalanalysis}
\begin{theorem}[Satz von Hahn-Banach; Fortsetzungsversion]  \index{Satz, Darstellungssatz von Riesz!von HahnBanach@-- von Hahn-Banach}
Sei $X$ ein normierter Raum und $U$ ein Untervektorraum.
Zu jedem stetigen linearen Funktional $u^*:U \to \R$ existiert dann ein stetiges lineares Funktional
$x^*: X \to \R$ mit $x^*\big|_U=u^*$ und $\|x^*\|=\|u^*\|$.
\end{theorem}

\begin{folgerung} \label{FolgerungAnnulator}
Seien $X$ ein normierter Raum, $U$ ein abgeschlossener Unterraum und $x \in X$, $x \not\in U$.
Dann existiert ein $x^* \in X^*$ mit $x^*\big|_U=0$ und $\langle x^*,x \rangle \not=0$.
\end{folgerung}

\begin{theorem}[Satz von Hahn-Banach; Trennungsversion]  \index{Satz, Darstellungssatz von Riesz!Trennungssatz@--, Trennungssatz}
Seien $X$ ein normierter Raum, $V_1,V_2 \subseteq X$ konvex, $V_1$ offen und $V_1 \cap V_2 = \emptyset$.
Dann existiert ein $x^* \in X^*$ mit
$$\langle x^*, v_1 \rangle < \langle x^*, v_2 \rangle \qquad \mbox{ f"ur alle } v_1 \in V_1, v_2 \in V_2.$$
\end{theorem}

Eine Abbildung $T$ hei"st offen\index{Abbildung, beschränkte!offene@--, offene}, wenn $T$ offene Mengen auf offene Mengen abbildet.

\begin{theorem}[Satz von der offenen Abbildung] \label{SatzOffeneAbbildung}  \index{Satz, Darstellungssatz von Riesz!von der offen@-- von der offenen Abbildung}
Es seien $X$, $Y$ Banachr"aume und $T \in L(X,Y)$ surjektiv.
Dann ist $T$ offen.
\end{theorem}

Seien $X,Y$ normierte R"aume und $T \in L(X,Y)$.
Der adjungierte Operator\index{adjungierter Operator} $T^*:Y^* \to X^*$ ist durch
$\langle T^*y^*, x \rangle = \langle y^*, Tx \rangle$
definiert. Offensichtlich folgt daraus $T^* \in L(Y^*,X^*)$.
Seien nun $U \subseteq X$ und $V \subseteq X^*$.
Wir definieren die Mengen
\begin{eqnarray*}
U^\perp &=& \{ x^* \in X^* \,|\, \langle x^*, x \rangle =0 \mbox{ f"ur alle } x \in U\}, \\
V_\perp &=& \{ x \in X \,|\, \langle x^*, x \rangle =0 \mbox{ f"ur alle } x^* \in V\}.
\end{eqnarray*}

\begin{lemma}[Satz vom abgeschlossenen Bild] \label{SatzAbgeschlossenesBild} \index{Satz, Darstellungssatz von Riesz!vom ab@-- vom abgeschlossenen Bild}
Seien $X$, $Y$ Banachr"aume, und es sei $T \in L(X,Y)$.
Dann gelten die "Aquivalenzen:
\begin{eqnarray*}
{\rm Im\,}T \mbox{ ist abgeschlossen } &\Leftrightarrow& {\rm Im\,}T=({\rm Ker\,}T^*)_\perp \\
\Leftrightarrow \; {\rm Im\,}T^* \mbox{ ist abgeschlossen} &\Leftrightarrow&  {\rm Im\,}T^*=({\rm Ker\,}T)^\perp.
\end{eqnarray*}
\end{lemma}

\begin{satz}[Fixpunktsatz von Weissinger] \label{SatzWeissinger} \index{Satz, Darstellungssatz von Riesz!Fixpunkt@--, Fixpunktsatz von Weissinger}
Es sei $U$ eine nichtleere abgeschlossene Teilmenge des Banachraumes $X$,
ferner $\displaystyle \sum_{n=1}^\infty a_n$ eine konvergente Reihe positiver Zahlen und
$A:U \to U$ eine Selbstabbildung von $U$ mit $\|A^n u-A^n v\| \leq a_n\|u-v\|$ f"ur alle $u,v \in U$ und $n \in \N$.
Dann besitzt $A$ genau einen Fixpunkt, d.\,h. es gibt genau ein $u \in U$ mit $Au=u$. \\
Dieser Fixpunkt ist Grenzwert der Iterationsfolge $u_n=Au_{n-1}, n=1,2,...,$ bei beliebigem Startwert $u_0 \in U$.
Schlie"slich gilt die Fehlerabsch"atzung
$$\|u-u_n\| \leq \|u_1-u_0\| \cdot \sum_{k=n}^\infty a_k.$$
\end{satz}

\begin{satz}[Satz über implizite Funktionen] \label{SatzImpliziteFunktionen}
Es seien $X$, $Y$ und $Z$ Banachräume, $U$ eine Umgebung des Punktes $(x_0,y_0) \in X \times Y$ und $F:U \to Z$
eine stetig differenzierbare Abbildung.
Ferner sei $F(x_0,y_0)=0$ und es sei die partielle Ableitung $F_y(x_0,y_0):Y \to Z$ ein linearer Homöomorphismus.
Dann gibt es auf einer Umgebung $U$ von $x_0$ eine stetig differenzierbare Abbildung $x \to y(x)$ nach $Y$ derart,
dass für jedes $x \in U$ die Gleichung $F\big(x,y(x)\big)=0$ und die Beziehung $y'(x)=-\big[F_y\big(x,y(x)\big)\big]^{-1} \circ F_x\big(x,y(x)\big)$ gilt.
\end{satz}

\subsection{Der Darstellungssatz von Riesz} \label{AnhangRiesz}
\begin{satz}[Rieszscher Darstellungssatz] \label{SatzRiesz}\index{Satz, Darstellungssatz von Riesz}\index{Darstellungssatz von Riesz}
Zu jedem $x^* \in C^*([t_0,t_1],\R)$ gibt es genau ein signiertes regul"ares Borelsches Ma"s $\mu$
auf der Borelschen $\sigma$-Algebra "uber $[t_0,t_1]$ mit
$$\langle x^*,x(\cdot) \rangle = \int_{t_0}^{t_1} x(t) \, d\mu(t), \qquad \|x^*\|= |\mu|([t_0,t_1]).$$
\end{satz}

\begin{folgerung} \label{FolgerungRiesz}
Jedes stetige lineare Funktional $x^* \in C^*([t_0,t_1],\R^n)$ l"asst sich auf eindeutige Weise in der Gestalt
$$\langle x^*,x(\cdot) \rangle = \langle a , x(t_0) \rangle + \sum_{k=1}^n \int_{t_0}^{t_1} x_k(t) \, d\mu_k(t)$$
darstellen, wobei $a \in \R^n$ gilt und $\mu_1(\cdot),...,\mu_n(\cdot)$ normalisierte Funktionen beschr"ankter Variation sind.
Oder "aquivalent:
$$\langle x^*,x(\cdot) \rangle =  \sum_{k=1}^n \int_{t_0}^{t_1} x_k(t) \, d\mu_k(t)$$
mit Funktionen $\mu_1(\cdot),...,\mu_n(\cdot)$ beschr"ankter Variation,
die mit eventueller Ausnahme des Punktes $t_0$ rechtsseitig stetig sind.
\end{folgerung}

\begin{folgerung} \label{FolgerungRiesz3}
Zu jedem $x^* \in C_0^*([t_0,t_1],\R)$ gilt der Rieszsche Darstellungssatz \ref{SatzRiesz},
wobei das signierte regul"are Borelsche Ma"s $\mu$ keinen atomaren Anteil in $t_0$ besitzt.
\end{folgerung}

\begin{lemma} \label{FolgerungRiesz2}
Zu jedem stetigen linearen Funktional $x^* \in C_1^*([t_0,t_1],\R^n)$ gibt es eindeutig bestimmte $a,b \in \R^n$ und
normalisierte Funktionen beschr"ankter Variation $\mu_1(\cdot),...,\mu_n(\cdot)$ derart,
dass $x^*$ folgende Darstellung besitzt:
$$\langle x^*,x(\cdot) \rangle
  = \langle a , x(t_0) \rangle + \langle b , \dot{x}(t_0) \rangle + \sum_{k=1}^n \int_{t_0}^{t_1} \dot{x}_k(t) \, d\mu_k(t).$$
\end{lemma}

\begin{satz}[Rieszscher Darstellungssatz] \label{SatzRieszC0}\index{Satz, Darstellungssatz von Riesz}\index{Darstellungssatz von Riesz}
Der Dualraum $C_0^*(\R_+,\R)$ des Raumes der stetigen Funktionen, die im Unendlichen verschwinden,
ist unter der Abbildung
$$\Lambda(\mu)x(\cdot)= \int_{\R_+} x(t) \, d\mu(t)$$
isometrisch isomorph zu den signierten regul"aren Borelschen Ma"se auf der Borelschen $\sigma$-Algebra "uber $\R_+$.
\end{satz}

An dieser Stelle erweitern wir die Darstellungss"atze um die Beispiele (h)--(k) in Abschnitt \ref{AbschnittNormRaum}.
Dazu leiten wir einen allgemeinen Darstellungssatz her:

\begin{lemma} \label{CorollaryDarstellung}
Es sei $\Omega$ eine nichtleere Menge und es sei $X_0=\{x_0(\cdot): \Omega \to \R^n\}$ ein Banachraum mit Norm $\|\cdot\|_{X_0}$.
Ferner sei $X$ der Raum derjenigen Elemente $x(\cdot)$,
die die eindeutige Darstellung $x(\cdot)=x_0(\cdot)+a$,
d.\,h. $x(\omega)=x_0(\omega)+a$ für alle $\omega \in \Omega$,
mit $x_0(\cdot) \in X_0$ und $a \in \R^n$ besitzen.
Weiterhin sei $X$ versehen mit einer $p$-Norm, d.\,h.
$$\|x(\cdot)\|^p_X = \|x_0(\cdot)\|^p_{X_0} + \|a\|^p \mbox{ f"ur } 1 \leq p < \infty, \qquad
  \|x(\cdot)\|_X = \|x_0(\cdot)\|_{X_0} + \|a\| \mbox{ f"ur } p=\infty.$$
Dann gibt es zu jedem $x^* \in X^*$ ein eindeutiges Paar $(x_0^*,\alpha) \in X_0^* \times \R^n$ mit
$$\langle x^*,x(\cdot) \rangle = \langle x_0^*, x_0(\cdot) \rangle + \alpha^T a \quad \mbox{f"ur alle } x(\cdot)=x_0(\cdot)+a \in X.$$
\end{lemma}

{\em Beweis:} Versehen mit einer $p$-Norm ist $X$ selbst ein Banachraum.
Wir betrachten die stetige lineare Abbildung $\Lambda$ mit $\Lambda x(\cdot) = x_0(\cdot)$, die $X$ auf $X_0$ abbildet.
Es gilt offenbar
$${\rm Ker\,}\Lambda=\{ x(\cdot) \in X \,|\, x(\cdot)=x_0(\cdot)+a, \mbox{ wobei } x_0(\omega)=0\; \forall \omega \in \Omega\}.$$
Sei $x^* \in X^*$.
Wir bezeichnen mit $\alpha_i$ den Wert von $x^*$ f"ur diejenige vektorwertige Funktion,
deren $i$-te Komponente identisch $1$ ist und alle weiteren identisch verschwinden.
Nun betrachten wir das Funktional $x_1^* \in X^*$, welches wie folgt definiert ist:
$$\langle x_1^*,x(\cdot) \rangle = \langle x^*,x(\cdot) \rangle - \alpha^T a, \quad x(\cdot)=x_0(\cdot)+a, \quad \alpha=(\alpha_1,...,\alpha_n).$$
Offensichtlich gilt $x_1^* \in ({\rm Ker\,}\Lambda)^\perp$.
Nach dem Satz vom abgeschlossenen Bild existiert ein Funktional $x_0^* \in X_0^*$ mit $x_1^* = \Lambda^* x_0^*$.
Es folgt
$$\langle x_1^*,x(\cdot) \rangle = \langle \Lambda^* x_0^*,x(\cdot) \rangle = \langle x_0^*, \Lambda x(\cdot) \rangle =  \langle x_0^*, x_0(\cdot) \rangle$$
f"ur alle $x(\cdot) \in X$.
Also erhalten wir f"ur $x^* \in X^*$ die Darstellung
$$\langle x^*,x(\cdot) \rangle = \langle x_1^*,x(\cdot) \rangle + \alpha^T a =  \langle x_0^*, x_0(\cdot) \rangle + \alpha^T a$$
mit $x_0^* \in X_0^*$ und $\alpha \in \R^n$.
Die Eindeutigkeit kann direkt nachgepr"uft werden. \hfill $\square$ \\[2mm]
Im Raum $C_{\lim}(\R_+,\R^n)$ unterscheiden sich die Normen $\|x(\cdot)\|_{\infty}=\|x_0(\cdot)+a\|_{\infty}$ und die in
Lemma \ref{CorollaryDarstellung} für $p=\infty$ definierte Norm $\|x(\cdot)\|_{C_{\lim}}=\|x_0(\cdot)\|_{\infty}+\|a\|$. 
Aber es gilt:

\begin{lemma}
Auf dem Raum $C_{\lim}(\R_+,\R^n)$ sind die Normen $\|\cdot\|_{C_{\lim}}$ und $\|\cdot\|_{\infty}$ "aquivalent.
\end{lemma}

{\em Beweis:} Es sei $x(\cdot) \in C_{\lim}(\R_+,\R^n)$ und $x(t)=x_0(t)+a$ mit $x_0(\cdot) \in C_0(\R_+,\R^n)$, $a \in \R^n$ 
die eindeutige Darstellung.
Sei $0 \leq \|x_0(\cdot)\|_\infty \leq \|a\|$.
Wegen $x(t) \to a$ f"ur $t \to \infty$ folgt
$$2\|x(\cdot)\|_\infty = 2\|x_0(\cdot)+a\|_\infty \geq 2 \|a\| \geq \|x_0(\cdot)\|_\infty + \|a\| = \|x(\cdot)\|_{C_{\lim}}.$$
Sei nun $0 \leq \|a\| \leq \|x_0(\cdot)\|_\infty$.
Dann gibt es ein $\lambda \in [0,1]$ mit $\|a\| = \lambda \|x_0(\cdot)\|_\infty$.
Also
\begin{eqnarray*}
\|x_0(\cdot)+a\|_\infty &\geq& \max \{\|x_0(\cdot)\|_\infty - \|a\|,\|a\|\}= \max \{ (1-\lambda) \|x_0(\cdot)\|_\infty, \lambda \|x_0(\cdot)\|_\infty\} \\
&\geq& \min_{\lambda \in [0,1]} \Big(\max\{1-\lambda,\lambda\}\Big)\cdot \|x_0(\cdot)\|_\infty = \frac{1}{2} \|x_0(\cdot)\|_\infty.
\end{eqnarray*}
Damit gelten die Ungleichungen $\|x(\cdot)\|_\infty \leq \|x(\cdot)\|_{C_{\lim}} \leq 4\|x(\cdot)\|_\infty$ auf $C_{\lim}(\R_+,\R^n)$. \hfill $\square$

\begin{folgerung}[Rieszscher Darstellungssatz für $C_{\lim}(\R_+,\R^n)$]
\label{FolgerungRieszClim} \index{Satz, Darstellungssatz von Riesz}\index{Darstellungssatz von Riesz}
Jedes stetige lineare Funktional $x^*$ auf $C_{\lim}(\R_+,\R^n)$ besitzt nach Lemma \ref{CorollaryDarstellung} die eindeutige Darstellung
$$\langle x^*(\cdot),x(\cdot) \rangle = \int_0^\infty \langle x_0(t), d\mu_0(t) \rangle + \alpha^T a,$$
wobei $\alpha \in \R^n$ gilt und $\mu_0$ ein signiertes regul"ares Borelsches Vektorma"s auf $\R_+$ ist. \\
Oder "aquivalent: Zu jedem $x^* \in C^*(\overline{\R}_+,\R^n)$ existiert ein eindeutig bestimmtes signiertes reguläres
Borelsches Ma"s $\mu$ "uber $\overline{\R}_+$ (Anhang \ref{AnhangMassR} und Definition \ref{DefinitionMassR}) mit
$$\langle x^*(\cdot),x(\cdot) \rangle =\int_0^\infty \langle x(t), d\mu(t) \rangle
  =\int_0^\infty \langle x(t), d\mu_0(t) \rangle+ x(\infty)^T \mu(\{\infty\}),\quad x(\infty)=a,$$
wobei $\mu$ die eindeutige Zerlegung $\mu=\mu_0+\mu(\{\infty\})$ mit einem signierten regul"aren Borelschen Vektorma"s $\mu_0$
auf $\R_+$ besitzt.
\end{folgerung}

Im Fall $1 \leq p < \infty$ lässt sich Lemma \ref{CorollaryDarstellung} unmittelbar auf die
Fälle $X_0=L_p(\R_+,\R^n)$ und $X_0=W^1_p(\R_+,\R^n)$ anwenden:

\begin{folgerung}
Jedes stetig lineare Funktional $x^*$ auf $L_{p,\lim}(\R_+,\R^n)$ besitzt die eindeutige Darstellung
$$\langle x^*(\cdot),x(\cdot) \rangle = \int_0^\infty \langle y_0(t), x_0(t) \rangle  \, dt + \alpha^T a,$$
wobei $\alpha \in \R^n$ und $y_0(\cdot) \in L_q(\R_+,\R^n)$ mit $1/p+1/q=1$ falls $p>1$ und $q=\infty$ falls $p=1$ gelten.
Zus"atzlich hat man $y(\cdot) = y_0(\cdot) +a \in L_{q,\lim}(\R_+,\R^n)$ f"ur $1<p<\infty$.
\end{folgerung}

\begin{folgerung}
Jedes stetig lineare Funktional $x^*$ auf $W^1_{p,\lim}(\R_+,\R^n)$ besitzt die eindeutige Darstellung
$$\langle x^*(\cdot),x(\cdot) \rangle = \int_0^\infty \big[\langle y_0(t), x_0(t) \rangle + \langle \dot{y}_0(t), \dot{x}_0(t) \rangle\big] \, dt + \alpha^T a,$$
wobei $\alpha \in \R^n$ und $y_0(\cdot) \in W^1_q(\R_+,\R^n)$ mit $1/p+1/q=1$ falls $p>1$ und $q=\infty$ falls $p=1$ gelten.
Zus"atzlich hat man $y(\cdot) = y_0(\cdot) +a \in W^1_{q,\lim}(\R_+,\R^n)$ f"ur $1<p<\infty$.
\end{folgerung}

Abschlie"send geben wir für $1 \leq p < \infty$ mit Hilfe von Lemma \ref{CorollaryDarstellung} das nachstehende Resultat
für konvergente $\ell_p$-Folgen an:

\begin{folgerung}
Jedes $x^*$ auf $\ell_{p,\lim}$ besitzt die eindeutige Darstellung
$$\langle x^*,x \rangle = \sum_{n \in N} y^0_n x^0_n+ \alpha a,$$
mit $\alpha \in \R^n$, $y^0 \in \ell_q$ und $1/p+1/q=1$ falls $p>1$ und $q=\infty$ falls $p=1$. 
Au"serdem ist $y=y^0+\alpha \in \ell_{q,\lim}$ f"ur $1<p<\infty$.
\end{folgerung}


\subsection{Der Satz von Ljusternik}
In diesem Abschnitt befassen wir uns mit dem fundamentalen Satz von Ljusternik (Ljusternik \cite{Ljusternik}).
Die vorliegende zusammenfassende Darstellung und die Verallgemeinerung ist Ioffe \& Tichomirov \cite{Ioffe} entnommen.
Eine vollst"andige Beweisf"uhrung ist wiederum in Ioffe \& Tichomirov \cite{Ioffe} zu finden.

\begin{definition}[Lokaler Tangentialkegel]
Sei $X$ ein Banachraum und $x_0 \in M \subseteq X$.
Mit $\mathscr{C}(M,x_0)$ bezeichnen wir die Menge aller Elemente $x \in X$,
zu denen ein $\varepsilon_0 > 0$ und eine Abbildung $r(\varepsilon): [0,\varepsilon_0] \to X$ mit den Eigenschaften
$$\lim_{\varepsilon \to 0^+} \frac{\| r(\varepsilon) \|}{\varepsilon} = 0
  \qquad\mbox{und}\qquad x_0 + \varepsilon x + r(\varepsilon) \in M \quad \mbox{ f"ur alle } \varepsilon \in [0,\varepsilon_0]$$
existieren. 
$\mathscr{C}(M,x_0)$ hei"st der lokale Tangentialkegel an die Menge $M$ im Punkt $x_0$.
\end{definition}

\begin{theorem}[Satz von Ljusternik] \label{SatzLjusternik}\index{Satz, Darstellungssatz von Riesz!von Ljusternik@-- von Ljusternik}
Es seien $X$ und $Y$ Banachr"aume,
$V$ eine Umgebung des Punktes $x_* \in X$ und $F$ eine Fr\'echet-differenzierbare Abbildungen der Menge $V$ in $Y$.
Wir setzen voraus,
die Abbildung $F$ sei regul"ar im Punkt $x_*$, d.\,h., es gelte
$${\rm Im\,} F'(x_*) = Y,$$
au"serdem sei ihre Ableitung in diesem Punkt in der gleichm"a"sigen Operatorentopologie des Raumes $L(X,Y)$ stetig.
Dann stimmt der lokale Tangentialkegel an die Menge \\ $M= \big\{ x \in V \,\big|\, F(x) = F(x_*) \big\}$
im Punkt $x_*$ mit dem Kern des Operators $F'(x_*)$ "uberein:
$$\mathscr{C}(M,x_*) = {\rm Ker\,} F'(x_*).$$
\end{theorem}

\begin{lemma} \label{LemmaVerallgemeinerterSatzLjusternik}
Es seien $X$ und $Y$ Banachr"aume und $\Lambda:X \to Y$ ein stetiger linearer Operator.
Wir setzen
$$C(\Lambda) := \sup_{y \not= 0} \left( \frac{\inf \{ \|x\| \,|\, x \in X, \Lambda x = y\} }{\|y\|} \right).$$
Wenn ${\rm Im\,} \Lambda = Y$ gilt, so ist $C(\Lambda) < \infty$.
\end{lemma}

\begin{theorem}[Der verallgemeinerte Satz von Ljusternik] \label{SatzVerallgemeinerterSatzLjusternik}
Es seien $X$ und $Y$ Banachr"aume,
$\Lambda: X \to Y$ ein stetiger linearer Operator und $F$ eine Abbildung einer gewissen Umgebung $V$ des Punktes $x_* \in X$ in $Y$.
Wir setzen voraus, es sei ${\rm Im\,}\Lambda = Y$ und es gebe eine Zahl $\delta > 0$ derart, dass erstens
$$\delta C(\Lambda) < \frac{1}{2}$$
und zweitens
\begin{equation} \label{SatzVerallgemeinerterSatzLjusternik1}
\| F(x) - F(x') - \Lambda(x-x') \| \leq \delta \|x-x'\|
\end{equation}
f"ur alle $x,x'$ aus $V$ gilt.
Dann existieren eine Umgebung $V' \subseteq V$ des Punktes $x_*$,
eine Zahl $K>0$ und eine Abbildung $\xi \to x(\xi)$ der Umgebung $V'$ in $X$
derart, dass die Beziehungen
$$F\big( \xi + x(\xi) \big) = F(x_*) \qquad\mbox{ und }\qquad \| x(\xi) \| \leq K \| F(\xi) - F(x_*) \|$$
f"ur alle $\xi \in V'$ erf"ullt sind.
\end{theorem}

%% file: III-Differentialgleichung.tex
\section{Differentialgleichungen, Volterrasche Integralgleichungen} \label{AnhangDGL}
Ein grundlegende Einf"uhrung zu gew"ohnlichen Differentialgleichungen liefert Heuser \cite{HeuserGD}.
Bei Differentialgleichungen mit messbaren rechten Seiten verweisen wir auf Filippov \cite{Filippov}. \\
Einführende Darstellungen zu Differentialgleichungen beruhen zumeist auf Gleichungen mit rechten Seiten,
die stetig in der Gesamtheit der Variablen sind. 
Da in Steuerungsproblemen die rechte Seite der Dynamik für einen Steuerungsprozess typischerweise nicht stetig in der
Zeitvariable ist, und somit die Rahmenbedingungen von der Standardtheorie abweichen,
geben wir für die benötigten Resultate die Beweise an.
Dabei folgen wir Ioffe \& Tichomirov \cite{Ioffe} in den Betrachtungen von Differentialgleichungen mit messbaren
rechten Seiten.
Au"serdem erweitern wir die Theorie um notwendige Ergebnisse zu Differentialgleichungen über dem unbeschr"ankten
Intervall $\R_+$ und um Ergebnisse zu Volterraschen Integralgleichungen 2.\,Art.

\subsection{Lineare Gleichungen}
Im Weiteren betrachten wir bezüglich $x(\cdot)$ lineare Integralgleichungssysteme der Form
\begin{eqnarray}
\label{LinDGL1} \zeta(t) &=& x(t) - \int_\tau^t \big[A(t,s)x(s)+a(s)\big] \, ds, \\
\label{LinDGL2} \zeta(t) &=& x(t) - \int_\tau^t \big[A(s)x(s)+a(s)\big] \, ds
\end{eqnarray}
und das lineare Differentialgleichungssystem
\begin{equation}\label{LinDGL3}
\dot{x}(t)=A(t)x(t)+a(t).
\end{equation}
Offensichtlich sind (\ref{LinDGL3}) ein Spezialfall von (\ref{LinDGL2}) und
(\ref{LinDGL2}) ein Spezialfall von (\ref{LinDGL1}).
Demzufolge basiert die eindeutige Lösbarkeit der Gleichungen (\ref{LinDGL2}) und (\ref{LinDGL3}) auf den Betrachtungen
zur Gleichung (\ref{LinDGL1}) und dem zentralen Ergebnis in Form von Lemma \ref{LemmaDGL1}. \\
Weiterhin beziehen sich die weiteren Ergebnisse, die aus Lemma \ref{LemmaDGL1} erzielt werden,
auf die eindeutige Lösbarkeit der Gleichungen (\ref{LinDGL1})--(\ref{LinDGL3}) im Rahmen der Räume 
$C([t_0,t_1],\R^n)$, $C_{\lim}(\R_+,\R^n)$ bzw. $L_\infty([t_0,t_1],\R^n)$.

\begin{bemerkung}{\rm \label{BemDGL}
Ist $\zeta(\cdot)$ in $\tau$ stetig,
so muss eine Lösung $x(\cdot)$ der Gleichung (\ref{LinDGL1}) in $t=\tau$ stetig sein und $x(\tau)=\zeta(\tau)$ gelten.
Daher erfüllt $x(\cdot)$ ebenso die Gleichung
$$\zeta(t) - \zeta(\tau) = x(t) - x(\tau) - \int_\tau^t \big[A(t,s)x(s)+a(s)\big] \, ds.$$
Folglich zieht die Lösbarkeit der Gleichung (\ref{LinDGL1}) die Existenz einer Lösung der Gleichung
$$\zeta_0(t) = x(t) - x(\tau) - \int_\tau^t \big[A(t,s)x(s)+a(s)\big] \, ds$$
mit einer in $t=\tau$ stetigen Funktion $\zeta_0(\cdot)$, für die $\zeta_0(\tau)=0$ gilt, nach sich.
\hfill $\square$}
\end{bemerkung}

\begin{lemma} \label{LemmaDGL1}
Es sei die Abbildung $A(t,s)$ für $t_0 \leq s, t\leq t_1$ in $t$ stetig und in $s$ messbar.
Ferner existiere eine über $[t_0,t_1]$ integrierbare Funktion $c(\cdot)$ derart,
dass $\|A(t,s)\| \leq c(s)$ für fast alle $(t,s)$ mit $t_0 \leq s, t\leq t_1$ gilt.
Weiter sei die Vektorfunktion $a(\cdot)$ "uber $[t_0,t_1]$ integrierbar.
Dann existiert zu jeder Vektorfunktion $\zeta(\cdot) \in C([t_0,t_1],\R^n)$ und jedem $\tau \in [t_0,t_1]$
eine eindeutig bestimmte Lösung $x(\cdot) \in C([t_0,t_1],\R^n)$ der Gleichung
$$\zeta(t)=x(t) - \int_\tau^t \big[A(t,s)x(s)+a(s)\big] \, ds \quad \mbox{ für alle } t \in [t_0,t_1].$$
\end{lemma}

{\bf Beweis} Wir werden im Folgenden zeigen,
dass die Fixpunktgleichung $x(\cdot) = T\big(x(\cdot)\big)$, wobei der Operator $T$ durch
$$x(\cdot) \to T\big(x(\cdot)\big), \quad
  T\big(x(\cdot)\big)(t) = \zeta(t) + \int_\tau^t \big[A(t,s)x(s)+a(s)\big] \, ds, \quad t \in [t_0,t_1],$$
gegeben wird, stets eine eindeutige L"osung besitzt.
Da $A(t,s)$ in $t$ stetig ist, bildet der Operator $T$ den Raum $C([t_0,t_1],\R^n)$ in sich ab.
Mit der integrierbaren Funktion $c(\cdot)$ in der Voraussetzung des Satzes seien zur abk"urzenden Schreibweise
$$C(t) = \int_\tau^t c(s) \, ds, \qquad c_0=\int_{t_0}^{t_1} c(s) \, ds.$$
Bei mehrfacher Anwendung des Operators $T$ ergeben sich f"ur $x_1(\cdot),x_2(\cdot) \in C([t_0,t_1],\R^n)$ die Beziehungen
\begin{eqnarray*}
\lefteqn{\big\| \big[T\big(x_1(\cdot) - x_2(\cdot)\big)\big](t) \big\|
         \leq \int_\tau^t c(s) \, ds \cdot \| x_1(\cdot) - x_2(\cdot) \|_\infty,} \\
\lefteqn{\big\| \big[T^2\big(x_1(\cdot) - x_2(\cdot)\big)\big](t) \big\|
         \leq \int_\tau^t c(s) \big\| \big[T\big(x_1(\cdot) - x_2(\cdot)\big)\big](s) \big\| \, ds} \\
&& \hspace*{10mm} \leq \int_\tau^t c(s) C(s) \, ds \cdot \big\| x_1(\cdot) - x_2(\cdot) \big\|_\infty
   = \frac{1}{2} C^2(t) \cdot \| x_1(\cdot) - x_2(\cdot) \|_\infty, \\
&& \hspace*{30mm}\vdots \\
\lefteqn{\big\| \big[T^m \big(x_1(\cdot) - x_2(\cdot)\big)\big](t) \big\|
         \leq \int_\tau^t c(s) \big\| \big[T^{m-1}\big(x_1(\cdot) - x_2(\cdot)\big)\big](s) \big\| \, ds} \\
&& \hspace*{10mm} \leq \int_\tau^t c(s) \frac{C^{m-1}(s)}{(m-1)!} \, ds \cdot \big\| x_1(\cdot) - x_2(\cdot) \big\|_\infty
       =\frac{C^m(t)}{m!} \cdot \| x_1(\cdot) - x_2(\cdot) \|_\infty.
\end{eqnarray*}
In der Topologie des Raumes $C([t_0,t_1],\R^n)$ gilt daher
$$\big\|T^m \big(x_1(\cdot) - x_2(\cdot)\big) \big\|_\infty \leq \frac{c_0^m}{m!} \cdot \| x_1(\cdot) - x_2(\cdot) \|_\infty.$$
Die Zahlen $a_m = c_0^m/m!$ liefern eine Folge, deren Reihe konvergiert.
Nach dem Fixpunktsatz von Weissinger (Satz \ref{SatzWeissinger}) existiert daher genau ein
$x(\cdot)$ mit $x(\cdot) = T x(\cdot)$.
\hfill $\blacksquare$

\begin{lemma} \label{LemmaDGL2}
Unter den Annahmen in Satz \ref{LemmaDGL1} gibt es zu jedem $\zeta(\cdot) \in L_\infty([t_0,t_1],\R^n)$
und jedem $\tau \in [t_0,t_1]$
eine eindeutig bestimmte Lösung $x(\cdot) \in L_\infty([t_0,t_1],\R^n)$ der Gleichung
$$\zeta(t)=x(t) - \int_\tau^t \big[A(t,s)x(s)+a(s)\big] \, ds \quad \mbox{ für fast alle } t \in [t_0,t_1].$$
\end{lemma}

{\bf Beweis} Zum Nachweis von Lemma \ref{LemmaDGL2} ist im Beweis von Lemma \ref{LemmaDGL1} der Operator $T$ als Abbildung
von $L_\infty([t_0,t_1],\R^n)$ in sich aufzufassen.
\hfill $\blacksquare$

\begin{lemma} \label{LemmaDGL3}
Es seien die Abbildung $t \to A(t)$ und die Vektorfunktion $a(\cdot)$ "uber $[t_0,t_1]$ integrierbar.
Dann existiert zu jedem $\zeta(\cdot) \in C([t_0,t_1],\R^n)$ und jedem $\tau \in [t_0,t_1]$
eine eindeutig bestimmte Lösung $x(\cdot) \in C([t_0,t_1],\R^n)$ der Gleichung
$$\zeta(t)=x(t) - \int_\tau^t \big[A(s)x(s)+a(s)\big] \, ds \quad \mbox{ für alle } t \in [t_0,t_1].$$
\end{lemma}

{\bf Beweis} Mit der Setzung $c(t) = \|A(t)\|$ folgt dies unmittelbar aus Lemma \ref{LemmaDGL1}.
\hfill $\blacksquare$

\begin{lemma} \label{LemmaDGL4}
Es seien die Abbildung $t \to A(t)$ und die Vektorfunktion $a(\cdot)$ "uber $\R_+$ integrierbar.
Dann existiert zu jedem $\zeta(\cdot) \in C_{\lim}(\R_+,\R^n)$ und jedem $\tau \in \overline{\R}_+$
eine eindeutig bestimmte Lösung $x(\cdot) \in C_{\lim}(\R_+,\R^n)$ der Gleichung
$$\zeta(t)=x(t) - \int_\tau^t \big[A(s)x(s)+a(s)\big] \, ds \quad \mbox{ für alle } t \in \overline{\R}_+.$$
\end{lemma}

{\bf Beweis} Wir ersetzen im Beweis von Satz \ref{LemmaDGL1} das Intervall $[t_0,t_1]$ durch $\R_+$.
Ferner seien der Operator $T$ wie dort definiert und
$$c(t) = \|A(t)\|, \qquad C(t) = \int_\tau^t c(s) \, ds, \qquad c_0 = \int_0^\infty c(s) \, ds.$$
Bei mehrfacher Anwendung des Operators $T$ ergibt sich f"ur $x_1(\cdot),x_2(\cdot) \in C_{\lim}(\R_+,\R^n)$:
$$\big\|T^m \big(x_1(\cdot) - x_2(\cdot)\big) \big\|_\infty \leq \frac{c_0^m}{m!} \cdot \| x_1(\cdot) - x_2(\cdot) \|_\infty.$$
Die Zahlen $a_m = c_0^m / m!$ liefern eine Folge, deren Reihe konvergiert.
Nach dem Fixpunktsatz von Weissinger existiert daher genau ein $x(\cdot)$ mit $x(\cdot) = T x(\cdot)$.
\hfill $\blacksquare$

\begin{lemma} \label{LemmaDGL5}
Es seien die Voraussetzungen in Lemma \ref{LemmaDGL3} erfüllt und $x(\cdot) \in C([t_0,t_1],\R^n)$ die eindeutige Lösung.
Dann hat $x(\cdot)$ eine verallgemeinerte Ableitung $\dot{x}(\cdot) \in L_1([t_0,t_1],\R^n)$ und
$\dot{x}(\cdot)$ genügt für fast alle $t \in [t_0,t_1]$ der Gleichung
$$\dot{x}(t)=A(t)x(t)+a(t), \qquad x(\tau)=\zeta.$$
\end{lemma}

{\bf Beweis} Die Behauptung folgt unmittelbar aus Lemma \ref{LemmaDGL3} mit $\zeta(t) \equiv \zeta$. \hfill $\blacksquare$

\begin{lemma} \label{LemmaDGL6}
Es seien die Voraussetzungen in Lemma \ref{LemmaDGL4} erfüllt und $x(\cdot) \in C_{\lim}(\R_+,\R^n)$ die eindeutige Lösung.
Dann besitzt $x(\cdot)$ eine verallgemeinerte Ableitung $\dot{x}(\cdot) \in L_1(\R_+,\R^n)$ und
$\dot{x}(\cdot)$ genügt für fast alle $t \in \R_+$ der Gleichung
$$\dot{x}(t)=A(t)x(t)+a(t), \qquad x(\tau)=\zeta, \quad x(\infty):=\lim_{t \to \infty} x(t).$$
\end{lemma}

{\bf Beweis} Die Behauptung folgt aus Lemma \ref{LemmaDGL4},
wenn wir $\zeta(t) \equiv \zeta$ setzen. \hfill $\blacksquare$


\subsection{Existenz und Eindeutigkeit, Abhängigkeit von Anfangsdaten} \label{AbschnittDGLAbhaengigkeit}
Wir widmen unser Interesse nun den nichtlinearen Differentialgleichungen der Form
\begin{equation}\label{DGL1}
\dot{x}(t) = \varphi\big(t,x(t)\big)
\end{equation}
und den nichtlinearen Volterraschen Integralgleichungen der Gestalt
\begin{equation}\label{DGL2}
x(t)=x(t_0)+ \int_{t_0}^t \varphi\big(t,s,x(s)\big) \, ds.
\end{equation}
Da sich die Differentialgleichung (\ref{DGL1}) in der äquivalenten Form
\begin{equation}\label{DGL3}
x(t)=x(t_0)+ \int_{t_0}^t \varphi\big(s,x(s)\big) \, ds
\end{equation}
als Spezialfall von (\ref{DGL2}) erweist, konzentrieren wir uns auf die Integralgleichung.
Dementsprechend sind für die Gleichung (\ref{DGL3}) die nachstehenden Resultate lediglich ohne der äußeren Zeitvariable $t$
zu formulieren und zu übernehmen. \\[2mm]
Wir treffen für die weiteren Untersuchungen folgende Festlegungen:
\begin{enumerate}
\item[(1)] Für $\gamma >0$ und $x(\cdot) \in C([a,b],\R^n)$ sei $\mathscr{R}$ der kompakte Umgebungsstreifen 
           $$\{(t,s,x) \in \R \times \R \times \R^n \,|\, a \leq s, t\leq b, \, \|x-x(t)\| \leq \gamma \}.$$
           Ferner bezeichne $U_\delta(x) \subset \R^n$ die offene Kugel $\{z \in \R^n \,|\, \|z-x\| < \delta\}$.
\item[(2)] Die Abbildung $\varphi(t,s,x)$ sei auf der Menge $\mathscr{R}$ in der Variable $x$ gleichmäßig stetig und
           gleichmäßig stetig differenzierbar, d.\,h. es gibt eine positive Konstante $L$ mit
           $$\|\varphi(t,s,x_1)-\varphi(t,s,x_2)\| \leq L \|x_1-x_2\|, \quad \|\varphi_x(t,s,x_1)-\varphi_x(t,s,x_2)\| \leq L \|x_1-x_2\|$$
           für alle $(t,s,x_1),\, (t,s,x_2) \in \mathscr{R}$.
\item[(3)] Die Abbildungen $\varphi(t,s,x)$ und $\varphi_x(t,s,x)$ seien auf der Menge $\mathscr{R}$ in der Variable $t$ stetig
           und in der Variable $s$ messbar.
\item[(4)] Die Abbildungen $\varphi(t,s,x)$ und $\varphi_x(t,s,x)$ seien über $\mathscr{R}$ beschränkt und es gelten
           $\|\varphi(t,s,x)\| \leq M,\, \|\varphi_x(t,s,x)\| \leq M$ über $\mathscr{R}$ mit der positiven Konstanten $M$.
\end{enumerate}

\begin{satz}[Lokaler Fixpunkt] \label{SatzFixpunktlokal}
Zu jedem inneren Punkt $(t_0,t_0,\xi)$ der Menge $\mathscr{R}$ gibt es Zahlen $\varepsilon >0$ und $\delta >0$ derart,
dass über dem Intervall $I=[t_0-\varepsilon,t_0+\varepsilon]$ die Abbildung
$$A: C(I,\R^n) \to C(I,\R^n), \quad
  [Ax(\cdot)](t)= \xi(t) + \int_{t_0}^t \varphi\big(t,s,x(s)\big) \, ds, \; t \in I,$$
für jedes $\xi(\cdot) \in C(I,\R^n)$ mit $\|\xi(\cdot)-\xi\|_\infty < \delta$ genau einen Fixpunkt $x_\xi(\cdot) \in C(I,\R^n)$ besitzt.
Unter diesen Voraussetzungen gilt ferner für die Fixpunkte $x_{\xi_1}(\cdot), \,x_{\xi_2}(\cdot) \in C(I,\R^n)$ der Abbildung $A$
zu den Funktionen $\xi_1(\cdot),\,\xi_2(\cdot) \in C(I,\R^n)$ die Ungleichung
$$\|x_{\xi_1}(\cdot)-x_{\xi_2}(\cdot)\|_\infty \leq  \frac{1}{1-\varepsilon L}\|\xi_1(\cdot)-\xi_2(\cdot)\|_\infty.$$
\end{satz}

{\bf Beweis} Wir wählen die positiven Zahlen $\varepsilon$ und $\delta$ so, dass 
\begin{equation} \label{DGLoffen}
\{(t,s,x) \in \R \times \R \times \R^n \,|\, a < t_0-\varepsilon \leq s, t\leq t_0+\varepsilon <b, \, \|x-\xi\| < 2\delta \}\subset \mathscr{R}
\end{equation}
gelten und weiterhin die Beziehungen $\varepsilon L < 1$ und $\varepsilon M <\delta$ erfüllt sind. \\
Wir betrachten im Raum $C(I,\R^n)$ die abgeschlossene und konvexe Teilmenge
$$K=\{ x(\cdot) \in C(I,\R^n) \,|\, \|x(t)-\xi\|\leq 2\delta \mbox{ für alle } t \in I = [t_0-\varepsilon,t_0+\varepsilon]\}.$$
Es sei $\xi(\cdot) \in C(I,\R^n)$ mit $\|\xi(\cdot)-\xi\|_\infty < \delta$.
Für $x(\cdot) \in K$ und $t \in I$ gilt dann
\begin{eqnarray*}
        \|[Ax(\cdot)](t) -\xi\|
&\leq&  \|[Ax(\cdot)](t) -\xi(t)\| + \|\xi(t)-\xi\| < \bigg\| \int_{t_0}^t \varphi\big(t,s,x(s)\big) \, ds \bigg\| + \delta \\
&\leq& \int_{t_0}^t M \, ds + \delta \leq \varepsilon M + \delta < 2\delta,
\end{eqnarray*}
d.\,h., die Abbildung $A$ bildet die Menge $K$ in sich ab.
Ferner gilt für $x_1(\cdot),\, x_2(\cdot) \in K$
\begin{eqnarray*}
    \|Ax_1(\cdot) - Ax_2(\cdot)\|_\infty
&=& \max_{t \in I} \bigg\| \int_{t_0}^t \big[\varphi\big(t,s,x_1(s)\big)-\varphi\big(t,s,x_2(s)\big)\big] \, ds \bigg\| \\
&\leq& \max_{t \in I} \int_{t_0}^t L \|x_1(s)-x_2(s)\| \, ds
       \leq \varepsilon L \cdot \big\| x_1(\cdot)- x_2(\cdot)\big\|_\infty,
\end{eqnarray*}
d.\,h., die Abbildung $A$ ist wegen $\varepsilon L < 1$ kontrahierend.
Nach dem Banachschen Fixpunktsatz besitzt die Abbildung $A$ einen eindeutigen Fixpunkt $x(\cdot) \in K$. \\ 
Es seien die Funktion $\xi_1(\cdot),\,\xi_2(\cdot)  \in C(I,\R^n)$ mit
$\|\xi_1(\cdot)-\xi\|_\infty \leq \delta$, $\|\xi_2(\cdot)-\xi\|_\infty \leq \delta$ und 
es seien $x_{\xi_1}(\cdot), \,x_{\xi_2}(\cdot) \in C(I,\R^n)$ die zugehörigen Fixpunkte der Abbildung $A$.
Dann gilt
\begin{eqnarray*}
    \|x_{\xi_1}(\cdot)-x_{\xi_2}(\cdot)\|_\infty
&=& \max_{t \in I}
    \bigg\| \xi_1(t)-\xi_2(t)+\int_{t_0}^t \big[\varphi\big(t,s,x_{\xi_1}(s)\big)-\varphi\big(t,s,x_{\xi_2}(s)\big)\big] \, ds \bigg\| \\
&\leq&  \|\xi_1(\cdot)-\xi_2(\cdot)\|_\infty + \varepsilon L \cdot \big\|x_{\xi_1}(\cdot)-x_{\xi_2}(\cdot)\big\|_\infty.
\end{eqnarray*}
Daraus ergibt sich abschlie"send
$$\big\|x_{\xi_1}(\cdot)-x_{\xi_2}(\cdot)\big\|_\infty \leq  \frac{1}{1-\varepsilon L}\|\xi_1(\cdot)-\xi_2(\cdot)\|_\infty.$$
Der Satz \ref{SatzFixpunktlokal} ist damit vollständig nachgewiesen. \hfill $\blacksquare$

\begin{satz}[Lokaler Existenz-, Eindeutigkeits- und Stetigkeitssatz] \label{SatzEElokal}
Zu jedem inneren Punkt $(t_0,t_0,\xi)$ der Menge $\mathscr{R}$ gibt es Zahlen $\varepsilon >0$ und $\delta >0$,
so dass zu jedem $z \in U_\delta(\xi)$ genau eine stetige L"osung $x_z(\cdot)$ der Gleichung (\ref{DGL2}) zur Anfangsbedingung $x(t_0)=z$
auf dem Intervall $I=[t_0-\varepsilon,t_0+\varepsilon]$ existiert. \\
Ist ferner $\{z_k\}$ eine Folge aus $U_\delta(\xi)$, die gegen $z \in U_\delta(\xi)$ konvergiert,
so gilt
$$\max_{t \in I} \big\|x_{z_k}(t)-x_z(t)\big\|_\infty \leq  \frac{1}{1-\varepsilon L}\|z_k-z\|$$
und $\{x_{z_k}(\cdot)\}$ konvergiert auf $I$ gleichm"a"sig gegen $x_z(\cdot)$.
\end{satz}

{\bf Beweis} Nach Satz \ref{SatzFixpunktlokal} können wir $\varepsilon > 0$ und $\delta >0$ so wählen,
dass die Abbildung 
$$[Ax(\cdot)](t)= z + \int_{t_0}^t \varphi\big(t,s,x(s)\big) \, ds, \quad t \in I=[t_0-\varepsilon,t_0+\varepsilon],$$
für jedes $z \in U_\delta(\xi)$ genau einen Fixpunkt $x_z(\cdot) \in C(I,\R^n)$ besitzt, d.\,h.
$$Ax_z(\cdot) = x_z(\cdot) \quad\Leftrightarrow\quad 
  x_z(t)=z + \int_{t_0}^t \varphi\big(t,s,x_z(s)\big) \, ds \mbox{ für alle } t \in I=[t_0-\varepsilon,t_0+\varepsilon].$$
Die Existenz und Eindeutigkeit der Lösung ist damit gezeigt.\\
Weiter sei $\delta >0$ so gewählt, dass (\ref{DGLoffen}) erfüllt ist.
Es seien nun $z_1,z_2 \in U_\delta(\xi)$.
Dann existieren die Lösungen $x_{z_1}(\cdot)$, $x_{z_2}(\cdot)$ der Gleichung (\ref{DGL2})
über $I$ und nehmen Werte in der Menge $U_{2\delta}(\xi)$ an.
Nach Satz \ref{SatzFixpunktlokal} ist die Ungleichung
$$\big\|x_{z_1}(\cdot)-x_{z_2}(\cdot)\big\|_\infty \leq  \frac{1}{1-\varepsilon L}\|z_1-z_2\|$$
erfüllt,
welche die gleichmäßige Konvergenz von $x_{z_k}(\cdot)$ gegen $x_z(\cdot)$ für $z_k \to z$ bedeutet.
Satz \ref{SatzEElokal} ist somit vollständig bewiesen. \hfill $\blacksquare$ \\[2mm]
Ist $x(\cdot)$ über $[t_0,t_1]$ eine Lösung der Gleichung (\ref{DGL3}),
so lässt sich $x(\cdot)$ formal durch 
$$x(t)=x(t_1)+ \int_{t_1}^t \varphi\big(s,x(s)\big) \, ds$$
zu einer Lösung über $[t_0,t_1+\varepsilon]$ fortsetzen.
Im Fall von Gleichung (\ref{DGL2}) führt
$$x(t)=x(t_1)+ \int_{t_1}^t \varphi\big(t,s,x(s)\big) \, ds$$
im Allgemeinen nicht zu der gewünschten Fortsetzung.
Stattdessen erfolgt eine Fortsetzung der Gleichung (\ref{DGL2}) an der Stelle $t=t_1$ nicht mit dem Anfangswert $x(t_1)$,
sondern mit einer Funktion $\xi(t)$, für die $\xi(t_1)=x(t_1)$ gilt.
Dieser Sachverhalt ist der Gegenstand der nachstehen Untersuchung:

\begin{lemma}[Lokale Fortsetzbarkeit] \label{LemmaDGLFortsetzung}
Es sei $x(\cdot)$ über $[t_0,t_1] \subset (a,b)$ eine Lösung von (\ref{DGL2}) und es sei $\big(t_1,t_1,x(t_1)\big)$
ein innerer Punkt der Menge $\mathscr{R}$.
Dann existiert eine Zahl $\varepsilon >0$ derart,
dass $x(\cdot)$ eindeutig zu einer stetigen Lösung der Gleichung (\ref{DGL2}) über dem Intervall $[t_0,t_1+\varepsilon]$ fortgesetzt werden kann.
Die Fortsetzung ergibt sich für $t \in [t_1,t_1+\varepsilon]$ als die eindeutige Lösung der Gleichung
\begin{equation}\label{DGLFortsetzung}
x(t) = \displaystyle \xi(t)+ \int_{t_1}^t \varphi\big(t,s,x(s)\big) \, ds \quad\mbox{mit}\quad
\xi(t)=x(t_0)  + \int_{t_0}^{t_1} \varphi\big(t,s,x(s)\big) \, ds.
\end{equation}
\end{lemma}

{\bf Beweis} Als Lösung der Gleichung (\ref{DGL2}) gilt für $x(\cdot)$ in $t=t_1$ die Beziehung
$$x(t_1)= x(t_0) + \int_{t_0}^{t_1} \varphi\big(t_1,s,x(s)\big) \, ds.$$
Damit bringen wir die Abbildung $t \to \xi(t)$ in die Form
\begin{equation} \label{DGLUmformung}
\xi(t) = x(t_1)+ \int_{t_0}^{t_1} \big[ \varphi\big(t,s,x(s)\big) - \varphi\big(t_1,s,x(s)\big) \big] \, ds.
\end{equation}
Nun sieht man unmittelbar,
dass $t \to \xi(t)$ für hinreichend kleine $\varepsilon >0$ auf $[t_1,t_1+\varepsilon]$ stetig ist und $\xi(t_1)=x(t_1)$ gilt.
Daher kann zu jedem $\delta >0$ ein $\varepsilon >0$ gewählt werden, dass $\|\xi(t)-x(t_1)\| < \delta$ über $[t_1,t_1+\varepsilon]$ erfüllt ist.
Nach Satz \ref{SatzFixpunktlokal} besitzt die Abbildung
$$[Ay(\cdot)](t)= \xi(t) + \int_{t_1}^t \varphi\big(t,s,y(s)\big) \, ds, \qquad t \in [t_1,t_1+\varepsilon],$$
für hinreichend kleine $\varepsilon > 0$ und $\delta >0$ genau einen Fixpunkt $y(\cdot) \in C([t_1,t_1+\varepsilon],\R^n)$.
Es gilt $y(t_1)=x(t_1)$.
Setzen wir zudem $y(t)=x(t)$ für $t \in [t_0,t_1]$,
so ist $y(\cdot) \in C([t_0,t_1+\varepsilon],\R^n)$ eine Fortsetzung der Lösung $x(\cdot)$ von (\ref{DGL2}).
Denn offenbar erfüllt $y(\cdot)$ für jedes $t \in [t_0,t_1]$ die Gleichung (\ref{DGL2}).
Für $t \in (t_1,t_1+\varepsilon]$ folgt
\begin{eqnarray*}
y(t) &=& \xi(t) + \int_{t_1}^t \varphi\big(t,s,y(s)\big) \, ds
         = x(t_0)  + \int_{t_0}^{t_1} \varphi\big(t,s,x(s)\big) \, ds + \int_{t_1}^t \varphi\big(t,s,y(s)\big) \, ds \\
     &=& x(t_0)  + \int_{t_0}^{t_1} \varphi\big(t,s,y(s)\big) \, ds + \int_{t_1}^t \varphi\big(t,s,y(s)\big) \, ds
         = x(t_0 ) + \int_{t_0}^t \varphi\big(t,s,y(s)\big) \, ds.
\end{eqnarray*}
Damit löst $y(\cdot)$ die Gleichung $\displaystyle y(t) = x(t_0 ) + \int_{t_0}^t \varphi\big(t,s,y(s)\big) \, ds$ für alle $t \in [t_0,t_1+\varepsilon]$. \\
Ist umgekehrt $\zeta(\cdot)$ eine fortgesetzte Lösung der Gleichung (\ref{DGL2}) von $x(\cdot)$,
d.\,h. $\zeta(t)=x(t)$ für alle $t \in [t_0,t_1]$,
so gilt für $t > t_1$:
$$\zeta(t) = x(t_0 ) + \int_{t_0}^t \varphi\big(t,s,\zeta(s)\big) \, ds
           = \underbrace{x(t_0 ) + \int_{t_0}^{t_1} \varphi\big(t,s,x(s)\big) \, ds}_{=\xi(t)} + \int_{t_1}^t \varphi\big(t,s,\zeta(s)\big) \, ds.$$
Damit wird $x(\cdot)$ nur durch die eindeutige Lösung von (\ref{DGLFortsetzung}) fortgesetzt. \hfill $\blacksquare$ 

\begin{satz}[Globaler Existenz-, Eindeutigkeits- und Stetigkeitssatz] \label{SatzEEglobal}
Über $[t_0,t_1]$ sei $x(\cdot)$ eine stetige Lösung der Gleichung (\ref{DGL2}) zum Anfangswert $x(t_0)=x_0$,
für die die Menge $\{(t,s,x) \in \R \times \R \times \R^n \,|\, t_0 \leq s, t\leq t_1, \, x=x(s) \}$ dem Inneren der Menge $\mathscr{R}$ angehört. \\
Dann existiert ein $\delta>0$ derart,
dass für jedes $z \in U_\delta(x_0)$ die Gleichung
\begin{equation} \label{AnhangDGLMax}
x(t) = z+ \int_{t_0}^t \varphi\big(t,s,x(s)\big) \, ds
\end{equation}
eine eindeutige stetige Lösung $x_z(\cdot)$ über $[t_0,t_1]$ besitzt und
die Funktionen $x_z(\cdot)$ gleichmäßig gegen $x(\cdot)$ über $[t_0,t_1]$ für $z \to x_0$ konvergieren.
\end{satz}

{\bf Beweis} Nach Satz \ref{SatzEElokal} gibt es ein $\delta >0$ und ein $\tau >t_0$ derart,
dass für jedes $z \in U_\delta(x_0)$ die Gleichung (\ref{AnhangDGLMax})
eine Lösung über $[t_0,\tau]$ besitzt und die Funktionen $x_z(\cdot)$ gleichmäßig gegen $x(\cdot)$ über $[t_0,\tau]$ für $z \to x_0$
konvergieren.
Im Weiteren sei $\overline{\tau}$ durch das Supremum über alle $\tau > t_0$ mit dieser Eigenschaft definiert.
D.\,h. es gibt ein $\delta >0$ derart,
dass für jedes $z \in U_\delta(x_0)$ eine Lösung $x_z(\cdot)$ von (\ref{AnhangDGLMax}) über $[t_0,\tau]$ existiert und
die Funktionen $x_z(\cdot)$ gleichmäßig gegen $x(\cdot)$ über $[t_0,\tau]$ für $z \to x_0$ konvergieren.
Zum Beweis des Satzes ist es ausreichend, die Beziehung $\overline{\tau}>t_1$ nachzuweisen. \\[1mm]
Angenommen, es ist $\overline{\tau} \leq t_1$.
Dann ist $\big(\overline{\tau},\overline{\tau},x(\overline{\tau}) \big)$ ein innerer Punkt der Menge $\mathscr{R}$.
Nach Satz \ref{SatzEElokal} lässt sich zu jedem $\varrho >0$ ein $\delta >0$ derart angeben,
dass $\|x_z(\cdot)-x(\cdot)\|_\infty \leq \varrho$ für alle $z \in U_\delta(x_0)$ gilt.
Insbesondere konvergiert $x_z(\cdot)$ gegen $x(\cdot)$ gleichmäßig für $z \to x_0$
und es ist $x_z(\overline{\tau}) \in U_\varrho\big(x(\overline{\tau})\big)$. \\

\begin{figure}[h]
	\centering
	\fbox{\includegraphics[height=5cm]{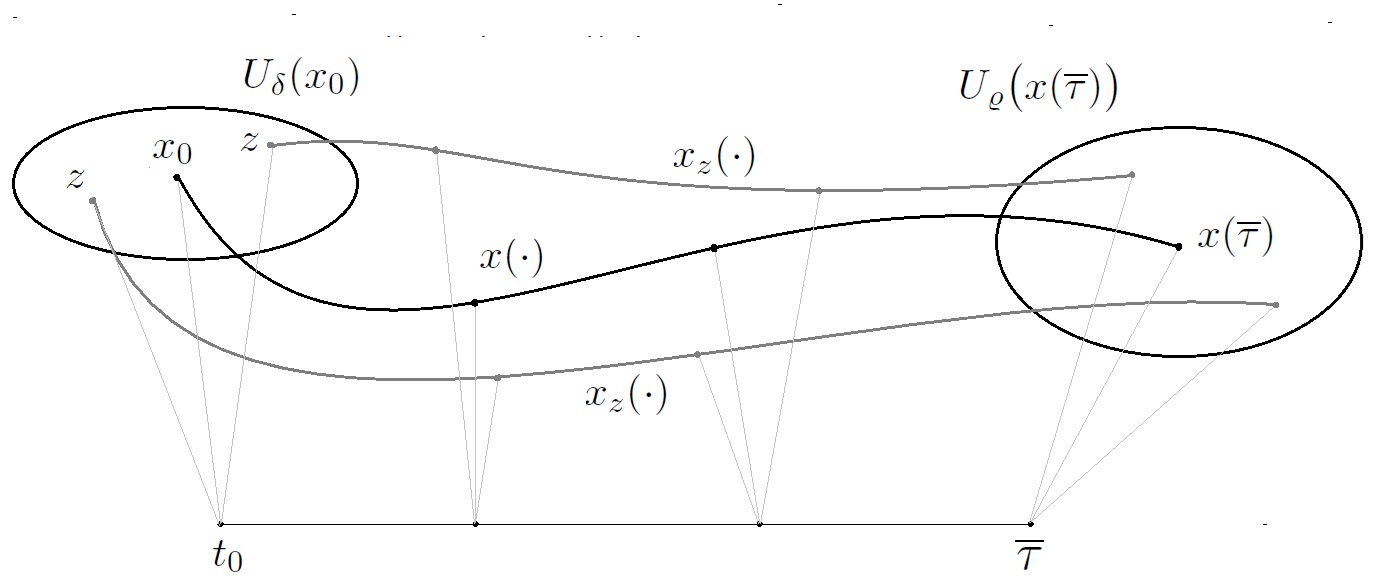}}
	\caption[Globale Existenz und Stetigkeit bei Differentialgleichungen]{Für jedes $z \in U_\delta(x_0)$
                                                                          gehört $x_z(\overline{\tau})$ zu $U_\varrho\big(x(\overline{\tau})\big)$.}
\end{figure}

Sei zunächst $x_z(\cdot)$ die stetige Lösung von (\ref{AnhangDGLMax}) über $[t_0,\overline{\tau}]$ zu einem festen $z \in U_\delta(x_0)$.
Dann gibt es nach Lemma \ref{LemmaDGLFortsetzung} eine Zahl $\varepsilon > 0$ und eine stetige Fortsetzung $y_z(\cdot)$ von $x_z(\cdot)$,
welche die Gleichung (\ref{AnhangDGLMax}) über $[t_0,\overline{\tau}+\varepsilon]$ löst.
Die Fortsetzung $y_z(\cdot)$ ist für $t \in [\overline{\tau},\overline{\tau}+\varepsilon]$ die Lösung der Gleichung
$$y_z(t) = \xi_z(t)+ \int_{\overline{\tau}}^t \varphi\big(t,s,y_z(s)\big)\big) \, ds \quad\mbox{mit}\quad
\xi_z(t)= z  + \int_{t_0}^{\overline{\tau}} \varphi\big(t,s,x_z(s)\big) \, ds.$$
Es ist nun zu zeigen, dass zu jedem $\varrho>0$ eine Kugel $U_{\delta'}(x_0) \subseteq U_\delta(x_0)$ und ein $\varepsilon > 0$ so angeben lassen,
dass für alle $z \in U_{\delta'}(x_0)$ sich die Lösungen $x_z(\cdot)$ durch $y_z(\cdot)$ auf $t \in [\overline{\tau},\overline{\tau}+\varepsilon]$
fortsetzen lassen und die Ungleichung
$$\max_{t \in [\overline{\tau},\overline{\tau}+\varepsilon]} \|y_z(t)-y(t) \| \leq \varrho$$
erfüllt ist.
Denn dann konvergiert $x_z(\cdot)$ zusammen mit der zugehörigen Fortsetzung $y_z(\cdot)$ gleichmäßig über $[t_0,\overline{\tau}+\varepsilon]$ gegen
die Lösung $x(\cdot)$ und deren Fortsetzung $y(\cdot)$.

\newpage
\begin{figure}[h]
	\centering
	\fbox{\includegraphics[height=5cm]{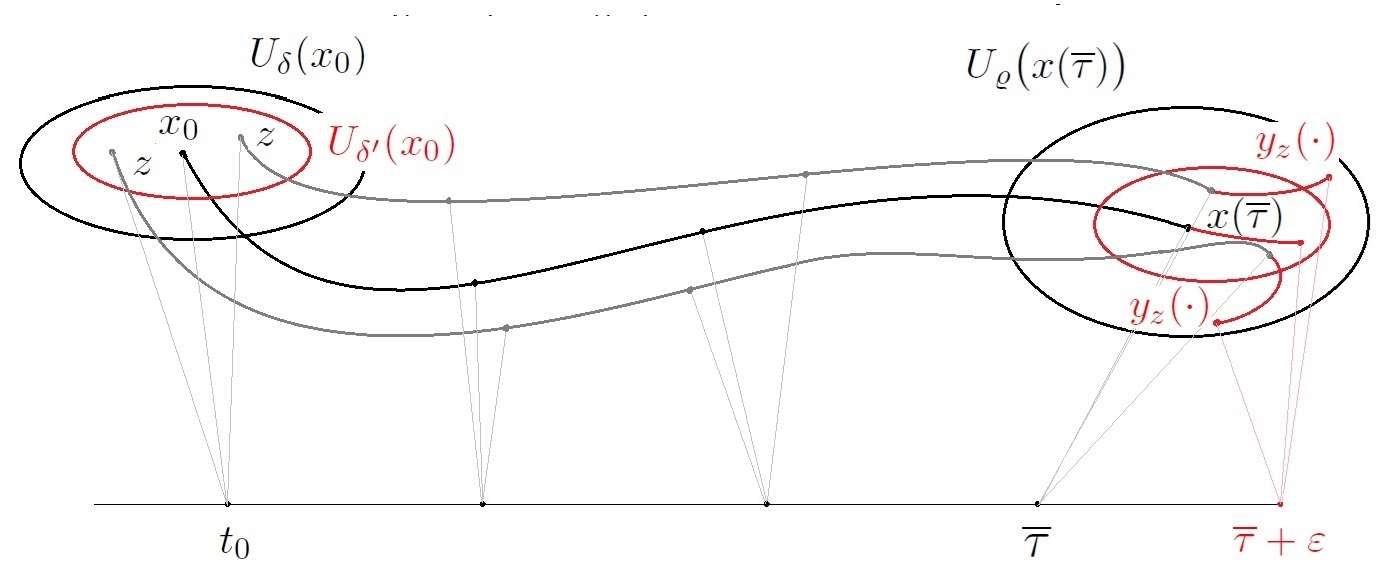}}
	\caption[Fortsetzungen einer nicht-maximalen Lösung]{Die Fortsetzungen $y_z(\cdot)$ von $x_z(\cdot)$.}
\end{figure}

Nach Satz \ref{SatzFixpunktlokal} existieren $\varepsilon >0$ und $\varrho >0$ derart,
dass die Abbildung $A$ für jedes stetige $\xi(\cdot)$ mit $\|\xi(t)-x(\overline{\tau})\| < \varrho$ für alle $t \in [\overline{\tau},\overline{\tau}+\varepsilon]$
genau einen Fixpunkt $x_\xi(\cdot)$ über $[\overline{\tau},\overline{\tau}+\varepsilon]$ besitzt.
Weiterhin wählen wir $\varepsilon$ mit $\varepsilon L <1$ und den Radius $\delta'$ der Kugel $U_{\delta'}(x_0)$ mit
$$\max_{t \in [t_0,\overline{\tau}]} \|x_z(t)-x(t) \| < \frac{1-\varepsilon L}{1 + 2L ( \overline{\tau}-t_0 ) } \varrho \quad
  \mbox{ für alle } z \in U_{\delta'}(x_0).$$
Die positiven Zahlen $\varepsilon, \, \varrho,\, \delta'$ seien im Weiteren so vorgegeben.
Ferner bezeichnen $\xi(\cdot)$ bzw. $\xi_z(\cdot)$ die Funktion,
mit deren Hilfe die Fortsetzungen für $t > \overline{\tau}$ gebildet werden:
$$\xi(t)= x(t_0)  + \int_{t_0}^{\overline{\tau}} \varphi\big(t,s,x(s)\big) \, ds,
  \quad \xi_z(t)= z  + \int_{t_0}^{\overline{\tau}} \varphi\big(t,s,x_z(s)\big) \, ds.$$
Für $z \in U_{\delta'}(x_0)$ ergeben sich mit den Darstellungen von $\xi(\cdot)$ und $\xi_z(\cdot)$ in der Form (\ref{DGLUmformung}):
\begin{eqnarray*}
\|\xi_z(t)-\xi(t)\|
&=& \bigg\|x_z(\overline{\tau}) - x(\overline{\tau})  
       + \int_{t_0}^{\overline{\tau}} \big[ \varphi\big(t,s,x_z(s)\big) - \varphi\big(t,s,x(s)\big) \big]\, ds \\
&&  + \int_{t_0}^{\overline{\tau}} \big[ \varphi\big(\overline{\tau},s,x_z(s)\big) - \varphi\big(\overline{\tau},s,x(s)\big) \big]\, ds \bigg\|\\
&\leq& \|x_z(\overline{\tau}) - x(\overline{\tau}) \| + 
       2L(\overline{\tau}-t_0) \cdot \max_{t \in [t_0,\overline{\tau}]} \|x_z(t)-x(t)\|.
\end{eqnarray*}
Mit Hilfe von Satz \ref{SatzFixpunktlokal} ergibt sich für die Fortsetzungen $y_z(\cdot), \, y(\cdot)$ nun der Abstand
\begin{eqnarray*}
\max_{t \in [\overline{\tau},\overline{\tau}+\varepsilon]} \|y_z(t)-y(t)\|
&\leq& \frac{1}{1-\varepsilon L} \cdot \max_{t \in [\overline{\tau},\overline{\tau}+\varepsilon]} \|\xi_z(t)-\xi(t)\| \\
&\leq& \frac{1 + 2L(\overline{\tau}-t_0)}{1-\varepsilon L} \cdot \max_{t \in [t_0,\overline{\tau}]} \|x_z(t)-x(t)\|
       < \varrho.
\end{eqnarray*}
Dies zeigt, dass sich zu jedem $\varrho >0$ positive Zahlen $\varepsilon$ und $\delta'$ so angeben lassen,
dass über $[t_0,\overline{\tau}+\varepsilon]$ die Ungleichung $\|y_z(t)-y(t)\| < \varrho$ für alle $z \in U_{\delta'}(x_0)$
erfüllt ist;
also die fortgesetzten Lösungen $x_z(\cdot)$ gleichmäßig gegen $x(\cdot)$ über $[t_0,\overline{\tau}+\varepsilon]$ für $z \to x_0$ konvergieren. 
Daher muss $\overline{\tau} > t_1$ gelten. 
Der Satz \ref{SatzEEglobal} ist damit bewiesen. \hfill $\blacksquare$

\newpage
\begin{satz}[Differenzierbarkeitssatz] \label{SatzDGLDifferenzierbarkeit}
Es seien $x_0(\cdot)$ und $x_z(\cdot)$ die Lösungen von (\ref{AnhangDGLMax}) zu den Anfangswerten $x_0$ bzw. $z \in U_\delta(x_0)$ in Satz \ref{SatzEEglobal}.
Weiter bezeichne $\Phi$ die Abbildung, die jedem Anfangswert $z \in U_\delta(x_0)$ die Lösung $x_z(\cdot)$ von (\ref{AnhangDGLMax}) zuordnet. \\
Dann ist $\Phi:U_\delta(x_0) \to C([t_0,t_1],\R^n)$ mit $\Phi(z)=x_z(\cdot)$ in $x_0$ stetig differenzierbar,
und für jedes $y \in \R^n$ genügt die Ableitung
$$\lim_{\lambda \to 0} \frac{\Phi(x_0+\lambda y)-\Phi(x_0)}{\lambda} = \Phi'(x_0)y=\xi_y(\cdot) \in C([t_0,t_1],\R^n)$$
der linearen Integralgleichung
$$\xi_y(t)= y + \int_{t_0}^t \varphi_x\big(t,s,x_0(s)\big) \xi_y(s) \, ds, \quad t \in [t_0,t_1].$$ 
\end{satz}

{\bf Beweis} Wir betrachten die Abbildung $F:\R^n \times C([t_0,t_1],\R^n) \to C([t_0,t_1],\R^n)$, die durch
$$\big[F\big(z,x(\cdot)\big)\big](t)= x(t)- z- \int_{t_0}^t \varphi\big(t,s,x(s)\big) \, ds = \big[G\big(x(\cdot)\big)\big](t) - z$$
definiert ist.
Wie im Beispiel \ref{DiffDynamik1} zeigt man,
dass die Abbildung $F$ in einer Umgebung des Punktes $x_0(\cdot)$ nach $x(\cdot)$ stetig differenzierbar ist und die Ableitung
$$\big[F_{x(\cdot)}\big(z,x(\cdot)\big)y(\cdot)\big](t)= y(t)- \int_{t_0}^t \varphi_x\big(t,s,x(s)\big) y(s) \, ds$$
besitrzt.
Es gilt $F\big(x_0,x_0(\cdot)\big) = 0$,
da $x_0(\cdot)$ Lösung (\ref{AnhangDGLMax}) zum Anfangswert $x_0$ ist.
Ferner stellt nach Lemma \ref{LemmaDGL1} der Operator $F_{x(\cdot)}\big(x_0,x_0(\cdot)\big)$ eine lineare homöomorphe Abbildung des Raumes $C([t_0,t_1],\R^n)$ in sich dar.
Gemä"s Satz \ref{SatzImpliziteFunktionen} über implizite Funktionen wird in einer hinreichend kleinen Umgebung $U_\delta(x_0)$ des Punktes $x_0$ eine
stetig differenzierbare Abbildung $z \to x_z(\cdot)$ in den Raum $C([t_0,t_1],\R^n)$ definiert,
für die $F\big(z,x_z(\cdot)\big)=0$ gilt.
Diese Abbildung ist in $x_0$ stetig differenzierbar. 
Die Bedingung $F\big(z,x_z(\cdot)\big)=0$ bedeutet
$$x_z(t)= z + \int_{t_0}^t \varphi\big(t,s,x_z(s)\big) \, ds,$$
d.\,h., dass diese Abbildung jedem Anfangswert $z \in U_\delta(x_0)$ die Lösung $x_z(\cdot)$ von (\ref{AnhangDGLMax}) zuordnet.
Demnach wird dadurch die Abbildung $\Phi$ beschrieben und $\Phi$ ist in $x_0$ stetig differenzierbar.
Ihre Ableitung ordnet nach Satz \ref{SatzImpliziteFunktionen} jedem $y \in \R^n$ die Vektorfunktion
$$\Phi'(x_0)y=\xi_y(\cdot) = - F^{-1}_{x(\cdot)}\big(x_0,x_0(\cdot)\big) \circ \big[ F_z\big(x_0,x_0(\cdot)\big)y \big] = -F^{-1}_{x(\cdot)}\big(x_0,x_0(\cdot)\big) y(\cdot)$$
zu, wobei $F_z\big(x_0,x_0(\cdot)\big)y$ die Funktion $y(\cdot)$ mit $y(t) \equiv y$ liefert.
Weiter lässt sich die Beziehung $\xi_y(\cdot) = F^{-1}_{x(\cdot)}\big(x_0,x_0(\cdot)\big) y(\cdot)$
in die Form
$$F_{x(\cdot)}\big(x_0,x_0(\cdot)\big)\xi_y(\cdot) =  y(\cdot) \quad\Leftrightarrow\quad
  \xi_y(t) - \int_{t_0}^t \varphi_x\big(t,s,x_0(s)\big) \xi_y(s) \, ds = y$$
bringen. Der Satz ist damit bewiesen. \hfill $\blacksquare$

%% file: IV-KonvexeAnalysis.tex
\section{Elemente der Konvexen Analysis} \label{AbschnittKonvexeAnalysis}
Bei der Zusammenstellung der grundlegenden Ergebnisse beschr"anken wir uns auf die Eigenschaften konvexer und lokalkonvexer Funktionen
nach Clarke \cite{Clarke}, Ioffe \& Tichomirov \cite{Ioffe} und Rockafellar \cite{Rockafellar}.
Im vorliegenden Rahmen stimmen die klassische Richtungsableitung und der Clarkesche Gradient "uberein.
Deswegen verweisen wir bez"uglich Lemma \ref{LemmaRichtungsableitung} und bez"uglich der Kettenregel \ref{SatzKettenregel}
auf \cite{Ioffe}.

\subsection{Das Subdifferential konvexer Funktionen}
Es seien $X,Y$ Banachr"aume.
Eine Funktion $f$ auf $X$ ist in der Konvexen Analysis eine Abbildung in die erweiterte reelle Zahlengerade, d.\,h.
$f:X \to \overline{\R} = [-\infty,\infty]$.
Der effektive Definitionsbereich
der Abbildung $f$ ist die Menge ${\rm dom\,}f = \{ x \in X | f(x) < \infty\}$. 
Die Funktion $f$ hei"st eigentlich\index{Funktion, absolutstetige!eigentliche@--, eigentliche},
falls ${\rm dom\,}f \not= \emptyset$ und $f(x) > -\infty$ f"ur alle $x \in X$ gelten. \\
Die eigentliche Funktion $f$ hei"st konvex\index{Funktion, absolutstetige!konvexe@--, konvexe},
wenn f"ur alle $x_1,x_2 \in X$ und alle $0 \leq \alpha \leq 1$ folgende Ungleichung gilt:
$$f\big(\alpha x_1 + (1-\alpha) x_2 \big) \leq \alpha f(x_1) + (1-\alpha) f(x_2).$$
Die Funktion $f$ hei"st homogen\index{Funktion, absolutstetige!homogene@--, homogene},
falls $f(0) = 0$ und $f(\lambda x) = \lambda f(x)$ f"ur alle $x \in X,\, \lambda > 0$ ist.
Eine eigentliche konvexe Funktion ist genau dann in einem Punkt stetig,
wenn sie auf einer Umgebung dieses Punktes nach oben beschr"ankt ist.
In diesem Fall ist das Innere des effektiven Definitionsbereichs nichtleer.
Ist andererseits eine homogene Funktion auf einer Umgebung des Nullpunktes stetig, so ist sie auf $X$ stetig. \\
Ist $f$ eine eigentliche konvexe Funktion auf $X$,
dann existiert in jedem Punkt der Menge ${\rm dom\,}f$ die klassische
Richtungsableitung\index{Ableitung, Fr\'echet-!Richtungs@--, Richtungs-}\label{Richtungsableitung},
d.\,h. f"ur alle $z \in X$ der Grenzwert
$$f'(x;z)=\lim_{\lambda \to 0^+} \frac{f(x + \lambda z) - f(x)}{\lambda}.$$
Sei $f$ eigentlich, konvex und in $x$ stetig.
Dann ist $f$ auf einer Umgebung des Punktes $x$ nach oben beschr"ankt,
in $x$ lokal Lipschitz-stetig und es existiert der Clarkesche Gradient\index{Clarkescher Gradient}\label{ClarkescherGradient}
$$f^\circ(x;z) = \limsup_{\substack{y \to x \\ \lambda \to 0^+}} \frac{f(y + \lambda z) - f(y)}{\lambda}.$$
Unter diesen Voraussetzungen ist au"serdem die Funktion $f$ ist im Punkt $x$ regul"ar
\index{Funktion, absolutstetige!reguläre@--, im Sinn der Konvexen Analysis reguläre}
im Sinn der Konvexen Analysis,
d.\,h. $f'(x;\cdot)=f^\circ(x;\cdot)$.
Das Subdifferential\index{Subdifferential} der eigentlichen konvexen Funktion $f$ besteht im Punkt $x$ aus allen Subgradienten
$x^* \in X^*$, d.\,h. \label{Subdifferential}
$$\partial f(x) = \{ x^* \in X^* | f(z) - f(x) \geq \langle x^*, z-x \rangle \mbox{ f"ur alle } z \in X \}.$$
F"ur eine eigentliche konvexe Funktion $f$ gilt $\partial f(x)=\partial f'(x;0)$ f"ur alle $x \in {\rm dom\,}f$.
Ist $f$ eine eigentliche homogene konvexe Funktion und $x \not= 0$, dann ist
$$\partial f(x) = \{ x^* \in \partial f(0) | f(x) = \langle x^*, x \rangle \}.$$
Mit $\partial_x f(x,y)$ bezeichnen wir das Subdifferential der Abbildung $x \to f(x,y)$.\label{partSubdifferential}

\subsection{Lokalkonvexe Funktionen}
Es seien $X,Y$ Banachr"aume.
Eine auf $X$ definierte Funktion $G$ hei"st im Punkt $x_0$ lokalkonvex\index{Funktion, absolutstetige!lokalkonvexe@--, lokalkonvexe},
wenn ihre Richtungsableitung in diesem Punkt existiert und $x \to G'(x_0;x)$ konvex ist.
Im Folgenden seien $g: X \to Y$ im Punkt $x_0 \in X$ Fr\'echet-differenzierbar und
$f: Y \to \overline{\R}$ eigentlich, konvex und im Punkt $g(x_0)$ stetig.

\begin{lemma} \label{LemmaRichtungsableitung}
Die Funktion $G: X \to \overline{\R}$, $G(x) = f\big( g(x)\big)$, besitzt in $x_0$ eine klassische Richtungsableitung, es gilt
$$G'(x_0;x) = f'\big(g(x_0);g'(x_0)x\big)$$
und die Richtungsableitung konvergiert bez"uglich jeder Richtung $x$
gleichm"a"sig\index{Funktion, absolutstetige!gleichmäßig@--, gleichmäßig differenzierbare}:
$$\bigg| \frac{G(x_0+\lambda z) - G(x_0)}{\lambda} - G'(x_0;x) \bigg| < \varepsilon \quad
  \mbox{f"ur alle } z \in U(x), \; 0<\lambda<\lambda_0.$$
Insbesondere folgt aus der gleichm"a"sigen Differenzierbarkeit bez"uglich aller Richtungen,
dass die Richtungsableitung der Abbildung $G$ eine stetige Funktion ist.
\end{lemma}

\begin{lemma} \label{LemmaClarke}
Die Funktion $G(x) = f\big(g(x)\big)$ ist in $x_0$ regul"ar.
\end{lemma}
{\bf Beweis} Nach Definition des $\limsup$ existieren Folgen $x_n \to x_0$ und $\lambda_n \to 0^+$ mit
$$\lim_{n \to \infty} \frac{G(x_n + \lambda_n x) - G(x_n)}{\lambda_n} = G^\circ(x_0;x).$$
Unter den getroffenen Voraussetzungen gilt $f'\big(g(x_0);g'(x_0)x\big) = f^\circ\big(g(x_0);g'(x_0)x\big)$ und wir erhalten 
\begin{eqnarray*}
G^\circ(x_0;x) &=& \lim_{n \to \infty} \frac{G(x_n + \lambda_n x) - G(x_n)}{\lambda_n}
                   = \lim_{n \to \infty} \frac{f\big(g(x_n + \lambda_n x)\big) - f\big(g(x_n)\big)}{\lambda_n} \\
               &\leq& \limsup_{\substack{y \to x \\ \lambda \to 0^+}} \frac{f\big(g(y + \lambda x)\big) - f\big(g(y)\big)}{\lambda} 
                   = f'\big(g(x_0);g'(x_0)x\big) = G'(x_0;x).
\end{eqnarray*}
Andererseits folgt unmittelbar die Relation
\begin{eqnarray*}
G'(x_0;x) &=& \lim_{\lambda \to 0^+} \frac{G(x_0 + \lambda x) - G(x_0)}{\lambda}
            \leq \limsup_{\substack{y \to x \\ \lambda \to 0^+}} \frac{G(y + \lambda x) - G(y)}{\lambda} = G^\circ(x_0;x).
\end{eqnarray*}
Beide Ungleichungen zeigen $G'(x_0;\cdot)=G^\circ(x_0;\cdot)$. \hfill $\blacksquare$

\begin{satz}[Kettenregel] \label{SatzKettenregel}
Es seien $g: X \to Y$ im Punkt $x_0 \in X$ Fr\'echet-differenzierbar und $f: Y \to \overline{\R}$ eigentlich,
konvex und im Punkt $g(x_0)$ stetig.
Dann ist die Funktion $G(x) = f\big( g(x)\big)$ im Punkt $x_0$ regul"ar und es gilt
$$\partial G(x_0) = g'^*(x_0) \partial f\big(g(x_0)\big).$$
\end{satz}


\subsection{Das Subdifferential konkreter Funktionen}
Im Folgenden bestimmen wir die Subdifferentiale der Abbildungen
$$f\big(x(\cdot)\big) = \max_{t \in I} x(t), \qquad G\big(x(\cdot)\big) = \max_{t \in I} g\big(t,x(t)\big)$$
im Raum stetiger Funktionen auf den Intervallen $I=[t_0,t_1]$, $I=\R_+$ bzw. $I = \overline{\R}_+$. \\
Insbesondere f"ur das unbeschr"ankte Intervall zeigt sich,
dass erst im Fall des Abschlusses $I = \overline{\R}_+$ das Subdifferential der Funktion $f$ eine passende Darstellung besitzt.
Demgegen"uber weist das Subdifferential von $f$ "uber $I=\R_+$ eine gewisse Unvollst"andigkeit auf.
Der Grund daf"ur ist der Rieszsche Darstellungssatz \ref{SatzRieszC0} f"ur die stetigen linearen Funktionale im Raum $C_0(\R_+, \R^n)$,
der im Gegensatz zum Darstellungssatz \ref{FolgerungRieszClim} f"ur $C_{\lim}(\R_+, \R^n)$ nicht das in $t=\infty$ konzentrierte
signierte Borelsche Ma"s enth"alt. 

\begin{beispiel} \label{SubdifferentialMaximum1}
{\rm Es sei $[t_0,t_1] \subset \R$ ein kompaktes Intervall.
Im Raum $C([t_0,t_1],\R)$ betrachten wir die Funktion
$$f\big(x(\cdot)\big) = \max_{t \in [t_0,t_1]} x(t).$$
Diese Funktion ist stetig, konvex und homogen.
Das Subdifferential $\partial f(0)$ besteht nach Definition aus denjenigen signierten regul"aren Borelschen Ma"sen $\mu$ auf $[t_0,t_1]$,
die der Bedingung
$$\max_{t \in [t_0,t_1]} x(t) \geq \int_{t_0}^{t_1} x(t) \, d\mu(t) \qquad \mbox{f"ur alle } x(\cdot) \in C([t_0,t_1],\R)$$
gen"ugen.
Hieraus folgt, dass das Ma"s $\mu$ nichtnegativ ist und $\|\mu\| = 1$ gilt.
Man erh"alt n"amlich f"ur $-x(t)$ nach obiger Ungleichung
$$\int_{t_0}^{t_1} x(t) \, d\mu(t) \geq - \max_{t \in [t_0,t_1]} \big( -x(t)\big) = \min_{t \in [t_0,t_1]} x(t).$$
Ist daher $x(t) \geq 0$ f"ur alle $t \in [t_0,t_1]$,
so ist auch $\displaystyle\int_{t_0}^{t_1} x(t) \, d\mu(t) \geq 0$ und $\mu$ ist positiv.
Weiterhin besitzen alle $\mu \in \partial f(0)$ die Norm $\|\mu\| = 1$,
denn nur dann gilt die Relation
$$\max_{t \in [t_0,t_1]} x(t) \geq \int_{t_0}^{t_1} x(t) \,d\mu(t) \geq \min_{t \in [t_0,t_1]} x(t)$$
f"ur konstante Funktionen $x(\cdot)$. \\
Offenbar ist auch umgekehrt, wenn $\mu \geq 0$ ist und Norm Eins hat, die Ungleichung
$$\max_{t \in [t_0,t_1]} x(t) \geq \int_{t_0}^{t_1} x(t) \, d\mu(t)$$
richtig.
Somit haben wir das Subdifferential der Funktion $f$ im Nullpunkt berechnet. \\[2mm]
F"ur $x(\cdot) \not= 0$ gilt f"ur die eigentliche homogene konvexe Funktion $f\big(x(\cdot)\big)$:
$$\partial f\big(x(\cdot)\big) = \big\{ x^* \in \partial f(0) \big| \langle x^*,x(\cdot) \rangle = f\big(x(\cdot)\big) \big\},$$
wobei $\partial f(0)$ die regul"aren Borelschen Ma"se mit $\|\mu\|=1$ enth"alt.
Wir zeigen, dass die Ma"se $\mu \in \partial f\big(x(\cdot)\big)$ auf den Mengen
$$T = \bigg\{ \tau \in [t_0,t_1] \,\bigg|\, x(\tau) = \max_{t \in [t_0,t_1]} x(t) \bigg\}$$
konzentriert sind:
Der K"urze halber bezeichne $M$ das Maximum von $x(t)$ auf $[t_0,t_1]$.
Die Menge $[t_0,t_1] \setminus T$ ist offen, da $T$ abgeschlossen ist.
Angenommen, es existiert eine me"sbare Menge $B \subseteq [t_0,t_1]$ mit $B \cap T = \emptyset$ und $\mu(B) > 0.$
Dann gilt
$$\int_B x(t)\, d\mu(t) < \int_B M\, d\mu(t) \Rightarrow \int_{[t_0,t_1] \setminus T} x(t)\, d\mu(t) < \int_{[t_0,t_1] \setminus T} M\, d\mu(t).$$
Damit folgt der Widerspruch
$$\max_{t\in [t_0,t_1]} x(t) = \int_{t_0}^{t_1} M\, d\mu(t) > \int_T M\, d\mu(t) + \int_{[t_0,t_1] \setminus T} x(t) \, d\mu(t)
   = \int_{t_0}^{t_1} x(t)\, d\mu(t).$$
Ist umgekehrt $\mu \geq 0$ und besitzt Norm Eins, dann geh"ort $\mu$ der Menge $\partial f(0)$ an.
Wenn $\mu$ zus"atzlich auf der Menge $T$ konzentriert ist, dann ist auch
$$\max_{t \in [t_0,t_1]} x(t) = \int_{t_0}^{t_1} x(t)\, d\mu(t)$$
erf"ullt.
Damit besteht das Subdifferential der Funktion $f$ in einem vom Nullpunkt verschiedenen Punkt $x(\cdot)$ aus den regul"aren
Borelschen Ma"sen auf $[t_0,t_1]$, deren Norm gleich Eins ist und die auf der Menge $T$ konzentriert sind. \hfill $\square$}
\end{beispiel}

\begin{beispiel} \label{SubdifferentialMaximum2}
{\rm Es sei $x_*(\cdot) \in C([t_0,t_1],\R)$ und es sei $g(t,x) : \R \times \R^n \to \R$
auf der Menge $V_\gamma = \{ (t,x) \in [t_0,t_1] \times \R^n \,|\, \|x-x_*(t)\| \leq \gamma\}$
stetig und bez"uglich $x$ stetig differenzierbar. \\
Unter diesen Voraussetzungen ist nach Beispiel \ref{DiffAbbildung} die Abbildung
$$\tilde{g}: C([t_0,t_1],\R^n) \to C([t_0,t_1],\R), \qquad \big[\tilde{g}\big(x(\cdot)\big)\big] (t) = g\big(t,x(t)\big),\quad t \in [t_0,t_1],$$
in $x_*(\cdot)$ Fr\'echet-differenzierbar und es gilt
$$\big[\tilde{g}'\big(x_*(\cdot)\big) x(\cdot)\big] (t) = \big\langle g_x \big(t,x_*(t)\big), x(t) \big\rangle, \qquad t \in [t_0,t_1].$$
Weiterhin ist die Funktion $f\big(x(\cdot)\big)$ im Beispiel \ref{SubdifferentialMaximum1}
auf $C([t_0,t_1],\R)$ stetig. \\[2mm]
Wir bestimmen nun in $x_*(\cdot)$ das Subdifferential der Funktion
$$G\big(x(\cdot)\big) = f\Big(g\big(x(\cdot)\big)\Big) =\max_{t \in [t_0,t_1]} g\big(t,x(t)\big).$$
Dazu k"onnen wir die Kettenregel (Satz \ref{SatzKettenregel}) anwenden:
$$\partial G\big(x_*(\cdot)\big) = \tilde{g}'^*\big(x_*(\cdot)\big) \partial f \big(\tilde{g}\big(x_*(\cdot)\big)\big).$$
F"ur $x^* \in \tilde{g}'^*\big(x_*(\cdot)\big) \partial f \big(\tilde{g}\big(x_*(\cdot)\big)\big)$ und
$x(\cdot) \in C([t_0,t_1],\R^n)$ gilt dann:
$$\big\langle x^*,x(\cdot) \big\rangle
   = \big\langle \tilde{g}'^*\big(x_*(\cdot)\big)\, \mu, x(\cdot) \big\rangle
   = \big\langle \mu , \tilde{g}'\big(x_*(\cdot)\big)x(\cdot) \big\rangle
   = \int_{t_0}^{t_1} \big\langle g_x \big(t,x_*(t)\big), x(t) \big\rangle \,d\mu(t).$$
Daher besteht das Subdifferential der Funktion $G$ im Punkt $x_*(\cdot)$ genau aus denjenigen stetigen linearen
Funktionalen $x^*$,
die die Darstellung
$$\big\langle x^*,x(\cdot) \big\rangle = \int_{t_0}^{t_1} \big\langle g_x \big(t,x_*(t)\big), x(t) \big\rangle \,d\mu(t)$$
besitzen,
wobei das regul"are Borelsches Ma"s $\mu$ auf
$T= \big\{ t \in [t_0,t_1] \,\big|\, g\big(t,x_*(t)\big) = G\big(x_*(\cdot)\big) \big\}$
konzentriert ist und $\|\mu\| =1$ gilt. \hfill $\square$}
\end{beispiel}

\begin{beispiel} \label{SubdifferentialMaximum3}
{\rm Wir betrachten im Raum $C_{\lim}(\R_+,\R)$ die Funktion
$$f\big(x(\cdot)\big) = \max_{t \in \overline{\R}_+} x(t).$$
Mit den gleichen Argumenten wie in Beispiel \ref{SubdifferentialMaximum1}
ist diese Funktion stetig, konvex und homogen.
Au"serdem besitzen nach dem Rieszschen Darstellungssatz in Folgerung \ref{FolgerungRieszClim} die stetigen linearen 
Funktional $x^*$ im Raum $C_{\lim}(\R_+,\R)$ die Darstellung
$$\langle x^*(\cdot),x(\cdot) \rangle = \int_0^\infty x(t) \, d\mu(t)=
  \int_0^\infty x(t) \, d\mu_0(t) + x(\infty) \mu(\{\infty\}),$$
wobei $\mu$ ein signiertes reguläres Borelsches Maß auf $\overline{\R}_+$ ist,
das nach Anhang \ref{AnhangMassR}, Definition \ref{DefinitionMassR} die eindeutige Zerlegung
$\mu=\mu_0+\mu(\{\infty\})$ mit einem signierten regul"aren Borelschen Ma"s $\mu_0$ auf $\R_+$
und einem signierten, in $t=\infty$ konzentrierten Borelschen Ma"s $\mu(\{\infty\})$ besitzt. \\
Somit besteht das Subdifferential $\partial f(0)$ aus denjenigen Borelschen Ma"sen $\mu$ auf $\overline{\R}_+$,
die nach dem Rieszschen Darstellungssatz die Zerlegung $\mu=\mu_0+\mu(\{\infty\})$ wie eben beschrieben besitzen,
die nichtnegativ sind und die die Totalvariation $\|\mu\| = 1$ besitzen. \\
Weiterhin besteht das Subdifferential der Funktion $f$ in einem vom Nullpunkt verschiedenen Punkt $x(\cdot)$ aus den Borelschen Ma"sen
$\mu=\mu_0+\mu_\infty$ auf $\overline{\R}_+$,
deren Norm gleich Eins ist und die auf der Menge
$T = \big\{ t \in \overline{\R}_+ \,\big|\, x(t) = f\big(x(\cdot)\big) \big\}$
konzentriert sind. \hfill $\square$}
\end{beispiel}

\begin{beispiel} \label{SubdifferentialMaximum4}
{\rm Es sei $x_*(\cdot) \in C_{\lim}(\R_+,\R)$ und es sei $g(t,x) : \R \times \R^n \to \R$
auf der Menge $V_\gamma = \{ (t,x) \in \overline{\R}_+ \times \R^n \,|\, \|x-x_*(t)\| \leq \gamma\}$
gleichm"a"sig stetig und bez"uglich $x$ gleichm"a"sig stetig differenzierbar.
Unter diesen Voraussetzungen ist nach Beispiel \ref{DiffAbbildung} die Abbildung
$$\tilde{g}: C_{\lim}(\R_+,\R^n) \to C_{\lim}(\R_+,\R), \qquad \big[\tilde{g}\big(x(\cdot)\big)\big] (t) = g\big(t,x(t)\big),\quad t \in \R_+,$$
in $x_*(\cdot)$ Fr\'echet-differenzierbar und es gilt
$$\big[\tilde{g}'\big(x_*(\cdot)\big) x(\cdot)\big] (t) = \big\langle g_x \big(t,x_*(t)\big), x(t) \big\rangle, \qquad t \in \R_+.$$
Weiterhin ist die Funktion $f$ im Beispiel \ref{SubdifferentialMaximum3} auf $C_{\lim}(\R_+,\R)$ stetig.
Auf die Funktion
$$G\big(x(\cdot)\big) = f\Big(g\big(x(\cdot)\big)\Big) =\max_{t \in \overline{\R}_+} g\big(t,x(t)\big).$$
wenden wir die Kettenregel an:
Wie im Beispiel \ref{SubdifferentialMaximum2} besteht 
das Subdifferential der Funktion $G$ im Punkt $x_*(\cdot)$ genau aus denjenigen stetigen linearen Funktionalen $x^*$,
die mit einem regulären Borelsches Maß $\mu=\mu_0+\mu(\{\infty\})$ auf $\overline{\R}_+$ die Darstellung
$$\big\langle x^*,x(\cdot) \big\rangle = \int_0^\infty \big\langle g_x \big(t,x_*(t)\big), x(t) \big\rangle \,d\mu_0(t) + 
                              \lim_{t \to \infty} \big\langle g_x \big(t,x_*(t)\big), x(t) \big\rangle \, \mu(\{\infty\})$$
besitzen,
wobei das positive Borelsche Ma"s $\mu$ auf $T= \big\{ t \in \overline{\R}_+ \,\big|\, g\big(t,x_*(t)\big) = G\big(x_*(\cdot)\big) \big\}$
konzentriert ist und $\|\mu\| =1$ gilt. \hfill $\square$}
\end{beispiel}

\begin{beispiel} \label{SubdifferentialMaximum5}
{\rm Im Raum $C_0(\R_+,\R)$ betrachten wir die Funktion $f\big(x(\cdot)\big) = \sup\limits_{t \in \R_+} x(t)$.
Sie ist stetig, homogen und konvex.
F"ur das Subdifferential $\partial f(0)$ liefert die Ungleichung
$$\sup_{t \in \R_+} x(t) \geq \int_0^\infty x(t) \, d\mu(t) \quad \mbox{ f"ur alle }x(\cdot) \in C_0(\R_+,\R),$$
dass das "uber $\R_+$ signierte regul"are Borelsche Ma"s $\mu$ nichtnegativ ist.
Aber im Gegensatz zu Beispiel \ref{SubdifferentialMaximum3} erhalten wir aus der Ungleichungskette 
$$\sup_{t \in \R_+} x(t) \geq \int_0^\infty x(t) \, d\mu(t) \geq \inf_{t \in \R_+} x(t)$$
lediglich $\|\mu\|\leq 1$, da die Funktionen $x(\cdot)$ im Unendlichen stets verschwinden. \\
Umgekehrt gilt, wenn $\mu \geq 0$ und $\|\mu\|\leq 1$, die Ungleichung
$$\sup_{t \in \R_+} x(t) \geq \sup_{t \in \R_+} x(t) \cdot \int_0^\infty d\mu(t) = \int_0^\infty \sup_{t \in \R_+} x(t) \,d\mu(t)
  \geq \int_0^\infty x(t) \,d\mu(t)$$
f"ur alle $x(\cdot) \in C_0(\R_+,\R)$.
Damit besteht das Subdifferential der Funktion $f$ in $x(\cdot)=0$ aus allen positiven, regul"aren Borelschen Ma"sen $\mu$ auf $\R_+$ mit $\|\mu\| \leq 1$. \\
F"ur $x(\cdot) \not= 0$ m"ussen wir die Unterscheidung treffen,
ob die Funktion $x(\cdot) \in C_0(\R_+,\R)$ ein Maximum besitzt oder nicht.
Die Funktion $x(\cdot)$ besitzt genau dann kein Maximum "uber $\R_+$,
wenn $x(t)<0$ f"ur alle $t \in \R_+$ gilt.
In diesem Fall erhalten wir die Entartung
$$0 = \sup_{t\in \R_+} x(t) = \int_{\R_+} x(t)\, d\mu(t) \qquad\Leftrightarrow\qquad \|\mu\|=0.$$
Nimmt die Funktion $x(\cdot)$ ihr Maximum "uber $\R_+$ an,
dann gilt $\|\mu\|=1$ und $\mu$ ist auf der Menge $T=\big\{ t \in \R_+ \,\big|\, x(t) = f\big(x(\cdot)\big) \big\}$ konzentriert. \hfill $\square$}
\end{beispiel}


%% file: V-Nadelvariation.tex
\section{Mehrfache Nadelvariationen nach Ioffe \& Tichomirov} \label{AnhangNV}
Die einfache Nadelvariation\label{AbschnittNV}\index{Nadelvariation, mehrfache}, 
auf deren Basis das Maximumprinzip f"ur die Aufgabe mit freiem Endpunkt im Abschnitt \ref{AbschnittPMPeinfach} hergeleitet wurde,
erweist sich als ungeeignet in einem Steuerungsproblem mit beidseitigen Randbedingungen.
Deswegen stellen wir in diesem Abschnitt mehrfache Nadelvariationen vor. \\[2mm]
Wir folgen bei der Darstellung mehrfacher Nadelvariationen Ioffe \& Tichomirov \cite{Ioffe}.
Als erstes geben wir eine Konstruktion geeigneter Tr"agerfamilien an,
auf denen die Nadelvariationen durchgef"uhrt werden.
Danach definieren wir die mehrfachen Nadelvariationen und untersuchen abschlie"send die Eigenschaften,
die die mehrfache Nadelvariation dem Steuerungsproblem "ubergibt.

%% file: V-1-Konstruktion.tex
\subsection{Die Konstruktion der Tr\"agerfamilien} \label{AbschnittMengenfamilie}
\begin{lemma} \label{LemmaMengenfamilie1}
Zu jeder auf $[t_0,t_1]$ definierten Treppenfunktion $y(\cdot), y(t) \in \R^n,$ und jedem $\delta > 0$ kann man eine einparametrige
Familie
$$\big\{ M(\alpha) \big\}_{0 \leq \alpha \leq 1} = \big\{ M(\alpha, y(\cdot),\delta) \big\}_{0 \leq \alpha \leq 1}$$
messbarer Teilmengen des Intervalls $[t_0,t_1]$ derart angeben, dass die Beziehungen
\begin{eqnarray*}
&&|M(\alpha)| = \alpha (t_1-t_0), \qquad M(\alpha') \subseteq M(\alpha), \qquad 0 \leq \alpha' \leq \alpha \leq 1, \\
&& \max_{t \in [t_0,t_1]}
  \left\| \int_{t_0}^t \Big[ \big( \chi_{M(\alpha)}(\tau) - \chi_{M(\alpha')}(\tau) \big) y(\tau)
                     - (\alpha-\alpha') y(\tau) \Big] \, d\tau \right\| \leq \delta |\alpha-\alpha'|
\end{eqnarray*}
gelten.
\end{lemma}

Der Grundgedanke der nachstehend beschriebenen Konstruktion l"asst sich am leichtesten in dem Fall verstehen,
dass die Vektorfunktion $y(\cdot)$ konstant ist: $y(t) \equiv C$.
Zerlegt man das Intervall $[t_0,t_1]$ in gleiche Teile $\Delta_i$, deren L"ange kleiner als $\delta/C$ ist,
so kann man als Mengenfamilie $\{M(\alpha)\}$ die Vereinigung derjenigen Intervall $\Delta_i(\alpha)$ nehmen,
deren linkes Ende jeweils mit dem linken Ende des Intervalls $\Delta_i$ und deren L"ange jeweils mit dem $\alpha$-ten Teil
der L"ange von $\Delta_i$ "ubereinstimmt. \\
Für $y(t) \equiv C$ bezeichnen zur Illustration $Y(\cdot)$ bzw. $Y_\alpha(\cdot)$ die Funktionen
$$Y(t)=\int_{t_0}^t y(\tau)\, d\tau,\qquad Y_\alpha(t)=\int_{t_0}^t \chi_{M(\alpha)}(\tau) y(\tau)\, d\tau, \qquad t \in [t_0,t_1].$$
Es sei $\delta>0$ gewählt.
Dann kann die Zerlegung des Intervalls $[t_0,t_1]$ in Teilintervalle $\Delta_i$ bestimmt werden.
Der Abbildung \ref{AbbMengenfamilie1} entnimmt man das Aussehen der Funktion $Y_\alpha(\cdot)$ bei Wahl des Parameters $\alpha \in [0,1]$:
Bei Verkleinerung des Parameters $\alpha$ werden die einzelnen Stufenabsätze immer breiter und die Höhenunterschiede
zwischen nacheinander folgender Stufen immer kleiner.
\begin{center}
\includegraphics[width=12cm]{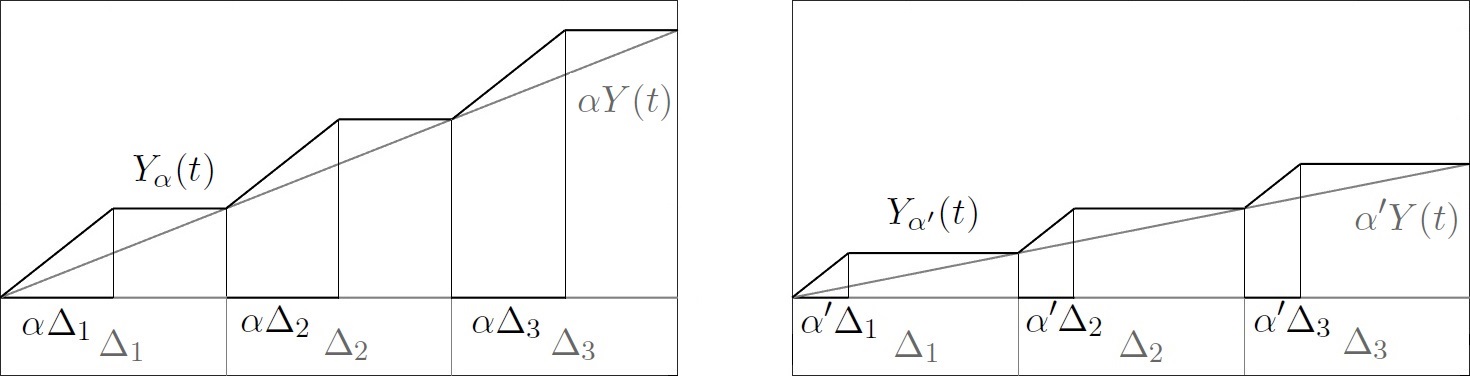}
\captionof{figure}[Mengenfamilien mehrfacher Nadelvariation -- Paramater $\alpha$]{Mengenfamilie $\{M(\alpha)\}$ bei Verkleinerung des Parameters $\alpha' < \alpha$.}
\label{AbbMengenfamilie1}
\end{center}

Wird andererseits ein $\delta'<\delta$ zu Beginn der Konstruktion der Mengenfamilien ausgewählt,
so zeigt Abbildung \ref{AbbMengenfamilie2} die entstehende feinere Approximation von $\alpha Y(\cdot)$ durch $Y_\alpha(\cdot)$. 
\begin{center}
\includegraphics[width=12cm]{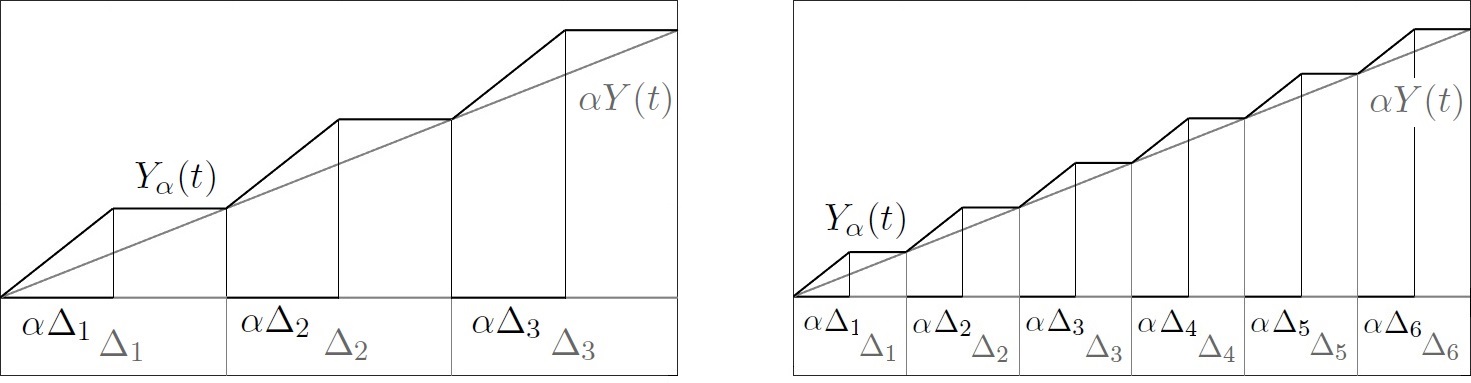}
\captionof{figure}[Mengenfamilien mehrfacher Nadelvariation -- Paramater $\delta$]{Mengenfamilie $\{M(\alpha)\}$ bei Verkleinerung der Feinheit $\delta' < \delta$.}
\label{AbbMengenfamilie2}
\end{center}

Wir erweitern diesen Gedankengang auf den Fall von Treppenfunktionen und führen den Beweis von Lemma \ref{LemmaMengenfamilie1}. \\[2mm]
{\bf Beweis} Da $y(\cdot)$ eine Treppenfunktion auf dem Intervall $[t_0,t_1]$ ist,
gibt es endlich viele $y_1,...,y_d \in \R^n$ und eine disjunkte Zerlegung $\{ A_1,...,A_d\}$ von $[t_0,t_1]$ mit
$$y(t) = \sum_{j=1}^d y_j \chi_{A_j}(t).$$
Ferner sei $C= \max\limits_j \|y_j\|$.
Wir zerlegen das Intervall $[t_0,t_1]$ in gleiche Teile $\Delta_1,...,\Delta_r$,
deren L"ange h"ochstens $\delta / (2C)$ ist und setzen
$$M_{ij}(\alpha) = \bigg\{ t \in (A_j \cap \Delta_i)
                      \,\bigg|\, \int_{t_0}^t \chi_{(A_j \cap \Delta_i)} (\tau) \, d\tau
                                                                < \alpha \cdot |A_j \cap \Delta_i| \bigg\}.$$
Dann ist $M(\alpha) = \displaystyle\bigcup_{i,j} M_{ij}(\alpha)$ die gesuchte Menge.
Denn es gilt
$$|M(\alpha)| = \sum_{i,j} |M_{ij}(\alpha)| = \alpha \sum_{i,j} |A_j \cap \Delta_i| = \alpha (t_1-t_0).$$
Ferner zieht die Ungleichung $\alpha \geq \alpha'$ die Beziehung $M_{ij}(\alpha') \subseteq M_{ij}(\alpha)$,
d.\,h. $M(\alpha') \subseteq M(\alpha)$, nach sich.
Schlie"slich gilt
$$\int_{M_{ij}(\alpha)} y(t) \, dt = \int_{t_0}^{t_1} \chi_{M_{ij}(\alpha)}(t) y(t) \, dt = \alpha y_j \cdot |A_j \cap \Delta_i|,$$
weil jeweils auf den Mengen $A_j$ die Vektorfunktion $y(\cdot)$ konstant und gleich $y_j$ ist.
Hieraus folgt in den Endpunkten der Intervalle $\Delta_i$:
\begin{eqnarray*}
\lefteqn{\int_{\Delta_i} \big( \chi_{M(\alpha)}(t) - \chi_{M(\alpha')}(t) \big) y(t) \, dt} \\
&=& \sum_j \left( \int_{M_{ij}(\alpha)} y(t) \, dt - \int_{M_{ij}(\alpha')} y(t) \, dt \right) 
    = (\alpha-\alpha') \sum_j y_j \cdot |A_j \cap \Delta_i| \\
&=& (\alpha-\alpha') \sum_j \int_{\Delta_i} \chi_{(A_j \cap \Delta_i)}(t) y(t) \, dt
    = (\alpha-\alpha') \int_{\Delta_i} y(t) \, dt,
\end{eqnarray*}
also, dass in den Endpunkten der Intervalle $\Delta_i$ die Werte der Integrale
$$\int_{t_0}^{t} \big( \chi_{M(\alpha)}(\tau) - \chi_{M(\alpha')}(\tau) \big) y(\tau) \, d\tau, \qquad 
  \int_{t_0}^{t} (\alpha-\alpha') y(\tau) \, d\tau$$
"ubereinstimmen.
F"ur $\Delta_i = [\tau_i,\tau_{i+1}]$ und $\tau_i < t < \tau_{i+1}$ ergibt sich
$$\left\| \int_{\tau_i}^t \Big[ \big( \chi_{M(\alpha)}(\tau) - \chi_{M(\alpha')}(\tau) \big) y(\tau)
                                 - (\alpha-\alpha') y(\tau) \Big] \, d\tau \right\|
  \leq 2 C |\alpha-\alpha'| \cdot |t-\tau_i| \leq \delta |\alpha-\alpha'|.$$
Damit ist Lemma \ref{LemmaMengenfamilie1} gezeigt. \hfill $\blacksquare$

\begin{lemma} \label{LemmaMengenfamilie2}
Zu jeder auf $[t_0,t_1]$ definierten beschr"ankten messbaren Vektorfunktion $y(\cdot)$ und jedem $\delta > 0$ gibt es eine
einparametrige Familie
$$\big\{ M(\alpha) \big\}_{0 \leq \alpha \leq 1} = \big\{ M(\alpha, y(\cdot),\delta) \big\}_{0 \leq \alpha \leq 1}$$
mit den Eigenschaften 
\begin{eqnarray*}
&&|M(\alpha)| = \alpha (t_1-t_0), \qquad M(\alpha') \subseteq M(\alpha), \qquad 0 \leq \alpha' \leq \alpha \leq 1, \\
&& \max_{t \in [t_0,t_1]}
  \left\| \int_{t_0}^t \Big[ \big( \chi_{M(\alpha)}(\tau) - \chi_{M(\alpha')}(\tau) \big) y(\tau)
                     - (\alpha-\alpha') y(\tau) \Big] \, d\tau \right\| \leq \delta |\alpha-\alpha'|.
\end{eqnarray*}
\end{lemma}

{\bf Beweis} Jede beschr"ankte messbare Vektorfunktion auf $[t_0,t_1]$ ist gleichm"a"siger Grenzwert von Treppenfunktionen.
D.\,h., es gibt eine Treppenfunktion $\tilde{y}(\cdot)$ mit
$$\sup_{t \in [t_0,t_1]} \|y(t)-\tilde{y}(t)\| \leq \frac{\delta}{t_1-t_0}.$$
Nach letztem Lemma kann man eine Familie $\big\{ M(\alpha, \tilde{y}(\cdot),\delta) \big\}_{0 \leq \alpha \leq 1}$ mit den entsprechenden
Eigenschaften f"ur $\tilde{y}(\cdot)$ angeben.
Ist dann $\alpha \geq \alpha'$, so gilt
$$\left\| \int_{t_0}^t \big( \chi_{M(\alpha)}(\tau) - \chi_{M(\alpha')}(\tau) \big) \cdot \big(y(\tau)-\tilde{y}(\tau)\big) d\tau \right\|
  \leq \frac{\delta}{t_1-t_0} |M(\alpha) \setminus M(\alpha')| \leq \delta |\alpha-\alpha'|$$
und
$$\left\| \int_{t_0}^t (\alpha-\alpha') \big(y(\tau) - \tilde{y}(\tau)\big) d\tau \right\| \leq 
  |\alpha-\alpha'| \int_{t_0}^t \big\|  \big(y(\tau)-\tilde{y}(\tau)\big) \big\| d\tau
  \leq \delta |\alpha-\alpha'|.$$
Zusammenfassend erhalten wir mit der Aussage des letzten Lemma f"ur Treppenfunktionen:
\begin{eqnarray*}
\lefteqn{\left\| \int_{t_0}^t \Big[ \big( \chi_{M(\alpha)}(\tau) - \chi_{M(\alpha')}(\tau) \big) y(\tau)
                     - (\alpha-\alpha') y(\tau) \Big] \, d\tau \right\|} \\
&\leq& \left\| \int_{t_0}^t \big( \chi_{M(\alpha)}(\tau) - \chi_{M(\alpha')}(\tau) \big) \cdot \big(y(\tau) - \tilde{y}(\tau)\big)\, d\tau
               \right\|
     + \left\| \int_{t_0}^t (\alpha-\alpha') \big(y(\tau) - \tilde{y}(\tau)\big) d\tau \right\|  \\
&&     + \left\| \int_{t_0}^t \Big[ \big( \chi_{M(\alpha)}(\tau) - \chi_{M(\alpha')}(\tau) \big) \tilde{y}(\tau)
                     - (\alpha-\alpha') \tilde{y}(\tau) \Big] \, d\tau \right\| \leq 3 \delta |\alpha-\alpha'|,
\end{eqnarray*}
also $M(\alpha) = M(\alpha, y(\cdot),3\delta)$. \hfill $\blacksquare$

\begin{lemma} \label{LemmaMengenfamilie3}
Es seien $y_i(\cdot),\, y_i : [t_0,t_1] \to \R^{n_i},\, i=1,...,d,$ beschr"ankte messbare Vektorfunktionen.
Dann gibt es zu jedem $\delta > 0$ einparametrige Mengenfamilien $M_1(\alpha),...,M_d(\alpha)$
messbarer Teilmengen des Intervalls $[t_0,t_1]$,
wobei der Parameter $\alpha$ Werte zwischen $0$ und $1/d$ annimmt derart, dass f"ur alle
$i=1,...,d,$ $0 \leq \alpha' \leq \alpha \leq 1/d$ und $i \not= i'$ gilt:
\begin{eqnarray}
&& \label{Nadelvariation1} 
   \hspace*{-15mm} |M_i(\alpha)| = \alpha (t_1-t_0), \qquad M_i(\alpha') \subseteq M_i(\alpha),
   \qquad M_i(\alpha) \cap M_{i'}(\alpha') = \emptyset,\\
&& \label{Nadelvariation2}
   \hspace*{-15mm} \max_{t \in [t_0,t_1]}
   \bigg\| \int_{t_0}^t \big( \chi_{M_i(\alpha)}(\tau) - \chi_{M_i(\alpha')}(\tau) \big) y_i(\tau) \, d\tau
   - (\alpha-\alpha') \int_{t_0}^t y_i(\tau) \, d\tau \bigg\| \leq \delta |\alpha-\alpha'|.
\end{eqnarray}
\end{lemma}

{\bf Beweis} Es sei $n=n_1 + ... + n_d$. Dann ist
$$t \to z(t) = \big( y_1(t),...,y_d(t) \big)$$
eine messbare beschr"ankte Abbildung des Intervalls $[t_0,t_1]$ in $\R^n$. \\
Wir w"ahlen eine Familie $\{M(\alpha)\}_{0 \leq \alpha \leq 1}$ von messbaren Teilmengen des Intervalls $[t_0,t_1]$,
die gemeinsam mit $z(\cdot)$ und $\delta$ den Bedingungen von Lemma \ref{LemmaMengenfamilie2} gen"ugen.
Weiterhin sei $0 \leq \alpha \leq 1/d$.
Setzen wir
$$M_i(\alpha) = M\big( (i-1)/d+\alpha \big) \setminus M\big( (i-1)/d \big), \qquad i=1,...,d,$$
so sind diese $M_i(\alpha)$ die gesuchten Familien:
F"ur $i \not=i'$, $0 \leq \alpha' \leq \alpha \leq 1/d$ sind die Mengen $M_i(\alpha)$ und $M_{i'}(\alpha')$
disjunkt,
es ist $|M_i(\alpha)| = \alpha (t_1 - t_0)$ und f"ur $0 \leq \alpha' \leq \alpha \leq 1/d$ gilt
$M_i(\alpha') \subseteq M_i(\alpha)$.
Genau wie im Beweis des letzten Lemma k"onnen wir schlie"sen:
$$\left\| \int_{t_0}^t \Big[ \big( \chi_{M_i(\alpha)}(\tau) - \chi_{M_i(\alpha')}(\tau) \big) z(\tau)
                     - (\alpha-\alpha') z(\tau) \Big] \, d\tau \right\| \leq 3 \delta |\alpha-\alpha'|, \quad i=1,...,d.$$
Die letzte Behauptung des Lemma folgt unmittelbar aus der Tatsache, dass f"ur jeden Vektor
$z=(y_1,...,y_d) \in \R^n$, $y_i \in \R^{n_{i}}$, die Ungleichung $\|y_i\| \leq \|z\|$ erf"ullt ist. \hfill $\blacksquare$

%% file: V-2-Definition.tex
\subsection{Die Definition der mehrfachen Nadelvariation} \label{AbschnittDefinitionNadelvariation}
Die Mengenfamilien $\big\{ M(\alpha) \big\}_{0 \leq \alpha \leq 1}$,
deren Konstruktion nach Ioffe \& Tichomirov \cite{Ioffe} im letzten Abschnitt angegeben wurde,
bilden die Basis zur Einführung der mehrfachen Nadelvariation einer Steuerung $u_*(\cdot)$ in Steuerungsproblemen.
Die Wahl der beschr"ankten messbaren Vektorfunktionen $y_i(\cdot)$ in Lemma \ref{LemmaMengenfamilie3} ist im Folgenden angegeben. \\[2mm]
Es seien $x_*(\cdot) \in C([t_0,t_1],\R^n)$, $u_*(\cdot),u_1(\cdot),...,u_d(\cdot) \in L_\infty\big([t_0,t_1],U\big)$.
Ferner bezeichne $V_\gamma$ den Umgebungsstreifen
$$V_\gamma= \{ (t,x) \in [t_0,t_1] \times \R^n\,|\, \|x-x_*(t)\|\leq \gamma\}.$$
Wir nehmen an, dass die Abbildungen $f(t,x,u)$, $\varphi(t,x,u)$ auf $V_\gamma \times \R^m$ stetig in der Gesamtheit der
Variablen und stetig differenzierbar bez"uglich $x$ sind.
D.\,h., dass die partiellen Ableitungen $f_x(t,x,u)$, $\varphi_x(t,x,u)$ auf $V_\gamma \times \R^m$ stetig in der 
Gesamtheit der Variablen sind.
Unter diesen Annahmen sind die Vektorfunktionen $y_i(\cdot)$,
\begin{eqnarray*}
y_i (t) &=& \Big( \varphi\big(t,x_*(t),u_i(t) \big) - \varphi\big(t,x_*(t),u_*(t) \big), \\
             & & \hspace*{20mm} f\big(t,x_*(t),u_i(t) \big) - f\big(t,x_*(t),u_*(t) \big) \Big), \quad i = 1,...,d,
\end{eqnarray*}
auf $[t_0,t_1]$ messbar und beschr"ankt, weil die Steuerungen $u_*(\cdot),u_1(\cdot),...,u_{d}(\cdot)$ messbar, beschr"ankt und $f, \varphi$ stetig sind.
Daher existieren nach Lemma \ref{LemmaMengenfamilie3} auf dem Intervall $[t_0,t_1]$
einparametrige Mengenfamilien $\big\{ M_i(\alpha) \big\}$ mit
\begin{eqnarray*}
&& |M_i(\alpha)| = \alpha (t_1-t_0), \qquad M_i(\alpha') \subseteq M_i(\alpha),
   \qquad M_i(\alpha) \cap M_{i'}(\alpha') = \emptyset,\\
&& \max_{t \in [t_0,t_1]}
   \bigg\| \int_{t_0}^t \big( \chi_{M_i(\alpha)}(\tau) - \chi_{M_i(\alpha')}(\tau) \big) y_i(\tau) \, d\tau
   - (\alpha-\alpha') \int_{t_0}^t y_i(\tau) \, d\tau \bigg\| \leq \delta |\alpha-\alpha'|
\end{eqnarray*}
f"ur alle $i=1,...,d,$ $0 \leq \alpha' \leq \alpha \leq 1/d$ und $i \not= i'$. \\[2mm]
Demnach ist auf dem $d-$dimensionalen Quader
$$Q^{d} = \Big\{ \alpha = (\alpha_1,...,\alpha_{d}) \in \R^{d} \,\Big|\,
                                0 \leq \alpha_i  \leq 1/d, \, i =1,...,d \Big\}$$
die Abbildung $\alpha \to u_\alpha(\cdot) \in L_\infty\big([t_0,t_1],U\big),$
$$u_\alpha(t) = u_*(t) + \sum_{i=1}^{d} \chi_{M_i(\alpha_i)}(t) \cdot \big( u_i(t)-u_*(t) \big),$$
wohldefiniert und wir nennen $u_\alpha(\cdot)$ eine mehrfache Nadelvariation von $u_*(\cdot)$.

%% file: V-3-Eigenschaften.tex
\subsection{Eigenschaften der mehrfachen Nadelvariation} \label{AbschnittEigenschaftenNadelvariation}
Es seien $x_*(\cdot) \in C([t_0,t_1],\R^n)$ und $u_*(\cdot),u_1(\cdot),...,u_d(\cdot) \in L_\infty\big([t_0,t_1],U\big)$.
Weiterhin seien $f(t,x,u)$, $\varphi(t,x,u)$ auf $V_\gamma \times \R^m$ stetig in der Gesamtheit der
Variablen und stetig differenzierbar bez"uglich $x$.
Wir definieren die Mengen
\begin{eqnarray*}
\Sigma(\Delta) &=& \left\{ \alpha = (\alpha_1,...,\alpha_{d}) \in \R^{d} \,\bigg|\,
                                \alpha_1,...,\alpha_{d} \geq 0,\, \sum_{i=1}^{d} \alpha_i \leq \Delta \right\}, \\
V(\sigma) &=& \{ x(\cdot) \in C([t_0,t_1],\R^n) \,\big|\, \|x(\cdot) - x_*(\cdot)\|_{\infty} \leq \sigma \}.
\end{eqnarray*}
Dann betrachten wir f"ur $t \in [t_0,t_1]$ die Abbildung
$$\Phi_1\big(x(\cdot),\alpha\big)(t) = \int_{t_0}^t \big[\varphi\big(\tau,x(\tau),u_\alpha(\tau)\big)
                                                 - \varphi\big(\tau,x_*(\tau),u_*(\tau)\big)\big] \, d\tau$$
und den linearen Operator
\begin{eqnarray*}
\Lambda_1\big(x(\cdot),\alpha\big)(t) &=& \int_{t_0}^t \Big[\varphi_x\big(\tau,x_*(\tau),u_*(\tau)\big) \big(x(\tau)-x_*(\tau)\big) \\
&&  \hspace*{5mm}+ \sum_{i=1}^{d} \alpha_i \cdot 
   \Big( \varphi\big(\tau,x_*(\tau),u_i(\tau)\big) - \varphi\big(\tau,x_*(\tau),u_*(\tau)\big)\Big)\Big] \, d\tau.
\end{eqnarray*}
Dabei bezeichnet $u_\alpha(\cdot)$ die mehrfache Nadelvariation
$$u_\alpha(t) = u_*(t) + \sum_{i=1}^d \chi_{M_i(\alpha_i)}(t) \cdot \big( u_i(t)-u_*(t) \big).$$
Dementsprechend bezieht sich der zweite Term im Operator $\Lambda_1$ auf eine Linearisierung bez"uglich
der Variation auf den Tr"agermengen. 
Ziel dieses Abschnitts ist es zu zeigen,
dass der Operator $\Phi_1$ zwar nicht stetig differenzierbar, aber wenigstens streng differenzierbar ist.
D.\,h., dass die Beziehung
\begin{eqnarray*}
&& \Big\| \big[\Phi_1\big(x(\cdot),\alpha\big)-\Phi_1\big(x'(\cdot),\alpha'\big)
          -\Lambda_1\big(x(\cdot),\alpha\big)-\Lambda_1\big(x'(\cdot),\alpha'\big)\big] (\cdot) \Big\|_\infty \\
&& \hspace*{20mm} \leq\delta \bigg( \|x(\cdot)-x'(\cdot)\|_\infty + \sum_{i=1}^d |\alpha_i - \alpha_i'| \bigg)
\end{eqnarray*}
f"ur alle $\alpha, \alpha' \in \Sigma(\Delta)$ und alle $x(\cdot), x'(\cdot) \in V(\sigma)$ gilt.
Au"serdem wir zeigen,
dass gleichzeitig f"ur die Abbildungen $\Phi_2$ und $\Lambda_2$,
\begin{eqnarray*}
\Phi_2\big(x(\cdot),\alpha\big) &=&  \int_{t_0}^{t_1} \big[f\big(t,x(t),u_\alpha(t)\big)-f\big(t,x(t),u_*(t)\big)\big] \, dt, \\
\Lambda_2\big(x(\cdot),\alpha\big) &=& \sum_{i=1}^{d} \int_{t_0}^{t_1} \alpha_i \cdot \Big( f\big(t,x(t),u_i(t)\big)
                                                    - f\big(t,x(t),u_*(t)\big)\Big) \, dt,
\end{eqnarray*}
die Relation
$$\Phi_2\big(x(\cdot),\alpha\big) - \Lambda_2\big(x(\cdot),\alpha\big) \leq \delta \sum_{i=1}^d \alpha_i$$
f"ur alle $\alpha \in \Sigma(\Delta)$ und alle $x(\cdot) \in V(\sigma)$ erf"ullt ist. \\[2mm]
Wir treffen einige Vorbereitungen: 
Aus den Eigenschaften (\ref{Nadelvariation1}) in Lemma \ref{LemmaMengenfamilie3}
der Mengenfamilien der mehrfachen Nadelvariation
$$u_\alpha(t) = u_*(t) + \sum_{i=1}^d \chi_{M_i(\alpha_i)}(t) \cdot \big( u_i(t)-u_*(t) \big)$$
gelten f"ur jede auf $\R \times \R^n \times \R^m$ definierte Vektorfunktion $h$ und alle $\alpha, \alpha' \in Q^d$:
\begin{eqnarray}
\lefteqn{h\big(t,x,u_\alpha(t)\big) - h\big(t,x,u_*(t)\big)} \nonumber\\
&=& \label{EigenschaftNadelvariation1}
   \sum_{i=1}^{d} \chi_{M_i(\alpha_i)}(t) \cdot \Big( h\big(t,x,u_i(t)\big) - h\big(t,x,u_*(t)\big) \Big), \\
\lefteqn{h\big(t,x,u_\alpha(t) \big) - h\big(t,x,u_{\alpha'}(t) \big)} \nonumber \\
&=& \label{EigenschaftNadelvariation2} \sum_{i=1}^{d} \big( \chi_{M_i({\alpha}_i)}(t) -  \chi_{M_i(\alpha'_i)}(t) \big) 
    \cdot \Big( h\big(t,x,u_i(t)\big) - h\big(t,x,u_*(t)\big) \Big).
\end{eqnarray}
Weil alle Steuerungen $u_*(\cdot), u_1(\cdot),...,u_d(\cdot)$ beschr"ankt sind,
sind ihre Werte f"ur fast alle $t \in [t_0,t_1]$ in einer kompakten Menge $U_1 \subset \R^m$ enthalten.
Nutzen wir dar"uber hinaus die Stetigkeit von $f$, $\varphi$ und $\varphi_x$ aus,
so k"onnen wir ein $\sigma > 0$ derart angeben, dass
\begin{eqnarray}
&& \label{EigenschaftNadelvariation3}
   \big\| \varphi(t,x,u) - \varphi(t,x_*(t),u) \big\| \leq \frac{\delta}{4(t_1-t_0)}, \\
&& \label{EigenschaftNadelvariation4}
   \big| f(t,x,u) - f(t,x_*(t),u) \big| \leq \frac{\delta}{8(t_1-t_0)}, \\
&& \label{EigenschaftNadelvariation5}
   \big\| \varphi(t,x,u) - \varphi(t,x',u) - \varphi_x (t,x_*(t),u) (x-x') \big\| \leq \frac{\delta \|x-x'\|}{2(t_1-t_0)}
\end{eqnarray}
f"ur fast alle $t \in [t_0,t_1]$ und alle $x,x',u$ mit $\|x-x_*(t)\| \leq \sigma$, $\|x'-x_*(t)\| \leq \sigma$,
$u \in U_1$ erf"ullt sind.
Ferner w"ahlen wir $\Delta \in (0,1/d]$ derart, dass f"ur fast alle $t \in [t_0,t_1]$
\begin{equation} \label{EigenschaftNadelvariation6}
\Delta \cdot (t_1-t_0) \cdot \max_{ u \in U_1} \big\| \varphi_x(t,x_*(t),u) \big\| \leq \frac{\delta}{4}
\end{equation}
gilt.
Offenbar ist $\Sigma(\Delta) \subseteq Q^{d}$ f"ur $\Delta \in (0,1/d]$.
Daher ist die mehrfache Nadelvariation $u_\alpha(\cdot)$ f"ur $\Delta \in (0,1/d]$
auf dem gesamten Simplex $\Sigma(\Delta)$ definiert.

\begin{lemma} \label{LemmaEigenschaftNadelvariation1}
Es existieren $\Delta \in (0,1/d]$ und $\sigma>0$ derart, dass
\begin{eqnarray}
&& \Big\| \big[\Phi_1\big(x(\cdot),\alpha\big)-\Phi_1\big(x'(\cdot),\alpha'\big)
          -\Lambda_1\big(x(\cdot),\alpha\big)+\Lambda_1\big(x'(\cdot),\alpha'\big)\big] (\cdot) \Big\|_\infty \nonumber \\
&& \label{EigenschaftNadelvariation7}
   \hspace*{20mm} \leq\delta \bigg( \|x(\cdot)-x'(\cdot)\|_\infty + \sum_{i=1}^d |\alpha_i - \alpha_i'| \bigg)
\end{eqnarray}
f"ur alle $\alpha, \alpha' \in \Sigma(\Delta)$ und alle $x(\cdot), x'(\cdot) \in V(\sigma)$ gilt.  
\end{lemma}

Ausf"uhrlich geschrieben ist die linke Seite in (\ref{EigenschaftNadelvariation7}) gleich
\begin{eqnarray}
&&\hspace*{-10mm} \max_{t \in [t_0,t_1]} \bigg\| \int_{t_0}^t \bigg[\varphi \big( \tau, x(\tau), u_\alpha(\tau) \big)
                            - \varphi \big( \tau, x'(\tau), u_{\alpha'}(\tau) \big)
               - \varphi_x \big( \tau, x_*(\tau), u_*(\tau) \big) \big( x(\tau) - x'(\tau) \big) \nonumber \\
&& \label{EigenschaftNV1}\hspace*{20mm} - \sum_{i=1}^d (\alpha_i- \alpha'_i) \cdot 
   \Big( \varphi \big( \tau, x_*(\tau), u_i(\tau) \big) - \varphi \big( \tau, x_*(\tau), u_*(\tau) \big) \Big) \bigg] d\tau \bigg\|.
\end{eqnarray}   
      
{\bf Beweis} Den Ausdruck (\ref{EigenschaftNV1}) k"onnen wir nach oben gegen
\begin{eqnarray*}
&& \int_{t_0}^{t_1} \Big\| \varphi \big(t,x(t),u_\alpha(t)\big) - \varphi \big(t,x'(t),u_\alpha(t)\big)
   - \varphi_x \big(t,x_*(t),u_\alpha(t)\big) \big(x(t)-x'(t)\big) \Big\| dt \\
& & + \int_{t_0}^{t_1} \Big\| \Big( \varphi_x \big(t,x_*(t),u_\alpha(t)\big) 
                                 - \varphi_x \big(t,x_*(t),u_*(t)\big) \Big) \big(x(t)-x'(t)\big) \Big\| dt \\
& & +\int_{t_0}^{t_1} \Big\| \varphi \big(t,x'(t),u_\alpha(t)\big) - \varphi \big(t,x'(t),u_{\alpha'}(t)\big) \\
& & \hspace*{15mm} - \varphi \big(t,x_*(t),u_\alpha(t)\big) + \varphi \big(t,x_*(t),u_{\alpha'}(t)\big) \Big\| dt \\
& &  + \max_{t \in [t_0,t_1]} \bigg\| \int_{t_0}^t \varphi \big(\tau,x_*(\tau),u_\alpha(\tau)\big)
                   - \varphi \big(\tau,x_*(\tau),u_{\alpha'}(\tau)\big) \\
& & \hspace*{15mm} - \sum_{i=1}^{d} (\alpha_i- \alpha'_i)
      \Big( \varphi \big(\tau,x_*(\tau),u_i(\tau)\big) - \varphi \big(\tau,x_*(\tau),u_*(\tau)\big) \Big) d\tau \bigg\|
\end{eqnarray*}
absch"atzen.
Der erste Summand ist nach (\ref{EigenschaftNadelvariation5}) kleiner oder gleich
$$\frac{\delta}{2} \|x(\cdot)-x'(\cdot)\|_\infty.$$
Aus (\ref{EigenschaftNadelvariation1}) folgt, dass der zweite Summand gleich
$$\int_{t_0}^{t_1} \bigg\| \sum_{i=1}^{d} \chi_{M_i(\alpha_i)}(t) 
         \cdot \Big( \varphi_x \big(t,x_*(t),u_i(t)\big) - \varphi_x \big(t,x_*(t),u_*(t)\big) \Big)
         \big( x(t) - x'(t) \big) \bigg\| \, dt$$
ist.
Mit (\ref{Nadelvariation1}), (\ref{EigenschaftNadelvariation6}) und $\alpha \in \Sigma(\Delta)$
f"allt dieser Ausdruck nicht gr"o"ser aus als
$$\int_{t_0}^{t_1} \bigg( \sum_{i=1}^{d} \alpha_i \bigg) \cdot 2 \max_{ u \in U_1} \big\| \varphi_x(t,x_*(t),u) \big\| \, dt \cdot
   \| x(\cdot) - x'(\cdot) \|_\infty \\
  \leq \frac{\delta}{2} \|x(\cdot)-x'(\cdot)\|_\infty.$$
Wegen (\ref{EigenschaftNadelvariation2}) ist der dritte Summand gleich
\begin{eqnarray*}
\lefteqn{\int_{t_0}^{t_1} \bigg\| \sum_{i=1}^{d} \big( \chi_{M_i(\alpha_i)}(t) - \chi_{M_i(\alpha'_i)}(t) \big) 
   \cdot \Big( \varphi \big(t,x'(t),u_i(t)\big) - \varphi \big(t,x'(t),u_*(t)\big) \Big)} \\
&& - \sum_{i=1}^{d} \big( \chi_{M_i(\alpha_i)}(t) - \chi_{M_i(\alpha'_i)}(t) \big)
   \cdot \Big( \varphi \big(t,x_*(t),u_i(t)\big) - \varphi \big(t,x_*(t),u_*(t)\big) \Big) \bigg\| \, dt.
\end{eqnarray*}
Gem"a"s (\ref{EigenschaftNadelvariation3}) ist dieses Integral nicht gr"o"ser als
$$2 \cdot \frac{\delta}{4(t_1-t_0)} \cdot \sum_{i=1}^{d} \int_{t_0}^{t_1}
    \big\| \chi_{M_i(\alpha_i)}(t) - \chi_{M_i(\alpha'_i)}(t) \big\| \, dt
  = \frac{\delta}{2} \sum_{i=1}^{d} | \alpha_i - \alpha'_i |.$$
Schlie"slich l"asst sich der vierte Summand mittels (\ref{EigenschaftNadelvariation2}) in die Form
\begin{eqnarray*}
&& \hspace*{-10mm} \max_{t \in [t_0,t_1]} \bigg\| \int_{t_0}^t \sum_{i=1}^{d} \bigg[
      \big( \chi_{M_i(\alpha_i)}(\tau) - \chi_{M_i(\alpha'_i)}(\tau) \big) \cdot
       \Big( \varphi \big(\tau,x_*(\tau),u_i(\tau)\big)- \varphi \big(\tau,x_*(\tau),u_*(\tau)\big) \Big) \\
&& \hspace*{40mm} - (\alpha_i - \alpha'_i) \cdot 
      \Big( \varphi \big(\tau,x_*(\tau),u_i(\tau)\big) - \varphi \big(\tau,x_*(\tau),u_*(\tau)\big) \Big)
      \bigg] d\tau \bigg\|     
\end{eqnarray*}
bringen.
Die Differenzen
$$\Big( \varphi \big(\tau,x_*(\tau),u_i(\tau)\big) - \varphi \big(\tau,x_*(\tau),u_*(\tau)\big) \Big)$$
enthalten die ersten $n$ Koordinaten der Vektorfunktionen 
\begin{eqnarray*}
y_i (t) &=& \Big( \varphi\big(t,x_*(t),u_i(t) \big) - \varphi\big(t,x_*(t),u_*(t) \big), \\
             & & \hspace*{20mm} f\big(t,x_*(t),u_i(t) \big) - f\big(t,x_*(t),u_*(t) \big) \Big)
\end{eqnarray*}
aus Abschnitt \ref{AbschnittDefinitionNadelvariation}.
Daher ist der obenstehende Ausdruck kleiner oder gleich
$$\max_{t \in [t_0,t_1]} \bigg\| \sum_{i=1}^{d} \int_{t_0}^t
           \Big[ \big( \chi_{M_i(\alpha_i)}(\tau) - \chi_{M_i(\alpha'_i)}(\tau) \big) y_i(\tau)
               - (\alpha_i - \alpha'_i) y_i(\tau) \Big] d\tau \bigg\|.$$
Mit (\ref{Nadelvariation2}) ist dieser Ausdruck nicht gr"o"ser als
$$\frac{\delta}{2} \sum_{i=1}^{d} | \alpha_i - \alpha'_i |.$$
Damit ist die Beziehung (\ref{EigenschaftNadelvariation7}) bewiesen. \hfill $\blacksquare$ 

\begin{lemma} \label{LemmaEigenschaftNadelvariation2}
Es existieren $\Delta \in (0,1/d]$ und $\sigma>0$ derart, dass
\begin{equation} \label{EigenschaftNadelvariation8}
\Phi_2\big(x(\cdot),\alpha\big) - \Lambda_2\big(x(\cdot),\alpha\big) \leq \delta \sum_{i=1}^d \alpha_i
\end{equation}
f"ur alle $\alpha \in \Sigma(\Delta)$ und alle $x(\cdot) \in V(\sigma)$ gilt.
\end{lemma}

{\bf Beweis} Beachten wir (\ref{EigenschaftNadelvariation1}),
so ist die linke Seite in (\ref{EigenschaftNadelvariation8}) gleich
\begin{eqnarray*}
& &  \sum_{i=1}^{d} \int_{t_0}^{t_1} \bigg[ \chi_{M_i(\alpha_i)}(t) \cdot
                \Big( f \big(t,x(t),u_i(t)\big) - f\big(t,x(t),u_*(t)\big)\Big) \\
& & \hspace*{35mm} - \alpha_i \cdot \Big( f\big(t,x(t),u_i(t)\big) - f\big(t,x(t),u_*(t)\big)\Big) \bigg] dt.
\end{eqnarray*}
Dieser Ausdruck ist weiterhin kleiner oder gleich der folgenden Summe:

\begin{eqnarray*}
&& \sum_{i=1}^{d} \bigg| \int_{t_0}^{t_1} \chi_{M_i(\alpha_i)}(t) \cdot
               \Big( f \big(t,x_*(t),u_i(t)\big) - f\big(t,x_*(t),u_*(t)\big)\Big) \\
& & \hspace*{35mm}
     - \alpha_i \cdot \Big( f \big(t,x_*(t),u_i(t)\big) - f\big(t,x_*(t),u_*(t)\big)\Big) dt \bigg| \\
& & + \, \sum_{i=1}^{d} \Bigg[ \int_{t_0}^{t_1} \chi_{M_i(\alpha_i)}(t)  \cdot
         \Big[ \big|  f \big(t,x(t),u_i(t)\big) - f\big(t,x_*(t),u_i(t)\big) \big| \\
& & \hspace*{35mm} + \big| f\big(t,x_*(t),u_*(t)\big) - f \big(t,x(t),u_*(t)\big) \big| \Big] dt \\
& & \hspace*{15mm} + \,  \alpha_i \int_{t_0}^{t_1} 
         \Big[      \big| f \big(t,x(t),u_i(t)\big) - f\big(t,x_*(t),u_i(t)\big) \big| \\
& & \hspace*{35mm} + \big| f \big(t,x(t),u_*(t)\big) - f\big(t,x_*(t),u_*(t)\big) \big| \Big] dt \Bigg].
\end{eqnarray*}
F"ur den ersten Summanden gilt die Relation (\ref{Nadelvariation2}), d.\,h.
\begin{eqnarray*}
&& \sum_{i=1}^{d} \bigg| \int_{t_0}^{t_1} \chi_{M_i(\alpha_i)}(t) \cdot
               \Big( f \big(t,x_*(t),u_i(t)\big) - f\big(t,x_*(t),u_*(t)\big)\Big) \\
&& \hspace*{25mm}
     - \alpha_i \cdot \Big( f \big(t,x_*(t),u_i(t)\big) - f\big(t,x_*(t),u_*(t)\big)\Big) dt \bigg|
   \leq \frac{\delta}{2} \sum_{i=1}^{d} \alpha_i.
\end{eqnarray*}
Mit (\ref{EigenschaftNadelvariation4}) folgt, dass die zweite Summe kleiner gleich
$$\sum_{i=1}^{d} \left( 2 \int_{t_0}^{t_1} \chi_{M_i(\alpha_i)}(t) \,dt
      + 2 \alpha_i \int_{t_0}^{t_1} dt \right) \cdot \frac{\delta}{8(t_1-t_0)}
     = \frac{\delta}{2} \sum_{i=1}^{d} \alpha_i$$
ist. Damit ist das Lemma bewiesen. \hfill $\blacksquare$

%% file: V-4-Multiprozesse.tex
\subsection{Mehrfache Nadelvariationen in Multiprozessen} \label{AnhangNVMP}
Wir gehen in diesem Abschnitt auf die Variation von Multiprozessen in Kapitel \ref{KapitelMultiprozess} ein.
Unser besonderes Augenmerk liegt dabei auf die Einbindung der $k$-fachen Zerlegungen.
\index{Nadelvariation, mehrfache!Multi@-- Multiprozesse} 
Im Weiteren verwenden wir die Bezeichnungen, die in Kapitel \ref{KapitelMultiprozess} eingef"uhrt wurden. \\
Es seien $x_*(\cdot) \in C([t_0,t_1],\R^n)$ und $u_*(\cdot),u_i(\cdot) \in L_\infty\big([t_0,t_1],U\big)$, $\mathscr{A}_*,\mathscr{A}_i \in \Zt$
f"ur $i=1,...,d$.
Ferner bezeichne $V_\gamma$ die Menge
$$V_\gamma= \{ (t,x) \in [t_0,t_1] \times \R^n\,|\, \|x-x_*(t)\|\leq \gamma\}.$$
Weiterhin seien $f_i(t,x,u^{i})$, $\varphi_i(t,x,u{i})$ auf $V_\gamma \times \R^{m_i}$ stetig in der Gesamtheit der
Variablen und stetig differenzierbar bez"uglich $x$. \\[2mm]
Unter diesen Voraussetzungen sind die $(n+1)$-dimensionalen Vektorfunktionen $y_i(\cdot)$,
\begin{eqnarray*}
y_i (t) &=& \Big( \chi_{\mathscr{A}_i}(t) \circ \varphi \big(t,x_*(t),u_i(t) \big)
                       - \chi_{\mathscr{A}_*}(t) \circ \varphi \big(t,x_*(t),u_*(t) \big), \\
             & & \hspace*{10mm} \chi_{\mathscr{A}_i}(t) \circ f \big(t,x_*(t),u_i(t) \big)
                       - \chi_{\mathscr{A}_*}(t) \circ f \big(t,x_*(t),u_*(t) \big) \Big), \quad i = 1,...,d,
\end{eqnarray*}
me"sbar und beschr"ankt.
Daher existieren nach Lemma \ref{LemmaMengenfamilie3} Mengenfamilien $\big\{ M_i(\alpha) \big\}$ mit
\begin{eqnarray*}
&& \hspace*{-15mm} |M_i(\alpha)| = \alpha (t_1-t_0), \qquad M_i(\alpha') \subseteq M_i(\alpha),
   \qquad M_i(\alpha) \cap M_{i'}(\alpha') = \emptyset,\\
&& \hspace*{-15mm} \max_{t \in [t_0,t_1]}
   \bigg\| \int_{t_0}^t \big( \chi_{M_i(\alpha)}(\tau) - \chi_{M_i(\alpha')}(\tau) \big) y_i(\tau) \, d\tau
   - (\alpha-\alpha') \int_{t_0}^t y_i(\tau) \, d\tau \bigg\| \leq \delta |\alpha-\alpha'|
\end{eqnarray*}
f"ur alle $i=1,...,d,$ $0 \leq \alpha' \leq \alpha \leq 1/d$ und $i \not= i'$. 
Wir beachten ferner,
dass f"ur Zerlegungen $\mathscr{A}, \mathscr{B} \in \Zt$ und f"ur eine messbare Menge $M \subseteq [t_0,t_1]$ die Funktion
$$\chi_{\mathscr{C}}(t) = \chi_{\mathscr{A}}(t) + \chi_M(t) \cdot \big(\chi_{\mathscr{B}}(t)-\chi_{\mathscr{A}}(t)\big)$$
stets mit einer Menge $\mathscr{C} \in \Zt$ korrespondiert.
Dann ist auf dem Quader
$$Q^d = \Big\{ \alpha = (\alpha_1,...,\alpha_d) \in \R^d \,\Big|\,
                                0 \leq \alpha_i  \leq 1/d, \, i =1,...,d \Big\}$$
die Abbildung $\alpha \to \big( u_\alpha(\cdot), \mathscr{A}_\alpha \big) \in L_\infty\big([t_0,t_1],U\big) \times \Zt$ wohldefiniert.
Dabei bezeichnen
$$u_\alpha(t) = u_*(t) + \sum_{i=1}^d \chi_{M_i(\alpha_i)}(t) \cdot \big( u_i(t)-u_*(t) \big)$$
und $\mathscr{A}_\alpha$ identifizieren wir durch die charakteristische Vektorfunktion
$$\chi_{\mathscr{A}_\alpha}(t) = \chi_{\mathscr{A}_*}(t)
    + \sum_{i=1}^d \chi_{M_i(\alpha_i)}(t) \cdot \big( \chi_{\mathscr{A}_i}(t) - \chi_{\mathscr{A}_*}(t) \big).$$
Aus den Eigenschaften der Mengenfamilien in der Setzung der verallgemeinerten Nadelvariation
folgt f"ur jede wohldefinierte Vektorfunktion $h=(h_1,...,h_k)$ und jedes $\alpha \in Q^d$
\begin{eqnarray}
\lefteqn{\chi_{\mathscr{A}_\alpha}(t) \circ h \big(t,x,u_\alpha(t) \big)
    - \chi_{\mathscr{A}_*}(t) \circ h \big(t,x,u_*(t) \big)} \nonumber \\
&=& \label{EigenschaftNadelvariationHybrid1}
    \sum_{i=1}^{d} \chi_{M_i(\alpha_i)}(t) \cdot \Big(  \chi_{\mathscr{A}_i}(t) \circ h \big(t,x,u_i(t) \big) 
    - \chi_{\mathscr{A}_*}(t) \circ h \big(t,x,u_*(t) \big) \Big).
\end{eqnarray}
Denn verwenden wir zuerst die Setzung von $u_\alpha(\cdot)$, so erhalten wir
\begin{eqnarray*}
\lefteqn{ \chi_{\mathscr{A}_\alpha}(t) \circ h \big(t,x,u_\alpha(t) \big) = \chi_{\mathscr{A}_\alpha}(t) \circ 
    h\bigg( t,x,u_*(t) + \sum_{i=1}^{d} \chi_{M_i(\alpha_i)}(t) \big( u_i(t) - u_*(t) \big) \bigg)  } \\
&=& \chi_{\mathscr{A}_\alpha}(t) \circ 
    \bigg( h\big(t,x,u_*(t)\big) +
    \sum_{i=1}^{d} \chi_{M_i(\alpha_i)}(t) \cdot \Big( h\big(t,x,u_i(t)\big) - h\big(t,x,u_*(t)\big) \Big) \bigg) \\
&=& \chi_{\mathscr{A}_\alpha}(t) \circ h\big(t,x,u_*(t)\big) + \chi_{\mathscr{A}_\alpha}(t) \circ
    \sum_{i=1}^{d} \chi_{M_i(\alpha_i)}(t) \cdot \Big( h\big(t,x,u_i(t)\big) - h\big(t,x,u_*(t)\big) \Big).
\end{eqnarray*}
Nutzen wir nun die Gestalt von $\chi_{\mathscr{A}_\alpha}(\cdot)$ aus und verwenden im zweiten Summanden,
dass die Mengenfamilien $\{M_i(\alpha_i)\}$ nach Konstruktion die Eigenschaft 
$$\bigg( \sum_{i=1}^{d} \chi_{M_i(\alpha_i)}(t) \bigg)^2 = \sum_{i=1}^{d} \chi^2_{M_i(\alpha_i)}(t)
  = \sum_{i=1}^{d} \chi_{M_i(\alpha_i)}(t)$$
besitzen, dann ergibt sich
\begin{eqnarray*}
\lefteqn{\chi_{\mathscr{A}_\alpha}(t) \circ h\big(t,x,u_*(t)\big) + \chi_{\mathscr{A}_\alpha}(t) \circ
    \sum_{i=1}^{d} \chi_{M_i(\alpha_i)}(t) \cdot \Big( h\big(t,x,u_i(t)\big) - h\big(t,x,u_*(t)\big) \Big)} \\
&=& \chi_{\mathscr{A}_*}(t) \circ h \big(t,x,u_*(t) \big) + 
    \sum_{i=1}^{d} \chi_{M_i(\alpha_i)}(t) \cdot \big( \chi_{\mathscr{A}_i}(t) - \chi_{\mathscr{A}_*}(t) \big)
    \circ h \big(t,x,u_*(t) \big) \\
& & + \sum_{i=1}^{d} \chi_{M_i(\alpha_i)}(t) \cdot \chi_{\mathscr{A}_*}(t) \circ
      \Big( h\big(t,x,u_i(t)\big) - h\big(t,x,u_*(t)\big) \Big) \\
& & + \sum_{i=1}^{d} \chi_{M_i(\alpha_i)}(t) \cdot \big( \chi_{\mathscr{A}_i}(t) - \chi_{\mathscr{A}_*}(t) \big)
    \circ \Big( h\big(t,x,u_i(t)\big) - h\big(t,x,u_*(t)\big) \Big).
\end{eqnarray*}
Daraus folgt nun nach Vereinfachung
\begin{eqnarray*}
\lefteqn{ \chi_{\mathscr{A}_\alpha}(t) \circ h \big(t,x,u_\alpha(t) \big) 
= \chi_{\mathscr{A}_*}(t) \circ h \big(t,x,u_*(t) \big) } \\
& & + \sum_{i=1}^{d} \chi_{M_i(\alpha_i)}(t) \cdot \Big(  \chi_{\mathscr{A}_i}(t) \circ h \big(t,x,u_i(t) \big) 
    - \chi_{\mathscr{A}_*}(t) \circ h \big(t,x,u_*(t) \big) \Big),
\end{eqnarray*}
und nach Subtraktion von $\chi_{\mathscr{A}_*}(t) \circ h \big(t,x,u_*(t) \big)$ die Beziehung (\ref{EigenschaftNadelvariationHybrid1}). \\
Betrachten wir die Verkn"upfung der Vektorfunktion $h=(h_1,...,h_k)$ mit den Zerlegungen $\mathscr{A}_\alpha, \mathscr{A}_{\alpha'} \in \Zt$,
$\alpha, \alpha' \in Q^{d}$, dann erhalten wir mit (\ref{EigenschaftNadelvariationHybrid1})
\begin{eqnarray}
&& \hspace*{-15mm} \chi_{\mathscr{A}_\alpha}(t) \circ h \big(t,x,u_\alpha(t) \big)
               - \chi_{\mathscr{A}_{\alpha'}}(t) \circ h \big(t,x,u_{\alpha'}(t) \big) \nonumber \\
&& \hspace*{-15mm} = \label{EigenschaftNadelvariationHybrid2}
    \sum_{i=1}^{d} \big( \chi_{M_i({\alpha}_i)}(t) -  \chi_{M_i(\alpha'_i)}(t) \big) \cdot
    \Big( \chi_{\mathscr{A}_i}(t) \circ h \big(t,x,u_i(t) \big) - \chi_{\mathscr{A}_*}(t) \circ h\big(t,x,u_*(t) \big) \Big).
\end{eqnarray}
Die Beziehungen (\ref{EigenschaftNadelvariationHybrid1}) und (\ref{EigenschaftNadelvariationHybrid2}) f"ur Multiprozesse entsprechen
gerade den Eigenschaften (\ref{EigenschaftNadelvariation1}), (\ref{EigenschaftNadelvariation2}).
Beachten wir nun au"serdem die Ungleichung
$$\|\chi_{\mathscr{A}}(t) \circ h(t,x,u)\| \leq \sum_{s=1}^k \|h^{s}(t,x,u^{s})\|,$$
dann ergeben sich auf die gleiche Weise wie im Abschnitt \ref{AbschnittEigenschaftenNadelvariation} f"ur hinreichend kleine Parameter $\Delta \in (0,1/d]$
auf den Mengen
\begin{eqnarray*}
\Sigma(\Delta) &=& \left\{ \alpha = (\alpha_1,...,\alpha_{d}) \in \R^{d} \,\bigg|\,
                                \alpha_1,...,\alpha_{d} \geq 0,\, \sum_{i=1}^{d} \alpha_i \leq \Delta \right\}, \\
V &=& \{ x(\cdot) \in C([t_0,t_1],\R^n) \,\big|\, \|x(\cdot) - x_*(\cdot)\|_{\infty} \leq \sigma \}
\end{eqnarray*}
die folgenden Aussagen:

\begin{lemma} \label{LemmaEigenschaftNadelvariationHybrid1}
F"ur alle $\alpha, \alpha' \in \Sigma(\Delta)$ und alle $x(\cdot), x'(\cdot) \in V$ gilt:
\begin{eqnarray*}
\lefteqn{\max_{t \in [t_0,t_1]} \bigg\| \int_{t_0}^t \bigg[
    \chi_{\mathscr{A}_\alpha}(\tau) \circ \varphi \big( \tau, x(\tau), u_\alpha(\tau) \big)
   - \chi_{\mathscr{A}_{\alpha'}}(\tau) \circ \varphi \big( \tau, x'(\tau), u_{\alpha'}(\tau) \big)} \\
&& \hspace*{20mm} - \chi_{\mathscr{A}_*}(\tau) \circ \varphi_x \big( \tau, x_*(\tau), u_*(\tau) \big) \big( x(\tau) - x'(\tau) \big) \\
&& \hspace*{10mm} - \sum_{i=1}^{d} (\alpha_i- \alpha'_i)
   \Big( \chi_{\mathscr{A}_i}(\tau) \circ \varphi \big( \tau, x_*(\tau), u_i(\tau) \big)
    - \chi_{\mathscr{A}_*}(\tau) \circ \varphi \big( \tau, x_*(\tau), u_*(\tau) \big) \Big) \bigg] d\tau \bigg\| \\
&& \leq \delta \left( \big\|x(\cdot) - x'(\cdot) \big\|_\infty + 
   \sum_{i=1}^{d} |\alpha_i- \alpha'_i| \right).
\end{eqnarray*}      
\end{lemma}

\begin{lemma} \label{LemmaEigenschaftNadelvariationHybrid2}
F"ur alle $\alpha \in \Sigma(\Delta)$ und alle $x(\cdot) \in V$ gilt:
\begin{eqnarray*}
\lefteqn{\int_{t_0}^{t_1} \bigg[ \chi_{\mathscr{A}_\alpha}(t) \circ f \big(t,x(t),u_\alpha(t)\big)
                                   - \chi_{\mathscr{A}_*}(t) \circ f\big(t,x(t),u_*(t)\big)} \\
&& \hspace*{5mm} - \sum_{i=1}^{d}
     \alpha_i \cdot \Big( \chi_{\mathscr{A}_i}(t) \circ f\big(t,x(t),u_i(t)\big) 
   -  \chi_{\mathscr{A}_*}(t) \circ f\big(t,x(t),u_*(t)\big)\Big) \bigg] \, dt \leq \delta \sum_{i=1}^{d} \alpha_i.
\end{eqnarray*}
\end{lemma}

%% file: V-5-UnendlicherHorizont.tex
\subsection{Mehrfache Nadelvariationen \"uber dem unendlichen Zeithorizont} \label{AnhangNVUH}
Der Beweis des Pontrjaginschen Maximumprinzips f"ur die Aufgabe mit unendlichem Zeithorizont
im Kapitel \ref{KapitelStrong} basiert ebenfalls auf einer mehrfachen Nadelvariation.
\index{Nadelvariation, mehrfache!Unendlich@-- unendlicher Zeithorizont} 
Die wesentlichen Grundlagen dazu haben wir in den vorhergehenden Abschnitten gelegt.
Allerdings werden wir mit zwei neuen Herausforderungen konfrontiert,
n"amlich mit dem unbeschr"ankten Zeitintervall $\R_+$ und mit der m"oglicherweise lokal unbeschr"ankten Dichtefunktion
$\omega(\cdot) \in L_1(\R_+,\R_+)$ im Integranden des Zielfunktionals.
Ein typischer Vertreter ist die Weibull-Verteilung $\omega(t)=t^{k-1}e^{-t^k}$ mit einem Formparameter $k \in (0,1)$. \\[2mm]
Ausgangspunkt f"ur die Betrachtung der mehrfachen Nadelvariation ist die Angabe geeigneter Tr"agerfamilien
(Lemma \ref{LemmaMengenfamilie3}) und die Eigenschaften,
die nach Lemma \ref{LemmaEigenschaftNadelvariation1} und Lemma \ref{LemmaEigenschaftNadelvariation2}
dem Steuerungsproblem "ubergeben werden.
Wir werden nun zeigen,
wie wir die wesentliche Argumentation f"ur das unbeschr"ankte Zeitintervall und f"ur einen lokal unbeschr"ankten Integranden
auf eine kompakte Menge $K$ zur"uckf"uhren k"onnen. \\[2mm]
Zu einer messbaren Menge $A$ bezeichnet $|A|$ das Lebesgue-Ma"s dieser Menge.
Sei $K \subset \R_+$ eine kompakte Menge mit $|K|>0$.
Dann existiert ein $T>0$ mit $K \subseteq [0,T]$.
Au"serdem gibt es zu jedem $n \in \N$ Zahlen $0=t_0 < t_1 <... < t_n =T$ mit $|[t_{i-1},t_i] \cap K| = |K| / n$, $i=1,...,n$.
Damit definieren die Mengen $\Delta_i= [t_{i-1},t_i] \cap K$ eine geordenete Zerlegung gleicher L"ange der Menge $K$.
Mit diesen Vorbereitungen ergibt sich das folgende Lemma:

\begin{lemma} \label{LemmaMengenfamilie4}
Es sei $K \subset \R_+$ kompakt und es seien $y_i(\cdot),\, y_i : K \to \R^{n_i},\, i=1,...,d,$ beschr"ankte messbare Vektorfunktionen.
Dann existieren zu jedem $\delta > 0$ einparametrige Mengenfamilien $M_1(\alpha),...,M_d(\alpha)$, $0\leq \alpha \leq 1/d$,
messbarer Teilmengen der Menge $K$ derart, dass f"ur alle
$i=1,...,d,$ $0 \leq \alpha' \leq \alpha \leq 1/d$ und $i \not= i'$ gilt:
\begin{eqnarray*}
&& \hspace*{-5mm} |M_i(\alpha)| = \alpha |K|, \qquad M_i(\alpha') \subseteq M_i(\alpha),
   \qquad M_i(\alpha) \cap M_{i'}(\alpha') = \emptyset,\\
&& \hspace*{-5mm} \max_{t \in K}
   \bigg\| \int_{[0,t] \cap K} \big( \chi_{M_i(\alpha)}(\tau) - \chi_{M_i(\alpha')}(\tau) \big) y_i(\tau) \, d\tau 
            - (\alpha-\alpha') \int_{[0,t] \cap K} y_i(\tau) \, d\tau \bigg\| \leq \delta |\alpha-\alpha'|.
\end{eqnarray*}
\end{lemma}

{\bf Beweis} Es sei  $\{\Delta_1,...,\Delta_r\}$ eine geordnete Zerlegung der Menge $K$ gleicher L"ange nicht gr"o"ser als $\delta / (2C)$.
Betrachten wir damit die Mengen
$$M_{ij}(\alpha) = \bigg\{ t \in (A_j \cap \Delta_i)
                      \,\bigg|\, \int_0^t \chi_{(A_j \cap \Delta_i)} (\tau) \, d\tau
                                                                < \alpha \cdot |A_j \cap \Delta_i| \bigg\},$$
dann ergibt sich der Nachweis unmittelbar aus dem Beweis von Lemma \ref{LemmaMengenfamilie1} . \hfill $\blacksquare$ \\[2mm]
F"ur die weiteren Betrachtungen geben wir einen "Uberblick "uber die ben"otigten Eigenschaften an die Aufgabe (\ref{PAUH1})--(\ref{PAUH5}),
die wir im Kapitel \ref{KapitelStrong} getroffen haben: \\
Es seien $x_*(\cdot) \in C_{\lim}(\R_+,\R^n)$ und $u_*(\cdot) \in L_\infty(\R_+,U)$.
Zu $x_*(\cdot)$ bezeichnet $V_\gamma$ die Menge
$$V_\gamma= \{ (t,x) \in \overline{\R}_+ \times \R^n \,|\, \|x-x_*(t)\| \leq \gamma\}.$$
Wir nehmen an,
dass es zu jeder kompakten Menge $U_1 \subseteq \R^m$ eine Zahl $\gamma>0$ derart gibt,
dass auf der Menge $V_\gamma \times U_1$ die Abbildungen
$f(t,x,u)$ und $\varphi(t,x,u)$ gleichm"a"sig stetig und gleichm"a"sig stetig differenzierbar bez"uglich $x$ sind.
Ferner seien f"ur das Paar $\big(x_*(\cdot),u_*(\cdot)\big)$ folgende Eigenschaften erf"ullt:
Es sind die Lebesgue-Integrale
$$\int_0^\infty \big\|\varphi\big(t,x_*(t),u_*(t)\big)\big\| \, dt, \qquad \int_0^\infty \big\|\varphi_x\big(t,x_*(t),u_*(t)\big)\big\| \, dt$$
endlich.
Au"serdem nehmen wir an, es existiert zu jedem $\delta>0$ ein $T>0$ mit
\begin{eqnarray}
&& \int_T^\infty \big\| \varphi\big(t,x(t),u_*(t)\big)-\varphi\big(t,x'(t),u_*(t)\big) - \varphi_x\big(t,x_*(t),u_*(t)\big)\big(x(t)-x'(t)\big) \big\| \, dt
   \nonumber \\
&& \label{AnhangNVUH1} \hspace*{20mm} \leq \delta \|x(\cdot)-x'(\cdot)\|_\infty
\end{eqnarray}
f"ur alle $x(\cdot), x'(\cdot) \in C_{\lim}(\R_+,\R^n)$ mit $\|x(\cdot)-x_*(\cdot)\|_\infty \leq \gamma$, $\|x'(\cdot)-x_*(\cdot)\|_\infty \leq \gamma$. \\
Ferner bezeichnet $\mathscr{U}$ die Menge aller $u(\cdot) \in L_\infty(\R_+,U)$,
die die Darstellung
$$u(t)=u_*(t) + \chi_M(t)\big(w(t)-u_*(t)\big)$$
mit $w(\cdot) \in L_\infty(\R_+,U)$ und einer me"sbar und beschr"ankten Menge $M \subset \R_+$ besitzen. \\[2mm]
Es seien $u_1(\cdot),...,u_d(\cdot) \in \mathscr{U}$ und $\delta >0$ gegeben.
Dann l"asst sich eine Zahl $T>0$ derart w"ahlen, dass
die Mengen $M_i$, die in den Darstellungen der Steuerungen $u_i(\cdot) \in \mathscr{U}$ auftreten, im Intervall $[0,T]$ enthalten sind und
die Relation (\ref{AnhangNVUH1}) mit $\delta/3$ erf"ullt ist. \\
Ferner existiert nach dem Satz von Lusin eine kompakte Menge $K \subseteq [0,T]$,
auf der die Funktion $\omega(\cdot)$ stetig ist und zudem f"ur $i=1,...,d$ folgende Relationen gelten:
\begin{eqnarray}
&& \hspace*{-1.5cm} \label{WahlT2}
   \int_{[0,T] \setminus K} \big\|\varphi\big(t,x_*(t),u_i(t)\big)-\varphi\big(t,x_*(t),u_*(t)\big)\big\| \, dt \leq \frac{\delta}{3}, \\
&& \hspace*{-1.5cm} \label{WahlT3}
   \int_{[0,T] \setminus K} \omega(t)\big|f\big(t,x_*(t),u_i(t)\big)-f\big(t,x_*(t),u_*(t)\big)\big| \, dt \leq \frac{\delta}{2}.
\end{eqnarray}
Wir betrachten im Folgenden f"ur $i=1,...,d$ die Funktionen $y_i(\cdot)$,
$$y_i (t) = \Big( \varphi\big(t,x_*(t),u_i(t) \big) - \varphi\big(t,x_*(t),u_*(t) \big),
                  \omega(t)\big[ f\big(t,x_*(t),u_i(t) \big) - f\big(t,x_*(t),u_*(t) \big)\big] \Big),$$
die nach Wahl von $K$ "uber dieser Menge messbar und beschr"ankt sind.
Mit den zugeh"origen Tr"agermengen $M_1(\alpha),...,M_d(\alpha) \subseteq K$ aus Lemma \ref{LemmaMengenfamilie4} definieren wir "uber $\R_+$ 
die Abbildung $\alpha \to u_\alpha(\cdot) \in \mathscr{U}$,
$\displaystyle u_\alpha(t) = u_*(t) + \sum_{i=1}^{d} \chi_{M_i(\alpha_i)}(t) \cdot \big( u_i(t)-u_*(t) \big)$.
Diese ist auf dem Quader
$\displaystyle Q^{d} = \{ \alpha = (\alpha_1,...,\alpha_{d}) \in \R^{d}
                                    \,|\, 0 \leq \alpha_i  \leq 1/d, \, i =1,...,d \}$
wohldefiniert. \\
Im Weiteren bezeichnen $\Sigma(\Delta)$ und $V(\sigma)$ die Mengen
\begin{eqnarray*}
\Sigma(\Delta) &=& \left\{ \alpha = (\alpha_1,...,\alpha_{d}) \in \R^{d} \,\bigg|\,
                                \alpha_1,...,\alpha_{d} \geq 0,\, \sum_{i=1}^{d} \alpha_i \leq \Delta \right\}, \\
V(\sigma) &=& \{ x(\cdot) \in C_{\lim}(\R_+,\R^n) \,\big|\, \|x(\cdot) - x_*(\cdot)\|_{\infty} \leq \sigma \}.
\end{eqnarray*}
Die Menge $V(\sigma)$ sei so gewählt,
dass die Eigenschaften (\ref{EigenschaftNadelvariation3})--(\ref{EigenschaftNadelvariation5}) über $[0,T]$ statt über $[t_0,t_1]$ gelten.
Auf diesen Mengen betrachten wir die Abbildungen
\begin{eqnarray*}
    \Phi_1\big(x(\cdot),\alpha\big)(t)
&=& \int_0^t \big[\varphi\big(\tau,x(\tau),u_\alpha(\tau)\big) - \varphi\big(\tau,x_*(\tau),u_*(\tau)\big)\big] \, d\tau, \quad t \in \R_+, \\
    \Lambda_1\big(x(\cdot),\alpha\big)(t)
&=& \int_0^t \bigg[\varphi_x\big(\tau,x_*(\tau),u_*(\tau)\big) \big(x(\tau)-x_*(\tau)\big) \\
& & \hspace*{7mm} + \sum_{i=1}^{d} \alpha_i \cdot 
    \Big( \varphi\big(\tau,x_*(\tau),u_i(\tau)\big) - \varphi\big(\tau,x_*(\tau),u_*(\tau)\big)\Big)\bigg] \, d\tau, \quad t \in \R_+
\end{eqnarray*}
und die Funktionale
\begin{eqnarray*}
    \Phi_2\big(x(\cdot),\alpha\big)
&=& \int_0^\infty \omega(t)\big[f\big(t,x(t),u_\alpha(t)\big)-f\big(t,x(t),u_*(t)\big)\big] \, dt, \\
    \Lambda_2\big(x(\cdot),\alpha\big)
&=& \sum_{i=1}^{d} \int_0^\infty \alpha_i \cdot \omega(t) \big[ f\big(t,x(t),u_i(t)\big)
                                                    - f\big(t,x(t),u_*(t)\big)\big] \, dt.
\end{eqnarray*}

\begin{lemma} \label{LemmaNVUH1}
Es existieren $\Delta \in (0,1/d]$ und $\sigma>0$ derart, dass
\begin{eqnarray}
&& \Big\| \big[\Phi_1\big(x(\cdot),\alpha\big)-\Phi_1\big(x'(\cdot),\alpha'\big)
          -\Lambda_1\big(x(\cdot),\alpha\big)-\Lambda_1\big(x'(\cdot),\alpha'\big)\big] (\cdot) \Big\|_\infty \nonumber \\
&& \label{NVUH1}
   \hspace*{20mm} \leq\delta \bigg( \|x(\cdot)-x'(\cdot)\|_\infty + \sum_{i=1}^d |\alpha_i - \alpha_i'| \bigg)
\end{eqnarray}
 f"ur alle $\alpha, \alpha' \in \Sigma(\Delta)$ und alle $x(\cdot), x'(\cdot) \in V(\sigma)$ gilt.
\end{lemma}

{\bf Beweis} Beachten wir, dass die Abbildung $\alpha \to u_\alpha(\cdot)$ nur f"ur $t \in K$ einflie"st,
so k"onnen wir die linke Seite in (\ref{NVUH1}) gegen folgenden Ausdruck nach oben absch"atzen:

\begin{eqnarray*}
& & \int_{[0,T] \setminus K}  \big\| \varphi\big(t,x(t),u_*(t)\big)-\varphi\big(t,x'(t),u_*(t)\big)
                                      -\varphi_x\big(t,x_*(t),u_*(t)\big)\big(x(t)-x'(t)\big)\big\| \, dt  \\
& & \hspace*{15mm} + \sum_{i=1}^d (\alpha_i- \alpha'_i) \cdot
    \int_{[0,T] \setminus K} \big\|\varphi\big(t,x_*(t),u_i(t)\big)-\varphi\big(t,x_*(t),u_*(t)\big)\big\| \, dt \\
&+& \max_{t \in [0,T]} \bigg\| \int_{[0,t] \cap K} \bigg[\varphi \big( \tau, x(\tau), u_\alpha(\tau) \big)
                            - \varphi \big( \tau, x'(\tau), u_{\alpha'}(\tau) \big) \\
& & \hspace*{30mm} - \varphi_x \big( \tau, x_*(\tau), u_*(\tau) \big) \big( x(\tau) - x'(\tau) \big) \\
& & \hspace*{30mm} - \sum_{i=1}^d (\alpha_i- \alpha'_i) \cdot 
    \Big( \varphi \big( \tau, x_*(\tau), u_i(\tau) \big) - \varphi \big( \tau, x_*(\tau), u_*(\tau) \big) \Big) \bigg] d\tau \bigg\| \\
&+& \int_T^\infty \big\| \varphi\big(t,x(t),u_*(t)\big)-\varphi\big(t,x'(t),u_*(t)\big)
                         - \varphi_x\big(t,x_*(t),u_*(t)\big)\big(x(t)-x'(t)\big) \big\| \, dt.
\end{eqnarray*}
Wegen (\ref{EigenschaftNadelvariation5}) und (\ref{WahlT2}) ist der erste Summand kleiner oder gleich
$$\frac{\delta}{3} \bigg( \|x(\cdot)-x'(\cdot)\|_\infty + \sum_{i=1}^d |\alpha_i - \alpha_i'| \bigg).$$
Nach Wahl der Zahl $T$ f"allt der letzte Summand kleiner gleich $\delta/3 \cdot \|x(\cdot)-x'(\cdot)\|_\infty$ aus. \\
Bez"uglich der kompakten Menge $K$ l"asst sich genauso wie im Beweis von Lemma \ref{LemmaEigenschaftNadelvariation1} ein
$\Delta \in (0,1/d]$ derart angeben, dass der mittlere Summand kleiner oder gleich
$$\frac{\delta}{3} \bigg( \|x(\cdot)-x'(\cdot)\|_\infty + \sum_{i=1}^d |\alpha_i - \alpha_i'| \bigg)$$
f"ur alle $\alpha, \alpha' \in \Sigma(\Delta)$ und alle $x(\cdot), x'(\cdot) \in V(\sigma)$ ausf"allt. \hfill $\blacksquare$

\begin{lemma} \label{LemmaNVUH2}
Es existieren $\Delta \in (0,1/d]$ und $\sigma>0$ derart,
dass die Ungleichung
\begin{equation} \label{NVUH2}
\Phi_2\big(x(\cdot),\alpha\big) - \Lambda_2\big(x(\cdot),\alpha\big) \leq \delta \sum_{i=1}^d \alpha_i
\end{equation}
f"ur alle $\alpha \in \Sigma(\Delta)$ und alle $x(\cdot) \in V(\sigma)$ gilt.
\end{lemma}

{\bf Beweis} Die linke Seite in (\ref{NVUH2}) ist gleich
$$\sum_{i=1}^{d} \int_0^\infty \Big(\chi_{M_i(\alpha_i)}(t) -\alpha_i \Big) \cdot
                \Big( \omega(t) \big[f \big(t,x(t),u_i(t)\big) - f\big(t,x(t),u_*(t)\big)\big]\Big) \, dt.$$
Beachten wir $M_i(\alpha_i) \subseteq K$, so folgt nach (\ref{WahlT3}) bei der Wahl der Menge $K$:
$$\sum_{i=1}^{d} \int_{[0,T] \setminus K} \alpha_i \cdot
                \Big( \omega(t) \big[f \big(t,x(t),u_i(t)\big) - f\big(t,x(t),u_*(t)\big)\big]\Big) \, dt
    \leq \frac{\delta}{2} \sum_{i=1}^d \alpha_i.$$
Nach Wahl von $T$ gilt $u_i(t)=u_*(t)$ f"ur alle $t >T$ und $i=1,...,d$.
Daher ist
$$\sum_{i=1}^{d} \int_T^\infty \alpha_i \cdot \Big( \omega(t) \big[f \big(t,x(t),u_i(t)\big) - f\big(t,x(t),u_*(t)\big)\big]\Big) \, dt = 0.$$
Ferner erhalten wir wie im Beweis von Lemma \ref{LemmaEigenschaftNadelvariation2} "uber der kompakten Menge $K$ die Relation
$$\sum_{i=1}^{d} \int_K \Big(\chi_{M_i(\alpha_i)}(t) -\alpha_i \Big) \cdot
                \Big( \omega(t) \big[f \big(t,x(t),u_i(t)\big) - f\big(t,x(t),u_*(t)\big)\big]\Big) \, dt
  \leq \frac{\delta}{2} \sum_{i=1}^d \alpha_i$$
f"ur alle $\alpha \in \Sigma(\Delta)$ und alle $x(\cdot) \in V(\sigma)$. \hfill $\blacksquare$

%% file: V-6-Integralgleichung.tex
\subsection{Mehrfache Nadelvariationen für Integralgleichungen} \label{AnhangNVIG}
Die Herausforderungen bei der Untersuchungen von mehrfachen Nadelvariationen für Integralgleichungen
werden durch verschiedene Elemente der Untersuchungen im Kapitel \ref{KapitelIGL} hervorgerufen:
\begin{enumerate}
\item[$\cdot$] Einerseits sehen wir uns durch die innere und die äußere Zeitvariable mit einem zweidimensionalen Zeitbereich und
               einer Familie von rechten Seiten in der Dynamik konfrontiert.
\item[$\cdot$] Andererseits führt die Methode der Substitution der Zeit in der Aufgabe mit freiem Anfangs- und Endzeitpunkt im Abschnitt \ref{AbschnittFreieZeitIGL}
               zu Zuständen der äußeren Zeit, \index{Zustände der äußeren Zeit}
               welche in die Dynamik und den Integranden des Zielfunktionals einfließen.
\end{enumerate}

In den vorherigen Abschnitten wurden die mehrfachen Nadelvariationen zu messbar und beschränkten Funktionen
auf Basis einparametriger Mengenfamilien $\big\{ M_i(\alpha) \big\}$ ausgearbeitet.
\index{Nadelvariation, mehrfache!Integral@-- Integralgleichungen} 
Im Fokus stand dabei vorrangig die rechte Seite der Dynamik (\ref{PMP2}) der Standardaufgabe (\ref{PMP1})--(\ref{PMP5}),
die in Integraldarstellung die Gestalt
\begin{equation}\label{AnhangNVIG1} x(t)=x(t_0)+\int_{t_0}^t \varphi\big(s,x(s),u(s)\big) \, ds\end{equation}
besitzt.
Im Rahmen der Steuerung von Integralgleichungen hat die Dynamik die Form
\begin{equation}\label{AnhangNVIG2}x(t)=x(t_0)+\int_{t_0}^t \varphi\big(t,s,x(s),u(s)\big) \, ds.\end{equation}
Während der Integrand in (\ref{AnhangNVIG1}) stets nur eine Abbildung $s \to\varphi\big(s,x(s),u(s)\big)$ enthält,
liegt in (\ref{AnhangNVIG2}) zu jedem $t \in [t_0,t_1]$ eine Familie
$s \to \big\{\varphi\big(\tau,s,x(s),u(s)\big) \,|\, \tau \in [t_0,t]\big\}$ vor.
Die Herausforderung besteht nun in dem Nachweis,
dass bezüglich der Abbildungsfamilie in der Dynamik (\ref{AnhangNVIG2}) ebenfalls einparametrige Mengenfamilien
$\big\{ M_i(\alpha) \big\}$ in der Form angegeben werden können,
dass eine Entsprechung von Lemma \ref{LemmaEigenschaftNadelvariation1} gültig ist -- und das unabhängig
vom Parameter $t \in [t_0,t_1]$.
Dazu werden wir im Folgenden die Abbildung $\varphi(t,s,x,u)$ bezüglich der Variable $t$ geeignet approximieren
und anschließend die Mengenfamilien $\big\{ M_i(\alpha) \big\}$ bezüglich der Variable $s$ ausarbeiten. \\[2mm]
Die Methode der Substitution der Zeit im Abschnitt \ref{AbschnittFreieZeitIGL} führt zu einem Zielfunktional und einer Dynamik,
welche sich in die folgende allgemeine Form eingliedern lassen:
\begin{eqnarray*}
&& \hspace*{-1cm} J\big(x(\cdot),u(\cdot)\big)  = \int_{t_0}^{t_1} f\big(s,x(s),x(t_0),x(t_1),u(s)\big) \, ds, \\
&& \hspace*{-1cm} x(t) = x(t_0) + \int_{t_0}^t \varphi\big(t,s,x(t),x(s),u(s)\big) \, ds, \quad t \in [t_0,t_1].\
\end{eqnarray*}
Darin fließen die Zustände $x(t_0)$, $x(t_1)$ in das Zielfunktional bzw. und der Zustand $x(t)$ in die Dynamik ein.
Deswegen müssen wir im Folgenden die Zustände der äußeren Zeit in die Betrachtungen aufnehmen. 

\newpage
Wir befassen uns also mit den Abbildungen
$f=f(s,x,x_0,x_1,u)$ und $\varphi=\varphi(t,s,\xi,x,u)$.
In Verallgemeinerung zu $V^{\mathcal{I}}_\gamma$ führen wir für $x_*(\cdot)$ die Menge $\mathcal{V}^{\mathcal{I}}_\gamma$ all derjenigen Punkte
$(t,s,\xi,x,x_0,x_1) \in \R \times \R \times \R^n \times \R^n \times \R^n \times \R^n$ ein,
für die die Ungleichungen
$$t_0 \leq s, t\leq t_1, \;\; \|\xi-x_*(t)\| \leq \gamma, \;\; \|x-x_*(s)\| \leq \gamma, \;\; \|x_0-x_*(t_0)\| \leq \gamma, \;\; \|x_1-x_*(t_1)\|$$
erfüllt sind.
Wir nehmen an, dass die Abbildungen $f(s,x,x_0,x_1,u)$ und $\varphi(t,s,\xi,x,u)$ auf $\mathcal{V}^{\mathcal{I}}_\gamma \times \R^m$
stetig in der Gesamtheit der Variablen und stetig differenzierbar bez"uglich der Zustandsvariablen $\xi,x$ und $x_0,x_1$ sind.
Daher sind für jede kompakte Menge $U_1 \subset \R^m$ die Abbildungen $\varphi$, $\varphi_\xi$ und $\varphi_x$ gleichmäßig stetig auf
$\mathcal{V}^{\mathcal{I}}_\gamma \times U_1$: 
Zu jedem $\varepsilon >0$ gibt es ein $\sigma >0$ mit 
$$\|\varphi(t,s,\xi,x,u) - \varphi(t',s',\xi',x',u')\| \leq \varepsilon$$
und
$$\|\varphi_\xi(t,s,\xi,x,u) - \varphi_\xi(t',s',\xi',x',u')\| \leq \varepsilon, 
   \quad \|\varphi_x(t,s,\xi,x,u) - \varphi_x(t',s',\xi',x',u')\| \leq \varepsilon$$
für alle $(t,s,\xi,x,u),\, (t',s',\xi',x',u') \in \mathcal{V}^{\mathcal{I}}_\gamma \times U_1$ mit $\|(t,s,\xi,x,u) -(t',s',\xi',x',u')\| \leq \sigma$.
Hier und im Folgenden bezeichnet $(t,s,\xi,x) \in \mathcal{V}^{\mathcal{I}}_\gamma$ die Einschränkung von
$(t,s,\xi,x,x_0,x_1) \in \mathcal{V}^{\mathcal{I}}_\gamma$ auf die ersten vier Komponenten. \\[2mm]
Es seien $x_*(\cdot) \in C([t_0,t_1],\R^n)$, $u_*(\cdot) \in L_\infty\big([t_0,t_1],U\big)$ und
$u_1(\cdot),...,u_d(\cdot) \in L_\infty\big([t_0,t_1],U\big)$.
Weil sämtliche Funktionen $u_*(\cdot), u_1(\cdot),...,u_d(\cdot)$ messbar und beschr"ankt sind,
sind ihre Werte fast überall in $[t_0,t_1]$ in einer kompakten Menge $U_1 \subset \R^m$ enthalten.
Wegen der gleichmäßigen Stetigkeit der Abbildungen $\varphi(t,s,\xi,x,u)$, $\varphi_\xi(t,s,\xi,x,u)$ und $\varphi_x(t,s,\xi,x,u)$
auf $\mathcal{V}^{\mathcal{I}}_\gamma \times U_1$ gibt es zu jedem $\delta >0$ eine Zerlegung $t_0=\tau_0<\tau_1 < ...<\tau_K=t_1$ von $[t_0,t_1]$
und eine Zahl $0< \sigma \leq \gamma$ derart,
dass für jedes $j=1,...,K$ die Ungleichungen
\begin{equation} \label{EigenschaftNadelvariationIGDelta}
\left. \begin{array}{l}
\|x_*(t)-x_*(\tau_j)\| \leq \sigma, \\[1mm]
\| \varphi(t,s,\xi,x,u) - \varphi(\tau_j,s,x_*(\tau_j),x,u) \| \leq \delta, \\[1mm]
\| \varphi_\xi(t,s,\xi,x,u) - \varphi_\xi(\tau_j,s,x_*(\tau_j),x,u) \| \leq \delta, \\[1mm]
\| \varphi_x(t,s,\xi,x,u) - \varphi_x(\tau_j,s,x_*(\tau_j),x,u) \| \leq \delta
\end{array}\right\}
\end{equation}
für alle $(t,s,\xi,x,u) \in \mathcal{V}^{\mathcal{I}}_\sigma \times U_1$ mit $t \in [\tau_{j-1},\tau_j]$ gelten. \\[2mm]
Die Zerlegung $t_0=\tau_0<\tau_1 < ...<\tau_K=t_1$ des Intervalls $[t_0,t_1]$ mit der Eigenschaft (\ref{EigenschaftNadelvariationIGDelta})
bezüglich der äußeren Zeitvariable $t$ bildet die Grundlage bei der Reduktion des zweidimensionalen Zeitbereichs auf das Zeitintervall $[t_0,t_1]$
der inneren Zeitvariable $s$. \\
Dieser Schritt erfolgt nun durch die Definition der Abbildungen $\psi_j$ in $t=\tau_j$:
\begin{equation} \label{EigenschaftNadelvariationIG1}
\psi_j(s,x,u) = \varphi(\tau_j,s,x_*(\tau_j),x,u) \quad \mbox{für } s \in [t_0,t_1].
\end{equation}

\begin{center}
\includegraphics[width=10cm]{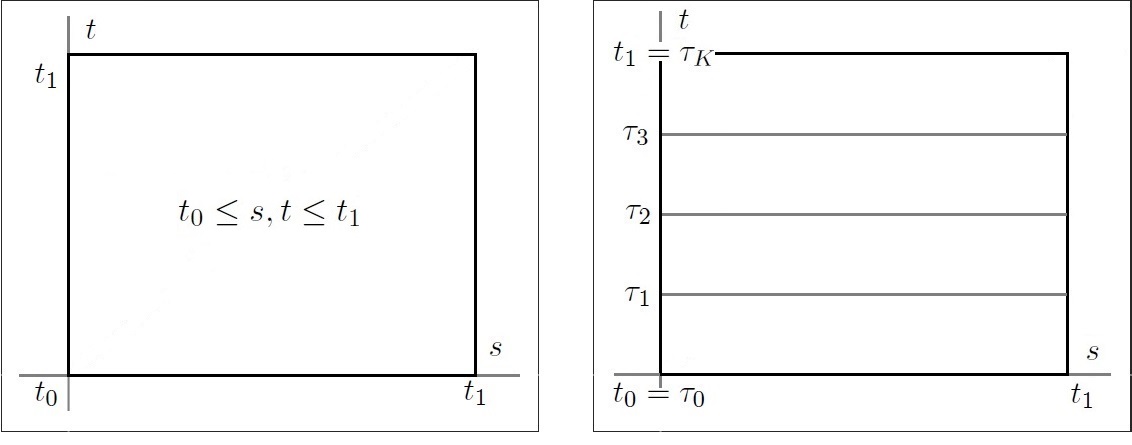}
\captionof{figure}[Zerlegung des Zeitbereiches $t_0 \leq s, t\leq t_1$]{Zerlegung des Zeitbereiches $t_0 \leq s, t\leq t_1$ in der $(s,t)-$Ebene.}
\label{AbbMengenfamilieInt1}
\end{center}

Die Abbildungen $\psi_j$ stellen die Einschränkungen von $\varphi$ über den horizontalen Linien $t=\tau_j$ innerhalb des zweidimensionalen Zeitbereichs
$\{(t,s) \in \R^2 \,|\, t_0 \leq s, t\leq t_1\}$ dar.
Wir werden nun zeigen,
dass die Abbildung $\psi_j$ eine geeignete Approximationen der Abbildung $\varphi$ über dem Rechteck
$[\tau_{j-1},\tau_j] \times [t_0,t_1]$ bildet.

\begin{center}
\includegraphics[width=5cm]{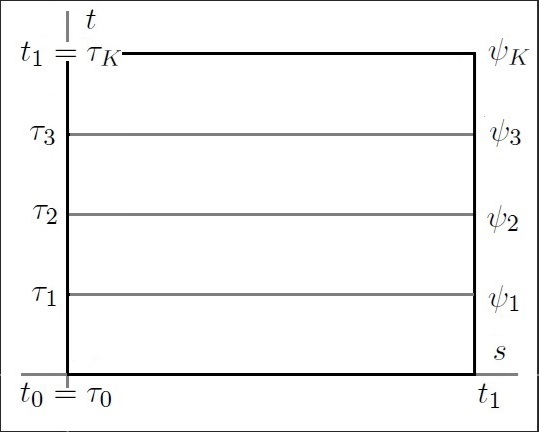}
\captionof{figure}[Abbildungen $\psi_j$ innerhalb der $(s,t)-$Ebene]{Abbildungen $\psi_j$ innerhalb der $(s,t)-$Ebene.}
\label{AbbMengenfamilieInt2}
\end{center}

Nach Wahl der Stellen $\tau_j \in [t_0,t_1]$ und der Zahl $\sigma >0$ liefert für $t \in [\tau_{j-1},\tau_j]$ die jeweilige Abbildung $\psi_j$ nach
(\ref{EigenschaftNadelvariationIGDelta}) eine $\delta-$Approximation von $\varphi$
für alle $(t,s,\xi,x,u) \in \mathcal{V}^{\mathcal{I}}_\sigma \times U_1$ mit $t \in [\tau_{j-1},\tau_j]$.
Daher ist für alle $(t,s,\xi,x,u) \in \mathcal{V}^{\mathcal{I}}_\sigma \times U_1$ die Ungleichung
\begin{eqnarray}
\lefteqn{\Big\| \varphi(t,s,\xi,x,u) - \sum_{j=1}^K \chi_{(\tau_{j-1},\tau_j]}(t) \cdot \varphi(\tau_j,s,x_*(\tau_j),x,u) \Big\|} \nonumber \\
&=& \label{EigenschaftNadelvariationIG2}
    \Big\| \varphi(t,s,\xi,x,u) - \sum_{j=1}^K \chi_{(\tau_{j-1},\tau_j]}(t) \cdot \psi_j(s,x,u) \Big\| \leq \delta
\end{eqnarray}
erfüllt. Die Stelle $t=t_0$ ist stillschweigend inbegriffen. \\[2mm]
Mit Hilfe der Abbildungen $\psi_j$ führen wir nun die gesuchten Mengenfamilien $\big\{ M_i(\alpha) \big\}$ bezüglich der Variable $s$ ein.
Dazu fassen wir die $\psi_j$ zu einer Abbildung $\psi$ zusammen:
$$\psi(s,x,u) = \big(\psi_1(s,x,u),...,\psi_K(s,x,u) \big), \quad (s,x,u) \in [t_0,t_1] \times \R^n \times \R^m.$$
Konstruieren wir nun die Mengenfamilien $\big\{ M_i(\alpha) \big\}$ bezüglich der Abbildung $\psi$,
so werden die Mengen $\big\{ M_i(\alpha) \big\}$ gleichmäßig für alle $\psi_j$ erzeugt (Abbildung \ref{AbbMengenfamilieInt3}) und besitzen
aufgrund $\|\psi_j\| \leq \|\psi\|$ gleichsam für alle $\psi_j$ die erforderlichen Eigenschaften.  

\begin{center}
\includegraphics[width=5cm]{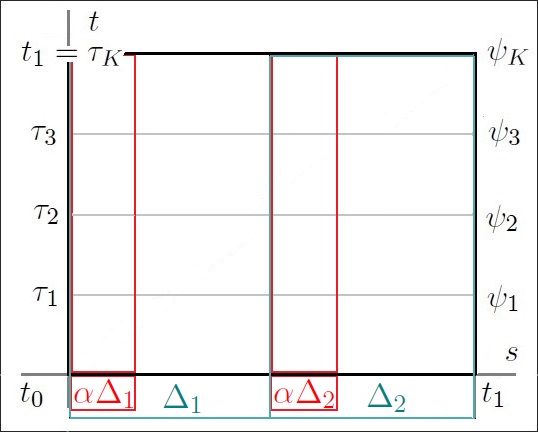}
\captionof{figure}[Gleichmäßige Einführung der Trägerfamilien bei Abbildungsfamilien]{Die gleichmäßige Einführung der Trägerfamilien für alle Abbildungen $\psi_j$.}
\label{AbbMengenfamilieInt3}
\end{center}

Es sei $\delta>0$ gegeben.
Ferner sei die Zahl $\sigma>0$ gewählt.
Außerdem sei bezüglich der Variable $t$ die Zerlegung von $[t_0,t_1]$ in $K$ Teilintervalle $[\tau_{j-1},\tau_j]$
gemäß (\ref{EigenschaftNadelvariationIGDelta}) erfolgt.
Entsprechend dem Abschnitt \ref{AbschnittDefinitionNadelvariation} betrachten wir zu den gegebenen Funktionen
$x_*(\cdot)$, $u_*(\cdot)$ und $u_1(\cdot),...,u_d(\cdot)$ nun für $i=1,...,d$ die messbaren und beschränkten Abbildungen
\begin{eqnarray*}
y_i (s) &=& \Big( \psi\big(s,x_*(s),u_i(s) \big) - \psi\big(s,x_*(s),u_*(s) \big), \\
             & & \hspace*{5mm} f\big(s,x_*(s),x_*(t_0),x_*(t_1),u_i(s) \big) - f\big(s,x_*(s),x_*(t_0),x_*(t_1),u_*(s) \big) \Big).
\end{eqnarray*}
Dann gibt es nach Lemma \ref{LemmaMengenfamilie3} über $[t_0,t_1]$ 
einparametrige Mengenfamilien $\big\{ M_i(\alpha) \big\}$ mit
\begin{eqnarray*}
&& |M_i(\alpha)| = \alpha (t_1-t_0), \qquad M_i(\alpha') \subseteq M_i(\alpha),
   \qquad M_i(\alpha) \cap M_{i'}(\alpha') = \emptyset,\\
&& \max_{t \in [t_0,t_1]}
   \bigg\| \int_{t_0}^t \big( \chi_{M_i(\alpha)}(s) - \chi_{M_i(\alpha')}(s) \big) y_i(s) \, ds
   - (\alpha-\alpha') \int_{t_0}^t y_i(s) \, ds \bigg\| \leq \frac{\delta}{2K} |\alpha-\alpha'|
\end{eqnarray*}
f"ur alle $i=1,...,d,$ $0 \leq \alpha' \leq \alpha \leq 1/d$ und $i \not= i'$.
Wegen der Setzung von $\psi$ gilt die letzte dieser Beziehungen ebenfalls für jedes $j \in \{1,...,K\}$ für die Abbildungen
\begin{eqnarray*}
y_{ij} (s) &=& \Big( \psi_j\big(s,x_*(s),u_i(s) \big) - \psi_j\big(s,x_*(s),u_*(s) \big), \\
             & & \hspace*{5mm} f\big(s,x_*(s),x_*(t_0),x_*(t_1),u_i(s) \big) - f\big(s,x_*(s),x_*(t_0),x_*(t_1),u_*(s) \big) \Big)
\end{eqnarray*}
mit den bezüglich $\psi$ angegebenen Mengenfamilien $\big\{ M_i(\alpha) \big\}$:
$$\max_{t \in [t_0,t_1]}
   \bigg\| \int_{t_0}^t \big( \chi_{M_i(\alpha)}(s) - \chi_{M_i(\alpha')}(s) \big) y_{ij}(s) \, ds
   - (\alpha-\alpha') \int_{t_0}^t y_{ij}(s) \, ds \bigg\| \leq \frac{\delta}{2K} |\alpha-\alpha'|.$$
Über einem Teilintervall $[s_0,s_1] \subseteq [t_0,t_1]$ ergibt sich außerdem
\begin{eqnarray*}
&& \max_{t \in [s_0,s_1]}
   \bigg\| \int_{s_0}^t \big( \chi_{M_i(\alpha)}(s) - \chi_{M_i(\alpha')}(s) \big) y_{ij}(s) \, ds
   - (\alpha-\alpha') \int_{s_0}^t y_{ij}(s) \, ds \bigg\| \\
&\leq& \max_{t \in [t_0,s_1]}
   \bigg\| \int_{t_0}^t \big( \chi_{M_i(\alpha)}(s) - \chi_{M_i(\alpha')}(s) \big) y_{ij}(s) \, ds
   - (\alpha-\alpha') \int_{t_0}^t y_{ij}(s) \, ds \bigg\| \\
&& + \max_{t \in [t_0,s_0]}
   \bigg\| \int_{t_0}^t \big( \chi_{M_i(\alpha)}(s) - \chi_{M_i(\alpha')}(s) \big) y_{ij}(s) \, ds
   - (\alpha-\alpha') \int_{t_0}^t y_{ij}(s) \, ds \bigg\|.
\end{eqnarray*}   
Daher ist über jedem Teilintervall $[s_0,s_1] \subseteq [t_0,t_1]$ die Ungleichung
\begin{equation} \label{EigenschaftNadelvariationIGTeilintervall}
\max_{t \in [s_0,s_1]}  \bigg\| \int_{s_0}^t \big( \chi_{M_i(\alpha)}(s) - \chi_{M_i(\alpha')}(s) \big) y_{ij}(s) \, ds
   - (\alpha-\alpha') \int_{s_0}^t y_{ij}(s) \, ds \bigg\| \leq \frac{\delta}{K} |\alpha-\alpha'|
\end{equation}
erfüllt.
Für die Zerlegung von $[t_0,t_1]$ in die Teilintervalle $[\tau_{j-1},\tau_j]$ ergibt sich schließlich
\begin{eqnarray*}
&& \max_{t \in [t_0,t_1]}
   \bigg\| \int_{t_0}^t \sum_{j=1}^K \chi_{(\tau_{j-1},\tau_j]}(t) \cdot
           \Big[\big( \chi_{M_i(\alpha)}(s) - \chi_{M_i(\alpha')}(s) \big) y_{ij}(s) - (\alpha-\alpha') \cdot y_{ij}(s)\Big] \, ds \bigg\| \\
&\leq& \sum_{j=1}^K \bigg( \max_{t \in [\tau_{j-1},\tau_j]} \bigg\|\int_{\tau_{j-1}}^t \big( \chi_{M_i(\alpha)}(s) - \chi_{M_i(\alpha')}(s) \big) y_{ij}(s) \, ds
      - (\alpha-\alpha') \int_{\tau_{j-1}}^t y_{ij}(s) \, ds \bigg\|\bigg).
\end{eqnarray*}
Mit (\ref{EigenschaftNadelvariationIGTeilintervall}) können wir die Ungleichung weiterführen und erhalten
$$\leq \sum_{j=1}^K \frac{\delta}{K} |\alpha-\alpha'| = \delta |\alpha-\alpha'|.$$
Zusammenfassend liefern die in (\ref{EigenschaftNadelvariationIG1}) definierten Abbildungen $\psi_j$ bezüglich der Variable $t$ 
eine $\delta-$Approximation der Abbildung $\varphi$ im Sinn von (\ref{EigenschaftNadelvariationIG2}).
Andererseits sind für die Abbildungen $\psi_j$ bezüglich der Variable $s$ mit den Mengenfamilien $\big\{ M_i(\alpha) \big\}$ die Beziehungen
\begin{equation} \label{EigenschaftNadelvariationIG3}
\left.\begin{array}{l}
 |M_i(\alpha)| = \alpha (t_1-t_0), \qquad M_i(\alpha') \subseteq M_i(\alpha),
   \qquad M_i(\alpha) \cap M_{i'}(\alpha') = \emptyset, \\[1mm]
 \displaystyle \max_{t \in [t_0,t_1]}
   \bigg\| \int_{t_0}^t \sum_{j=1}^K \chi_{(\tau_{j-1},\tau_j]}(s) \cdot 
           \big( \chi_{M_i(\alpha)}(s) - \chi_{M_i(\alpha')}(s) \big) y_{ij}(s) \, ds \\[1mm]
 \displaystyle \hspace*{25mm}
   - (\alpha-\alpha') \int_{t_0}^t \sum_{j=1}^K \chi_{(\tau_{j-1},\tau_j]}(s) \cdot y_{ij}(s) \, ds \bigg\|
    \leq \delta |\alpha-\alpha'|
\end{array} \right\}
\end{equation}
f"ur alle $i=1,...,d,$ $0 \leq \alpha' \leq \alpha \leq 1/d$ und $i \not= i'$ erfüllt.
Die Vorbereitungen bezüglich der Trägerfamilien $\big\{ M_i(\alpha) \big\}$ sind erbracht und wir können nun die geforderten Eigenschaften
entsprechend von Lemma \ref{LemmaEigenschaftNadelvariation1} und von Lemma \ref{LemmaEigenschaftNadelvariation2} angehen.

\newpage
Im Weiteren bezeichnen $Q^{d}$, $\Sigma(\Delta)$ und $V(\sigma)$ die Mengen aus den Abschnitten
\ref{AbschnittDefinitionNadelvariation} und \ref{AbschnittEigenschaftenNadelvariation}. 
Weil alle Steuerungen $u_*(\cdot), u_1(\cdot),...,u_d(\cdot)$ beschr"ankt sind,
sind ihre Werte f"ur fast alle $t \in [t_0,t_1]$ in einer kompakten Menge $U_1 \subset \R^m$ enthalten.
Nutzen wir dar"uber hinaus die Stetigkeit von $f$, $\varphi$, $\varphi_\xi$ und $\varphi_x$ aus,
so k"onnen wir ein $\sigma > 0$ derart angeben, dass
\begin{eqnarray}
&& \label{EigenschaftNadelvariationIG4}
   \hspace*{-10mm} \big\| \varphi(t,s,\xi,x,u) - \varphi(t,s,x_*(t),x_*(s),u) \big\| \leq \frac{\delta}{4(t_1-t_0)}, \\
&& \label{EigenschaftNadelvariationEinschubf}
   \hspace*{-10mm} \big| f(t,x,x_0,x_1,u) - f(t,x_*(t),x_*(t_0),x_*(t_1),u) \big| \leq \frac{\delta}{8(t_1-t_0)}, \\
&& \hspace*{-10mm} \big\| \varphi(t,s,\xi,x,u) - \varphi(t,s,\xi',x',u) - \varphi_\xi (t,s,x_*(t),x_*(s),u) (\xi-\xi') \nonumber \\
&& \label{EigenschaftNadelvariationIG5}
   \hspace*{15mm} - \varphi_x (t,s,x_*(t),x_*(s),u) (x-x') \big\| \leq \frac{\delta (\|\xi-\xi'\| + \|x-x'\|)}{4(t_1-t_0)}
\end{eqnarray}
für alle $(t,s,\xi,x,x_0,x_1,u) \in \mathcal{V}^{\mathcal{I}}_\sigma \times U_1$ erfüllt sind.
Ferner w"ahlen wir $\Delta \in (0,1/d]$ mit
\begin{equation} \label{EigenschaftNadelvariationIG6}
\left.\begin{array}{l}
\displaystyle \Delta \cdot \max_{ u \in U_1} \big\| \varphi_\xi(t,s,x_*(t),x_*(s),u) \big\| \leq \frac{\delta}{4(t_1-t_0)}\\[3mm]
\displaystyle \Delta \cdot \max_{ u \in U_1} \big\| \varphi_x(t,s,x_*(t),x_*(s),u) \big\| \leq \frac{\delta}{4(t_1-t_0)}
\end{array}\right\}
\quad \mbox{für fast alle } (t,s).
\end{equation}

Au"serdem kann eine Zerlegung $t_0=\tau_0<\tau_1 < ...<\tau_K=t_1$ von $[t_0,t_1]$ gewählt werden,
so dass für $j=1,...,K$ die Ungleichung
\begin{equation}
\|x_*(t)-x_*(\tau_j)\| \leq \sigma \quad \mbox{ für alle } t \in [\tau_{j-1},\tau_j]
\end{equation} 
erfüllt ist und $\varphi$ nach (\ref{EigenschaftNadelvariationIG2})  durch die Abbildungen $\psi_j$
in (\ref{EigenschaftNadelvariationIG1}) mit der Güte
\begin{equation} \label{EigenschaftNadelvariationIG7}
\Big\| \varphi(t,s,\xi,x,u) - \sum_{j=1}^K \chi_{(\tau_{j-1},\tau_j]}(t) \cdot \psi_j(s,x,u) \Big\| \leq \frac{\delta}{16(t_1-t_0)}
\end{equation}
für alle $(t,s,\xi,x,u) \in \mathcal{V}^{\mathcal{I}}_\sigma \times U_1$ approximiert wird. \\[2mm]
Wir betrachten im Weiteren f"ur $t \in [t_0,t_1]$ die Abbildung
$$\Phi_1\big(x(\cdot),\alpha\big)(t) = \int_{t_0}^t \big[\varphi\big(t,s,x(t),x(s),u_\alpha(s)\big)
                                                   - \varphi\big(t,s,x_*(t),x_*(s),u_*(s)\big)\big] \, ds$$
und den linearen Operator
\begin{eqnarray*}
\Lambda_1\big(x(\cdot),\alpha\big)(t) &=& \int_{t_0}^t \Big[\varphi_\xi\big(t,s,x_*(t),x_*(s),u_*(s)\big) \big(x(t)-x_*(t)\big) \\
&& \hspace*{20mm} + \varphi_x\big(t,s,x_*(t),x_*(s),u_*(s)\big) \big(x(s)-x_*(s)\big) \\
&& \hspace*{-5mm}+ \sum_{i=1}^{d} \alpha_i \cdot 
   \Big( \varphi\big(t,s,x_*(t),x_*(s),u_i(s)\big) - \varphi\big(t,s,x_*(t),x_*(s),u_*(s)\big)\Big)\Big] \, ds.
\end{eqnarray*}

\begin{lemma} \label{LemmaEigenschaftNadelvariationIG1}
Es existieren $\Delta \in (0,1/d]$ und $\sigma>0$ derart, dass
\begin{eqnarray}
&& \Big\| \big[\Phi_1\big(x(\cdot),\alpha\big)-\Phi_1\big(x'(\cdot),\alpha'\big)
          -\Lambda_1\big(x(\cdot),\alpha\big)+\Lambda_1\big(x'(\cdot),\alpha'\big)\big] (\cdot) \Big\|_\infty \nonumber \\
&& \label{EigenschaftNadelvariationIG8}
   \hspace*{20mm} \leq\delta \bigg( \|x(\cdot)-x'(\cdot)\|_\infty + \sum_{i=1}^d |\alpha_i - \alpha_i'| \bigg)
\end{eqnarray}
f"ur alle $\alpha, \alpha' \in \Sigma(\Delta)$ und alle $x(\cdot), x'(\cdot) \in V(\sigma)$ gilt.  
\end{lemma}
  
{\bf Beweis} Genauso wie im Beweis von Lemma \ref{LemmaEigenschaftNadelvariation1} ergeben sich die einzelnen Summanden
\begin{eqnarray*}
&& \max_{t \in [t_0,t_1]} 
   \int_{t_0}^t \Big\| \varphi \big(t,s,x(t),x(s),u_\alpha(s)\big) - \varphi \big(t,s,x'(t),x'(s),u_\alpha(s)\big) \\
& & \hspace*{30mm} - \varphi_\xi \big(t,s,x_*(t),x_*(s),u_\alpha(s)\big) \big(x(t)-x'(t)\big) \\
& & \hspace*{30mm} - \varphi_x \big(t,s,x_*(t),x_*(s),u_\alpha(s)\big) \big(x(s)-x'(s)\big) \Big\| ds \\
& & + \max_{t \in [t_0,t_1]}
   \int_{t_0}^t \Big\| \Big( \varphi_\xi \big(t,s,x_*(t),x_*(s),u_\alpha(s)\big) \\
& & \hspace*{30mm}  - \varphi_\xi \big(t,s,x_*(t),x_*(s),u_*(s)\big) \Big) \big(x(t)-x'(t)\big) \Big\| ds \\
& & + \max_{t \in [t_0,t_1]}
   \int_{t_0}^t \Big\| \Big( \varphi_x \big(t,s,x_*(t),x_*(s),u_\alpha(s)\big) \\
& & \hspace*{30mm}  - \varphi_x \big(t,s,x_*(t),x_*(s),u_*(s)\big) \Big) \big(x(s)-x'(s)\big) \Big\| ds \\
& & + \max_{t \in [t_0,t_1]}
   \int_{t_0}^t \Big\| \varphi \big(t,s,x'(t),x'(s),u_\alpha(s)\big) - \varphi \big(t,s,x'(t),x'(s),u_{\alpha'}(s)\big) \\
& & \hspace*{30mm} - \varphi \big(t,s,x_*(t),x_*(s),u_\alpha(s)\big) + \varphi \big(t,s,x_*(t),x_*(s),u_{\alpha'}(s)\big) \Big\| ds
\end{eqnarray*}

und ferner der Summand von eigentlichem Interesse
\begin{eqnarray}
&& \hspace*{-20mm} + \max_{t \in [t_0,t_1]} \bigg\| \int_{t_0}^t \varphi \big(t,s,x_*(t),x_*(s),u_\alpha(s)\big)
                   - \varphi \big(t,s,x_*(t),x_*(s),u_{\alpha'}(s)\big) \nonumber \\
&& \hspace*{-10mm} \label{EigenschaftNadelvariationIG9} - \sum_{i=1}^{d} (\alpha_i- \alpha'_i)
      \Big( \varphi \big(t,s,x_*(t),x_*(s),u_i(s)\big) - \varphi \big(t,s,x_*(t),x_*(s),u_*(s)\big) \Big) ds \bigg\|.
\end{eqnarray}
Unter Beachtung der Eigenschaften (\ref{EigenschaftNadelvariation1}) und (\ref{EigenschaftNadelvariation2})
sind die ersten vier Summanden nach (\ref{EigenschaftNadelvariationIG4})--(\ref{EigenschaftNadelvariationIG6}) kleiner oder gleich
\begin{equation} \label{EigenschaftNadelvariationIG10}
\frac{\delta}{2} \|x(\cdot)-x'(\cdot)\|_\infty + \frac{\delta}{4} \|x(\cdot)-x'(\cdot)\|_\infty
+\frac{\delta}{4} \|x(\cdot)-x'(\cdot)\|_\infty + \frac{\delta}{2} \sum_{i=1}^{d} | \alpha_i - \alpha'_i |.
\end{equation}
Den Summanden (\ref{EigenschaftNadelvariationIG9}) 
können wir mit Hilfe der Abbildungen $\psi_j$ durch folgende Ausdrücke nach oben abschätzen:
\begin{eqnarray*}
& &  \max_{t \in [t_0,t_1]} \bigg\| \int_{t_0}^t \Big[\Big(\varphi \big(t,s,x_*(t),x_*(s),u_\alpha(s)\big)
                            - \sum_{j=1}^K \chi_{(\tau_{j-1},\tau_j]}(t) \cdot \psi_j\big(s,x_*(s),u_\alpha(s)\big) \Big) \\
& & \hspace*{10mm}
       - \Big(\varphi \big(t,s,x_*(t),x_*(s),u_{\alpha'}(s)\big)
            - \sum_{j=1}^K \chi_{(\tau_{j-1},\tau_j]}(t) \cdot \psi_j\big(s,x_*(s),u_{\alpha'}(s)\big)\Big)\Big] ds \bigg\| \\
&+& \max_{t \in [t_0,t_1]} \bigg\| \int_{t_0}^t  \sum_{i=1}^{d} (\alpha_i- \alpha'_i)
      \Big( \varphi \big(t,s,x_*(t),x_*(s),u_i(s)\big) \\
& & \hspace*{50mm} - \sum_{j=1}^K \chi_{(\tau_{j-1},\tau_j]}(t) \cdot \psi_j \big(s,x_*(s),u_i(s)\big) \Big) \\
& &  - \sum_{i=1}^{d} (\alpha_i- \alpha'_i)
      \Big( \varphi \big(t,s,,x_*(t),x_*(s),u_*(s)\big)
            - \sum_{j=1}^K \chi_{(\tau_{j-1},\tau_j]}(t) \cdot \psi_j \big(s,x_*(s),u_*(s)\big) \Big) ds \bigg\| \\
&+&  \max_{t \in [t_0,t_1]} \bigg\| \int_{t_0}^t \sum_{j=1}^K \chi_{(\tau_{j-1},\tau_j]}(t) \cdot \Big[
             \psi_j \big(s,x_*(s),u_\alpha(s)\big) - \psi_j \big(s,x_*(s),u_{\alpha'}(s)\big) \\
& & \hspace*{30mm} - \sum_{i=1}^{d} (\alpha_i- \alpha'_i)
      \Big( \psi_j\big(s,x_*(s),u_i(s)\big) - \psi_j\big(s,x_*(s),u_*(s)\big) \Big)\Big] ds \bigg\|.
\end{eqnarray*}

In dieser Summe fallen nach (\ref{EigenschaftNadelvariationIG7}) die ersten beiden Summanden kleiner oder gleich
\begin{equation} \label{EigenschaftNadelvariationIG11}
\frac{2\delta}{16} \sum_{i=1}^{d} | \alpha_i - \alpha'_i |+\frac{2\delta}{16} \sum_{i=1}^{d} | \alpha_i - \alpha'_i |
  =\frac{\delta}{4} \sum_{i=1}^{d} | \alpha_i - \alpha'_i |
\end{equation}
aus (bezüglich dem ersten Summanden ist (\ref{EigenschaftNadelvariation2}) zu beachten).
Der dritte Summand wurde bereits in den vorbereitenden Betrachtungen diskutiert.
Das Ergebnis ist die Version von Lemma \ref{LemmaEigenschaftNadelvariation1} für die Abbildungen $\psi_j$
in Form der Beziehung (\ref{EigenschaftNadelvariationIG3}).
Demnach lassen sich die Mengenfamilien $\{M_i(\alpha)\}$ so angeben,
dass der dritte Summand nicht größer als
\begin{equation} \label{EigenschaftNadelvariationIG12}
\frac{\delta}{4} \sum_{i=1}^{d} | \alpha_i - \alpha'_i |
\end{equation}
ausfällt.
Zusammenfassend zeigen die Abschätzungen (\ref{EigenschaftNadelvariationIG10})--(\ref{EigenschaftNadelvariationIG12}),
dass die linke Seite in (\ref{EigenschaftNadelvariationIG8})
nicht größer als
$$\delta \bigg(\|x(\cdot)-x'(\cdot)\|_\infty + \sum_{i=1}^{d} | \alpha_i - \alpha'_i |\bigg)$$
ist.
Das Lemma \ref{LemmaEigenschaftNadelvariationIG1} ist damit bewiesen. \hfill $\blacksquare$

\newpage
Nachdem wir die aufwendige Untersuchung der Dynamik abgeschlossen haben widmen wir uns dem Zielfunktional
$$J\big(x(\cdot),u(\cdot)\big)  = \int_{t_0}^{t_1} f\big(s,x(s),x(t_0),x(t_1),u(s)\big) \, ds,$$
in welchem die äußeren Zustände $x(t_0)$ und $x(t_1)$ zusätzlich auftreten.
Mit dem Funktional $J$ verbinden wir die Abbildungen
\begin{eqnarray*}
\Phi_2\big(x(\cdot),\alpha\big) &=&  \int_{t_0}^{t_1} \big[f\big(t,x(t),x(t_0),x(t_1),u_\alpha(t)\big)-f\big(t,x(t),x(t_0),x(t_1),u_*(t)\big)\big] \, dt, \\
\Lambda_2\big(x(\cdot),\alpha\big) &=& \sum_{i=1}^{d} \int_{t_0}^{t_1} \alpha_i \cdot \Big( f\big(t,x(t),x(t_0),x(t_1),u_i(t)\big) \\
                                   & &\hspace*{35mm} - f\big(t,x(t),x(t_0),x(t_1),u_*(t)\big)\Big) \, dt.
\end{eqnarray*}

\begin{lemma} \label{LemmaEigenschaftNadelvariationIG2}
Es existieren $\Delta \in (0,1/d]$ und $\sigma>0$ derart, dass
\begin{equation} \label{EigenschaftNadelvariationIG13}
\Phi_2\big(x(\cdot),\alpha\big) - \Lambda_2\big(x(\cdot),\alpha\big) \leq \delta \sum_{i=1}^d \alpha_i
\end{equation}
f"ur alle $\alpha \in \Sigma(\Delta)$ und alle $x(\cdot) \in V(\sigma)$ gilt.
\end{lemma}

{\bf Beweis} Beachten wir (\ref{EigenschaftNadelvariation1}),
so ist die linke Seite in (\ref{EigenschaftNadelvariationIG13}) gleich
\begin{eqnarray*}
& &  \sum_{i=1}^{d} \int_{t_0}^{t_1} \bigg[ \chi_{M_i(\alpha_i)}(t) \cdot
                \Big( f \big(t,x(t),x(t_0),x(t_1),u_i(t)\big) - f\big(t,x(t),x(t_0),x(t_1),u_*(t)\big)\Big) \\
& & \hspace*{25mm} - \alpha_i \cdot \Big( f\big(t,x(t),x(t_0),x(t_1),u_i(t)\big) - f\big(t,x(t),x(t_0),x(t_1),u_*(t)\big)\Big) \bigg] dt.
\end{eqnarray*}
Im Vergleich zum Beweis von Lemma \ref{LemmaEigenschaftNadelvariation2}
müssen im weiteren Vorgehen lediglich die äußeren Zustände $x(t_0)$ und $x(t_1)$ ergänzt werden.
Ersetzen wir die entsprechenden Verweise auf (\ref{Nadelvariation2}) durch (\ref{EigenschaftNadelvariationIG3})
und auf (\ref{EigenschaftNadelvariation4}) durch (\ref{EigenschaftNadelvariationEinschubf}),
so ergibt sich die Gültigkeit von Lemma \ref{LemmaEigenschaftNadelvariationIG2}
auf die gleiche Weise wie Lemma \ref{LemmaEigenschaftNadelvariation2}. \hfill $\blacksquare$

%% file: VI-Extremalprinzip.tex
\section{Zur Theorie der Extremalaufgaben} \label{KapitelExtremalaufgaben}
Zur Ermittlung notwendiger Optimalit"atsbedingungen f"ur Aufgaben der Klassischen Variationsrechnung und anschlie"send f"ur Probleme der Optimalen Steuerung
zeichnete sich zu Beginn des letzten Jahrhunderts ab,
dass eine allgemeinere Methodik erarbeitet werden musste (vgl. Plail \cite{Plail}).
Daraus entwickelte sich das Verst"andnis von abstrakten Optimierungsaufgaben unter Nebenbedingungen in Funktionenr"aumen
und in Folge dessen die Herausforderung des Nachweises des fundamentalen Prinzips von Lagrange in einem unendlich-dimensionalen Rahmen.
Die ben"otigten Werkzeuge lieferte die zu dieser Zeit parallele Entwicklung der Funktionalanalysis. \\
Grundlegende Beitr"age zu einer Verallgemeinerung lieferten die Chicagoer Schule um Bolza (ein Sch"uler von Karl Weierstra"s), Bliss und McShane.
Eine allgemeine Entwicklung einer Theorie der Extremalaufgaben erfolgte ma"sgeblich in den 1960er- und 1970er-Jahren.
Diesbez"uglich geben wir die Arbeiten von
Dubovickii \& Milyutin \cite{DuboMil}, Girsanov \cite{Girsanov}, Halkin \cite{Halkin2}, Halkin \& Neustadt \cite{HalkinNeustadt},
Ioffe \& Tichomirov \cite{Ioffe}, Kurcyusz \cite{Kurcyusz} und Neustadt \cite{Neustadt} an. \\[2mm]
In der weiteren Darstellung seien $\mathscr{X},\mathscr{Z}$ Banachr"aume, $\mathscr{U}$ eine nichtleere Menge,
$\mathscr{K} \subseteq \mathscr{Z}$ ein abgeschlossener konvexer Kegel mit Spitze im Ursprung $0 \in \mathscr{Z}$ und mit nichtleerem Inneren.
Ferner seien $J$ ein Funktional und $\mathscr{F}$, $G$ bzw. $G_j$, $j=1,...,l$, gewisse Abbildungen.
Die Aufgabenklassen, die wir im Folgenden behandeln werden, k"onnen wir damit wie folgt katalogisieren:
\begin{enumerate}
\item[(1)] Glatte Aufgaben mit Gleichungsnebenbedingungen:
           $$J(x) \to \inf; \quad \mathscr{F}(x)=0, \quad x \in \mathscr{X}.$$
\item[(2)] Glatte Aufgaben mit Gleichungs- und Ungleichungsnebenbedingungen:
           $$J(x) \to \inf; \quad \mathscr{F}(x)=0, \quad G(x) \in \mathscr{K}, \quad x \in \mathscr{X}.$$
\item[(3)] Schwaches lokales Minimum unter Gleichungs- und Ungleichungsnebenbedingungen:
           $$J(x,u) \to \inf; \quad \mathscr{F}(x,u)=0, \quad G(x) \in \mathscr{K}, \quad x \in \mathscr{X},\; u \in \mathscr{U}.$$
           Dabei ist $\mathscr{U}$ eine nichtleere konvexe Teilmenge eines normierten Raumes. 
\item[(4)] Starkes lokales Minimum unter Gleichungs- und Ungleichungsnebenbedingungen:
           $$J(x,u) \to \inf; \quad \mathscr{F}(x,u)=0, \quad G_j(x) \leq 0, \; j=1,...,l, \quad x \in \mathscr{X},\; u \in \mathscr{U}.$$
\end{enumerate}
Den Kern bei der Herleitung des Extremalprinzips,
d.\,h. beim Nachweis der G"ultigkeit des Lagrangeschen Prinzips,
bildet in den einzelnen Aufgaben (2)--(4) das Dubovickii-Milyutin-Schema.
Die Vorgehensweise l"asst sich dabei wie folgt abk"urzend skizzieren:
\begin{enumerate}
\item[(i)] Zu einem lokalen Minimum $x_* \in \mathscr{X}$ bezeichne $\mathscr{C}$ die Menge derjenigen $(\eta_0,y,z)$ mit
           $$\eta_0 > J'(x_*)x, \quad y = \mathscr{F}'(x_*)x, \quad z = G(x_*) + G'(x_*)x.$$
           Die Menge $\mathscr{C}$ ist konvex und besitzt ein nichtleeres Inneres.
\item[(ii)] W"urde nun der Ursprung dem Inneren der Menge $\mathscr{C}$ angeh"oren,
            so gibt es ein $\overline{x} \in \mathscr{X}$ und ferner nach dem Satz von Ljusternik Elemente
            $x_\varepsilon=x_*+ \varepsilon \overline{x}+r(\varepsilon)$ mit
            $$J(x_\varepsilon)<J(x_*), \quad \mathscr{F}(x_\varepsilon)=0, \quad G(x_\varepsilon) \in \mathscr{K},
              \quad x_\varepsilon \to x_* \mbox{ f"ur } \varepsilon \to 0;$$
            im Widerspruch zur lokalen Optimalit"at von $x_*$.
\item[(iii)] Nach dem Trennungssatz f"ur konvexe Mengen lassen sich die Menge $\mathscr{C}$ und der Ursprung trennen.
             Demnach existiert ein nichttriviales stetiges lineares Funktional (nichttriviale Lagrangesche Multiplikatoren),
             mit dem (denen) das Lagrangesche Prinzip gilt.
\end{enumerate}
Im Gegensatz zu den Aufgaben (1) und (2) unterscheiden wir in den Aufgaben (3) und (4) die Gr"o"sen $x$ und $u$.
Dies erm"oglicht einerseits die Unterscheidung zwischen schwacher und starker lokaler Optimalit"at.
Andererseits kann in Steuerungsproblemen explizit dem eigenst"andigen Charakter einer Steuerung Rechnung getragen werden.
Insbesondere gibt es in der Aufgabe (4) keine Annahmen "uber die Struktur der Menge $\mathscr{U}$.
Die Erweiterung des Dubovickii-Milyutin-Schema f"ur starke lokale Minimalstellen ist Ioffe \& Tichomirov \cite{Ioffe} entnommen. \\[2mm]
In den Aufgaben (2)--(4) behandeln wir die Ungleichungsrestriktionen in den Varianten $G(x) \in \mathscr{K}$ und $G_j(x) \leq 0$ f"ur $j=1,...,l$.
Im einfachen Fall der Ungleichungen $x_j(t) \leq 0$ f"ur alle $t \in [t_0,t_1]$ und $j=1,...,l$ bezeichnet $\mathscr{K}$ den Kegel
$$\mathscr{K}= \{ x(\cdot) \in C([t_0,t_1], \R^l) \,|\, x(t) \preceq 0 \Leftrightarrow x_j(t) \leq 0 \mbox{ f"ur alle } t \in [t_0,t_1], \; j=1,...,l\}.$$
Andererseits k"onnen wir die Ungleichungen $x_j(t) \leq 0$ f"ur alle $t \in [t_0,t_1]$ mit Hilfe der regul"ar lokalkonvexen Funktionen
$$G_j\big(x(\cdot)\big) = \max_{t \in [t_0,t_1]} x(t)$$
behandeln.
Im Rahmen der Extremalprinzipien ergeben sich in beiden Varianten entsprechende komplement"are Schlupfbedingungen.
Die Herleitungen der Extremalprinzipien setzen umfassende Kenntnisse der vorhergehenden Teile des Anhanges voraus.
Insbesondere stellt die abstrakte Abbildung in der Voraussetzung (C) von Theorem \ref{SatzExtremalprinzipStark} die konstruierte mehrfache
Nadelvariation mit den entsprechenden Tr"agermengen im Anhang \ref{AnhangNV} dar.

%% file: VI-1-GlattesExtremalprinzip.tex
\subsection{Glatte Aufgaben mit Gleichungsnebenbedingungen}
In diesem Abschnitt seien $\mathscr{X}$, $\mathscr{Y}$ Banachr"aume, $J$ ein Funktional auf $\mathscr{X}$,
$\mathscr{F}$ eine Abbildung des Raumes $\mathscr{X}$ in den Raum $\mathscr{Y}$. \\[2mm]
Wir betrachten im Weiteren die Extremalaufgabe
\begin{equation}\label{EAAnhang1}
J(x) \to \inf; \quad \mathscr{F}(x)=0, \quad x \in \mathscr{X}.
\end{equation}
Der Punkt $x \in \mathscr{X}$ ist ein zul"assiges Element der Aufgabe (\ref{EAAnhang1}),
falls die Nebenbedingung $\mathscr{F}(x)=0$ erf"ullt ist.
Ein zul"assiges Element $x_* \in \mathscr{X}$ hei"st ein lokales Minimum der Extremalaufgabe (\ref{EAAnhang1}),
wenn ein $\varepsilon >0$ derart existiert,
dass f"ur alle zul"assigen Punkte $x \in \mathscr{X}$ mit $\|x-x_*\|_{\mathscr{X}} \leq \varepsilon$ die Ungleichung $J(x_*) \leq J(x)$ gilt. \\[2mm]
Auf $\mathscr{X} \times \R \times \mathscr{Y}^*$ definieren wir zur Aufgabe (\ref{EAAnhang1}) die Lagrange-Funktion
$$\mathscr{L}(x,\lambda_0,y^*) = \lambda_0 J(x) + \langle y^*, \mathscr{F}(x) \rangle.$$

\begin{theorem}[Extremalprinzip] \label{SatzExtremalprinzipGlatt}
\index{Extremalprinzip!lokal@-- lokales Minimum}
Sei $x_*$ ein zul"assiges Element der Aufgabe (\ref{EAAnhang1}).
\begin{enumerate}
\item[(A)] Wir nehmen an, dass
           \begin{enumerate}
           \item[(A$_1$)] die Funktion $J(x)$ im Punkt $x_*$ Fr\'echet-differenzierbar ist;
           \item[(A$_2$)] die Abbildung $\mathscr{F}(x)$ im Punkt $x_*$ stetig Fr\'echet-differenzierbar ist.
           \end{enumerate}
\item[(B)] Weiterhin sei ${\rm Im\,}\mathscr{F}'(x_*)$ abgeschlossen.
\end{enumerate}
Ist $x_*$ eine lokale Minimalstelle der Aufgabe (\ref{EAAnhang1}),
dann existieren nicht gleichzeitig verschwindende Lagrangesche Multiplikatoren $\lambda_0 \geq 0$ und $y^* \in \mathscr{Y}^*$ derart,
dass die Lagrange-Funktion bez"uglich $x$ in $x_*$ einen station"aren Punkt besitzt, d.\,h.
\begin{equation}\label{SatzEPglatt1}
0 = \mathscr{L}_x(x_*,\lambda_0,y^*).
\end{equation}
Gilt au"serdem ${\rm Im\,}\mathscr{F}'(x*)=\mathscr{Y}$, so ist $\lambda_0 \not=0$ und man kann $\lambda_0=1$ setzen.
\end{theorem}

{\bf Beweis im entarteten Fall:}
Es sei ${\rm Im\,}\mathscr{F}'(x*)$ ein echter abgeschlossener Unterraum von $\mathscr{Y}$.
Dann existiert nach Folgerung \ref{FolgerungAnnulator} ein nichttriviales Funktional $y^* \in \mathscr{Y}^*$ mit
$$0=\langle y^*, \mathscr{F}'(x_*)x \rangle$$
f"ur alle $x \in \mathscr{X}$.
Setzen wir au"serdem $\lambda_0=0$,
dann gilt die Beziehung (\ref{SatzEPglatt1}) des Extremalprinzips \ref{SatzExtremalprinzipGlatt}. \hfill $\blacksquare$ \\[2mm]
{\bf Beweis im regul"aren Fall:}
Im regul"aren Fall gen"ugt die Abbildung $\mathscr{F}(x)$ im Punkt $x_*$ den Voraussetzungen
des Satzes von Ljusternik (Theorem \ref{SatzLjusternik}).
Demzufolge stimmt im Punkt $x_*$ der lokale Tangentialkegel an die Menge $M = \{ x \,|\, \mathscr{F}(x)=0 \}$
mit dem Kern des Operators $\mathscr{F}'(x_*)$ "uberein.
F"ur $x_0 \in {\rm Ker\,}\mathscr{F}'(x_*)$ existiert daher eine Variation
$$x(\varepsilon) = x_* + \varepsilon x_0 + r(\varepsilon), \qquad
  \lim_{\varepsilon \to 0} \frac{\|r(\varepsilon)\|_{\mathscr{X}} }{\varepsilon}=0,$$
derart, dass f"ur alle $|\varepsilon| \leq \varepsilon_0$, $\varepsilon_0>0$, die Gleichung $\mathscr{F}(x_\varepsilon)=0$ gilt.
Dabei d"urfen wir das Intervall $[0,\varepsilon_0]$ in der Definition des lokalen Tangentialkegels auf $[-\varepsilon_0,\varepsilon_0]$ 
erweitern,
da ${\rm Ker\,}\mathscr{F}'(x_*)$ ein linearer Teilraum ist. \\[2mm]
Die Funktion $\varphi(\varepsilon)=J\big(x(\varepsilon)\big)$ besitzt in $\varepsilon=0$ ein lokales Minimum.
Daher gilt
$$\varphi'(0)=J'(x_*)x_0=0.$$
Diese Beziehung muss f"ur alle $x \in {\rm Ker\,}\mathscr{F}'(x_*)$ gelten.
Deswegen ist
$$J'(x_*) \in \big({\rm Ker\,}\mathscr{F}'(x_*)\big)^\perp.$$
Nach dem Satz vom abgeschlossenen Bild (Theorem \ref{SatzAbgeschlossenesBild}) existiert ein $y^* \in \mathscr{Y}^*$ mit
$$0 = J'(x_*) +{\mathscr{F}'}^*(x_*)y^*.$$
Damit ist der regul"are Fall des Extremalprinzips \ref{SatzExtremalprinzipGlatt} mit $\lambda_0=1$ gezeigt. \hfill $\blacksquare$ 


\subsection{Glatte Aufgaben mit Ungleichungsnebenbedingungen}
Es seien $\mathscr{X}$, $\mathscr{Y}$, $\mathscr{Z}$ Banachr"aume, $J$ ein Funktional auf $\mathscr{X}$,
$\mathscr{F}$ eine Abbildung des Raumes $\mathscr{X}$ in den Raum $\mathscr{Y}$
und $G$ eine Abbildung des Raumes $\mathscr{X}$ in den Raum $\mathscr{Z}$.
Weiterhin sei $\mathscr{K} \subseteq \mathscr{Z}$ ein abgeschlossener, konvexer Kegel mit Spitze im Ursprung $0 \in Z$ und mit nichtleerem Inneren.
Au"serdem bezeichnen wir mit $\mathscr{K}^*$ den dualen Kegel
$$\mathscr{K}^*=\{z^* \in \mathscr{Z}^* \,|\, \langle z^*,z \rangle \leq 0 \mbox{ f"ur alle } z \in \mathscr{K} \}.$$

Wir betrachten im Weiteren die Extremalaufgabe
\begin{equation}\label{EAAnhang11}
J(x) \to \inf; \quad \mathscr{F}(x)=0, \quad G(x) \in \mathscr{K}, \quad x \in \mathscr{X}.
\end{equation}
Der Punkt $x \in \mathscr{X}$ ist ein zul"assiges Element der Aufgabe (\ref{EAAnhang11}),
falls die Nebenbedingungen $\mathscr{F}(x)=0$ und $G(x) \in \mathscr{K}$ erf"ullt sind.
Ein zul"assiges Element $x_* \in \mathscr{X}$ hei"st ein lokales Minimum der Extremalaufgabe (\ref{EAAnhang11}),
wenn ein $\varepsilon >0$ derart existiert,
dass f"ur alle zul"assigen Punkte $x \in \mathscr{X}$ mit $\|x-x_*\|_{\mathscr{X}} \leq \varepsilon$ die Ungleichung $J(x_*) \leq J(x)$ gilt. \\[2mm]
Auf $\mathscr{X} \times \R \times \mathscr{Y}^* \times \mathscr{Z}^*$ definieren wir zur Aufgabe (\ref{EAAnhang11}) die Lagrange-Funktion
$$\mathscr{L}(x,\lambda_0,y^*, z^*) = \lambda_0 J(x) + \langle y^*, \mathscr{F}(x) \rangle +  \langle z^*, G(x) \rangle.$$

\begin{theorem}[Extremalprinzip] \label{SatzExtremalprinzipGlatt2}
\index{Extremalprinzip!lokal@-- lokales Minimum}
Sei $x_*$ ein zul"assiges Element der Aufgabe (\ref{EAAnhang11}).
\begin{enumerate}
\item[(A)] Wir nehmen an, dass
           \begin{enumerate}
           \item[(A$_1$)] die Funktion $J(x)$ im Punkt $x_*$ Fr\'echet-differenzierbar ist;
           \item[(A$_2$)] die Abbildung $\mathscr{F}(x)$ im Punkt $x_*$ stetig Fr\'echet-differenzierbar ist;
           \item[(A$_3$)] die Abbildung $G(x)$ im Punkt $x_*$ Fr\'echet-differenzierbar ist.
           \end{enumerate}
\item[(B)] Weiterhin sei ${\rm Im\,}\mathscr{F}'(x_*)$ abgeschlossen.
\end{enumerate}
Ist $x_*$ eine lokale Minimalstelle der Aufgabe (\ref{EAAnhang11}),
dann existieren nicht gleichzeitig verschwindende Lagrangesche Multiplikatoren
$\lambda_0 \geq 0$, $y^* \in \mathscr{Y}^*$ und $z^* \in \mathscr{Z}^*$ derart,
dass folgende Bedingungen gelten:
\begin{enumerate}
\item[(a)] Die Lagrange-Funktion besitzt bez"uglich $x$ in $x_*$ einen station"aren Punkt, d.\,h.
           \begin{equation}\label{SatzEPglatt21}
           0 = \mathscr{L}_x(x_*,\lambda_0,y^*,z^*);
           \end{equation}   
\item[(b)] Die komplement"aren Schlupfbedingungen gelten, d.\,h.
           \begin{equation}\label{SatzEPglatt22}
           0 = \langle z^*, G(x_*) \rangle, \qquad z^* \in \mathscr{K}^*.
           \end{equation}
\end{enumerate}
\end{theorem}

{\bf Beweis im entarteten Fall:}
Es sei ${\rm Im\,}\mathscr{F}'(x*)$ ein echter abgeschlossener Unterraum von $\mathscr{Y}$.
Dann existiert nach Folgerung \ref{FolgerungAnnulator} ein nichttriviales Funktional $y^* \in \mathscr{Y}^*$ mit
$$0=\langle y^*, \mathscr{F}'(x_*)x \rangle$$
f"ur alle $x \in \mathscr{X}$.
Setzen wir au"serdem $\lambda_0=0$ und $z^*=0$,
dann gelten die Beziehungen (\ref{SatzEPglatt21}) und (\ref{SatzEPglatt22})
des Extremalprinzips \ref{SatzExtremalprinzipGlatt2}. \hfill $\blacksquare$ \\[2mm]
{\bf Beweis im regul"aren Fall:} Wir nehmen $J(x_*)=0$ an.
Im Weiteren sei $\mathscr{C}$ die Menge derjenigen
$(\eta_0,y,z) \in \R \times \mathscr{Y} \times \mathscr{Z}$ mit folgender Eigenschaft:
Zu jedem Element existieren ein $x \in \mathscr{X}$ und ein $\eta \in {\rm int\,}\mathscr{K}$ mit
$$\eta_0 > J'(x_*)x, \quad y = \mathscr{F}'(x_*)x, \quad z = G(x_*) + G'(x_*)x - \eta.$$
Zum Beweis des Extremalprinzips \ref{SatzExtremalprinzipGlatt2} gen"ugt es im regul"aren Fall die Beziehungen
$${\rm int\,}\mathscr{C} \not= \emptyset, \qquad 0 \not\in {\rm int\,}\mathscr{C}$$
nachzuweisen.
Die Menge $\mathscr{C}$ ist offensichtlich konvex.
Wenn also diese Beziehungen gelten,
dann folgt die Existenz eines nichttrivialen Funktionals $(\lambda_0,y^*,z^*) \in \R \times \mathscr{Y}^* \times \mathscr{Z}^*$,
das die Menge ${\rm int\,}\mathscr{C}$ vom Ursprung trennt,
d.\,h. ein solches Funktional, mit dem die Ungleichung
$$0 \leq \lambda_0 \eta_0 + \langle y^*,y \rangle + \langle z^*,z \rangle$$
f"ur alle $(\eta_0,y,z) \in\mathscr{C}$ gilt.
Wir erhalten daraus, dass f"ur alle $x \in \mathscr{X}$ und alle $\eta \in \mathscr{K}$ die Beziehung
$$\langle z^*,\eta \rangle \leq \lambda_0 J'(x_*)x + \langle y^*, \mathscr{F}'(x_*)x \rangle + \langle z^*,G(x_*)+G'(x_*)x \rangle$$
erf"ullt ist.
Betrachten wir diese Ungleichung zun"achst f"ur $x=0$,
dann folgt
$$\langle z^*,\eta \rangle \leq \langle z^*,G(x_*) \rangle$$
f"ur alle $\eta \in {\rm int\,}\mathscr{K}$.
Da $\mathscr{K}$ ein abgeschlossener Kegel mit Spitze in Null ist, ergeben sich daraus die Beziehungen
$\langle z^*,z \rangle \leq 0$ f"ur alle $z \in \mathscr{K}$ und $\langle z^*,G(x_*) \rangle = 0$.
Damit ist (\ref{SatzEPglatt22}) gezeigt.
In Folge dessen ergibt sich f"ur alle $x \in \mathscr{X}$ die Ungleichung
$$0 \leq \lambda_0 J'(x_*)x + \langle y^*, \mathscr{F}'(x_*)x \rangle + \langle z^*,G'(x_*)x \rangle = \mathscr{L}_x(x_*,\lambda_0,y^*,z^*), x \rangle,$$
aus welcher sich die Beziehung (\ref{SatzEPglatt21}) ableitet. \\[2mm]
Wir zeigen nun ${\rm int\,}\mathscr{C} \not= \emptyset$:
Es sei $S=\{x \in \mathscr{X}\,|\, \|x\|_{\mathscr{X}} <1\}$ die offene Einheitskugel des Raumes $\mathscr{X}$.
Dann ist die Menge $\mathscr{F}'(x_*)S$ nach dem Satz von der offenen Abbildung offen.
Ferner besitzt der Kegel $\mathscr{K}$ ein nichtleeres Inneres.
Schlie"slich sei $c_0=\|J_x(x_*)\|$ und $\R_{c_0} = \{ a \in \R \,|\, a > c_0\}$.
Dann besitzt die Menge $\mathscr{C}_0= \R_{c_0} \times \mathscr{F}'(x_*)S \times (G(x_*)-{\rm int\,}\mathscr{K})$
ein nichtleeres Inneres und es gilt $\mathscr{C}_0 \subseteq \mathscr{C}$.
Demzufolge ist ${\rm int\,}\mathscr{C} \not= \emptyset$. \\[2mm]
Wir nehmen nun an, es ist $0 \in {\rm int\,}\mathscr{C}$.
Dann existieren Elemente $\overline{x} \in \mathscr{X}$, $\overline{\eta} \in {\rm int\,}\mathscr{K}$ und eine Zahl $c>0$ mit
$$-c > J'(x_*)\overline{x}, \quad 0 = \mathscr{F}'(x_*)\overline{x}, \quad 0 = G(x_*)+ G'(x_*)\overline{x} - \overline{\eta}.$$
Da die Abbildung $\mathscr{F}(x)$ im Punkt $x_*$ regul"ar ist und $\overline{x}$ dem Kern des Operators $\mathscr{F}'(x_*)$ angeh"ort,
existieren nach dem Satzes von Ljusternik (Theorem \ref{SatzLjusternik}) eine Zahl $\varepsilon_0>0$ und eine Abbildung
des Intervalls $[0,\varepsilon_0]$ in $\mathscr{X}$ derart,
dass
$$\varepsilon \to x_\varepsilon = x_* + \varepsilon \overline{x} + r(\varepsilon), \qquad \lim_{\varepsilon \to 0} \frac{\|r(\varepsilon)\|_{\mathscr{X}}}{\varepsilon}=0$$
gelten und f"ur alle $\varepsilon \in [0,\varepsilon]$ die Gleichung $\mathscr{F}(x_\varepsilon)=0$ erf"ullt ist. \\
Da die Abbildung $G$ im Punkt $x_*$ Fr\'echet-differenzierbar ist,
ergibt sich
$$G(x_\varepsilon) = G(x_*) + \varepsilon G'(x_*)\overline{x} + r(\varepsilon).$$
Wir beachten,
dass $\mathscr{K}$ konvex ist,
sowie $G(x_*) \in \mathscr{K}$ und $G(x_*)+ G'(x_*)\overline{x} = \overline{\eta} \in {\rm int\,}\mathscr{K}$ gelten.
Damit ergibt sich
$$G(x_\varepsilon) = (1-\varepsilon) G(x_*) + \varepsilon [G(x_*)+ G'(x_*)\overline{x}] + r(\varepsilon) 
                   = (1-\varepsilon) G(x_*) + \varepsilon [\overline{\eta} + r(\varepsilon) / \varepsilon].$$
Mit einem hinreichend kleinen $\varepsilon_0$ erhalten wir daraus
$G(x_\varepsilon) \in {\rm int\,}\mathscr{K}$ f"ur alle $\varepsilon \in [0,\varepsilon_0]$. \\
Da das Funktional $J$ im Punkt $x_*$ Fr\'echet-differenzierbar ist, gilt
$$J(x_\varepsilon) = J(x_*) + \varepsilon J'(x_*)\overline{x} + o(\varepsilon) < J(x_*) - \varepsilon c + o(\varepsilon).$$
Die Beziehungen zeigen,
dass f"ur hinreichend kleine $\varepsilon >0$ die Elemente $x_\varepsilon$ zul"assig in der Aufgabe (\ref{EAAnhang11}) sind.
Au"serdem ist $J(x_\varepsilon) < J(x_*)$ f"ur hinreichend kleine $\varepsilon >0$.
Wegen $x_\varepsilon \to x_*$ f"ur $\varepsilon \to 0$ bedeutet dies,
dass der Punkt $x_*$ im Widerspruch zur Voraussetzung kein schwaches lokales Minimum sein k"onnte.
Daher war die Annahme $0 \in {\rm int\,}\mathscr{C}$ falsch. \hfill $\blacksquare$

%% file: VI-2-SchwachesExtremalprinzip.tex
\subsection{Ein Extremalprinzip f\"ur ein schwaches lokales Minimum} \label{AbschnittEPschwach}
In der folgenden Darstellung des Extremalprinzips f"ur ein schwaches lokales Minimum seien $\mathscr{X}$, $\mathscr{Y}$, $\mathscr{Z}$ Banachr"aume,
$\mathscr{W}$ ein normierter Raum, $\mathscr{U}$ eine konvexe Teilmenge des Raumes $\mathscr{W}$,
$\mathscr{K} \subseteq \mathscr{Z}$ ein abgeschlossener konvexer Kegel mit Spitze im Ursprung und es sei ${\rm int\,}\mathscr{K} \not= \emptyset$.
Weiterhin seien $J$ ein Funktional auf $\mathscr{X} \times \mathscr{U}$,
$\mathscr{F}$ eine Abbildung des Produktes $\mathscr{X} \times \mathscr{U}$ in den Raum $\mathscr{Y}$ und $G:\mathscr{X} \to \mathscr{Z}$. \\[2mm]
Unter diesen Angaben betrachten wir in diesem Abschnitt die Extremalaufgabe
\begin{equation}\label{EAAnhang2}
J(x,u) \to \inf; \quad \mathscr{F}(x,u)=0, \quad G(x) \in \mathscr{K}, \quad x \in \mathscr{X},\; u \in \mathscr{U},\; \mathscr{U} \mbox{ konvex}.
\end{equation}
Der Punkt $(x,u)$ ist ein zul"assiges Element der Aufgabe (\ref{EAAnhang2}),
falls s"amtliche Nebenbedingungen erf"ullt sind.
Der Punkt $(x_*,u_*)$ hei"st ein schwaches lokales Minimum\index{Minimum, schwaches lokales!Extremal@-- Extremalaufgabe}
der Extremalaufgabe (\ref{EAAnhang2}), wenn ein $\varepsilon >0$ derart existiert,
dass f"ur alle zul"assigen Paare $(x,u)$ mit $\|(x,u)-(x_*,u_*)\|_{\mathscr{X} \times \mathscr{W}} \leq \varepsilon$
die Ungleichung $J(x_*,u_*) \leq J(x,u)$ gilt. \\[2mm]
Auf $\mathscr{X} \times \mathscr{W} \times \R \times \mathscr{Y}^* \times \mathscr{Z}^*$ definieren wir zur Aufgabe (\ref{EAAnhang2}) die
Lagrange-Funktion
$$\mathscr{L}(x,u,\lambda_0,y^*,z^*)
  = \lambda_0 J(x,u) + \langle y^*, \mathscr{F}(x,u) \rangle + \langle z^*, G(x) \rangle.$$
Au"serdem bezeichnet $\mathscr{K}^*$ den dualen Kegel
$$\mathscr{K}^*=\{z^* \in \mathscr{Z}^* \,|\, \langle z^*,z \rangle \leq 0 \mbox{ f"ur alle } z \in \mathscr{K} \}.$$

\begin{theorem}[Extremalprinzip] \label{SatzExtremalprinzipSchwach}
\index{Extremalprinzip!Schwach@-- schwaches lokales Minimum}
Sei $(x_*,u_*)$ ein zul"assiges Element der Aufgabe (\ref{EAAnhang2}).
\begin{enumerate}
\item[(A)] Wir nehmen an, dass der Punkt $(x_*,u_*)$ eine Umgebung $V$ mit folgenden Eigenschaften besitzt:
           \begin{enumerate}
           \item[(A$_1$)] Die Funktion $J(x,u)$ ist im Punkt $(x_*,u_*)$ Fr\'echet-differenzierbar;
           \item[(A$_2$)] Die Abbildung $\mathscr{F}(x,u)$ ist auf der Umgebung $V$ Fr\'echet-differenzierbar
                          und im Punkt $(x_*,u_*)$ stetig Fr\'echet-differenzierbar.
           \item[(A$_3$)] Die Abbildung $G(x)$ ist im Punkt $x_*$ Fr\'echet-differenzierbar.
           \end{enumerate}
\item[(B)] Weiterhin setzen wir voraus, dass der Operator $\mathscr{F}_x(x_*,u_*)$ eine endliche Kodimension besitzt.
\end{enumerate}
Ist dann $(x_*,u_*)$ schwache lokale Minimalstelle der Aufgabe (\ref{EAAnhang2}),
so ist f"ur die Aufgabe (\ref{EAAnhang2}) das Lagrangesche Prinzip g"ultig,
d.\,h., es existieren nicht gleichzeitig verschwindende Lagrangesche Multiplikatoren $\lambda_0 \geq 0$, $y^* \in \mathscr{Y}^*$
und $z^* \in \mathscr{Z}^*$ derart,
dass folgende Bedingungen gelten:
\begin{enumerate}
\item[(a)] Die Lagrange-Funktion besitzt bez"uglich $x$ in $x_*$ einen station"aren Punkt, d.\,h.
           \begin{equation}\label{SatzEPschwach1}
           0 = \mathscr{L}_x(x_*,u_*,\lambda_0,y^*,z^*);
           \end{equation}   
\item[(b)] Die Lagrange-Funktion erf"ullt bez"uglich $u$ in $u_*$ die Variationsungleichung
           \begin{equation}\label{SatzEPschwach2}
           0 \leq \langle \mathscr{L}_u(x_*,u_*,\lambda_0,y^*,z^*), u-u_* \rangle \qquad \mbox{f"ur alle } u \in \mathscr{U};
           \end{equation}
\item[(c)] Die komplement"aren Schlupfbedingungen gelten, d.\,h.
           \begin{equation}\label{SatzEPschwach3}
           0 = \langle z^*, G(x_*) \rangle, \qquad z^* \in \mathscr{K}^*.
           \end{equation}
\end{enumerate}
\end{theorem}

Wir werden im Weiteren folgende Bezeichnungen verwenden:
$$L_0={\rm Im\,}\mathscr{F}_x(x_*,u_*) \subseteq \mathscr{Y},$$
die Wertemenge des stetigen linearen Operators $\mathscr{F}_x(x_*,u_*)$;
$$B = L_0 + \mathscr{F}_u(x_*,u_*)(\mathscr{U}-u_*),$$
die Gesamtheit derjenigen $y \in \mathscr{Y}$, zu denen es ein $x \in \mathscr{X}$ und ein $u \in \mathscr{U}$ gibt mit
$$y=\mathscr{F}_x(x_*,u_*)x+\mathscr{F}_u(x_*,u_*)(u-u_*).$$
Au"serdem bezeichnet $L= {\rm lin\,} B$ die lineare H"ulle der Menge $B$.
Nach Voraussetzung hat der Teilraum $L_0$ eine endliche Kodimension,
so dass $L_0$ und $L$ abgeschlossene Teilr"aume von $\mathscr{Y}$ sind.
Die Annahme, dass ${\rm Im\,}\mathscr{F}_x(x_*,u_*)$ abgeschlossen sei, w"are im folgenden Beweis nicht ausreichend.

\begin{lemma} \label{LemmaEPschwach}
Es sei $L=\mathscr{Y}$.
Dann ist ${\rm int\,}B \not= \emptyset$.
Ist au"serdem $0 \in {\rm int\,}B$, so existieren $u_1,...,u_m \in \mathscr{U}$ derart, dass
f"ur $z_j=\pi\big(\mathscr{F}_u(x_*,u_*)(u_j-u_*)\big)$ folgende Beziehungen gelten:
$${\rm lin\,}\{z_1,...,z_m\}= \mathscr{Y}/L_0, \qquad z_1+...+z_m=0.$$
Dabei bezeichnet $\pi:\mathscr{Y} \to \mathscr{Y}/L_0$ die kanonische Abbildung,
d.\,h. es gilt $\pi y_1=\pi y_2$ genau dann, wenn $y_1-y_2 \in L_0$ ist.
\end{lemma}

{\bf Beweis:} Da die Kodimension der Menge $L_0$ endlich ist,
ist der Quotientenraum $\mathscr{Y}/L_0$ endlichdimensional.
Die Menge $B$ ist offensichtlich konvex.
Deshalb ist auch $\pi(B)$ konvex.
Da die lineare H"ulle der Menge $B$ mit $\mathscr{Y}$ "ubereinstimmt,
f"allt die lineare H"ulle der Menge $\pi(B)$ mit $\mathscr{Y}/L_0$ zusammen.
Wegen der Konvexit"at von $\pi(B)$ muss daher die Menge $\pi(B)$ ein nichtleeres Inneres besitzen.
Au"serdem ist die Beziehung $\pi^{-1}\big(\pi(B)\big) = B$ erf"ullt.
Weil $\pi$ eine stetige Abbildung ist, bedeutet dies, dass ${\rm int\,}B \not= \emptyset$ gilt. \\
Weil $\mathscr{Y}/L_0$ endlichdimensional ist,
existieren im Fall $0 \in {\rm int\,}B$ offenbar endlich viele Punkte $z_1,...,z_m$ aus $\pi(B)$ mit den geforderten Eigenschaften.
Wir w"ahlen dazu einfach $z_j$ als die Ecken eines hinreichend kleinen Würfels entsprechender Dimension, dessen Mittelpunkt im Ursprung liegt.
Nach Definition der Menge $\pi(B)$ gibt es ferner Elemente $u_j \in \mathscr{U}$ mit
$z_j=\pi\big(\mathscr{F}_u(x_*,u_*)(u_j-u_*)\big)$. \hfill $\blacksquare$ \\[2mm]
Wir wenden uns nun dem Beweis des Extremalprinzips \ref{SatzExtremalprinzipSchwach} zu.
Der Beweis ist in drei F"alle aufgeteilt, n"amlich zwei entarteten und einem nichtentarteten. \\[2mm]
{\bf Beweis im ersten entarteten Fall:} Es sei $L \not= \mathscr{Y}$.
Dann existiert nach Folgerung \ref{FolgerungAnnulator} ein nichttriviales Funktional $y^* \in \mathscr{Y}^*$ mit
$0=\langle y^*, \mathscr{F}_x(x_*,u_*)x+\mathscr{F}_u(x_*,u_*)(u-u_*) \rangle$
f"ur alle $x \in \mathscr{X}$ und $u \in \mathscr{U}$.
Setzen wir au"serdem $\lambda_0=0$ und $z^*=0$,
dann gelten die Bedingungen des Extremalprinzips \ref{SatzExtremalprinzipSchwach}. \hfill $\blacksquare$ \\[2mm]
{\bf Beweis im zweiten entarteten Fall:} Es sei $L = \mathscr{Y}$ und $0 \not\in {\rm int\,}B$.
Da die konvexe Menge $B$ in diesem Fall ein nichtleeres Inneres besitzt,
existiert nach dem Trennungssatz ein nichttriviales $y^* \in \mathscr{Y}^*$,
das die Menge $B$ und das Nullelement trennt.
Dies bedeutet,
dass f"ur alle $x \in \mathscr{X}$ und $u \in \mathscr{U}$ die Ungleichung
$$0\leq \langle y^*, \mathscr{F}_x(x_*,u_*)x+\mathscr{F}_u(x_*,u_*)(u-u_*) \rangle$$
erf"ullt ist.
Setzen wir hierin $u=u_*$, so erhalten wir $0\leq \langle y^*, \mathscr{F}_x(x_*,u_*)x \rangle$ f"ur alle $x \in \mathscr{X}$.
Folglich gilt $\mathscr{F}_x^*(x_*,u_*)y^*=0$.
F"ur $x=0$ ergibt sich f"ur alle $u \in \mathscr{U}$ die Ungleichung $0\leq \langle y^*, \mathscr{F}_u(x_*,u_*)(u-u_*) \rangle$.
Wie im vorhergehenden Fall sind daher $\lambda_0=0$, $z^*=0$ und $y^*$ ein System gesuchter Multiplikatoren. \hfill $\blacksquare$ \\[2mm]
{\bf Beweis im regul"aren Fall:} Es sei $L = \mathscr{Y}$ und $0 \in {\rm int\,}B$.
Wir nehmen $J(x_*,u_*)=0$ an.
Im Weiteren sei $\mathscr{C}$ die Menge derjenigen
$(\eta_0,y,z) \in \R \times \mathscr{Y} \times \mathscr{Z}$ mit folgender Eigenschaft:
Zu jedem Element existieren ein $x \in \mathscr{X}$, ein $u \in \mathscr{U}$ und ein $\eta \in {\rm int\,}\mathscr{K}$ mit
$$\eta_0 > J'(x_*,u_*)(x,u-u_*), \quad y = \mathscr{F}'(x_*,u_*)(x,u-u_*), \quad z = G(x_*) + G'(x_*)x - \eta.$$
Die Menge $\mathscr{C}$ ist konvex, da nach Voraussetzung die Menge $\mathscr{U}$ konvex ist.
Zum Beweis des Extremalprinzips \ref{SatzExtremalprinzipSchwach}  im regul"aren Fall gen"ugt es ${\rm int\,}\mathscr{C} \not= \emptyset$ und
$0 \not\in {\rm int\,}\mathscr{C}$ nachzuweisen.
Denn falls diese Beziehungen gelten,
dann folgt die Existenz eines nichttrivialen Funktionals $(\lambda_0,y^*,z^*) \in \R \times \mathscr{Y}^* \times \mathscr{Z}^*$,
das die Menge ${\rm int\,}\mathscr{C}$ vom Ursprung trennt,
d.\,h. ein solches Funktional,
dass $0 \leq \lambda_0 \eta_0 + \langle y^*,y \rangle + \langle z^*,z \rangle$ f"ur alle $(\eta_0,y,z) \in\mathscr{C}$ gilt.
Beachten wir $J(x_*,u_*)=0$ und $\mathscr{F}(x_*,u_*)=0$,
so erhalten wir daraus, dass f"ur alle $x \in \mathscr{X}$, $u \in \mathscr{U}$ und alle $\eta \in \mathscr{K}$ die Beziehung
\begin{eqnarray*} 
       \langle z^*,\eta \rangle
&\leq& \lambda_0 [J(x_*,u_*) + J'(x_*,u_*)(x,u-u_*)] \\
&    & + \langle y^*, \mathscr{F}(x_*,u_*) + \mathscr{F}'(x_*,u_*)(x,u-u_*) \rangle + \langle z^*,G(x_*)+G'(x_*)x \rangle
\end{eqnarray*}
erf"ullt ist.
Betrachten wir diese Ungleichung zun"achst f"ur $x=0$ und $u=u_*$,
dann folgt $\langle z^*,\eta \rangle \leq \langle z^*,G(x_*) \rangle$ f"ur alle $\eta \in {\rm int\,}\mathscr{K}$.
Da $\mathscr{K}$ ein abgeschlossener Kegel mit Spitze im Ursprung ist, ergeben sich daraus die Beziehungen
$\langle z^*,z \rangle \leq 0$ f"ur alle $z \in \mathscr{K}$ und $\langle z^*,G(x_*) \rangle = 0$.
Damit ist (\ref{SatzEPschwach3}) gezeigt.
Weiter ergibt sich f"ur alle $x \in \mathscr{X}$ und $u \in \mathscr{U}$ die Ungleichung
$0 \leq \langle \mathscr{L}_x(x_*,u_*,\lambda_0,y^*,z^*), x \rangle + \langle \mathscr{L}_u(x_*,u_*,\lambda_0,y^*,z^*),u-u_* \rangle$.
Betrachten wir nun nacheinander $u=u_*$ und $x=0$,
so kommen wir zu den Beziehungen (\ref{SatzEPschwach1}), (\ref{SatzEPschwach2}) des Extremalprinzips \ref{SatzExtremalprinzipSchwach}. \\[2mm]
Wir zeigen ${\rm int\,}\mathscr{C} \not= \emptyset$:
Es seien $u_1,...,u_m$ die Elemente mit den entsprechenden Eigenschaften in Lemma \ref{LemmaEPschwach},
$U_0={\rm conv\,} \{u_1,...,u_m\} \subseteq \mathscr{U}$, $S=\{x \in \mathscr{X}\,|\, \|x\| <1\}$ und
$$B_0=\mathscr{F}_x(x_*,u_*)S + \mathscr{F}_u(x_*,u_*)(U_0-u_*).$$
Die Menge $B_0$ ist konvex und besitzt ein nichtleeres Inneres,
da $\pi(B_0)$ die Punkte $z_1,...,z_m$ (vgl. Lemma \ref{LemmaEPschwach}) enth"alt.
Au"serdem ist nach dem Satz von der offenen Abbildung die Menge $\mathscr{F}_x(x_*,u_*)S$ offen in $L_0$.
Ferner seien
$$c_0 = \|J_x(x_*,u_*)\| +\max_{j=1,...,m} J_u(x_*,u_*)(u_j-u_*), \qquad \R_{c_0} = \{ a \in \R \,|\, a > c_0\}.$$
Dann besitzt die Menge $\mathscr{C}_0= \R_{c_0} \times B_0 \times (G(x_*)-{\rm int\,}\mathscr{K})$
ein nichtleeres Inneres und es gilt $\mathscr{C}_0 \subseteq \mathscr{C}$.
Demzufolge ist ${\rm int\,}\mathscr{C} \not= \emptyset$. \\[2mm]
Wir zeigen $0 \not\in {\rm int\,}\mathscr{C}$:
Wir nehmen nun an, es ist $0 \in {\rm int\,}\mathscr{C}$.
Dann existieren Elemente $\overline{x} \in \mathscr{X}$, $\overline{u} \in \mathscr{U}$, $\overline{\eta} \in {\rm int\,}\mathscr{K}$
und eine Zahl $c>0$ mit
\begin{equation} \label{BeweisWEPr2}
-c > J'(x_*,u_*)(\overline{x},\overline{u}-u_*), \quad 0 = \mathscr{F}'(x_*,u_*)(\overline{x},\overline{u}-u_*), \quad
0 = G(x_*)+ G'(x_*)\overline{x} - \overline{\eta}.
\end{equation}
Angenommen,
die Beziehungen (\ref{BeweisWEPr2}) seien erf"ullt.
Es sei $\varrho >0$ fest gew"ahlt.
Ferner seien $u_1,...,u_m \in \mathscr{U}$ die Elemente in Lemma \ref{LemmaEPschwach}.
Nach Voraussetzung (A$_2$) ist $\mathscr{F}$ auf der Umgebung $V$ des Punktes $(x_*,u_*)$ stetig.
Daher ist in einer Umgebung des Punktes $(x_*,0,0) \in \mathscr{X} \times \R \times \R^m$ durch
$$\Phi(x,\alpha_0,\alpha) = \mathscr{F}\bigg(x,u_*+\alpha_0(\overline{u}-u_*) + \varrho\sum_{j=1}^m \alpha_j (u_j-u_*)\bigg),
  \quad \alpha=(\alpha_1,...,\alpha_m),$$
eine Abbildung in den Raum $\mathscr{Y}$ definiert.
Dabei gilt $\Phi(x_*,0,0) = 0$.
Die Bedingung (A$_2$) liefert,
dass $\Phi$ auf einer Umgebung des Punktes $(x_*,0,0)$ Fr\'echet-differenzierbar und im Punkt $(x_*,0,0)$ stetig Fr\'echet-differenzierbar mit der Ableitung
$$\Phi'(x_*,0,0)(x,\alpha_0,\alpha) = \mathscr{F}'(x_*,u_*)\big(x,\alpha_0(\overline{u}-u_*)\big) + \varrho\sum_{j=1}^m \alpha_j \mathscr{F}_u(x_*,u_*)(u_j-u_*)$$
ist.
Ferner enth"alt die Wertemenge des linearen Operators $\Phi'(x_*,0,0)$ die Menge $B_0$ und stimmt folglich mit dem ganzen Raum $\mathscr{Y}$ "uberein.
Au"serdem gibt es nach Wahl der Elemente $u_j$ in Lemma \ref{LemmaEPschwach} ein $x' \in \mathscr{X}$ mit
$$\mathscr{F}_x(x_*,u_*)x' + \sum_{j=1}^m \mathscr{F}_u(x_*,u_*)(u_j-u_*) =0.$$
Damit folgt zusammen mit (\ref{BeweisWEPr2}) die Bedingung $\Phi'(x_*,0,0)(\overline{x} + \varrho x',1,1)=0$,
d.\,h., dass der Punkt $(\overline{x}+ \varrho x',1,1)$ dem Kern des Operators $\Phi'(x_*,0,0)$ angeh"ort. \\
Nach dem Satz von Ljusternik (Theorem \ref{SatzLjusternik}) existieren $\varepsilon_0 >0$ und eine Abbildung
$\varepsilon \to \big(x(\varepsilon),\alpha_0(\varepsilon),\alpha(\varepsilon)\big)$ des Intervalls $[0,\varepsilon_0]$ in den Raum
$\mathscr{X} \times \R \times \R^m$ derart, dass
$\big\|\big(x(\varepsilon),\alpha_0(\varepsilon),\alpha(\varepsilon)\big)\big\|_{\mathscr{X} \times \R \times \R^m} \to 0$
f"ur $\varepsilon \to 0^+$ gilt und au"serdem
\begin{equation} \label{BeweisWEPr5}
\Phi\big(x_* + \varepsilon [\overline{x} + \varrho x' + x(\varepsilon)],\varepsilon[1 + \alpha_0(\varepsilon)],
         \varepsilon[1 + \alpha(\varepsilon)]\big) =\Phi(x_*,0,0)
\end{equation}
f"ur alle $\varepsilon \in [0,\varepsilon_0]$ erf"ullt ist.
Wir setzen der K"urze halber $\tilde{x}(\varepsilon) = x_* + \varepsilon [\overline{x} + \varrho x' + x(\varepsilon)]$ und
$$\tilde{u}(\varepsilon) = u_*+\varepsilon[1 + \alpha_0(\varepsilon)](\overline{u}-u_*)
                           + \varepsilon\varrho\sum_{j=1}^m [1+ \alpha_j(\varepsilon)] (u_j-u_*).$$
F"ur alle $\varepsilon \in [0,\varepsilon_0]$ mit hinreichend kleinem $\varepsilon_0>0$ stellt $\tilde{u}(\varepsilon)$  eine
Konvexkombination der Elemene $u_*,\overline{u},u_1,...,u_m$ dar, d.\,h. $\tilde{u}(\varepsilon) \in \mathscr{U}$,
und es folgt unmittelbar aus (\ref{BeweisWEPr5}):
\begin{equation} \label{BeweisWEPr6}
\mathscr{F}\big(\tilde{x}(\varepsilon),\tilde{u}(\varepsilon)\big) = 0.
\end{equation}
Da die Abbildung $G$ im Punkt $x_*$ Fr\'echet-differenzierbar ist,
ergibt sich
$$G\big(\tilde{x}(\varepsilon)\big) = G(x_*) + \varepsilon G'(x_*)(\overline{x}+\varrho x') + r(\varepsilon), \quad
  \lim_{\varepsilon \to 0^+} \frac{\|r(\varepsilon)\|_{\mathscr{Z}}}{\varepsilon}=0.$$
Wir beachten,
dass $\mathscr{K}$ konvex ist,
sowie $G(x_*) \in \mathscr{K}$ und $G(x_*)+ G'(x_*)\overline{x} = \overline{\eta} \in {\rm int\,}\mathscr{K}$ nach (\ref{BeweisWEPr2}) gelten.
Es sei $\varrho >0$ mit $\overline{\eta} + \varrho G'(x_*) x' \in {\rm int\,}\mathscr{K}$ gew"ahlt. 
Damit erhalten wir f"ur alle $\varepsilon \in (0,\varepsilon_0]$ mit einem hinreichend kleinen $\varepsilon_0$:
\begin{eqnarray}
G\big(\tilde{x}(\varepsilon)\big) &=& (1-\varepsilon) G(x_*) + \varepsilon [G(x_*)+ G'(x_*)\overline{x}]
                                      +\varepsilon \varrho  G'(x_*) x' + r(\varepsilon) \nonumber \\
                                  &=& (1-\varepsilon) G(x_*) + \varepsilon [\overline{\eta} + \varrho G'(x_*) x' + r(\varepsilon) / \varepsilon]
\label{BeweisWEPr7}                    \in {\rm int\,}\mathscr{K}.
\end{eqnarray}
Zus"atzlich zu der Bedingung $\overline{\eta} + \varrho G'(x_*) x' \in {\rm int\,}\mathscr{K}$ w"ahlen wir $\varrho>0$ so,
dass mit der Zahl $c >0$ in (\ref{BeweisWEPr2}) die Relation
$$\varrho \bigg[ J_x(x_*,u_*)x' + \sum_{j=1}^m J_u(x_*,u_*)(u_j-u_*) \bigg] \leq  \frac{c}{2}$$
erf"ullt ist.
Da das Funktional $J$ im Punkt $(x_*,u_*)$ Fr\'echet-differenzierbar ist, ergibt sich
\begin{eqnarray*}
    J\big(\tilde{x}(\varepsilon),\tilde{u}(\varepsilon)\big)
&=& J(x_*,u_*) + \varepsilon J'(x_*,u_*)(\overline{x},\overline{u}-u_*) \\
& &   + \varepsilon \varrho \bigg[ J_x(x_*,u_*)x' + \sum_{j=1}^m J_u(x_*,u_*)(u_j-u_*) \bigg] + o(\varepsilon).
\end{eqnarray*}
Mit (\ref{BeweisWEPr2}) und nach Wahl von $\varrho>0$ erhalten wir f"ur alle $\varepsilon \in [0,\varepsilon_0]$:
\begin{equation} \label{BeweisWEPr8}
J\big(\tilde{x}(\varepsilon),\tilde{u}(\varepsilon)\big) \leq J(x_*,u_*) - \varepsilon c + \varepsilon \frac{c}{2} + o(\varepsilon).
\end{equation}
Die Beziehungen (\ref{BeweisWEPr6}) und (\ref{BeweisWEPr7}) zeigen,
dass f"ur hinreichend kleine $\varepsilon >0$ das Paar $\big(\tilde{x}(\varepsilon),\tilde{u}(\varepsilon)\big)$ zul"assig in der Aufgabe (\ref{EAAnhang2}) ist.
Au"serdem ist $J\big(\tilde{x}(\varepsilon),\tilde{u}(\varepsilon)\big) < J(x_*,u_*)$ nach (\ref{BeweisWEPr8}) f"ur hinreichend kleine $\varepsilon >0$.
Wegen $\tilde{x}(\varepsilon) \to x_*$ und $\tilde{u}(\varepsilon) \to u_*$ f"ur $\varepsilon \to 0$ bedeutet dies,
dass der Punkt $(x_*,u_*)$ im Widerspruch zur Voraussetzung kein schwaches lokales Minimum sein k"onnte. \hfill $\blacksquare$

%% file: VI-3-StarkesExtremalprinzip.tex
\subsection{Ein Extremalprinzip f\"ur ein starkes lokales Minimum}
Es seien $\mathscr{X}$, $\mathscr{Y}$ Banachr"aume und $\mathscr{U}$ eine beliebige Menge.
Ferner seien $J$ ein Funktional auf $\mathscr{X} \times \mathscr{U}$,
$\mathscr{F}$ eine Abbildung des Produktes $\mathscr{X} \times \mathscr{U}$ in den Raum $\mathscr{Y}$ und $G_j:\mathscr{X} \to \R$ f"ur $j=1,...,l$.
Unter diesen Angaben betrachten wir die Extremalaufgabe
\begin{equation}\label{EAAnhang3}
J(x,u) \to \inf; \quad \mathscr{F}(x,u)=0, \quad G_j(x) \leq 0, \; j=1,...,l, \quad x \in \mathscr{X},\; u \in \mathscr{U}.
\end{equation}
Der Punkt $(x,u)$ ist ein zul"assiges Element der Aufgabe (\ref{EAAnhang3}),
falls s"amtliche Nebenbedingungen erf"ullt sind.
Ein Punkt $(x_*,u_*)$ hei"st ein starkes lokales Minimum\index{Minimum, starkes lokales!Extremal@-- Extremalaufgabe}
der Extremalaufgabe (\ref{EAAnhang3}),
wenn ein $\varepsilon >0$ derart existiert,
dass f"ur alle zul"assigen Paare $(x,u)$ mit $\|x-x_*\|_{\mathscr{X}} \leq \varepsilon$ die Ungleichung $J(x_*,u_*) \leq J(x,u)$ gilt. \\
Auf $\mathscr{X} \times \mathscr{U} \times \R \times \mathscr{Y}^* \times \R^l$ definieren wir zur Aufgabe (\ref{EAAnhang3}) die Lagrange-Funktion
$$\mathscr{L}(x,u,\lambda_0,y^*,\lambda_1,...,\lambda_l)
  = \lambda_0 J(x,u) + \langle y^*, \mathscr{F}(x,u) \rangle + \sum_{j=1}^l  \lambda_j G_j(x).$$
Au"serdem bezeichnet in diesem Abschnitt $\Sigma(\Delta)$ das folgende $d$-dimensionale Simplex:
$$\Sigma(\Delta) = \bigg\{ \alpha = (\alpha_1,...,\alpha_d) \in \R^d \,\bigg|\, \alpha_1,...,\alpha_d \geq 0,\, \sum_{i=1}^{d} \alpha_i \leq \Delta \bigg\}.$$

\begin{theorem}[Extremalprinzip] \label{SatzExtremalprinzipStark}
\index{Extremalprinzip!Stark@-- starkes lokales Minimum}
Sei $(x_*,u_*)$ ein zul"assiges Element der Aufgabe (\ref{EAAnhang3}).
\begin{enumerate}
\item[(A)] Wir nehmen an, dass der Punkt $x_*$ eine Umgebung $V$ besitzt mit:
           \begin{enumerate}
           \item[(A$_1$)] F"ur jedes $u \in U$ ist $x \to J(x,u)$ im Punkt $x_*$ Fr\'echet-differenzierbar;
           \item[(A$_2$)] F"ur jedes $u \in U$ ist $x \to \mathscr{F}(x,u)$ im Punkt $x_*$ G\^ateaux-differenzierbar;
           \item[(A$_3$)] Die Funktionen $G_j(x)$ sind im Punkt $x_*$ lokalkonvex und bez"uglich jeder Richtung gleichm"a"sig differenzierbar.
           \end{enumerate}
\item[(B)] Weiterhin setzen wir voraus, dass der Operator $\mathscr{F}_x(x_*,u_*)$ eine endliche Kodimension besitzt.
\item[(C)] Au"serdem nehmen wir an,
           dass zu jedem endlichen System von Punkten $u_1,...,u_d$ aus $\mathscr{U}$ und zu jedem $\delta>0$ eine Zahl $\Delta >0$ und eine Abbildung
           $u:\Sigma(\Delta) \to \mathscr{U}$ derart existieren, dass $u(0)=u_*$ gilt und f"ur alle $x,x'$ aus $V'$ und alle $\alpha,\alpha'$
           aus $\Sigma(\Delta)$ folgende Ungleichungen erf"ullt sind:
           \begin{eqnarray*}
           && \hspace*{-5mm} \bigg\| \mathscr{F}\big(x,u(\alpha)\big) - \mathscr{F}\big(x',u(\alpha')\big)-\mathscr{F}_x(x_*,u_*)(x-x') \\
           && -\sum_{k=1}^d (\alpha_k - \alpha_k') \big( \mathscr{F}(x_*,u_k)-\mathscr{F}(x_*,u_*)\big) \bigg\|_{\mathscr{Y}} 
              \leq \delta \bigg( \|x-x'\|_{\mathscr{X}} + \sum_{k=1}^d |\alpha_k - \alpha_k'|\bigg), \\
           && \hspace*{-5mm} J\big(x,u(\alpha)\big) - J(x,u_*) - \sum_{k=1}^d \alpha_k \big( J(x,u_k)-J(x,u_*)\big)
              \leq \delta \bigg( \|x-x_*\|_{\mathscr{X}} + \sum_{k=1}^d \alpha_k \bigg).
           \end{eqnarray*}
\end{enumerate}
Ist dann $(x_*,u_*)$ starke lokale Minimalstelle der Aufgabe (\ref{EAAnhang3}),
so existieren nicht gleichzeitig verschwindende Lagrangesche Multiplikatoren $\lambda_0 \geq 0,...,\lambda_l \geq 0$ und $y^* \in \mathscr{Y}^*$ derart,
dass folgende Bedingungen gelten:
\begin{enumerate}
\item[(a)] Die Lagrange-Funktion besitzt bez"uglich $x$ in $x_*$ einen station"aren Punkt, d.\,h.
           \begin{equation}\label{SatzEPstark1}
           0 \in \partial_x \mathscr{L}(x_*,u_*,\lambda_0,y^*,\lambda_1,...,\lambda_l);
           \end{equation}   
\item[(b)] Die Lagrange-Funktion erf"ullt bez"uglich $u$ in $u_*$ die Minimumbedingung
           \begin{equation}\label{SatzEPstark2}
           \mathscr{L}(x_*,u_*,\lambda_0,y^*,\lambda_1,...,\lambda_l)
           = \min_{u \in \mathscr{U}} \mathscr{L}(x_*,u,\lambda_0,y^*,\lambda_1,...,\lambda_l);
           \end{equation}
\item[(c)] Die komplement"aren Schlupfbedingungen gelten, d.\,h.
           \begin{equation}\label{SatzEPstark3}
           0 = \lambda_j G_j(x_*) \qquad\mbox{f"ur } j=1,...,l.
           \end{equation}
\end{enumerate}
\end{theorem}

Im Weiteren seien $L_0={\rm Im\,}\mathscr{F}_x(x_*,u_*)$, $B = L_0 + {\rm conv\,}\mathscr{F}(x_*,\mathscr{U})$ und $L= {\rm lin\,} B$.
Zum Beweis des Extremalprinzips \ref{SatzExtremalprinzipStark} gehen wir "ahnlich wie im Abschnitt \ref{AbschnittEPschwach} vor.
Die entarteten F"alle ergeben sich dabei v"ollig analog.

\begin{lemma} \label{LemmaEPstark}
Es sei $L=\mathscr{Y}$.
Dann ist ${\rm int\,}B \not= \emptyset$.
Ist au"serdem $0 \in {\rm int\,}B$, so gibt es $x' \in \mathscr{X}$, $u_1,...,u_m \in \mathscr{U}$ und $\gamma_j >0$ mit
$\mathscr{F}_x(x_*,u_*)x' + \gamma_1\mathscr{F}(x_*,u_1) + ...+ \gamma_m \mathscr{F}(x_*,u_m) = 0$ und
${\rm lin\,} \big(L_0 \cup \{ \mathscr{F}(x_*,u_1),...,\mathscr{F}(x_*,u_m)\}\big) =\mathscr{Y}$.
\end{lemma}

{\bf Beweis:} Der erste Teil, n"amlich 
$\mathscr{F}_x(x_*,u_*)x' + \gamma_1\mathscr{F}(x_*,u_1) + ...+ \gamma_m \mathscr{F}(x_*,u_m) = 0$,
ergibt sich ebenso wie im Beweis von Lemma \ref{LemmaEPschwach}.
Ist $L=\mathscr{Y}$ und $0 \in {\rm int\,}B$,
so gibt es $x' \in \mathscr{X}$ und $y_1,...,y_r \in {\rm conv\,}\mathscr{F}(x_*,\mathscr{U})$ mit (Lemma \ref{LemmaEPschwach})
$$\mathscr{F}_x(x_*,u_*)x' + y_1 + ... + y_r =0, \quad {\rm lin\,} (L_0 \cup \{ y_1,...,y_r\}) =\mathscr{Y}.$$
Weiterhin existieren zu jedem $j=1,...,r$ Zahlen $\gamma_{js}>0$ und Punkte $u_{js} \in \mathscr{U}$ mit
$$y_j = \sum_{s=1}^{m_j} \gamma_{js} \mathscr{F}(x_*,u_{js}), \qquad \sum_{s=1}^{m_j} \gamma_{js} =1.$$
Damit sind $\gamma_{js}>0$ die gesuchten Zahlen und $u_{js}$ die gesuchten Elemente. \hfill $\blacksquare$ \\[2mm]
{\bf Beweis in den entarteten F"allen:} Ist $L \not= \mathscr{Y}$,
so existiert ein nichttriviales Funktional $y^* \in \mathscr{Y}^*$, das dem Annulator des Teilraumes $L$ angeh"ort,
und $\lambda_0=...=\lambda_l=0$ und $y^*$ sind die gesuchten Multiplikatoren. \\[2mm]
Ist $L = \mathscr{Y}$ und $0 \not\in {\rm int\,}B$, 
dann besitzt die konvexe Menge $B$ ein nichtleeres Inneres.
Nach dem Trennungssatz existiert ein nichttriviales $y^* \in \mathscr{Y}^*$,
das die Menge $B$ und den Ursprung im Raum $\mathscr{Y}$ trennt.
Daher sind $\lambda_0=...=\lambda_l=0$ und $y^*$ ein System gesuchter Multiplikatoren. \hfill $\blacksquare$ \\[2mm]
{\bf Beweis im regul"aren Fall:} Es sei $L = \mathscr{Y}$ und $0 \in {\rm int\,}B$.
Wir nehmen $J(x_*,u_*)=0$ an.
Au"serdem seien $G_j(x_*)=0$ f"ur $j=1,...,l'$ und $G_j(x_*)<0$ f"ur $j=l'+1,...,l$.
Im Weiteren bezeichnet $\mathscr{C}$ die Menge derjenigen
$(\eta_0,y,\eta_1,...,\eta_{l'}) \in \R \times \mathscr{Y} \times \R^{l'}$ mit folgender Eigenschaft:
Zu jedem Element existieren ein $x \in \mathscr{X}$ und ein $u \in \mathscr{U}$ mit
\begin{eqnarray*}
&&\hspace*{-5mm} \eta_0 > J_x(x_*,u_*)x + J(x_*,u)-J(x_*,u_*), \quad y = \mathscr{F}_x(x_*,u_*)x + \mathscr{F}(x_*,u)-\mathscr{F}(x_*,u_*), \\
&&\hspace*{-5mm} \eta_j > G_j(x_*)+G'_j(x_*;x), \quad j=1,...,l'.
\end{eqnarray*}
Zum Beweis des Extremalprinzips \ref{SatzExtremalprinzipStark} gen"ugt es im regul"aren Fall ${\rm int\,}({\rm conv\,}\mathscr{C}) \not= \emptyset$
und $0 \not\in {\rm conv\,}\mathscr{C}$ nachzuweisen.
Denn falls diese Beziehungen richtig sind,
dann existiert ein nichttriviales Funktional $(\lambda_0,y^*,\lambda) \in \R \times \mathscr{Y}^* \times \R^{l'}$,
das die Menge ${\rm conv\,}\mathscr{C}$ vom Urpsrung trennt.
Dies bedeutet, dass die Ungleichung $0 \leq \lambda_0 \eta_0 + \langle y^* , y\rangle + \lambda_1 \eta_1 + ... +\lambda_{l'} \eta_{l'}$
f"ur alle $(\eta_0,y,\eta_1,...,\eta_{l'}) \in \mathscr{C} \subseteq {\rm conv\,}\mathscr{C}$ gilt.
In dieser Ungleichung folgt durch die Betrachtung $\eta_j \to \infty$ sofort $\lambda_j\geq 0$ f"ur $j=0,...,l'$ und,
dass f"ur alle $x \in \mathscr{X}$, $u \in \mathscr{U}$ die Ungleichung
\begin{eqnarray*}
0 &\leq &\lambda_0 [J_x(x_*,u_*)x + J(x_*,u)-J(x_*,u_*)] \\
&& + \langle y^*, \mathscr{F}_x(x_*,u_*)x + \mathscr{F}(x_*,u)-\mathscr{F}(x_*,u_*) \rangle +\sum_{j=1}^{l'} \lambda_j [G_j(x_*)+G_j'(x_*;x)]
\end{eqnarray*}
erf"ullt ist.
Setzen wir $\lambda_{l'+1}=...=\lambda_l=0$, dann gelten die Schlupfbedingungen (\ref{SatzEPstark3}).
Betrachten wir weiterhin diese Ungleichung nacheinander f"ur $x=0$ bzw. f"ur $u=u_*$,
so folgen (\ref{SatzEPstark1}) und (\ref{SatzEPstark2}). \\[2mm]
Wir zeigen ${\rm int\,}({\rm conv\,}\mathscr{C}) \not= \emptyset$:
Es seien $u_1,...,u_m$ die Elemente mit den entsprechenden Eigenschaften in Lemma \ref{LemmaEPstark}.
Im Weiteren bezeichnen $S=\{x \in \mathscr{X}\,|\, \|x\|_{\mathscr{X}} <1\}$ und
$$Y_0={\rm conv\,} \{\mathscr{F}(x_*,u_1),...,\mathscr{F}(x_*,u_m)\}, \qquad B_0=\mathscr{F}_x(x_*,u_*)S + Y_0.$$
Die Menge $B_0$ ist konvex und besitzt ein nichtleeres Inneres, da ${\rm lin\,}B_0=\mathscr{Y}$.
Au"serdem ist die Menge $\mathscr{F}_x(x_*,u_*)S$ nach dem Satz von der offenen Abbildung offen.
Wir setzen
$$c_1 = \|J_x(x_*,u_*)\| +\max_{j=1,...,m} J(x_*,u_j), \qquad c_2 = \max_{j=1,...,l'} \sup_{x \in S} |G'_j(x_*;x)|.$$
Die Zahl $c_2$ ist endlich, da die Richtungsableitung nach Lemma \ref{LemmaRichtungsableitung} eine stetige Funktion ist.
Ferner seien $c_0 = \max \{c_1,c_2\}$ und $\R_{c_0} = \{ a \in \R \,|\, a > c_0\}$.
Dann besitzt die Menge $\mathscr{C}_0= \R_{c_0} \times B_0 \times \R_{c_0}^{l'} \subseteq {\rm conv\,}\mathscr{C}$ ein nichtleeres Inneres.
Demzufolge ist ${\rm int\,}({\rm conv\,}\mathscr{C}) \not= \emptyset$. \\[2mm]
Wir zeigen $0 \not\in {\rm conv\,}\mathscr{C}$:
Es sei $0 \in {\rm conv\,}\mathscr{C}$.
Dann existieren Vektoren $\overline{\xi}_k \in \mathscr{X}$, $\overline{w}_k \in \mathscr{U}$
und $\overline{a}_1>0,...,\overline{a}_r>0$ mit $\overline{a}_1 +...+ \overline{a}_r=1$ derart,
dass mit $\overline{\xi} = \overline{a}_1 \overline{\xi}_1 +...+ \overline{a}_r \overline{\xi}_r$ die Beziehungen
\begin{eqnarray*}
&& 0 > J_x(x_*,u_*)\overline{\xi} + \sum_{k=1}^r \overline{a}_k \big[J(x_*,\overline{w}_k)-J(x_*,u_*)\big], \\
&& 0 = \mathscr{F}_x(x_*,u_*)\overline{\xi} + \sum_{k=1}^r \overline{a}_k \big[\mathscr{F}(x_*,\overline{w}_k)-\mathscr{F}(x_*,u_*)\big]
\end{eqnarray*}
erf"ullt sind und wegen der Konvexit"at der Abbildungen $x \to G_j'(x_*;x)$ folgende Ungleichungen f"ur $j=1,...,l'$ gelten:
$$0 > G_j(x_*) +  \sum_{k=1}^r \overline{a}_k G_j'(x_*;\overline{\xi}_k) \geq G_j(x_*)+ G_j'(x_*;\overline{\xi}).$$
Nach Lemma \ref{LemmaEPstark} existieren $x' \in \mathscr{X}$, $u_1,...,u_m \in \mathscr{U}$ und Zahlen $\gamma_1>0,...,\gamma_m >0$ mit
$$\mathscr{F}_x(x_*,u_*)x' + \gamma_1\mathscr{F}(x_*,u_1) + ...+ \gamma_m \mathscr{F}(x_*,u_m) = 0.$$
Da die Abbildungen $x \to J_x(x_*,u_*)x$, $x \to G_j(x_*;x)$ stetig und $J(x_*,u_j)$ endlich sind,
gibt es positive Zahlen $\varrho$ und $c>0$ mit
\begin{eqnarray*}
&&\hspace*{-7mm} -4c > J_x(x_*,u_*)\overline{\xi} + \sum_{k=1}^r \overline{a}_k \big[J(x_*,\overline{w}_k)-J(x_*,u_*)\big]
       +\varrho \bigg[ J_x(x_*,u_*)x' + \sum_{j=1}^m \gamma_j J(x_*,u_j) \bigg]\!, \\
&&\hspace*{-4mm} 0 = \mathscr{F}_x(x_*,u_*)\overline{\xi} + \sum_{k=1}^r \overline{a}_k \big[\mathscr{F}(x_*,\overline{w}_k)-\mathscr{F}(x_*,u_*)\big]
       +\varrho \bigg[ \mathscr{F}_x(x_*,u_*)x' + \sum_{j=1}^m \gamma_j \mathscr{F}(x_*,u_j) \bigg]\!, \\
&&\hspace*{-6mm} -c > G_j(x_*)+ G_j'(x_*;\overline{\xi}+\varrho x'), \quad j=1,...,l'.
\end{eqnarray*}
Wir f"uhren durch $\overline{x}=\overline{\xi}+\varrho x'$, $\overline{u}_1=\overline{w}_1,...,\overline{u}_d=u_m$
und $\overline{\alpha}_1=\overline{a}_1,...,\overline{\alpha}_d=\gamma_m$ eine vereinheitlichte Bezeichnungsweise ein.
Dann gelten im Fall $0 \in {\rm conv\,}\mathscr{C}$ die Beziehungen
\begin{eqnarray}
\label{BeweisSEPr2} -4c &>& J_x(x_*,u_*)\overline{x} + \sum_{k=1}^d \overline{\alpha}_k [J(x_*,\overline{u}_k) - J(x_*,u_*)], \\
\label{BeweisSEPr3} 0 &=& \mathscr{F}_x(x_*,u_*)\overline{x} + \sum_{k=1}^d \overline{\alpha}_k [\mathscr{F}(x_*,\overline{u}_k) - \mathscr{F}(x_*,u_*)], \\
\label{BeweisSEPr4} -c &>& G_j(x_*)+ G_j'(x_*;\overline{x}), \quad j=1,...,l',
\end{eqnarray}
und au"serdem, da die Elemente $u_j$ aus Lemma \ref{LemmaEPstark} in den $\overline{u}_k$ enthalten sind,
\begin{equation} \label{BeweisSEPr1}
{\rm lin\,} \big(L_0 \cup \{ \mathscr{F}(x_*,\overline{u}_1),...,\mathscr{F}(x_*,\overline{u}_d)\}\big) =\mathscr{Y}.
\end{equation}
Wir nehmen an,
die Beziehungen (\ref{BeweisSEPr2})--(\ref{BeweisSEPr4}) seien erf"ullt.
Es sei $\delta >0$ fest gew"ahlt.
Dann kann man f"ur $\overline{u}_1,...,\overline{u}_d$ und $\delta$ eine Umgebung $V'$ des Punktes $x_*$,
eine Zahl $\Delta >0$ und eine Abbildung $u: \Sigma(\Delta) \to \mathscr{U}$ w"ahlen,
die der Voraussetzung (C) gen"ugen. \\
F"ur jedes $a \in \R$ setzen wir $a^+=\max\{a,0\}$, $a^-=\min\{a,0\}$ und verwenden f"ur jeden Vektor
$\alpha = (\alpha_1,...,\alpha_m) \in \R^m$ die Bezeichnungen $\alpha^+=(\alpha_1^+,...,\alpha_m^+)$, $\alpha^-=(\alpha_1^-,...,\alpha_m^-)$. \\
Damit ist in einer Umgebung des Punktes $(x_*,0) \in \mathscr{X} \times \R^d$ durch
$$\Phi(x,\alpha) = \mathscr{F}\big(x,u(\alpha^+)\big) + \sum_{k=1}^d \alpha_k^- \big(\mathscr{F}(x_*,\overline{u}_k) - \mathscr{F}(x_*,u_*)\big)$$
eine Abbildung in den Raum $\mathscr{Y}$ definiert.
Offenbar gilt dabei $\Phi(x_*,0) = \mathscr{F}(x_*,u_*) = 0$. \\
Schlie"slich bezeichnen wir mit $\Lambda$ die Abbildung
$$\Lambda(x,\alpha) = \mathscr{F}_x(x_*,u_*)x+ \sum_{k=1}^d \alpha_k \big(\mathscr{F}(x_*,\overline{u}_k) - \mathscr{F}(x_*,u_*)\big).$$
Die Abbildung $\Lambda$ definiert einen stetigen linearen Operator aus $\mathscr{X} \times \R^d$ in $\mathscr{Y}$,
da nach Voraussetzung (A$_2$) der Operator $\mathscr{F}_x(x_*,u_*)$ stetig und linear ist.
Aufgrund der Voraussetzung (C) gilt daher f"ur alle $(x,\alpha)$ und $(x',\alpha')$ aus einer Umgebung des Punktes
$(x_*,0) \in \mathscr{X} \times \R^d$ die Ungleichung
\begin{eqnarray}
&& \|\Phi(x,\alpha)-\Phi(x',\alpha')-\Lambda(x,\alpha)+\Lambda(x',\alpha') \|_{\mathscr{Y}} \nonumber \\
&& =\bigg\| \mathscr{F}\big(x,u(\alpha^+)\big) - \mathscr{F}\big(x',u({\alpha'}^+)\big)-\mathscr{F}_x(x_*,u_*)(x-x') \nonumber \\
&& \hspace*{20mm} -\sum_{k=1}^d (\alpha_k^+ - {\alpha_k'}^+) \big( \mathscr{F}(x_*,\overline{u}_k)-\mathscr{F}(x_*,u_*)\big) \bigg\|_{\mathscr{Y}}
   \nonumber \\
&& \label{BeweisSEPr5}
   \leq \delta \bigg( \|x-x'\|_{\mathscr{X}} + \sum_{k=1}^d |\alpha^+_k - {\alpha^+_k}'|\bigg)
   \leq \delta \bigg( \|x-x'\|_{\mathscr{X}} + \sum_{k=1}^d |\alpha_k - \alpha_k'|\bigg).
\end{eqnarray}
Aus der Definition des Operators $\Lambda$ und (\ref{BeweisSEPr1}) folgt ${\rm Im\,}\Lambda=\mathscr{Y}$.
Gem"a"s Lemma \ref{LemmaVerallgemeinerterSatzLjusternik} ist die Zahl $C(\Lambda)< \infty$.
W"ahlen wir nun $\delta >0$ so, dass $\delta C(\Lambda)< 1/2$ erf"ullt ist,
so gen"ugen die Abbildung $\Phi$ und der Operator $\Lambda$ dem verallgemeinerten Satz von Ljusternik. \\
Nach (\ref{BeweisSEPr3}) geh"ort der Vektor $(\overline{x},\overline{\alpha})$ mit
$\overline{\alpha}=(\overline{\alpha}_1,...,\overline{\alpha}_d)$ dem Kern des Operators $\Lambda$ an.
Daher existieren nach dem verallgemeinerten Satz von Ljusternik eine Zahlen $\varepsilon_0 >0$, $K>0$ und Abbildungen
$\varepsilon \to \big(x(\varepsilon),\alpha(\varepsilon)\big) = \big(x(\varepsilon),\alpha_1(\varepsilon),...,\alpha_d(\varepsilon)\big)$
des Intervalls $[0,\varepsilon_0]$ in den Raum $\mathscr{X} \times \R^d$ derart, dass
\begin{equation} \label{BeweisSEPr6}
\Phi\big(x_* + \varepsilon \overline{x} + x(\varepsilon),\varepsilon \overline{\alpha} + \alpha(\varepsilon)\big) = \Phi(x_*,0)
\end{equation}
f"ur alle $\varepsilon \in [0,\varepsilon_0]$ und,
da $(\overline{x},\overline{\alpha})$ dem Kern des Operators $\Lambda$ angeh"ort,
\begin{equation} \label{BeweisSEPr6-1}
\|x(\varepsilon)\|_{\mathscr{X}} + \sum_{k=1}^d |\alpha_k(\varepsilon)|
\leq K \|\Phi(x_*+\varepsilon \overline{x}, \varepsilon \overline{\alpha}) - \Phi(x_*,0)
         - \varepsilon \Lambda(\overline{x},\overline{\alpha})\|_{\mathscr{Y}}
\end{equation}
gelten.
Aus (\ref{BeweisSEPr6-1}) ergibt sich zusammen mit (\ref{BeweisSEPr5}) die Ungleichung
\begin{equation} \label{BeweisSEPr7}
       \|x(\varepsilon)\|_{\mathscr{X}} + \sum_{k=1}^d |\alpha_k(\varepsilon)|
\leq \varepsilon K \delta \bigg( \|\overline{x}\|_{\mathscr{X}} + \sum_{k=1}^d \overline{\alpha}_k \bigg).
\end{equation}
Hieraus ergibt sich insbesondere $x(\varepsilon) \to 0$ und $\alpha(\varepsilon) \to 0$ f"ur $\varepsilon \to 0$.
Wir nehmen jetzt ferner an, dass $\delta >0$ neben $\delta C(\Lambda)< 1/2$ noch folgende Ungleichung erf"ullt:
\begin{equation}\label{BeweisSEPr8}
K \delta \bigg( \|\overline{x}\|_{\mathscr{X}} + \sum_{k=1}^d \overline{\alpha}_k \bigg) < \min\{\overline{\alpha}_1,...,\overline{\alpha}_d\}.
\end{equation}
Ein solches $\delta$ existiert, da $\overline{\alpha}_1>0,...,\overline{\alpha}_d>0$ gelten.
Aus (\ref{BeweisSEPr7}) und (\ref{BeweisSEPr8}) folgt sofort $\varepsilon \overline{\alpha}_k + \alpha_k(\varepsilon)>0$ f"ur alle
$\varepsilon \in (0,\varepsilon_0]$ und damit, dass $u\big(\varepsilon \overline{a} + \alpha(\varepsilon)\big) \in \mathscr{U}$ gilt.
Setzen wir $\tilde{x}(\varepsilon) = x_* + \varepsilon \overline{x} + x(\varepsilon)$ und
$\tilde{u}(\varepsilon) = u\big( \varepsilon \overline{\alpha} + \alpha(\varepsilon) \big)$,
so folgt wegen (\ref{BeweisSEPr6})
\begin{equation} \label{BeweisSEPr9}
\mathscr{F}\big(\tilde{x}(\varepsilon),\tilde{u}(\varepsilon)\big) = 0.
\end{equation}
Da die Abbildungen $G_j(x)$ im Punkt $x_*$ gleichm"a"sig differenzierbar sind,
l"asst sich (vgl. Lemma \ref{LemmaRichtungsableitung}) eine Zahl $\sigma >0$ derart angeben,
dass f"ur $j=1,...,l'$ die Ungleichungen
$$G_j(x_* + \varepsilon x) \leq G_j(x_*) + \varepsilon \big( G_j'(x_*;\overline{x})+c\big)$$
f"ur alle $0 \leq \varepsilon \leq \sigma$ und alle $x \in \mathscr{X}$ mit $\|x-\overline{x}\|_{\mathscr{X}} \leq \sigma$ gelten.
Beachten wir (\ref{BeweisSEPr4}), sowie $G_j(x_*)=0$ f"ur $j=1,...,l'$ und $G_j(x_*)<0$ f"ur $j=l'+1,...,l$,
dann ergeben sich
\begin{equation} \label{BeweisSEPr10}
\left. \begin{array}{ll}
\hspace*{9mm} G_j\big(\tilde{x}(\varepsilon)\big) \leq G_j(x_*) - \varepsilon c +o(\varepsilon), & \quad j=1,...,l', \\[1mm]
\lim\limits_{\varepsilon \to 0^+} G_j\big(\tilde{x}(\varepsilon)\big) = G_j(x_*)<0, & \quad j=l'+1,...,l.
\end{array} \right\}
\end{equation}
Nach Voraussetzung (A$_1$) ist die Abbildung $x \to J(x,u)$ f"ur jedes $u \in \mathscr{U}$ im Punkt $x_*$ Fr\'echet-differenzierbar.
Daher gibt es f"ur alle $x$ eine Abbildung $o(\varepsilon)$ mit
$$J(x_* + \varepsilon x, u_*) = J(x_*,u_*) + \varepsilon J_x(x_*,u_*) \overline{x} + o(\varepsilon)
                                 + \varepsilon J_x(x_*,u_*) (x-\overline{x}).$$
Daher lassen sich $\sigma >0$ und $\varepsilon_0>0$ so w"ahlen,
dass mit $c>0$ in (\ref{BeweisSEPr2}) die Ungleichung
\begin{equation} \label{BeweisSEPr11}
J(x_* +  \varepsilon x, u_*) \leq J(x_*,u_*) + \varepsilon \big(J_x(x_*,u_*)\overline{x} + c\big)
\end{equation}
f"ur alle $\varepsilon \in [0,\varepsilon_0]$ und alle $x \in \mathscr{X}$ mit $\|x- \overline{x}\|_{\mathscr{X}} \leq \sigma$ gilt. \\
Wir spezifizieren erneut die Zahl $\delta$.
Neben $\delta C(\Lambda)< 1/2$ sei $\delta>0$ so gew"ahlt, dass
\begin{eqnarray}
&& \label{BeweisSEPr12} K \delta \bigg( \|\overline{x}\|_{\mathscr{X}} + \sum_{k=1}^d \overline{\alpha}_k \bigg)
                        < \min\{\overline{\alpha}_1,...,\overline{\alpha}_d,\sigma\}, \\
&& \label{BeweisSEPr13} K \delta \bigg( \|\overline{x}\|_{\mathscr{X}} + \sum_{k=1}^d \overline{\alpha}_k \bigg)
                        \cdot \max_{1\leq k\leq d} |J(x_*,\overline{u}_k)-J(x_*,u_*)| < c, \\
&& \label{BeweisSEPr14} (\delta + K \delta^2) \bigg( \|\overline{x}\|_{\mathscr{X}} + \sum_{k=1}^d \overline{\alpha}_k \bigg) < c.
\end{eqnarray}
Aus der zweiten Ungleichung in Voraussetzung (C) ergibt sich f"ur $\varepsilon \in [0,\varepsilon_0]$:
\begin{eqnarray}
       J\big(\tilde{x}(\varepsilon),\tilde{u}(\varepsilon)\big)
&\leq& J\big(\tilde{x}(\varepsilon),u_*\big)
       + \sum_{k=1}^d \big(\varepsilon\overline{\alpha}_k + \alpha_k(\varepsilon)\big)
                      \Big(J\big(\tilde{x}(\varepsilon),\overline{u}_k\big)-J\big(\tilde{x}(\varepsilon),u_*\big)\Big) \nonumber \\
&    & \label{BeweisSEPr15} + \delta \bigg( \|\varepsilon \overline{x} + x(\varepsilon)\|_{\mathscr{X}}
                               + \sum_{k=1}^d \big(\varepsilon\overline{\alpha}_k + \alpha_k(\varepsilon)\big) \bigg).
\end{eqnarray}
Wir zeigen, dass aus dieser Ungleichung die Relation
\begin{equation} \label{BeweisSEPr16}
J\big(\tilde{x}(\varepsilon),\tilde{u}(\varepsilon)\big) \leq J(x_*,u_*) - \varepsilon c + o(\varepsilon)
\end{equation}
f"ur alle hinreichend kleine $\varepsilon >0$ folgt.
Wegen (\ref{BeweisSEPr7}) und (\ref{BeweisSEPr12}) ist $\varepsilon^{-1} \|x(\varepsilon)\|_{\mathscr{X}} \leq \sigma$.
Daher ist $J(\tilde{x}(\varepsilon), u_*) \leq J(x_*,u_*) + \varepsilon \big(J_x(x_*,u_*)\overline{x} + c\big)$ nach (\ref{BeweisSEPr11})
f"ur $0\leq \varepsilon \leq \varepsilon_0$ erf"ullt.
Nach Voraussetzung (A$_1$) und $\alpha_k(\varepsilon) \to 0$ f"ur $\varepsilon \to 0$ gilt
\begin{eqnarray*}
\lefteqn{\sum_{k=1}^d \big(\varepsilon\overline{\alpha}_k + \alpha_k(\varepsilon)\big)
   \Big(J\big(\tilde{x}(\varepsilon),\overline{u}_k\big)-J(x_*,\overline{u}_k)-J\big(\tilde{x}(\varepsilon),u_*\big) + J(x_*,u_*)\Big)} \\
&=& \sum_{k=1}^d \big(\varepsilon\overline{\alpha}_k + \alpha_k(\varepsilon)\big)
    \Big(J_x(x_*,\overline{u}_k)(\varepsilon \overline{x})-J_x\big(x_*,u_*)(\varepsilon \overline{x}) + o(\varepsilon)\Big) = o(\varepsilon).
\end{eqnarray*}
Beachten wir ferner die Beziehungen (\ref{BeweisSEPr7}), (\ref{BeweisSEPr12}) und (\ref{BeweisSEPr13}),
so erhalten wir
\begin{eqnarray*}
\lefteqn{\sum_{k=1}^d \big(\varepsilon\overline{\alpha}_k + \alpha_k(\varepsilon)\big)
   \Big(J\big(\tilde{x}(\varepsilon),\overline{u}_k\big)-J\big(\tilde{x}(\varepsilon),u_*\big)\Big)} \\
&=&  \varepsilon \sum_{k=1}^d \overline{\alpha}_k \big(J(x_*,\overline{u}_k)-J(x_*,u_*)\big)
   + \sum_{k=1}^d \alpha_k(\varepsilon) \big(J(x_*,\overline{u}_k)-J(x_*,u_*)\big) \\
&& + \sum_{k=1}^d \big(\varepsilon\overline{\alpha}_k + \alpha_k(\varepsilon)\big)
   \Big(J\big(\tilde{x}(\varepsilon),\overline{u}_k\big)-J(x_*,\overline{u}_k)-J\big(\tilde{x}(\varepsilon),u_*\big) + J(x_*,u_*)\Big) \\
&\leq&  \varepsilon \bigg( \sum_{k=1}^d \overline{\alpha}_k \big(J(x_*,\overline{u}_k)-J(x_*,u_*)\big) + c \bigg) + o(\varepsilon).
\end{eqnarray*}
Schlie"slich gilt gem"a"s (\ref{BeweisSEPr7}) und (\ref{BeweisSEPr14})
\begin{eqnarray*}
\delta \bigg( \|\varepsilon \overline{x} + x(\varepsilon)\|_{\mathscr{X}}
                               + \sum_{k=1}^d \big(\varepsilon\overline{\alpha}_k + \alpha_k(\varepsilon)\big) \bigg)
\!\!\!&\leq&\!\!\! \varepsilon \delta \bigg(\|\overline{x}\|_{\mathscr{X}} + \frac{\|x(\varepsilon)\|_{\mathscr{X}}}{\varepsilon}
                                +\sum_{k=1}^d \Big( \overline{\alpha}_k + \frac{|\alpha_k(\varepsilon)|}{\varepsilon}\Big) \bigg) \\
\!\!\!&\leq&\!\!\! \varepsilon (\delta + K \delta^2) \bigg( \|\overline{x}\|_{\mathscr{X}} + \sum_{k=1}^d \overline{\alpha}_k \bigg) < \varepsilon c.
\end{eqnarray*}
Zusammen ergibt sich aus diesen Ungleichungen
$$J\big(\tilde{x}(\varepsilon),\tilde{u}(\varepsilon)\big)
 \leq J(x_*,u_*) + \varepsilon \bigg( J_x(x_*,u_*)\overline{x} + \sum_{k=1}^d \overline{\alpha}_k \big(J(x_*,\overline{u}_k) - J(x_*,u_*)\big) +3c \bigg)
       + o(\varepsilon).$$
Beachten wir darin (\ref{BeweisSEPr2}), dann erhalten wir unmittelbar (\ref{BeweisSEPr16}). \\[2mm]
Die Beziehungen (\ref{BeweisSEPr9}) und (\ref{BeweisSEPr10}) zeigen,
dass das Paar $\big(\tilde{x}(\varepsilon),\tilde{u}(\varepsilon)\big)$ f"ur hinreichend kleine $\varepsilon >0$ zul"assig in der
Aufgabe (\ref{EAAnhang3}) ist.
Weiterhin folgt $J\big(\tilde{x}(\varepsilon),\tilde{u}(\varepsilon)\big) < J(x_*,u_*)$ aus (\ref{BeweisSEPr16}).
Wegen $\tilde{x}(\varepsilon) \to x_*$ f"ur $\varepsilon \to 0$ bedeutet dies,
dass der Punkt $(x_*,u_*)$ im Widerspruch zur Voraussetzung kein starkes lokales Minimum sein k"onnte. \hfill $\blacksquare$

%% file: Manuskript.bbl
\begin{thebibliography}{...}
\bibitem{Arnold} Arnold,\,L.: Wachstumstheorie. Franz Vahlen, München (1997).              
\bibitem{AseKry} Aseev,\,S.M., Kryazhimskii,\,A.V.: The Pontryagin Maximum Principle and Optimal Economic Growth Problems.
                 Proc. Steklov Inst. Math., 257, 1--255 (2007).
\bibitem{AseVel} Aseev,\,S.M., Veliov,\,V.M.:
                 Maximum Principle for Infinite-Horizon Optimal Control Problems with Dominating Discount.
                 Dynamics of Continuous, Discrete and Impulsive Systems, Series B, Vol. 19 1-2 (2012).
\bibitem{AseVel2} Aseev,\,S.M., Veliov,\,V.M.: Needle variations in infinite-horizon optimal control.
                  Variational and Optimal Control Problems on Unbounded Domains,
                  eds. G.\,Wolansky, A.J.\,Zaslavski, Amer. Math. Soc. Contemporary Mathematics, 619 (2014).
\bibitem{AseVel3} Aseev,\,S.M., Veliov,\,V.M.:
                  Maximum principle for infinite-horizon optimal control problems under weak regularity
                  assumptions. Trudy Inst. Mat. i Mekh. UrO RAN, 20, no.3, 41--57 (2014).
\bibitem{Bakke} Bakke,\,V.L.:
                A maximum principle for an optimal control problem with integral constraints, JOTA, 13, 32--55 (1974).
\bibitem{Barro} Barro,\,R.J., Sala-i-Martin,\,X.: Economic growth. McGraw Hill, New York (1995).
\bibitem{Bauer} Bauer,\,H.: Ma"s- und Integrationstheorie. 2.\,Auflage. Walter de Gruyter \& Co., Berlin-New York (1992).
\bibitem{Bellman} Bellman,\,R.E.: Dynamic Programming. Princeton University Press (1957).
\bibitem{Bock} Bock,\,H.G., Longman,\,R.W.:
               Computation of Optimal Controls on Disjoint Control Sets for Minimum Energy Subway Operations. 
               Proc. AAS Symp. on Engineering Science and Mechanics, Taiwan (1981). 
\bibitem{Boltjanski} Boltjanski,\,W.G.: Mathematische Methoden der Optimalen Steuerung.
                     Akademische Verlagsgesellschaft Geest \& Portig K.-G., Leipzig (1971). 
\bibitem{Bonnans} J. F. Bonnans,\,J.\,F.,  De la Vega,\,C.: Optimal control of state constrained integral equations.
                  Set-Valued Var. Analysis 18, 307--326 (2010).
\bibitem{Brodskii} Brodskii,\,Yu.I.:
                   Necessary Conditions for a Weak Extremum in Optimal Control Problems on an Infinite Time Interval.
                   Mat. Sb., Vol. 105(147), no. 3, 371--388 (1978).
\bibitem{Carlson} Carlson,\,D.A.: 
                  An elementary proof of the maximum principle for optimal control problems governed by a Volterra integral
                  equation. JOTA, 54, 43--61 (1987).
\bibitem{CarHauLei} Carlson,\,D.A., Haurie,\,A.B., Leizarowitz,\,A.: Infinite Horizon Optimal Control.
                    Springer-Verlag, New York, Berlin, Heidelberg (1991).
\bibitem{Clarke} Clarke,\,F.: Optimization and Nonsmooth Analysis. John Wiley \& Sons, New York (1983).
\bibitem{ClarkeVinter} Clarke,\,F.H., Vinter,\,R.B.: Optimal multiprocesses.
                       SIAM J. Control and Optimization, 27, pp. 1072--1091 (1989).
\bibitem{ClarkeVinterII} Clarke,\,F.H., Vinter,\,R.B.: Applications of Optimal Multiprocesses.
                         SIAM J. Control and Optimization, 27, pp. 1048--1071 (1989).
\bibitem{Dmitruk} Dmitruk,\,A.V., Kaganovich,\,A.M.:
                  The Hybrid Maximum Principle is a consequence of Pontryagin Maximum Principle.
                  Systems \& Control Letters, Vol. 57, no. 11, p. 964--970 (2008).
\bibitem{DmitrukOsmo1} Dmitruk,\,A.V., Osmolovskii,\,N.P.:
                       Necessary conditions for a weak minimum in optimal control problems with integral equations
                       subject to state and mixed constraints.
                       SIAM J. on Control and Optimization, 52, 3437--3462 (2014).
\bibitem{DmitrukOsmo2} Dmitruk,\,A.V., Osmolovskii,\,N.P.:
                       Necessary conditions for a weak minimum in optimal control problems with integral equations on a
                       variable time interval.
                       Discret. Continuous Dyn. Systems. 35(9), 4323--4343 (2015).
\bibitem{Dockner} Dockner\,E., Feichtinger,\,G., Mehlmann,\,A.:
                  Noncooperative Solutions for a Differential Game Model of Fishery.
                  Journal of Economic Dynamics and Control 13, 1--20 (1989).
\bibitem{Dorfman} Dorfman,\,R.: An economic interpretation of optimal control theory.
                  AER 59, 817--831 (1969).
\bibitem{DuboMil} Dubovitskii,\,A.Ya., Milyutin,\,A.A.: Extremum problems in the presence of restrictions.
                  USSR Computational Mathematics and Mathematical Physics, Volume 5, Issue 3, 1--80 (1965).
\bibitem{Elstrodt} Elstrodt,\,J.: Ma"s- und Integrationstheorie. Springer, Berlin, (1996).
\bibitem{Feichtinger} Feichtinger,\,G., Hartl,\,R.F.: Optimale Kontrolle "okonomischer Prozesse.
                      de Gruyter; Berlin - New York (1986).
\bibitem{Filippov} Filippov,\,A.F.: Differential Equations with Discontinuous Right-Hand Sides.
                   Nauka, Moscow, (1985), (Kluwer, Dordrecht, 1988).
\bibitem{Galbraith} Galbraith,\,G.N., Vinter,\,R.B.: Optimal control of hybrid systems with an infinite set of discrete states.
                    Journal of Dynamical and Control Systems, 9, pp. 563--584 (2003).
\bibitem{Gale} Gale,\,D.: On Optimal Development in a Multi-Sector Economy. RES 34, 1--18 (1967).
\bibitem{Gamkrelidze} Gamkrelidze,\,R.V.: Theorie zeitoptimaler Prozesses in linearen Systemen.
                      Bull. Acad. Sci. URSS, Ser. Math. 22 (1958). 
\bibitem{Garavello} Garavello,\,M., Piccoli,\,B.: Hybrid Necessary Principle.
                    SIAM J. Control Optim., 43, pp. 1867--1887 (2005).
\bibitem{Girsanov} Girsanov,\,I.V.: Lectures on Mathematical Theory of Extremum Problems.
                Volume 67 of Lecture Notes in Economics and Mathematical Systems. Springer, Berlin-Heidelberg-New York, (1972).
\bibitem{GoMa} Göllmann,\,L., Maurer,\,H.: Optimal control problems with time delays: Two case studies in biomedicine.
                    Mathematical Biosciences \& Engineering 15(5), 1137--1154 (2018).
\bibitem{GoMa2} Göllmann,\,L., Kern,\,D., Maurer,\,H.: Optimal control problems with delays in state and
                control and mixed control-state constraints, Optimal Control Applications and Methods 30, 341--365 (2009).
\bibitem{Grass} Grass,\,D., Caulkins,\,J.P., Feichtinger,\,G., Tragler,\,G., Behrens,\,D.A.:
                     Optimal Control of Ninlinear Processes: With Applications in Drugs, Corruption and Terror.
                     Springer, Berlin, (2008).
\bibitem{Halkin2} Halkin,\,H.: Nonlinear convex programming in an infinite dimensional space.
                  Mathematical theory of control, Academic Press, New York, 10--25 (1967).
\bibitem{Halkin} Halkin,\,H.: Necessary conditions for optimal control problems with infinite horizons.
                 Econometrica 42, 267--272 (1974).
\bibitem{HalkinNeustadt} Halkin,\,H., Neustadt,\,L.W.: General necessary conditions for optimization problems.
                    Proceedings of the National Academy of Sciences of the United States of America, 56 (4), 1066--1071 (1966).
\bibitem{Hartl} Hartl,\,R.F., Sethi,\,S.P., Vickson,\,R.G.: A survey of the maximum principles for optimal control problems
                with state constraints. SIAM Review 37 (2), 181--218 (1995).
\bibitem{Hestenes} Hestenes,\,M.: Calculus of variations and the optimal control theory.
                   Wiley, New-York-London (1966).
\bibitem{HeuserI} Heuser,\,H.: Lehrbuch der Analysis, Teil 1. 11.\,Auflage. Teubner Stuttgart 1994. 
\bibitem{HeuserII} Heuser,\,H.: Lehrbuch der Analysis, Teil 2. 11.\,Auflage. Teubner Stuttgart-Leipzig-Wiesbaden 2000. 
\bibitem{HeuserGD} Heuser,\,H.: Gew"ohnliche Differentialgleichungen. 3.\,Auflage. Teubner Stuttgart 1995.
\bibitem{HeuserFA} Heuser,\,H.: Funktionalanalysis. 3.\,Auflage. B.\,G. Teubner Stuttgart 1992.
\bibitem{Hritonenko} Hritonenko,\,N., Yatsenko,\,Y.: Mathematical Modeling in Economics, Ecology, and the Environment. Kluwer Dordrecht (1999)
\bibitem{Isaacs0} Isaacs,\,R.P.: Games of Pursuit. Paper No. P-257, RAND Corporation, Santa Monica (1951).
\bibitem{Isaacs} Isaacs,\,R.P.: Differential Games. John Wiley \& Sons, New York (1965).
\bibitem{Ioffe} Ioffe,\,A.D., Tichomirov,\,V.M.: Theorie der Extremalaufgaben.
                VEB Deutscher Verlag der Wissenschaften Berlin, (1979).
\bibitem{Kamien} Kamien,\,M.I., Schwartz,\,N.L.: Dynamic Optimization:
                 The Calculus of Variations and Optimal Control in Economics and Management.
                 North-Holland, New York (1981).
\bibitem{Kurcyusz} Kurcyusz,\,S.: On the Existence and Nonexistence of Lagrange Multipliers in Banach Spaces.
                   Journal of Optimization Theory and Applications, 20 (1), 81--110 (1976).
\bibitem{Lancaster} Lancaster,\,K.: The dynamic inefficiency of capitalism. Journal of Political Economy 81, 1092--1109 (1973).
\bibitem{Ljusternik} Ljusternik,\,L.A.: On Conditional Extremums of Functionals. Mat. Sb. 41, 390--401 (1934).
\bibitem{Mangasarian} Mangasarian,\,O.L.: Sufficient Conditions for the Optimal Control of Nonlinear Systems.
                      SIAM J. Control 4, 139--152 (1966).
\bibitem{Maurer} Maurer,\,H.: Skript zur Vorlesung Optimale Steuerungsprozesse.
                 https://www.uni-muenster.de/AMM/num/Vorlesungen/maurer/Skript{\%}20Opt{\%}20Steuer.pdf
\bibitem{McShane} McShane,\,E.J.: On Multipliers for Lagrange Problems. American Journal of Mathematics 61, 809--819 (1939).
\bibitem{Michel} Michel,\,P.:
                 On the Transversality Condition in Infinite Horizon Optimal Problems. Econometrica 50, 975--985 (1982).
\bibitem{Natanson} Natanson,\,I.P.: Theorie der Funktionen einer reellen Ver\"anderlichen.
                   Verlag Harri Deutsch, Thun-Frankfurt/Main, (1981).
\bibitem{Neustadt} Neustadt,\,L.W.: Optimization: A Theory of Necessary Conditions. Princeton, New Jersey (1976).
\bibitem{PeschPlail} Pesch,\,H.J., Plail,\,M.: The Maximum Principle of optimal control:
                     A history of ingenious ideas and missed opportunities. Control and Cybernetics 38 (2009).
\bibitem{Plail} Plail,\,M.: Die Entwicklung der optimalen Steuerungen. Vandenhoeck \& Ruprecht, G"ottingen (1998).
\bibitem{Pontrjagin} Pontrjagin,\,L.S., Boltjanskij,\,V.G., Gamkrelidze,\,R.V., Mis\v{c}enko,\,E.F.:
                     Mathematische Theorie optimaler Prozesse. R.\,Oldenbourg, M\"unchen-Wien, (1964).
\bibitem{Ramsey} Ramsey\,F.P.: A Mathematical Theory of Saving. Econ. J. 38, 543--559 (1928).
\bibitem{Rihan} Rihan,\,F., Abdelrahman,\,D.H., Al-Maskari,\,F., Ibrahim,\,F.,Abdeen,\,M.A.:
                Delay differential model for tumour-immune-response with chemoimmunotherapy and optimal control.
                Computational and Mathematical Methods in Medicine, Hindawi Publishing Corporation,
                Vol. 2014, Article ID 982978, (2014).
\bibitem{Rockafellar} Rockafellar,\,R.T.: Convex Analysis. Princeton University Press, Princeton, New Jersey (1970).
\bibitem{Rudin} Rudin,\,W.: Real and Complex Analysis. McGraw Hill, 3rd ed., (1987).
\bibitem{Seierstad} Seierstad,\,A., Syds\ae ter,\,K.: Optimal Control Theory with Economic Applications.
                    North-Holland Amsterdam-New York-Oxford-Tokyo, (1987).
\bibitem{Sethi} Sethi,\,S.P., Thompson,\,G.L.: Optimal Control Theory. Applications to Management Science and Economics.
                Kluwer, Boston-Dordrecht-London, 2nd ed. (2000).
\bibitem{Shaikh04} Shaikh,\,M.S.: Optimal control of hybrid systems: Theory and algorithms. Ph.D. dissertation,
                   Department of Electrical and Computer Engineering, McGill University, Montreal, Canada, May 2004.
\bibitem{Shaikh07} Shaikh,\,M.S., Caines,\,P.\,E.: On the hybrid optimal control problem: Theory and algorithms.
                   IEEE Transactions on Automatic Control, 52, pp. 1587--1603, 2007.
\bibitem{Sussmann} Sussmann,\,H.J.: A nonsmooth hybrid maximum principle. Stability and Stabilization of Nonlinear Systems.
                   Lecture Notes in Control and Information Sciences 246, Springer, pp. 325--354, London 1999.
\bibitem{TauchnitzDiss} Tauchnitz,\,N.: Das Pontrjaginsche Maximumprinzip f"ur eine Klasse hybrider Steuerungsprobleme mit
                        Zustandsbeschr"ankungen und seine Anwendung.
                        PhD Thesis, Cottbus (2010).
\bibitem{Vinokurov} Vinokurov,\;V.R.: Optimal control of processes described by integral equations.
                    SIAM J. on Control, 7, 324--355 (1969).
\bibitem{Vinter} Vinter,\,R.B.: Optimal Control. Birkhauser (2000).
\bibitem{VonWeiz} Von Weizsaecker,\,C.C.: Existence of optimal programs of accumulation for an infinite time horizon.
                  Review of Economic Studies 32, 85--92 (1965).
\bibitem{Werner} Werner,\,D.: Funktionalanalysis. Springer-Verlag Berlin-Heidelberg-New York, (1995).
\bibitem{Young} Young,\,L.C.: Generalized curves and the existence of an attained absolute minimum in the calculus
                of variations.
                C. R. Acad. Soc. Sci. et Lettr., Varsovie 3:30, 212--234 (1937).
\end{thebibliography}
